\documentclass[12pt]{article}

\usepackage[utf8]{inputenc}
\usepackage{natbib}
\RequirePackage{amsmath,amssymb}
\usepackage{epsfig}
\usepackage{graphicx}
\usepackage{epstopdf,color}
\usepackage[usenames,dvipsnames,svgnames]{xcolor}
\usepackage{multirow,enumerate}
\usepackage{listings}
\usepackage{caption}
\RequirePackage[colorlinks,citecolor=blue,urlcolor=blue]{hyperref}
\usepackage{float,threeparttable}
\usepackage{subcaption}

\usepackage{ulem}

\allowdisplaybreaks[4]
\usepackage{xr}
\usepackage[a4paper]{geometry}

\newcommand{\blind}{0}

\newtheorem{thm}{Theorem}
\newtheorem{lem}{Lemma}

\newtheorem{example}{Example}

\def\ba{{\boldsymbol a}}

\def\bb{{\boldsymbol b}}

\def\bE{{\mathbb E}}

\def\cF{{\cal F}}

\def\cG{{\cal G}}

\def\Gn{{\mathbb G}_n}

\def\bH{\boldsymbol{H}}

\def\tH{\tilde{\bH}}

\def\bI{{\boldsymbol I}}
\def\bJ{{\boldsymbol J}}
\def\hJ{{\hat J}}
\def\hK{{\hat K}}
\def\tK{{\tilde K}}

\def\rP{\boldsymbol{\mathrm P}}

\def\Pn{\mathbb{P}_{n}}

\def\bR{\mathbb{R}}

\def\tT{\tilde{T}}

\def\bv{{\boldsymbol v}}
\def\bV{{\boldsymbol V}}

\def\tV{{\tilde{V}}}

\def\bW{{\boldsymbol W}}

\def\bx{{\boldsymbol{x}}}

\def\bX{{\boldsymbol X}}
\def\tbX{\tilde{{\boldsymbol X}}}

\def\tX{\widetilde{X}}

\def\bZ{{\boldsymbol Z}}

\def\eps{\epsilon}

\def\real{\mathop{{\rm I}\kern-.2em\hbox{\rm R}}\nolimits}

\def\trans{^{\rm T}}

\def\htheta{\hat{\theta}}

\def\trace{\mbox{trace}}

\def\balpha{\boldsymbol \alpha}

\def\tbalpha{\tilde{\boldsymbol \alpha}}
\def\hbalpha{\hat{\boldsymbol \alpha}}
\def\btheta{\boldsymbol \theta}
\def\bbeta{\boldsymbol \beta}

\def\tbbeta{\tilde{\boldsymbol \beta}}

\def\bmu{\boldsymbol \mu}
\def\tmu{\tilde{\mu} }

\def\tpsi{\tilde{\psi}}

\def\hbtheta{\hat{\boldsymbol \theta}}

\def\bxi{\boldsymbol \xi}

\def\bzeta{\boldsymbol{\zeta}}

\def\bvarsigma{\boldsymbol \varsigma}

\def\Var{\hbox{Var}}
\def\Cov{\hbox{Cov}}

\def\bone{{\boldsymbol 1}}
\def\bzero{{\boldsymbol 0}}

\newtheorem{lemma}{Lemma}
\newtheorem{remark}{Remark}
\newtheorem{theorem}{Theorem}

\newcommand{\E}{\textrm{E}}

\newcommand{\pkonv}{\stackrel{p}{\rightarrow}}
\newcommand{\lpkonv}{\stackrel{p}{\longrightarrow}}

\newcommand{\askonv}{\stackrel{a.s.}{\longrightarrow}}
\newcommand{\lkonv}{\stackrel{\cal{L}}{\longrightarrow}}

\def\pvalue{\mbox {p-value}}

\def\Gp{\mbox {grouping parameter}}
\def\Gv{\mbox {grouping variable}}
\def\GDp{\mbox {grouping difference parameter}}

\begin{document}
\def\spacingset#1{\renewcommand{\baselinestretch}%
{#1}\small\normalsize} \spacingset{1}
\captionsetup[figure]{name={Fig.},labelsep=period}

\date{}

\if0\blind{
	\title{\bf Efficient subgroup testing in change-plane models}
	\author{
		Xu Liu \\
		School of Statistics and Management,\\
		Shanghai University of Finance and Economics, Shanghai 200433, China\\
		Jian Huang \\
		Department of Applied Mathematics,\\
		The Hong Kong Polytechnic University, Hung Hom, Kowloon, Hong Kong, China\\
		Yong Zhou \\
		Academy of Statistics and Interdisciplinary Sciences,\\
		East China Normal University, Shanghai 200062, China\\
		Feipeng Zhang\\
		School of Economics and Finance, \\
		Xi'an Jiaotong University, Xi'an 710049, China\\
		Panpan Ren\thanks{Corresponding author. Email address: panpanren@stu.sufe.edu.cn} \\
		School of Statistics and Management,\\
		Shanghai University of Finance and Economics, Shanghai 200433, China\\
	}
	\maketitle
} \fi

\if1\blind{
	\bigskip
	\bigskip
	\bigskip
	\begin{center}
		{\LARGE\bf Efficient subgroup testing in change-plane models}
	\end{center}
	\medskip
} \fi

\bigskip
\begin{abstract}
	Considered here is a hypothesis test for the coefficients in the change-plane regression models to detect the existence of a change plane. 
  The test that is considered is from the class of test problems in which some parameters are not identifiable under the null hypothesis. The classic exponential average tests do not work well in practice. To overcome this drawback, a novel test statistic is proposed by taking the weighted average of the squared score test statistic (WAST) over the grouping parameter's space, which has a closed form from the perspective of the conjugate priors when an appropriate weight is chosen. The WAST significantly improves the power in practice, particularly in cases where the number of the grouping variables is large. The asymptotic distributions of the WAST are derived under the null and alternative hypotheses. The approximation of critical value by the bootstrap method is investigated, which is theoretically guaranteed. Furthermore, the proposed test method is naturally extended to the generalized estimating equation (GEE) framework, as well as multiple change planes that can test if there are three or more subgroups. The WAST performs well in simulation studies, and its performance is further validated by applying it to real datasets.
\end{abstract}

\noindent%
{\it Keywords:}  Nonstandard tests; Quantile regression; Subgroup detection; U-statistics
\vfill

\spacingset{1.45} 
\section{Introduction}
\label{sec:intro}

Subgroup analysis has emerged as a powerful tool in diverse fields covering public management, clinical medicine and the physical and social sciences. 
In public management, the effect of a certain policy may be heterogeneous, which means that the same policy can have distinct impacts on different subpopulations. 
Ignoring a population's heterogeneity in classical models may cause biased estimates, and identifying subgroups of individuals is the core problem. 
Procedures for estimating and identifying subgroups have been extensively studied based on various models \citep{2000Assmann,hansen2000sample,2011Subgroup, 2012Estimating, 2013Effectively, 2018The, 2018Oracle, 2020Threshold, 2021Multithreshold, 2021Subgroup, 2022Single, 2022Inferece}. Before identifying subgroups, an appropriate statistic is required to test for their existence. 
However, this test problem has received relatively little attention, especially when there are large number of grouping variables 
that can divide the population into subgroups with different associations between predictors and responses. 
The existing literature for the testing problem is outlined below, focusing on several classical change-plane regression models. Therefore, herein, we delve into subgroup testing in the change-plane regression models, which belongs to the class of hypothesis tests when some parameters are not identifiable under the null hypothesis.

To motivate our work, we look at the data from the 2013 China Health and Retirement Longitudinal Study (CHARLS). It is widely acknowledged that population aging is a significant global social issue, and China has been facing this challenge since 2001.  In response, the Chinese government has established a social old-age insurance system to alleviate the burden on adult children. 
As is known, child support and social pension are two major sources of income for rural elderly in China. 
Does the pension income from the New Rural Social Pension Insurance (NRSPI) effectively alleviate the financial burden on adult children of the elderly? 
Many scholars have offered different answers to this question \citep{zhang2019beneficiaries,nikolov2019private,ko2021chipping}. 
Addressing this question can help find out the shortcomings of the current pension insurance system and make improvements. 
Therefore, it is crucial to identify whether there is a subpopulation in which pension income has a distinct impact on support income of the elderly. 
This paper considers a dataset from CHARLS in Shanxi Province. 
The support income of the elderly in Shanxi Province, which is considered as the response variable, displays heavy-tailed and asymmetry characteristics (see Figure D.1 in the Supplementary Material), prompting us to develop a quantile change-plane regression model. 

For the different types of data and goals, there are many change-plane regressions considered in the literature.
Here, we give some common and useful examples.
\begin{example}\label{exam:qr}{\bf Quantile regression. }
	\citet{2022Single} and \citet{2018Oracle,2011Testing} considered single-index thresholding in quantile change-plane regression:
	\begin{align}\label{model:qr}
		Q_{Y_i}(\tau|\tbX_i, \bX_i, \bZ_i) = \tbX_i\trans\balpha(\tau)+\bX_i\trans\bbeta(\tau)\bone(\bZ_i\trans\btheta(\tau)\geq0),
	\end{align}
	where $Q_Y(\tau|\cdot)$ denotes the $\tau$th conditional quantile of the response $Y$ given covariates $\tbX$, $\bX$, and $\bZ$, $\tau\in(0, 1)$ is the quantile level of interest, $\bX_i$ is a subset of $\tbX_i$, and $\bone(\cdot)$ denotes the indicator function.
	For fixed $\btheta$, the score function with respect to $\bbeta$ at $\bbeta=\bzero$ is
	\begin{align*}
		\psi(\bV_i, \hbalpha, \bzero, \btheta)=
		\bone(\bZ_i\trans\btheta(\tau)\geq0)\left[\bone(Y_i-\tbX_i\trans\hbalpha(\tau)\leq 0)-\tau\right]\bX_i,
	\end{align*}
	where $\hbalpha=\hbalpha(\tau)$ is the estimate of $\balpha=\balpha(\tau)$ when $\bbeta=\bzero$ according to
	\begin{align*}
		\sum_{i=1}^{n}\psi_1(\bV_i, \balpha) = \sum_{i=1}^{n}\left[\bone(Y_i-\tbX_i\trans\balpha(\tau) \leq 0)-\tau\right]\tbX_i=\bzero.
	\end{align*}

\end{example}

\begin{example}\label{exam:glm} {\bf Logistic regression model.} \citet{2021Threshold} studied the logistic change-plane regression model, which belongs to generalized linear models (GLMs) with the likelihood
	\begin{align}\label{model:glm}
		\begin{split}
			\exp\left[ Y_i\mu_i - c(\mu_i)\right]a(Y_i), ~\mbox{with}~
			\mu_i =& \tbX_i\trans\balpha + \bX_i\trans\bbeta\bone(\bZ_i\trans\btheta\geq0),
		\end{split}
	\end{align}
	where $\mu_i$ is the canonical parameter, and $a(\cdot)$ and $c(\cdot)$ are known functions. 
For fixed $\btheta$, the score function with respect to $\bbeta$ at $\bbeta=\bzero$ is
	\begin{align*}
		\psi(\bV_i, \hbalpha, \bzero, \btheta)= \bone(Z_i\geq\theta)\left[Y_i-\partial c(\tbX_i\trans\hbalpha)/\partial\mu_i\right]\bX_i,
	\end{align*}
	where $\hbalpha$ is the estimate of $\balpha$ when $\bbeta=\bzero$ according to
	\begin{align*}
		\sum_{i=1}^{n}\psi_1(\bV_i, \balpha) = \sum_{i=1}^{n}\left[Y_i-\partial c(\tbX_i\trans\balpha)/\partial\mu_i\right]\tbX_i=\bzero.
	\end{align*}
\end{example}

\begin{example}\label{exam:prob}{\bf Probit regression model.}
	\citet{2011Testing} considered the probit change-plane regression model with likelihood
	\begin{align}\label{model:prob}
		\Phi(\tbX_i\trans\balpha+\bX_i\trans\bbeta\bone(Z_i\geq\theta))^{Y_i}\Phi(-\tbX_i\trans\balpha-\bX_i\trans\bbeta\bone(Z_i\geq\theta))^{1-Y_i},
	\end{align}
	where $\Phi(\cdot)$ is the cumulative distribution function of the standard normal distribution. 
For fixed $\btheta$, the score function with respect to $\bbeta$ at $\bbeta=\bzero$ is
	\begin{align*}
		\psi(\bV_i, \hbalpha, \bzero, \btheta)=
		\bone(Z_i\geq\theta)\left[Y_i\frac{\phi(\tbX_i\trans\hbalpha)}{\Phi(\tbX_i\trans\hbalpha)} - (1-Y_i)\frac{\phi(-\tbX_i\trans\hbalpha)}{\Phi(-\tbX_i\trans\hbalpha)}\right]\bX_i,
	\end{align*}
	where $\phi(\cdot)$ is the density of the standard normal distribution, and $\hbalpha$ is the estimate of $\balpha$ when $\bbeta=\bzero$ according to
	\begin{align*}
		\sum_{i=1}^{n}\psi_1(\bV_i, \balpha) = \sum_{i=1}^{n}\left[Y_i\frac{\phi(\tbX_i\trans\balpha)}{\Phi(\tbX_i\trans\balpha)} - (1-Y_i)\frac{\phi(-\tbX_i\trans\balpha)}{\Phi(-\tbX_i\trans\balpha)}\right]\tbX_i=\bzero.
	\end{align*}
	
\end{example}

\begin{example}\label{exam:semi} {\bf Semiparametric model.}
	\citet{2017Change} considered the semiparametric change-plane model
	\begin{align}\label{model:semi}
		Y_i = \gamma(\bZ_i) +\beta X_i \bone(\bZ_i\trans \btheta\geq0)+\eps_i,
	\end{align}
	where $\gamma(\bZ)$ is an unknown baseline mean function for patients in treatment $X=0$, $X_i\in\{0, 1\}$, and $\bE(\eps_i|X_i, \bZ_i)$=0. 
For fixed $\btheta$, the score function with respect to $\bbeta$ at $\bbeta=\bzero$ is
	\begin{align*}
		\psi(\bV_i, \hbalpha, \bzero, \btheta)=
		\bone(\bZ_i\trans \btheta\geq0)\left[X_i-\pi(\bZ_i, \hbalpha_1)\right]\left[Y_i-\gamma(\bZ_i, \hbalpha_2)\right],
	\end{align*}
	where $\hbalpha=(\hbalpha_1\trans, \hbalpha_2\trans)\trans$ is the estimate of $\balpha=(\balpha_1\trans, \balpha_2\trans)\trans$ when $\bbeta=\bzero$ according to
	\begin{align*}
		\sum_{i=1}^{n}\psi_1(\bV_i, \balpha)
		=\sum_{i=1}^{n}\left(
		\begin{array}{c}
			D_{\balpha_1}(\bZ_i)\{X_i-\pi(\bZ_i, \balpha_1)\}\\
			D_{\balpha_2}(\bZ_i)\{Y_i-\gamma(\bZ_i, \balpha_2)\}
		\end{array}
		\right)=\bzero
	\end{align*}
	with $D_{\balpha_1}(\bZ_i)=[\pi(\bZ_i, \balpha_1)(1-\pi(\bZ_i, \balpha_1))]^{-1}\partial \pi(\bZ_i, \balpha_1)/\partial \balpha_1$ and $D_{\balpha_2}(\bZ_i)=\partial \gamma(\bZ_i, \balpha_2)/\partial \balpha_2$. Following \citet{2017Change}, we emphasize that either the baseline mean function $\gamma(\bZ, \balpha_2)$ or the propensity model $\pi(\bZ_i, \balpha_1)$ is specified correctly, but not necessarily both.
\end{example}

The goal, subgroup or change-plane test, considered herein is to test whether $\bbeta$ is $\bzero$ or not, i.e.,
\begin{align}\label{test}
H_0: \bbeta=\bzero\quad \textrm{versus} \quad H_1: \bbeta\neq \bzero,
\end{align}
where the nuisance parameter $\btheta$ is only identifiable under the alternative hypothesis.
For convenience of expression, we call $\bZ$ the $\Gv$, $\btheta$ the $\Gp$, $\bX$ the grouping difference variable, and $\bbeta$ the $\GDp$. The grouping difference parameter measures the effect difference between two subgroups. $\tbX_i\trans\balpha$ is the baseline model, $\tbX$ is the baseline variable, and $\balpha$ is called the effect of the baseline group.

From examples 1--4, it is natural to extend these classical models to the generalized estimating equation (GEE) framework.
Let the observed data consist of $\{\bV_i=(Y_i,\tbX_i,\bX_i,\bZ_i), i=1, \cdots, n\}$, which are $n$ independent and identically distributed copies of $\bV=(Y,\tbX,\bX,\bZ)$.
Consider the GEEs $\bE[\bone(\bZ\trans\btheta\geq 0)\psi_0(\bV, \balpha, \bbeta,\btheta)]=\bzero$ and $\bE\psi_1(\bV, \balpha)=\bzero$ and the corresponding empirical estimating functions
\begin{align}\label{model1a}
\begin{split}
\Psi_{n}(\balpha, \bbeta, \btheta)=&\sum_{i=1}^{n} \bone(\bZ_i\trans\btheta\geq 0)\psi_0(\bV_i, \balpha, \bbeta,\btheta), \\
\Psi_{1n}(\balpha)=&\sum_{i=1}^{n} \psi_1(\bV_i, \balpha),
\end{split}
\end{align}
where $\balpha=(\alpha_1, \cdots, \alpha_r)\trans\in\Theta_{\alpha}\subseteq\mathbb{R}^{r}$, $\bbeta=(\beta_1, \cdots, \beta_p)\trans\in\Theta_{\beta}\subseteq\mathbb{R}^p$, and $\btheta=(\theta_1, \cdots, \theta_q)\trans\in\Theta_{\theta}\subseteq\mathbb{R}^q$ are unknown parameters. Here, $\psi_0(\bV, \balpha, \bzero,\btheta)$
is independent of both $\bbeta$ and $\btheta$. 
For ease of expression, we denote $\psi(\bV, \balpha, \bbeta, \btheta)=\bone(\bZ\trans\btheta\geq 0)\psi_0(\bV, \balpha, \bbeta,\btheta)$ and $\psi_0(\bV, \balpha)=\psi_0(\bV, \balpha, \bzero,\btheta)$.
In many hypothesis testing problems, $\btheta$ is a nuisance parameter, such as in the examples above. The functions~\eqref{model1a} for considering testing problem~\eqref{test} are flexible, with many existing models included as special cases, such as examples 1--4.

In the ordinary testing problem in which there is no identifiability problem under either the null or alternative hypothesis, the classic Wald- and score-type tests perform well. When nuisance parameters are present, \citet{1999Estimation} and \citet{1994Optimal, 1995Admissibility} studied the optimal tests based on the weighted average power criterion, which was originally introduced by \citet{1943Tests} when studying the likelihood ratio test under regularity conditions. \citet{1977Hypothesis, 1987Hypothesis, 2002Hypothesis} investigated the supremum of the squared score test statistic for mixture models, which is a special case of the optimal tests. \citet{2007Inference, 2009On, 2016Subgroup}, and \citet{2017Subgroup} considered nonstandard tests for semiparametric models with the application of censoring data. \citet{2020Testing} studied a score-type test to detect the existence of social network dependence with time-to-event data.

All the aforementioned methods require taking the supremum of a classic Wald- or score-type test statistic. 
The large feasible region of calculating the supremum may cause these tests to not perform well, and they are time-consuming computationally.
To address these limitations, this paper introduces a novel test statistic based on the weighted average of the squared score test statistic (WAST) over the nuisance parametric space. Unlike the aforementioned tests, the proposed WAST has a closed form from the perspective of the conjugate priors when an appropriate weight is chosen.
It is a by-product that the proposed WAST reduces the computational cost dramatically because a closed-form expression is obtained after integrating out the nuisance parameter. Extensive simulation studies with finite samples show that the WAST exhibits 
(i) type-\uppercase\expandafter{\romannumeral1} errors that are close to the nominal significance level, (ii) higher power compared to its competitors in practice.  
The limiting distributions of the test statistic are established under the null and alternative hypotheses, and the approximation of critical value by the bootstrap method is provided, as are the theoretical guarantees. 

The remainder of this paper is organized as follows. In Section~\ref{sec:test_stat}, we introduce the novel test statistic, WAST, and establish its limiting distributions under the null and alternative hypotheses. In Section~\ref{sec:examples}, we investigate the choice of the weight function based on several useful examples that are all special cases of the considered GEEs, with detailed computations included for each example. In Section~\ref{sec:simulations}, we report the results of simulation studies conducted to evaluate the finite-sample performance of the proposed test statistic and to compare its power with competitors. In Section~\ref{sec:case_studies}, we apply the proposed test method to a real data set. Finally, in Section~\ref{sec:conclusion}, we conclude with remarks and further extensions. The Assumptions are provided in the Appendix, the proofs and additional simulation studies are included in the Supplementary
and in Supplementary C.3, we further extend the proposed test procedure to the GEEs with multiple change planes. An R package named ``wast'' is available at \url{https://github.com/xliusufe/wast}.

\section{Test Statistics}\label{sec:test_stat}

We first give some useful notation for convenience of expression. For a vector $\bv\in \bR^{d}$ and a squared matrix $A=(a_{ij})\in\mathbb{R}^{d\times d}$, denote by $\|\bv\|$ the Euclidean norm of $\bv$, by $\|\bv\|_{p}=(\sum_{i=1}^{d}|v_i|^p)^{1/p}$ the $p$-norm of $\bv$; $\|\bv\|_{A}^2=\sum_{i, j}a_{ij}v_iv_j$ and $\bv^{\otimes2}=\bv\bv\trans$. Denote by $\|A\|_p=\sup\{\|A\bx\|_p: \bx\in\bR^{d}, \|\bx\|_p=1\}$ the induced operator norm for a matrix $A=(a_{ij})\in\bR^{m\times d}$. 
Denote by $\rP$ the ordinary probability measure such that $\bE f=\int f d\rP$ for any measurable function $f$. Let $L^p(Q)$ be the space of all measurable functions $f$ such that $\|f\|_{Q, p}:=(Q|f|^p)^{1/p}<\infty$, where $p\in[1, \infty)$ and $(Q|f|^p)^{1/p}$ denotes the essential supremum when $p=\infty$. Let $N(\eps, \cF, \|\cdot\|_{Q, 2})$ be an $\eps$-covering number of $\cF$ with respect to the $L^2(Q)$ seminorm $\|\cdot\|_{Q, 2}$, where $\cF$ is a class of measure functions and $Q$ is finite.

\subsection{Preliminaries}\label{sec:sst}
We are interested in the testing problem \eqref{test}.
For any known $\btheta\in\Theta_{\theta}$, we consider a squared score test statistic for testing $\bbeta=\bzero$, i.e.,
\begin{align}\label{score_test0}
	\tT_n(\btheta) = n^{-1}\|\Psi_n(\hbalpha, \bzero, \btheta)\|^2_{\tV(\btheta)^{-1}},
\end{align}
where $\hbalpha$ is an estimate of $\balpha$ based on the empirical GEE $\Psi_{1n}(\balpha)=0$ and $\tV(\btheta)=\\n^{-1}\sum_{i=1}^{n}\{\psi(\bV_i, \hbalpha, \bzero, \btheta) - \hat{K}(\btheta)\hat{J}\psi_{1}(\bV_i, \hbalpha)\}^{\otimes2}$. Here, $\hK(\btheta)$ and $\hJ$ are consistent estimators of $K(\btheta)=\partial\bE\psi(\bV, \balpha_0, \bzero, \btheta)/\partial\balpha\trans\in\mathbb{R}^{p\times r}$ and $J=\left[\partial\bE\psi_1(\bV, \balpha_0)/\partial\balpha\trans\right]^{-1}\in\mathbb{R}^{r\times r}$, respectively.

\begin{lemma}\label{lem1}
	If Assumptions~{(A1)--(A5)} in the Appendix hold, then for any fixed $\btheta\in \Theta_\theta$, $\tT_n(\btheta)$ converges in distribution to a $\chi^2$ one with $p$ degrees of freedom under $H_0$ as $n \rightarrow \infty$.
\end{lemma}

\begin{remark}
	Because $\tT_n(\btheta)$ includes an unknown parameter $\btheta$, Lemma~\ref{lem1} cannot be used directly in practice. However, it essentially demonstrates that the asymptotic distributions are the same for any nuisance parameter $\btheta\in\Theta_\theta$. Thus, the supremum and the weighted average of $\tT_n(\btheta)$ achieve correct type-\uppercase\expandafter{\romannumeral1} errors. We construct the proposed test statistic by making use of this point.
\end{remark}

Because equations~\eqref{model1a} does not depend on $\btheta$ when $\bbeta = 0$, the parameter $\btheta$ is identifiable only under the alternative hypothesis. Lemma~\ref{lem1} suggests that $\tT_n(\btheta)$ for any $\btheta\in\Theta_{\theta}$ is a test statistic that can control type-\uppercase\expandafter{\romannumeral1} errors. However, $\tT_n(\btheta)$ may suffer a loss of power for some $\btheta$, making it difficult to determine $\btheta$ in practice. To pursue the largest power theoretically, a good choice is to take the supremum of the squared score test statistic over the parametric space $\Theta_{\theta}$ (SST), i.e.,
\begin{align}\label{score_test1}
	\tT_n = \sup_{\btheta\in\Theta_{\theta}} \left\{n^{-1}\|\Psi_n(\hbalpha, \bzero, \btheta)\|^2_{\tV(\btheta)^{-1}}\right\}.
\end{align}
The test statistic $\tT_n$ has been widely studied in the literature for the case in which $\psi_0(\bV, \balpha, \bbeta,\btheta)$ is one-dimensional,
and the dimension of the $\Gv$ $\bZ$ is small; see \citet{1992Optimal, 1993Tests, 1994Optimal, 1977Hypothesis, 1987Hypothesis, 2002Hypothesis, 2009On,2017Change} and \citet{Shen2020}. We summarize their results under the null and alternative hypotheses in the GEE framework in Lemmas~A.4 and A.5 in the Supplementary Material.

\subsection{Weighted average of squared score test statistic}\label{sec:wast}
The SST considered in Section~\ref{sec:sst} may sacrifice power in practice for pursuing theoretically the largest power, especially for large $q$ and $p$. To overcome this drawback and take advantage of Lemma~\ref{lem1}, it is natural to average $\tT_n(\btheta)$ over $\Theta_{\theta}$. In practice, we can take a weighted average over $\Theta_{\theta}$.

Based on the above arguments, we propose the WAST based on $\tT_n(\btheta)$ by taking integral after removing both the inverse of covariance and cross-interaction terms:
\begin{align}\label{statistic1}
	T_n = \frac{1}{n(n-1)}\sum_{i\neq j}\rho_{ij}
\end{align}
with
\begin{align*}
	\rho_{ij} = \psi_0(\bV_i, \hbalpha)\trans \psi_0(\bV_j, \hbalpha)\int_{\btheta\in\Theta_{\theta}}\bone(\bZ_i\trans\btheta\geq 0) \bone(\bZ_j\trans\btheta\geq 0)w(\btheta)d\btheta,
\end{align*}
where $w(\btheta)$ is a weight satisfying $w(\btheta)\geq0$ for all $\btheta\in\Theta_{\theta}$ and $\int_{\btheta\in\Theta_{\theta}}w(\btheta)d\btheta=1$.

Another explanation for $w(\btheta)$ comes from a Bayesian perspective, in which the weight $w(\btheta)$ can be treated as the density of the prior of the $\Gp$. Herein, we focus mainly on testing the existence of subgroups instead of estimating the $\Gp$, so there is no requirement to calculate the posterior distribution, which usually makes implementation relatively easy.

In general, there is no closed-form expression for $\rho_{ij}$, but it usually has a simple formulation if we choose a specific weight $w(\btheta)$. 
The choice of weight can affect the computation of the test statistic since the numerical integration over $\bR^{q}$ is required. 
Fortunately, there is a closed-form expression for an appropriate weight from the perspective of the conjugate priors, which is quite attractive in practice. 
\citet{FU2019} and \citet{FU2022} have discussed a similar weight that facilitates the derivation of the closed-form expression. We discuss the calculation of $\rho_{ij}$ in detail in Section~\ref{sec:examples}.

\citet{1994Optimal} studied the average Lagrange multiplier (LM), denoted by $\mbox{ATM}_n = \int_{\btheta\in\Theta_{\theta}}\tT_n(\btheta)w(\btheta)d\btheta$, which is a special case of
\begin{align*}
	\mbox{Exp-LM}_n=(1+c)^{-p/2}\int\exp\left(\frac{c}{2(1+c)}\mbox{LM}_n(\btheta)\right)dw(\btheta),
\end{align*}
where $c$ is a scalar constant. However, the main drawback of this approach is that it is hard to find a closed-form expression for $\mbox{ATM}_n$, which may lose power in practice and incurs a expensive computational cost for large $p$. By contrast, the proposed WAST $T_n$ has a different form from that of $\mbox{ATM}_n$, and it is usually easy to achieve a closed-form representation by choosing an appropriate weight $w(\btheta)$ (see Section~\ref{sec:examples} for details).

\subsection{Asymptotic analysis}
This subsection is devoted to exploring the asymptotic properties of the proposed test statistic $T_{n}$. Before establishing its asymptotic distribution, we introduce some additional notations. Let the kernel of a U-statistic under the null hypothesis be
\begin{align}\label{kernel0}
	\begin{split}
		h(\bV_i, \bV_j)
		=& \int_{\btheta\in\Theta_{\theta}}\psi(\bV_i, \balpha_0, \bzero, \btheta)\trans \psi(\bV_j, \balpha_0, \bzero, \btheta) w(\btheta)d\btheta\\
		& - \psi_1(\bV_i, \balpha_0)\trans K_{j}- K_{i}\trans \psi_1(\bV_j, \balpha_0) + \psi_1(\bV_i, \balpha_0)\trans H\psi_1(\bV_j, \balpha_0),
	\end{split}
\end{align}
where
\begin{align*}
	H=\int_{\btheta\in\Theta_{\theta}}J\trans K(\btheta)\trans K(\btheta)J w(\btheta)d\btheta ~\mbox{and}~
	K_i = \int_{\btheta\in\Theta_{\theta}}J\trans K(\btheta)\trans\psi(\bV_i, \balpha_0, \bzero, \btheta) w(\btheta)d\btheta.
\end{align*}

\begin{theorem}\label{thm11}
	If Assumptions~{(A1)--(A5)} in the Appendix hold, then under the null hypothesis, we have
	\begin{align*}
		nT_n -\mu_0 \lkonv \zeta,
	\end{align*}
	where $\lkonv$ denotes convergence in distribution, $\mu_0=-2\bE[\psi_1(\bV_1, \balpha_0)\trans K_{1}]
	+\bE[\psi_1(\bV, \balpha_0)\trans H\psi_1(\bV, \balpha_0)]$,
	$\zeta$ is a random variable of the form $\zeta=\sum_{j=1}^{\infty}\lambda_{j}(\chi^2_{1j}-1)$, and $\chi^2_{11}, \chi^2_{12}, \cdots$ are independent $\chi^2_{1}$ variables, i.e., $\zeta$ has the characteristic function
	\begin{align*}
		\bE\left[e^{it\zeta}\right]=\prod_{j=1}^{\infty}(1-2it\lambda_{j})^{-1/2}e^{-it\lambda_{j}}.
	\end{align*}
	Here, $i=\sqrt{-1}$ standards for the imaginary unit,
	and $\{\lambda_{j}\}$ are the eigenvalues of the kernel $h(\bv_1, \bv_2)$ under $f(\bv, \balpha_0, \bzero, \btheta_0)$, i.e., they are the solutions of $\lambda_{j}g_{j}(\bv_2)=\int_{0}^{\infty}h(\bv_1, \bv_2)g_{j}(\bv_1)\\ f(\bv_1, \balpha_0, \bzero, \btheta_0)d\bv_1$ for nonzero $g_{j}$, where $f(\bv, \balpha, \bbeta, \btheta)$ is the density of $\bV$ with parameters $\balpha$, $\bbeta$, and $\btheta$.
\end{theorem}

Next, we investigate the power performance of the proposed test statistic under the local alternative.
\begin{theorem}\label{thm12}
	If Assumptions~{(A1)--(A6)} in the Appendix hold, then under the local alternative hypothesis $H_1$, that is, $\bbeta = n^{-1/2}\bxi$ with a fixed vector $\bxi\in\Theta_{\beta}$, we have
	\begin{align*}
		nT_n-\mu_0 \lkonv \zeta_1,
	\end{align*}
	where $\mu_0$ is defined in Theorem~\ref{thm11}, $\zeta_1$ is a random variable of the form $\zeta_1=\sum_{j=1}^{\infty}\lambda_{j}(\chi^2_{1j}(\mu_{aj})-1)$, and $\chi^2_{11}(\mu_{a1}), \chi^2_{12}(\mu_{a2}), \cdots$ are independent noncentral $\chi^2_{1}$ variables, i.e., $\zeta_1$ has the characteristic function
	\begin{align*}
		\bE\left[e^{it\zeta_1}\right]=\prod_{j=1}^{\infty}(1-2it\lambda_{j})^{-1/2}
		\exp\left(-it\lambda_j+\frac{it\lambda_{j}\mu_{aj}}{1-2it\lambda_j}\right).
	\end{align*}
	Here, $\{\lambda_{j}\}$ are the eigenvalues of the kernel $h(\bv_1, \bv_2)$ defined in \eqref{kernel0} under $f(\bv, \balpha_0, \bzero, \btheta_0)$, i.e., they are the solutions of $\lambda_{j}g_{j}(\bv_2)=\int_{0}^{\infty}h(\bv_1, \bv_2)g_{j}(\bv_1)f(\bv_1, \balpha_0, \bzero, \btheta_0)d\bv_1$ for nonzero $g_{j}$, and each noncentrality parameter of $\chi^2_{1j}(\mu_{aj})$ is
	\begin{align*}
		\mu_{aj} = \bE\left[\phi_j(\bV_{0})\bxi\trans\partial \log(f(\bV_{0}, \balpha_0, \bzero, \btheta_0))/\partial\bbeta\right], \quad j=1, 2, \cdots,
	\end{align*}
	where $\{\phi_j(\bv)\}$ denotes orthonormal eigenfunctions corresponding to the eigenvalues $\{\lambda_j\}$, and $\bV_0$ 
	is generated from the null distribution $f(\bv, \balpha_0, \bzero, \btheta_0)$.
\end{theorem}

Theorem~\ref{thm12} implies that the power function of $nT_n-\mu_0$ can be theoretically approximated by $F_{\zeta_1}$, where $F_{\zeta_1}$ is the cumulative distribution function of $\zeta_1$. The proof of Theorem~\ref{thm12} shows that $0< \bE h(\bV_1, \bV_2)=\sum_{j=1}^{\infty}\lambda_j\mu_{aj}+o(1)$ under the local alternative hypothesis. The additional mean $\sum_{j=1}^{\infty}\lambda_j\mu_{aj}$ under $H_1$ can be viewed as a measure of the difference between $H_0$ and $H_1$. Because the asymptotic distribution is not common, it is difficult to use in practice; instead, in Section~\ref{sec:critical_value}, we provide a bootstrap method for calculating the critical value or $\pvalue$.
Theorems~\ref{thm11}--\ref{thm12} show that the bias $\mu_0$ arises from estimating $\balpha$, which means that $\mu_0$ is zero if $\balpha$ is known. In fact, if $\balpha$ is known, then $T_n$ is exactly a U-statistic with zero mean under the null hypothesis.

We extend the proposed test method to the GEEs with multiple change planes. Due to limit of space, we provide asymptotic distributions of the proposed statistic and simulation results in Appendix C.3 of the Supplementary Material. 

\subsection{Computation of critical value}\label{sec:critical_value}
In Section \ref{sec:wast}, we deduced the asymptotic distributions for the WAST $T_n$ under the null and alternative hypotheses. However, because these asymptotic distributions are not as common as normal, $\chi^2$, and F distributions, it is difficult to calculate the critical value directly. Instead, many related articles recommend the strategy of approximating the limiting null distribution of the test statistic. \citet{2017Change} proposed a perturbed test statistic for the SST
\begin{align*}
	\tT^*_n=\sup_{\btheta\in\Theta_{\theta}}\frac{1}{n}\left\|\sum_{i=1}^{n}\left[\nu_i\{\psi(\bV_i, \hbalpha, \bzero, \btheta) + \hat{K}(\btheta)\hat{J}\psi_{1}(\bV_i, \hbalpha)\}\right]\right\|^2_{\tV(\btheta)^{-1}},
\end{align*}
where $\tV(\btheta)$ is defined in \eqref{score_test0}, and $\{\nu_i,i=1,\cdots,n\}$ are $n$ independent and identically distributed copies of a random variable $\nu$ generated from $\mathcal{N}(0,1)$ that is independent of $\bV$.

The critical value $C_{\alpha}$ can be calculated by the empirical distribution of $\tT^*_n$ according to replicates, say $\{\tT_n^{*, j},j=1,\cdots,N\}$, where $C_{\alpha}$ is an upper quantile of the empirical distribution of $\tT^*_n$ and $N$ is a large integer. The null hypothesis is rejected for an $\alpha$-level test when $\tT_n>C_{\alpha}$, and the $\pvalue$ is approximated by $N^{-1}\sum_{j=1}^{N}\bone(\tT_n>\tT_n^{*, j})$.

For the WAST, we recommend the parametric bootstrap method to approximate the critical value or $\pvalue$. For each bootstrap sample $\{\bV_i^{*, b}, i=1, \cdots, n\}$, generated as described in Section 3, where $b=1, \cdots, B$, we calculate $T_n^{*, b}$ by
\begin{align}\label{statistic1*}
	T_n^{*, b} = \frac{1}{n(n-1)}\sum_{i\neq j}\int_{\btheta\in\Theta_{\theta}}\psi(\bV_i^{*, b}, \hbalpha^{*, b}, \bzero, \btheta)\trans\psi(\bV_j^{*, b}, \hbalpha^{*, b}, \bzero, \btheta)w(\btheta)d\btheta,
\end{align}
where $B$ is a large integer and $\hbalpha^{*, b}$ is the estimate of $\balpha$ under the null hypothesis. The critical value can be approximated by the upper quantile of $\{T_n^{*, b}, b=1,\cdots,B\}$. The $\pvalue$ is approximated by $B^{-1}\sum_{b=1}^{B}\bone(T_n>T_n^{*, b})$. In Section~\ref{sec:examples}, we discuss the parametric bootstrap method using several useful examples.

If the density $f(\bv; \balpha, \bbeta, \btheta)$ of $\bV$ is given except for the unknown parameters $\balpha$, $\bbeta$, and $\btheta$, and bootstrap samples $\{\bV_i^*, i=1, \cdots, n\}$ are generated from the distribution $f(\bv; \hbalpha, \bzero, \btheta_0)$, then Theorem~\ref{thm13} provides the asymptotic consistency of the bootstrap distribution of the test statistic, where $\hbalpha$ is the estimate of $\balpha$ under the null hypothesis. 

\begin{theorem}\label{thm13}
	If Assumptions~{(A1)--(A5)} and {(A7)} in the Appendix hold, then under the null hypothesis, we have
	\begin{align*}
		\sup_{x\in\bR}\left|\rP^*(nT_n^*\leq x)-\rP(nT_n\leq x)\right| \lpkonv 0,
	\end{align*}
	where $\lpkonv$ denotes convergence in probability and $\rP^*$ denotes the probability under the bootstrap procedure.
\end{theorem}

\section{Choice of the Weight}\label{sec:examples}
In Section~\ref{sec:test_stat}, we have theoretically established the asymptotic distributions of $T_n$ under the null and alternative hypotheses in the GEE framework. Here, we discuss in detail the closed form of $T_n$ in practice by taking an appropriate weight. We also consider the bootstrap method for each example.

We first consider the quantile regression in example~\ref{exam:qr} of Section~\ref{sec:intro}. By Lemma~A.1 in the Supplementary Material and taking the weight $w(\btheta)$ as the multivariate Gaussian density with mean $\bmu$ and covariance $\Sigma$, $\rho_{ij}$ in the WAST \eqref{statistic1} reduces to
\begin{align}\label{rho_qr}
	\begin{split} \rho_{ij}=\omega_{ij}\bX_i\trans\bX_j\left[\bone(Y_i-\tbX_i\trans\hbalpha(\tau)\leq0)-\tau\right]\left[\bone(Y_j-\tbX_j\trans\hbalpha(\tau)\leq0)-\tau\right],
	\end{split}
\end{align}
where
\begin{align*}
	\omega_{ij} = \int_{-\infty}^{\frac{\bZ_j\trans\bmu}{\|\bZ_j\|_{\Sigma}}}\phi(z)\Phi\left(\frac{1}{\sqrt{1-\varrho_{ij}^2}}\left(\frac{\bZ_i\trans\bmu}{\|\bZ_i\|_{\Sigma}}-\varrho_{ij} z\right)\right)dz, \quad \mbox{if}~ i\neq j,
\end{align*}
and $\varrho_{ij} = \bZ_i\trans\Sigma\bZ_j(\|\bZ_i\|_{\Sigma}\|\bZ_j\|_{\Sigma})^{-1}$.

From a Bayesian perspective, $\bmu$ and $\Sigma$ can be treated as priors and chosen empirically. Although there is no closed-form expression of $\omega_{ij}$ for $i\neq j$ in general, we can approximate $\omega_{ij}$ by the Monte Carlo method. More specifically, we approximate $\omega_{ij}$ by
\begin{align}\label{wij_approx} \omega_{ij}\approx N^{-1}\sum_{k=1}^N\left[\Phi\left(\frac{1}{\sqrt{1-\varrho_{ij}^2}}\left(\frac{\bZ_i\trans\bmu}{\|\bZ_i\|_{\Sigma}}-\varrho_{ij} z_k\right)\right)\bI\left(z_k\leq \|\bZ_j\|_{\Sigma}^{-1}\bZ_j\trans\bmu\right)\right],
\end{align}
where $z_k$ follows $\mathcal{N}(0, 1)$ and $N$ is a large integer.
In practice, we can simply choose $\bmu=\bzero$ and $\Sigma=\bI$, whereby Lemma~A.1 in the Supplementary Material implies that $\omega_{ij}$ reduces to
\begin{align}\label{wij}
	\omega_{ij}=\frac{1}{4}+\frac{1}{2\pi}\arctan\left(\frac{\varrho_{ij}}{\sqrt{1-\varrho_{ij}^2}}\right) \quad \mbox{if}~ i\neq j.
\end{align}

If $\theta\in[0, 1]$, then we can take $w(\theta)$ as being the density of the beta distribution Beta($\lambda_1$, $\lambda_2$) with shape parameters $\lambda_1, \lambda_2>0$. In this case, we find that
\begin{align*}
	\omega_{ij} = B(\min\{Z_i, Z_j\};\lambda_1, \lambda_2),
\end{align*}
where $B(\cdot;\lambda_1, \lambda_2)$ is the cumulative distribution function of Beta($\lambda_1$, $\lambda_2$).
If $\theta\in\mathbb{R}$, then we can take $w(\theta) = \phi(\theta;\mu, \sigma^2)$ as the density of the Gaussian distribution with mean $\mu$ and variance $\sigma^2$, which implies that
\begin{align*}
	\omega_{ij} = \Phi(\min\{Z_i, Z_j\};\mu, \sigma^2),
\end{align*}
where $\Phi(\cdot;\mu, \sigma^2)$ is the cumulative distribution function of $\mathcal{N}(\mu, \sigma^2)$.
Alternatively, we can write $\check{\btheta}=(-\theta, 1)\trans$ and $\check{\bZ}=(1, Z)\trans$, which imply that $\bone(Z\geq\theta)$ is equal to $\bone(\check{\bZ}\trans\check{\btheta}\geq0)$. This allows us to use the bivariate normal density as the weight, which results in $\omega_{ij}$ in \eqref{wij}.

As emphasized in Section 2.2, the choice of the weight affects the computation of the test statistic since the numerical integration over $\bR^q$. 
We conduct numerical studies based on GLMs, probit, quantile and semiparametric models to investigate the sensitivity of the weight's choice by comparing $\omega_{ij}$ with the closed-form in \eqref{wij} and the approximated $\omega_{ij}$ in \eqref{wij_approx}, see details in Appendix B.5 in the Supplemetary Material. From the numerical results, compared with the approximated $\omega_{ij}$ in \eqref{wij_approx}, the test statistic with $\omega_{ij}$ in \eqref{wij} has higher power uniformly, and only take 10\% time computationally when $N=10000$. This strongly recommends one to use the WAST with the closed-form in \eqref{wij}.

For the bootstrap procedure, as with the wild bootstrap of \citet{feng2011}, we independently generate a random variable $\nu$ satisfying the following conditions: (i) there are two positive constants $c_1$ and $c_2$ satisfying $\sup\{\nu\in\cG:\nu\leq 0\}=-c_1$ and $\inf\{\nu\in\cG:\nu\geq 0\}= c_2$, where $\cG$ is the support of $\nu$; (ii) the distribution $G$ of $\nu$ satisfies $\int_{0}^{\infty}\nu^{-1}g(\nu)d\nu=-\int_{-\infty}^{0}\nu^{-1}g(\nu)d\nu=1/2$ and $\bE[|\nu|]<\infty$, where $g(\nu)$ is the density function of $\nu$; (iii) the $\tau$th quantile of $\nu$ is zero. We obtain the bootstrapped samples $\{(Y_i^*, \tbX_i, \bX_i, \bZ_i), i=1, \cdots, n\}$ and calculate $T_n^*$ given in \eqref{statistic1*}, where $Y_i^*=\tbX_i\hbalpha(\tau)+\nu_i|\tilde{r}_i|$ and $\tilde{r}_i=Y_i-\tbX_i\hbalpha(\tau)$, with $\hbalpha(\tau)$ being the estimate of $\balpha(\tau)$ based on $\Psi_n(\balpha)=\bzero$ under the null hypothesis. In practice, one can generate $\nu$ following a discrete distribution with $\rP(\nu=2(1-\tau))=1-\tau$ and $\rP(\nu=-2\tau)=\tau$. 

Consider the GLMs in example \ref{exam:glm}. $\rho_{ij}$ in the proposed test statistic becomes
\begin{align}\label{rho_glm}
	\rho_{ij} = \omega_{ij}\bX_i\trans\bX_j\left[Y_i-\partial c(\tbX_i\trans\hbalpha)/\partial\mu_i\right]\left[Y_j-\partial c(\tbX_j\trans\hbalpha)/\partial\mu_j\right],
\end{align}
where the weight $w(\btheta)$ can be chosen as in the case of quantile regression. We generate the bootstrap samples from GLM with the density $f(\bV_i, \hbalpha,\bzero,\btheta_0)$, which is the density under the null hypothesis by replacing $\balpha$ with $\hbalpha$. We obtain the bootstrap samples $\{(Y_i^*, \tbX_i, \bX_i, \bZ_i), i=1, \cdots, n\}$ and calculate $T_n^*$ given in \eqref{statistic1*}. The consistency of the test statistic $T_n^*$ is given in Theorem~\ref{thm13}.

We now consider probit regression model in example~\ref{exam:prob}, in which the test statistic is
\begin{align}\label{rho_probit}
	\rho_{ij} = \omega_{ij}\psi_0(\bV_i, \hbalpha)\psi_0(\bV_j, \hbalpha),
\end{align}
where $\psi_0(\bV_i, \balpha) = Y_i\bX_i\frac{\phi(\tbX_i\balpha)}{\Phi(\tbX_i\balpha)} - (1-Y_i)\bX_i\frac{\phi(-\tbX_i\balpha)}{\Phi(-\tbX_i\balpha)}$ and $\omega_{ij}$ can be chosen as in GLMs.
The bootstrap samples can be drawn from the probit regression model under the null hypothesis by replacing $\balpha$ with $\hbalpha$. We obtain the bootstrap samples $\{(Y_i^*, \tbX_i, \bX_i, \bZ_i), i=1, \cdots, n\}$ and calculate $T_n^*$ given in \eqref{statistic1*}, where $Y_i^*= \bone(\nu_i\leq \tbX_i\trans\hbalpha)$ and $\nu_i$ is independent of $\bV_i$ and normally distributed with zero mean and unit variance. The consistency of the test statistic $T_n^*$ is given in Theorem~\ref{thm13}.

For the semiparametric models in example \ref{exam:semi}, we have
\begin{align}\label{rho_dr}
	\begin{split}
		\rho_{ij}=\omega_{ij}\left[A_i-\pi(\bZ_i, \hbalpha_1)\right]\left[Y_i-\gamma(\bZ_i, \hbalpha_2)\right]\left[A_j-\pi(\bZ_j, \hbalpha_1)\right]\left[Y_j-\gamma(\bZ_j, \hbalpha_2)\right],
	\end{split}
\end{align}
where $\omega_{ij}$ is the same as that in the quantile regression of example~\ref{exam:qr}. To calculate the critical value or p-value of $T_n$, we can use the wild bootstrap strategy to generate bootstrap samples \citep{1986Wu, 1988Liu, 1993Mammen}. More specifically, we obtain the bootstrap samples as $\{(Y_i^*, A_i, \bZ_i), i=1, \cdots, n\}$ with $Y_i^*=\gamma(\bZ_i, \hbalpha_2)+\nu_i\tilde{r}_i$ and calculate $T_n^*$ given in \eqref{statistic1*}, where $\tilde{r}_i=Y_i-\gamma(\bZ_i, \hbalpha_2)$ and $\nu_i$ is independent of $(Y_i, A_i, \bZ_i)$ and normally distributed with zero mean and unit variance. 

\section{Simulation Studies}\label{sec:simulations}

\subsection{Performance for GLMs}\label{GLM} 
In this section, we demonstrate the finite-sample performance of the proposed WAST. We report here simulation results from several scenarios involving different numbers of $\bZ$ and distributions of $\Gv$ $\bZ$ for GLMs and quantile regression, probit regression and semiparametric models. Due to limit of space, we give simulation results of probit regression and semiparametric models for high- or low-dimensional $\bZ$, as well as multiple change planes, in the Appendix B and C of the Supplementary Material. As mentioned in Section \ref{sec:examples}, we choose Gaussian density with zero mean and covariance $\bI$ as the weight in all simulations.

For fair comparison, we use the same settings as considered by \citet{2021Threshold}.
Consider the GLM in model~(\ref{model:glm}) with the canonical parameter
\begin{align*}
	\mu = \tbX\trans\balpha + \bX\trans\bbeta\bone(\bZ\trans\btheta\geq0),
\end{align*}
where $\bbeta = \kappa\bone_p$ and $\alpha_2=\cdots=\alpha_r=\log(1.4)$. For Gaussian and Poisson families,
we set $\alpha_1 = 0.5$,  and for binomial family $\alpha_1$ is chosen so that the proportion of cases in the data set is ${1}/{3}$ on average under the null hypothesis. $\theta_2,\cdots,\theta_{q}$ are equally spaced numbers from $-1$ to $1$, and $\theta_1$ is chosen as the negative of the 0.65 percentile of $Z_2\theta_2+\cdots+Z_q\theta_q$, which means that $\bZ\trans\btheta$ divides the population into two groups with 0.35 and 0.65 observations, respectively. For the predictor, we generate $(v_2, \cdots,v_{\max\{r,p\}}, Z_2,\cdots,Z_q)\trans$ from a multivariate normal distribution with mean $\bzero$ and covariance $\Sigma=(\tilde{\rho}_{ij})$, where $\tilde{\rho}_{ij}=1$ if $i=j$ and $\tilde{\rho}_{ij}=\rho$ otherwise. Here, we consider both $\rho=0$ and $\rho=0.3$. We set $\tX_j=\bone(v_j>0)$ with $j=2,\cdots,r$, $X_k = \bone(v_k>0)$ with $k=2,\cdots,p$, and $\tX_1=1$, $X_1=1$, and $Z_1=1$.

We evaluate the power under a sequence of alternative models indexed by $\kappa$, i.e., $H_1^{\kappa}: \bbeta = \kappa\bone_p$ with $\kappa$ equally spaced in the range $(0,\kappa_{\max}]$. We set the sample size as $n=(300, 600)$, and we set 1000 repetitions and 1000 bootstrap samples. For comparison, we consider three methods here: (i) the proposed WAST, (ii) the SST, and (iii) the approximated supremum of likelihood-ratio test statistic (SLRT) given by \citet{2021Threshold}. We calculate both SLRT and SST over $\{\btheta^{(k)}=(\theta^{(k)}_1,\cdots,\theta^{(k)}_q)\trans: k=1,\cdots,K\}$. Let $\btheta^{(k)}_{-1}=\tilde{\btheta}^{(k)}_{-1}/\|\tilde{\btheta}^{(k)}_{-1}\|$, where $\btheta^{(k)}_{-1}=(\theta^{(k)}_2,\cdots,\theta^{(k)}_q)\trans$ and $\tilde{\btheta}^{(k)}_{-1}=(\tilde{\theta}^{(k)}_2,\cdots,\tilde{\theta}^{(k)}_q)\trans$, and $\tilde{\btheta}^{(k)}_{-1}$ is drawn independently from the multivariate normal distribution $\mathcal{N}(\bzero_{q-1},\bone_{q-1})$. For each $\theta_1^{(k)}$, $k=1,\cdots,K$, we select it by an equal grid search in the range from the lower tenth percentile to the upper tenth percentile of the data points of $\{\theta^{(k)}_2Z_{i2}+\cdots+\theta^{(k)}_qZ_{iq}\}_{i=1}^n$, which is the same as that in \citet{2021Threshold}. Here, we set $K=1000$. 
And the critical value of the test statistic SST is calculated with the resampling method used in \cite{2017Change}. The results for each scenario are summarized based on 1000 repetitions.

\autoref{table_size_m} lists the type-\uppercase\expandafter{\romannumeral1} errors ($\kappa=0$) with the nominal significance level of 0.05. As can be seen, the sizes of the proposed method are close to 0.05. For the SST, its sizes are far from 0.05 in most of the scenarios, and it performs the worst of these three methods. For the SLRT, its sizes are close to 0.05. The type-\uppercase\expandafter{\romannumeral1} errors of the Gaussian and Poisson families and other settings for different parameter dimensions can be found in Table C.1 in the Supplementary Material, from which we can observe that the proposed method WAST has a good control over the type-\uppercase\expandafter{\romannumeral1} error rate.

\autoref{fig_binomial_m} and Figures C.1-C.3 in the Appendix C.1 of the Supplementary Material show the powers for all the scenarios considered in \autoref{table_size_m} and Table C.1. As can be seen, the powers increase as the sample size $n$ increases, which is verified by the asymptotic theory. For the Gaussian, binomial and Poisson families, the power of the proposed WAST increases faster than those of the SST and SLRT, and this superiority is particularly evident when $n$ is small or $p$ and $q$ increase. In summary, the proposed test statistic shows very competitive performance.

We consider different distributions of grouping variable $\bZ$ in Appendix C.2 of the Supplementary Material. From Tables C.2-C.3 and Figures C.4-C.9, we 
have similar conclusion as these for GLMs with $\bZ$ generated from a normal distribution with zero mean and unit variance. 
For large numbers of $\bZ$, we analyze numerically performance of the WAST in Appendices B.3-B.4 of the Supplementary Material, from which  the WAST still overperforms the competitors 
as scenarios with small numbers of grouping variables.

\begin{table}
	\centering
	\caption{\label{table_size_m} Type-\uppercase\expandafter{\romannumeral1} errors of the WAST, SST, and SLRT  based on resampling for binomial families with different numbers of $\bZ$ and $\rho=0$. The nominal significant level is 0.05.}
	\resizebox{\textwidth}{!}{
		\begin{threeparttable}
			\begin{tabular}{l*{15}{c}}\\
				\hline
				\multirow{2}{*}{$n$} & \multicolumn{3}{c}{$(r,p,q)=(2,2,3)$} && \multicolumn{3}{c}{$(6,6,3)$} && \multicolumn{3}{c}{$(2,2,11)$} && \multicolumn{3}{c}{$(6,6,11)$} \\
				\cline{2-4} \cline{6-8} \cline{10-12} \cline{14-16}
				& WAST& SST & SLRT && WAST& SST & SLRT && WAST& SST & SLRT && WAST& SST & SLRT \\
				\cline{2-16}
				300 & 0.053 & 0.109 & 0.052 && 0.056 & 0.081 & 0.055 && 0.046 & 0.124 & 0.053 && 0.052 & 0.069 & 0.049 \\
				600 & 0.050 & 0.082 & 0.054 && 0.045 & 0.110 & 0.067 && 0.057 & 0.089 & 0.065 && 0.036 & 0.099 & 0.041 \\
				\hline
			\end{tabular}
		\end{threeparttable}
	}
\end{table}

\begin{figure}[!ht]
	\begin{center}
		\includegraphics[scale=0.285]{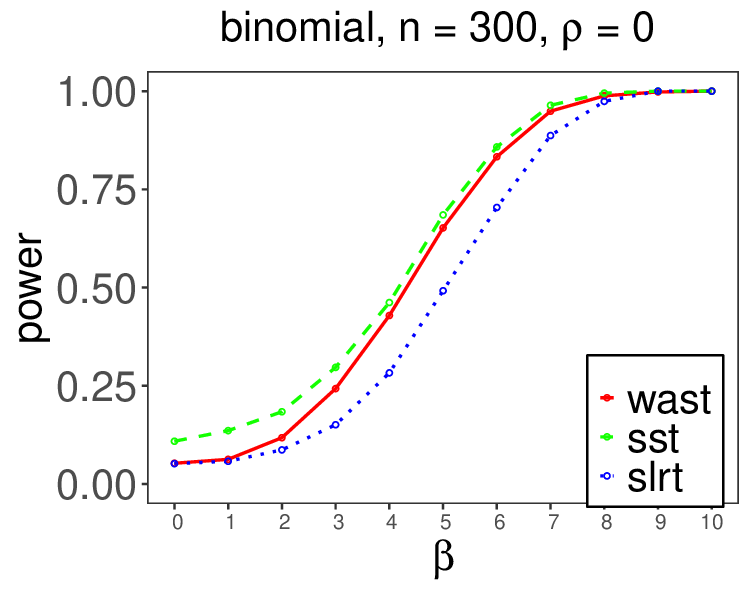}
		\includegraphics[scale=0.285]{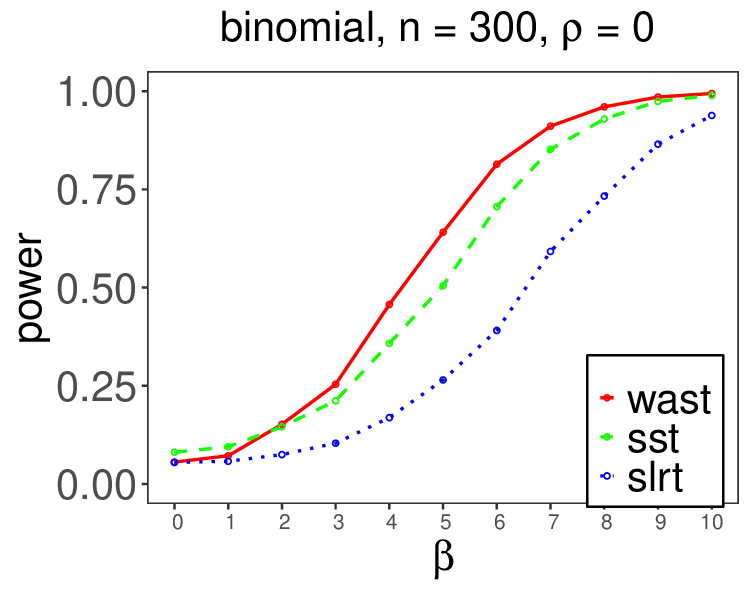}
		\includegraphics[scale=0.285]{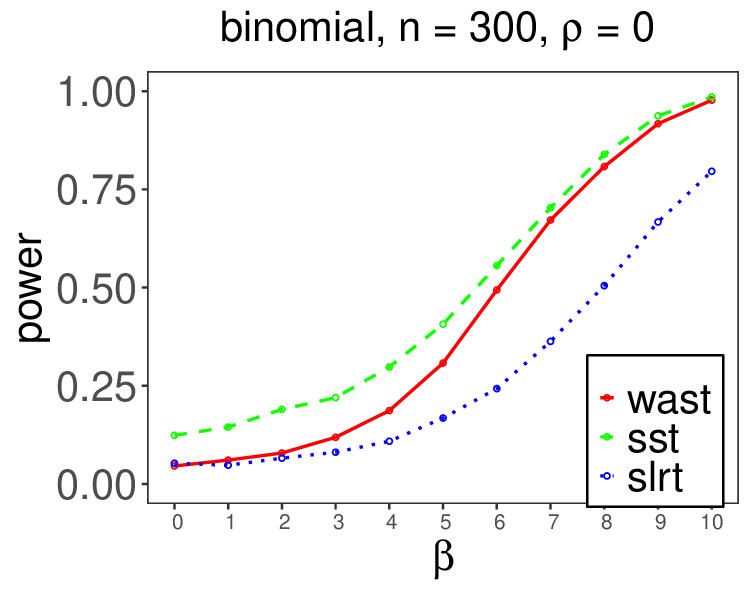}
		\includegraphics[scale=0.285]{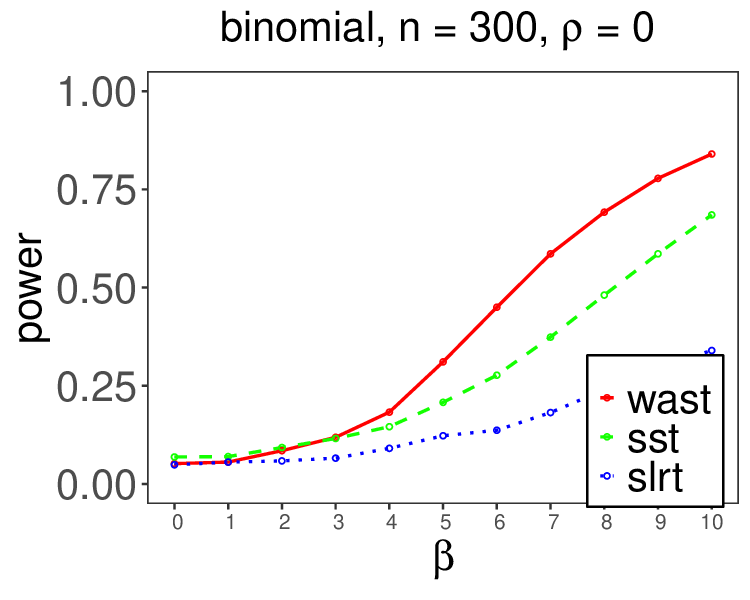}
		\includegraphics[scale=0.285]{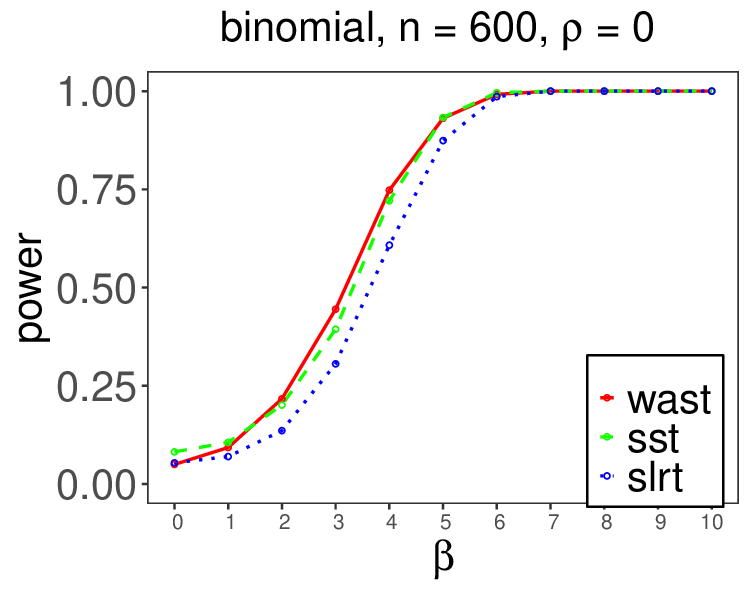}
		\includegraphics[scale=0.285]{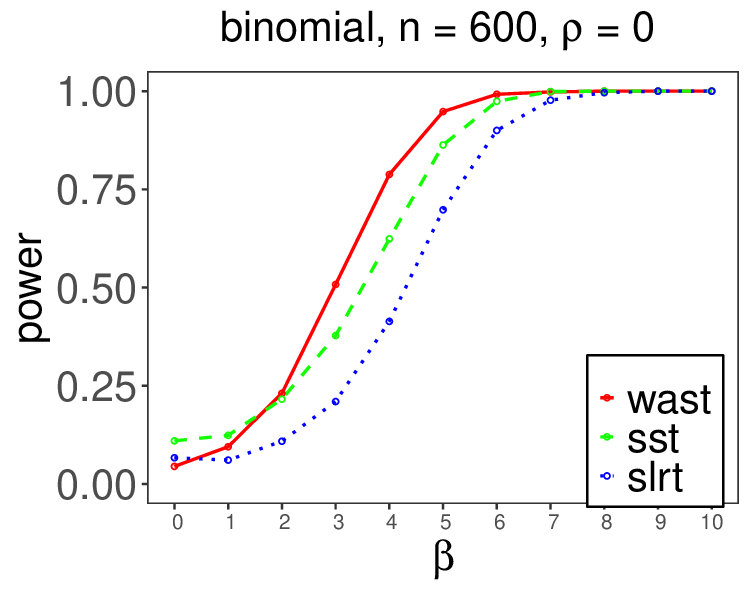}
		\includegraphics[scale=0.285]{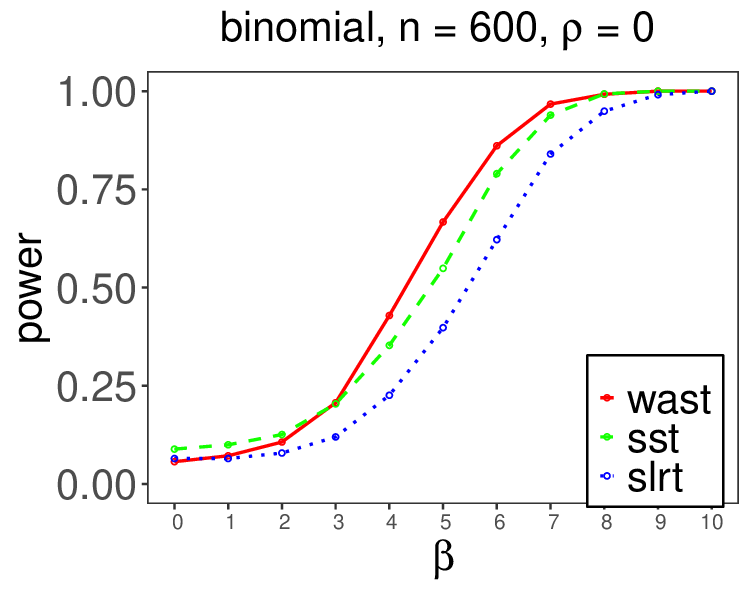}
		\includegraphics[scale=0.285]{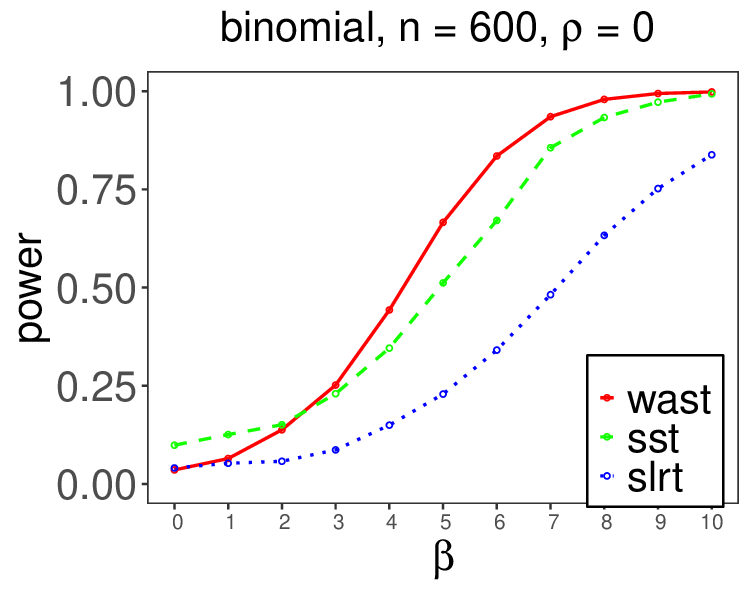}
		\caption{Powers of testing logistic regression by WAST (red solid line), SST (green dashed line), and SLRT (blue dotted line) for $n=(300,600)$. From left to right, each column depicts the powers for the cases $(r,p,q)=(2,2,3)$, $(6,6,3)$, $(2,2,11)$, and $(6,6,11)$.}
		\label{fig_binomial_m}
	\end{center}
\end{figure}

\subsection{Performance for quantile regression}
For fair comparison, we consider the settings as similar to the existing literature \citet{2011Testing,2017Change}, which are different from the GLMs in Section \ref{GLM}.

Consider the quantile regression model~(\ref{model:qr}),
\begin{align*}
	Y_i = 0.5+\tX_i+\bX_i\trans\bbeta(\tau)\bone(\bZ_i\trans\btheta(\tau)\geq0) + \eps_i,
\end{align*}
where $\tX_i$ is generated from standard normal distribution, $\bX_i$ from $p$-variate normal distribution with mean $\bzero$ and variance $\sqrt{2}I$, and $Z_{1i}=1$ and $(Z_{2i},\cdots,Z_{qi})\trans$ from $(q-1)$-variate standard normal distribution. The error $\eps_i$ is generated from standard normal distribution. 
We set $\bbeta=(1,\cdots,1)\trans$ under $H_0$, and $(\theta_2, \cdots,\theta_q)\trans = (1,2,\cdots,2)\trans$ under $H_1$, where $\theta_1$ is chosen as the negative of the 0.65 percentile of $Z_2\theta_2+\cdots+Z_q\theta_q$, which means that $\bZ\trans\btheta$ divides the population into two groups with 0.35 and 0.65 observations, respectively.

For all models with change plane analysis, we evaluate the power under a sequence of alternative models indexed by $\kappa$, that is $H_1^{\kappa}: \bbeta^{\kappa}=\kappa\bbeta^*$ with $\kappa = i/10$ for semiparametric model and $\kappa=i/20$ for others, $i=1,\cdots,10$, where $\bbeta^*=(1,\cdots,1)\trans$. We set sample size $n=(200, 600)$, 1000 repetitions and 1000 bootstrap samples, and report in \autoref{fig_qr} and Figures C.11-C.12 in the Appendices C.4 of the Supplementary Material the performance for both the WAST and SST. We calculate SST over $\{\btheta^{(k)}=(\theta^{(k)}_1,\cdots,\theta^{(k)}_q)\trans: k=1,\cdots,K\}$ with the number of threshold values $K=2000$. Let $\btheta^{(k)}_{-1}=\tilde{\btheta}^{(k)}_{-1}/\|\tilde{\btheta}^{(k)}_{-1}\|$, where $\btheta^{(k)}_{-1}=(\theta^{(k)}_2,\cdots,\theta^{(k)}_q)\trans$ and $\tilde{\btheta}^{(k)}_{-1}=(\tilde{\theta}^{(k)}_2,\cdots,\tilde{\theta}^{(k)}_q)\trans$, and $\tilde{\btheta}^{(k)}_{-1}$ is drawn independently from $(r-1)$-variate standard normal distribution.  For each $\theta_1^{(k)}$, $k=1,\cdots,K$, we select it by equal grid search in the range from the lower 10th percentile to upper 10th percentile of the data points of $\{\theta^{(k)}_2Z_{2i}+\cdots+\theta^{(k)}_qZ_{qi}\}_{i=1}^n$, which is same as that in \cite{2020Threshold}. Here we consider two combinations of $(p,q)=(1, 3), (5, 5)$. The power results of another two combinations of $(p,q)= (10, 10), (50, 20)$ can be found in Figures C.13-C.14 in Appendix C.4 of the Supplementary Material.

Type-\uppercase\expandafter{\romannumeral1} errors ($\kappa=0$) with sample size $n=(200, 600)$ are listed in \autoref{table_size1_m} and Table C.5 in Appendix C.4 of the Supplementary Material. We can see from these tables that the size of the proposed WAST are close to the nominal significance level $0.05$, but for most scenarios the size of the SST are much smaller than 0.05.
\autoref{fig_qr} and Figures C.11-C.14 in the Appendix C.4 of the Supplementary Material indicate that powers become greater as sample size $n$ increases, which are verified by the asymptotic theory. The WAST has comparable power with 
the SST for the quantile regression, but the size of the SST is much less than the nominal level 0.05 in \autoref{table_size1_m}. To save space here, we report the performance of the WAST 
for probit and semiparametric models in Appendix C.4 of the Supplementary Material, 
and the performance of the WAST for the large numbers of $\bZ$ is analysed in Appendices B.1-B.2 of the
Supplementary Material. We can observe that the WAST works as well in scenarios with a large number of grouping variables as it does in situations with a small number of grouping variables. 

\begin{table}[htp!]
	\def~{\hphantom{0}}
	\centering
    \footnotesize
	\caption{Type \uppercase\expandafter{\romannumeral1} errors of the proposed WAST and SST based on resampling for quantile regression with $\tau=0.5$ and $0.7$. The nominal significant level is 0.05.
}
	\begin{threeparttable}
		\begin{tabular}{lccccccccccc}
			\hline
			\multirow{3}{*}{$(p,q)$} & \multicolumn{5}{c}{$\tau=0.5$} && \multicolumn{5}{c}{$\tau=0.7$} \\
			\cline{2-6} \cline{8-12}
			& \multicolumn{2}{c}{$n=200$} && \multicolumn{2}{c}{$n=600$} && \multicolumn{2}{c}{$n=200$} && \multicolumn{2}{c}{$n=600$} \\
			\cline{2-3} \cline{5-6} \cline{8-9} \cline{11-12}
			& WAST& SST && WAST& SST && WAST& SST && WAST& SST\\
			\cline{2-12}
			$(1,3)$ & 0.060 & 0.039  && 0.053 & 0.034 && 0.044 & 0.020 && 0.053 & 0.037\\
			$(5,5)$ & 0.049 & 0.005  && 0.047 & 0.033 && 0.048 & 0.005 && 0.045 & 0.024\\
			\hline
		\end{tabular}
	\end{threeparttable}
	\label{table_size1_m}
\end{table}

\begin{figure}[!ht]
	\begin{center}
		\includegraphics[scale=0.285]{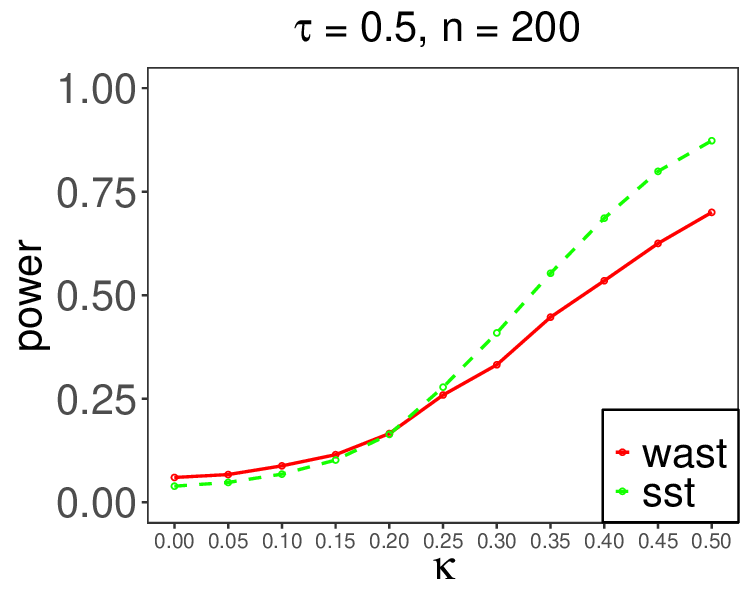}
		\includegraphics[scale=0.285]{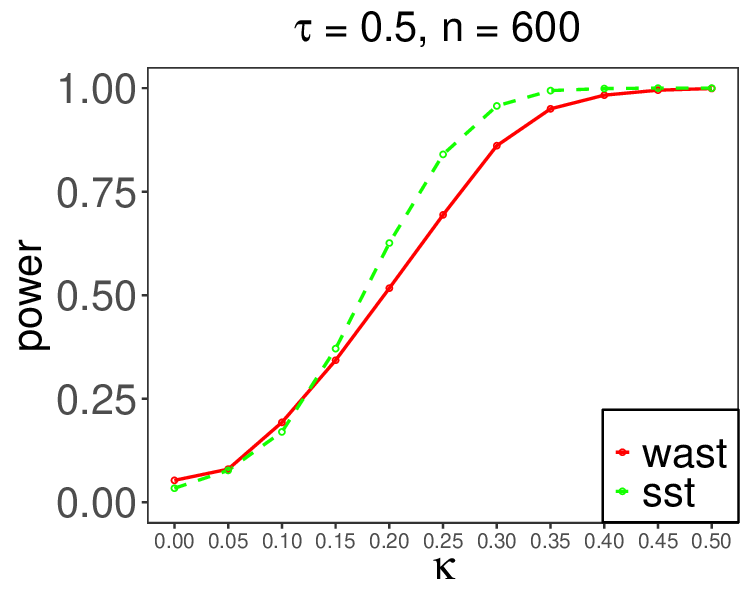}
		\includegraphics[scale=0.285]{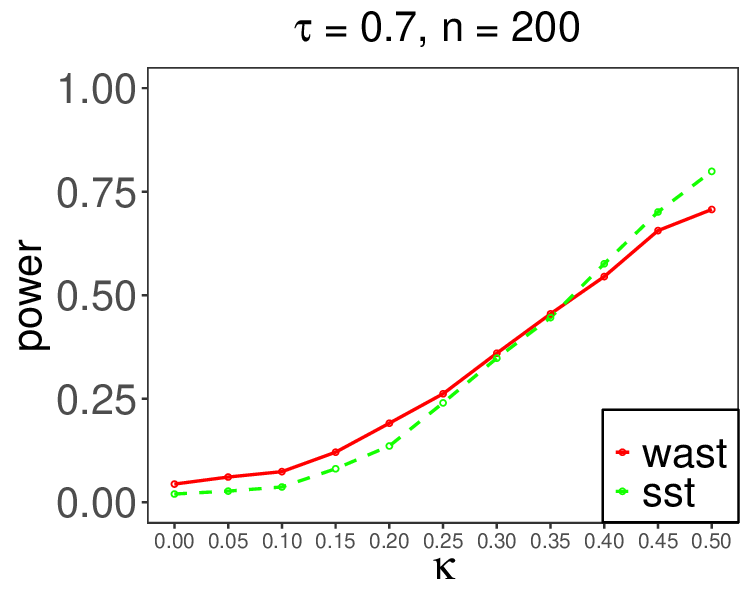}
		\includegraphics[scale=0.285]{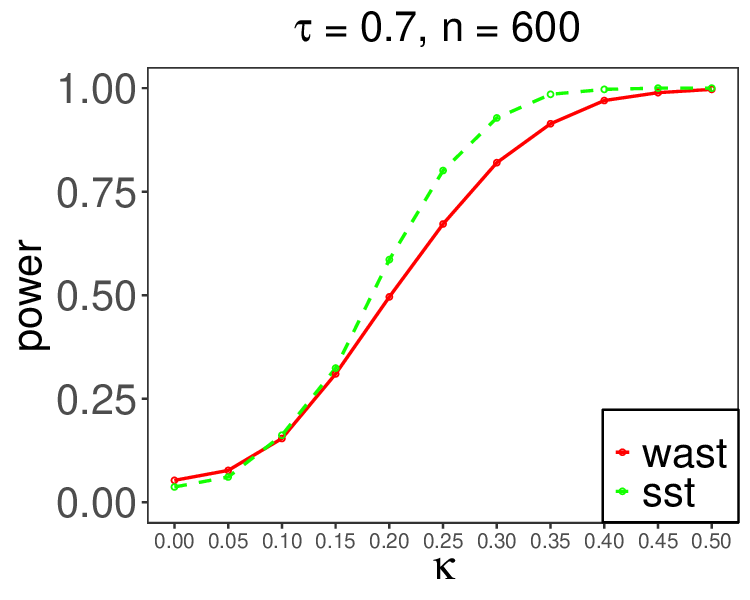}
		\includegraphics[scale=0.285]{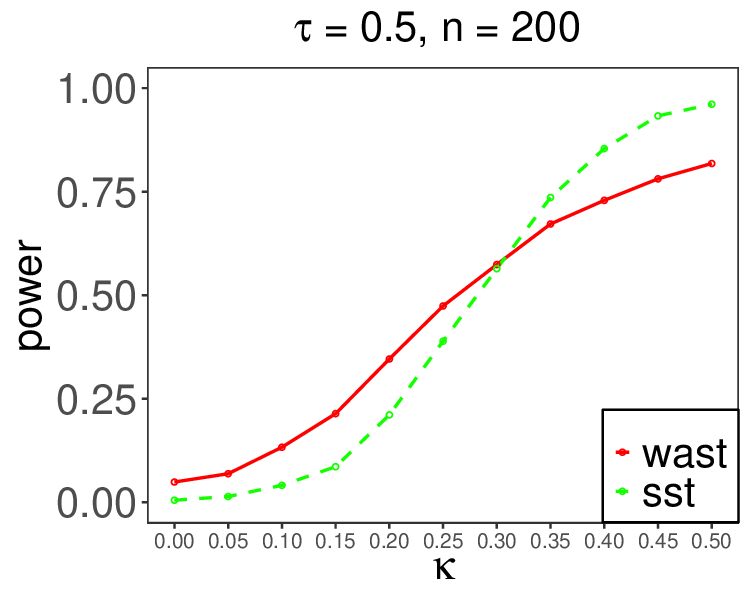}
		\includegraphics[scale=0.285]{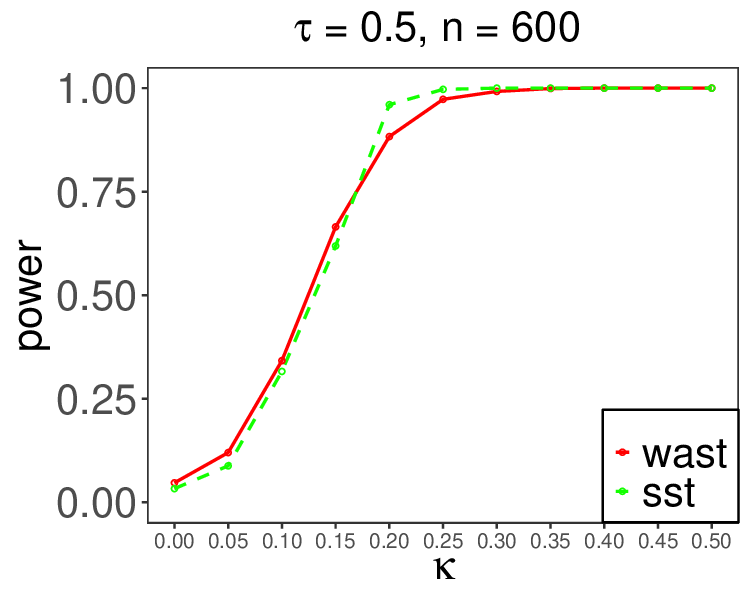}
		\includegraphics[scale=0.285]{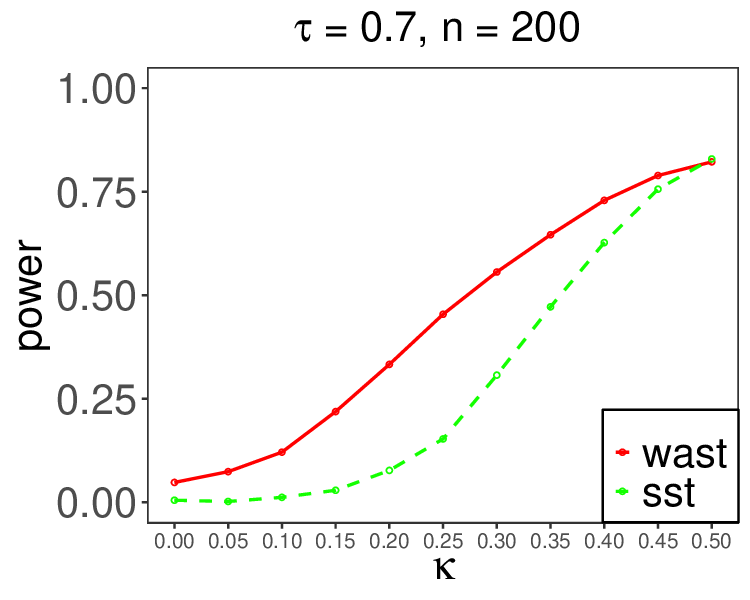}
		\includegraphics[scale=0.285]{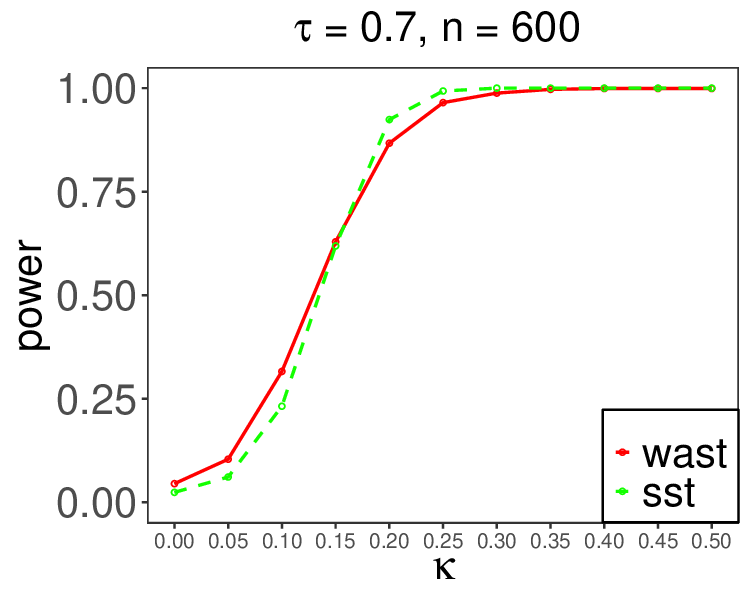}
		\caption{\it Powers of test statistic for quantile regression with $\tau=0.5$ and $\tau=0.7$ and with different numbers of $\bZ$ by the proposed WAST (red solid line) and SST (green dashed line). From top to bottom, each row depicts the powers for $(p,q)=(1,3)$, and $(p,q)=(5,5)$.}
		\label{fig_qr}
	\end{center}
\end{figure}

\section{Empirical Application}\label{sec:case_studies}


In this section, we demonstrate our proposed method by analyzing a data set obtained from China Health and Retirement Longitudinal Study (CHARLS) in 2013. This nationally representative data set focuses on individuals aged 45 and above in China, and includes comprehensive demographic, health, and economic information. CHARLS aims to investigate the primary health and economic adaptations associated with the country's rapid population aging. 
It is well known that the pension income from the society and financial support from adult children are two main sources of income for the elderly. 
Our analysis aims to determine if there are specific subgroups in which pension income affects the elderly's support income differently. The subjects of our study are rural seniors aged 60 and older in Shanxi Province. After data-cleaning, the data set in Shanxi Province includes information from 183 individuals.

We treat the support income (1000 CNY/year) for the elderly as the response variable $Y$. The baseline variable $\tbX$ includes 
pension income (10,000 CNY/year), the number of children with an annual income of more than 10,000 CNY, and emotional comfort (whether the elderly contact with their children at least every half month: Yes = 1, No = 0). 
To identify the different effects of pension income on support income is of particular interest, therefore, the grouping difference variable $\bX$ is taken as the pension income. 
The grouping variables $\bZ$ contain some individual characteristics of elderly people, including gender (male = 1, female = 0), age, and living arrangement (whether they live with children: Yes = 1, No = 0).

Figure D.1 in the Supplementary Material exhibits the histogram of support incomes for the elderly in Shanxi Province. 
It can be noticed that the support income has a heavy-tailed distribution. 
In view of this, we consider quantile change-plane regression for different $\tau\in (0,1)$:
\begin{align}
	Q_{Y_i}(\tau|\tbX_i, \bX_i, \bZ_i) = \tbX_i\trans\balpha(\tau)+\bX_i\trans\bbeta(\tau)\bone(\bZ_i\trans\btheta(\tau)\geq0),
\end{align}
where $Y_i, \tbX_i, \bX_i, \bZ_i$ are sample values of $Y, \tbX, \bX, \bZ$, respectively, and the first element of $\tbX_i$ is 1 which corresponds to the intercept.

Calculated based on 5000 bootstrap samples, the p-values of the proposed WAST are 0.021 if $\tau=0.25$, $<0.0001$ if $\tau=0.5$, and 0.265 if $\tau=0.75$. 
This demonstrates that for those who have low and median support incomes, we have strong evidence for rejecting $H_0$, i.e., there exist two subgroups wherein the pension income has a different impact on support income of the elderly. 
Furthermore, it suggests that the pension income does not have a subgroup effect on the support income for individuals with high support incomes.

Next, we estimate parameters $(\balpha,\bbeta,\btheta)$ according to the test results. We first randomly choose 1000 $\btheta$s.  For every given $\btheta$, 
the estimator of $\balpha$ and $\bbeta$ can be obtained from classic quantile regression. 
Then the estimated grouping parameter $\hbtheta$ is obtained by minimizing the empirical quantile loss functions over the given 1000 $\btheta$s, and consequently, one of subgroups is given by $\{i:\bZ_i\trans\hat{\btheta}\geq0\}$.

\autoref{table_CHARLS} lists the estimators for $\tau=0.25$ and $0.5$. According to $\htheta$ in \autoref{table_CHARLS}, one of subgroups, denoted by $G^+=\{i:\bZ^{T}_{i}\hat{\btheta}(\tau)\geq 0\}$, is $$G^+=\{i:-0.051*gender+0.014*age-0.023*living~ arrangement\ge 0.998\}.$$ There are 50 individuals in the subgroup $G^+$. Based on these results, we can see that for elderly people in the subgroup $G^+$, the pension income has a positive impact on the support income. 
Conversely, the effect is reversed in another subgroup. 

\begin{table}[htp!]
        \centering
	\caption{Estimates if the null hypothesis has been rejected.}
	\begin{threeparttable}
		\begin{tabular}{ccccccccccccc}
			\hline
			$\tau$&& \multicolumn{4}{c}{$\hat{\balpha}$}  && \multicolumn{4}{c}{$\hat{\btheta}$} && \multicolumn{1}{c}{$\hat{\bbeta}$}\\
            \cline{1-1} \cline{3-6} \cline{8-11} \cline{13-13}
			0.25 && 0.061 & -0.085 & 0.499 & 0.193  && -0.998 & -0.051 & 0.014 & -0.023 && 0.219\\
			0.5 && 0.283 & 0.087 & 1.170 & 0.217  && -0.998 & -0.051 & 0.014 & -0.023 && 1.226\\
			\hline
		\end{tabular}
	\end{threeparttable}
	\label{table_CHARLS}
\end{table}

All elderly people aged over 77 years fall within the estimated subgroup $G^{+}$, and thus their support income will increase as their pension income increases. 
There may be two reasons for this. Firstly, the living expenses of these individuals may be higher. Therefore, their pension income may be insufficient to cover their cost of living, necessitating the need for additional support income. 
Secondly, pension income can increase the unearned income of the elderly, allowing them to reduce working hours, have more leisure time, and consequently increase intergenerational care time. As a result, children may increase their financial support for their parents as a form of reciprocity. 

The estimator is $\hat{\balpha}=(1.019, -0.288, 2.125, 0.694)$ when $\tau=0.75$, which indicates that for individuals with high support income, the pension income exerts a negative impact on the support income.

\section{Conclusion}\label{sec:conclusion}

Herein, we have considered a class of change-plane tests for heterogeneous estimating equations. A novel test statistic, WAST, is introduced to test whether the change plane exists. In a comparison with existing test procedures, including the supremum of the squared score test from the framework of the average exponential LM test \citep{1993Tests, 1994Optimal, 2009On, 2017Change}, the proposed test statistic significantly improves the power in practice, particularly when the number of grouping variables is large. As a by-product, the WAST reduces the computational cost dramatically because it has a closed form when an appropriate weight is used. Asymptotic distributions are derived under the null and alternative hypotheses, and we extend the proposed method to models with multiple change planes in the GEE framework, as studied by \citet{2021Multithreshold}.
Many classic and useful models in the GEE framework are considered herein, such as quantile regression \citep{2022Single, 2018Oracle, 2011Testing}, GLMs \citep{2021Threshold}, probit regression \citep{2011Testing}, and semiparametric models \citep{2017Change}. Case studies are analyzed to further demonstrate the performance of the proposed test procedure. From the simulation studies and case studies, we conclude that compared with several competitors, the proposed WAST performs pretty well not only for low-dimensional scenarios but also for situations involving a large number of grouping variables.

We provide the hypothesis test to detect the existence of subgroups, but it is left blank how to partition the subjects, particularly in the setting with large numbers of grouping variables, which needs the estimate of parameters, especially of the grouping parameter $\btheta$. It is a possible strategy for the estimation of the grouping parameter to penalize the loss function with $\|\btheta\|_1$ under the assumption of sparsity.

\section*{Acknowledgements}
This work was supported by National Natural Science Foundation of China (12271329) and Program for Innovative Research Team of SUFE.

\section*{Disclosure statement}
The authors report there are no competing interests to declare.

\section*{Supplementary Material}\label{SM}
Appendix A includes proofs of Theorems \ref{thm11}--\ref{thm13}
and related Lemmas. Appendices B and C present additional simulation studies to illustrate the performance of the proposed test statistic. 
Appendix D shows the histogram of the support incomes for the elderly in Shanxi Province and additional analysis results of another real dataset.

\section*{Appendix}
We, herein, provide some assumptions required to establish the asymptotic properties.
\begin{enumerate}[({A}1)]
	\item $\bE\psi_1(\bV, \balpha_0)=0$ for some $\balpha_0\in\Theta_{\alpha}$. $\bE\psi_1(\bV, \hbalpha)\rightarrow \bzero$ implies that $\hbalpha-\balpha_0=o(1)$ for any sequence $\{\hbalpha\}\in\Theta_{\alpha}$.
	\item The maps $\balpha \mapsto \bE\psi(\bV, \balpha, \bbeta, \btheta)$ and $\balpha \mapsto \bE\psi_1(\bV, \balpha)$ are twice continuously differentiable for any $\bbeta\in\Theta_{\beta}$ and $\btheta\in\Theta_{\theta}$. 
	The map $\bbeta\mapsto \bE\psi(\bV, \balpha_0, \bbeta, \btheta)$ is continuously differentiable for any $\btheta\in\Theta_{\theta}$ and $\int_{\btheta\in\Theta_{\theta}}\partial\bE\psi(\bV, \balpha_0, \bbeta, \btheta)/\partial\bbeta w(\btheta)d\btheta=O(1)$.
	\item $\max_{j=1, \cdots, p}\bE\left[\psi^{(j)}(\bV, \balpha, \bzero, \btheta)-\psi^{(j)}(\bV, \balpha', \bzero, \btheta)\right]^2\leq L\|\balpha-\balpha'\|^2$ for any $\balpha\in \Theta_{\alpha}$, $\balpha'\in\Theta_{\alpha}$, and $\btheta\in\Theta_{\theta}$, where $\psi^{(j)}(\bV, \balpha, \bbeta, \btheta)$ is the $j$th component of $\psi(\bV, \balpha, \bbeta, \btheta)$ and $L$ is a constant. $\max_{j=1, \cdots, r}\bE\left[\psi_1^{(j)}(\bV, \balpha)-\psi_1^{(j)}(\bV, \balpha')\right]^2\leq L_{1}\|\balpha-\balpha'\|^2$ for any $\balpha\in \Theta_{\alpha}$, and $\balpha'\in\Theta_{\alpha}$, where $\psi_1^{(j)}(\bV, \balpha)$ is the $j$th component of $\psi_1(\bV, \balpha)$ and $L_1$ is a constant.
	\item The classes of functions $\cF^{(j)}=\{\bV \mapsto \psi^{(j)}(\bV, \balpha, \bbeta, \btheta): \balpha\in \Theta_{\alpha}, \bbeta\in\Theta_{\beta}, \btheta\in\Theta_{\theta}\}$ and $\cF_1^{(j)}=\{\bV \mapsto \psi_1^{(j)}(\bV, \balpha): \balpha\in \Theta_{\alpha}\}$ are pointwise measurable and satisfy the uniform entropy condition $\max_{j=1, \cdots, p}\log\sup_Q N(\eps\|F_j\|_{Q, 2}, \cF^{(j)}, \|\cdot\|_{Q, 2}) \leq C_0\log(e/\eps)$, where $F_{j}=\sup_{f\in \cF^{(j)}}|f|$, $\eps\in(0, 1]$ and $C_0\geq1$ is a constant. The classes $\cF^{(j)}$ and $\cF_1^{(j)}$ have measurable envelope functions $F_j\geq \sup_{f\in \cF^{(j)}}|f|$ satisfying $\bE F_j^s(\bV)\leq C$ and $F_{1j}\geq\sup_{f\in\cF_1^{(j)}}|f|$ satisfying $\bE F_{1j}^s(\bV)\leq C_1$, respectively, where $s\geq 4$.
	
	\item $J=\left[\partial\bE\psi_1(\bV, \balpha_0)/\partial\balpha\trans\right]^{-1}\in\mathbb{R}^{r\times r}$ is a finite and positive-definite deterministic matrix, and $\sup_{\btheta\in\Theta_{\theta}, \bbeta\in\Theta_{\beta}}\|K(\bbeta, \btheta)\|_{\infty}\leq C_K$, $\|J\|_2\leq C_c$, where
	$K(\bbeta, \btheta)=\partial\bE\psi(\bV, \balpha_0, \bbeta, \btheta)/\partial\balpha\in\mathbb{R}^{p\times r}$ is continuous with respect to $\bbeta\in\Theta_{\beta}$. Denote $K(\btheta)=K(\bzero, \btheta)$. $\bE\|\psi_1(\bV, \balpha)\|^2$ and $\bE\|\psi(\bV, \balpha, \bzero, \btheta)\|^2$ are bounded for any $\balpha\in\mathcal{U}_{\alpha_0}$ and $\btheta\in\Theta_{\theta}$, where $\mathcal{U}_{\alpha_0}$ is a neighborhood of $\balpha_0$.

	\item There is a positive function $b(\bv, \bxi)$ of $\bv$ relying on $\balpha_0$, $\btheta_0$ such that
	\begin{align*}
	\left|\bxi\trans\frac{\partial f(\bv; \balpha_0, r_n\bxi, \btheta_0)\partial\bbeta}{f(\bv; \balpha_0, \bzero, \btheta_0)}\right|\leq b(\bv, \bxi),
	\end{align*}
	and $\bE[b(\bV, \bxi)^2]$, $\lambda_{\max}\left(\bE[b(\bV, \bxi)\psi_1(\bV, \balpha_0)^{\otimes2}]\right)$, and $\lambda_{\max}\left(\bE\left[b(\bV, \bxi)\psi(\bV, \balpha_0,\bbeta,\btheta)^{\otimes2}\right]\right)$ 
	are bounded by $C_{\bxi}$, and for all $k$, $\bE[\phi_k(\bV)^2b(\bV, \bxi)]$ is bounded by $C_{\bxi}$, where $\bxi\in\Theta_{\beta}$, $r_n=o(1)$, $C_{\bxi}>0$ is a constant relying on $\bxi$, $\phi_k(\cdot)$ is defined in Theorem~\ref{thm11}, and $\bV$ is generated from the null distribution with density $f(\bv; \balpha_0, \bzero, \btheta_0)$.
	
	\item There is a positive function $b_1(\bv, \bxi_1)$ of $\bv$ relying on $\balpha_0$, $\btheta_0$ such that
	\begin{align*}
	\left|\bxi_1\trans \frac{\partial f(\bv; \balpha_0+r_n\bxi_1, \bzero, \btheta_0)/\partial\balpha}{f(\bv; \balpha_0, \bzero, \btheta_0)}\right|\leq b_1(\bv, \bxi_1),
	\end{align*}
	and $\bE[b_1(\bV, \bxi_1)^2]$, $\lambda_{\max}\left(\bE[b_1(\bV, \bxi_1)\psi_1(\bV, \balpha_0)^{\otimes2}]\right)$, 
	and $\lambda_{\max}\left(\bE\left[b_1(\bV, \bxi)\psi(\bV, \balpha_0,\bbeta,\btheta)^{\otimes2}\right]\right)$
	are bounded by $C_f(\bxi_1)$, and for all $j, k$, $\bE[\phi_k(\bV)^2b_1(\bV, \bxi_1)]$ and $\bE[|\phi_k(\bV)\phi_j(\bV)|b_1(\bV, \bxi_1)]$ are bounded by $C_f(\bxi_1)$, where $\bxi_1\in\Theta_{\alpha}$, $r_n=o(1)$, $C_f(\bxi_1)>0$ is a constant relying on $\bxi_1$, $\phi_k(\cdot)$ is defined in Theorem~\ref{thm11}, and $\bV$ is generated from the null distribution with density $f(\bv; \balpha_0, \bzero, \btheta_0)$. Assume that, for $\bxi_2\in\Theta_{\alpha}$ and any $\tbalpha\in\Theta_{\alpha}$ satisfying $\|\tbalpha-\balpha_0\|^2=O(n^{-1})$,
	\begin{align*}
	\bE\left[\left(\bxi_2\trans\frac{\partial h(\bV_1, \bV_2, \tbalpha)}{\partial\balpha}\right)^2(1+b_1(\bV_1, \bxi_1))(1+b_1(\bV_2, \bxi_1))\right]
	<\infty,
	\end{align*}
	where $h(\bV_1, \bV_2, \balpha)$ equals $h(\bV_1, \bV_2)$ defined in \eqref{kernel0} with $\balpha_0$ replaced by $\balpha$, and $\bV_1$ and $\bV_2$ are independently generated from the null distribution with density $f(\bv, \balpha_0, \bzero, \btheta_0)$.
\end{enumerate}

Assumptions~(A1)--(A5) are mild conditions for establishing the asymptotic properties of statistic $\tT_n$. Assumption~(A1) is an identifiability condition for the GEE; see Theorem~2.10 of \citet{kosorok2008}. Assumptions~(A2) and (A3) are imposed for the asymptotic properties of the estimate of $\balpha$ and the test statistic $\tT_n$. Assumption~(A4) involves entropy conditions and is analogous to the assumption imposed by \citet{2011Testing}. Assumption~(A5) is mild and easily verified in practice, and is used to derive the asymptotic distribution of the proposed statistic $T_n$.
Assumptions~(A6) and (A7) basically impose the Lipschitz condition, which is similar to the Lipschitz condition in Assumption~(A2).

\setcounter{section}{0}\setcounter{equation}{0}\def\theequation{A.\arabic{equation}}\def\thesection{A\arabic{section}}
\def\thesection{Appendix \Alph{section}}

\section{Proof of Theorems}\label{proof_th12}
\def\thelem{A.\arabic{lem}}

We first give some useful notation for convenience of expression. For a vector $\bv\in \bR^{d}$ and a squared matrix $A=(a_{ij})\in\mathbb{R}^{d\times d}$, denote by $\|\bv\|$ the Euclidean norm of $\bv$, by $\trace(A)=\sum_{i=1}^{d}a_{ii}$ the trace of $A$; $\|\bv\|_{A}^2=\sum_{i, j}a_{ij}v_iv_j$ and $\bv^{\otimes2}=\bv\bv\trans$. Denote by $\|A\|_p=\sup\{\|A\bx\|_p: \bx\in\bR^{d}, \|\bx\|_p=1\}$ the induced operator norm for a matrix $A=(a_{ij})\in\bR^{m\times d}$.

Denote by $\rP$ the ordinary probability measure such that $\bE f=\int f d\rP$ for any measurable function $f$, by $\Pn$ the empirical measure of a sample of random elements from $\rP$ such that $\Pn f=n^{-1}\sum_{i=1}^{n}f(\bV_i)$, and by $\Gn$ the empirical process indexed by a class $\cF$ of measurable functions such that $\Gn f=\sqrt{n}(\Pn-\rP)f$ for any $f\in \cF$. Let $L^p(Q)$ be the space of all measurable functions $f$ such that $\|f\|_{Q, p}:=(Q|f|^p)^{1/p}<\infty$, where $p\in[1, \infty)$ and $(Q|f|^p)^{1/p}$ denotes the essential supremum when $p=\infty$. Let $N(\eps, \cF, \|\cdot\|_{Q, 2})$ be an $\eps$-covering number of $\cF$ with respect to the $L^2(Q)$ seminorm $\|\cdot\|_{Q, 2}$, where $\cF$ is a class of measure functions and $Q$ is finite discrete.

\begin{lem}\label{lemma:1a}
	Let $\bX$ be normally distributed with mean $\bmu$ and covariance $\Sigma$. Denote by $\phi(z)$ and $\Phi(z)$ the density function and cumulative distribution function of the standard normal distribution, respectively. For any nonzero vectors $\ba\in\mathbb{R}^p$ and $\bb\in\mathbb{R}^p$ and any numbers $\eta_a\in\mathbb{R}$ and $\eta_b\in\mathbb{R}$, we have
	\begin{align*}	
		\bE[\bone(\ba\trans\bX\geq\eta_a)]=&	\Phi\left(\frac{\ba\trans\bmu-\eta_a}{(\ba\trans\Sigma\ba)^{1/2}}\right),\\ \bE[\bone(\ba\trans\bX\geq\eta_a)\bone(\bb\trans\bX\geq\eta_b)]=&\int_{\frac{\eta_b-\bb\trans\bmu}{(\bb\trans\Sigma\bb)^{1/2}}}^{+\infty}\phi(z)\Phi\left(\frac{1}{\sqrt{1-\rho^2}}\left(\frac{\ba\trans\bmu-\eta_a}{(\ba\trans\Sigma\ba)^{1/2}}+\rho z\right)\right)dz,
	\end{align*}
	where $\rho = \frac{\ba\trans\Sigma\bb\\}{(\ba\trans\Sigma\ba\bb\trans\Sigma\bb)^{1/2}}$, and $\bone(\cdot)$ denotes the indicator function.
	Moreover, when $\ba\trans\bmu-\eta_a=\bb\trans\bmu-\eta_b=0$, we have
	\begin{align*}
		\bE[\bone(\ba\trans\bX\geq\eta_a)]=&\frac{1}{2},\\
		\bE[\bone(\ba\trans\bX\geq\eta_a)\bone(\bb\trans\bX\geq\eta_b)]=&\frac{1}{4}+\frac{1}{2\pi}\arctan\left(\frac{\rho}{\sqrt{1-\rho^2}}\right).
	\end{align*}	
\end{lem}
\noindent{\bf Proof of Lemma \ref{lemma:1a}.} Let $\mu_1 = \ba\trans\bmu$, $\mu_2 = \bb\trans\bmu$, $\sigma_1^2 = \ba\trans\Sigma\ba$, $\sigma_2^2 = \bb\trans\Sigma\bb$.
Noting that the random variables $Z_1 = \ba\trans\bX-\mu_1 \sim N(0, \sigma^2_1)$ and $Z_2 = \bb\trans\bX-\mu_2\sim N(0, \sigma^2_2)$, we have $Z=(Z_1,Z_2)\trans\sim N(\bzero, \Gamma)$ with
\begin{align*}
	\Gamma = \left(
	\begin{array}{cc}
		\ba\trans\Sigma\ba & \ba\trans\Sigma\bb\\
		\ba\trans\Sigma\bb	& \bb\trans\Sigma\bb
	\end{array}
	\right),
\end{align*}
which implies that 
\begin{align*}
	\bE[\bone(\ba\trans\bX\geq\eta_a)]
	=& \bE[\bone(Z_1\geq\eta_a-\mu_1)]\\
	=&\int_{\eta_a-\mu_1}^{+\infty}\frac{1}{\sqrt{2\pi}\sigma_1}\exp\left\{-\frac{z_1^2}{2\sigma_1^2}\right\}dz_1\\
	=&\Phi\left(\frac{\mu_1-\eta_a}{\sigma_1}\right),
\end{align*}
and
\begin{align*}
	&\bE[\bone(\ba\trans\bX\geq\eta_a)\bone(\bb\trans\bX\geq\eta_b)] \\
	=& \bE[\bone(Z_1\geq\eta_a-\mu_1)\bone(Z_2\geq\eta_b-\mu_2)] \\
	=&\int_{\eta_b-\mu_2}^{+\infty}\int_{\eta_a-\mu_1}^{+\infty}\frac{1}{2\pi \sigma_1\sigma_2\sqrt{1-\rho^2}}\exp\left\{-\frac{1}{2}Z\trans\Gamma^{-1}Z\right\}dz_1dz_2\\
	=&\int_{\eta_b-\mu_2}^{+\infty}\int_{\eta_a-\mu_1}^{+\infty}\frac{1}{2\pi \sigma_1\sigma_2\sqrt{1-\rho^2}}\exp\left\{-\frac{1}{2(1-\rho^2)}\left(\frac{z_1}{\sigma_1}-\rho\frac{z_2}{\sigma_2}\right)^2-\frac{z_2^2}{2\sigma_2^2}\right\}dz_1dz_2\\
	=&\int_{\eta_b-\mu_2}^{+\infty}\frac{1}{\sqrt{2\pi}\sigma_2}
	\exp\left\{-\frac{z_2^2}{2\sigma_2^2}\right\}\Phi\left(-\frac{1}{\sqrt{1-\rho^2}}\left(\frac{\eta_a-\mu_1}{\sigma_1}-\frac{\rho}{\sigma_2}z_2\right)\right)dz_2\\
	=&\int_{\frac{\eta_b-\mu_2}{\sigma_2}}^{+\infty}\phi(z)\Phi\left(\frac{1}{\sqrt{1-\rho^2}}\left(\frac{\mu_1-\eta_a}{\sigma_1}+\rho z\right)\right)dz.
\end{align*}

When $\ba\trans\bmu-\eta_a=\bb\trans\bmu-\eta_b=0$, it is straightforward that
\begin{align*}
	\bE[\bone(\ba\trans\bX\geq\eta_a)\bone(\bb\trans\bX\geq\eta_b)]
	=&\int_{0}^{+\infty}\phi(z)\Phi\left(\frac{\rho}{\sqrt{1-\rho^2}}z\right)dz\\
	=&\frac{1}{4}+\frac{1}{2\pi}\arctan\left(\frac{\rho}{\sqrt{1-\rho^2}}\right).
\end{align*} \hfill $\square$

\begin{lem}\label{lem0}
	If Assumptions {\rm (A1)-(A4)} hold, we have
	\begin{align*}
		\hbalpha-\balpha_0
		=& - J\Pn\psi_1(\bV, \balpha_0)+o_p(n^{-1/2}),
	\end{align*}
	where $\balpha_0$ is given in Assumption (A1), and $J=\left(\partial\rP\psi_1(\bV, \balpha_0)/\partial\balpha\trans\right)^{-1}\in\mathbb{R}^{r\times r}$.	
\end{lem}
\noindent{\bf Proof of Lemma \ref{lem0}.}
Note that
\begin{align*}
	\rP[\psi_1(\bV, \hbalpha)- \psi_1(\bV, \balpha_0)]
    =& \Pn\psi_1(\bV, \hbalpha) -\Pn\psi_1(\bV, \balpha_0)\\
	 &- (\Pn-\rP)\left[\psi_1(\bV, \hbalpha)-\psi_1(\bV, \balpha_0)\right].
\end{align*}
It suffices to prove that
\begin{align}\label{B01}
	(\Pn-\rP)\left[\psi_1(\bV, \hbalpha)-\psi_1(\bV, \balpha_0)\right] = o_p(n^{-1/2}),
\end{align}
and
\begin{align}\label{B02}
	\rP[\psi_1(\bV, \hbalpha)- \psi_1(\bV, \balpha_0)] = J^{-1}(\hbalpha-\balpha_0)+o(n^{-1/2}).
\end{align}
	
	Let $\mathcal{U}_{\alpha_0}=\{\balpha\in\Theta_{\alpha}:\|\balpha-\balpha_0\|=o(1)\}$ and $\cF_{1\alpha}=\cup_{1\leq j\leq r}\cF_{1\alpha}^{(j)}$, where $\cF_{1\alpha}^{(j)}=\{\psi_1^{(j)}(\bV, \balpha)-\psi_1^{(j)}(\bV, \balpha_0):\balpha\in\mathcal{U}_{\alpha_0}\}$.
        By Assumptions (A3)-(A4) and Theorem 5.1 and 5.2 of \cite{chernozhukov2014gaussian} with $\sigma^2=\|F_{1\alpha}\|_{\rP,2}^2$ and $M_n=\max_{1\leq i \leq n}F_{1\alpha}(\bV_i)$, we have that, with probability $1-o(1)$, 
	\begin{align*}
		\sup_{f\in\cF_{1\alpha}}\|\Gn f\|\lesssim 2\sqrt{C_0r}L_{1}^{1/2}\|\balpha-\balpha_0\|+4C_0rM_nn^{-1/2}=o(1),
	\end{align*}
	which implies that (\ref{B01}) holds.

	It follows from the central limit theorem
	\begin{align}\label{B03}
		\Pn\psi_1(\bV, \balpha_0) = O_p(n^{-1/2}),
	\end{align}
	which combining (\ref{B01}) indicates that $\rP[\psi_1(\bV, \hbalpha)- \psi_1(\bV, \balpha_0)] = O(n^{-1/2})$ with probability $1-o(1)$, and consequently, we have, by Assumption (A1), $\hbalpha-\balpha_0=o(1)$ with probability $1-o(1)$.
	
	For any $\balpha\in \mathcal{U}_{\alpha_0}$, we have, by Assumption (A2) and Taylor's theorem,
	\begin{align*}
		\rP[\psi_1(\bV, \balpha)- \psi_1(\bV, \balpha_0)]
		= \frac{\partial\rP\psi_1(\bV, \balpha_0)}{\partial\balpha\trans}(\balpha-\balpha_0) + O(\|\balpha-\balpha_0\|^2),
	\end{align*}
	which indicates that, with probability $1-o(1)$,
	\begin{align*}
		\rP[\psi_1(\bV, \hbalpha)- \psi_1(\bV, \balpha_0)]
		= J^{-1}(\hbalpha-\balpha_0)(1+o(1)).
	\end{align*}
	Thus, combining (\ref{B01}) and (\ref{B03}) we have $\hbalpha-\balpha_0=O_p(n^{-1/2})$, and that (\ref{B02}) holds. \hfill $\square$

    \begin{lem}\label{lems1}
        If Assumptions~{(A1)--(A5)} hold, then for any $\bbeta\in\Theta_{\beta}$ and $\btheta\in\Theta_{\theta}$, we have
        \begin{align*}
        \Pn\psi(\bV, \hbalpha, \bbeta, \btheta)
        =&\Pn\psi(\bV, \balpha_0, \bbeta, \btheta) - K(\bbeta, \btheta)J\Pn\psi_{1}(\bV, \balpha_0)+o_p(n^{-1/2}),
        \end{align*}
        where $\balpha_0$ is given in Assumption~(A1), $K(\bbeta, \btheta)=\partial\rP\psi(\bV, \balpha_0, \bbeta, \btheta)/\partial\balpha\in\mathbb{R}^{p\times r}$, and $J=\{\partial\rP\psi_1(\bV, \balpha_0)/\partial\balpha\}^{-1}\in\mathbb{R}^{r\times r}$.
        \end{lem}
        \noindent{\bf Proof of Lemma \ref{lems1}.}
            Note that
            \begin{multline*}
                \rP[\psi(\bV, \hbalpha, \bbeta, \btheta)- \psi(\bV, \balpha_0, \bbeta, \btheta)] = \Pn\psi(\bV, \hbalpha, \bbeta, \btheta) -\Pn\psi(\bV, \balpha_0, \bbeta, \btheta)\\
                - (\Pn-\rP)\left[\psi(\bV, \hbalpha, \bbeta, \btheta)-\psi(\bV, \balpha_0, \bbeta, \btheta)\right].
            \end{multline*}
            It suffices to prove that
            \begin{align}\label{A01}
                (\Pn-\rP)\left[\psi(\bV, \hbalpha, \bbeta, \btheta)-\psi(\bV, \balpha_0, \bbeta, \btheta)\right] = o_p(n^{-1/2}),
            \end{align}
            and
            \begin{align}\label{A02}
                \rP[\psi(\bV, \hbalpha, \bbeta, \btheta)- \psi(\bV, \balpha_0, \bbeta, \btheta)] = - K(\btheta, \bbeta)J\Pn\psi_{1}(\bV, \balpha_0)+o(n^{-1/2}).
            \end{align}
            
            We first show that (\ref{A01}) holds.
            Let $\mathcal{U}_{\alpha_0}=\{\balpha\in\Theta_{\alpha}:\|\balpha-\balpha_0\|=O(n^{-1/2})\}$ and $\cF_{\alpha}=\cup_{1\leq j\leq p}\cF_{\alpha}^{(j)}$, where $\cF_{\alpha}^{(j)}=\{\psi^{(j)}(\bV, \balpha, \bbeta, \btheta)-\psi^{(j)}(\bV, \balpha_0, \bbeta, \btheta):\balpha\in\mathcal{U}_{\alpha_0},\bbeta\in\Theta_{\beta}, \btheta\in\Theta_{\theta}\}$. We define the envelope function of $\cF_{\alpha}$ as
            \begin{align*}
                F_{\alpha} = \max_{j=1,\cdots,p}F_{\alpha}^{(j)},
            \end{align*}
            where
            \begin{align*}
                F_{\alpha}^{(j)} = \sup_{\balpha\in\mathcal{U}_{\alpha_0}, \btheta\in\Theta_{\theta}}|\psi^{(j)}(\bV, \balpha, \bbeta, \btheta)-\psi^{(j)}(\bV, \balpha_0, \bbeta, \btheta)|.
            \end{align*}
            Assumption (A3) implies that
            \begin{align*}
                \|F_{\alpha}\|_{\rP, 2} \leq \|F\|_{\rP,2}\leq L^{1/2}\|\balpha-\balpha_0\|,
            \end{align*}
            and by the fact that $\cF_{\alpha}^{(j)}\subseteq \cF^{(j)}-\cF^{(j)}$,
            \begin{align*}
                \log\sup_Q N(\eps\|F_{\alpha}\|_{Q,2}, \cF_{\alpha}, \|\cdot\|_{Q,2})
                \leq& 2p \max_{j=1,\dots,p}\log\sup_Q N\left(\frac{\eps}{2}\|F_j\|_{Q,2}, \cF^{(j)}, \|\cdot\|_{Q,2}\right)\\
                \leq& 2C_0p\log(2e/\eps),
            \end{align*}
            which results in that, for $\delta\in (0,1)$,
            \begin{align*}
                \bJ(\delta,\cF_{\alpha}, F_{\alpha})
                =& \int_{0}^{\delta}\sup_Q\sqrt{1+\log N(\eps\|F_{\alpha}\|_{Q,2}, \cF_{\alpha}, \|\cdot\|_{Q,2})}d\eps\\
                \leq&   \int_{0}^{\delta}\sqrt{1+2C_0p\log(2e/\eps)}d\eps\\
                =&2e\int_{2e/\delta}^{\infty}\frac{1}{\eps^2}\sqrt{1+2C_0p\log(\eps)}d\eps.
            \end{align*}
            According to the fact that $1+\log(\eps)\geq1$ for all $\eps\in[1, \infty)$ and
            the integration by parts, we have
            \begin{align*}
                \int_{2e/\delta}^{\infty}\frac{1}{\eps^2}\sqrt{1+\log(\eps)}d\eps
                =&-\frac{1}{\eps}\sqrt{1+\log(\eps)}\bigg|_{\eps=2e/\delta}^{\infty}+\frac{1}{2}\int_{2e/\delta}^{\infty}\frac{1}{\eps^2\sqrt{1+\log(\eps)}}d\eps\\
                \leq&\frac{\delta}{2e}\sqrt{1+\log(2e/\delta)}+\frac{1}{2}\int_{2e/\delta}^{\infty}\frac{1}{\eps^2}\sqrt{1+\log(\eps)}d\eps,
            \end{align*}
            which implies that
            \begin{align*}
                \int_{2e/\delta}^{\infty}\frac{1}{\eps^2}\sqrt{1+\log(\eps)}d\eps
                \leq\frac{\delta}{e}\sqrt{1+\log(2e/\delta)},
            \end{align*}
            and consequently
            \begin{align*}
                \bJ(\delta,\cF_{\alpha}, F_{\alpha}) \leq \delta \sqrt{2C_0p(1+\log(2e/\delta))}.
            \end{align*}
            Therefore, by Assumption (A4) and Theorem 5.1 and 5.2 of \cite{chernozhukov2014gaussian} with $\sigma^2=\|F_{\alpha}\|_{\rP,2}^2$ and $M_n=\max_{1\leq i \leq n}F_{\alpha}(\bV_i)$, we have that, with probability $1-o(1)$,
            \begin{align*}
                \sup_{f\in\cF_{\alpha}}\|\Gn f\|\lesssim 2\sqrt{C_0p}L^{1/2}\|\balpha-\balpha_0\|+4C_0pM_nn^{-1/2}=o(1),
            \end{align*}
            which implies that, with probability $1-o(1)$,
            \begin{align*}
                \sup_{f\in\cF_{\alpha}}\|(\Pn-\rP)f\| = o(n^{-1/2})
            \end{align*}
            because of the fact that $\hbalpha \in \mathcal{U}_{\alpha_0}$ with probability $1-o(1)$ by Lemma \ref{lem0}, and consequently (\ref{A01}) holds.
            
            We, now show (\ref{A02}) holds. For any $\balpha\in \mathcal{U}_{\alpha_0}$, we have, by Assumption (A2) and Taylor's theorem,
            \begin{align*}
                \rP[\psi(\bV, \balpha, \bbeta, \btheta)- \psi(\bV, \balpha_0, \bbeta, \btheta)]
                = \frac{\partial\rP\psi(\bV, \balpha_0, \bbeta, \btheta)}{\partial\balpha\trans}(\balpha-\balpha_0) + O(\|\balpha-\balpha_0\|^2),
            \end{align*}
            which combining Lemma \ref{lem0} imply that, with probability $1-o(1)$,
            \begin{align*}
                \rP[\psi(\bV, \hbalpha, \bbeta, \btheta)- \psi(\bV, \balpha_0, \bbeta, \btheta)]
                = - K(\btheta,\bbeta)J\Pn\psi_1(\bV,\alpha_0) + O(n^{-1}),
            \end{align*}
            and consequently (\ref{A02}) holds. \hfill $\square$

            \noindent{\bf Proof of Lemma 1.}
By applying the central limit theorem to i.i.d. sample $\{\bV_i\}_{i=1}^n$, we have that under $H_0$
	\begin{align*}
		n^{1/2}\left[\Pn\psi(\bV, \balpha_0, \bzero, \btheta)	- K(\btheta)J\Pn\psi_{1}(\bV, \balpha_0)\right]\lkonv N(\bzero, \Gamma(\btheta)),	
	\end{align*}
	for any fixed $\btheta\in\Theta_{\theta}$, where $\Gamma(\btheta)=\rP\{\psi(\bV, \balpha_0, \bzero, \btheta) - K(\btheta)J\psi_{1}(\bV, \balpha_0)\}^{\otimes2}$. 
 By Lemma \ref{lems1}, 
 it is straightforward that for any fixed $\btheta\in\Theta_{\theta}$,
	\begin{align*}
		\tT_n(\btheta)
		= \|\Pn\psi(\bV, \balpha_0, \bzero, \btheta) - K(\btheta)J\Pn\psi_{1}(\bV, \balpha_0)\|_{\tV(\btheta)^{-1}}^2+o_p(1)
		\lkonv\chi_p^2.
	\end{align*}
	To complete the proof, it suffices to show that $\tV(\btheta)\pkonv \Gamma(\btheta)$ for any fixed $\btheta\in\Theta_{\theta}$. In fact, it is easy to show that $\tV(\btheta)\pkonv \Gamma(\btheta)$ for any fixed $\btheta\in\Theta_{\theta}$ by the law of large numbers.	\hfill $\square$

    \begin{lem}\label{lem2}
        If Assumptions~{(A1)--(A5)} hold, then $\tT_n$ converges in distribution to $\sup_{\btheta\in\Theta_{\theta}}\cG^2(\btheta)$ under $H_0$ as $n \rightarrow \infty$, where $\{\cG(\btheta); \btheta\in\Theta_{\theta}\}$ is a zero-mean Gaussian process with the covariance function
        \begin{align*}
        \Sigma(\btheta_1, \btheta_2) = \rP\psi_*(\bV, \balpha_0, \bzero, \btheta_1)\psi_*(\bV, \balpha_0, \bzero, \btheta_2)\trans,
        \end{align*}
        indexed by $\btheta_1\in\Theta_{\theta}$ and $\btheta_2\in\Theta_{\theta}$, where $\psi_*(\bV, \balpha_0, \bzero, \btheta)=\Gamma(\btheta)^{-1/2}\{\psi(\bV, \balpha_0, \bzero, \btheta) - K(\btheta)J\psi_{1}(\bV, \balpha_0)\}$
        and $\Gamma(\btheta)=\rP\{\psi(\bV, \balpha_0, \bzero, \btheta) - K(\btheta)J\psi_{1}(\bV, \balpha_0)\}^{\otimes2}$, with $K(\btheta)$ and $J$ defined in Assumption~(A5).
        \end{lem}
        \noindent{\bf Proof of Lemma \ref{lem2}.}
        
            Define a class of measurable functions $\cF_{2}=\cup_{1\leq j\leq p}\cF_2^{(j)} $ and its envelope function $F_{2}=\max_{j=1,\cdots,p}\sup_{f\in\cF_{2}^{(j)}}|f|$, where $\cF_2^{(j)}=\{\bV \mapsto \psi_*^{(j)}(\bV, \balpha_0, \bzero, \btheta): \bbeta\in\Theta_{\beta}, \btheta\in\Theta_{\theta}\}$. Assumption (A2) and (A4) imply that
            \begin{align*}
                \rP F_{2}^2(\bV)\leq (C+C_K^2C_c^2C_1).
            \end{align*}
            and by the fact that $\cF_{2}^{(j)}\subseteq \cF^{(j)}+\underbrace{\cF_1^{(j)}+\cdots+\cF_1^{(j)}}_{r}$, and
            \begin{align*}
                \log\sup_Q N(\eps\|F_{2}\|_{Q,2}, \cF_{2}, \|\cdot\|_{Q,2})
                \leq& p\max_{j=1,\dots,p}\log\sup_Q N\left(\frac{\eps}{2}\|F_j\|_{Q,2}, \cF^{(j)}, \|\cdot\|_{Q,2}\right)\\
                & +rp\max_{j=1,\dots,p}\log\sup_Q N\left(\frac{\eps}{2}\|F_{1j}\|_{Q,2}, \cF_1^{(j)}, \|\cdot\|_{Q,2}\right)\\
                \leq& 2C_0rp\log(2e/\eps),
            \end{align*}
            which results in that, 
         along the line of Lemma \ref{lems1}, 
         for $\delta\in (0,1)$,
            \begin{align*}
                \bJ(\delta,\cF_{2}, F_{2})
                =& \int_{0}^{\delta}\sup_Q\sqrt{1+\log N(\eps\|F_{2}\|_{Q,2}, \cF_{2}, \|\cdot\|_{Q,2})}d\eps\\
                \leq& \delta \sqrt{2C_0rp(1+\log(2e/\delta))}< \infty.
            \end{align*}
            By Theorem 2.5.2 of \cite{van1996weak}, $\cF_2$ is $\rP$-Donsker, which means that $\Gn\psi_*(\bV, \balpha_0, \bzero, \btheta)$ converges weekly to a mean zero Gaussian process with the covariance
            \begin{align*}
                \Sigma(\btheta_1, \btheta_2) = \rP\psi_*(\bV, \balpha_0, \bzero, \btheta_1)\psi_*(\bV, \balpha_0, \bzero, \btheta_2)\trans.
            \end{align*}
         By Lemma \ref{lems1}, 
         it remains to show that $\tV(\btheta)\pkonv \Gamma(\btheta)$ for any fixed $\btheta\in\Theta_{\theta}$. In fact, it is seen that $\tV(\btheta)\pkonv \Gamma(\btheta)$ for any fixed $\btheta\in\Theta_{\theta}$ by the law of large numbers. \hfill $\square$
   
         \begin{lem}\label{lem3}
            If Assumptions~{(A1)--(A5)} hold, then $\tT_n$ converges in distribution to $\sup_{\btheta\in\Theta_{\theta}}\cG_{r_n}^2(\btheta)$ under $H_{1n}$ as $n \rightarrow \infty$, where $\{\cG_{r_n}(\btheta); \btheta\in\Theta_{\theta}\}$ is a Gaussian process with the mean function
            \begin{align*}
            \mu(\btheta) = \Gamma(\btheta)^{-1/2}\frac{\partial\rP_0\psi(\bV, \balpha_0, \bzero, \btheta)}{\partial\bbeta\trans}\bxi
            \end{align*}
            and the covariance function $\Sigma(\btheta_1, \btheta_2)$ as defined in Lemma \ref{lem2}. Here, $\rP_0$ is the probability measure under the null hypothesis.
            \end{lem}
            \noindent{\bf Proof of Lemma \ref{lem3}.}
            
                Under the local alternative, as in Lemma \ref{lem2}, we can show that
                \begin{align*}
                    (\Pn-\rP)\left[\psi(\bV, \hbalpha, \bbeta, \btheta)-\psi(\bV, \balpha_0, \bbeta, \btheta)\right] = o_p(n^{-1/2}),
                \end{align*}
                and consequently,
                \begin{align}\label{A21}
                \begin{split}
                    \Pn\psi(\bV, \hbalpha, \bzero, \btheta)
                =& \Pn\psi(\bV, \balpha_0, \bzero, \btheta)\\
                 &+\rP[\psi(\bV, \hbalpha, \bzero, \btheta)-\psi(\bV, \balpha_0, \bzero, \btheta)]+o_p(n^{-1/2}).
                \end{split}
                \end{align}
                By the differentiability of $\rP\psi(\bV, \balpha, \bbeta, \btheta)$ with respect to $\bbeta$ in Assumption (A2), we have, for any $\balpha\in\Theta_{\alpha}$ and $\btheta\in\Theta_{\theta}$,
                \begin{align*}
                    \rP\psi(\bV, \balpha, \bbeta, \btheta) =\rP\psi(\bV, \balpha, \bzero, \btheta)+\bbeta\trans\frac{\partial \rP\psi(\bV, \balpha, \bzero, \btheta)}{\partial\bbeta} + O(\|\bbeta\|^2).
                \end{align*}
                Thus, for any $\balpha$ satisfying $\balpha-\balpha_0=O(n^{-1/2})$, we have, by Assumption (A2) and Taylor's theorem,
                \begin{align*}
                    \rP[\psi(\bV, \balpha, \bzero, \btheta)- \psi(\bV, \balpha_0, \bzero, \btheta)]
                    =&\rP\psi(\bV, \balpha, \bbeta, \btheta) -\bbeta\trans\frac{\partial \rP\psi(\bV, \balpha, \bzero, \btheta)}{\partial\bbeta}\\
                     &-\rP\psi(\bV, \balpha_0, \bbeta, \btheta)+\bbeta\trans\frac{\partial \rP\psi(\bV, \balpha_0, \bzero, \btheta)}{\partial\bbeta}+ O(\|\bbeta\|^2)\\
                    =& K(\bbeta,\btheta)(\balpha-\balpha_0)+\bbeta\trans\frac{\partial^2\rP[\psi(\bV, \balpha_0, \bzero, \btheta)}{\partial\bbeta\partial\balpha\trans}(\balpha-\balpha_0)\\
                     &+ O(\|\balpha-\balpha_0\|^2)+O(\|\bbeta\|^2),
                \end{align*}
                which combining Lemma \ref{lem0} implies that, with probability $1-o(1)$,
                \begin{align*}
                    \rP[\psi(\bV, \hbalpha, \bzero, \btheta)- \psi(\bV, \balpha_0, \bzero, \btheta)]
                    =& K(\bbeta,\btheta)J\Pn\psi_1(\bV,\alpha_0) + O(n^{-1}),
                \end{align*}
                which implies that
                \begin{align*}
                    \Pn\psi(\bV, \hbalpha, \bzero, \btheta)
                    =&\Pn\{\psi(\bV, \balpha_0, \bzero, \btheta)+ K(\bbeta,\btheta)J\psi_{1}(\bV, \balpha_0)\}+o_p(n^{-1/2})\\
                    =&\Pn\{\psi(\bV, \balpha_0, \bbeta, \btheta)+ K(\bbeta,\btheta)J\psi_{1}(\bV, \balpha_0)\}\\
                     &+\Pn\{\psi(\bV, \balpha_0, \bzero, \btheta)-\psi(\bV, \balpha_0, \bbeta, \btheta)\}+o_p(n^{-1/2}).
                \end{align*}
                Let $\rP_0$ denote the probability measure under the null. Since under the local alternative,
                \begin{align*}
                \mu(\btheta)
                =&\Gamma(\btheta)^{-1/2}\frac{\partial\rP_0\psi(\bV, \balpha_0, \bzero, \btheta)}{\partial\bbeta\trans}\bxi,
                \end{align*}
                we have, along the line of Lemma \ref{lems1}, 
                $\Gn\psi_*(\bV, \balpha_0, \bzero, \btheta)$ converges weekly to a Gaussian process with mean function $\mu(\btheta)$ and the covariance
                \begin{align*}
                    \Sigma(\btheta_1, \btheta_2) = \rP\psi_*(\bV, \balpha_0, \bzero, \btheta_1)\psi_*(\bV, \balpha_0, \bzero, \btheta_2)\trans.
                \end{align*}
                Therefore, we complete the proof by (\ref{A21}). \hfill $\square$

                \begin{lem}\label{lem:H}
                    Let \begin{align*}
                        H = \int_{\btheta\in\Theta_\theta}J\trans K(\btheta)\trans K(\btheta)J w(\btheta) d\btheta,
                    \end{align*}
                    where $J$ and $K(\btheta)$ are defined in Assumption (A5), and $w(\btheta)$ is a weight satisfying $w(\btheta)\geq0$ for all $\btheta\in\Theta_{\theta}$ and $\int_{\btheta\in\Theta_{\theta}} w(\btheta)d\btheta=1$.
                    If Assumption (A5) holds, we have
                    \begin{align*}
                        \lambda_{\max}(H)\leq\lambda_{\max}^2(J)\lambda_{\max}\left[\int_{\btheta\in\Theta_\theta} \Sigma_{K}(\btheta) w(\btheta) d\btheta\right],
                    \end{align*}
                    and
                    \begin{align*}
                        \trace(H)\leq \lambda_{\max}^2(J)\trace\left(\int_{\btheta\in\Theta_\theta} \Sigma_{K}(\btheta) w(\btheta) d\btheta\right),
                    \end{align*}
                    where $\Sigma_{K}(\btheta)=K(\btheta)\trans K(\btheta)$.
                    Furthermore, by Assumption (A5), $\lambda_{\max}(H)$ and $\trace(H)$ are bounded.
                \end{lem}
                \noindent{\bf Proof of Lemma \ref{lem:H}.} It is straightforward to complete the proof by utilizing Assumption (A5). \hfill $\square$
 
                \begin{lem}\label{lem:h2}
                    If Assumption (A4) holds, we have for any fixed $\bxi\in\Theta_{\beta}$,
                    \begin{align*}
                        \rP[h(\bV_1,\bV_2)^2]	\leq&4\int_{\btheta\in\Theta_\theta}\lambda_{\max}(\Sigma_\psi(\btheta))\trace(\Sigma_\psi(\btheta))w(\btheta)d\btheta\\
                        &+4\lambda_{\max}^2(\Sigma_{\psi_1})\trace(H^2)\\
                        &+8\lambda_{\max}(\Sigma_{\psi_1})\lambda_{\max}^2(J)\int_{\btheta\in\Theta_\theta}\lambda_{\max}(\Sigma_\psi)\trace\left(\Sigma_{K}(\btheta)\right)w(\btheta)d\btheta,
                    \end{align*}
                    where $h(\bV_1,\bV_2)$ is defined in (11) in the main paper, and $\bV_1$ and $\bV_2$ are generated from the null distribution with density $f(\bv,\balpha_0,\bzero,\btheta_0)$ and are independent of each other.
                    Furthermore, if Assumption (A5) holds, $\bE[h(\bV_1,\bV_2)^2]$ is bounded.
                    
                    Here we define matrices as following,
                    \begin{align*}
                        \Sigma_\psi(\btheta)=&\rP\left[\psi(\bV, \balpha_0, \bzero, \btheta)^{\otimes2}\right],\\
                        \Sigma_{\psi_1}(\btheta)=&\rP\left[\psi_1(\bV, \balpha_0)^{\otimes2}\right],\\
                        \Sigma_{K}(\btheta)=&K(\btheta)\trans K(\btheta).
                    \end{align*}
                \end{lem}
                \noindent{\bf Proof of Lemma \ref{lem:h2}.} We deal with each term in $h(\bV_1,\bV_2)$ separately. It is seen that
                \begin{align*}
                    &\rP\left[\left(\int_{\btheta\in\Theta_\theta}\psi(\bV_i, \balpha_0, \bzero,\btheta)\trans \psi(\bV_j, \balpha_0, \bzero,\btheta) w(\btheta)d\btheta\right)^2\right]\\
                    \leq&\int_{\btheta\in\Theta_\theta}\rP\left[\left(\psi(\bV_i, \balpha_0, \bzero,\btheta)\trans \psi(\bV_j, \balpha_0, \bzero,\btheta)\right)^2\right] w(\btheta)d\btheta\\
                    \leq&\int_{\btheta\in\Theta_\theta}\lambda_{\max}(\Sigma_\psi(\btheta))\rP\left[\left\|\psi(\bV_i, \balpha_0, \bzero,\btheta)\right\|^2\right]w(\btheta)d\btheta\\
                    \leq&\int_{\btheta\in\Theta_\theta}\lambda_{\max}(\Sigma_\psi(\btheta))\trace(\Sigma_\psi(\btheta))w(\btheta)d\btheta.
                \end{align*}
                By the definition of $K_i$, we have
                \begin{align*}
                    \rP\|K_i\|^2
                    \leq&\int_{\btheta\in\Theta_\theta}\trace\left\{\rP\left[ J\trans K(\btheta)\trans \psi(\bV_i, \balpha_0, \bzero,\btheta)\psi(\bV_i, \balpha_0, \bzero,\btheta)\trans K(\btheta)J\right]\right\} w(\btheta)d\btheta\\
                    \leq&\int_{\btheta\in\Theta_\theta}\lambda_{\max}(\Sigma_\psi(\btheta))\trace\left(J\trans K(\btheta)\trans K(\btheta)J\right) w(\btheta)d\btheta\\
                    =&\lambda_{\max}^2(J)\int_{\btheta\in\Theta_\theta}\lambda_{\max}(\Sigma_\psi)\trace\left(\Sigma_{K}(\btheta)\right)w(\btheta)d\btheta,
                \end{align*}
                and consequently,
                \begin{align*}
                    \rP\left[(\psi_1(\bV_i, \balpha_0)\trans K_{j})^2\right]
                    =&\rP\left[\left(K_{j}\trans \Sigma_{\psi_1}K_{j}\right)\right] \\
                    \leq&\lambda_{\max}(\Sigma_{\psi_1})\lambda_{\max}^2(J)\int_{\btheta\in\Theta_\theta}\lambda_{\max}(\Sigma_\psi)\trace\left(\Sigma_{K}(\btheta)\right)w(\btheta)d\btheta.
                \end{align*}
                And similarly, it can be obtained that
                \begin{align*}
                    \rP\left[\left(\psi_1(\bV_i, \balpha_0)\trans H \psi_1(\bV_j, \balpha_0)\right)^2\right]
                    \leq&\lambda_{\max}^2(\Sigma_{\psi_1})\trace(H^2).
                \end{align*}
                From the foregoing considerations, it is seen that
                \begin{align*}
                    \rP[h(\bV_1,\bV_2)^2]
                    \leq&4\rP\tpsi_{ij}(\balpha_0)+4\rP\left[\psi_1(\bV_i, \balpha_0)\trans H\psi_1(\bV_j, \balpha_0)\right]^2
                    +8\rP\left[\psi_1(\bV_i, \balpha_0)\trans K_{j}\right]^2\\
                    \leq&4\int_{\btheta\in\Theta_\theta}\lambda_{\max}(\Sigma_\psi(\btheta))\trace(\Sigma_\psi(\btheta))w(\btheta)d\btheta\\
                    &+4\lambda_{\max}^2(\Sigma_{\psi_1})\trace(H^2)\\
                    &+8\lambda_{\max}(\Sigma_{\psi_1})\lambda_{\max}^2(J)\int_{\btheta\in\Theta_\theta}\lambda_{\max}(\Sigma_\psi)\trace\left(\Sigma_{K}(\btheta)\right)w(\btheta)d\btheta.
                \end{align*}
                By Assumption (A5) and Lemma \ref{lem:H}, it is straightforward to see that $\bE[h(\bV_1,\bV_2)^2]$ is bounded. \hfill $\square$
 
                \begin{lem}\label{lem:h2b2}
                    If Assumption (A4) -- (A6) holds, we have for any fixed $\bxi\in\Theta_{\beta}$,
                    \begin{align*}
                        \rP[h(\bV_1,\bV_2)^2(1+b(\bV_1,\bxi))(1+b(\bV_2,\bxi))]
                    \end{align*}
                    is bounded,
                    where $b(\bv,\bxi)$ is defined in Assumption (A6), $h(\bV_1,\bV_2)$ is defined in (11) in the main paper, and $\bV_1$ and $\bV_2$ are generated from the null distribution with density $f(\bv,\balpha_0,\bzero,\btheta_0)$ and are independent of each other.
                    
                \end{lem}
                \noindent{\bf Proof of Lemma \ref{lem:h2b2}.} 	We defined $\Sigma_{\psi,b}$ and $\Sigma_{\psi_1,b}$ as following,
                \begin{align*}
                    \Sigma_{\psi,b}=&\rP\left[b(\bV_i,\bxi)\psi(\bV, \balpha_0, \bzero, \btheta)^{\otimes2}\right],\\
                    \Sigma_{\psi_1,b}=&\rP\left[b(\bV_i,\bxi)\psi_1(\bV_i, \balpha_0)^{\otimes2}\right].
                \end{align*}
                Along the lines of Lemma \ref{lem:h2}, we have by Assumption (A6),
                \begin{align*}
                    &\rP\left[\left(\int_{\btheta\in\Theta_\theta}\psi(\bV_i, \balpha_0, \bzero,\btheta)\trans \psi(\bV_j, \balpha_0, \bzero,\btheta) w(\btheta)d\btheta\right)^2b(\bV_i,\bxi)\right]\\
                    \leq&\int_{\btheta\in\Theta_\theta}\rP\left[\left(\psi(\bV_i, \balpha_0, \bzero,\btheta)\trans \psi(\bV_j, \balpha_0, \bzero,\btheta)\right)^2b(\bV_i,\bxi)\right] w(\btheta)d\btheta\\
                    \leq&\int_{\btheta\in\Theta_\theta}\lambda_{\max}(\Sigma_{\psi})\rP\left[\left\|\psi(\bV_i, \balpha_0, \bzero,\btheta)\right\|^2b(\bV_i,\bxi)\right]w(\btheta)d\btheta\\
                    \leq&\int_{\btheta\in\Theta_\theta}\lambda_{\max}(\Sigma_{\psi}(\btheta))\trace(\Sigma_{\psi,b}(\btheta)) w(\btheta)d\btheta.
                \end{align*}
                By the definition of $K_i$, we have
                \begin{align*}
                    \rP&\left[\|K_i\|^2b(\bV_i,\bxi)\right]\\
                    \leq&\int_{\btheta\in\Theta_\theta}\trace\left\{\rP\left[ J\trans K(\btheta)\trans \psi(\bV_i, \balpha_0, \bzero,\btheta)\psi(\bV_i, \balpha_0, \bzero,\btheta)\trans K(\btheta)Jb(\bV_i,\bxi)\right]\right\} w(\btheta)d\btheta\\
                    \leq&\int_{\btheta\in\Theta_\theta}\lambda_{\max}(\Sigma_{\psi,b}(\btheta))\trace\left(J\trans K(\btheta)\trans K(\btheta)J\right) w(\btheta)d\btheta\\
                    =&\lambda_{\max}^2(J)\int_{\btheta\in\Theta_\theta}\lambda_{\max}(\Sigma_{\psi,b}(\btheta))\trace\left(\Sigma_K(\btheta)\right) w(\btheta)d\btheta,
                \end{align*}
                and consequently,
                \begin{align*}
                    \rP&\left[(\psi_1(\bV_i, \balpha_0)\trans K_{j})^2b(\bV_i,\bxi)\right]\\
                    =&\rP\left[K_{j}\trans \Sigma_{\psi_1,b}K_{j}\right] \\
                    \leq&\lambda_{\max}(\Sigma_{\psi_1,b})\lambda_{\max}^2(J)\int_{\btheta\in\Theta_\theta}\lambda_{\max}(\Sigma_{\psi}(\btheta))\trace\left(\Sigma_K(\btheta)\right) w(\btheta)d\btheta,
                \end{align*}
                and
                \begin{align*}
                    \rP&\left[(\psi_1(\bV_j, \balpha_0)\trans K_{i})^2b(\bV_i,\bxi)\right]\\
                    =&\rP\left[K_{i}\trans \Sigma_{\psi_1}K_{i}b(\bV_i,\bxi)\right] \\
                    \leq&\lambda_{\max}(\Sigma_{\psi_1})\lambda_{\max}^2(J)\int_{\btheta\in\Theta_\theta}\lambda_{\max}(\Sigma_{\psi,b}(\btheta))\trace\left(\Sigma_K(\btheta)\right) w(\btheta)d\btheta.
                \end{align*}
                And similarly, it can be obtained that
                \begin{align*}
                    \rP\left[\left(\psi_1(\bV_i, \balpha_0)\trans H\psi_1(\bV_j, \balpha_0)\right)^2b(\bV_i,\bxi)\right]
                    \leq&\lambda_{\max}(\Sigma_{\psi_1,b})\lambda_{\max}(\Sigma_{\psi_1})\trace(H^2).
                \end{align*}
                From the foregoing considerations, it is seen that by Assumption (A7),
                \begin{align}\label{eq:h2b}
                    \begin{split}
                        \rP&[h(\bV_i,\bV_j)^2b(\bV_i,\bxi)]\\
                        \leq&4\rP\left[\tpsi_{ij}(\balpha_0)b(\bV_i,\bxi)\right]\\
                        &+4\rP\left[(\psi_1(\bV_i, \balpha_0)\trans H\psi_1(\bV_j, \balpha_0))^2 b(\bV_i,\bxi)\right]\\
                        &+4\rP\left[(\psi_1(\bV_i, \balpha_0)\trans K_{j})^2 b(\bV_i,\bxi)\right]+4\rP\left[(\psi_1(\bV_j, \balpha_0)\trans K_{i})^2 b(\bV_i,\bxi)\right]\\
                        \leq&4\int_{\btheta\in\Theta_\theta}\lambda_{\max}(\Sigma_{\psi}(\btheta))\trace(\Sigma_{\psi,b}(\btheta)) w(\btheta)d\btheta\\
                        &+4\lambda_{\max}(\Sigma_{\psi_1,b})\lambda_{\max}^2(J)\int_{\btheta\in\Theta_\theta}\lambda_{\max}(\Sigma_{\psi}(\btheta))\trace\left(\Sigma_K(\btheta)\right) w(\btheta)d\btheta\\
                        &+4\lambda_{\max}(\Sigma_{\psi_1})\lambda_{\max}^2(J)\int_{\btheta\in\Theta_\theta}\lambda_{\max}(\Sigma_{\psi,b}(\btheta))\trace\left(\Sigma_K(\btheta)\right) w(\btheta)d\btheta\\
                        &+4\lambda_{\max}(\Sigma_{\psi_1,b})\lambda_{\max}(\Sigma_{\psi_1})\trace(H^2).
                    \end{split}
                \end{align}
                
                Furthermore, we have
                \begin{align*}
                    \rP&\left[\tpsi_{ij}(\balpha_0)b(\bV_i,\bxi)b(\bV_i,\bxi)b(\bV_j,\bxi)\right]\\
                    =&\rP\left[\left(\int_{\btheta\in\Theta_\theta}\psi(\bV_i, \balpha_0, \bzero,\btheta)\trans \psi(\bV_j, \balpha_0, \bzero,\btheta) w(\btheta)d\btheta\right)^2b(\bV_i,\bxi)b(\bV_j,\bxi)\right]\\
                    \leq&\int_{\btheta\in\Theta_\theta}\rP\left[\left(\psi(\bV_i, \balpha_0, \bzero,\btheta)\trans \psi(\bV_j, \balpha_0, \bzero,\btheta)\right)^2b(\bV_i,\bxi)b(\bV_j,\bxi)\right] w(\btheta)d\btheta\\
                    \leq&\int_{\btheta\in\Theta_\theta}\lambda_{\max}(\Sigma_{\psi,b}(\btheta))\rP\left[\left\|\psi(\bV_i, \balpha_0, \bzero,\btheta)\right\|^2b(\bV_i,\bxi)\right]w(\btheta)d\btheta\\
                    \leq&\int_{\btheta\in\Theta_\theta} \lambda_{\max}(\Sigma_{\psi,b}(\btheta))\trace(\Sigma_{\psi,b}(\btheta))w(\btheta)d\btheta,
                \end{align*}
                and
                \begin{align*}
                    \rP&\left[(\psi_1(\bV_i, \balpha_0)\trans K_{j})^2b(\bV_i,\bxi)b(\bV_j,\bxi)\right]\\
                    =&\rP\left[K_{j}\trans \Sigma_{\psi_1,b}K_{j}b(\bV_j,\bxi)\right] \\
                    \leq&\lambda_{\max}(\Sigma_{\psi_1,b})\lambda_{\max}^2(J)\int_{\btheta\in\Theta_\theta}\lambda_{\max}(\Sigma_{\psi,b}(\btheta))\trace\left(\Sigma_K(\btheta)\right) w(\btheta)d\btheta,
                \end{align*}
                and
                \begin{align*}
                    \rP\left[\left(\psi_1(\bV_i, \balpha_0)\trans H\psi_1(\bV_j, \balpha_0)\right)^2b(\bV_i,\bxi)b(\bV_j,\bxi)\right]
                    \leq&\lambda_{\max}^2(\Sigma_{\psi_1,b})\trace(H^2).
                \end{align*}
                
                Therefore, we have
                \begin{align*}
                    \rP[&h(\bV_i,\bV_j)^2b(\bV_i,\bxi)b(\bV_j,\bxi)]\\
                    \leq&4\rP\left[\tpsi_{ij}(\balpha_0)b(\bV_i,\bxi)b(\bV_i,\bxi)b(\bV_j,\bxi)\right]\\
                    &+4\rP\left[(\psi_1(\bV_i, \balpha_0)\trans H\psi_1(\bV_j, \balpha_0))^2 b(\bV_i,\bxi)b(\bV_j,\bxi)\right]\\
                    &+4\rP\left[(\psi_1(\bV_i, \balpha_0)\trans K_{j})^2 b(\bV_i,\bxi)b(\bV_j,\bxi)\right]+4\rP\left[(\psi_1(\bV_j, \balpha_0)\trans K_{i})^2 b(\bV_i,\bxi)b(\bV_j,\bxi)\right]\\
                    \leq&4\int_{\btheta\in\Theta_\theta} \lambda_{\max}(\Sigma_{\psi,b}(\btheta))\trace(\Sigma_{\psi,b}(\btheta))w(\btheta)d\btheta\\
                    &+8\lambda_{\max}(\Sigma_{\psi_1,b})\lambda_{\max}^2(J)\int_{\btheta\in\Theta_\theta}\lambda_{\max}(\Sigma_{\psi,b}(\btheta))\trace\left(\Sigma_K(\btheta)\right) w(\btheta)d\btheta\\
                    &+4\lambda_{\max}^2(\Sigma_{\psi_1,b})\trace(H^2),
                \end{align*}
                which combining with (\ref{eq:h2b}) and Lemmas \ref{lem:H} and \ref{lem:h2} implies that
                \begin{align*}
                    \rP&[h(\bV_1,\bV_2)^2(1+b(\bV_i,\bxi))(1+b(\bV_j,\bxi))]
                \end{align*}
                is bounded, and consequently completes the proof. \hfill $\square$
 
                \begin{lem}\label{lem:h2b12}
                    If Assumptions (A4) -- (A5) and (A7) hold, we have for any fixed $\bxi_1\in\Theta_{\alpha}$,
                    \begin{align*}
                        \rP[h(\bV_1,\bV_2)^2(1+b_1(\bV_1,\bxi_1))(1+b_1(\bV_2,\bxi_1))]
                    \end{align*}
                    is bounded,
                    where $b_1(\bv,\bxi_1)$ is defined in Assumption (A7), $h(\bV_1,\bV_2)$ is defined in (11) in the main paper, and $\bV_1$ and $\bV_2$ are generated from the null distribution with density $f(\bv,\balpha_0,\bzero,\btheta_0)$ and are independent of each other.	
                \end{lem}
                \noindent{\bf Proof of Lemma \ref{lem:h2b12}.}	We defined $\Sigma_{\psi,b_1}$ and $\Sigma_{\psi_1,b_1}$ as following,
                \begin{align*}
                    \Sigma_{\psi,b_1}=&\rP\left[b_1(\bV_i,\bxi_1)\psi(\bV_i, \balpha_0, \bzero,\btheta)^{\otimes2}\right],\\
                    \Sigma_{\psi_1,b_1}=&\rP\left[b_1(\bV_i,\bxi_1)\psi_1(\bV_i, \balpha_0)^{\otimes2}\right].
                \end{align*}
                It is straightforward to complete the proof by the same arguments in Lemma \ref{lem:h2b2}. \hfill $\square$
     
                \vfill
\noindent{\bf Proof of Theorem 1.}
For easy to express, we give some additional notations as follows.
\begin{align*}
	K =& \int_{\btheta\in\Theta_{\theta}} K(\btheta)w(\btheta) d\btheta,\\
	\tpsi_{i} =& \int_{\btheta\in\Theta_{\theta}}\psi(\bV_i, \balpha_0, \bzero, \btheta) w(\btheta)d\btheta,\\
    \tpsi_{ij}(\balpha) =& \int_{\btheta\in\Theta_{\theta}}\psi(\bV_i, \balpha, \bzero,\btheta)\trans \psi(\bV_j, \balpha, \bzero,\btheta) w(\btheta)d\btheta.
\end{align*}
By Lemma \ref{lems1}, we have
\begin{align*}
	\Pn\psi(\bV, \hbalpha, \bzero,\btheta)
	=&\Pn\psi(\bV, \balpha_0, \bzero,\btheta)	- K(\btheta)J\Psi_{1n}(\balpha_0)+o_p(n^{-1/2}),
\end{align*}
and
\begin{align*}
	\frac{1}{n}\sum_{i\neq j}\tpsi_{ij}(\hbalpha)
	=&\frac{1}{n}\sum_{i\neq j}\int_{\btheta\in\Theta_{\theta}}\psi(\bV_i, \hbalpha, \bzero,\btheta)\trans \psi(\bV_j, \hbalpha, \bzero,\btheta) w(\btheta)d\btheta\\
	=&\frac{1}{n}\sum_{i\neq j}\left[\tpsi_{ij}(\balpha_0) - \frac{1}{n}\Psi_{1n}(\balpha_0)\trans K_{i}- \frac{1}{n}K_{j}\trans \Psi_{1n}(\balpha_0) + \frac{1}{n^2}\Psi_{1n}(\balpha_0)\trans H\Psi_{1n}(\balpha_0)\right]\\
	&+\frac{1}{n}\sum_{i\neq j}\left[\tpsi_{i}	- \frac{1}{n}KJ\Psi_{1n}(\balpha_0)\right] o_p(n^{-1/2})\\
	&+\frac{1}{n}\sum_{i\neq j}\left[\tpsi_{j}-\frac{1}{n}KJ\Psi_{1n}(\balpha_0)\right] o_p(n^{-1/2})+o_p(n^{-1}).
\end{align*}
By notations, we have, for any $i\neq j$,
\begin{align*}
	\bE[K_{i}] = \bzero\quad\mbox{and}\quad \bE[\tpsi_{ij}(\balpha_0)] = 0.
\end{align*}
It can be shown by the central limit theorem (CLT) that
\begin{align*}
	\frac{1}{n}\sum_{i\neq j}\left[\tpsi_{i}	- \frac{1}{n}KJ\Psi_{1n}(\balpha_0)\right]
	=&\frac{n-1}{n}\sum_{i=1}^n\left[\tpsi_{i}	- \frac{1}{n}KJ\Psi_{1n}(\balpha_0)\right]\\
	=&O_p(n^{1/2}),
\end{align*}
which implies that
\begin{align*}
	T_n =&\frac{1}{n(n-1)}\sum_{i\neq j}\tpsi_{ij}(\hbalpha)\\
	=&\frac{1}{n(n-1)}\sum_{i\neq j}\left[\tpsi_{ij}(\balpha_0) - \frac{1}{n}\Psi_{1n}(\balpha_0)\trans K_{i}- \frac{1}{n}K_{j}\trans \Psi_{1n}(\balpha_0) + \frac{1}{n^2}\Psi_{1n}(\balpha_0)\trans H\Psi_{1n}(\balpha_0)\right]\\
    &+o_p(n^{-1})\\
	\equiv&T_{n1}+T_{n2}+T_{n3}+o_p(n^{-1}),
\end{align*}
where
\begin{align*}
	T_{n1} =& \frac{1}{n(n-1)}\sum_{i\neq j}\tpsi_{ij}(\balpha_0),\\
	T_{n2} =& -\frac{2}{n^2}\sum_{i=1}^n \Psi_{1n}(\balpha_0)\trans K_{i},\\
	T_{n3} =& \frac{1}{n^2}\Psi_{1n}(\balpha_0)\trans H\Psi_{1n}(\balpha_0).
\end{align*}
By some tedious calculations we have that
\begin{align*}
	\rP T_{n2}
	= -\frac{2}{n^2}\sum_{i=1}^{n}\rP[\psi_1(\bV_i, \balpha_0)\trans K_{i}]
	= -\frac{2}{n}\rP[\psi_1(\bV_1, \balpha_0)\trans K_{1}],
\end{align*}
and
\begin{align*}
	\rP T_{n3}
	=\frac{1}{n^2}\sum_{k=1}^n\rP[\psi_1(\bV_k, \balpha_0)\trans H\psi_1(\bV_k, \balpha_0)]
	=\frac{1}{n}\rP[\psi_1(\bV_1, \balpha_0)\trans H\psi_1(\bV_1, \balpha_0)],
\end{align*}
which implies that
\begin{align*}
	\rP\left[T_{n1}+T_{n2}+T_{n3}\right]
	=&-\frac{2}{n}\rP\left[\psi_1(\bV_1, \balpha_0)\trans K_{1}\right]
	+\frac{1}{n}\rP\left[\psi_1(\bV_1, \balpha_0)\trans H\psi_1(\bV_1, \balpha_0)\right]\\
    &+O(n^{-2})\\
	\equiv& \frac{\mu_0}{n} + O(n^{-2}).
\end{align*}

Thus, we can rewrite $T_n$ as
\begin{align*}
	T_n - \frac{\mu_0}{n} = \tT_n + R_n +o_p(n^{-1}),
\end{align*}
where
\begin{align*}
	\tT_n = \frac{1}{n(n-1)}\sum_{i\neq j}h(\bV_i,\bV_j)
\end{align*}
with the kernel given in (11) in the main paper
\begin{align*}
	h(\bV_i,\bV_j)
    =& \tpsi_{ij}(\balpha_0) - \psi_1(\bV_i, \balpha_0)\trans K_{j}\\
    &- K_{i}\trans \psi_1(\bV_j, \balpha_0) + \psi_1(\bV_i, \balpha_0)\trans H\psi_1(\bV_j, \balpha_0)
\end{align*}
and
\begin{align*}
	R_n =& -\frac{\mu_0}{n} - \frac{2}{n^2}\sum_{i=1}^n\psi_1(\bV_i, \balpha_0)\trans K_{i}
	+\frac{1}{n^2}\sum_{i=1}^n\psi_1(\bV_i, \balpha_0)\trans H\psi_1(\bV_i, \balpha_0)\\
	&+\frac{2}{n^2(n-1)}\sum_{i\neq j}^n\psi_1(\bV_i, \balpha_0)\trans K_{j}
	-\frac{1}{n^2(n-1)}\sum_{i\neq j}^n\psi_1(\bV_i, \balpha_0)\trans H\psi_1(\bV_j, \balpha_0).
\end{align*}
We are now going to verify that $R_n=o_p(n^{-1})$. In fact, we have, by the law of large numbers,
\begin{align*}
	-\frac{2}{n^2}\sum_{i=1}^n\psi_1(\bV_i, \balpha_0)\trans K_{i}
	+\frac{1}{n^2}\sum_{i=1}^n\psi_1(\bV_i, \balpha_0)\trans H\psi_1(\bV_i, \balpha_0)-\frac{\mu_0}{n} = O_p(n^{-3/2}),
\end{align*}
and
\begin{align*}
	\frac{2}{n^2(n-1)}\sum_{i\neq j}^n\psi_1(\bV_i, \balpha_0)\trans K_{j}
	-\frac{1}{n^2(n-1)}\sum_{i\neq j}^n\psi_1(\bV_i, \balpha_0)\trans H\psi_1(\bV_j, \balpha_0)
	=o_p(n^{-1}).
\end{align*}

Therefore, to complete the proof, it remains to show that
\begin{align*}
	n\tT_n \lkonv  \zeta,
\end{align*}
where $\zeta$ is a random variable of the form $\zeta=\sum_{j=1}^{\infty}\lambda_{j}(\chi^2_{1j}-1)$, and $\chi^2_{11},\chi^2_{12},...$ are independent $\chi^2_{1}$ variables, that is, $\zeta$ has characteristic function
\begin{align*}
	\bE\left[e^{it\zeta}\right]=\prod_{j=1}^{\infty}(1-2it\lambda_{j})^{-1/2}e^{-it\lambda_{j}}.
\end{align*}
In fact, this can be shown by Theorem in 5.5.2 of \cite{Serfling1980} by checking its conditions.
It can be shown that
\begin{align*}
	g_1(v_1) = \bE[h(\bV_1,\bV_2)|\bV_1=v_1]=0,
\end{align*}
and by Lemma \ref{lem:h2},
\begin{align*}
	\Var(h(\bV_1,\bV_2)) = \bE[h(\bV_1,\bV_2)^2]=\nu^2<\infty.
\end{align*} \hfill $\square$

\vfill
\noindent{\bf Proof of Theorem 2.}

For easy to express, we use the same notations $K$, $\tpsi_{i}$, $K_i$, $H$ and $\tpsi_{ij}(\balpha)$ in the proof of Theorem 1.
By Lemma \ref{lems1}, we have
\begin{align*}
	\Pn\psi(\bV, \hbalpha, \bzero,\btheta)
	=&\Pn\psi(\bV, \balpha_0, \bzero,\btheta)	- K(\btheta)J\Psi_{1n}(\balpha_0)+o_p(n^{-1/2}),
\end{align*}
and
\begin{align*}
	\frac{1}{n}\sum_{i\neq j}\tpsi_{ij}(\hbalpha)
	=&\frac{1}{n}\sum_{i\neq j}\int_{\btheta\in\Theta_{\theta}}\psi(\bV_i, \hbalpha, \bzero,\btheta)\trans \psi(\bV_j, \hbalpha, \bzero, \btheta) w(\btheta)d\btheta\\
	=&\frac{1}{n}\sum_{i\neq j}\left[\tpsi_{ij}(\balpha_0) - \frac{1}{n}\Psi_{1n}(\balpha_0)\trans K_{i}- \frac{1}{n}K_{j}\trans \Psi_{1n}(\balpha_0) + \frac{1}{n^2}\Psi_{1n}(\balpha_0)\trans H\Psi_{1n}(\balpha_0)\right]\\
	&+\frac{1}{n}\sum_{i\neq j}\left[\tpsi_{i}	- \frac{1}{n}KJ\Psi_{1n}(\balpha_0)\right] o_p(n^{-1/2})\\
	&+\frac{1}{n}\sum_{i\neq j}\left[\tpsi_{j}-\frac{1}{n}KJ\Psi_{1n}(\balpha_0)\right] o_p(n^{-1/2})+o_p(n^{-1}).
\end{align*}
Under the alternative hypothesis and by Assumption (A2), we have, for any $i\neq j$,
\begin{align*}
	\rP[\tpsi_{i}]
	=& \int_{\btheta\in\Theta_\theta}\rP[\psi(\bV_i, \balpha_0, \bzero,\btheta)] w(\btheta)d\btheta\\
	=& \int_{\btheta\in\Theta_\theta}\rP[\psi(\bV_i, \balpha_0, \bbeta, \btheta)] w(\btheta)d\btheta \\
     &- \bbeta\trans\int_{\btheta\in\Theta_\theta}\partial\rP[ \psi(\bV_i, \balpha_0, \bbeta, \btheta)/\partial\bbeta] w(\btheta)d\btheta(1+o(1))\\
	=&O(n^{-1/2}),
\end{align*}
and
\begin{align*}
	\rP[\psi_1(\bV_i,\balpha_0)]
	= \bbeta\trans\partial\rP[\psi_1(\bV_i,\balpha_0)]/\partial\bbeta(1+o(1))
	=O(n^{-1/2}).
\end{align*}
It can be shown by the law of large numbers that
\begin{align*}
	\frac{1}{n}\sum_{i\neq j}\left[\tpsi_{i}	- \frac{1}{n}KJ\Psi_{1n}(\balpha_0)\right]
	=\frac{n-1}{n}\sum_{i=1}^n\left[\tpsi_{i}	- \frac{1}{n}KJ\Psi_{1n}(\balpha_0)\right]
	=O_p(n^{1/2}),
\end{align*}
which implies that
\begin{align*}
	T_n =&\frac{1}{n(n-1)}\sum_{i\neq j}\tpsi_{ij}(\hbalpha)\\
	=&\frac{1}{n(n-1)}\sum_{i\neq j}\left[\tpsi_{ij}(\balpha_0) - \frac{1}{n}\Psi_{1n}(\balpha_0)\trans K_{i}- \frac{1}{n}K_{j}\trans \Psi_{1n}(\balpha_0) + \frac{1}{n^2}\Psi_{1n}(\balpha_0)\trans H\Psi_{1n}(\balpha_0)\right]\\
 &+o_p(n^{-1})\\
	\equiv&T_{n1}+T_{n2}+T_{n3}+o_p(n^{-1}),
\end{align*}
where
\begin{align*}
	T_{n1} =& \frac{1}{n(n-1)}\sum_{i\neq j}\tpsi_{ij}(\balpha_0),\\
	T_{n2} =& -\frac{2}{n^2}\sum_{i=1}^n \Psi_{1n}(\balpha_0)\trans K_{i},\\
	T_{n3} =& \frac{1}{n^2}\Psi_{1n}(\balpha_0)\trans H\Psi_{1n}(\balpha_0).
\end{align*}
By some tedious calculations we have that
\begin{align*}
	\rP T_{n2}
	= -\frac{2}{n^2}\sum_{i=1}^{n}\rP[\psi_1(\bV_i, \balpha_0)\trans K_{i}]
	= -\frac{2}{n}\rP[\psi_1(\bV_1, \balpha_0)\trans K_{1}],
\end{align*}
and
\begin{align*}
	\rP T_{n3}
	=\frac{1}{n^2}\sum_{k=1}^n\rP[\psi_1(\bV_k, \balpha_0)\trans H\psi_1(\bV_k, \balpha_0)]
	=\frac{1}{n}\rP[\psi_1(\bV_1, \balpha_0)\trans H\psi_1(\bV_1, \balpha_0)],
\end{align*}
which implies that
\begin{align*}
	\rP\left[T_{n2}+T_{n3}\right]
	=&-\frac{2}{n}\rP[\psi_1(\bV_1, \balpha_0)\trans K_{1}]
	+\frac{1}{n}\rP[\psi_1(\bV_1, \balpha_0)\trans H\psi_1(\bV_1, \balpha_0)]+O(n^{-2})\\
	\equiv& \frac{\mu_0}{n} + O(n^{-2}).
\end{align*}

Thus, we can rewrite $T_n$ as
\begin{align*}
	T_n - \frac{\mu_0}{n} = \tT_n + R_n +o_p(n^{-1}),
\end{align*}
where
\begin{align*}
	\tT_n = \frac{1}{n(n-1)}\sum_{i\neq j}h(\bV_i,\bV_j)
\end{align*}
with the kernel defined in (11) in the main paper
\begin{align*}
	h(\bV_i,\bV_j) = \tpsi_{ij}(\balpha_0) - \psi_1(\bV_i, \balpha_0)\trans K_{j}- K_{i}\trans \psi_1(\bV_j, \balpha_0) + \psi_1(\bV_i, \balpha_0)\trans H\psi_1(\bV_j, \balpha_0)
\end{align*}
and
\begin{align*}
	R_n =& -\frac{\mu_0}{n} - \frac{2}{n^2}\sum_{i=1}^n\psi_1(\bV_i, \balpha_0)\trans K_{i}
	+\frac{1}{n^2}\sum_{i=1}^n\psi_1(\bV_i, \balpha_0)\trans H\psi_1(\bV_i, \balpha_0)\\
	&+\frac{2}{n^2(n-1)}\sum_{i\neq j}^n\psi_1(\bV_i, \balpha_0)\trans K_{j}
	-\frac{1}{n^2(n-1)}\sum_{i\neq j}^n\psi_1(\bV_i, \balpha_0)\trans H\psi_1(\bV_j, \balpha_0).
\end{align*}
We are now going to verify that $R_n=o_p(n^{-1})$. In fact, we have, by the law of large numbers,
\begin{align*}
	-\frac{2}{n^2}\sum_{i=1}^n\psi_1(\bV_i, \balpha_0)\trans K_{i}
	+\frac{1}{n^2}\sum_{i=1}^n\psi_1(\bV_i, \balpha_0)\trans H\psi_1(\bV_i, \balpha_0)-\frac{\mu_0}{n} = O_p(n^{-3/2}),
\end{align*}
and
\begin{align*}
	\frac{2}{n^2(n-1)}\sum_{i\neq j}^n\psi_1(\bV_i, \balpha_0)\trans K_{j}
	-\frac{1}{n^2(n-1)}\sum_{i\neq j}^n\psi_1(\bV_i, \balpha_0)\trans H\psi_1(\bV_j, \balpha_0)
	=o_p(n^{-1}).
\end{align*}

Therefore, we complete the proof if we can show that
\begin{align*}
	n\tT_n \lkonv \zeta_1,
\end{align*}
and consequently by Theorem 9.5.2 of \cite{Resnick2019}, it suffices to show that for each $t\in\bR$,
\begin{align*}
	\left|\rP\left[e^{it (n-1)\tT_{n}}\right]-\rP\left[e^{it \zeta_1}\right]\right|\rightarrow 0,
\end{align*}
as $n$ goes to infinity.

Let $\{\bV_{0i},i=1,\cdots,n\}$ be independent and identically distributed sample for each $\bV_{0i}$ from the null distribution with density $f(\bv,\balpha_0,\bzero,\btheta_0)$. We show this along the lines of Theorem in Section 5.5.2 of \cite{Serfling1980}. As shown in Theorem in Section 5.5.2 of \cite{Serfling1980}, we rewrite, by the representation of $h(\bv_1,\bv_2)$ (Page 1087 in \cite{Linear_part2}),
\begin{align*}
	h(\bv_1,\bv_2) = \sum_{k=1}^{\infty}\lambda_{k}\phi_k(\bv_1)\phi_k(\bv_2),
\end{align*}
where $\{\phi_k(\cdot)\}$ and $\{\lambda_k\}$ have following properties that
\begin{align*}
	\rP[\phi_j(\bV_{01})\phi_k(\bV_{01})] = \left\{
	\begin{array}{ll}
		1, & \text{ if } j = k, \\
		0,& \text{ otherwise},
	\end{array}
	\right.
\end{align*}
\begin{align*}
	\rP[\phi_k(\bV_{01})]=0, \quad \mbox{for all } k,
\end{align*}
and
\begin{align*}
	\rP[h(\bV_{01},\bV_{02})^2]=\sum_{k=1}^{\infty}\lambda_k^2<\infty.
\end{align*}
Let
\begin{align*}
	T_{nK} = \frac{1}{n}\sum_{i\neq j}\sum_{k=1}^{K}\lambda_{k}\phi_k(\bV_i)\phi_k(\bV_j).
\end{align*}
According to the fact that $|e^{it}-1|\leq |t|$, we have
\begin{align*}
	\left|\rP\left[e^{it (n-1)\tT_n}\right]-\rP\left[e^{it T_{nK}}\right]\right|
	\leq& \rP\left[\left|e^{it (n-1)\tT_n}-e^{it T_{nK}}\right|\right]\\
	\leq& |t|\rP\left[\left|(n-1)\tT_n-T_{nK}\right|\right]\\
	\leq&|t|\left[\rP\left((n-1)\tT_n-T_{nK}\right)^2\right]^{1/2}.
\end{align*}
Next we claim that
\begin{align}\label{eq:tTn-Tn0}
	\rP\left((n-1)\tT_n-T_{nK}\right)^2\leq 4(2C_{\xi}^2+1)\sum_{k=K+1}^{\infty}\lambda_{k}^2+\left(\sum_{k=K+1}^{\infty}\lambda_{k}\mu_{ak}^2\right)^2.
\end{align}
In fact, $(n-1)\tT_n-T_{nK}$ has basically the form of a U-statistic, that is,
\begin{align*}
	(n-1)\tT_n-T_{nK}=\frac{2}{n}\binom{n}{2}U_{nK},
\end{align*}
where
\begin{align*}
	U_{nK}=\frac{1}{n(n-1)}\sum_{i\neq j}g_K(\bV_i,\bV_j)
\end{align*}
with
\begin{align*}
	g_K(\bV_i,\bV_j)
	=&h(\bV_i,\bV_j)-\sum_{k=1}^{K}\lambda_{k}\phi_k(\bV_i)\phi_k(\bV_j)\\
	=&\sum_{k=K+1}^{\infty}\lambda_{k}\phi_k(\bV_i)\phi_k(\bV_j).
\end{align*}

Before showing the inequality (\ref{eq:tTn-Tn0}), we derive following useful inequalities. 
It is seen that
\begin{align*}
	\rP[\psi(\bV_i,\balpha_0,  0, \btheta)]
	=&\rP\left[\psi(\bV_{0i},\balpha_0,  0, \btheta)\frac{f(\bV_{0i},\balpha_0,\bbeta,\btheta_0)}{f(\bV_{0i},\balpha_0,\bzero,\btheta_0)}\right],
\end{align*}
and consequently, by  Assumption (A6) and Cauchy-Schwarz inequality, there is a constant $C_{\bxi}$ such that
\begin{align*}
	\|\rP\psi(\bV_i,\balpha_0,  0, \btheta)\|
    \leq\frac{1}{\sqrt{n}}C_{\bxi}\left(\rP\|\psi(\bV_{0i},\balpha_0,  0, \btheta)\|^2\right)^{1/2}.
\end{align*}
By Assumption (A5), it is followed from the above that there is a constant $C_1>0$ such that
\begin{align*}
	\rP[\tpsi_{ij}(\balpha_0)]
	=&\int_{\btheta\in\Theta_\theta}\left\|\rP[\psi(\bV_i,\balpha_0,  0, \btheta)]\right\|^2w(\btheta) d\btheta\\
	\leq&\frac{1}{n}C_{\bxi}^2\int_{\btheta\in\Theta_\theta}\rP\|\psi(\bV_{0i},\balpha_0,  0, \btheta)\|^2w(\btheta) d\btheta\\
    =&\frac{1}{n}C_1C_{\bxi}^2,
\end{align*}
and
\begin{align*}
	\rP[K_i]
	=&\int_{\btheta\in\Theta_\theta}J\trans K(\btheta)\trans\rP[\psi(\bV_i,\balpha_0,  0, \btheta)]w(\btheta) d\btheta,
\end{align*}
and by Assumption (A3), there is a constant $C_3>0$ such that

\begin{align*}
	\|\rP[K_i]\|^2
	\leq&\int_{\btheta\in\Theta_\theta}\|J\trans K(\btheta)\trans\rP[\psi(\bV_i,\balpha_0,  0, \btheta)]\|^2w(\btheta) d\btheta\\
	\leq&\lambda_{\max}^2(J)\int_{\btheta\in\Theta_\theta}\lambda_{\max}\left(K(\btheta) K(\btheta)\trans\right)\|\rP\left[\psi(\bV_i,\balpha_0,  0, \btheta)\right]\|^2w(\btheta) d\btheta\\
	\leq&\frac{1}{n}C_{\bxi}^2\lambda_{\max}^2(J)\int_{\btheta\in\Theta_\gamma}\lambda_{\max}\left(K(\btheta) K(\btheta)\trans\right)\rP\|\psi(\bV_{0i},\balpha_0,  0, \btheta)\|^2w(\btheta) d\btheta\\
	\leq&\frac{1}{n}C_3C_{\bxi}^2.
\end{align*}
As the same arguments, we have
\begin{align*}
	\rP[\psi_1(\bV_i,\balpha_0)]
	=\rP\left[\psi_1(\bV_{0i},\balpha_0)\frac{f(\bV_{0i},\balpha_0,\bbeta,\btheta_0)}{f(\bV_{0i},\balpha_0,\bzero,\btheta_0)}\right],
\end{align*}
and there are constant $C_2>0$,
\begin{align*}
	\|\rP[\psi(\bV_i,\balpha_0)]\|
	\leq&\frac{1}{\sqrt{n}}C_{\bxi}\left\{\rP\left[\|\psi_1(\bV_{0i},\balpha_0)\|^2\right]\right\}^{1/2}\\
	\leq&\frac{1}{\sqrt{n}}C_2C_{\bxi}
\end{align*}

From the foregoing considerations, we have
\begin{align}\label{eq:hxy}
	\begin{split}
		|\rP[h(\bV_i,\bV_j)]|
		=&\int_{\btheta\in\Theta_\gamma}\left\|\rP[\psi(\bV_i,  0, \btheta)]\right\|^2w(\btheta) d\btheta\\
		&+\rP\left[\psi_1(\bV_i, \balpha_0)\trans K_{j}+K_{i}\trans \psi_1(\bV_j, \balpha_0)\right]\\
        &+ \rP\left[\psi_1(\bV_i, \balpha_0)\trans H\psi_1(\bV_j, \balpha_0)\right]\\
		\leq&\frac{1}{n}C_1C_{\bxi}^2+2\frac{1}{n}C_3C_{\bxi}^2+\frac{1}{n}C_2C_{\bxi}^2\\
		=&O(n^{-1}),
	\end{split}
\end{align}
and by the mean value theorem and combing the dominated convergence theorem and Assumption (A6),
\begin{align*}
	\rP[\phi_k(\bV_i)]
	=&\rP\left[\phi_k(\bV_{0i})\frac{f(\bV_{0i},\balpha_0,\bbeta,\btheta_0)}{f(\bV_{0i},\balpha_0,\bzero,\btheta_0)}\right]\\
	=&\frac{1}{n^{1/2}}\rP\left[\phi_k(\bV_{0i})\frac{1}{f(\bV_{0i},\balpha_0,\bzero,\btheta_0)}\frac{\partial f(\bV_{0i},\balpha_0,\tbbeta,\btheta_0)}{\partial\bbeta\trans}\bxi\right]\\
	=&\frac{1}{n^{1/2}}\rP\left[\phi_k(\bV_{0i})\frac{1}{f(\bV_{0i},\balpha_0,\bzero,\btheta_0)}\frac{\partial f(\bV_{0i},\balpha_0,\bzero,\btheta_0)}{\partial\bbeta\trans}\bxi\right](1+o(1))\\
	=&\frac{1}{n^{1/2}}\mu_{ak}+o(n^{-1/2}),
\end{align*}
where $\tbbeta$ is between $\bzero$ and $\bbeta$,
\begin{align*}
	\left|\rP[\phi_k(\bV_i)]\right|
	\leq&\frac{1}{\sqrt{n}}C_{\bxi}\left(\rP\left[\phi_k(\bV_{0i})^2\right]\right)^{1/2}\\
	=&\frac{1}{\sqrt{n}}C_{\bxi},
\end{align*}
and
\begin{align*}
	\rP[\phi_k(\bV_i)^2]
	=&\rP\left[\phi_k(\bV_{0i})^2\frac{f(\bV_{0i},\balpha_0,\bbeta,\btheta_0)}{f(\bV_{0i},\balpha_0,\bzero,\btheta_0)}\right]\\
	\leq&1+\frac{1}{\sqrt{n}}\rP\left[\phi_k(\bV_{0i})^2b(\bV_{0i},\bxi)\right]\\
	=&1+O(n^{-1/2}).
\end{align*}
Since by Lemma \ref{lem:h2b2} combining Assumption (A6),
\begin{align}\label{eq:gk2}
\begin{split}
	\rP[g_K(\bV_i,\bV_j)^2]
	=&\int\int g_K(\bv_1,\bv_2)^2f(\bv_1,\balpha_0,\bbeta,\btheta_0)f(\bv_2,\balpha_0,\bbeta,\btheta_0)d\bv_1d\bv_2\\
	\leq&\rP[g_K(\bV_{0i},\bV_{0j})^2(1+n^{-1/2}b(\bV_{0i},\bxi))(1+n^{-1/2}b(\bV_{0j},\bxi))]\\
	\leq&\rP[g_K(\bV_{0i},\bV_{0j})^2(1+n^{-1/2}b(\bV_{0i},\bxi))^2]\\
	=&\sum_{k=K+1}^{\infty}\lambda_{k}^2\rP[\phi_k(\bV_{0i})^2(1+n^{-1/2}b(\bV_{0i},\bxi))^2]\rP[\phi_k(\bV_{0j})^2]\\
	=&\sum_{k=K+1}^{\infty}\lambda_{k}^2\rP[\phi_k(\bV_{0i})^2(1+n^{-1/2}b(\bV_{0i},\bxi))^2]\\
	\leq&(1+3Cn^{-1/2})\sum_{k=K+1}^{\infty}\lambda_{k}^2\\
	<&\infty,
\end{split}
\end{align}
where $C$ is a constant such that $\rP[\phi_k(\bV_i)^2b(\bV_{0i},\bxi)^2]\leq C$,
thus, it is straightforward to obtain that by the dominated theorem and for large $n$,
\begin{align*}
	\rP[g_K(\bV_i,\bV_j)]
	=&\sum_{k=K+1}^{\infty}\lambda_{k}(\rP[\phi_k(\bV_i)])^2\\
	=&\frac{1}{n}\left(\sum_{k=K+1}^{\infty}\lambda_{k}\mu_{ak}^2\right)(1+o(1)),
\end{align*}
and
\begin{align*}
	\rP[g_K(\bV_i,\bV_j)^2]
	=&\int\int g_K(\bv_1,\bv_2)^2f(\bv_1,\balpha_0,\bbeta,\btheta_0)f(\bv_2,\balpha_0,\bbeta,\btheta_0)d\bv_1d\bv_2\\
	\leq&\rP[g_K(\bV_{0i},\bV_{0j})^2](1+o(1))\\
	=&(1+o(1))\sum_{k=K+1}^{\infty}\lambda_{k}^2\\
	\leq&2\sum_{k=K+1}^{\infty}\lambda_{k}^2,
\end{align*}
which implies that
\begin{align*}
	\Var(g_K(\bV_i,\bV_j))
	\leq\rP[g_K(\bV_i,\bV_j)^2]
	\leq2\sum_{k=K+1}^{\infty}\lambda_{k}^2.
\end{align*}
Let $b(\bv)=\rP[g_K(\bV_i,\bV_j)|\bV_j=\bv]$. Noting that
\begin{align*}
	b(\bv)=&\rP[g_K(\bV_i,\bV_j)|\bV_j=\bv]\\
	=&\sum_{k=K+1}^{\infty}\lambda_{k}\phi_k(\bv)\rP[\phi_k(\bV_i)],
\end{align*}
we have, along same lines as (\ref{eq:gk2}), for $n$ large enough,
\begin{align*}
	\rP\left[b(\bV_i)^2\right]
	=&\rP\left[b(\bV_{0i})^2\frac{f(\bV_{0i},\balpha_0,\bbeta,\theta_0)}{f(\bV_{0i},\balpha_0,\bzero,\theta_0)}\right]\\
	\leq&\rP\left[b(\bV_{0i})^2(1+n^{-1/2}b(\bV_{0i},\bxi))\right]\\
	=&\rP\left[\left(\E[g_K(\bV_i,\bV_j)|\bV_i=\bV_{0i}]\right)^2(1+n^{-1/2}b(\bV_{0i},\bxi))\right]\\
	\leq&\rP\left[\rP[g_K(\bV_i,\bV_j)^2(1+n^{-1/2}b(\bV_{0i},\bxi))|\bV_i=\bV_{0i}]\right]\\
	=&\rP\left[g_K(\bV_{0i},\bV_j)^2(1+n^{-1/2}b(\bV_{0i},\bxi))\right]\\ =&\rP\left[g_K(\bV_{0i},\bV_{0j})^2(1+n^{-1/2}b(\bV_{0i},\bxi))
      \frac{f(\bV_{0j},\balpha_0,\bbeta,\theta_0)}{f(\bV_{0j},\balpha_0,\bzero,\theta_0)}\right]\\
	\leq&\rP\left[g_K(\bV_{0i},\bV_{0j})^2(1+n^{-1/2}b(\bV_{0i},\bxi))(1+n^{-1/2}b(\bV_{0j},\bxi))\right]\\
	<&\infty,
\end{align*}
and consequently by the dominated theorem,
\begin{align*}
	\Var(\rP[g_K(\bV_i,\bV_j)|\bV_i])
	\leq&\rP\left[b(\bV_i)^2\right]\\
	=&\rP\left[b(\bV_{0i})^2\right](1+o(1))\\
	=&(1+o(1))\sum_{k=K+1}^{\infty}\lambda_{k}^2(\rP[\phi_k(\bV_j)])^2\\
	\leq&\frac{2}{n}C_{\xi}^2\sum_{k=K+1}^{\infty}\lambda_{k}^2.
\end{align*}
Lemma A in Section 5.2.1 of \cite{Serfling1980} states that
\begin{align*}
	\Var(U_{nK})
	\leq&\binom{n}{2}^{-1}\left(4C_{\xi}^2\sum_{k=K+1}^{\infty}\lambda_{k}^2+2\sum_{k=K+1}^{\infty}\lambda_{k}^2\right)\\
	\leq&2(2C_{\xi}^2+1)\binom{n}{2}^{-1}\sum_{k=K+1}^{\infty}\lambda_{k}^2,
\end{align*}
which implies that
\begin{align*}
	\rP\left((n-1)\tT_n-T_{nK}\right)^2
	=&(n-1)^2\left(\Var(U_{nK})+(\rP[U_{nK}])^2\right)\\
	\leq&4(2C_{\xi}^2+1)\sum_{k=K+1}^{\infty}\lambda_{k}^2+(n-1)^2(\E[U_{nK}])^2,
\end{align*}
and consequently (\ref{eq:tTn-Tn0}) holds.

It is immediately obtained that $n\rP[h(\bV_i,\bV_j)]<+\infty$ by (\ref{eq:hxy}), which implies that the series
\begin{align*}
	\sum_{k=1}^{\infty}\lambda_{k}\mu_{ak}^2=n\rP[h(\bV_i,\bV_j)]<\infty,
\end{align*}
and consequently that
\begin{align*}
	\lim_{K\rightarrow\infty}\sum_{k=K+1}^{\infty}\lambda_{k}\mu_{ak}^2=0.
\end{align*}
For fixed $t$ and given $\eps>0$, choose and fix $K$ large enough such that
\begin{align*}
	2|t|\left(\sum_{k=K+1}^{\infty}\lambda_{k}^2\right)^{1/2}&<\eps, \quad \mbox{and},\\
	|t|\sum_{k=K+1}^{\infty}\lambda_{k}\mu_{ak}^2&<\eps.
\end{align*}
Then it can been shown that for all $n$ sufficiently large,
\begin{align*}
	|t|^2(\rP[U_{nK}])^2=|t|^2\left(\sum_{k=K+1}^{\infty}\lambda_{k}\mu_{ak}^2\right)^2<\eps^2,
\end{align*}
and immediately,
\begin{align}\label{eq:tTn-Tn10}
	\left|\rP\left[e^{it (n-1)\tT_n}\right]-\rP\left[e^{it T_{nK}}\right]\right|<(2C_{\bxi}^2+2)^{1/2}\eps.
\end{align}

We now are going to show that
\begin{align}\label{eq:tTnk0}
	T_{nK}\lkonv Y_K=\sum_{k=1}^{K}\lambda_k(W_k^2-1),
\end{align}
where $\{W_k,k=1,\cdots,K\}$ are independent and each $W_k$ is normally distributed with mean $\mu_{ak}$ and unit variance.
We can rewrite
\begin{align*}
	T_{nK}= \sum_{k=1}^{K}\lambda_k(W_{nk}^2-Z_{nk}),
\end{align*}
where
\begin{align*}
	W_{nk}= \frac{1}{\sqrt{n}}\sum_{i=1}^{n}\phi_k(\bV_i),
\end{align*}
and
\begin{align*}
	Z_{nk}= \frac{1}{n}\sum_{i=1}^{n}\phi_k^2(\bV_i).
\end{align*}
Recalling the previous arguments, it can be obtained that
\begin{align*}
	\rP[	W_{nk}]= \frac{1}{\sqrt{n}}\sum_{i=1}^{n}\frac{1}{\sqrt{n}}\mu_{ak}(1+o(1))=\mu_{ak}(1+o(1)),
\end{align*}
\begin{align*}
	\Var(W_{nk})= \frac{1}{n}\sum_{i=1}^{n}\Var(\phi_k(\bV_i))=1+O(n^{-1/2})-\frac{1}{n}\mu_{ak}^2\rightarrow1,
\end{align*}
and if $j\neq k$,
\begin{align*}
	\Cov(W_{nj}, W_{nk})
	=& \rP[W_{nj}W_{nk}] - \rP[	W_{nj}]\rP[	W_{nk}]\\
	=&\frac{1}{n}\sum_{i=1}^{n}\rP[\phi_j(\bV_i)\phi_k(\bV_i)]+\frac{1}{n}\sum_{i\neq l}^{n}\rP[\phi_j(\bV_i)\phi_k(\bV_l)]-\mu_{aj}\mu_{ak}(1+o(1))\\
	=&\frac{1}{n}\sum_{i=1}^{n}\rP[\phi_j(\bV_i)\phi_k(\bV_i)]+\frac{1}{n}\sum_{i\neq l}^{n}\frac{1}{n}\mu_{aj}\mu_{ak}(1+o(1))-\mu_{aj}\mu_{ak}(1+o(1))\\
	=&\rP[\phi_j(\bV_i)\phi_k(\bV_i)]-\frac{1}{n}\mu_{aj}\mu_{ak}(1+o(1))\\
	=&\rP\left[\phi_j(\bV_{0i})\phi_k(\bV_{0i})\right](1+n^{-1/2}C_{\bxi})-\frac{1}{n}\mu_{aj}\mu_{ak}(1+o(1))\\
	\rightarrow&0,
\end{align*}
which implies that if $j\neq k$,
\begin{align*}
	\left|\Cov(W_{nj}, W_{nk})\right|\rightarrow0,
\end{align*}
and consequently,
\begin{align*}
	\Cov(W_{nj}, W_{nk})\rightarrow\left\{
	\begin{array}{ll}
		1, & \text{ if } j = k, \\
		0,& \text{ otherwise}.
	\end{array}
	\right.
\end{align*}
Therefore, by Lindeberg-Levy central limit theorem, we have
\begin{align*}
	(W_{n1},\cdots, W_{nK})\trans \lkonv N(\bmu_a, \bone_{K\times K}),
\end{align*}
where $\bmu_a=(\mu_{a1},\cdots,\mu_{aK})\trans$. Since $\rP[\phi_k^2(\bV_i)]\rightarrow1$, it is seen that by the strong law of large numbers,
\begin{align*}
	Z_{nk}\askonv 1,
\end{align*}
which implies that
\begin{align*}
	(Z_{n1},\cdots,Z_{nK})\trans\askonv (1,\cdots,1)\trans.
\end{align*}
Thus, (\ref{eq:tTnk0}) holds and consequently,
\begin{align}\label{eq:tTn-Tn20}
	\left|\rP\left[e^{it T_{nK}}\right]-\rP\left[e^{it Y_{K}}\right]\right|<\eps,
\end{align}
for all n sufficiently large.

Noting that by (\ref{eq:tTn-Tn10}) and (\ref{eq:tTn-Tn20}), for all $n$ large enough,
\begin{align*}
	\left|\rP\left[e^{it (n-1)\tT_{n}}\right]-\rP\left[e^{it \zeta_1}\right]\right|
	\leq &\left|\rP\left[e^{it (n-1)\tT_{n}}\right]-\rP\left[e^{it T_{nK}}\right]\right|
	+\left|\rP\left[e^{it T_{nK}}\right]-\rP\left[e^{it Y_{K}}\right]\right|\\
	&+\left|\rP\left[e^{it Y_{K}}\right]-\rP\left[e^{it \zeta_1}\right]\right|\\
	<&(2C_{\bxi}^2+2)^{1/2}\eps+ \left|\rP\left[e^{it Y_{K}}\right]-\rP\left[e^{it \zeta_1}\right]\right|.
\end{align*}
Therefore, to complete the proof, it suffices to establish that for all $K$ sufficiently large,
\begin{align}\label{eq:Yk0}
	\left|\rP\left[e^{it Y_{K}}\right]-\rP\left[e^{it \zeta_1}\right]\right|< \eps\left\{\rP[(W_1^2-1)^2]\right\}^{1/2}.
\end{align}
Let the random variable $W_1$, $W_2$, $\cdots$ be defined on a common probability space and denote by $\zeta_1$ the limit in mean square of $Y_K$ as $K\rightarrow\infty$. Then, by the definition of $Y_K$, we have
\begin{align*}
	\left|\rP\left[e^{it Y_{K}}\right]-\rP\left[e^{it \zeta_1}\right]\right|
	\leq& |t|\left[\rP\left(Y_{K}-\zeta_1\right)^2\right]^{1/2}\\
	\leq&|t|\left[\rP\left(W_1^2-1\right)^2\right]^{1/2}\left(\sum_{k=K+1}^{\infty}\lambda_k^2\right)^{1/2}\\
	<&\eps\left[\rP\left(W_1^2-1\right)^2\right]^{1/2},
\end{align*}
which arrives at (\ref{eq:Yk0}), and then completes the proof. \hfill $\square$

\vfill
\noindent{\bf Proof of Theorem 3.}
It is straightforward to complete the proof if we have $nT_n^*-\mu_0 \lkonv\zeta$, where $\zeta$ is defined in Theorem 1. Thus, it suffices to show that
\begin{align}\label{Tns}
	nT_n^*-\mu_0 \lkonv \zeta.
\end{align}
Since the above asymptotic distribution is derived given the observed sample $\{\bV_i, i=1,\cdots,n\}$, we can treat $\hbalpha$ as a deterministic sequence satisfying  $\hbalpha-\balpha_0=\frac{1}{\sqrt{n}}\bxi_1$ for fixed $\bxi_1\in\Theta_{\alpha}$, where we use the fact that $\hbalpha$ is $\sqrt{n}$-consistent by Lemma \ref{lem0}.

Denote by $\bE^*$ the expectation generated by the resampling from $\bV$.
We give additional notations as follows.
\begin{align*}
	K^*(\btheta)=&\frac{\partial \bE^*\psi(\bV^*, \hbalpha, \bzero, \btheta)}{\partial\balpha},\\
	K^* =& \int_{\btheta\in\Theta_\gamma} K^*(\btheta)w(\btheta) d\btheta,\\
	H^*=&\int_{\btheta\in\Theta_\gamma}J\trans K^*(\btheta)\trans K^*(\btheta)J w(\btheta)d\btheta,\\
	\tpsi_{i}^* =& \int_{\btheta\in\Theta_\gamma}\psi(\bV_i^*, \hbalpha^*, \bzero, \btheta) w(\btheta)d\btheta,\\
	\tK_i^* =& \int_{\btheta\in\Theta_\gamma}J\trans K^*(\btheta)\trans S(\bV_i^*, \hbalpha^*, \bzero,\btheta) w(\btheta)d\btheta,\\
	K_i^* =& \int_{\btheta\in\Theta_\gamma}J\trans K(\btheta)\trans\psi(\bV_i^*, \hbalpha^*, \bzero,\btheta) w(\btheta)d\btheta,\\
	\tpsi_{ij}^*(\balpha) =& \int_{\btheta\in\Theta_\gamma}\psi(\bV_i^*, \balpha, \bzero,\btheta)\trans \psi(\bV_j^*, \balpha, \bzero,\btheta) w(\btheta)d\btheta.
\end{align*}

By Lemma \ref{lems1}, we have
\begin{align*}
	\Pn\psi(\bV_i^*, \hbalpha^*, \bzero,\btheta)
	=&\Pn\psi(\bV_i^*, \hbalpha, \bzero,\btheta)	- K^*(\btheta)J\Pn\psi_1\bV_i^*,\hbalpha)+o_{p^*}(n^{-1/2}),
\end{align*}
and
\begin{align*}
	\frac{1}{n}\sum_{i\neq j}\tpsi_{ij}^*(\hbalpha^*)
	=&\frac{1}{n}\sum_{i\neq j}\int_{\btheta\in\Theta_\gamma}\psi(\bV_i^*, \hbalpha^*, \bzero,\btheta)\trans \psi(\bV_j^*, \hbalpha^*, \bzero,\btheta) w(\btheta)d\btheta\\
	=&\frac{1}{n}\sum_{i\neq j}\left[\tpsi_{ij}^*(\hbalpha) - \frac{1}{n}\Psi_{0n}(\hbalpha)\trans \tK_{i}^*- \frac{1}{n}\tK_{j}^{*\trans}\Psi_{1n}(\hbalpha)  + \frac{1}{n^2}\Psi_{1n}(\hbalpha)\trans H^*\Psi_{1n}(\hbalpha)\right]\\
	&+\frac{1}{n}\sum_{i\neq j}\left[\tpsi_{i}^* - \frac{1}{n}\tK^*J\Psi_{1n}(\hbalpha)\right] o_{p^*}(n^{-1/2})\\
	&+\frac{1}{n}\sum_{i\neq j}\left[\tpsi_{j}^*-\frac{1}{n}\tK^*J\Psi_{1n}(\hbalpha)\right] o_{p^*}(n^{-1/2})+o_{p^*}(n^{-1}).
\end{align*}
By above notations, we have, for any $i\neq j$,
\begin{align*}
	\bE^*[\tK_{i}^*] = \bzero\quad\mbox{and}\quad \bE^*[\tpsi_{ij}^*(\hbalpha)] = 0.
\end{align*}
It can be shown by the law of large numbers applying to bootstrap sample $\{\bV_i^*,i=1, \cdots,n\}$ that
\begin{align*}
	\frac{1}{n}\sum_{i\neq j}\left[\tpsi_{i}^*	- \frac{1}{n}\tK^*J\Psi_{1n}(\hbalpha)\right]
	=\frac{n-1}{n}\sum_{i=1}^n\left[\tpsi_{i}^*	- \frac{1}{n}\tK^*J\Psi_{1n}(\hbalpha)\right]
	=O_{p^*}(n^{1/2}),
\end{align*}
which implies that
\begin{align*}
	T_n^* =&\frac{1}{n(n-1)}\sum_{i\neq j}\tpsi_{ij}^*(\hbalpha^*)\\
	=&\frac{1}{n(n-1)}\sum_{i\neq j}\left[\tpsi_{ij}^*(\hbalpha) - \frac{1}{n}\Psi_{1n}(\hbalpha)\trans \tK_{i}^*- \frac{1}{n}\tK_{j}^{*\trans} \Psi_{1n}(\hbalpha) + \frac{1}{n^2}\Psi_{1n}(\hbalpha)\trans H^*\Psi_{1n}(\hbalpha)\right]\\
	 &+o_p(n^{-1})\\
	\equiv&T_{n1}^*+T_{n2}^*+T_{n3}^*+o_{p^*}(n^{-1}),
\end{align*}
where
\begin{align*}
	T_{n1}^* =& \frac{1}{n(n-1)}\sum_{i\neq j}\tpsi_{ij}^*(\hbalpha),\\
	T_{n2}^* =& -\frac{2}{n^2}\sum_{i=1}^n \Psi_{1n}(\hbalpha)\trans \tK_{i}^*,\\
	T_{n3}^* =& \frac{1}{n^2}\Psi_{1n}(\hbalpha)\trans H^*\Psi_{1n}(\hbalpha).
\end{align*}
Furthermore, we have that
\begin{align*}
	\bE^* [T_{n2}^*]
	= -\frac{2}{n^2}\sum_{i=1}^{n}\bE^*[\psi_1(\bV^*, \hbalpha)\trans \tK_{i}^*]
	= -\frac{2}{n}\bE^*[\psi_1(\bV, \hbalpha)\trans K_{1}^*],
\end{align*}
and
\begin{align*}
	\bE^*[ T_{n3}^*]
	=\frac{1}{n^2}\sum_{k=1}^n\bE^*[\psi_1(\bV^*, \hbalpha)\trans H^*\psi_1(\bV_k, \hbalpha)]
	=\frac{1}{n}\bE^*[\psi_1(\bV^*, \hbalpha)\trans H^*\psi_1(\bV^*, \hbalpha)],
\end{align*}
which implies that
\begin{align*}
	\bE^*\left[T_{n1}^*+T_{n2}^*+T_{n3}^*\right]
	=&-\frac{2}{n}\bE^*[\psi_1(\bV^*, \hbalpha)\trans \tK_{1}^*]
	+\frac{1}{n}\bE^*[\psi_1(\bV^*, \hbalpha)\trans H^*\psi_1(\bV^*, \hbalpha)]+O(n^{-2})\\
	\equiv& \frac{\mu_0(\hbalpha)}{n} + O(n^{-2}).
\end{align*}
By the continuity of $\mu_0(\hbalpha)$, we have
\begin{align*}
\bE^*\left[T_{n1}^*+T_{n2}^*+T_{n3}^*\right] = \frac{\mu_0(\hbalpha)}{n}(1+o_p(1)).
\end{align*}

By the Assumption (A7) and the dominated convergence theorem, it is seen that for all $\btheta\in\Theta_{\gamma}$,
\begin{align*}
	K^*(\btheta)
	=&\partial\int \psi(\bv, \hbalpha, \bzero, \btheta)f(\bv,\hbalpha,\bzero,\btheta)d\bv/\partial\balpha\\
	=&(1+O(n^{-1/2}))\partial \int \psi(\bv, \balpha_0, \bzero, \btheta)f(\bv,\balpha_0,\bzero,\btheta)d\bv/\partial\balpha\\
	=&K(\btheta)(1+O(n^{-1/2})),
\end{align*}
and consequently,
\begin{align*}
	K^*=K(1+o(1)),\\
	H^*=H(1+o(1)).
\end{align*}
By the same arguments, we have
\begin{align*}
	\mu_0(\hbalpha) = \mu_0(1+o(1)).
\end{align*}
Thus, we can rewrite $T_n^*$ as
\begin{align*}
	T_n^* - \frac{\mu_0}{n}(1+o(1)) = \tT_n^*(1+o_{p^*}(1)) + R_n^* +o_{p^*}(n^{-1}),
\end{align*}
where
\begin{align*}
	\tT_n^* = \frac{1}{n(n-1)}\sum_{i\neq j}h(\bV_i^*,\bV_j^*,\hbalpha),
\end{align*}
\begin{align*}
	R_n^* =& -\frac{\mu_0(\hbalpha)}{n} - \frac{2}{n^2}\sum_{i=1}^n\psi_1(\bV_i^*, \hbalpha)\trans \tK_{i}^*
	+\frac{1}{n^2}\sum_{i=1}^n\psi_1(\bV_i^*, \hbalpha)\trans H^*\psi_1(\bV_i^*, \hbalpha)\\
	&+\frac{2}{n^2(n-1)}\sum_{i\neq j}^n\psi_1(\bV_i^*, \hbalpha)\trans \tK_{j}^*
	-\frac{1}{n^2(n-1)}\sum_{i\neq j}^n\psi_1(\bV_i^*, \hbalpha)\trans H^*\psi_1(\bV_j^*, \hbalpha),
\end{align*}
and
\begin{align*}
	h(\bV_i^*,\bV_j^*,\hbalpha)
	=& \tpsi_{ij}^*(\hbalpha) - \psi_1(\bV_i^*, \hbalpha)\trans K_{j}^*\\
	 &- (K_{i}^*)\trans \psi_1(\bV_j^*, \hbalpha) + \psi_1(\bV_i^*, \hbalpha)\trans H \psi_1(\bV_j^*, \hbalpha).
\end{align*}

It can be verified that $R_n^*=o_{p^*}(n^{-1})$. In fact, we have, by the law of large numbers applying to bootstrap sample $\{\bV_i^*,i=1, \cdots,n\}$,
\begin{align*}
	-\frac{2}{n^2}\sum_{i=1}^n\psi_1(\bV_i^*, \hbalpha)\trans \tK_{i}^*
	+\frac{1}{n^2}\sum_{i=1}^n\psi_1(\bV_i^*, \hbalpha)\trans H^*\psi_1(\bV_i^*, \hbalpha)-\frac{\mu_0(\hbalpha)}{n} = O_{p^*}(n^{-3/2}),
\end{align*}
and
\begin{align*}
	\frac{2}{n^2(n-1)}\sum_{i\neq j}^n\psi_1(\bV_i^*, \hbalpha)\trans \tK_{j}^*
	-\frac{1}{n^2(n-1)}\sum_{i\neq j}^n\psi_1(\bV_i^*, \hbalpha)\trans H^*\psi_1(\bV_j^*, \hbalpha)
	=o_{p^*}(n^{-1}).
\end{align*}

Next, we are going to show that
\begin{align*}
	n\tT_n^* \lkonv \zeta,
\end{align*}
where $\zeta$ is a random variable of the form $\zeta=\sum_{j=1}^{\infty}\lambda_{j}(\chi^2_{1j}-1)$, and $\chi^2_{11},\chi^2_{12},\cdots$ are independent $\chi^2_{1}$ variables, that is, $\zeta$ has characteristic function
\begin{align*}
	\rP\left[e^{it\zeta}\right]=\prod_{j=1}^{\infty}(1-2it\lambda_{j})^{-1/2}e^{-it\lambda_{j}}.
\end{align*}
Here $\{\lambda_{j}\}$ are the eigenvalues of kernel $h(\bv_1,\bv_2)$ under $f(\bv,\balpha_0,\bzero,\btheta_0)$, that is, they are the solutions of $\lambda_{j}g_{j}(\bv_2)=\int_{0}^{\infty}h(\bv_1,\bv_2)g_{j}(\bv_1)f(\bv_1)d\bv_1$ for non-zero $g_{j}$.

Thus, by Theorem 9.5.2 of \cite{Resnick2019}, it remains to show that for each $t\in\bR$,
\begin{align}\label{eq:charast}
	\left|\rP\left[e^{it (n-1)\tT_{n}^*}\right]-\rP\left[e^{it \zeta}\right]\right|\rightarrow 0,
\end{align}
as $n$ goes to infinity.

Let $\{\bV_{0i},i=1,\cdots,n\}$ be independent and identically distributed sample for each $\bV_{0i}$ from the null distribution with density $f(\bv,\balpha_0,\bzero,\btheta_0)$, where $\{\bV_{0i},i=1,\cdots,n\}$ are independent of the observed sample $\{\bV_{i},i=1,\cdots,n\}$. As shown in Theorem in Section 5.5.2 of \cite{Serfling1980}, we rewrite, by the representation of $h(\bv_1,\bv_2)$ (Page 1087 in \cite{Linear_part2}),
\begin{align*}
	h(\bv_1,\bv_2) = \sum_{k=1}^{\infty}\lambda_{k}\phi_k(\bv_1)\phi_k(\bv_2),
\end{align*}
where $\{\phi_k(\cdot)\}$ and $\{\lambda_k\}$ have following properties that
\begin{align*}
	\rP[\phi_j(\bV_{01})\phi_k(\bV_{01})] = \left\{
	\begin{array}{ll}
		1, & \text{ if } j = k, \\
		0,& \text{ otherwise},
	\end{array}
	\right.
\end{align*}
\begin{align*}
	\rP[\phi_k(\bV_{01})]=0, \quad \mbox{for all } k,
\end{align*}
and
\begin{align*}
	\rP[h(\bV_{01},\bV_{02})^2]=\sum_{k=1}^{\infty}\lambda_k^2<\infty.
\end{align*}

Let
\begin{align*}
	T_{nK} = \frac{1}{n}\sum_{i\neq j}\sum_{k=1}^{K}\lambda_{k}\left(\phi_k(\bV_i^*)-\rP[\phi_k(\bV_i^*)]\right)\left(\phi_k(\bV_j^*)-\rP[\phi_k(\bV_j^*)]\right).
\end{align*}
First, we are going to show that for any fixed $K$,
\begin{align}\label{eq:tTnks}
	T_{nK}\lkonv Y_K=\sum_{k=1}^{K}\lambda_k(W_k^2-1),
\end{align}
where $\{W_k,k=1,\cdots,K\}$ are independent and each $W_k$ is normally distributed with mean $\mu_{ak}$ and variance $1$.
We can rewrite
\begin{align*}
	T_{nK}= \sum_{k=1}^{K}\lambda_k(W_{nk}^2-Z_{nk}),
\end{align*}
where
\begin{align*}
	W_{nk}= \frac{1}{\sqrt{n}}\sum_{i=1}^{n}\left(\phi_k(\bV_i^*)-\rP[\phi_k(\bV_i^*)]\right),
\end{align*}
\begin{align*}
	Z_{nk}= \frac{1}{n}\sum_{i=1}^{n}\left(\phi_k(\bV_i^*)-\rP[\phi_k(\bV_i^*)]\right)^2.
\end{align*}
By the mean value theorem and combing the dominated convergence theorem and Assumption (A7), it is seen that
\begin{align}\label{eq:phi1}
	\begin{split}
	\rP[\phi_k(\bV_i^*)]
	=&\rP\left[\phi_k(\bV_{0i})\frac{f(\bV_{0i},\hbalpha,\bzero,\btheta_0)}{f(\bV_{0i},\balpha_0,\bzero,\btheta_0)}\right]\\	=&\frac{1}{n^{1/2}}\rP\left[\phi_k(\bV_{0i})\bxi_1\trans\frac{\partial f(\bV_{0i},\hbalpha',\bzero,\btheta_0)/\partial\balpha}{f(\bV_{0i},\balpha_0,\bzero,\btheta_0)}\right]\\
=&\frac{1}{n^{1/2}}\rP\left[\phi_k(\bV_{0i})\bxi_1\trans\frac{\partial f(\bV_{0i},\balpha_0,\bzero,\btheta_0)/\partial\balpha}{f(\bV_{0i},\balpha_0,\bzero,\btheta_0)}\right](1+o(1))\\
	=&\frac{1}{n^{1/2}}\mu_{ak}+o(n^{-1/2}),
	\end{split}
\end{align}
where $\hbalpha'$ is between $\bzero$ and $\hbalpha$, and then, we have
\begin{align*}
	\left|\rP[\phi_k(\bV_i^*)]\right|
	\leq\frac{1}{\sqrt{n}}C_f(\bxi_1),
\end{align*}
and
\begin{align}\label{eq:phi2}
	\begin{split}
	\rP[\phi_k(\bV_i^*)^2]
	=&\rP\left[\phi_k(\bV_{0i})^2\frac{f(\bV_{0i},\hbalpha,\bzero,\btheta_0)}{f(\bV_{0i},\balpha_0,\bzero,\btheta_0)}\right]\\
	\leq&1+\frac{1}{\sqrt{n}}\bE\left[\phi_k(\bV_{0i})^2b_1(\bV_{0i},\bxi_1)\right]\\
	=&1+\frac{1}{\sqrt{n}}C_f(\bxi_1).
	\end{split}
\end{align}
Thus, it can be obtained that
\begin{align*}
	\rP[W_{nk}]= 0,
\end{align*}
\begin{align*}
	\Var(W_{nk})= \frac{1}{n}\sum_{i=1}^{n}\Var(\phi_k(\bV_i^*))=1+O(n^{-1/2})-\frac{1}{n}\mu_{ak}^2\rightarrow1,
\end{align*}
and if $j\neq k$, by assumption (A7),
\begin{align*}
	\Cov(W_{nj}, W_{nk})
	=& \rP[W_{nj}W_{nk}] - \rP[	W_{nj}]\rP[	W_{nk}]\\
	=&\frac{1}{n}\sum_{i=1}^{n}\left\{\rP[\phi_j(\bV_i^*)\phi_k(\bV_i^*)]-\rP[\phi_j(\bV_i^*)]\rP[\phi_k(\bV_i^*)]\right\}\\
	=&\rP[\phi_j(\bV_i^*)\phi_k(\bV_i^*)]-\rP[\phi_j(\bV_i^*)]\rP[\phi_k(\bV_i^*)]\\ =&\frac{1}{n^{1/2}}\rP\left[\phi_j(\bV_{0i})\phi_k(\bV_{0i})\bxi_1\trans\frac{\partial f(\bV_{0i},\hbalpha,\bzero,\btheta_0)/\partial\balpha}{f(\bV_{0i},\balpha_0,\bzero,\btheta_0)}\right]\\
	&-\rP[\phi_j(\bV_i^*)]\rP[\phi_k(\bV_i^*)]\\
	\leq&\frac{1}{n^{1/2}}\rP\left[|\phi_j(\bV_{0i})\phi_k(\bV_{0i})|b_1(\bV_{0i},\bxi_1)\right]+\frac{1}{n}C_f(\bxi_1)\\
	=&O(n^{-1/2})+O(n^{-1})\\
	\rightarrow&0,
\end{align*}
which implies that if $j\neq k$,
\begin{align*}
	\left|\Cov(W_{nj}, W_{nk})\right|\rightarrow0,
\end{align*}
and consequently,
\begin{align*}
	\Cov(W_{nj}, W_{nk})\rightarrow\left\{
	\begin{array}{ll}
		1, & \text{ if } j = k, \\
		0,& \text{ otherwise}.
	\end{array}
	\right.
\end{align*}
Therefore, by Lindeberg-Levy central limit theorem, we have
\begin{align*}
	(W_{n1},\cdots, W_{nK})\trans \lkonv N(\bzero_K, \bone_{K\times K}).
\end{align*}
Since $\rP[\phi_k^2(\bV_i^*)]\rightarrow1$ and $(\rP[\phi_k(\bV_i^*)])^2\rightarrow0$, it is seen that by the strong law of large numbers,
\begin{align*}
	Z_{nk}\askonv 1,
\end{align*}
which implies that
\begin{align*}
	(Z_{n1},\cdots,Z_{nK})\trans\askonv (1,\cdots,1)\trans.
\end{align*}
Thus, (\ref{eq:tTnks}) holds and consequently,
\begin{align}\label{eq:tTns-Tns2}
	\left|\rP\left[e^{it T_{nK}}\right]-\rP\left[e^{it Y_{K}}\right]\right|<\eps,
\end{align}
for all $n$ sufficiently large.

If we have (we will show this later) that for fixed $t$ and given $\eps>0$, choose and fix $K$ large enough such that for all $n$ sufficiently large,
\begin{align}\label{eq:tTns-Tns1}
	\left|\rP\left[e^{it (n-1)\tT_n^*}\right]-\rP\left[e^{it T_{nK}}\right]\right|<6(C_f(\bxi_1)+1)\eps,
\end{align}
which combing with (\ref{eq:tTns-Tns2}) implies that, for all $n$ large enough,
\begin{align*}
	\left|\rP\left[e^{it (n-1)\tT_{n}^*}\right]-\rP\left[e^{it \zeta}\right]\right|
	\leq &\left|\rP\left[e^{it (n-1)\tT_{n}^*}\right]-\rP\left[e^{it T_{nK}}\right]\right|
	+\left|\rP\left[e^{it T_{nK}}\right]-\rP\left[e^{it Y_{K}}\right]\right|\\
	&+\left|\rP\left[e^{it Y_{K}}\right]-\rP\left[e^{it \zeta}\right]\right|\\
	<&(2C_f^2(\bxi_1)+2)^{1/2}\eps+ \left|\rP\left[e^{it Y_{K}}\right]-\rP\left[e^{it \zeta}\right]\right|.
\end{align*}
Therefore, to obtain the inequality (\ref{eq:charast}), it suffices to establish that for all $K$ sufficiently large,
\begin{align}\label{eq:Yks}
	\left|\rP\left[e^{it Y_{K}}\right]-\rP\left[e^{it \zeta}\right]\right|< \eps\left\{\rP[(W_1^2-1)^2]\right\}^{1/2}.
\end{align}
Let the random variable $W_1$, $W_2$, $\cdots$ be defined on a common probability space and denote by $\zeta$ the limit in mean square of $Y_K$ as $K\rightarrow\infty$. Then according to the fact that $|e^{it}-1|\leq |t|$, we have by the definition of $Y_K$,
\begin{align*}
	\left|\rP\left[e^{it Y_{K}}\right]-\rP\left[e^{it \zeta}\right]\right|
	\leq& |t|\left[\rP\left(Y_{K}-\zeta\right)^2\right]^{1/2}\\
	\leq&|t|\left[\rP\left(W_1^2-1\right)^2\right]^{1/2}\left[\sum_{k=K+1}^{\infty}\lambda_k^2\right]^{1/2}\\
	<&\eps\left[\rP\left(W_1^2-1\right)^2\right]^{1/2},
\end{align*}
which arrives at (\ref{eq:Yks}). Therefore, the inequality (\ref{eq:charast}) holds.

The proof will be completed if we claim that (\ref{eq:tTns-Tns1}) holds.
To facilitate the expression, we give some additional notations below. Let
\begin{align*}
	\tT_{nK} = \frac{1}{n}\sum_{i\neq j}\sum_{k=1}^{K}\lambda_{k}\phi_k(\bV_i^*)\phi_k(\bV_j^*),
\end{align*}
\begin{align*}
	R_{1nK} = (n-1)\sum_{k=1}^{K}\lambda_{k}(\rP[\phi_k(\bV_i^*)])^2,
\end{align*}
\begin{align*}
	R_{2nK} = \frac{n-1}{n}\sum_{i=1}^n\sum_{k=1}^{K}\lambda_{k}\rP[\phi_k(\bV_j^*)]\phi_k(\bV_i^*),
\end{align*}
\begin{align*}
	T_{0n}^* = \frac{1}{n}\sum_{i\neq j}h(\bV_i^*,\bV_j^*),
\end{align*}
and
\begin{align*}
	\tT_{0n}^* = \frac{1}{n}\sum_{i\neq j}\rP[h(\bV_i^*,\bV_j^*)|\bV_i^*).
\end{align*}
Thus, $(n-1)\tT_n^*-T_{nK}$ can be decomposed as follows,
\begin{align*}
	(n-1)\tT_n^*-T_{nK}
	=&(n-1)\tT_n^*-T_{0n}^*+2\tT_{0n}^*-\rP[T_{0n}^*]\\
	 &+T_{0n}^*-\tT_{nK}\\
	 &+\rP[T_{0n}^*]-R_{1nK}\\
	 &+2R_{2nK}-2\tT_{0n}^*\\
\equiv&\Pi_1+\Pi_2+\Pi_3+\Pi_4,
\end{align*}
where
\begin{align*}
\Pi_1=&\frac{1}{n}\sum_{i\neq j}\big\{h(\bV_i^*,\bV_j^*,\hbalpha)-h(\bV_i^*,\bV_j^*)\\
 &+\rP[h(\bV_i^*,\bV_j^*)|\bV_j^*]+\rP[h(\bV_i^*,\bV_j^*)|\bV_i^*]-\rP[h(\bV_i^*,\bV_j^*)]\big\},
\end{align*}
and
\begin{align*}
	\Pi_2=&T_{0n}^*-\tT_{nK},\\
	\Pi_3=&\rP[T_{0n}^*]-R_{1nK},\\
	\Pi_4=&2R_{2nK}-2\tT_{0n}^*.
\end{align*}
We will deal with $\Pi_j$ one by one as follows, $j=1,2,3,4$.

Before showing the upper bounds, we derive following useful inequalities. By the mean value theorem, it is seen that
\begin{align*}
	\rP[\psi(\bV_i^*,\balpha_0,  0, \btheta)]
	=\rP\left[\psi(\bV_{0i},\balpha_0,  0, \btheta)\frac{f(\bV_{0i},\hbalpha,\bzero,\btheta_0)}{f(\bV_{0i},\balpha_0,\bzero,\btheta_0)}\right],
\end{align*}
and consequently, by  Assumption (A7) and Cauchy-Schwarz inequality, there are constant $C_1>0$ such that
\begin{align*}
	\|\rP[\psi(\bV_i^*,\balpha_0,  0, \btheta)]\|
	\leq&\frac{C_f(\bxi_1)}{\sqrt{n}}\left\{\rP\left[\|\psi(\bV_{0i},\balpha_0,  0, \btheta)\|^2\right]\right\}^{1/2}\\
	\leq&\frac{1}{\sqrt{n}}C_1C_f(\bxi_1).
\end{align*}
And it is followed from above inequality that
\begin{align*}
	\rP[\tpsi_{ij}^*(\balpha_0)]
	=&\int_{\btheta\in\Theta_\theta}\left\|\rP[\psi(\bV_i^*,\balpha_0,  0, \btheta)]\right\|^2w(\btheta) d\btheta\\
	\leq&\frac{1}{n}C_1^2C_f^2(\bxi_1),
\end{align*}
and
\begin{align*}
	\rP[K_i^*]
	=&\int_{\btheta\in\Theta_\theta}J\trans K(\btheta)\trans\rP[\psi(\bV_i^*,\balpha_0,  0, \btheta)]w(\btheta) d\btheta,
\end{align*}
and by Assumption (A3), there is a constant $C_3>0$ such that
\begin{align*}
	\|\rP[K_i^*]\|^2
	\leq&\int_{\btheta\in\Theta_\theta}\|J\trans K(\btheta)\trans\rP[\psi(\bV_i^*,\balpha_0,  0, \btheta)]\|^2w(\btheta) d\btheta\\
	\leq&\int_{\btheta\in\Theta_\theta}\lambda_{\max}\left(\trace\left(K(\btheta)JJ\trans K(\btheta)\trans\right)\right)\|\rP\left[\psi(\bV_i^*,\balpha_0,  0, \btheta)\right]\|^2w(\btheta) d\btheta\\
	\leq&\frac{1}{n}C_1^2C_f^2(\bxi_1)\int_{\btheta\in\Theta_\theta}\lambda_{\max}\left(\trace\left(K(\btheta)JJ\trans K(\btheta)\trans\right)\right)w(\btheta) d\btheta\\
	\leq&\frac{1}{n}C_1^2C_3^2C_f^2(\bxi_1).
\end{align*}
As the same arguments, we have
\begin{align*}
	\rP[\psi_1(\bV_i^*,\balpha_0)]
	=\rP\left[\psi_1(\bV_{0i},\balpha_0)\frac{f(\bV_{0i},\hbalpha,\bzero,\btheta_0)}{f(\bV_{0i},\balpha_0,\bzero,\btheta_0)}\right],
\end{align*}
and there are constant $C_2>0$,
\begin{align*}
	\|\rP[\psi_1(\bV_i^*,\balpha_0)]\|
	\leq&\frac{C_f(\bxi_1)}{\sqrt{n}}\left\{\rP\left[\|\psi_1(\bV_{0i},\balpha_0)\|^2\right]\right\}^{1/2}\\
	\leq&\frac{1}{\sqrt{n}}C_2C_f(\bxi_1).
\end{align*}
And consequently, we have
\begin{align}\label{eq:hxys}
	\begin{split}
		|\rP[h(\bV_i^*,\bV_j^*)]|
		=&\bigg|\int_{\btheta\in\Theta_\theta}\left\|\rP[\psi(\bV_i^*, \balpha_0, 0, \btheta)]\right\|^2w(\btheta) d\btheta\\
		&+\rP\left[\psi_1(\bV_i^*, \balpha_0)\trans K_{j}^*+ \psi_1(\bV_j^*, \balpha_0)\trans K_i^*\right]\\
        &+ \rP\left[\psi_1(\bV_i^*, \balpha_0)\trans H\psi_1(\bV_j^*, \balpha_0)\right]\bigg|\\
		\leq&\frac{1}{n}\left[C_1^2+2\frac{1}{\sqrt{n}}C_1C_2C_3+C_2^2\lambda_{\max}(H)\right]C_f^2(\bxi_1)\\
		=&O(n^{-1}).
	\end{split}
\end{align}

{\bf To bound $\Pi_1$}. Rewrite
\begin{align*}
	\Pi_1=(n-1)U_{1n},
\end{align*}
where
\begin{align*}
	U_{1n}=\frac{1}{n(n-1)}\sum_{i\neq j}\check{h}(\bV_i^*,\bV_j^*),
\end{align*}
and
\begin{align*}
	\check{h}(\bV_i^*,\bV_j^*)
	=&h(\bV_i^*,\bV_j^*,\hbalpha)-h(\bV_i^*,\bV_j^*)\\
	&+\rP[h(\bV_i^*,\bV_j^*)|\bV_j^*)+\rP[h(\bV_i^*,\bV_j^*)|\bV_i^*)-\rP[h(\bV_i^*,\bV_j^*)].
\end{align*}
Since $U_{1n}$ is a U-statistic with  $\rP[\check{h}(\bV_i^*,\bV_j^*)|\bV_j^*]=0$, we have, by Lemma A in Section 5.2.1 of \cite{Serfling1980},
\begin{align*}
	\rP[\Pi_1^2]
	=&(n-1)^2\Var(U_{1n})\\
	=&(n-1)^2\binom{n}{2}^{-1}\Var(\check{h}(\bV_i^*,\bV_j^*))\\
	\leq&2\Var(\check{h}(\bV_i^*,\bV_j^*))\\
	=&2\rP[\check{h}(\bV_i^*,\bV_j^*)^2]\\
	\leq&4\rP\left[\left(h(\bV_i^*,\bV_j^*,\hbalpha)-h(\bV_i^*,\bV_j^*)+\rP[h(\bV_i^*,\bV_j^*)]\right)^2\right]\\
	&+16\rP\left[\left(\rP[h(\bV_i^*,\bV_j^*)|\bV_i^*)-\rP[h(\bV_i^*,\bV_j^*)]\right)^2\right].
\end{align*}
According to Assumption (A7) and the mean value theorem, it is seen that for $n$ large enough,
\begin{align*}
	\rP&\left[\left|h(\bV_i^*,\bV_j^*,\hbalpha)-h(\bV_i^*,\bV_j^*)\right|^2\right]\\ =&\rP\left[\left|h(\bV_{0i},\bV_{0j},\hbalpha)-h(\bV_{0i},\bV_{0j})\right|^2\frac{f(\bV_{0i},\hbalpha,\bzero,\btheta_0)}{f(\bV_{0i},\balpha_0,\bzero,\btheta_0)}\frac{f(\bV_{0j},\hbalpha,\bzero,\btheta_0)}{f(\bV_{0j},\balpha_0,\bzero,\btheta_0)}\right]\\
	\leq&\frac{1}{n}\rP\left[\left(\bxi_1\trans\frac{\partial h(\bV_{0i},\bV_{0j},\hbalpha)}{\partial\balpha}\right)^2(1+b_1(\bV_{0i},\bxi_1))(1+b_1(\bV_{0j},\bxi_1))\right],
\end{align*}
where $\hbalpha$ is between $\bzero$ and $\hbalpha$. By utilizing Assumption (A7), there is a constant $C_{\bxi_1}>0$ such that
\begin{align*}
	\rP\left[\left(\bxi_1\trans\frac{\partial h(\bV_{0i},\bV_{0j},\hbalpha)}{\partial\balpha}\right)^2\right]\leq& C_{\bxi_1},\\
	\rP\left[\left(\bxi_1\trans\frac{\partial h(\bV_{0i},\bV_{0j},\hbalpha)}{\partial\balpha}\right)^2b_1(\bV_{0i},\bxi_1)\right]\leq& C_{\bxi_1},\\
	\rP\left[\left(\bxi_1\trans\frac{\partial h(\bV_{0i},\bV_{0j},\hbalpha)}{\partial\balpha}\right)^2b_1(\bV_{0i},\bxi_1)b_1(\bV_{0j},\bxi_1)\right]\leq& C_{\bxi_1},
\end{align*}
which implies that
\begin{align*}
	\rP\left[\left|h(\bV_i^*,\bV_j^*,\hbalpha)-h(\bV_i^*,\bV_j^*)\right|^2\right]=O(n^{-1}).
\end{align*}

By the foregoing arguments, we have
\begin{align*}
	\rP[h(\bV_i^*,\bV_j^*)|\bV_i^*]
	=&\int_{\btheta\in\Theta_\theta}\rP[\psi(\bV_j^*,\balpha_0,  0, \btheta)]\trans \psi(\bV_i^*,\balpha_0,  0, \btheta) w(\btheta) d\btheta\\
	&+\psi_1(\bV_i^*, \balpha_0)\trans\rP[K_{j}^*]+ \rP[\psi_1(\bV_j^*, \balpha_0)]\trans K_i^*\\
    &+ \psi_1(\bV_i^*, \balpha_0)\trans H\rP[\psi_1(\bV_j^*, \balpha_0)],
\end{align*}
and consequently, for $n$ large enough,
\begin{align*}
	\rP&\left[\left(\rP[h(\bV_i^*,\bV_j^*)|\bV_i^*)]\right)^2\right]\\
	\leq&3\int_{\btheta\in\Theta_\theta}\rP\left[\{\rP[\psi(\bV_j^*,\balpha_0,  0, \btheta)]\trans \psi(\bV_i^*,\balpha_0,  0, \btheta)\}^2\right] w(\btheta) d\btheta\\
	&+6\rP\left[\{[\psi_1(\bV_i^*, \balpha_0)]\trans\rP[K_{j}^*]\}^2\right]\\
	&+ 3\rP\left[\{[\psi_1(\bV_i^*, \balpha_0)]\trans H\rP[\psi_1(\bV_j^*, \balpha_0)]\}^2\right]\\
	\leq&3\int_{\btheta\in\Theta_\theta}\|\rP[\psi(\bV_j^*,\balpha_0,  0, \btheta)]\|^2 \rP\left[\|\psi(\bV_i^*,\balpha_0,  0, \btheta)\|^2\right] w(\btheta) d\btheta\\
	&+6\rP\left[\|\psi_1(\bV_i^*, \balpha_0)\|^2\right]\|\rP[K_{j}^*]\|^2\\
	&+ 3\lambda_{\max}(H)^2\rP\left[\|\psi_1(\bV_i^*, \balpha_0)\|^2\right]\|\rP\left[\psi_1(\bV_j^*, \balpha_0)\right]\|^2\\
	\leq&3C_1^2C_f^2(\bxi_1)\frac{1}{n}\left(1+\frac{1}{\sqrt{n}}C_f(\bxi_1)C_M\right)\\
	&+6C_1^2C_3^2C_f^2(\bxi_1)\frac{1}{n}C_M+3C_2^2C_f^2(\bxi_1)\frac{1}{n}C_M\\
	=&O(n^{-1}),
\end{align*}
where $C_M$ is a constant.

Combining the above arguments with inequality (\ref{eq:hxys}), we have
\begin{align}\label{eq:pi1}
	\begin{split}
	\rP[\Pi_1^2]
	\leq&8\rP\left[\left\{h(\bV_i^*,\bV_j^*,\hbalpha)-h(\bV_i^*,\bV_j^*)\right\}^2\right]+8\left(\rP[h(\bV_i^*,\bV_j^*)]\right)^2\\
	&+32\rP\left[\left\{\rP[h(\bV_i^*,\bV_j^*)|\bV_i^*)\right\}^2\right]+32\left(\rP[h(\bV_i^*,\bV_j^*)]\right)^2\\
	=&O(n^{-1}).
	\end{split}
\end{align}

{\bf To bound $\Pi_2$}.
Noting that $T_{0n}^*-\tT_{nK}$ has basically the form of a U-statistic, that is,
\begin{align*}
	T_{0n}^*-\tT_{nK}=\frac{2}{n}\binom{n}{2}U_{nK},
\end{align*}
where
\begin{align*}
	U_{nK}=\frac{1}{n(n-1)}\sum_{i\neq j}g_K(\bV_i^*,\bV_j^*)
\end{align*}
with
\begin{align*}
	g_K(\bV_i^*,\bV_j^*)
	=&h(\bV_i^*,\bV_j^*)-\sum_{k=1}^{K}\lambda_{k}\phi_k(\bV_i^*)\phi_k(\bV_j^*)\\
	=&\sum_{k=K+1}^{\infty}\lambda_{k}\phi_k(\bV_i^*)\phi_k(\bV_j^*).
\end{align*}
Next we claim that
\begin{align}\label{eq:pi2}
	\rP\left(T_{0n}^*-\tT_{nK}\right)^2\leq 4(2C_f^2(\bxi_1)+1)\sum_{k=K+1}^{\infty}\lambda_{k}^2.
\end{align}

From inequality (\ref{eq:phi1}) and (\ref{eq:phi2}), it is straightforward to obtain that for large $n$,
\begin{align*}
	\rP[g_K(\bV_i^*,\bV_j^*)]
	=&\sum_{k=K+1}^{\infty}\lambda_{k}(\rP[\phi_k(\bV_i^*)])^2\\
	=&\frac{1}{n}\left(\sum_{k=K+1}^{\infty}\lambda_{k}\tmu_{ak}^2\right)(1+o(1)),
\end{align*}
where
\begin{align*}
	\tmu_{ak}=\rP\left[\phi_k(\bV_{0i})\bxi_2\trans\frac{\partial h(\bV_1,\bV_2, \hbalpha)}{\partial\balpha}\right].
\end{align*}
Since by Lemma \ref{lem:h2b12} combining with Assumption (A7) and Cauchy-Schwarz inequality,
\begin{align}\label{eq:gk2a}
	\begin{split}
	\rP&[g_K(\bV_i^*,\bV_j^*)^2]\\
	=&\int\int g_K(\bv_1,\bv_2)^2f(\bv_1,\hbalpha,\bzero,\btheta_0)f(\bv_2,\hbalpha,\bzero,\btheta_0)\bv_1d\bv_2\\
	\leq&\rP[g_K(\bV_{0i},\bV_{0j})^2(1+n^{-1/2}b_1(\bV_{0i},\bxi_1))(1+n^{-1/2}b_1(\bV_{0j},\bxi_1))]\\
	\leq&\rP[g_K(\bV_{0i},\bV_{0j})^2(1+n^{-1/2}b_1(\bV_{0i},\bxi_1))^2]\\
	=&\sum_{k=K+1}^{\infty}\lambda_{k}^2\rP[\phi_k(\bV_{0i})^2(1+n^{-1/2}b_1(\bV_{0i},\bxi_1))^2]\rP[\phi_k(\bV_{0j})^2]\\
	=&\sum_{k=K+1}^{\infty}\lambda_{k}^2\rP[\phi_k(\bV_{0i})^2(1+n^{-1/2}b_1(\bV_{0i},\bxi_1))^2]\\
	\leq&(1+3Cn^{-1/2})\sum_{k=K+1}^{\infty}\lambda_{k}^2\\	
	<&\infty,
	\end{split}
\end{align}
where $C$ is a constant such that $\rP[\phi_k(\bV_i)^2b_1(\bV_{0i},\bxi_1)^2]\leq C$,
thus, it is straightforward to obtain that by the dominated theorem and for large $n$,
\begin{align*}
	\rP[g_K(\bV_i,\bV_j)^2]
	=&\int\int g_K(\bv_1,\bv_2)^2f(\bv_1,\hbalpha,\bzero,\btheta_0)f(\bv_2,\hbalpha,\bzero,\btheta_0)d\bv_1d\bv_2\\
	\leq&\rP[g_K(\bV_{0i},\bV_{0j})^2](1+o(1))\\
	=&(1+o(1))\sum_{k=K+1}^{\infty}\lambda_{k}^2\\
	\leq&2\sum_{k=K+1}^{\infty}\lambda_{k}^2,
\end{align*}
which implies that
\begin{align*}
	\Var(g_K(\bV_i^*,\bV_j^*))
	\leq\rP[g_K(\bV_i^*,\bV_j^*)^2]
	\leq2\sum_{k=K+1}^{\infty}\lambda_{k}^2.
\end{align*}
Let $b(\bv)=\rP[g_K(\bV_i^*,\bV_j^*)|\bV_j^*=\bv]$. Noting that
\begin{align*}
	b(\bv)=&\rP[g_K(\bV_i^*,\bV_j^*)|\bV_j^*=\bv]\\
	=&\sum_{k=K+1}^{\infty}\lambda_{k}\phi_k(\bv)\rP[\phi_k(\bV_i^*)],
\end{align*}
we have as similar to (\ref{eq:gk2a}), for $n$ large enough,
\begin{align*}
	\rP\left[b(\bV_i^*)^2\right]
	=&\rP\left[b(\bV_{0i})^2\frac{f(\bV_{0i},\hbalpha,\bzero,\btheta_0)}{f(\bV_{0i},\balpha_0,\bzero,\btheta_0)}\right]\\
	\leq&\rP\left[b(\bV_{0i})^2(1+n^{-1/2}b_1(\bV_{0i},\bxi_1))\right]\\
	=&\rP\left[\left(\rP[g_K(\bV_i,\bV_j)|\bV_i=\bV_{0i}]\right)^2(1+n^{-1/2}b_1(\bV_{0i},\bxi_1))\right]\\
	\leq&\rP\left[\rP[g_K(\bV_i,\bV_j)^2(1+n^{-1/2}b_1(\bV_{0i},\bxi_1))|\bV_i=\bV_{0i}]\right]\\
	=&\rP\left[g_K(\bV_{0i},\bV_j)^2(1+n^{-1/2}b_1(\bV_{0i},\bxi_1))\right]\\
	=&\rP\left[g_K(\bV_{0i},\bV_{0j})^2(1+n^{-1/2}b_1(\bV_{0i},\bxi_1))\frac{f(\bV_{0j},\hbalpha,\bzeta,\btheta_0)}{f(\bV_{0j},\balpha_0,\bzero,\btheta_0)}\right]\\
	\leq&\rP\left[g_K(\bV_{0i},\bV_{0j})^2(1+n^{-1/2}b_1(\bV_{0i},\bxi_1))(1+n^{-1/2}b_1(\bV_{0j},\bxi_1))\right]\\
	<&\infty,
\end{align*}
and consequently by the dominated theorem,
\begin{align*}
	\Var(\rP[g_K(\bV_i^*,\bV_j^*)|\bV_i^*])
	\leq&\rP\left[b(\bV_i^*)^2\right]\\
	=&\rP\left[b(\bV_{0i})^2\right](1+o(1))\\
	=&(1+o(1))\sum_{k=K+1}^{\infty}\lambda_{k}^2(\rP[\phi_k(\bV_j)])^2\\
	\leq&\frac{2}{n}C_f^2(\bxi_1)\sum_{k=K+1}^{\infty}\lambda_{k}^2.
\end{align*}
Lemma A in Section 5.2.1 of \cite{Serfling1980} states that
\begin{align*}
	\Var(U_{nK})
	\leq&\binom{n}{2}^{-1}\left(4C_f^2\sum_{k=K+1}^{\infty}\lambda_{k}^2+2\sum_{k=K+1}^{\infty}\lambda_{k}^2\right)\\
	\leq&2(2C_f^2(\bxi_1)+1)\binom{n}{2}^{-1}\sum_{k=K+1}^{\infty}\lambda_{k}^2,
\end{align*}
which implies that
\begin{align*}
	\rP\left(\tT_n-\tT_{nK}\right)^2
	=&(n-1)^2\left(\Var(U_{nK})+(\rP[U_{nK}])^2\right)\\
	\leq&4(2C_f^2(\bxi_1)+1)\sum_{k=K+1}^{\infty}\lambda_{k}^2,
\end{align*}
and consequently (\ref{eq:pi2}) holds.

{\bf To bound $\Pi_3$}. We now show that for any large $n$,
\begin{align}\label{eq:pi3}
	\rP[|\Pi_3|]=\rP[|T_{0n}^*-R_{1nK}|]\rightarrow0,
\end{align}
as $K$ goes infinity.
It is immediately obtained that $n\rP[h(\bV_i^*,\bV_j^*)]<+\infty$ by (\ref{eq:hxys}), which implies that the series
\begin{align*}
	(n-1)\sum_{k=1}^{\infty}\lambda_{k}\left(\rP[\phi_k(\bV_i^*)]\right)^2=(n-1)\rP[h(\bV_i^*,\bV_j^*)]<\infty,
\end{align*}
and consequently that
\begin{align}\label{eq:phi3}
	\lim_{K\rightarrow\infty}(n-1)\sum_{k=K+1}^{\infty}\lambda_{k}\left(\rP[\phi_k(\bV_i^*)]\right)^2=0.
\end{align}
It follows immediately that for any large $n$,
\begin{align*}
	\bE[|T_{0n}^*-R_{1nK}|]
	=(n-1)\sum_{k=K+1}^{\infty}\lambda_{k}\left(\rP[\phi_k(\bV_i^*)]\right)^2\rightarrow0,
\end{align*}
as $K$ goes infinity.

{\bf To bound $\Pi_4$}. By the definition of $R_{2nK}$ and $\tT_{0n}^*$, it is seen that
\begin{align*}
	\frac{1}{n-1}\left(\tT_{0n}^*-R_{2nK}\right)
	=\frac{1}{n}\sum_{i=1}^{n}g_{1K}(\bV_i^*),
\end{align*}
where
\begin{align*}
	g_{1K}(\bV_i^*)
	=\sum_{k=K+1}^{\infty}\lambda_{k}\rP[\phi_k(\bV_j^*)]\phi_k(\bV_i^*).
\end{align*}
Invoking the properties of $\phi_k(\cdot)$, we have, by the inequality (\ref{eq:phi1}) and (\ref{eq:phi2}),
\begin{align*}
	\rP[g_{1K}(\bV_i^*)^2]
	=&\sum_{k=K+1}^{\infty}\lambda_{k}^2\left(\rP[\phi_k(\bV_j^*)]\right)^2\rP[\phi_k(\bV_i^*)^2]\\
	\leq&C_f^2(\bxi_1)\frac{1}{n}\sum_{k=K+1}^{\infty}\lambda_{k}^2,
\end{align*}
which implies that
\begin{align}\label{eq:pi4}
	\begin{split}
	\rP\left[\left(\tT_{0n}^*-R_{2nK}\right)^2\right]
	=&\frac{(n-1)^2}{n^2}\sum_{i=1}^{n}\Var(g_{1K}(\bV_i^*)^2)+\frac{(n-1)^2}{n}\left(\rP[g_{1K}(\bV_i^*)]\right)^2\\
	\leq&C_f^2(\bxi_1)\frac{(n-1)^2}{n^2}\sum_{k=K+1}^{\infty}\lambda_{k}^2+\frac{(n-1)^2}{n}\left(\rP[g_{1K}(\bV_i^*)]\right)^2\\
	\leq&C_f^2(\bxi_1)\sum_{k=K+1}^{\infty}\lambda_{k}^2+n\left(\rP[g_{1K}(\bV_i^*)]\right)^2.
	\end{split}
\end{align}

We have obtained the upper bounds for all $\Pi_j$, $j=1,2,3,4$. For fixed $t$ and given $\eps>0$, choose and fix $K$ large enough such that
\begin{align*}
	|t|\left(\sum_{k=K+1}^{\infty}\lambda_{k}^2\right)^{1/2}&<\eps,\quad \mbox{and, }\\
	|t|\sum_{k=K+1}^{\infty}\lambda_{k}\mu_{ak}^2&<\eps.
\end{align*}
According to (\ref{eq:phi1}), it is seen that
\begin{align*}
	\rP[g_{1K}(\bV_i^*)]
	=&\sum_{k=K+1}^{\infty}\lambda_{k}\left(\rP[\phi_k(\bV_i^*)]\right)^2\\
	=&\frac{1}{n}\sum_{k=K+1}^{\infty}\lambda_{k}\mu_{ak}^2+o(n^{-1}).
\end{align*}
Thus, by (\ref{eq:pi2}), (\ref{eq:pi3}) and (\ref{eq:pi4}), we have
\begin{align*}
	|t|\rP\left[|\Pi_2|\right]\leq&|t|\left(\rP\left[\Pi_2^2\right]\right)^{1/2}\leq2(2C_f^2(\bxi_1)+1)^{1/2}\eps\\
	|t|\rP\left[|\Pi_3|\right]\leq&\eps\\
	|t|\rP\left[|\Pi_4|\right]\leq&|t|\left(\rP\left[\Pi_4^2\right]\right)^{1/2}\leq 2(C_f(\bxi_1)+1)\eps\\
\end{align*}
According to the fact that $|e^{it}-1|\leq |t|$, we have, by the above arguments and inequality (\ref{eq:pi1}) and for all $n$ sufficiently large,
\begin{align*}
	\left|\rP\left[e^{it (n-1)\tT_n^*}\right]-\rP\left[e^{it T_{nK}}\right]\right|
	\leq& \rP\left[\left|e^{it (n-1)\tT_n^*}-e^{it T_{nK}}\right|\right]\\
	\leq& |t|\rP\left[\left|(n-1)\tT_n^*-T_{nK}\right|\right]\\
	\leq& |t|\left\{\rP[|\Pi_1|]+\rP[|\Pi_2|]+\rP[|\Pi_3|]+\rP[|\Pi_4|]\right\}\\
	<& 6(C_f(\bxi_1)+1)\eps,
\end{align*}
which establishes (\ref{eq:tTns-Tns1}) and then completes the whole proof. \hfill $\square$

\section{Simulation Studies for large \texorpdfstring{$q$}{}}\label{simulations_largeq}
\def\thetable{B.\arabic{table}}
\def\thefigure{B.\arabic{figure}}

\subsection{Change plane analysis with sparse \texorpdfstring{$\bZ$}{} for quantile, probit and semiparametric models}\label{simulation_ee_sparse}

For fair comparison, in \ref{simulation_ee_sparse}--\ref{simulation_glm_dense}, we consider the settings as similar to the existing literature except for large $q$ \citep{2011Testing,2017Change}.

We consider the probit regression model in the main paper
\begin{align*}
f(\bV_i, \balpha, \bbeta, \btheta)=\Phi(0.5+\tX_{1i}+\bX_i\trans\bbeta\bone(\bZ_i\trans\btheta\geq0))^{Y_i}+\Phi(-0.5-\tX_{1i}-\bX_i\trans\bbeta\bone(\bZ_i\trans\btheta\geq0))^{1-Y_i}.
\end{align*}
We generate $X_{1i}$ independently from binomial distribution $b(1,0.5)$, $\bX_i$ from $p$-variate normal distribution with mean $\bzero$ and variance $\sqrt{2}I$, and $Z_{1i}=1$ and $(Z_{2i},\cdots,Z_{qi})\trans$ from $(q-1)$-variate standard normal distribution.
We set $\bbeta=(1,\cdots,1)\trans$ under $H_0$. Under $H_1$, $\theta_2,\cdots,\theta_{11}$ are generated randomly from uniform distribution $U(1, 3)$ and $\theta_{12},\cdots,\theta_{q}$ are zero, where $\theta_1$ is chosen as the negative of the 0.65 percentile of $Z_2\theta_2+\cdots+Z_q\theta_q$, which means that $\bZ\trans\btheta$ divides the population into two groups with 0.35 and 0.65 observations, respectively.

Consider the quantile regression model in the main paper,
\begin{align*}
Y_i = 0.5+\tX_i+\bX_i\trans\bbeta(\tau)\bone(\bZ_i\trans\btheta(\tau)\geq0) + \eps_i,
\end{align*}
where $\tX_i$ is generated from standard normal distribution, $\bX_i$ and $\bZ_i$ are generated as same as these in probit regression model. The error $\eps_i$ is generated from $t_3$ distribution.
The parameters $\bbeta$ under $H_0$ and $\btheta$ under $H_1$ are same as these in probit model.

Consider the semiparametric model in the main paper,
\begin{align*}
Y_i = \gamma(\tbX_i,\balpha_2) +\bX_i\trans\bbeta A_i \bone(\bZ_i\trans \btheta\geq0)+\eps_i.
\end{align*}
As the same settings provided by \cite{2017Change}, we consider two settings for the propensity score model $\pi(\tbX_i,\balpha_1)$, denoted by P1 and P2:
\begin{enumerate}[({P}1)]
	\item $\pi(\tbX_i,\balpha_1)=0.5$;
	\item $\pi(\tbX_i,\balpha_1)=\frac{\exp(0.5\tX_{2i}+0.5\tX_{3i})}{1+\exp(0.5\tX_{2i}+0.5\tX_{3i})}$,
\end{enumerate}
and two baseline mean functions for $\gamma(\tbX_i,\balpha_2)$, denoted by B1 and B2:
\begin{enumerate}[({B}1)]
	\item $\gamma(\tbX_i,\balpha_2) = 1+0.5\tX_{2i}+\tX_{3i}^2$;
	\item $\gamma(\tbX_i,\balpha_2) = 1+\sin(\tX_{2i}+\tX_{3i})$.
\end{enumerate}
We calculate the WAST and SST in the main paper by fitting a linear model for baseline function $\gamma(\tbX_i,\balpha_2)$ and a logistic regression for  the propensity score model $\pi(\tbX_i,\balpha_1)$. Therefore, two baseline models $\gamma(\tbX_i,\balpha_2)$ are misspecified. We generate $v_{1i}$ independently from binomial distribution $b(1,0.5)$, $v_{2i}$ independently from uniform distribution $U(-1,1)$, and $X_{ji}$ independently from uniform distribution $U(-1,1)$ , where $i=1,\cdots,n$ and $j=1,\cdots,p$. We set $\tbX_i=(1, v_{1i}, v_{2i})\trans$ and generate $\bZ_i$ as same as these in probit regression model.
The parameters $\bbeta$ under $H_0$ and $\btheta$ under $H_1$ are same as these in probit model.

\begin{figure}[!ht]
	\begin{center}
		\includegraphics[scale=0.3]{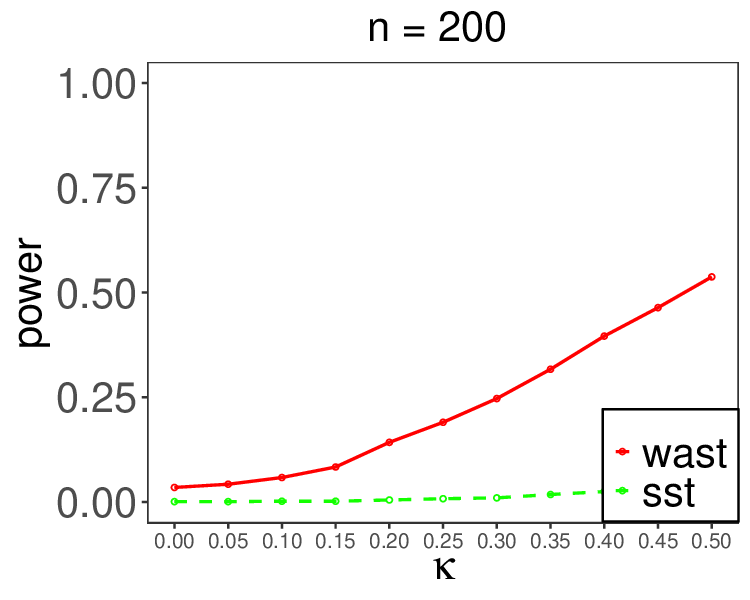}
		\includegraphics[scale=0.3]{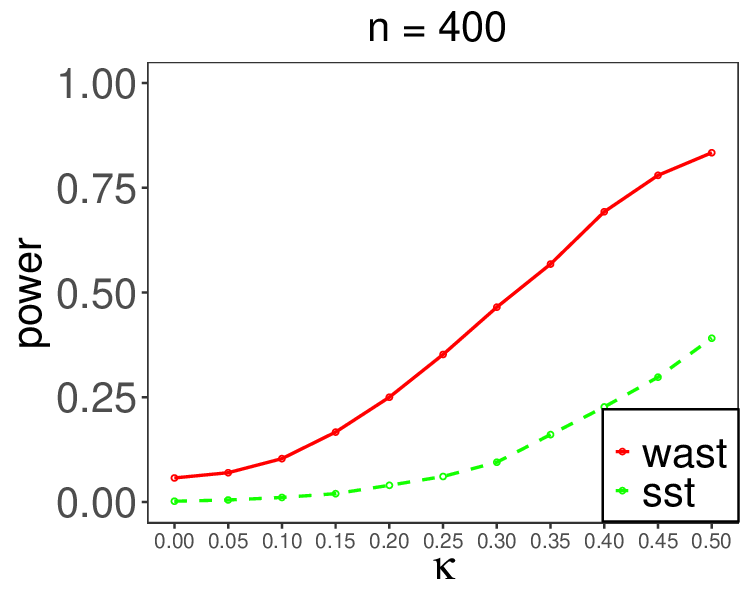}
		\includegraphics[scale=0.3]{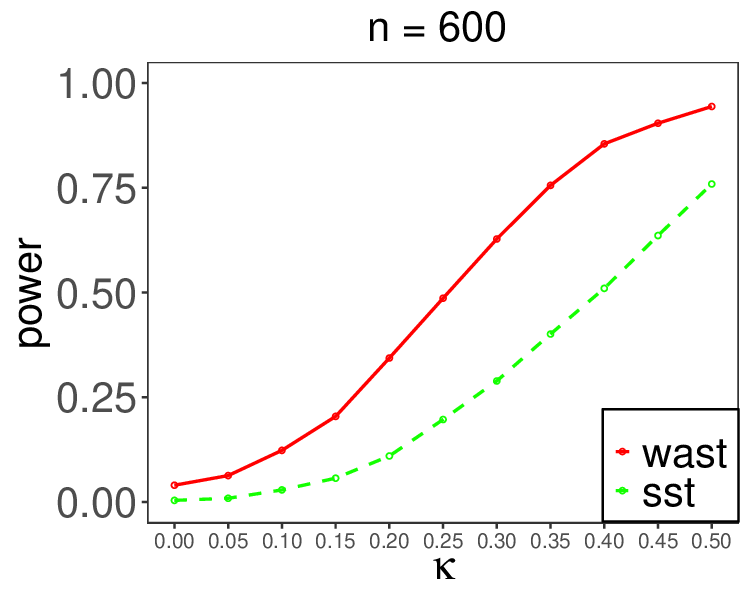}  \\
		\includegraphics[scale=0.3]{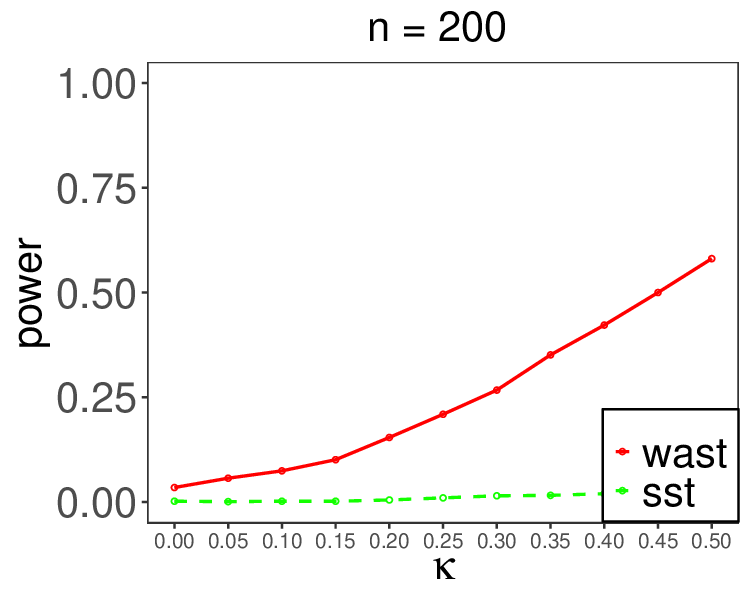}
		\includegraphics[scale=0.3]{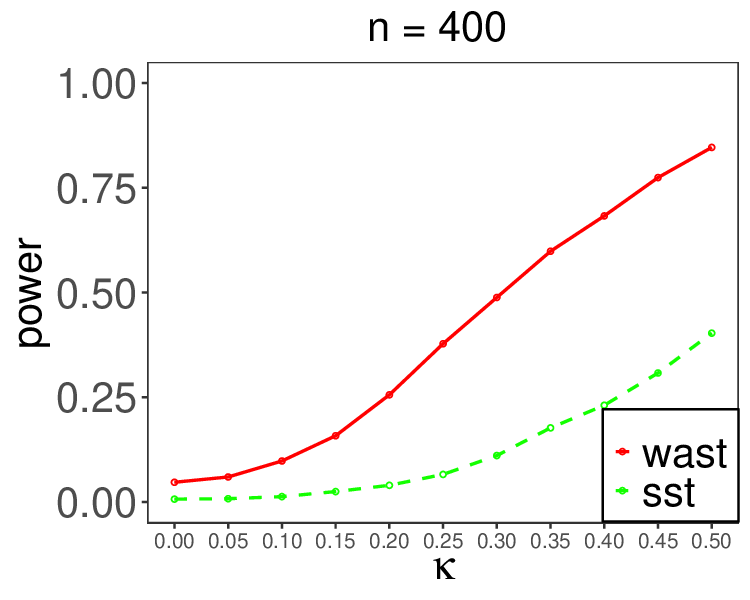}
		\includegraphics[scale=0.3]{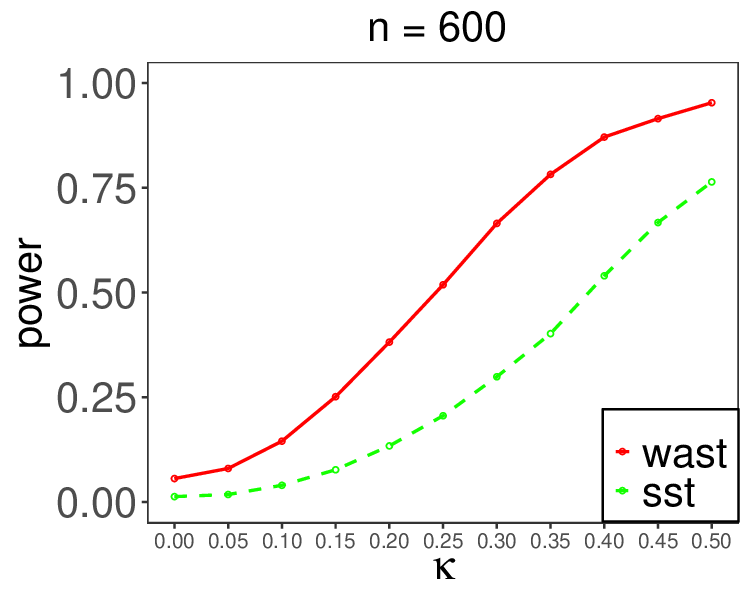}  \\
		\includegraphics[scale=0.3]{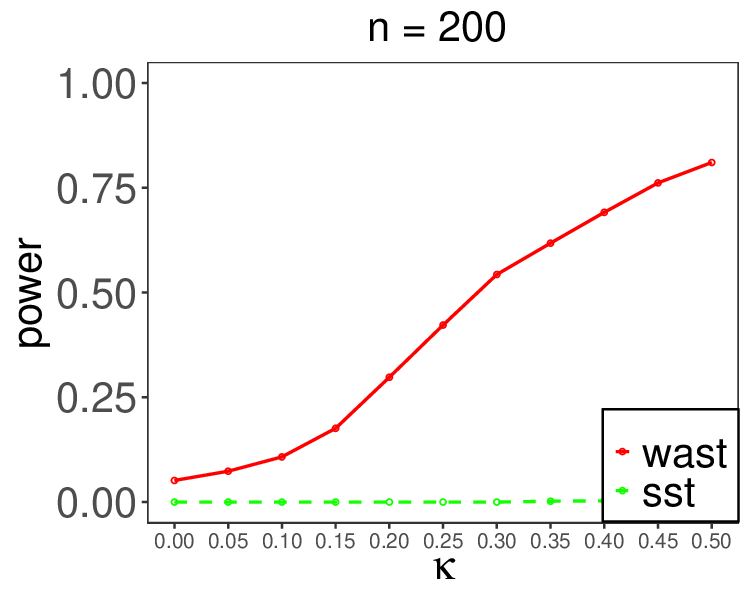}
		\includegraphics[scale=0.3]{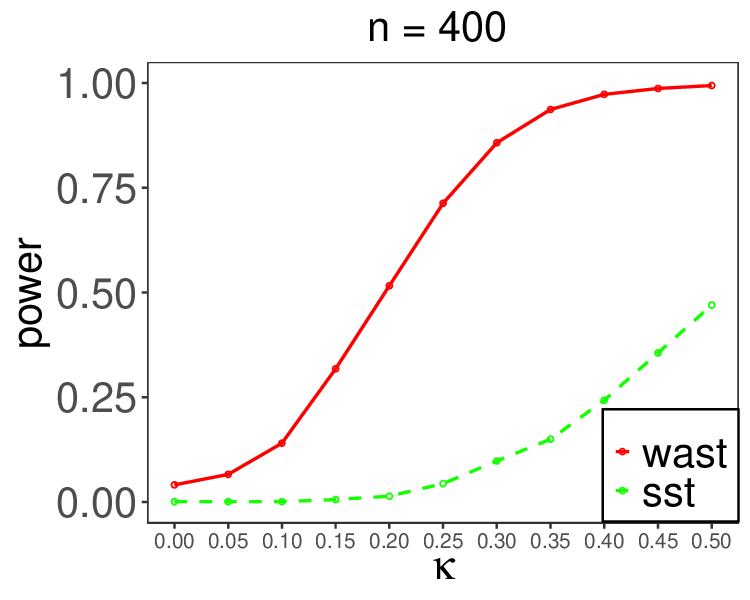}
		\includegraphics[scale=0.3]{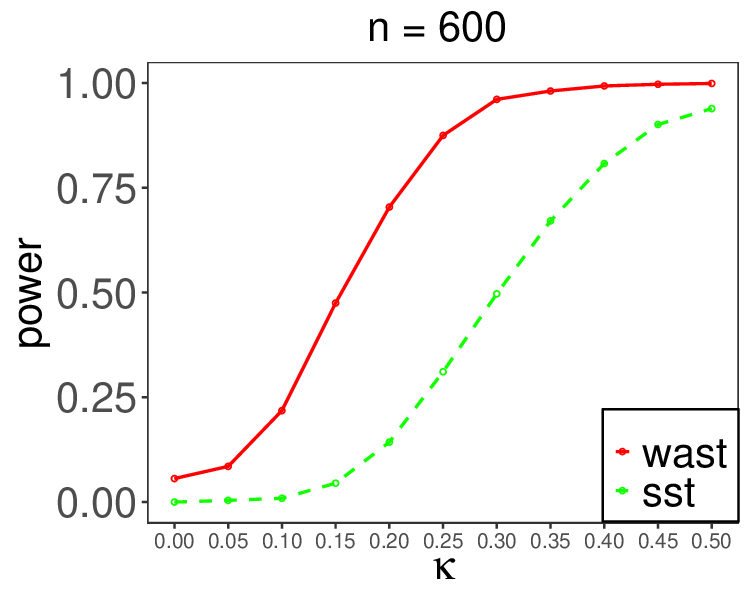}  \\
		\includegraphics[scale=0.3]{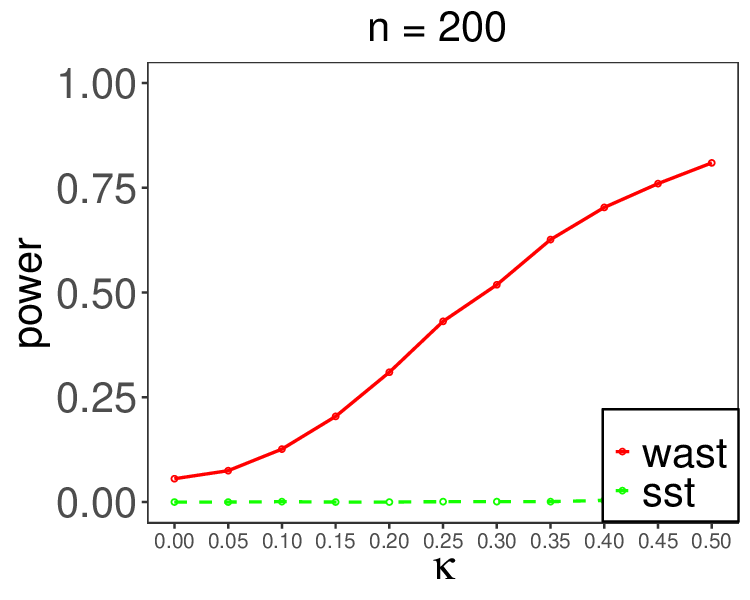}
		\includegraphics[scale=0.3]{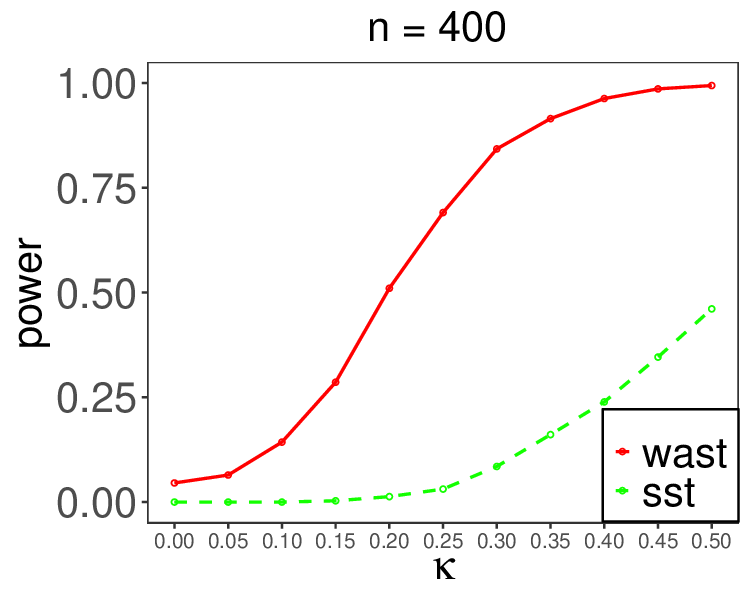}
		\includegraphics[scale=0.3]{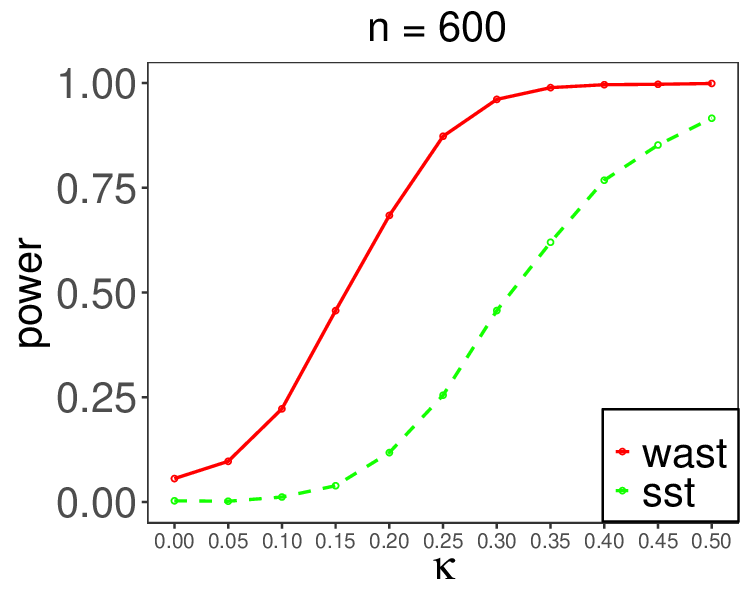}  \\
		\includegraphics[scale=0.3]{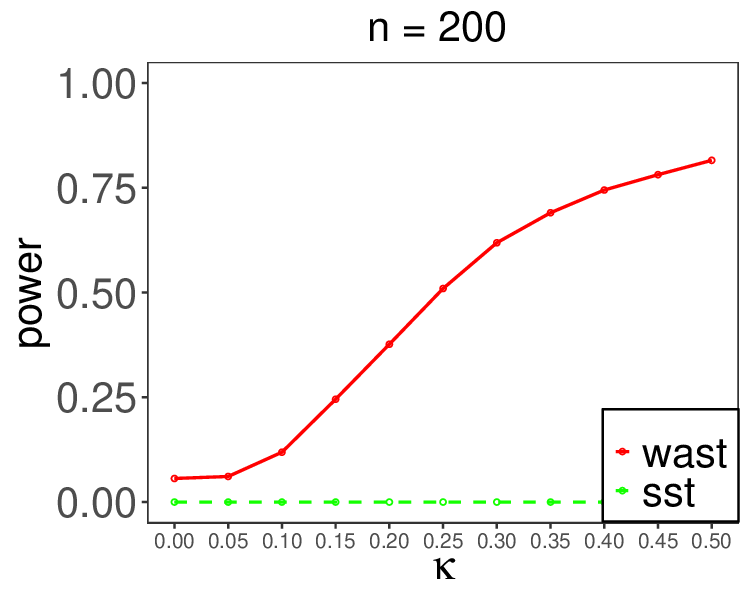}
		\includegraphics[scale=0.3]{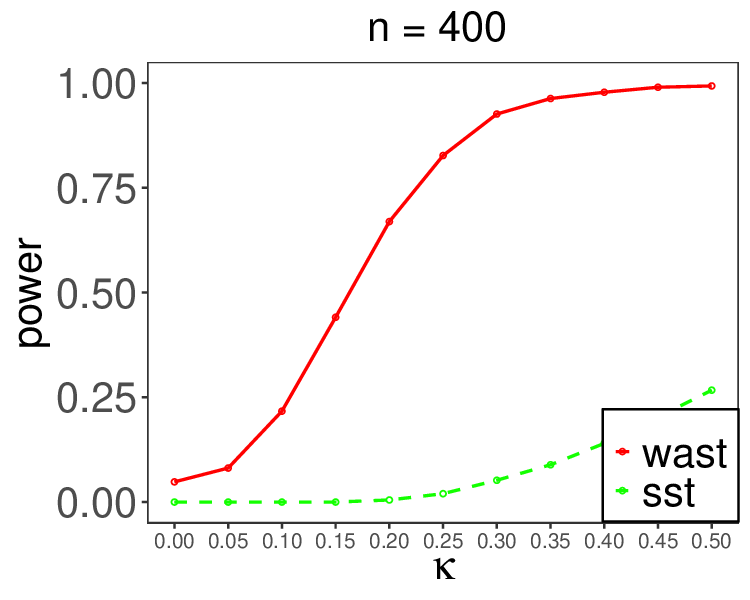}
		\includegraphics[scale=0.3]{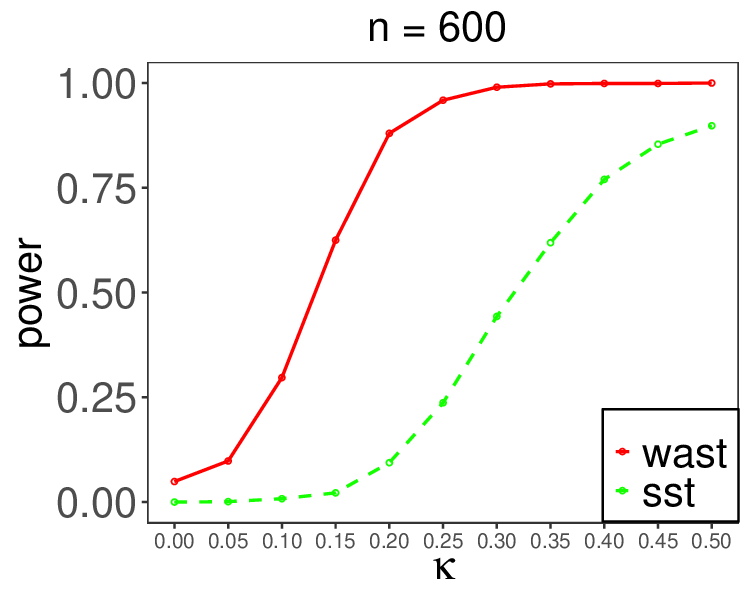} \\
		\includegraphics[scale=0.3]{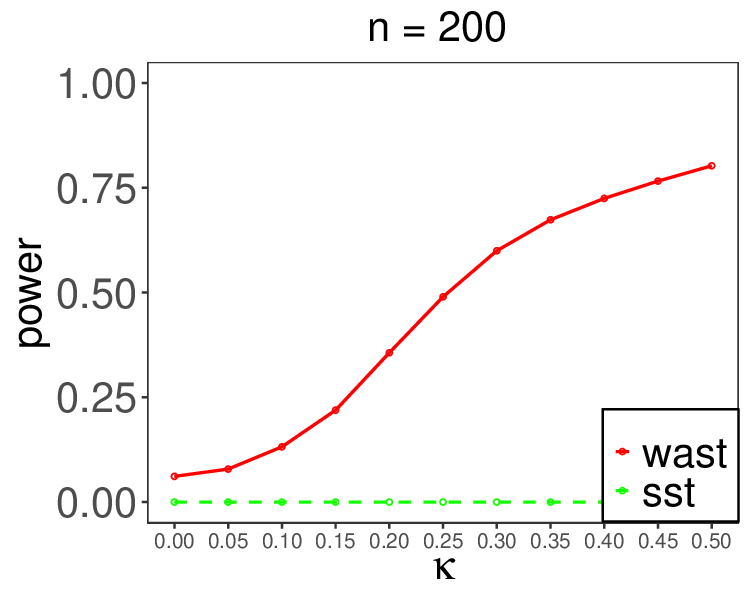}
		\includegraphics[scale=0.3]{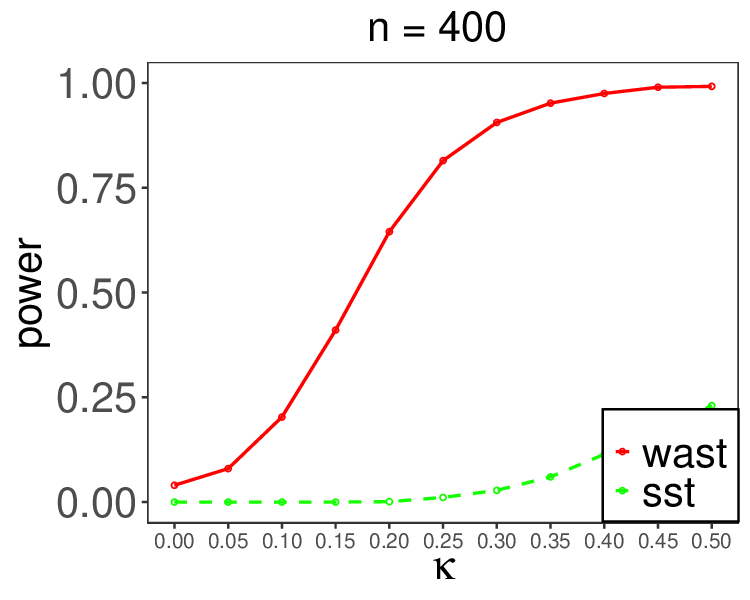}
		\includegraphics[scale=0.3]{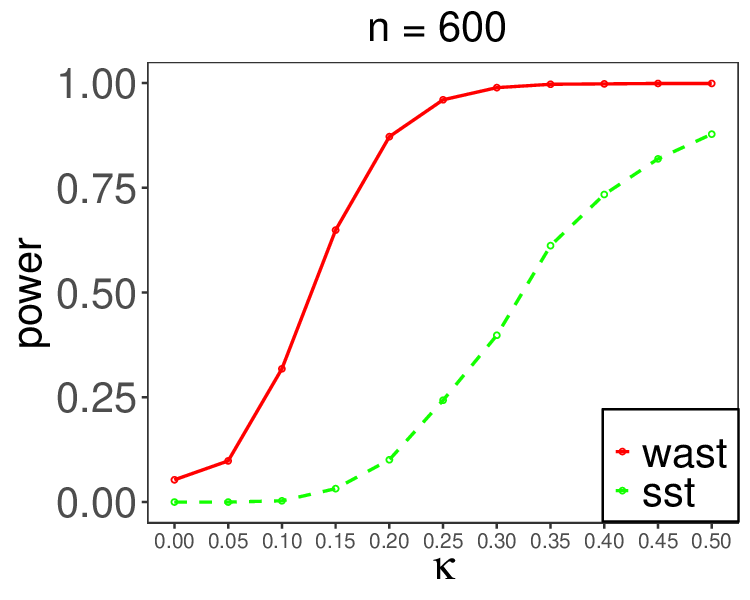}
		\caption{\it Powers of test statistic for probit model with large numbers of sparse $\bZ$ by the proposed WAST (red solid line) and SST (green dashed line). From top to bottom, each row depicts the powers for $(p,q)=(2,100)$, $(p,q)=(2,500)$, $(p,q)=(6,100)$, $(p,q)=(6,500)$, $(p,q)=(11,100)$, and $(p,q)=(11,500)$.}
		\label{fig_probit_sparse}
	\end{center}
\end{figure}

\begin{figure}[!ht]
	\begin{center}
		\includegraphics[scale=0.3]{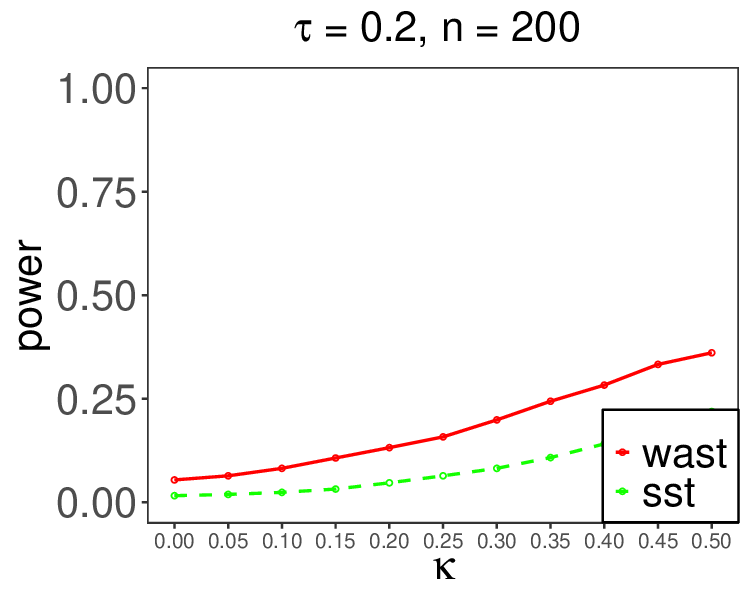}
		\includegraphics[scale=0.3]{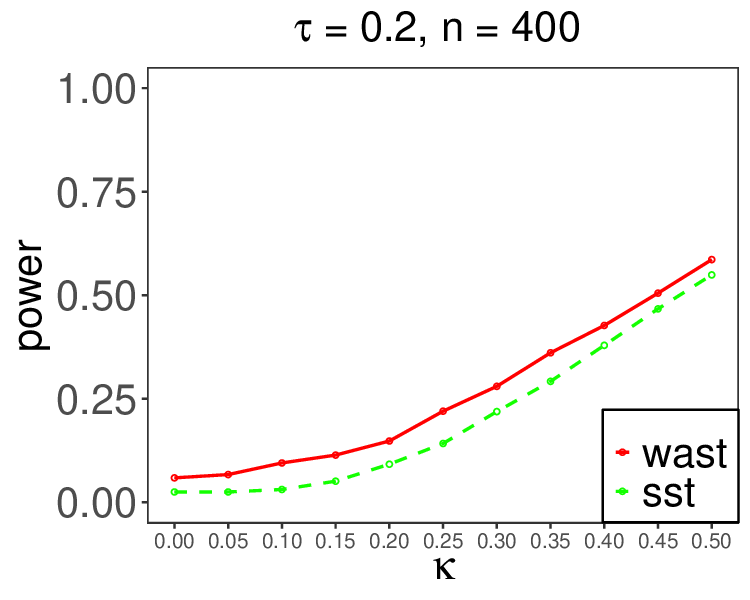}
		\includegraphics[scale=0.3]{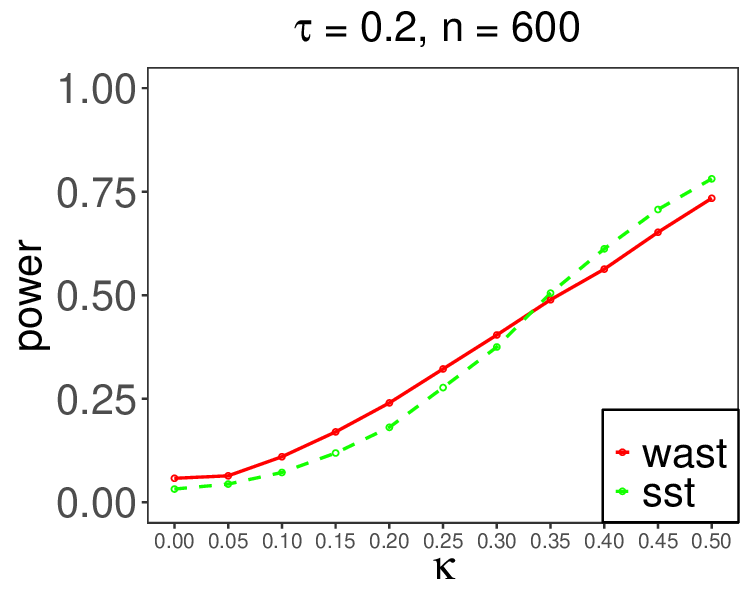}  \\
		\includegraphics[scale=0.3]{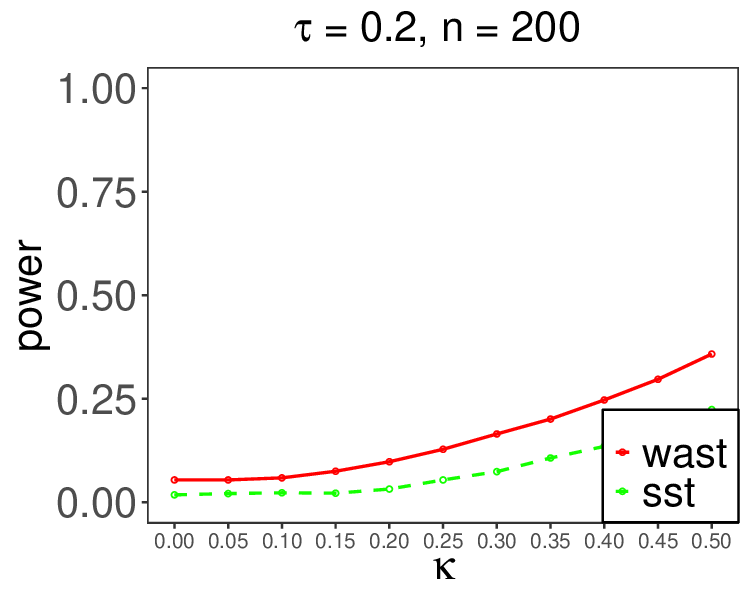}
		\includegraphics[scale=0.3]{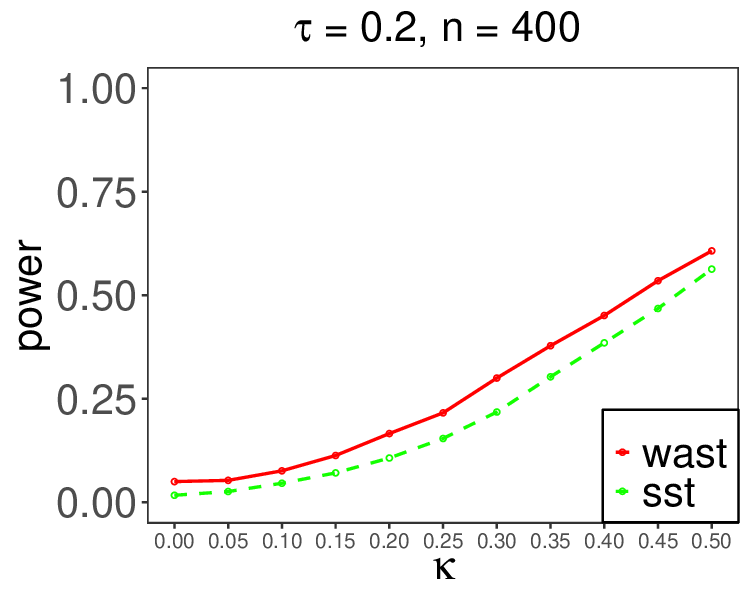}
		\includegraphics[scale=0.3]{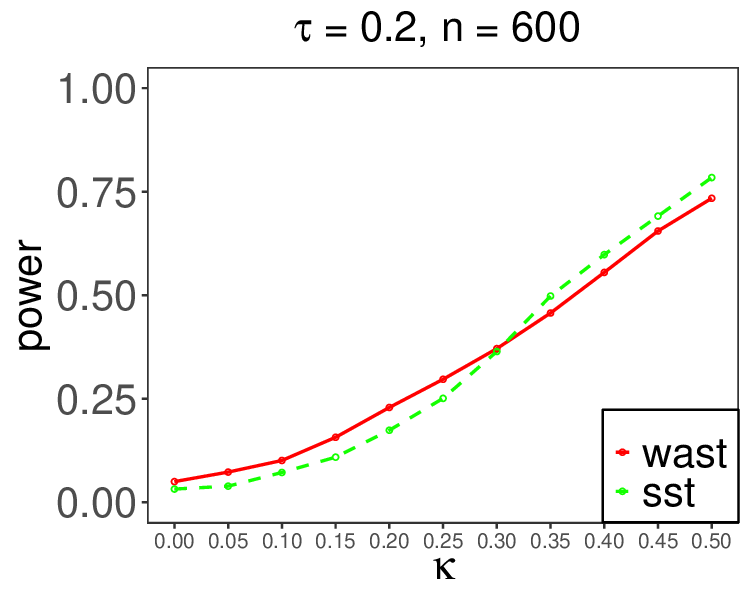}  \\
		\includegraphics[scale=0.3]{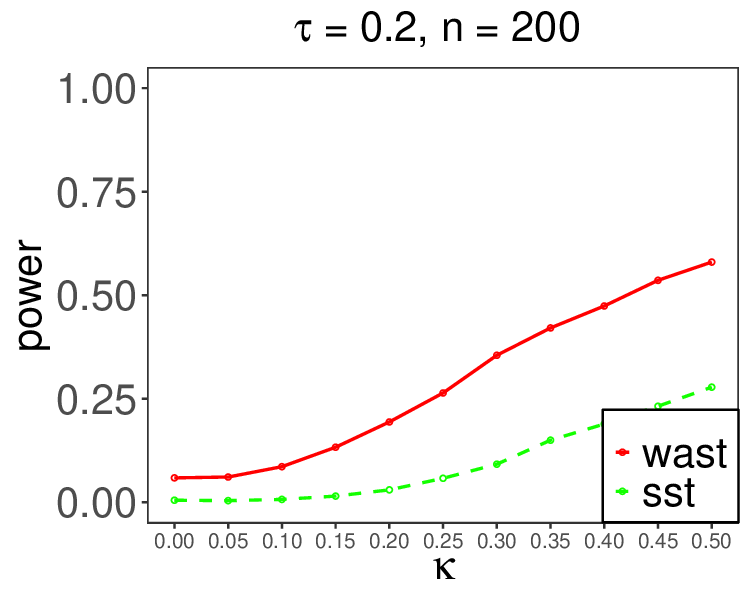}
		\includegraphics[scale=0.3]{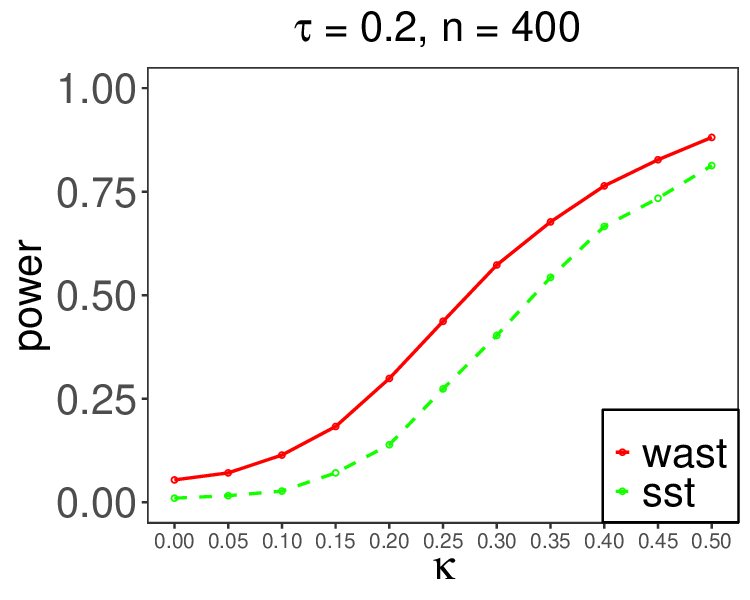}
		\includegraphics[scale=0.3]{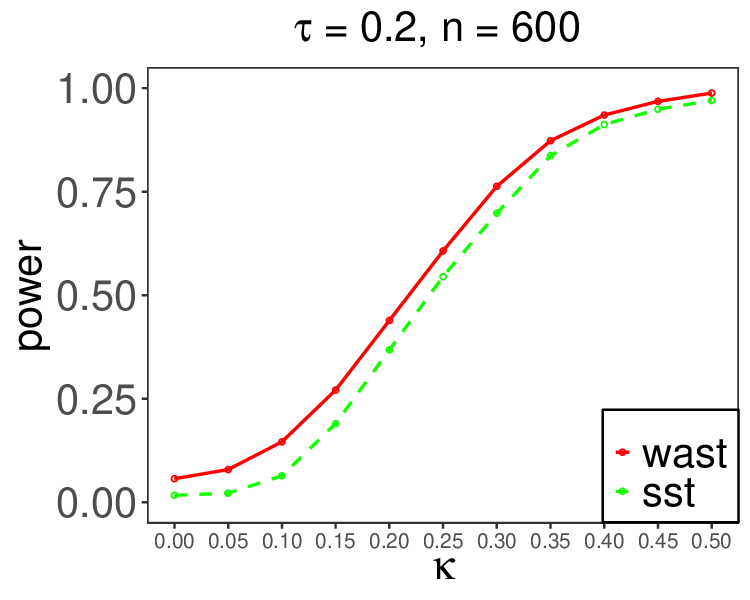}  \\
		\includegraphics[scale=0.3]{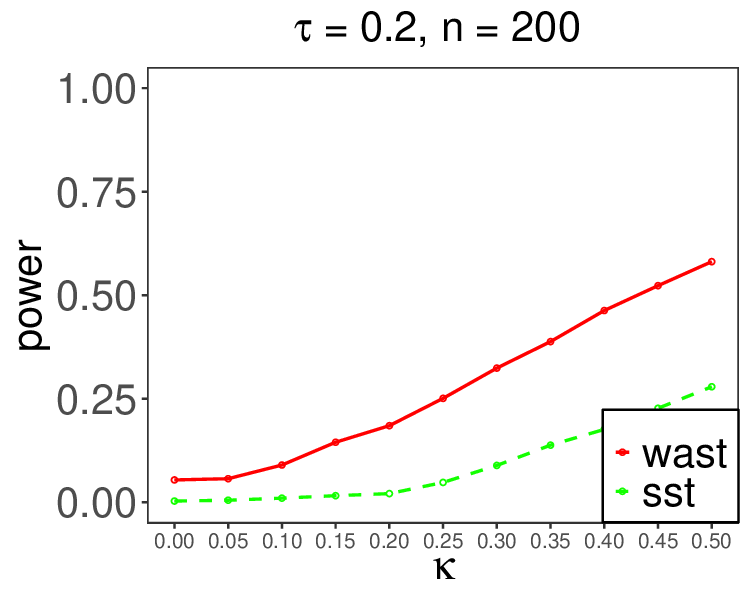}
		\includegraphics[scale=0.3]{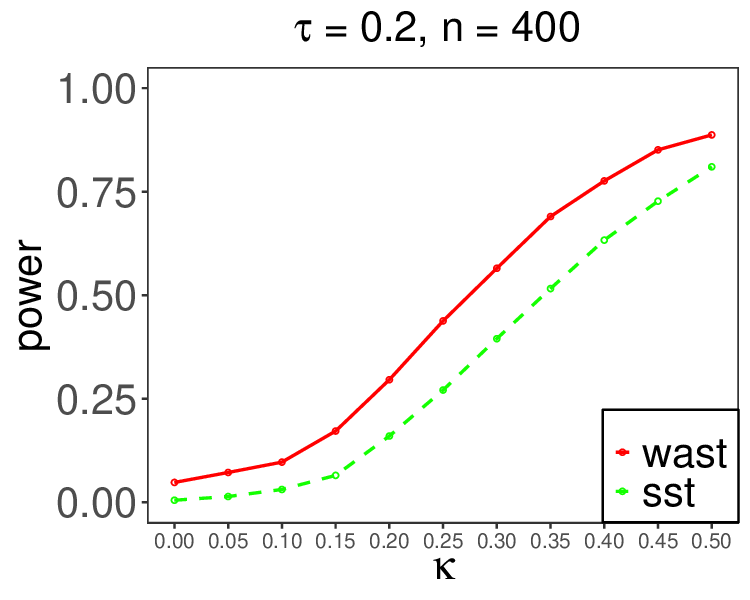}
		\includegraphics[scale=0.3]{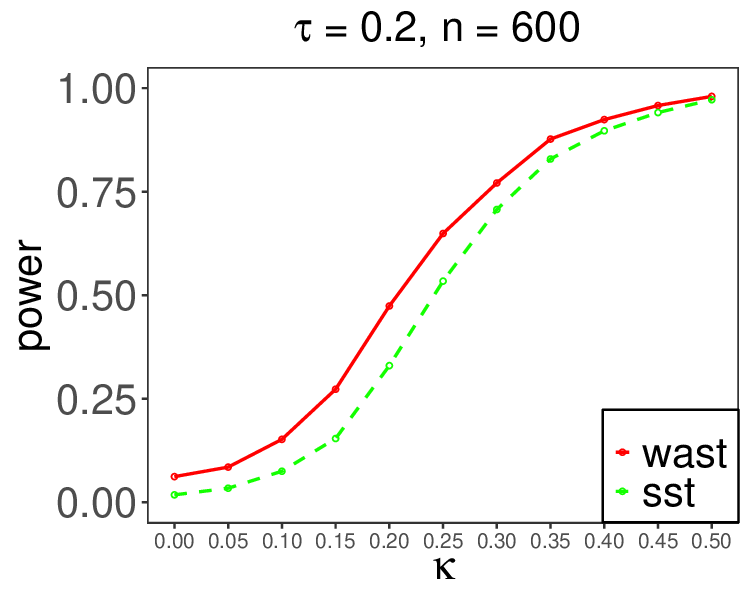}  \\
		\includegraphics[scale=0.3]{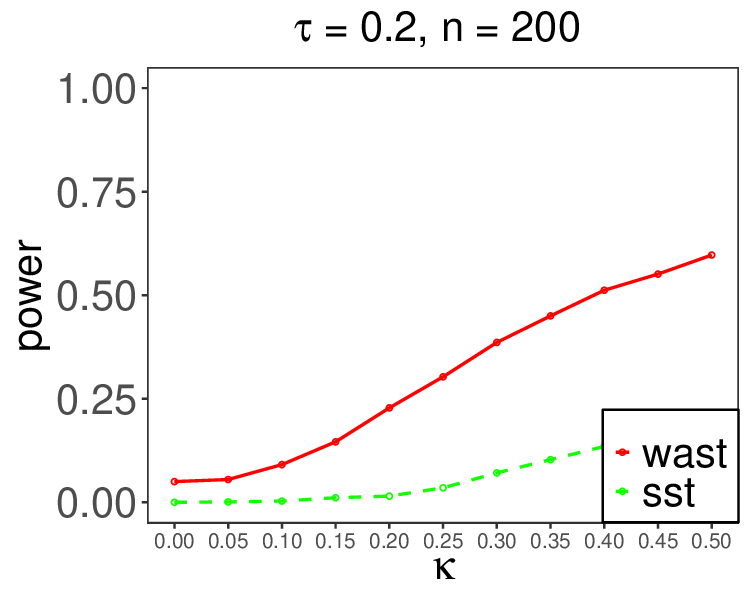}
		\includegraphics[scale=0.3]{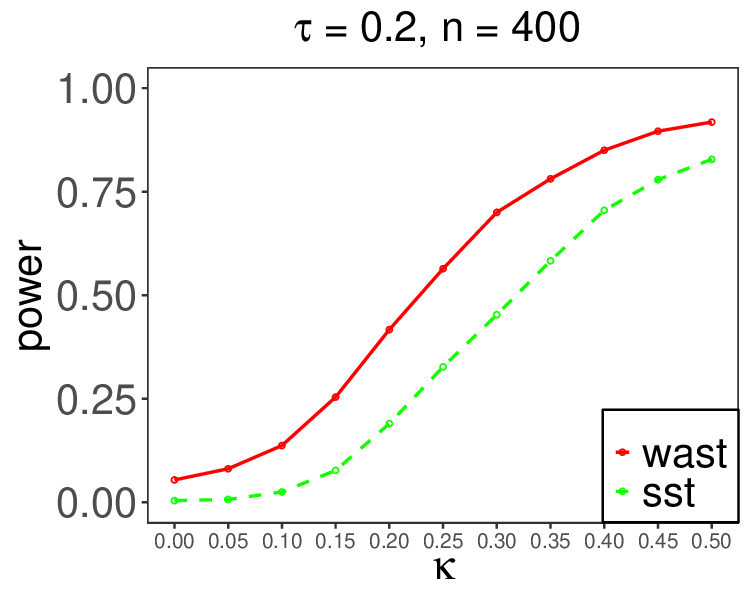}
		\includegraphics[scale=0.3]{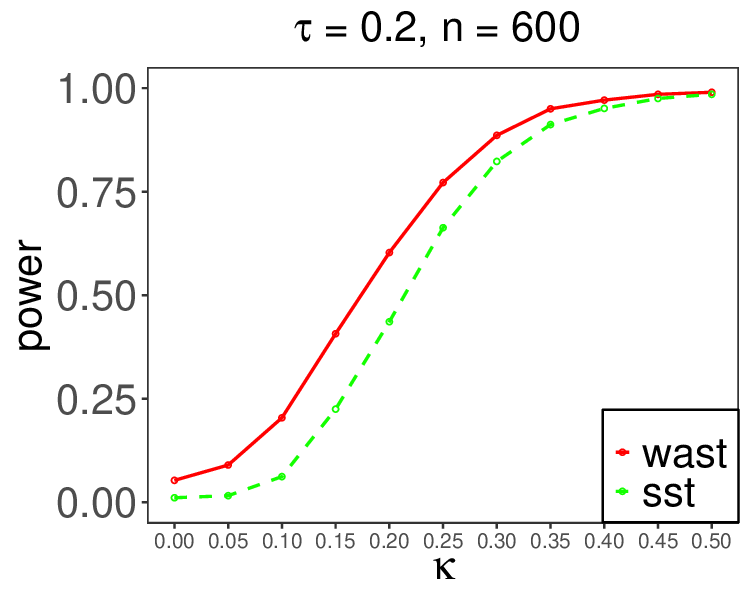} \\
		\includegraphics[scale=0.3]{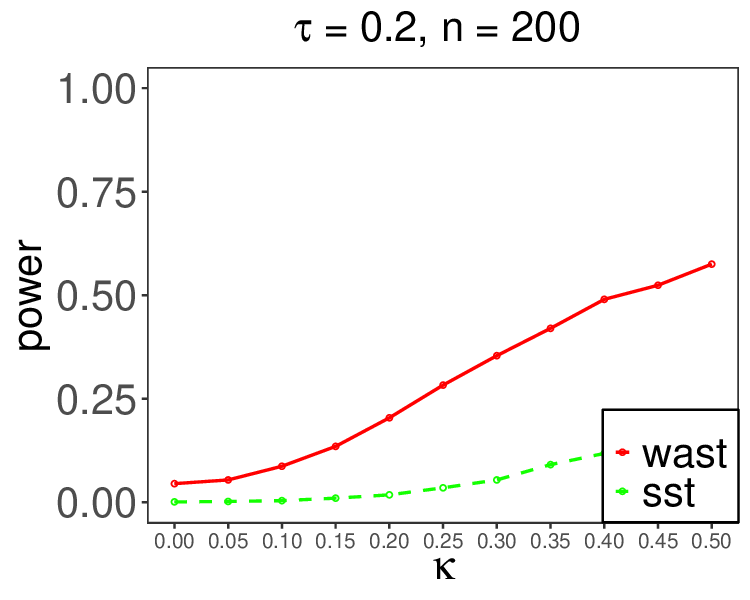}
		\includegraphics[scale=0.3]{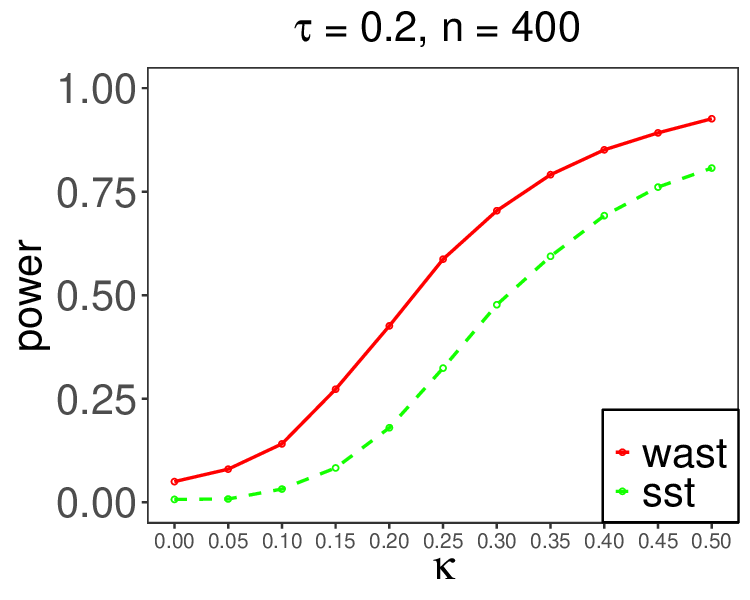}
		\includegraphics[scale=0.3]{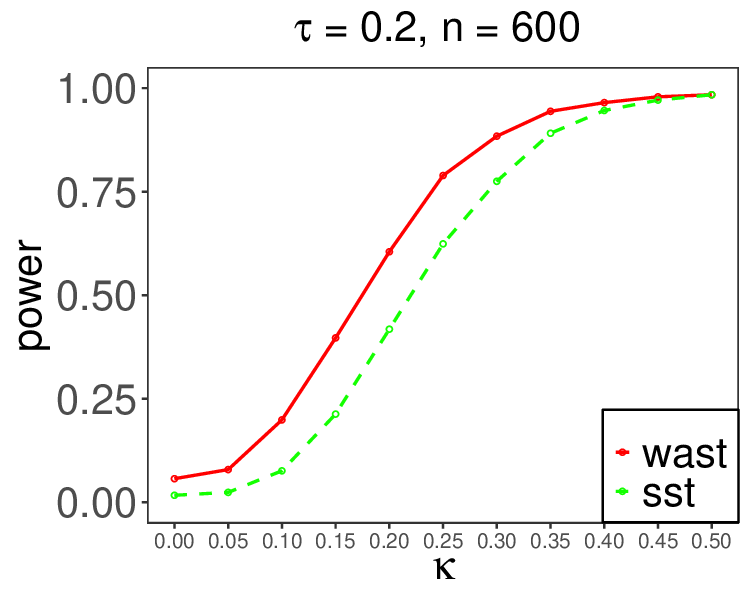}
		\caption{\it Powers of test statistic for quantile regression with $\tau=0.2$ and with large numbers of sparse $\bZ$ by the proposed WAST (red solid line) and SST (green dashed line). From top to bottom, each row depicts the powers for $(p,q)=(2,100)$, $(p,q)=(2,500)$, $(p,q)=(6,100)$, $(p,q)=(6,500)$, $(p,q)=(11,100)$, and $(p,q)=(11,500)$.}
		\label{fig_qr20_sparse}
	\end{center}
\end{figure}

\begin{figure}[!ht]
	\begin{center}
		\includegraphics[scale=0.3]{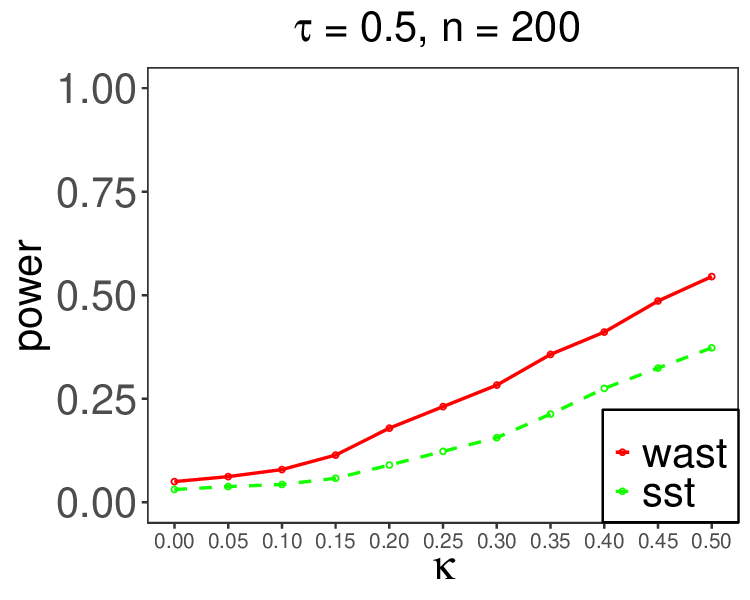}
		\includegraphics[scale=0.3]{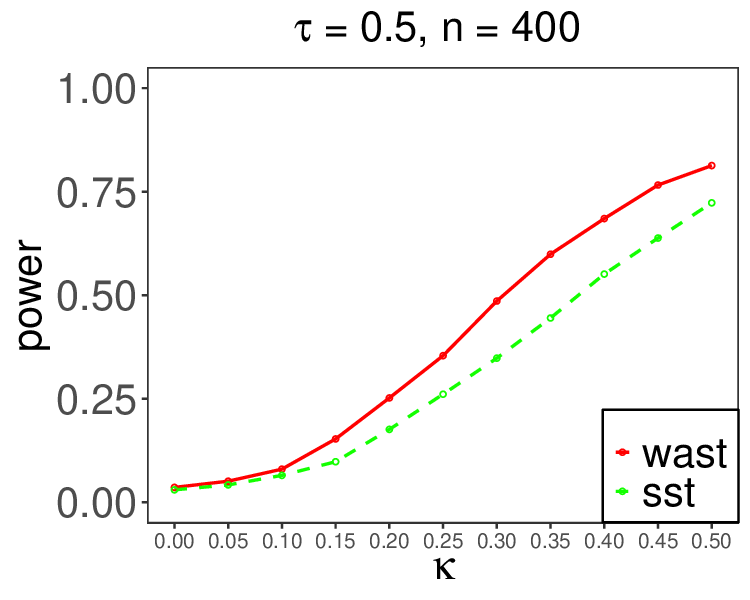}
		\includegraphics[scale=0.3]{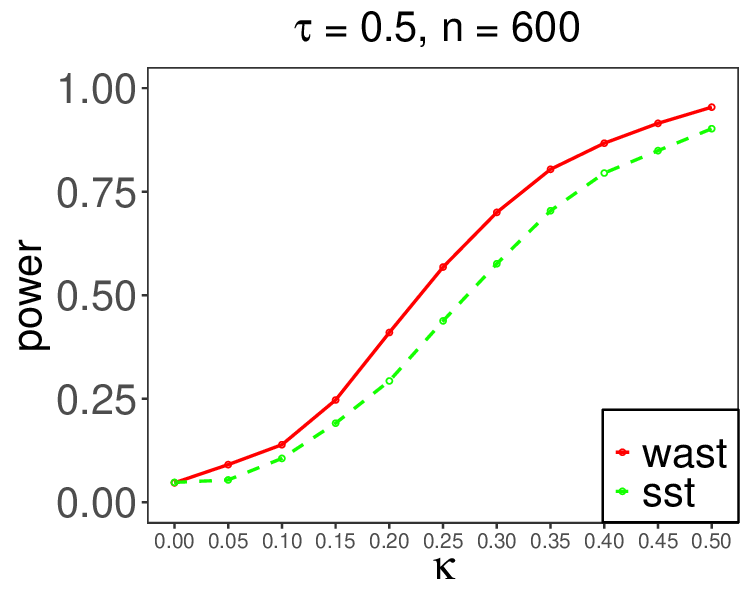}     \\
		\includegraphics[scale=0.3]{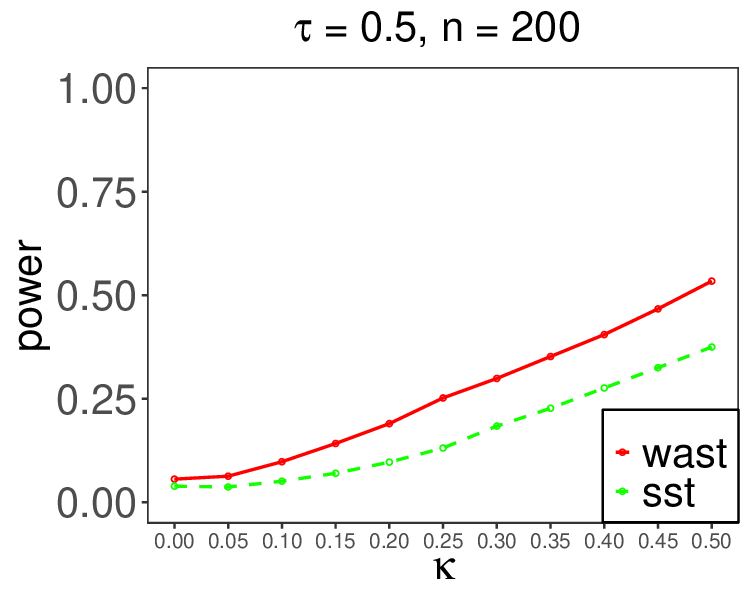}
		\includegraphics[scale=0.3]{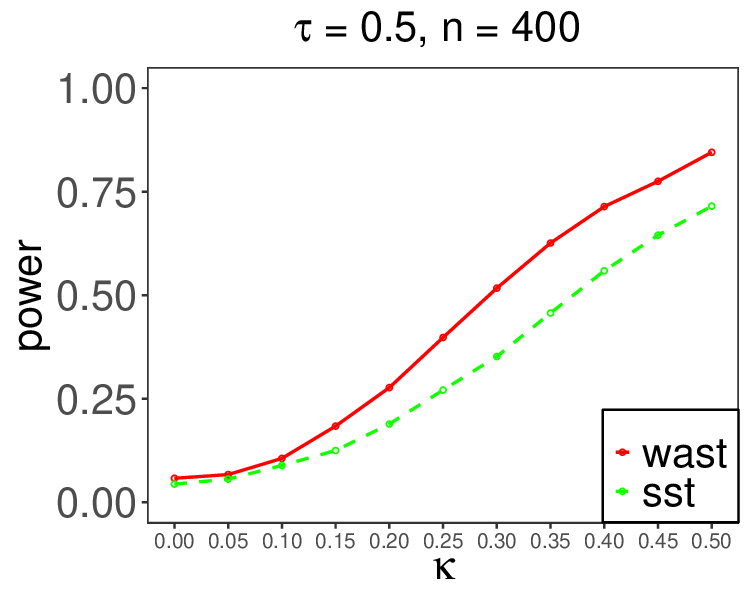}
		\includegraphics[scale=0.3]{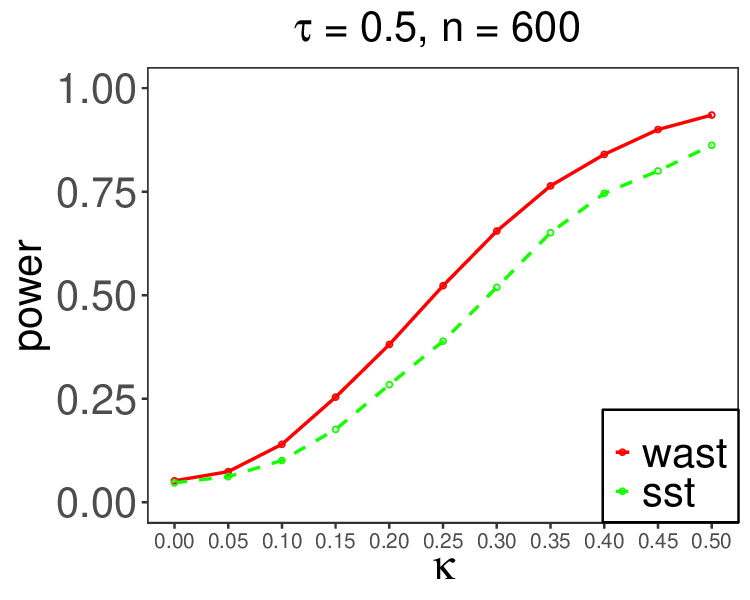}     \\
		\includegraphics[scale=0.3]{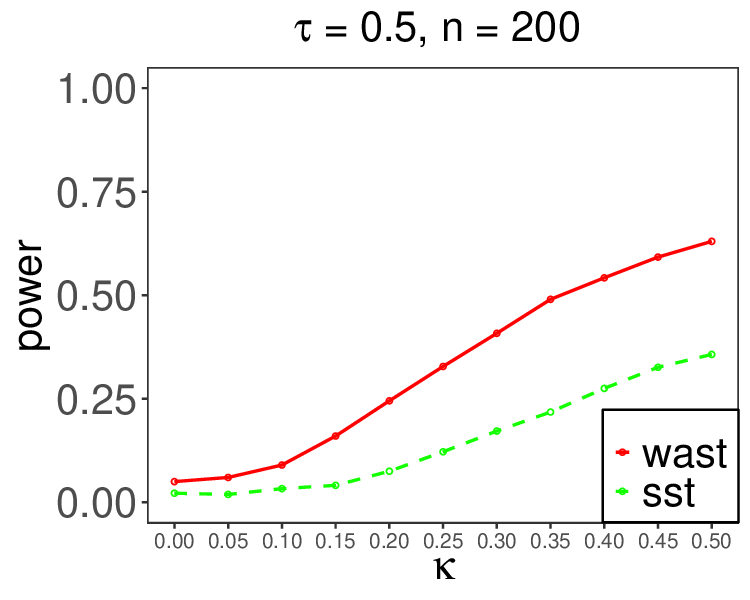}
		\includegraphics[scale=0.3]{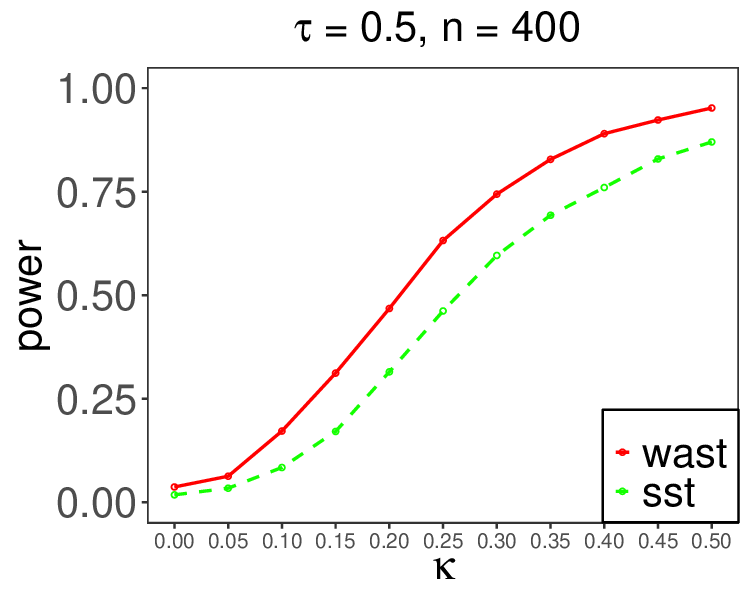}
		\includegraphics[scale=0.3]{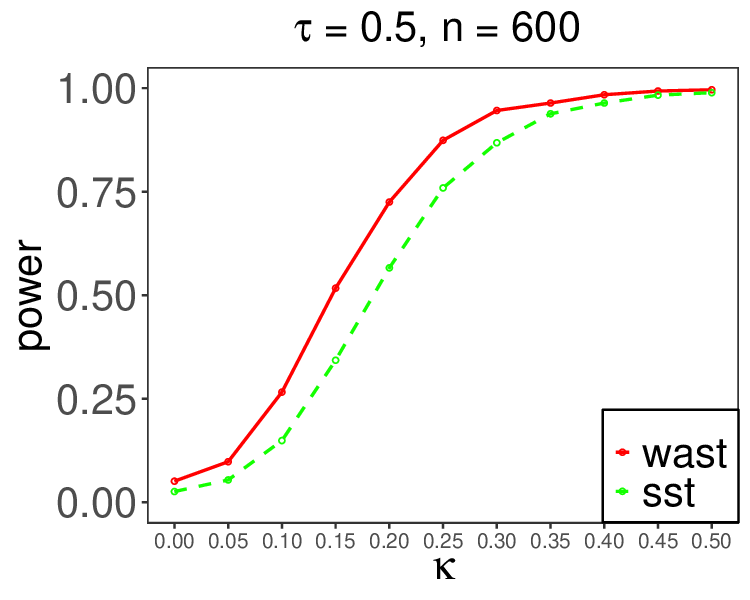}     \\
		\includegraphics[scale=0.3]{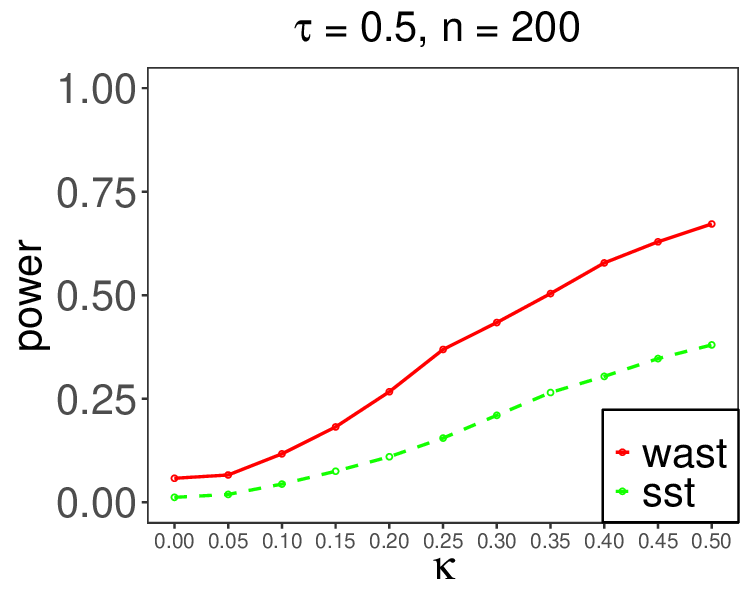}
		\includegraphics[scale=0.3]{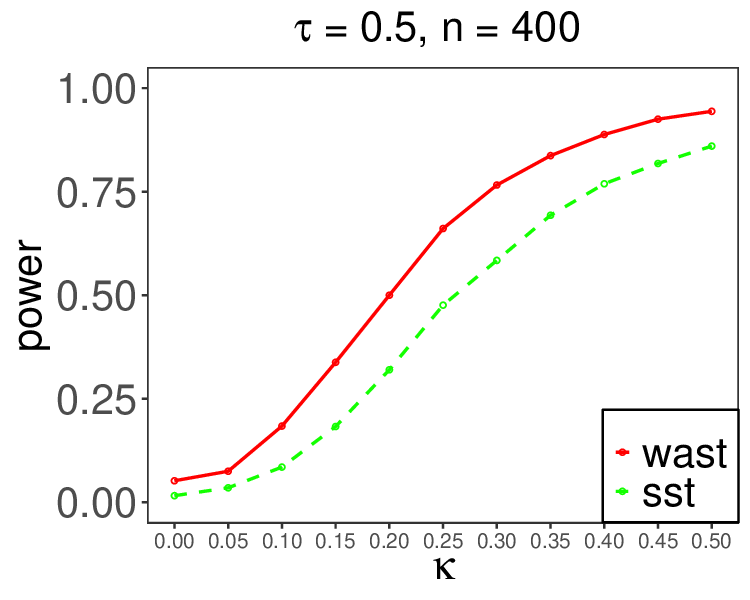}
		\includegraphics[scale=0.3]{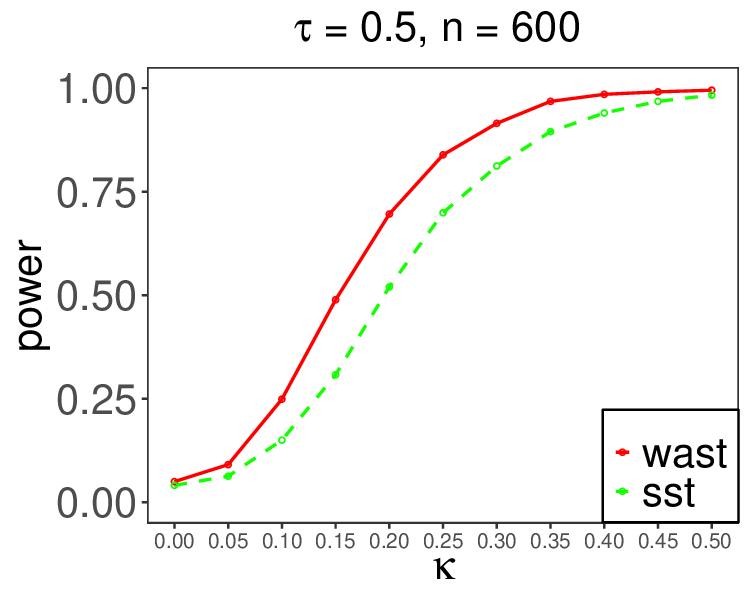}     \\
		\includegraphics[scale=0.3]{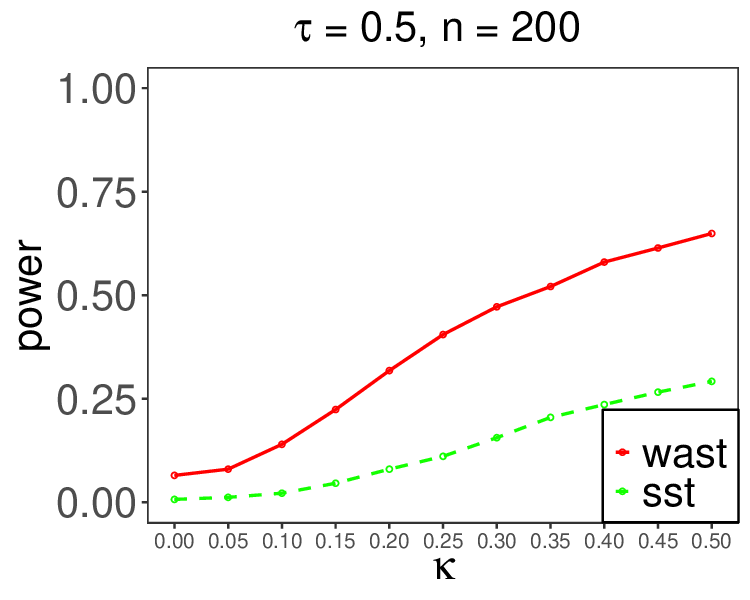}
		\includegraphics[scale=0.3]{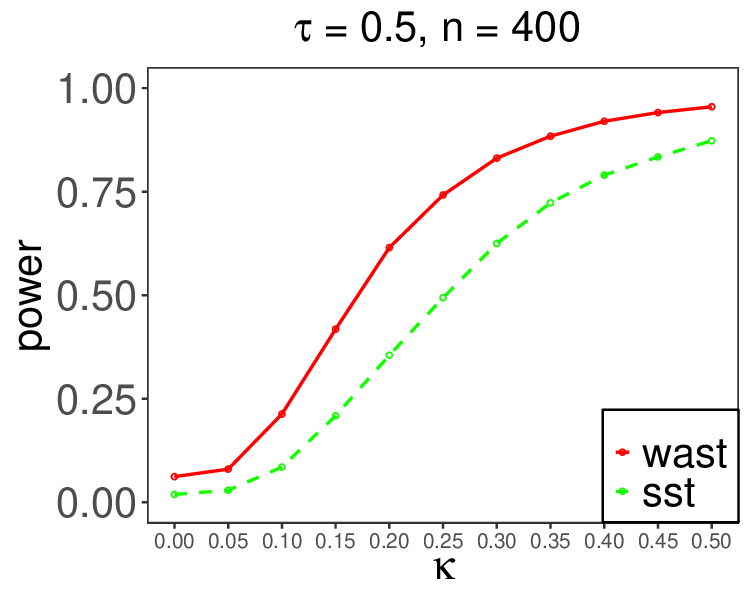}
		\includegraphics[scale=0.3]{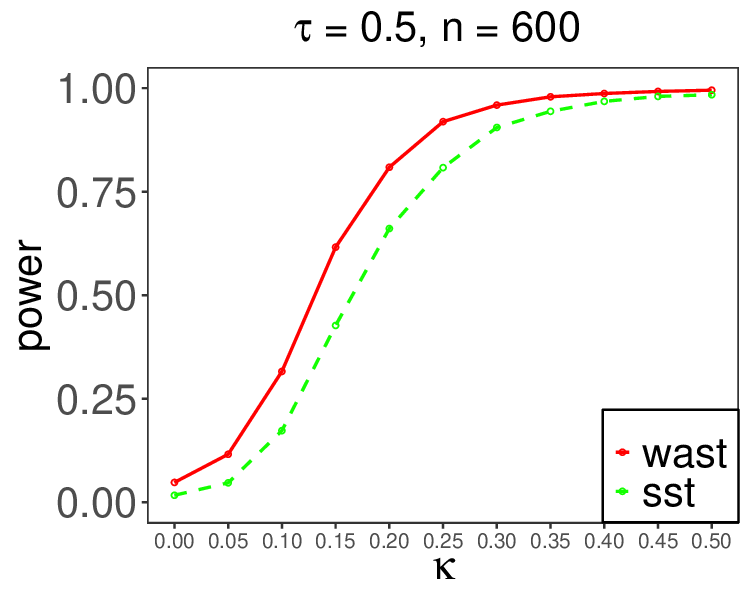}    \\
		\includegraphics[scale=0.3]{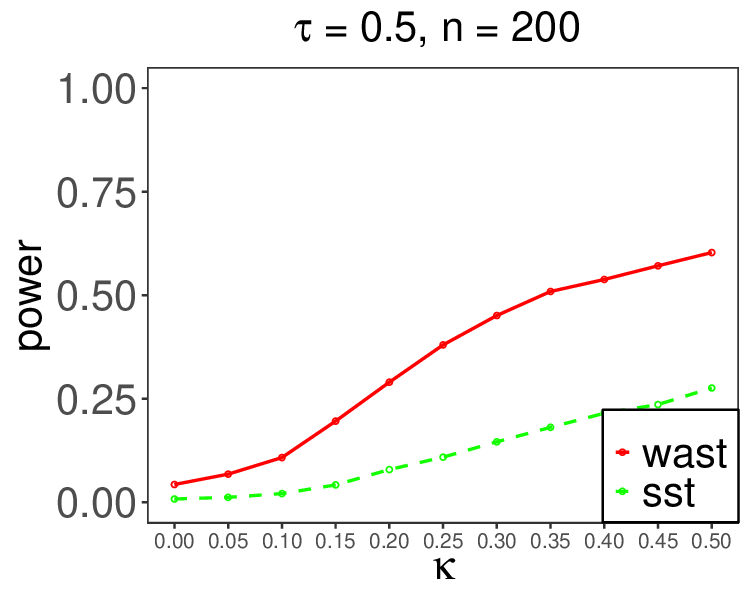}
		\includegraphics[scale=0.3]{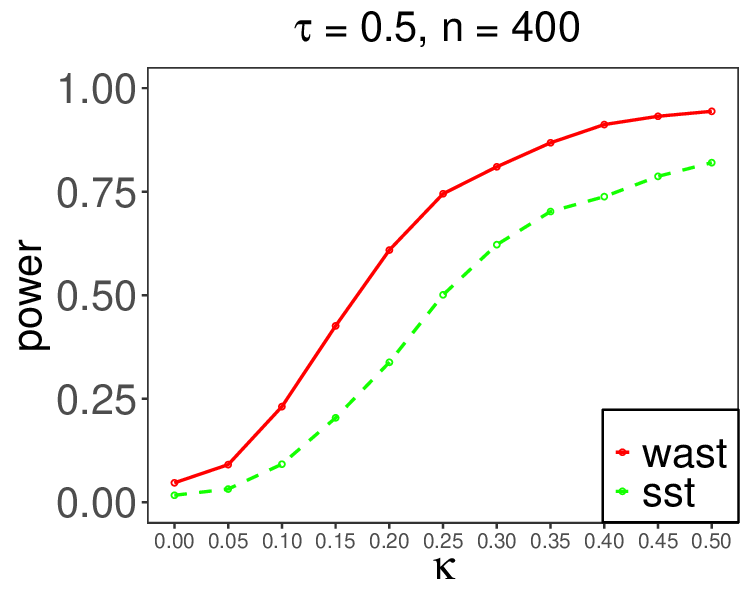}
		\includegraphics[scale=0.3]{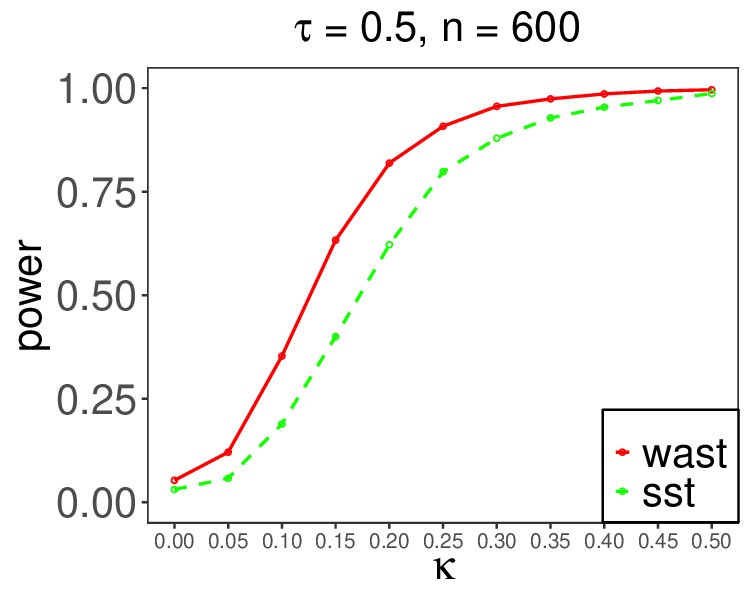}
		\caption{\it Powers of test statistic for quantile regression with $\tau=0.5$ and with large numbers of sparse $\bZ$ by the proposed WAST (red solid line) and SST (green dashed line). From top to bottom, each row depicts the powers for $(p,q)=(2,100)$, $(p,q)=(2,500)$, $(p,q)=(6,100)$, $(p,q)=(6,500)$, $(p,q)=(11,100)$, and $(p,q)=(11,500)$.}
		\label{fig_qr50_sparse}
	\end{center}
\end{figure}

\begin{figure}[!ht]
	\begin{center}
		\includegraphics[scale=0.3]{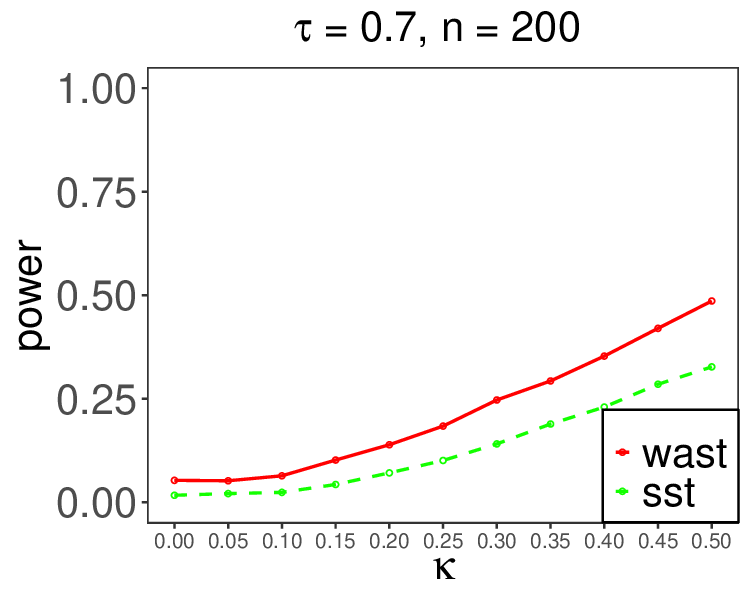}
		\includegraphics[scale=0.3]{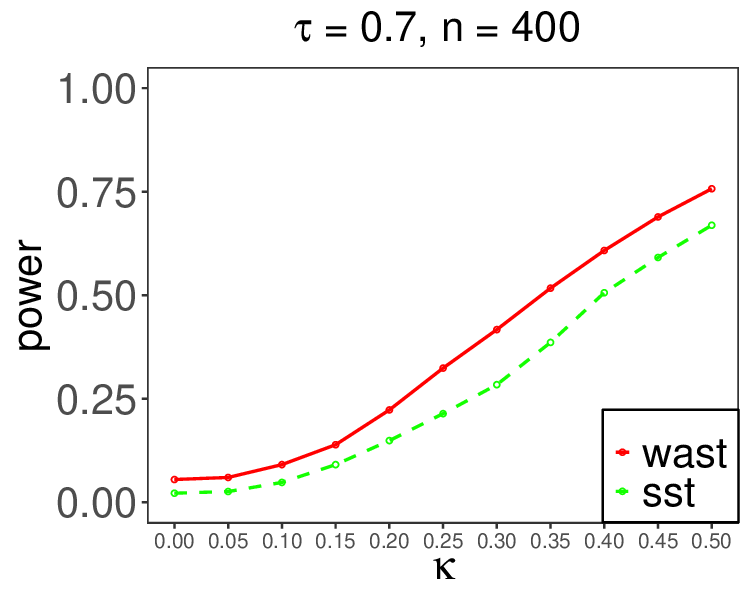}
		\includegraphics[scale=0.3]{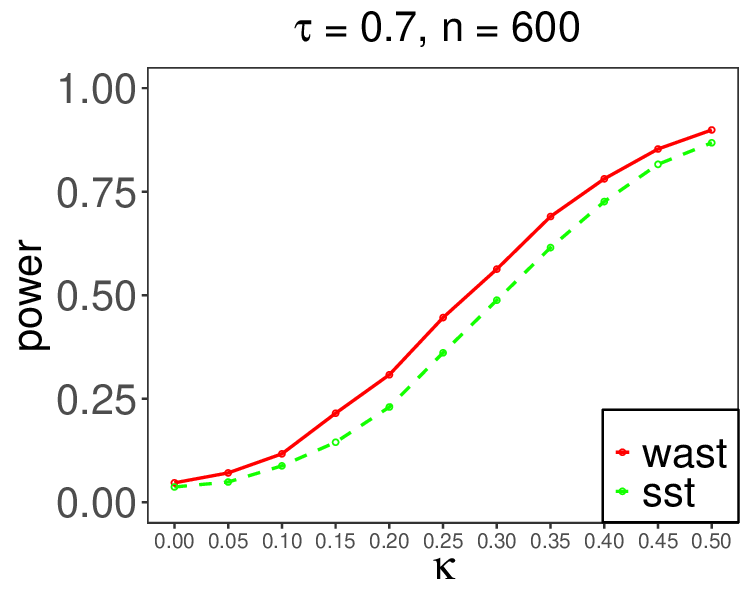}     \\
		\includegraphics[scale=0.3]{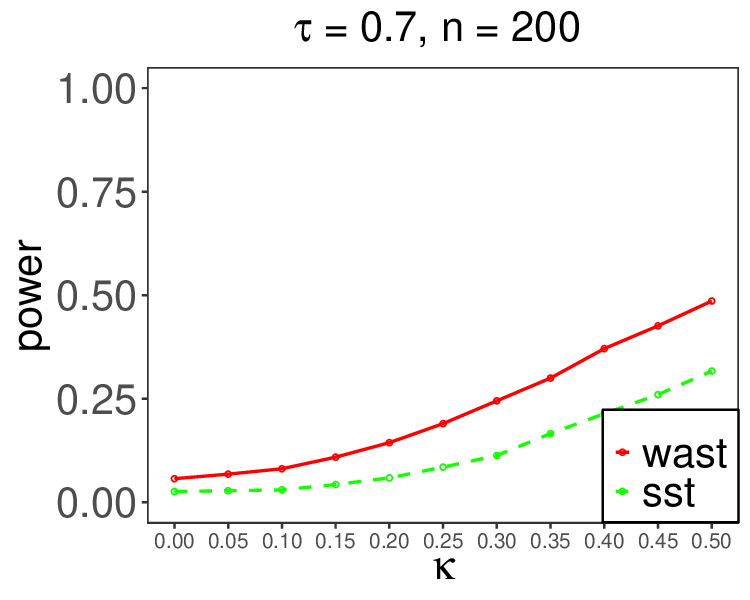}
		\includegraphics[scale=0.3]{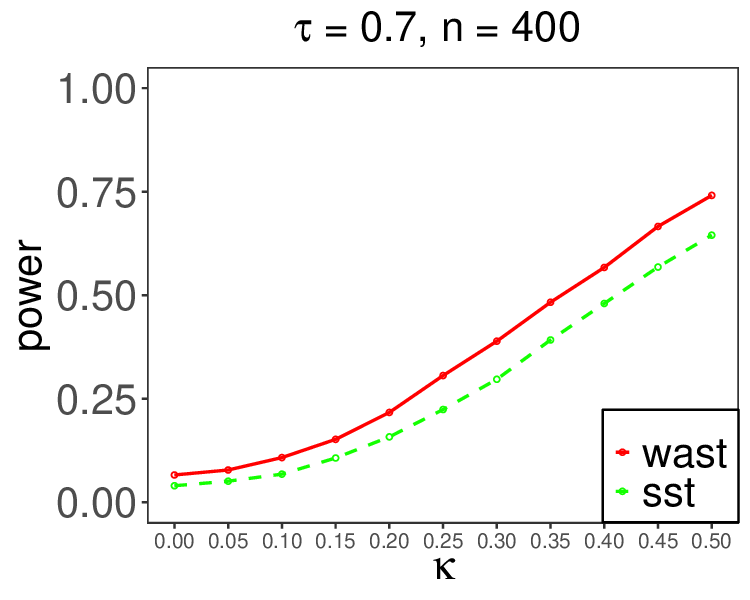}
		\includegraphics[scale=0.3]{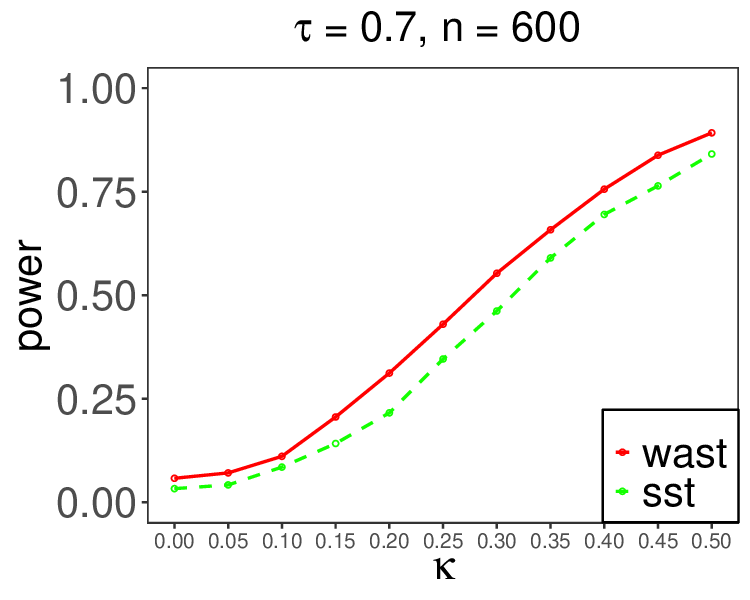}     \\
		\includegraphics[scale=0.3]{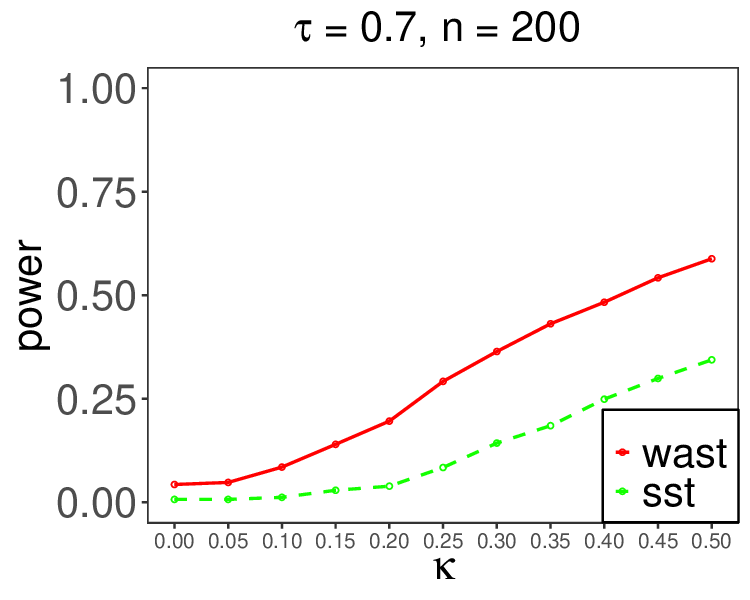}
		\includegraphics[scale=0.3]{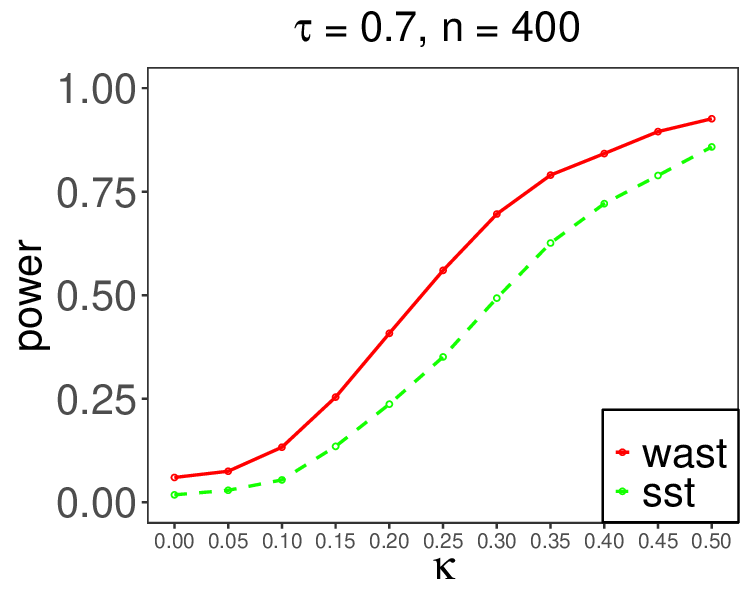}
		\includegraphics[scale=0.3]{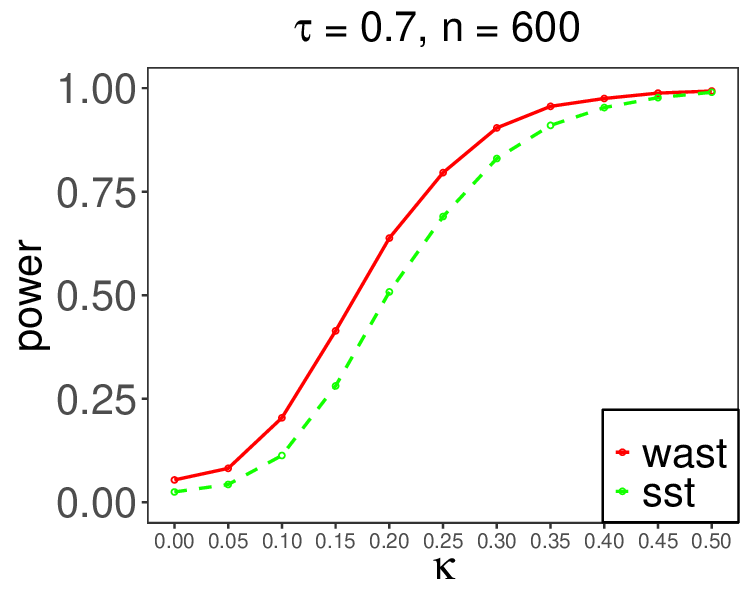}     \\
		\includegraphics[scale=0.3]{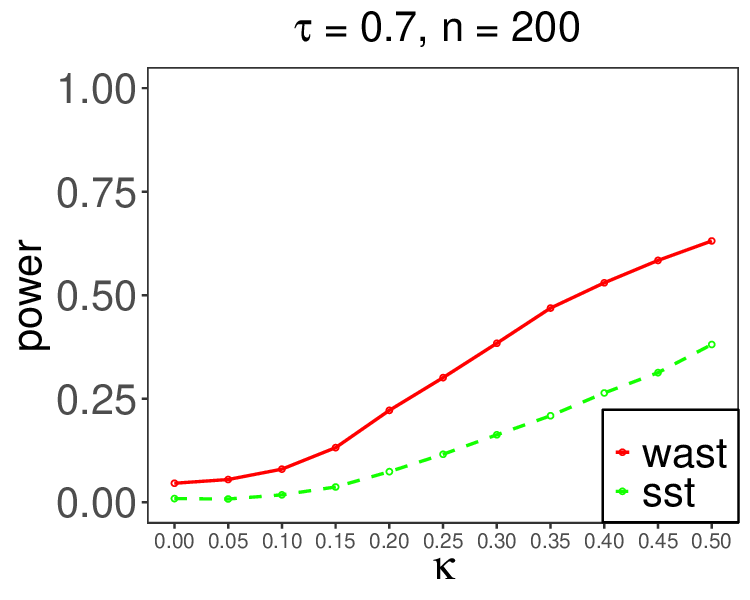}
		\includegraphics[scale=0.3]{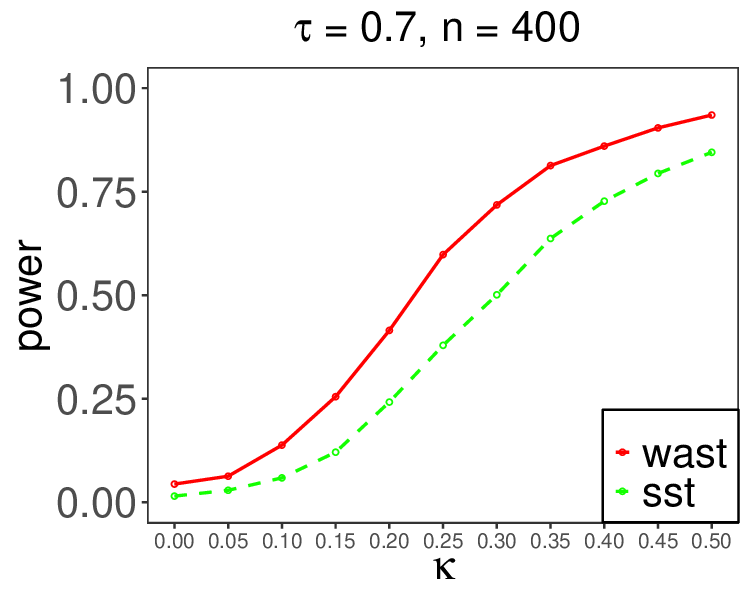}
		\includegraphics[scale=0.3]{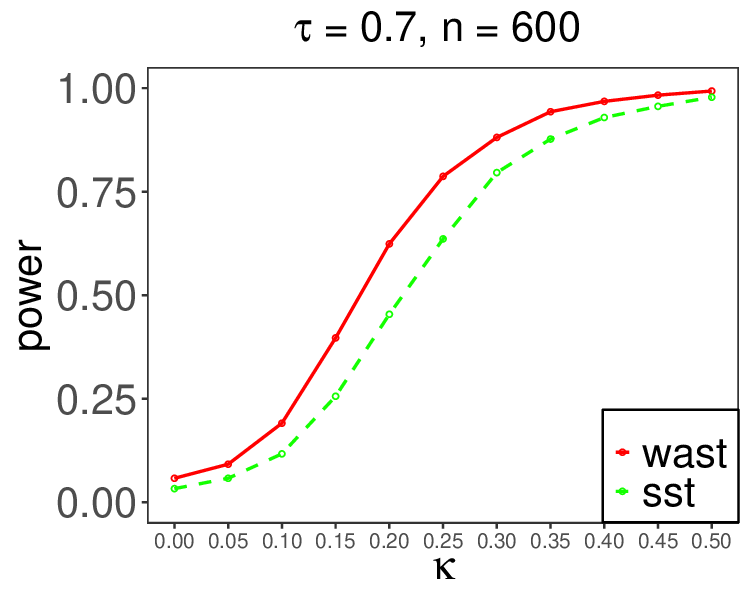}     \\
		\includegraphics[scale=0.3]{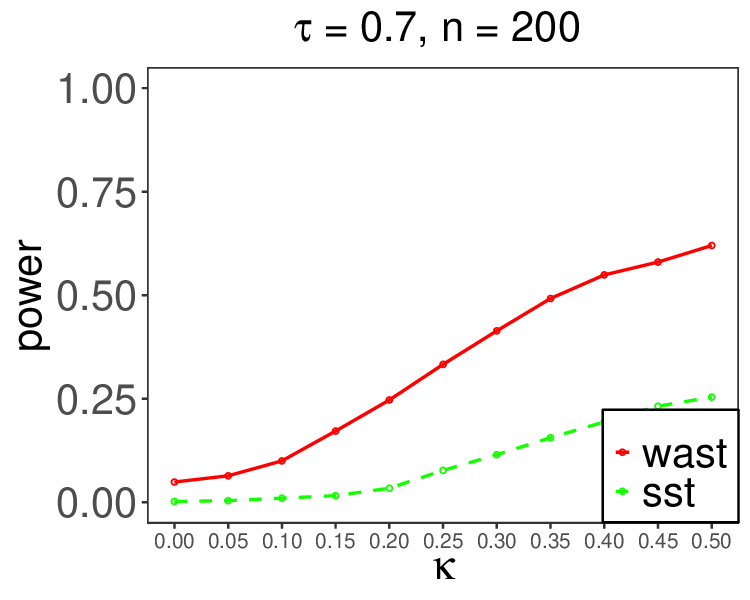}
		\includegraphics[scale=0.3]{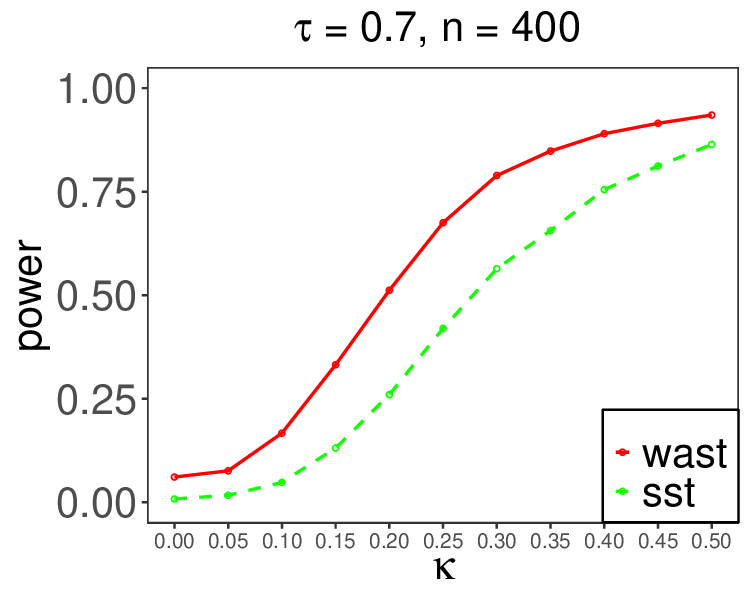}
		\includegraphics[scale=0.3]{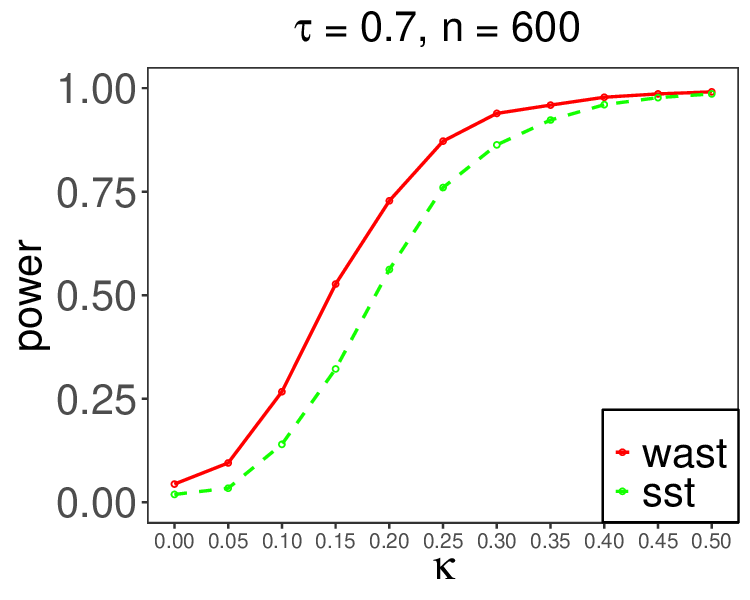}    \\
		\includegraphics[scale=0.3]{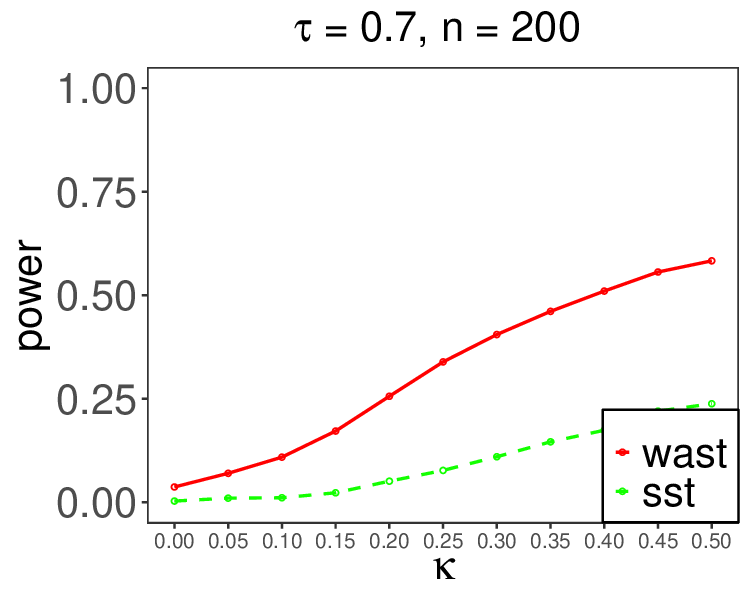}
		\includegraphics[scale=0.3]{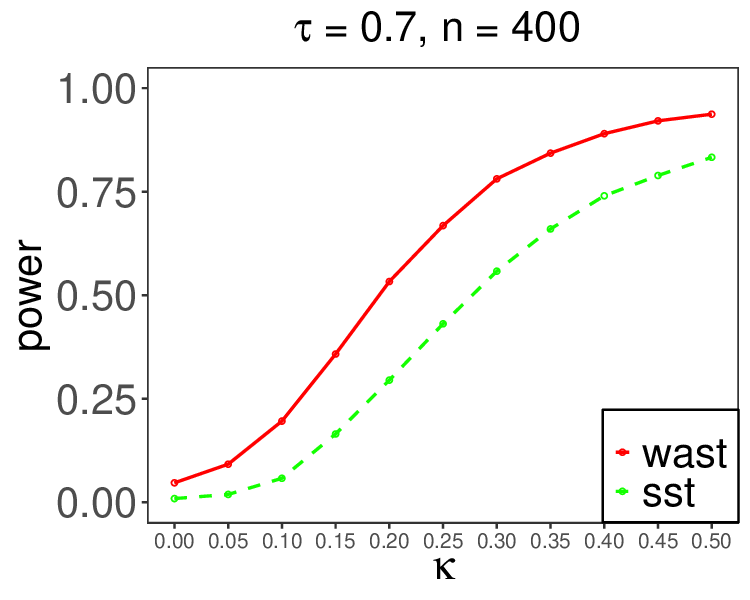}
		\includegraphics[scale=0.3]{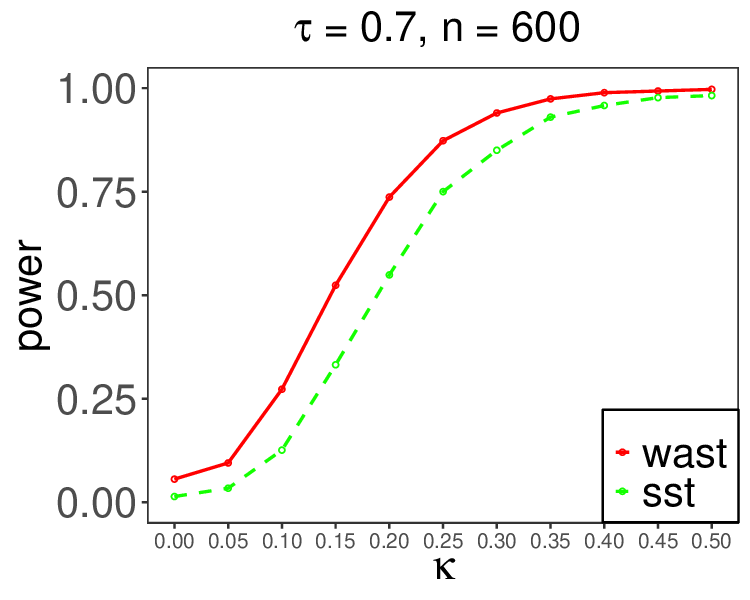}
		\caption{\it Powers of test statistic for quantile regression with $\tau=0.7$ and with large numbers of sparse $\bZ$ by the proposed WAST (red solid line) and SST (green dashed line). From top to bottom, each row depicts the powers for $(p,q)=(2,100)$, $(p,q)=(2,500)$, $(p,q)=(6,100)$, $(p,q)=(6,500)$, $(p,q)=(11,100)$, and $(p,q)=(11,500)$.}
		\label{fig_qr70_sparse}
	\end{center}
\end{figure}

\begin{figure}[!ht]
	\begin{center}
		\includegraphics[scale=0.3]{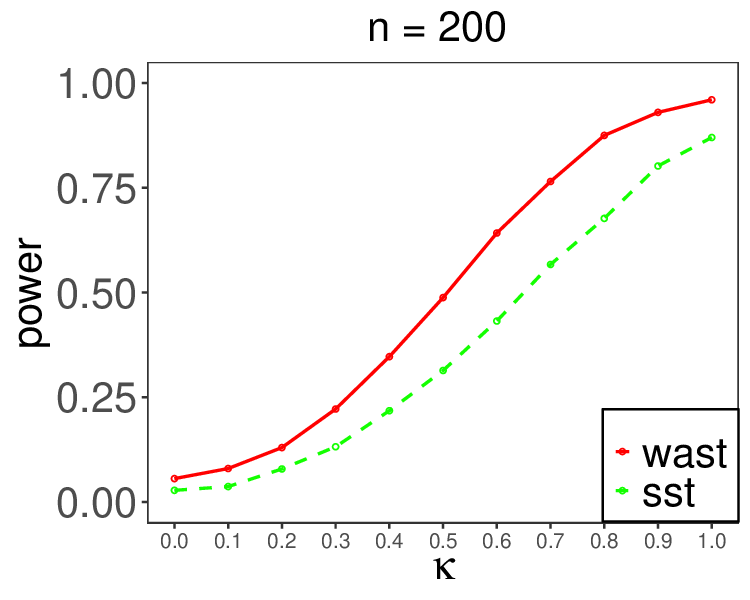}
		\includegraphics[scale=0.3]{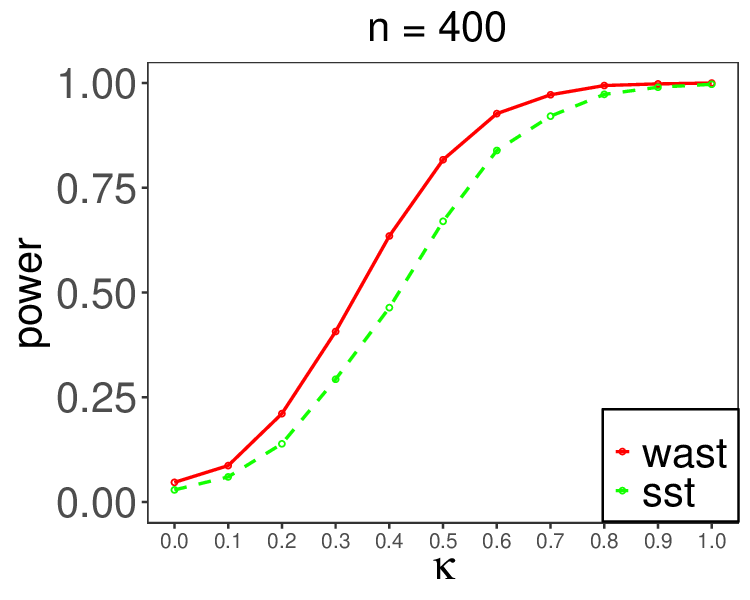}
		\includegraphics[scale=0.3]{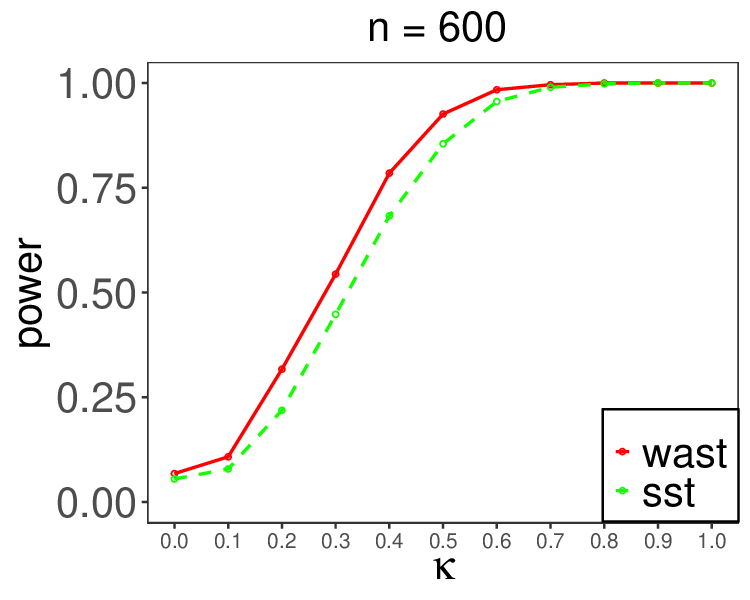}      \\
		\includegraphics[scale=0.3]{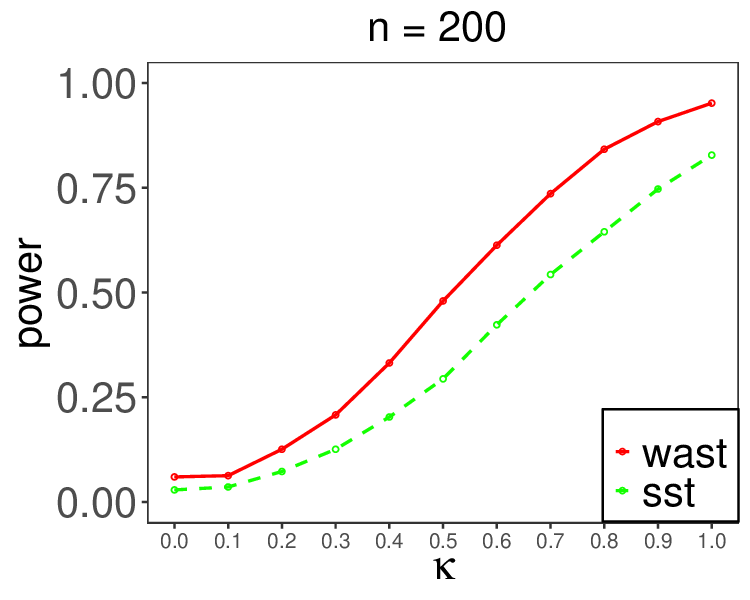}
		\includegraphics[scale=0.3]{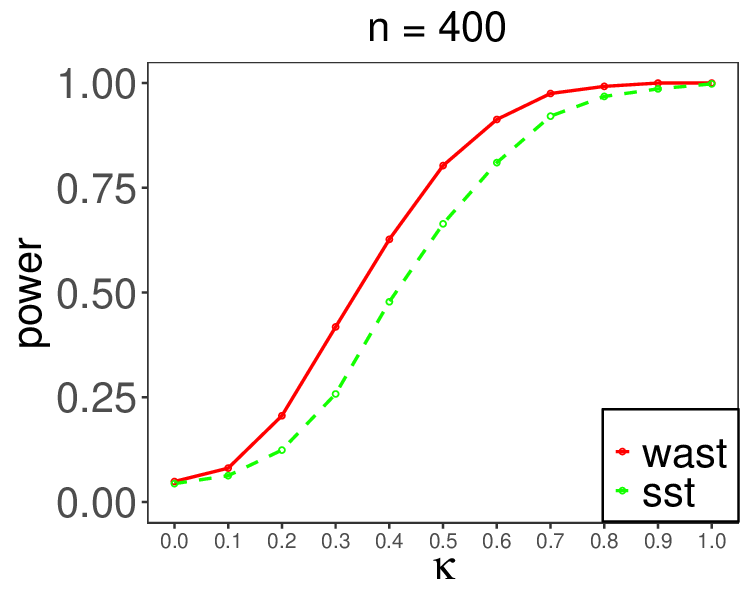}
		\includegraphics[scale=0.3]{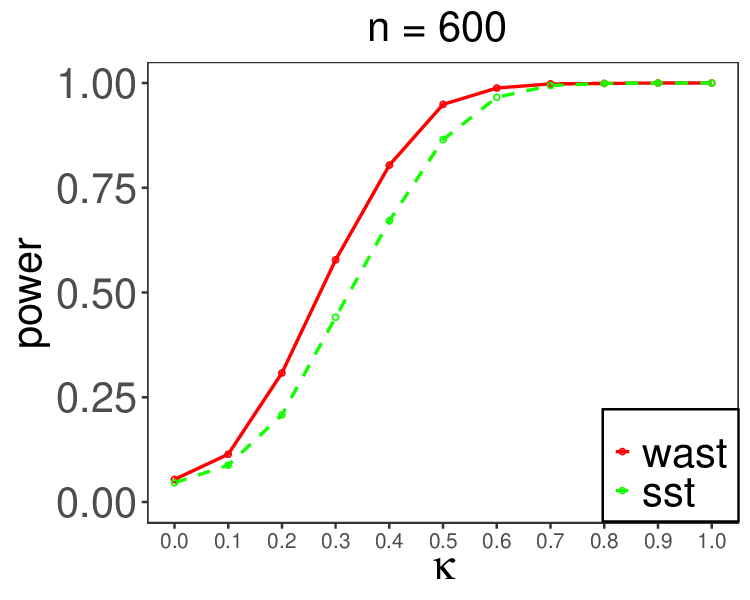}      \\
		\includegraphics[scale=0.3]{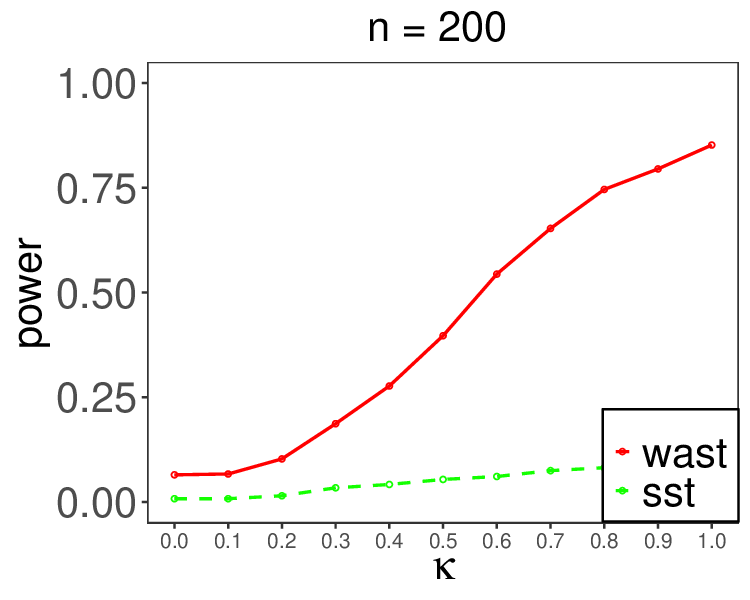}
		\includegraphics[scale=0.3]{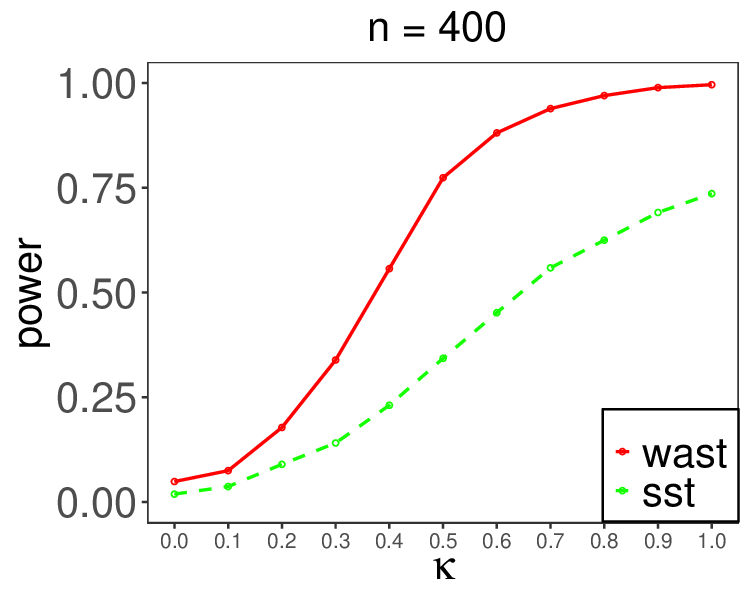}
		\includegraphics[scale=0.3]{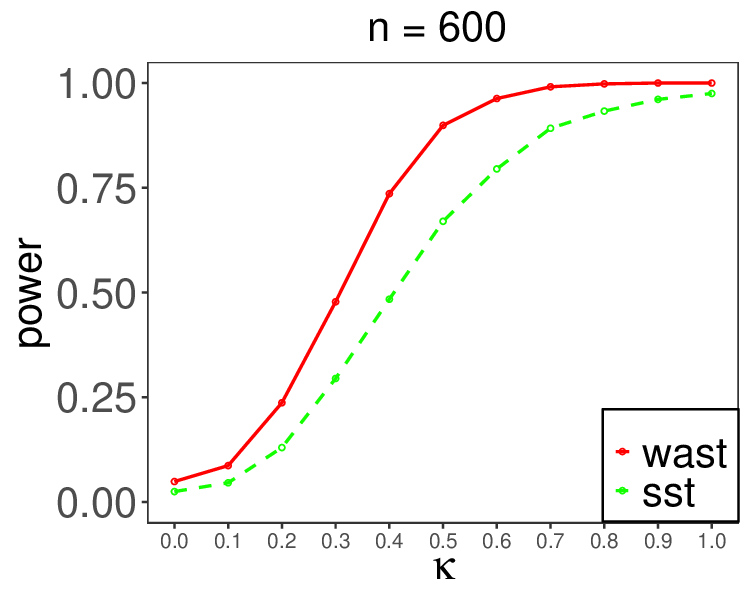}      \\
		\includegraphics[scale=0.3]{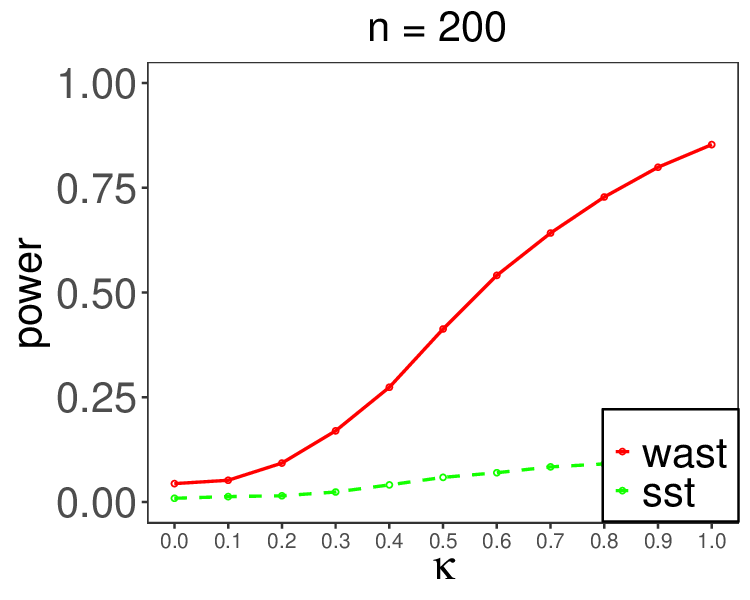}
		\includegraphics[scale=0.3]{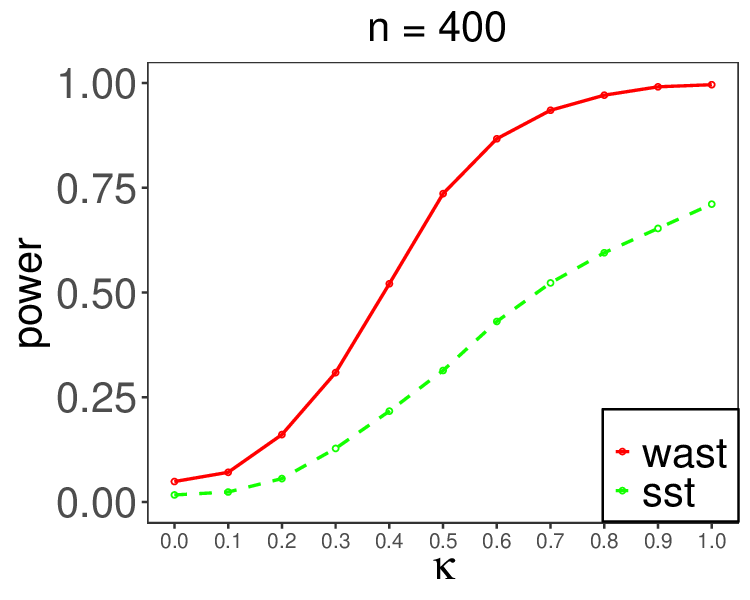}
		\includegraphics[scale=0.3]{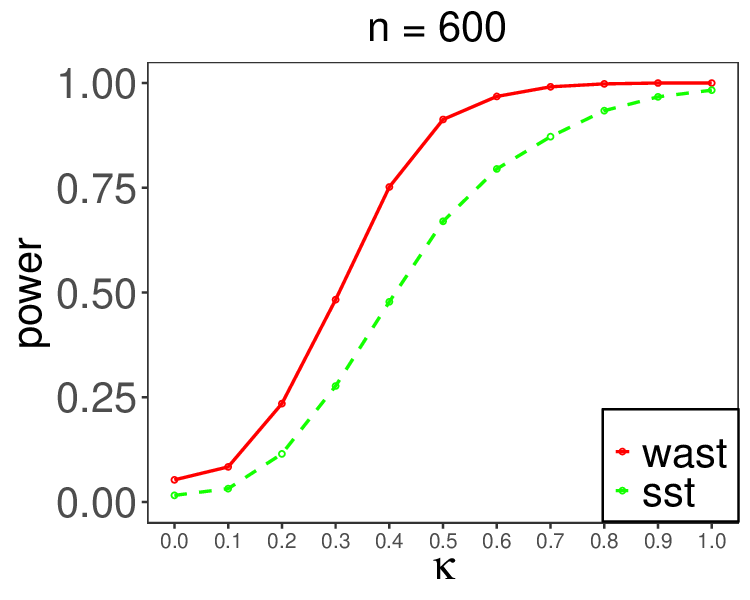}      \\
		\includegraphics[scale=0.3]{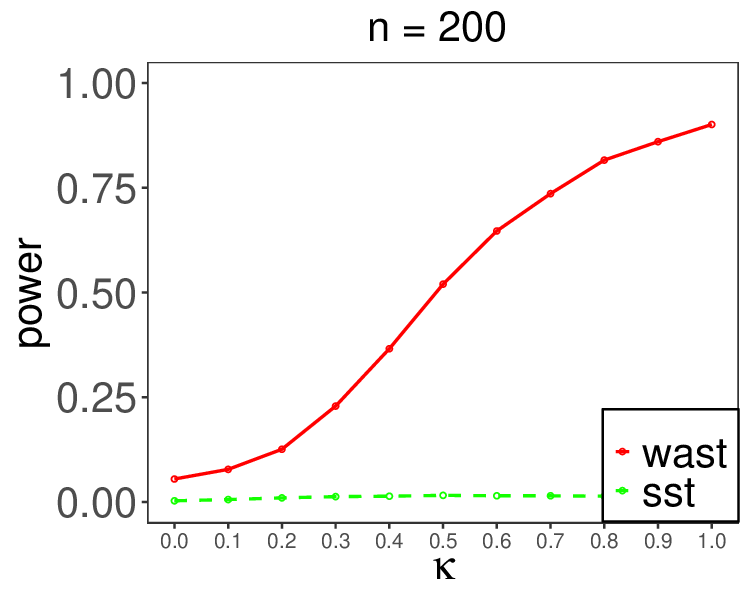}
		\includegraphics[scale=0.3]{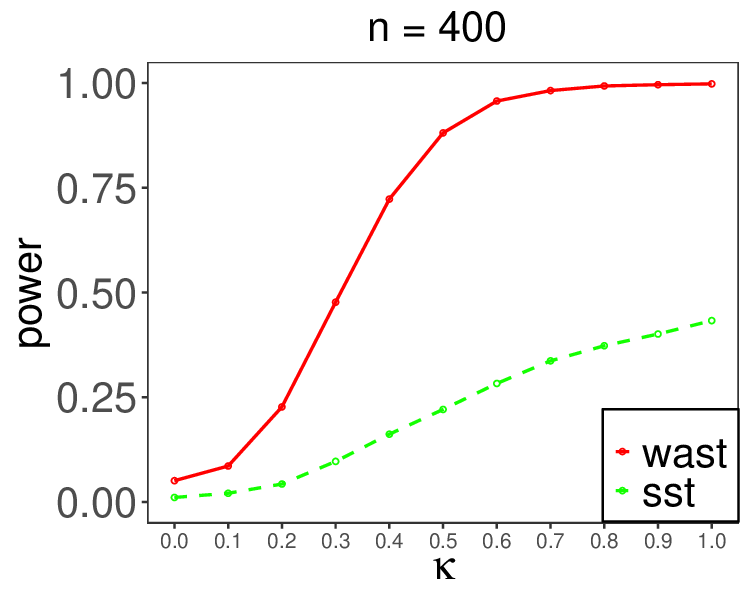}
		\includegraphics[scale=0.3]{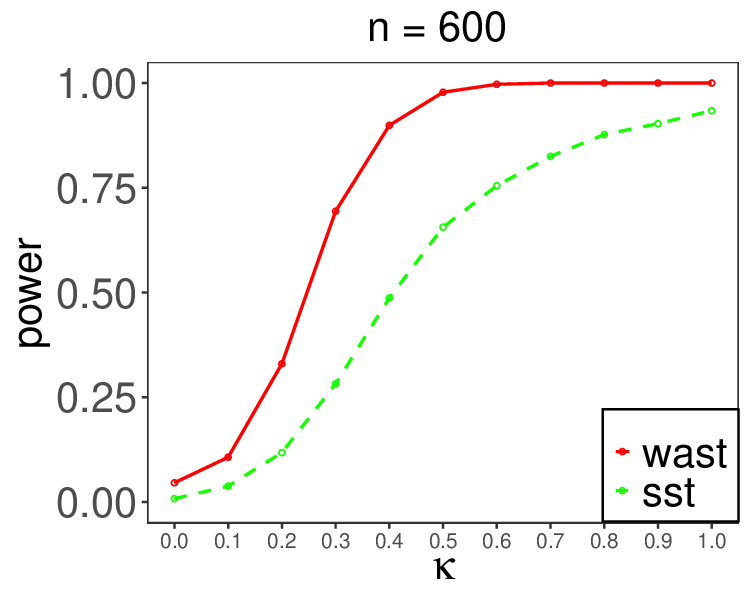}     \\
		\includegraphics[scale=0.3]{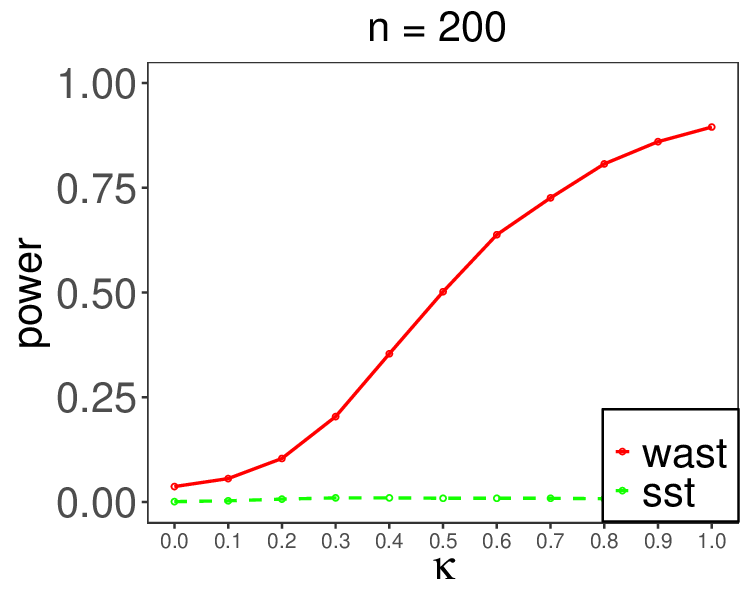}
		\includegraphics[scale=0.3]{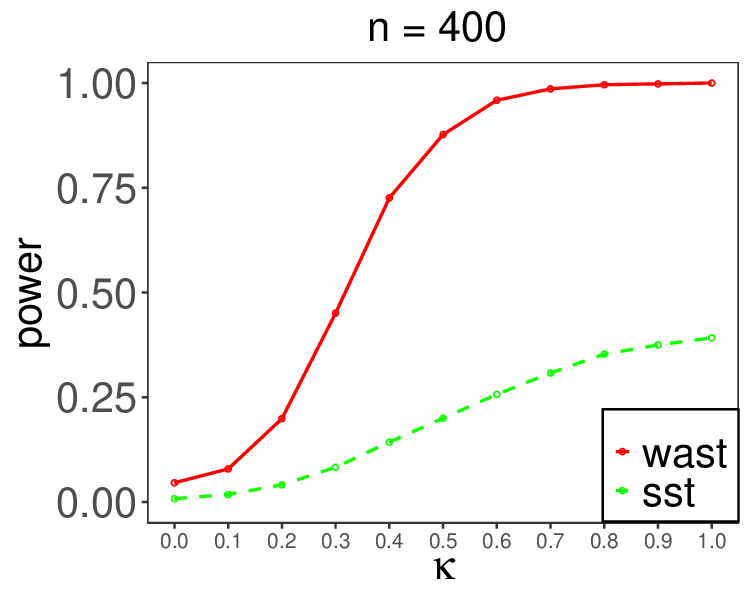}
		\includegraphics[scale=0.3]{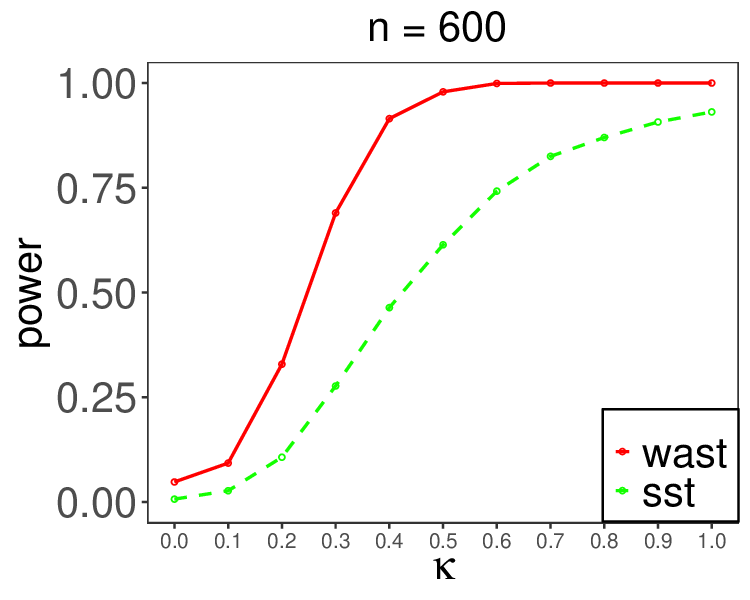}
		\caption{\it Powers of test statistic for semiparametric model (B1+P1) with large numbers of sparse $\bZ$ by the proposed WAST (red solid line) and SST (green dashed line). From top to bottom, each row depicts the powers for $(p,q)=(2,100)$, $(p,q)=(2,500)$, $(p,q)=(6,100)$, $(p,q)=(6,500)$, $(p,q)=(11,100)$, and $(p,q)=(11,500)$.}
		\label{fig_semiparam1_sparse}
	\end{center}
\end{figure}

\begin{figure}[!ht]
	\begin{center}
		\includegraphics[scale=0.3]{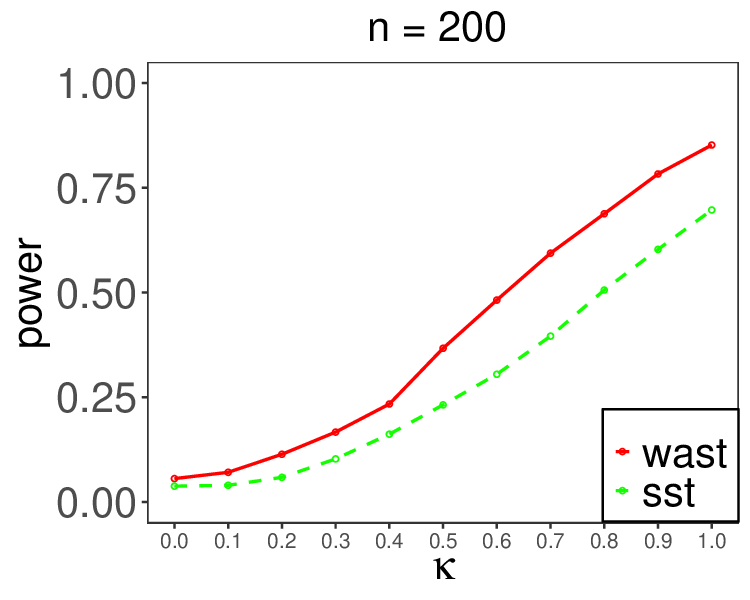}
		\includegraphics[scale=0.3]{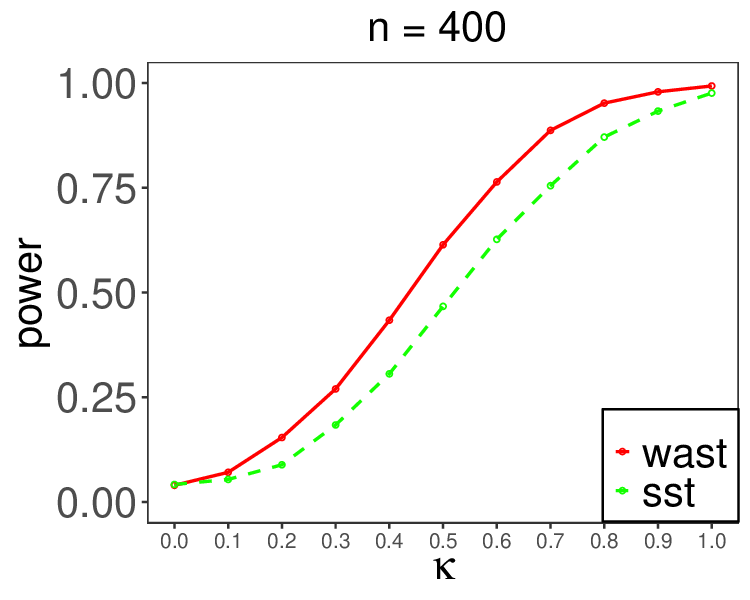}
		\includegraphics[scale=0.3]{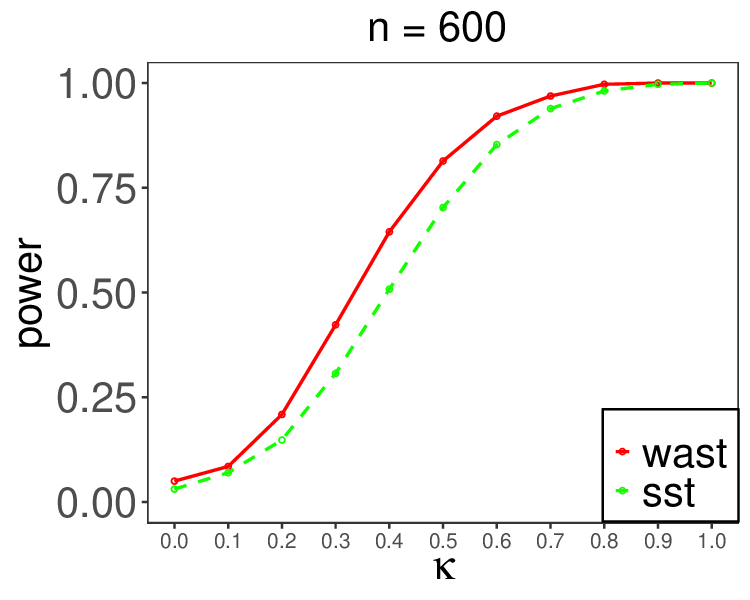}     \\
		\includegraphics[scale=0.3]{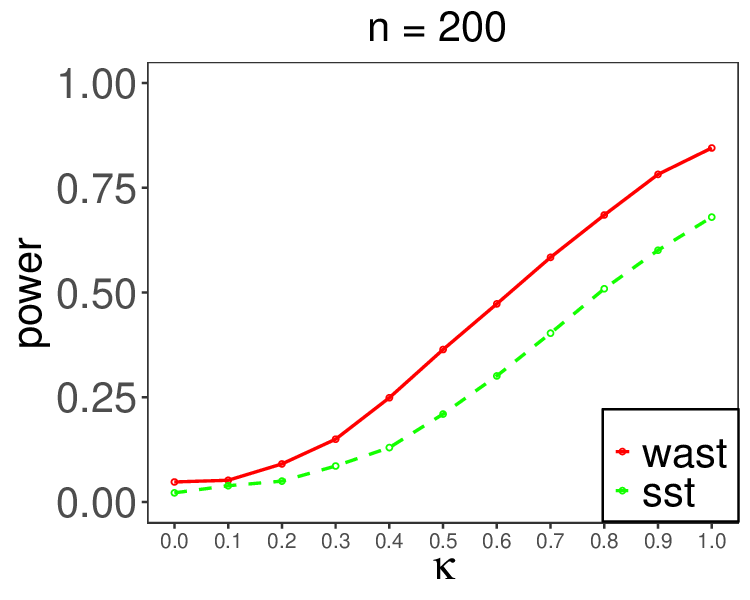}
		\includegraphics[scale=0.3]{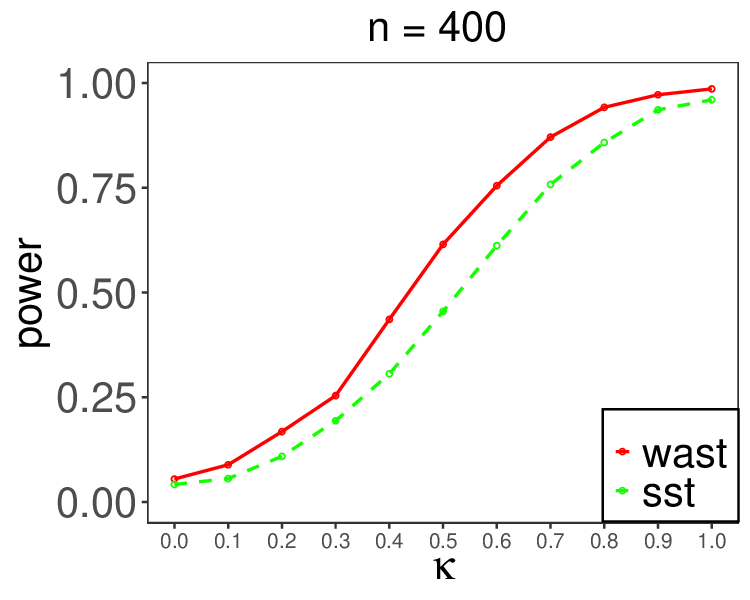}
		\includegraphics[scale=0.3]{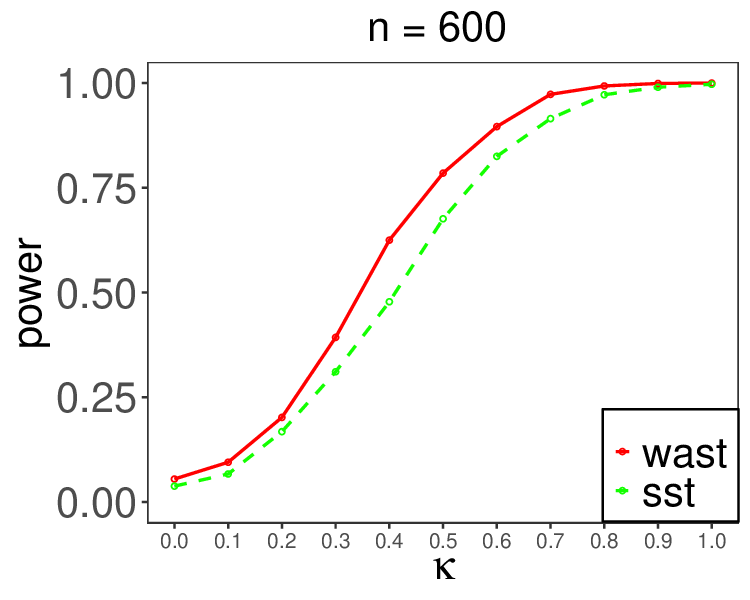}     \\
		\includegraphics[scale=0.3]{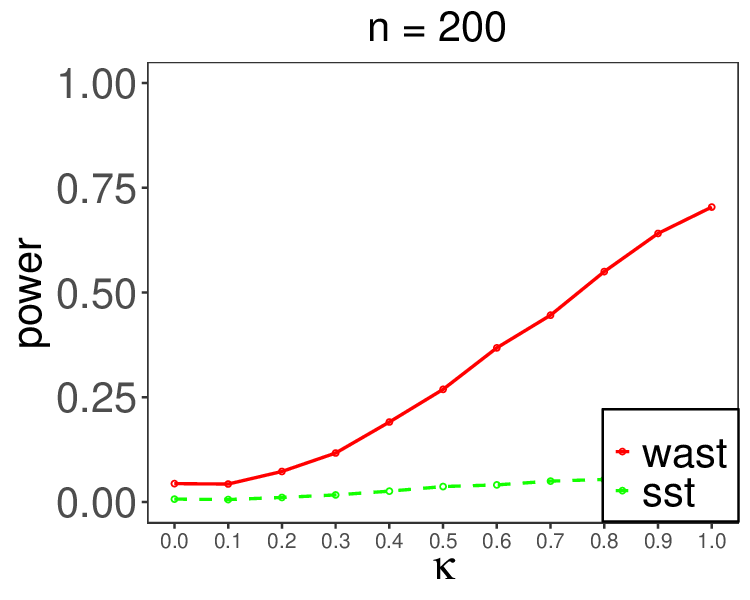}
		\includegraphics[scale=0.3]{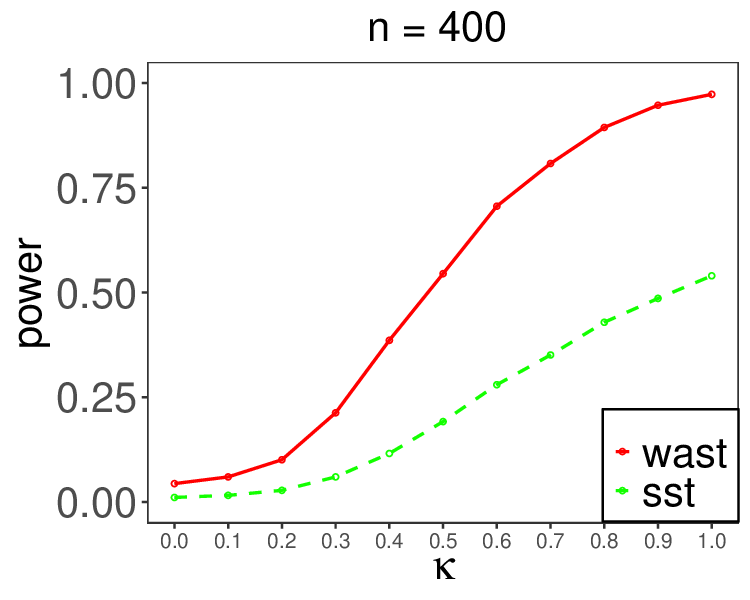}
		\includegraphics[scale=0.3]{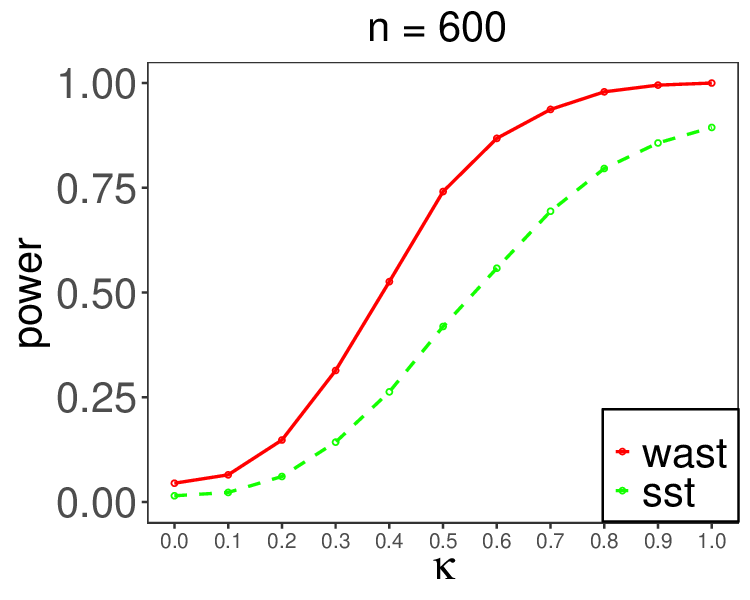}     \\
		\includegraphics[scale=0.3]{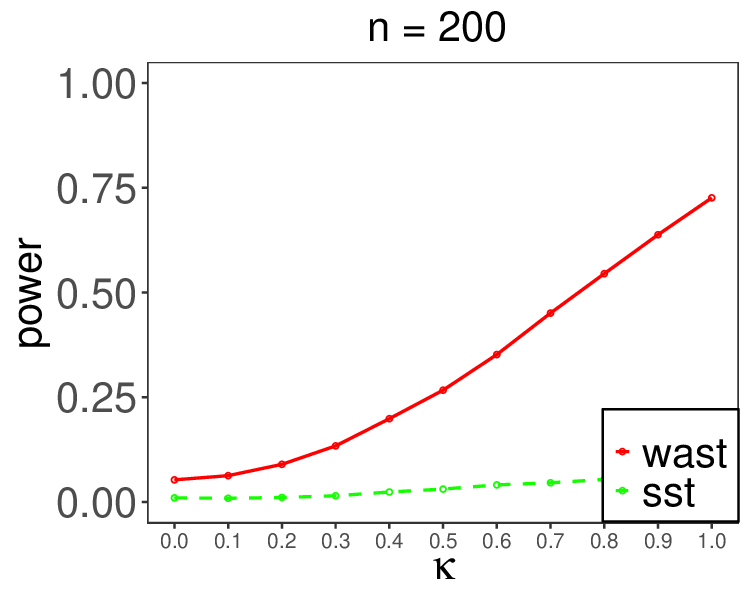}
		\includegraphics[scale=0.3]{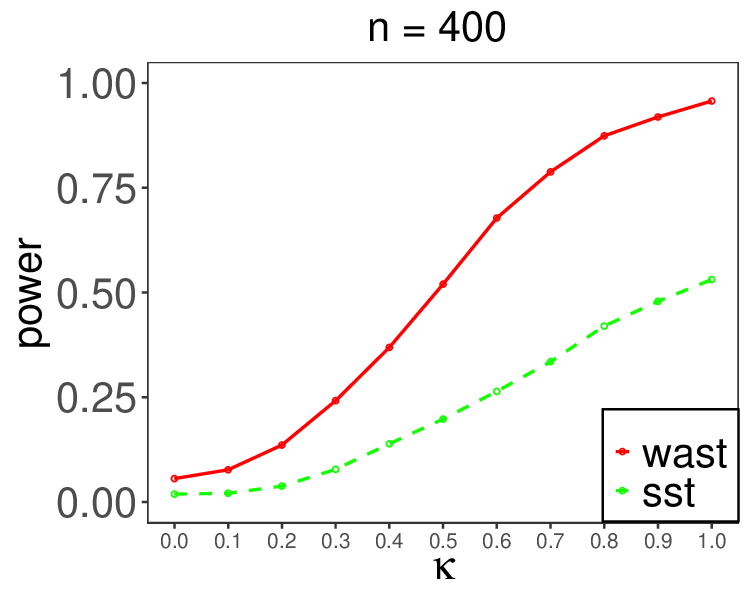}
		\includegraphics[scale=0.3]{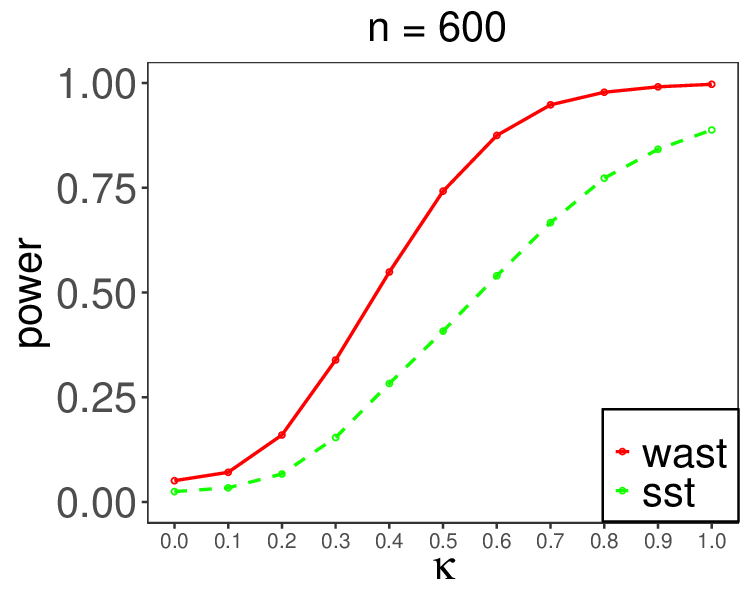}     \\
		\includegraphics[scale=0.3]{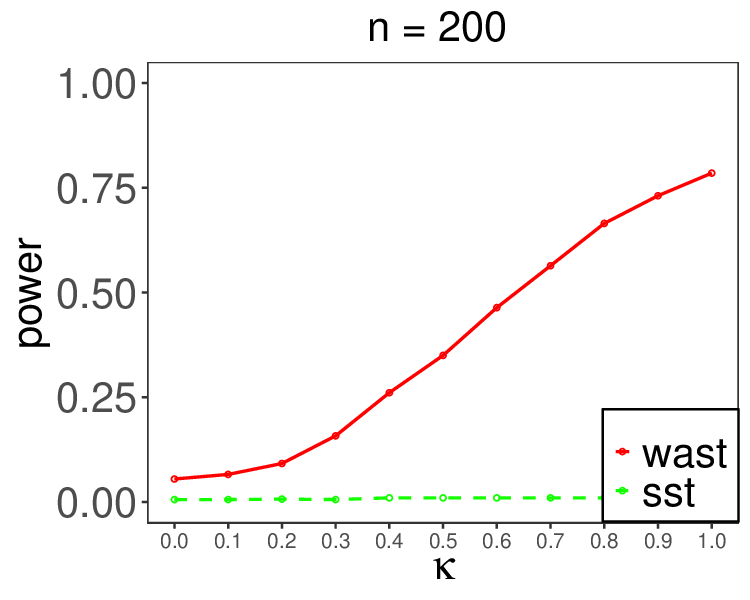}
		\includegraphics[scale=0.3]{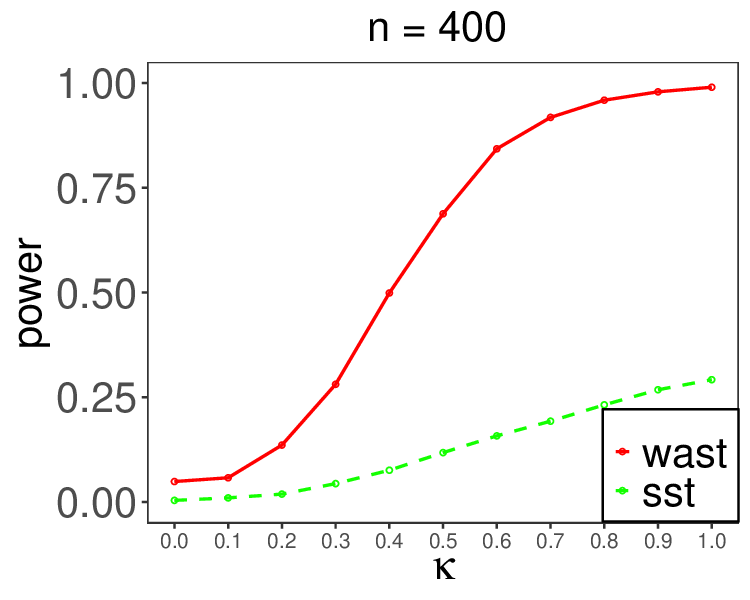}
		\includegraphics[scale=0.3]{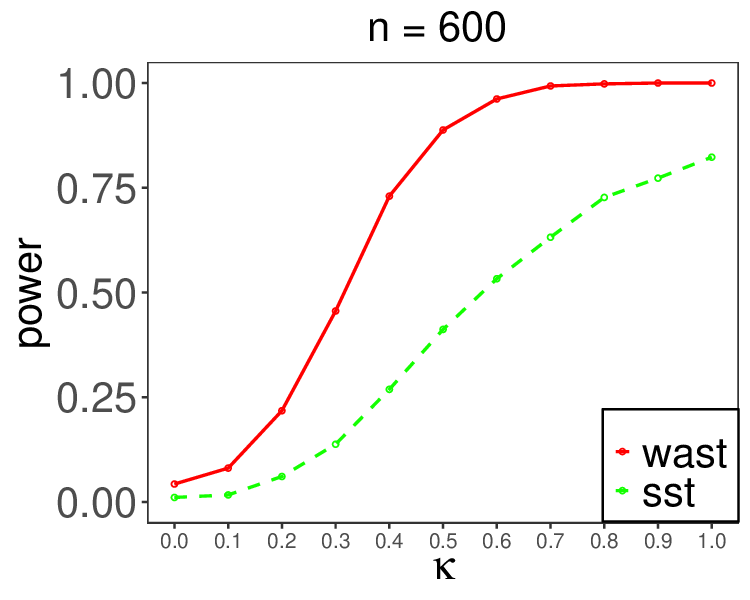}    \\
		\includegraphics[scale=0.3]{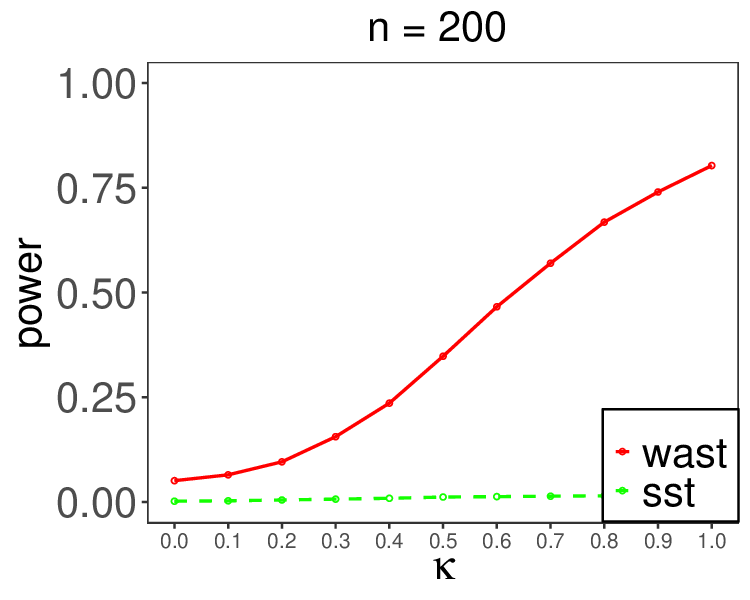}
		\includegraphics[scale=0.3]{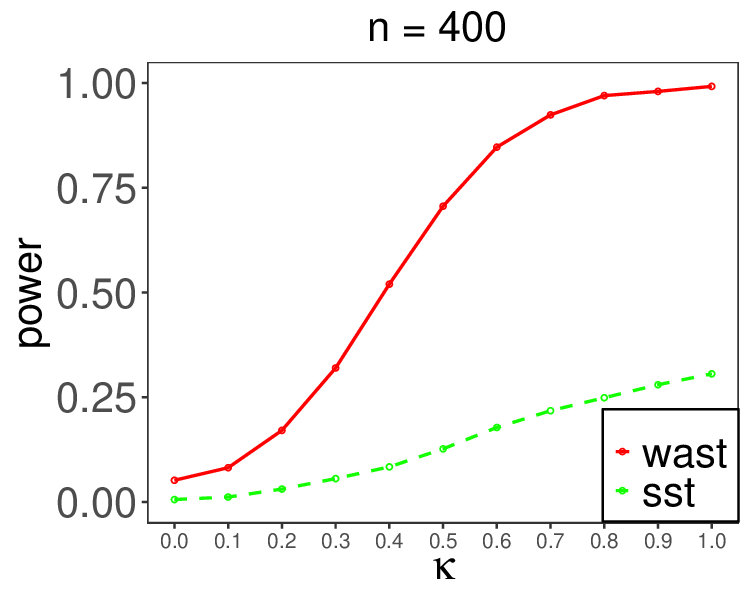}
		\includegraphics[scale=0.3]{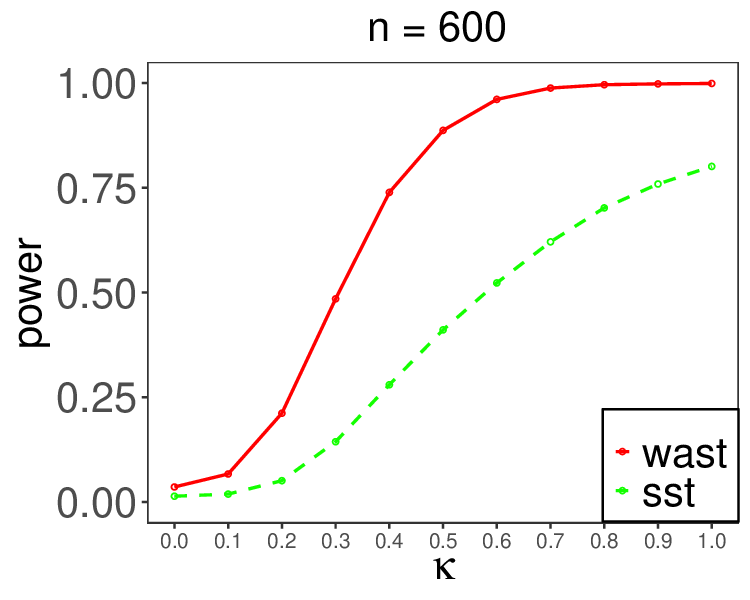}
		\caption{\it Powers of test statistic semiparametric model (B1+P1) with large numbers of sparse $\bZ$ by the proposed WAST (red solid line) and SST (green dashed line). From top to bottom, each row depicts the powers for $(p,q)=(2,100)$, $(p,q)=(2,500)$, $(p,q)=(6,100)$, $(p,q)=(6,500)$, $(p,q)=(11,100)$, and $(p,q)=(11,500)$.}
		\label{fig_semiparam2_sparse}
	\end{center}
\end{figure}

For all models with change plane analysis, we evaluated the power under a sequence of alternative models indexed by $\kappa$, that is $H_1^{\kappa}: \bbeta^{\kappa}=\kappa\bbeta^*$ with $\kappa = i/10$ for semiparametric model and $\kappa=i/20$ for others, $i=1,\cdots,10$, where $\bbeta^*=(1,\cdots,1)\trans$. We set sample size $n=(200, 400, 600)$ for quantile, probit and semiparametric models, 1000 repetitions and 1000 bootstrap samples, and report in Figure \ref{fig_probit_sparse}-\ref{fig_semiparam2_sparse} the performance for the WAST and SST. We calculate the SST over $\{\btheta^{(k)}=(\theta^{(k)}_1,\cdots,\theta^{(k)}_q)\trans: k=1,\cdots,K\}$ with the number of threshold values $K=1000$. Let $\btheta^{(k)}_{-1}=\tilde{\btheta}^{(k)}_{-1}/\|\tilde{\btheta}^{(k)}_{-1}\|$, where $\btheta^{(k)}_{-1}=(\theta^{(k)}_2,\cdots,\theta^{(k)}_q)\trans$ and $\tilde{\btheta}^{(k)}_{-1}=(\tilde{\theta}^{(k)}_2,\cdots,\tilde{\theta}^{(k)}_q)\trans$, and $\tilde{\btheta}^{(k)}_{-1}$ is drawn independently form $(r-1)$-variate standard normal distribution.  For each $\theta_1^{(k)}$, $k=1,\cdots,K$, we selected it by equal grid search in the range from the lower 10th percentile to upper 10th percentile of the data points of $\{\theta^{(k)}_2Z_{2i}+\cdots+\theta^{(k)}_qZ_{qi}\}_{i=1}^n$, which is same as that in \cite{2020Threshold}. Here we consider four combinations of $(p,q)=(2, 100), (2, 500), (6, 100), (6, 500), (11, 100), (11, 500)$.

Type \uppercase\expandafter{\romannumeral1} errors ($\kappa=0$) for probit, quantile and semiparametric models are listed in Table \ref{table_size_ee_sparse}. We can see from Table \ref{table_size_ee_sparse} that the size of the proposed WAST are close to the nominal significance level $0.05$, but for most scenarios the size of the SST are much smaller than 0.05.
Figure \ref{fig_probit_sparse}-\ref{fig_semiparam2_sparse} indicate that powers become greater as sample size $n$ increases, which are verified the asymptotic theory. The proposed WAST has comparable power with the SST for the semiparametric model, but the size of the SST is much less than the nominal level 0.05.

\begin{table}[htp!]
	\def~{\hphantom{0}}
	\tiny
	\caption{Type \uppercase\expandafter{\romannumeral1} errors of the proposed WAST and SST based on resampling for probit regression model (ProbitRE), quantile regression (QuantRE) and semiparametric model (SPMoldel) with large numbers of sparse $\bZ$. The nominal significant level is 0.05.
	}
	\resizebox{\textwidth}{!}{
		\begin{threeparttable}
			\begin{tabular}{llcccccccc}
				\hline
				\multirow{2}{*}{Model}&\multirow{2}{*}{$(p,q)$}
				&\multicolumn{2}{c}{ $n=200$} && \multicolumn{2}{c}{ $n=400$} && \multicolumn{2}{c}{ $n=600$} \\
				\cline{3-4} \cline{6-7} \cline{9-10}
				&&WAST& SST && WAST& SST && WAST& SST \\
				\cline{3-10}
				ProbitRE &$(2,100)$         & 0.035 & 0.001 && 0.058 & 0.002 && 0.040 & 0.004 \\
				&$(2,500)$                  & 0.035 & 0.002 && 0.047 & 0.007 && 0.056 & 0.013 \\
				&$(6,100)$                  & 0.052 & 0.000 && 0.041 & 0.001 && 0.056 & 0.000 \\
				&$(6,500)$                  & 0.056 & 0.000 && 0.046 & 0.000 && 0.056 & 0.003 \\
				&$(11,100)$                 & 0.056 & 0.000 && 0.048 & 0.000 && 0.049 & 0.000 \\
				&$(11,500)$                 & 0.061 & 0.000 && 0.040 & 0.000 && 0.053 & 0.000 \\
				[1 ex]
				SPModel &$(2,100)$          & 0.056 & 0.028 && 0.047 & 0.029 && 0.068 & 0.055 \\
				(B1+P1)&$(2,500)$           & 0.060 & 0.029 && 0.049 & 0.044 && 0.054 & 0.046 \\
				&$(6,100)$                  & 0.065 & 0.008 && 0.049 & 0.019 && 0.049 & 0.025 \\
				&$(6,500)$                  & 0.044 & 0.009 && 0.049 & 0.017 && 0.053 & 0.016 \\
				&$(11,100)$                 & 0.055 & 0.003 && 0.051 & 0.011 && 0.046 & 0.008 \\
				&$(11,500)$                 & 0.037 & 0.001 && 0.046 & 0.008 && 0.048 & 0.007 \\
				[1 ex]
				SPModel &$(2,100)$          & 0.056 & 0.038 && 0.040 & 0.042 && 0.050 & 0.031 \\
				(B2+P2)&$(2,500)$           & 0.048 & 0.022 && 0.055 & 0.042 && 0.055 & 0.038 \\
				&$(6,100)$                  & 0.044 & 0.007 && 0.044 & 0.011 && 0.045 & 0.015 \\
				&$(6,500)$                  & 0.053 & 0.010 && 0.056 & 0.019 && 0.051 & 0.025 \\
				&$(11,100)$                 & 0.055 & 0.006 && 0.049 & 0.004 && 0.043 & 0.011 \\
				&$(11,500)$                 & 0.051 & 0.002 && 0.052 & 0.006 && 0.036 & 0.014 \\
				[1 ex]
				QuantRE &$(2,100)$          & 0.054 & 0.016 && 0.059 & 0.025 && 0.058 & 0.032 \\
				($\tau=0.2$)&$(2,500)$      & 0.054 & 0.018 && 0.050 & 0.017 && 0.050 & 0.032 \\
				&$(6,100)$                  & 0.059 & 0.005 && 0.054 & 0.010 && 0.057 & 0.017 \\
				&$(6,500)$                  & 0.054 & 0.003 && 0.048 & 0.005 && 0.062 & 0.018 \\
				&$(11,100)$                 & 0.050 & 0.000 && 0.054 & 0.004 && 0.053 & 0.011 \\
				&$(11,500)$                 & 0.045 & 0.001 && 0.050 & 0.007 && 0.057 & 0.017 \\
				[1 ex]
				QuantRE &$(2,100)$          & 0.050 & 0.031 && 0.036 & 0.030 && 0.047 & 0.048 \\
				($\tau=0.5$)&$(2,500)$      & 0.056 & 0.039 && 0.058 & 0.044 && 0.052 & 0.047 \\
				&$(6,100)$                  & 0.050 & 0.022 && 0.037 & 0.018 && 0.051 & 0.026 \\
				&$(6,500)$                  & 0.058 & 0.012 && 0.052 & 0.016 && 0.050 & 0.041 \\
				&$(11,100)$                 & 0.065 & 0.007 && 0.062 & 0.019 && 0.048 & 0.017 \\
				&$(11,500)$                 & 0.043 & 0.008 && 0.047 & 0.017 && 0.053 & 0.031 \\
				[1 ex]
				QuantRE &$(2,100)$          & 0.053 & 0.017 && 0.055 & 0.022 && 0.047 & 0.037 \\
				($\tau=0.7$)&$(2,500)$      & 0.057 & 0.026 && 0.066 & 0.040 && 0.058 & 0.033 \\
				&$(6,100)$                  & 0.043 & 0.007 && 0.060 & 0.018 && 0.054 & 0.025 \\
				&$(6,500)$                  & 0.046 & 0.009 && 0.044 & 0.015 && 0.058 & 0.033 \\
				&$(11,100)$                 & 0.049 & 0.002 && 0.061 & 0.008 && 0.044 & 0.019 \\
				&$(11,500)$                 & 0.037 & 0.003 && 0.047 & 0.009 && 0.056 & 0.014 \\
				\hline
			\end{tabular}
		\end{threeparttable}
	}
	\label{table_size_ee_sparse}
\end{table}

\subsection{Change plane analysis with dense \texorpdfstring{$\bZ$}{} for quantile, probit and semiparametric models}\label{simulation_ee_dense}

For the sparse $\bZ$, we generate $\theta_2,\cdots,\theta_q$ from the uniform distribution $U(0,3)$. Other settings are same as these in Section \ref{simulation_ee_sparse}.

\begin{figure}[!ht]
	\begin{center}
		\includegraphics[scale=0.3]{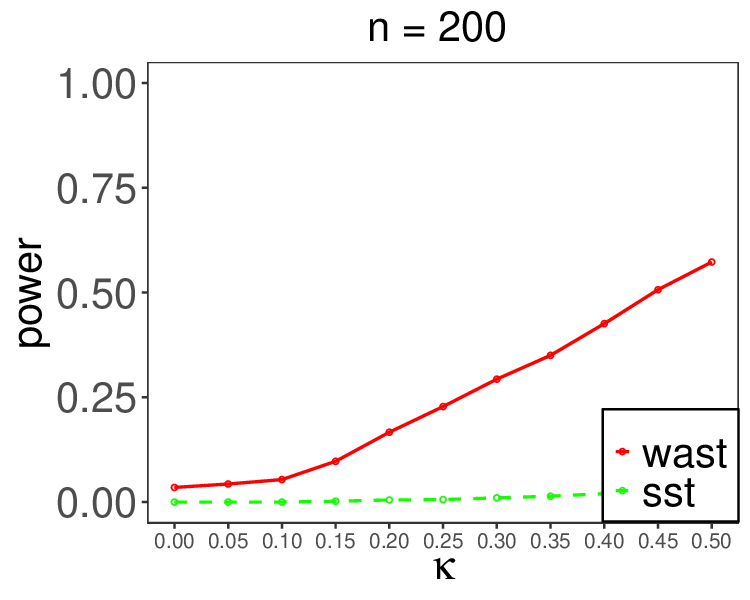}
        \includegraphics[scale=0.3]{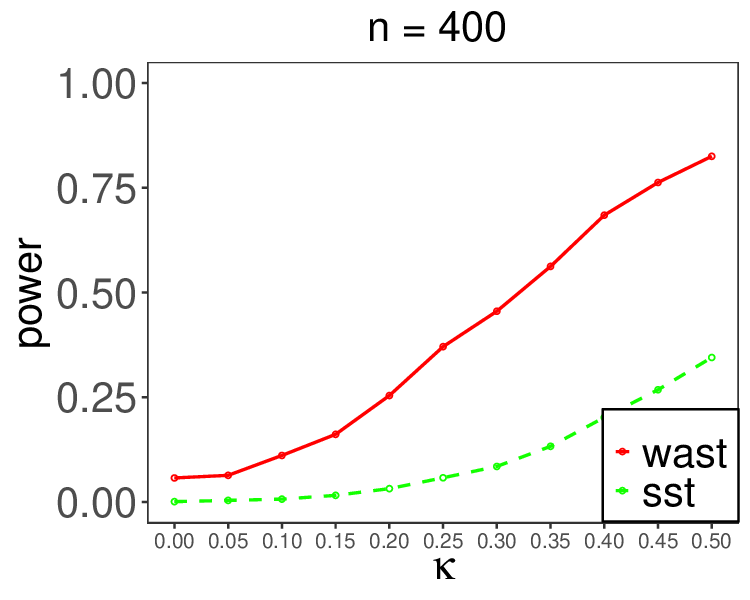}
        \includegraphics[scale=0.3]{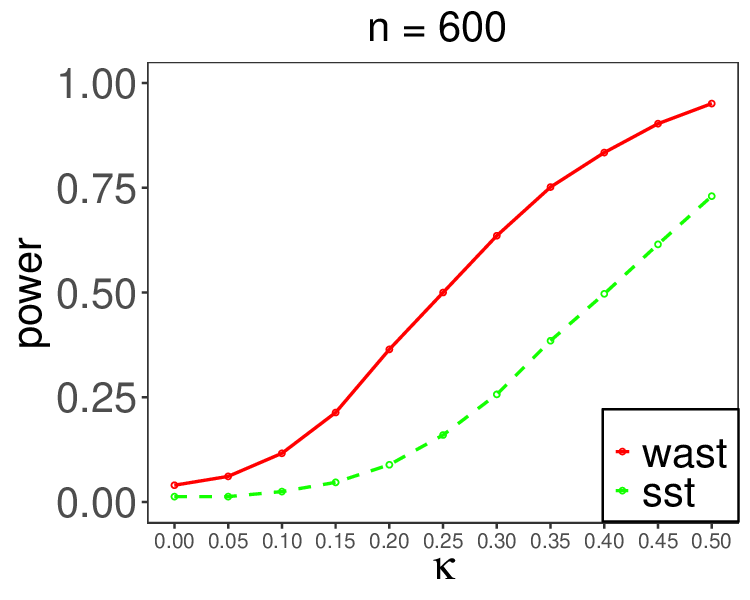}  \\
		\includegraphics[scale=0.3]{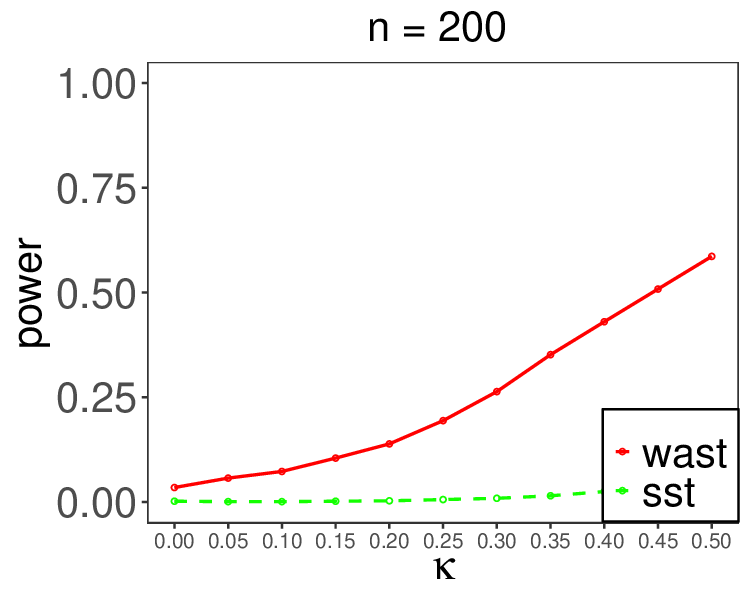}
        \includegraphics[scale=0.3]{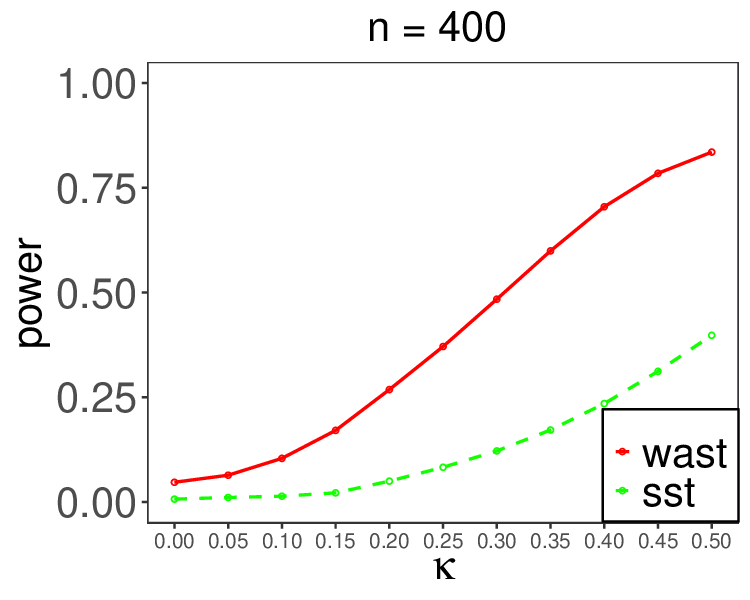}
        \includegraphics[scale=0.3]{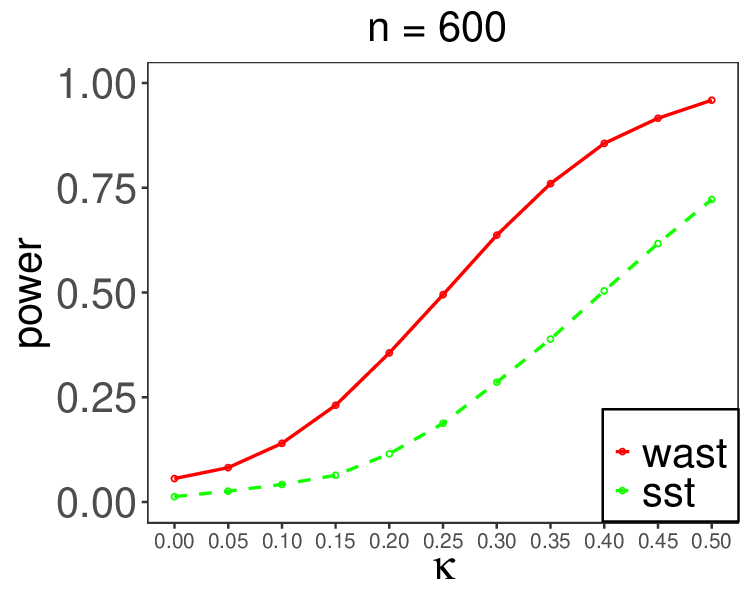}  \\
		\includegraphics[scale=0.3]{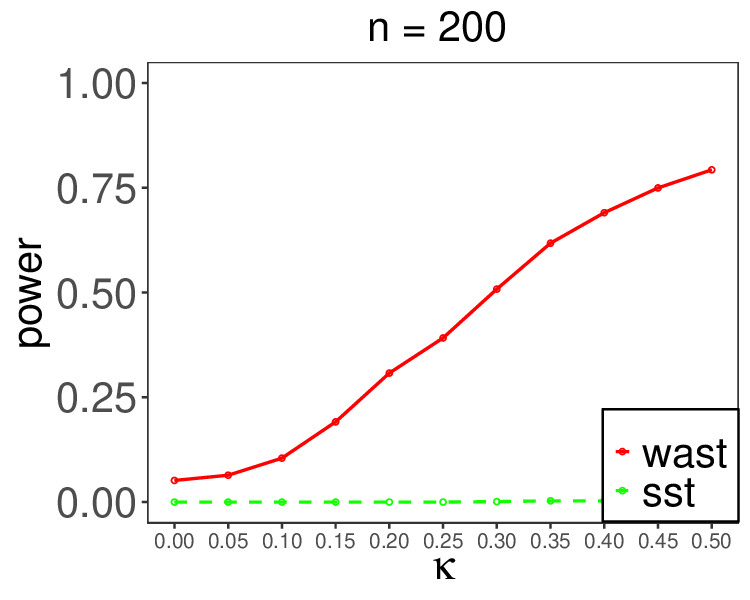}
        \includegraphics[scale=0.3]{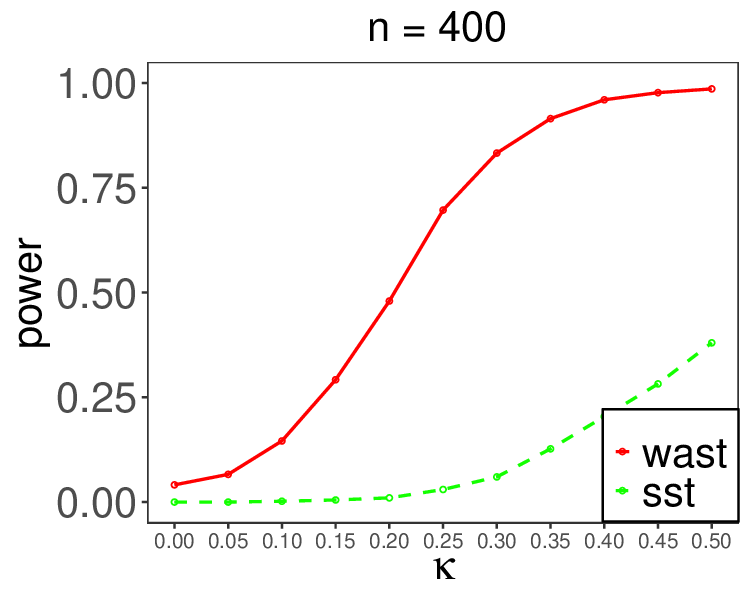}
        \includegraphics[scale=0.3]{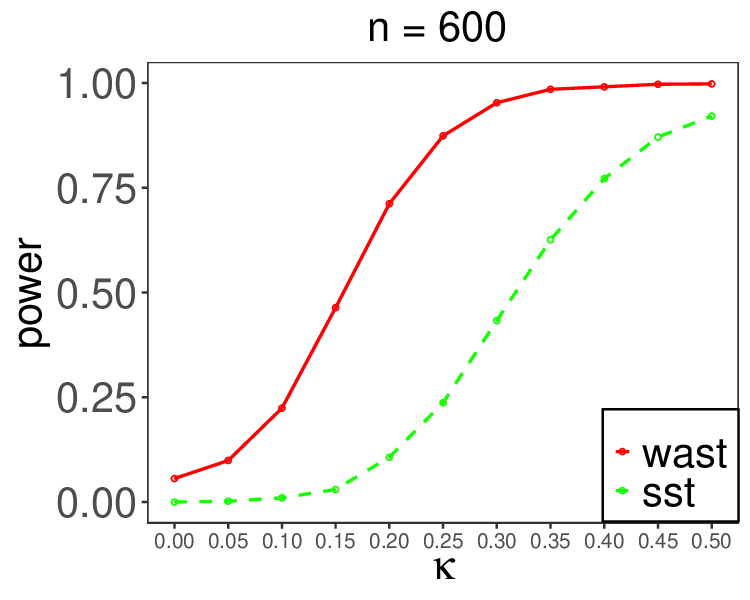}  \\
		\includegraphics[scale=0.3]{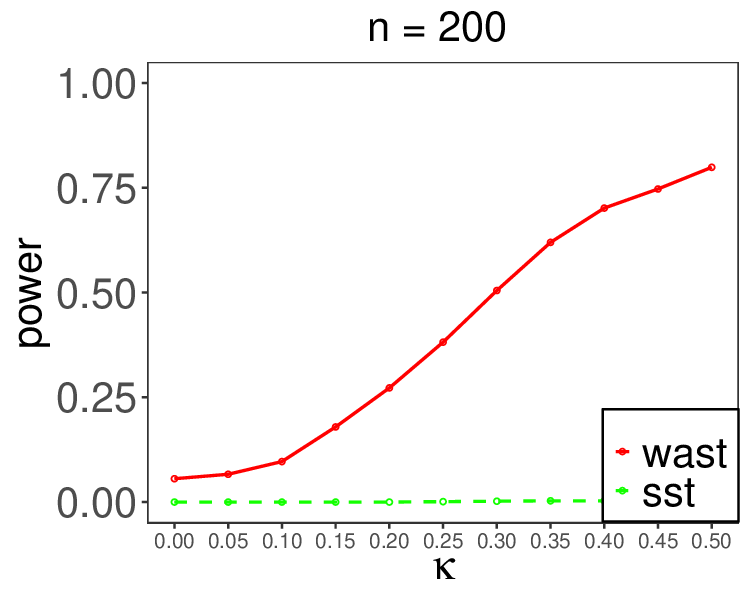}
        \includegraphics[scale=0.3]{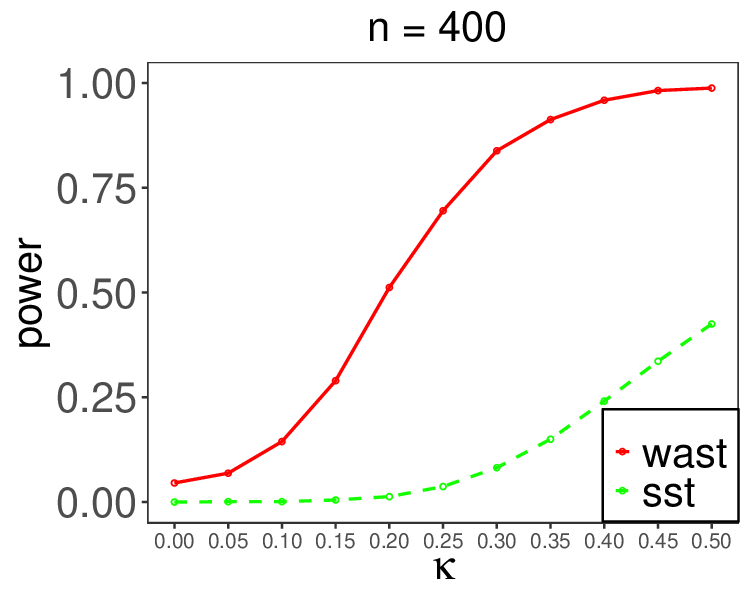}
        \includegraphics[scale=0.3]{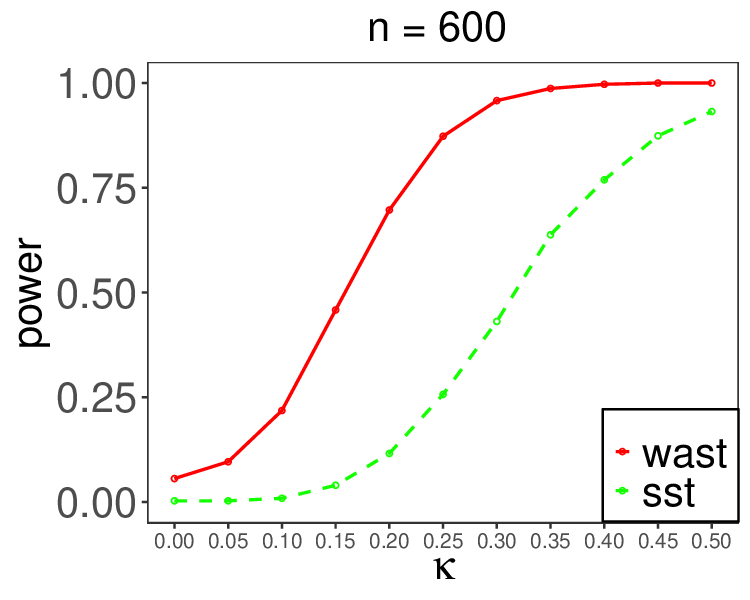}  \\
		\includegraphics[scale=0.3]{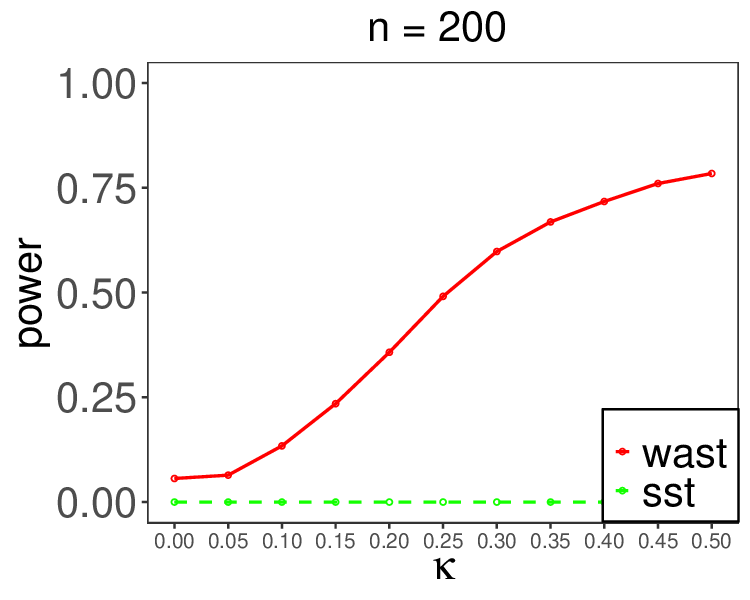}
        \includegraphics[scale=0.3]{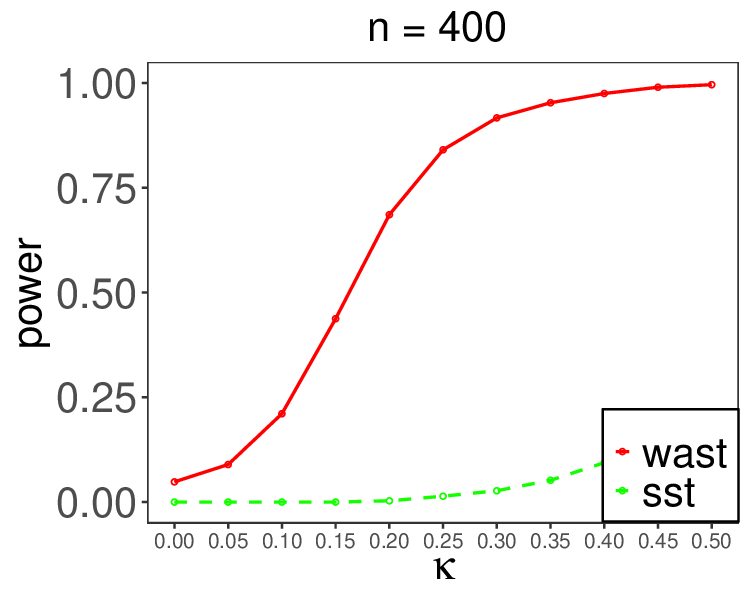}
        \includegraphics[scale=0.3]{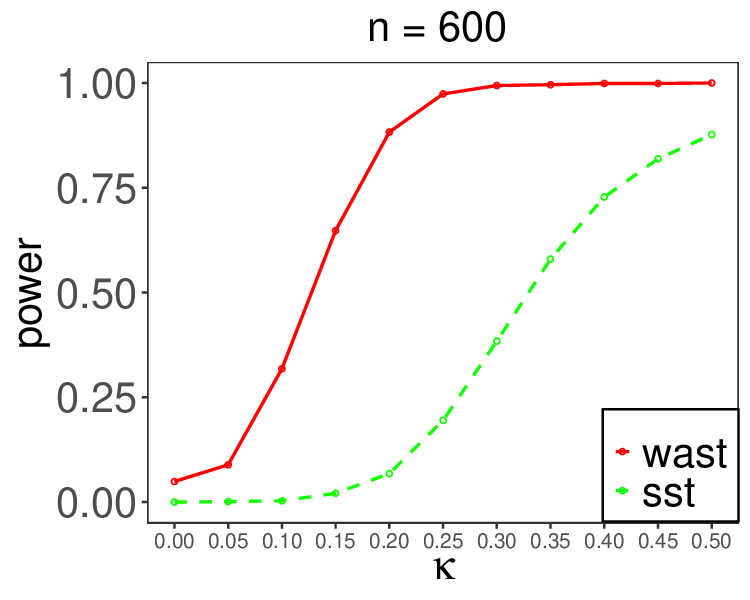} \\
		\includegraphics[scale=0.3]{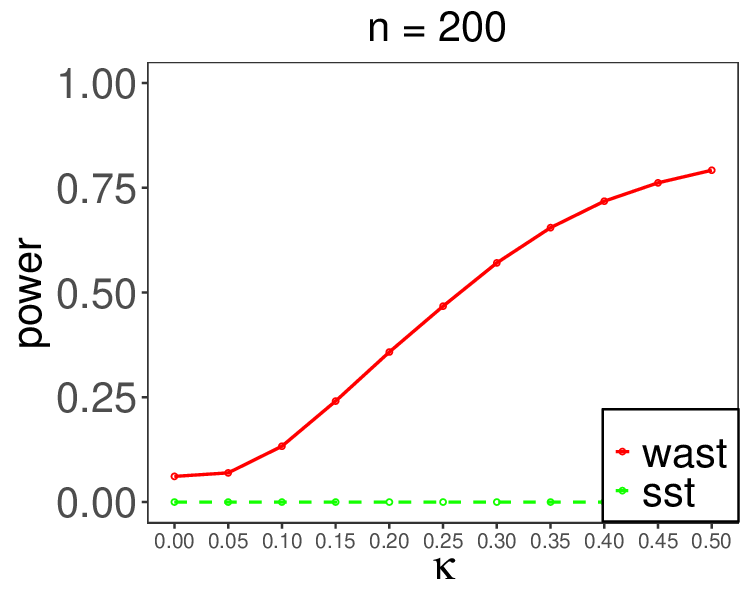}
        \includegraphics[scale=0.3]{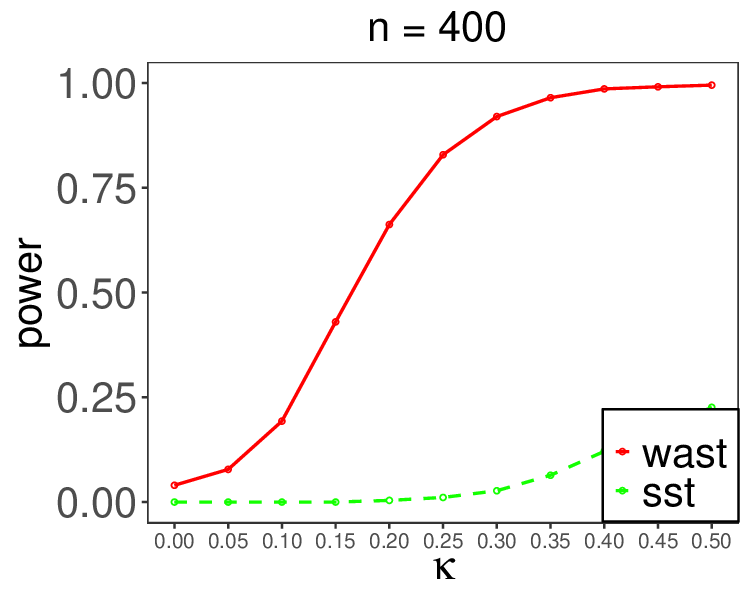}
        \includegraphics[scale=0.3]{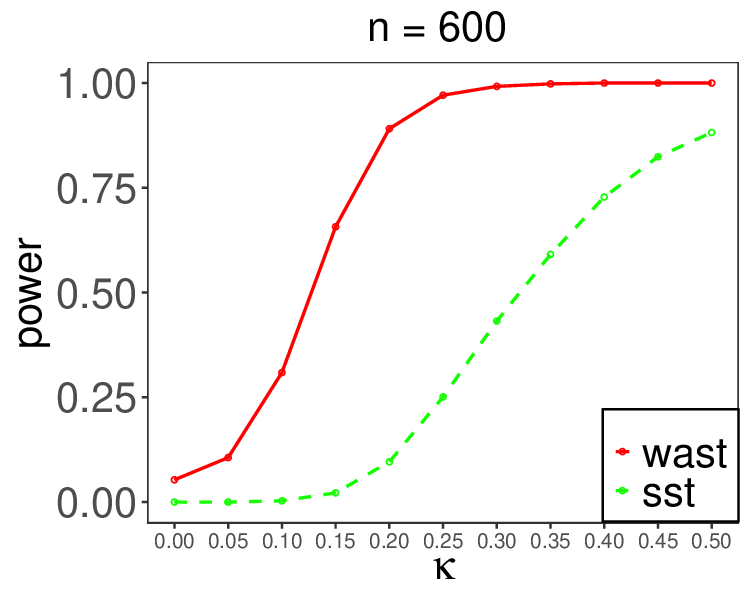}
		\caption{\it Powers of test statistic for probit model with large numbers of dense $\bZ$ by the proposed WAST (red solid line) and SST (green dashed line). From top to bottom, each row depicts the powers for $(p,q)=(2,100)$, $(p,q)=(2,500)$, $(p,q)=(6,100)$, $(p,q)=(6,500)$, $(p,q)=(11,100)$, and $(p,q)=(11,500)$.}
		\label{fig_probit_dense}
	\end{center}
\end{figure}

\begin{figure}[!ht]
	\begin{center}
		\includegraphics[scale=0.3]{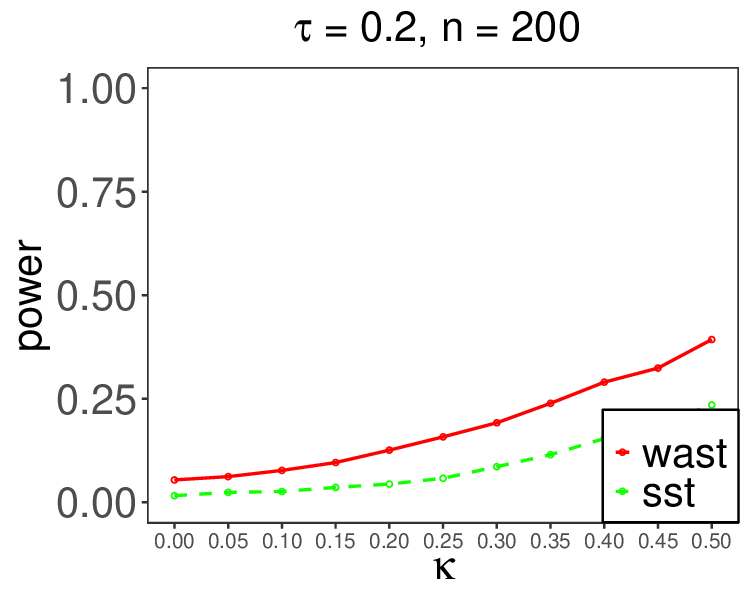}
		\includegraphics[scale=0.3]{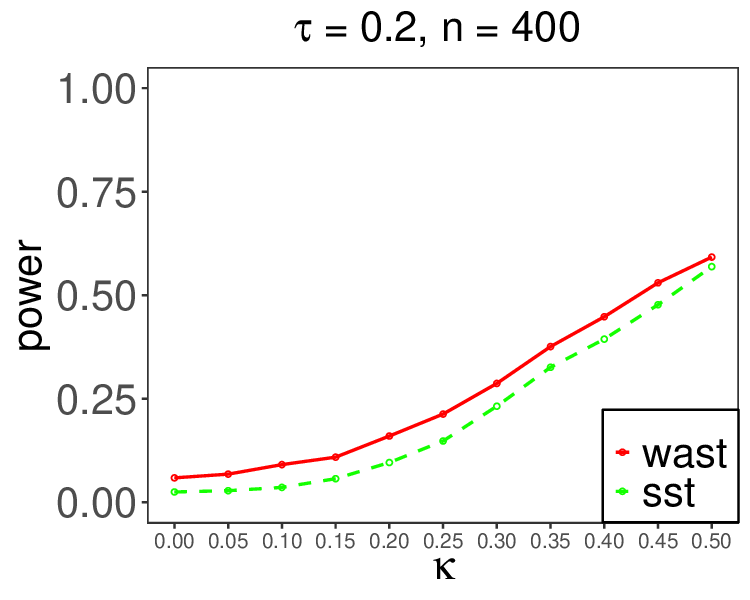}
		\includegraphics[scale=0.3]{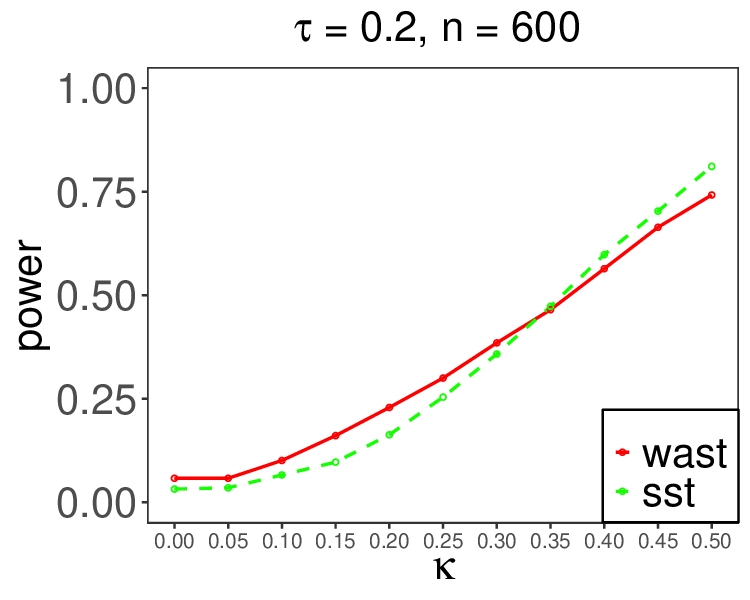}  \\
		\includegraphics[scale=0.3]{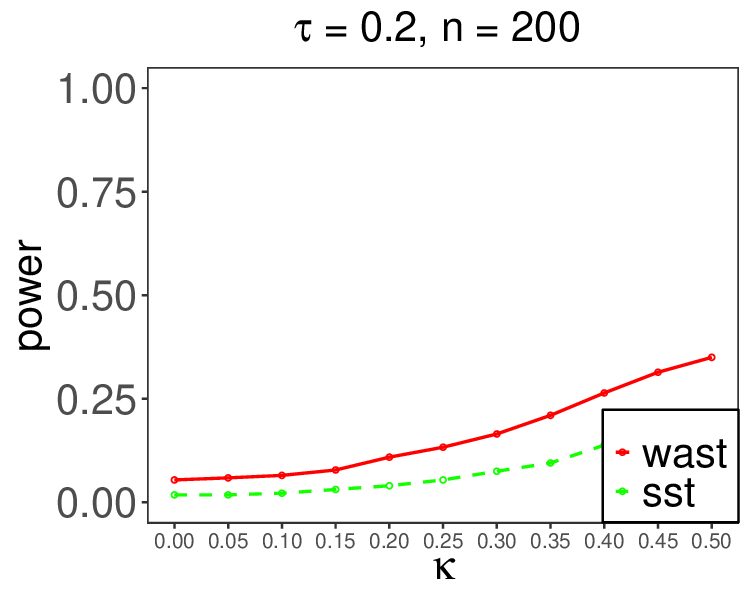}
		\includegraphics[scale=0.3]{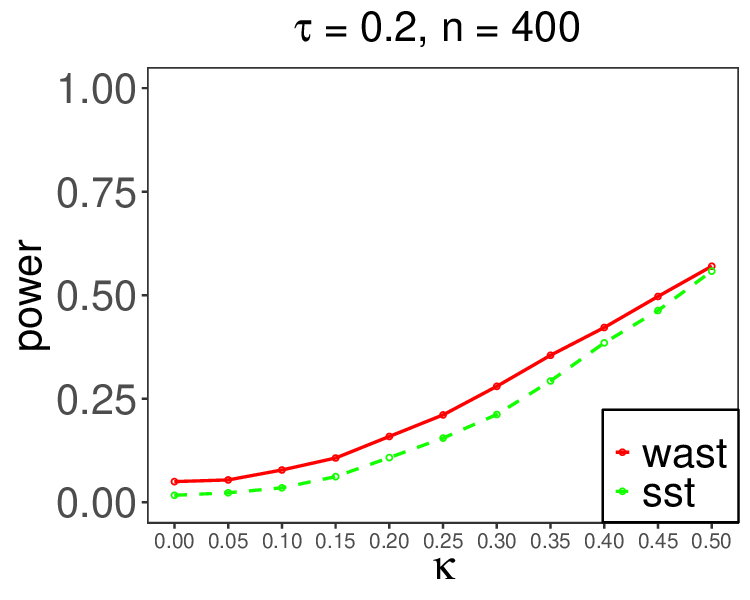}
		\includegraphics[scale=0.3]{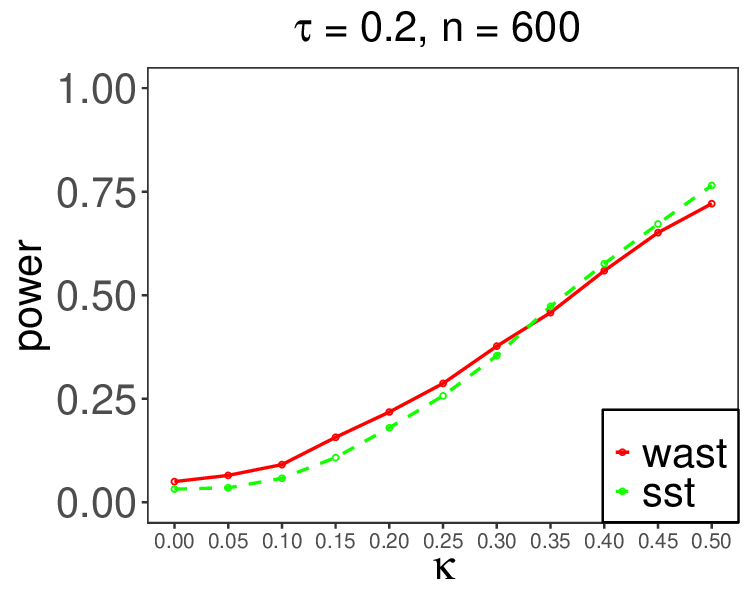}  \\
		\includegraphics[scale=0.3]{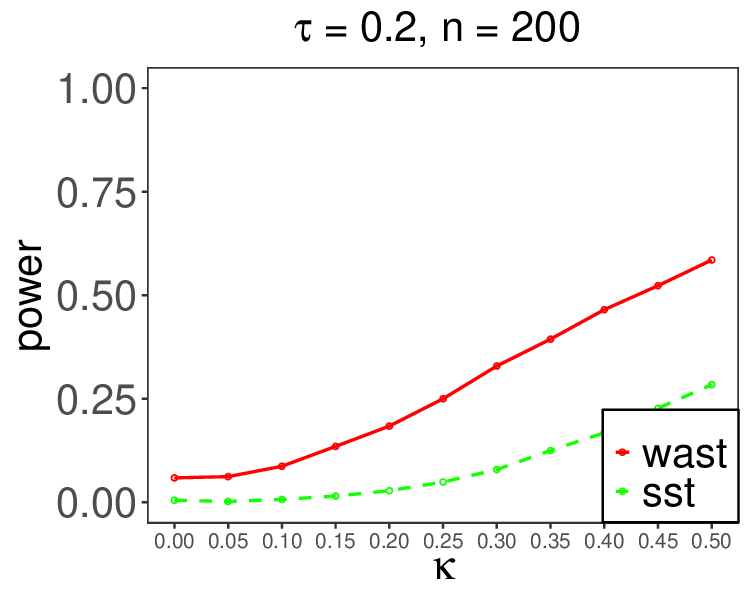}
		\includegraphics[scale=0.3]{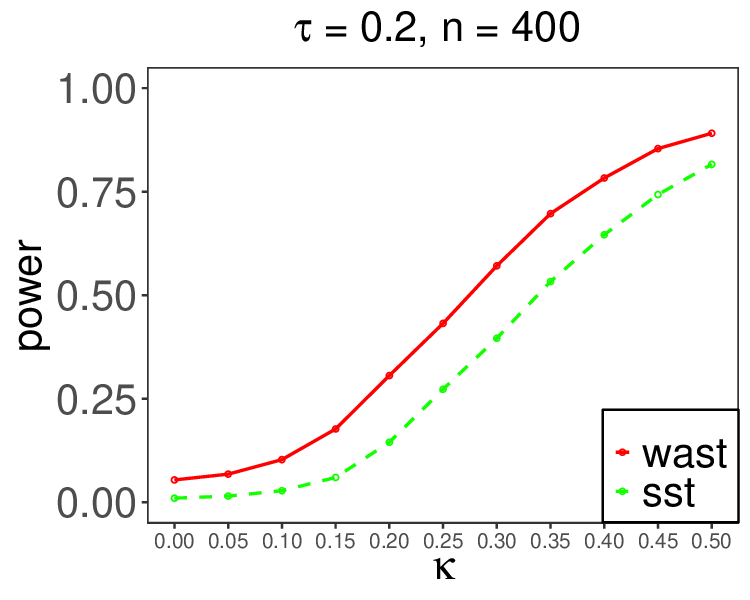}
		\includegraphics[scale=0.3]{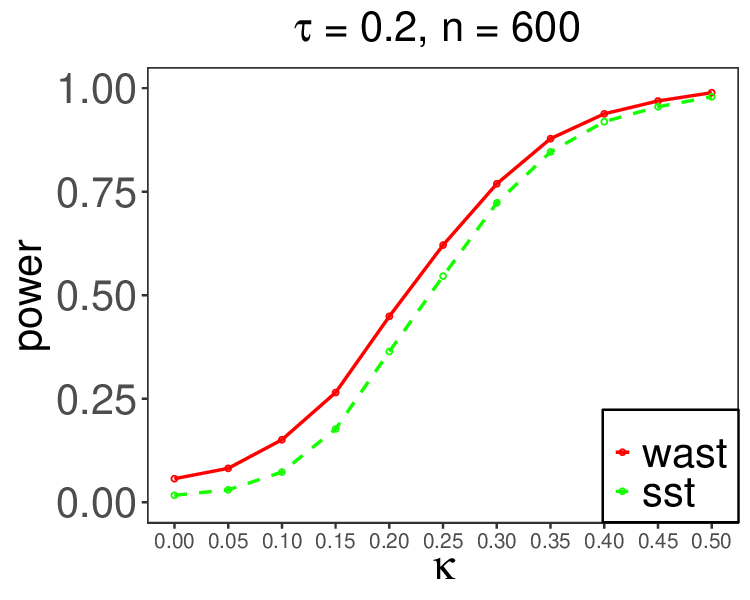}  \\
		\includegraphics[scale=0.3]{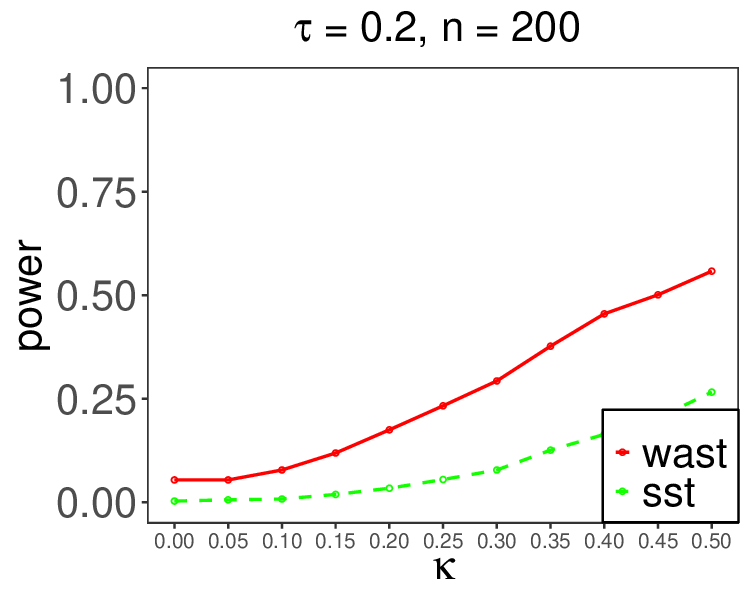}
		\includegraphics[scale=0.3]{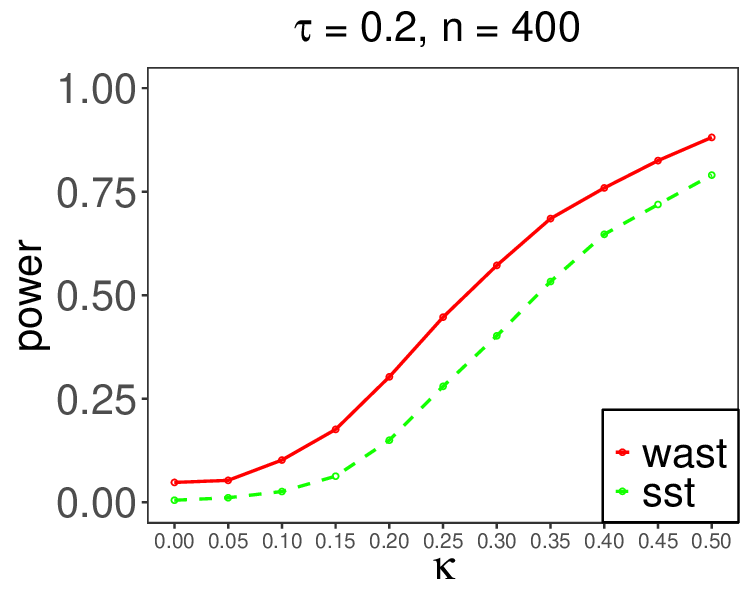}
		\includegraphics[scale=0.3]{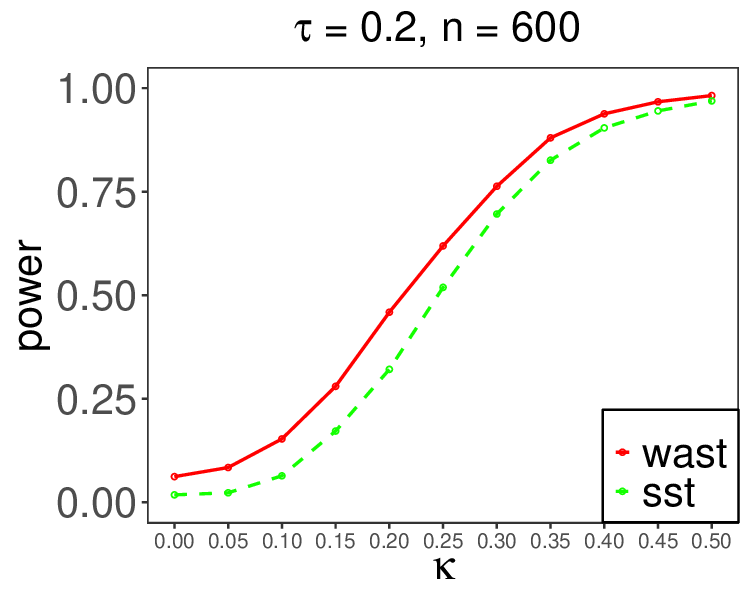}  \\
		\includegraphics[scale=0.3]{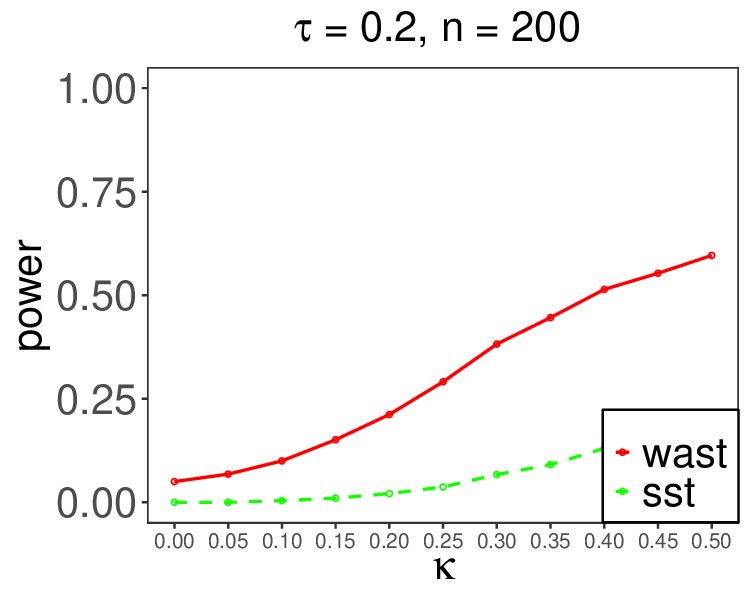}
		\includegraphics[scale=0.3]{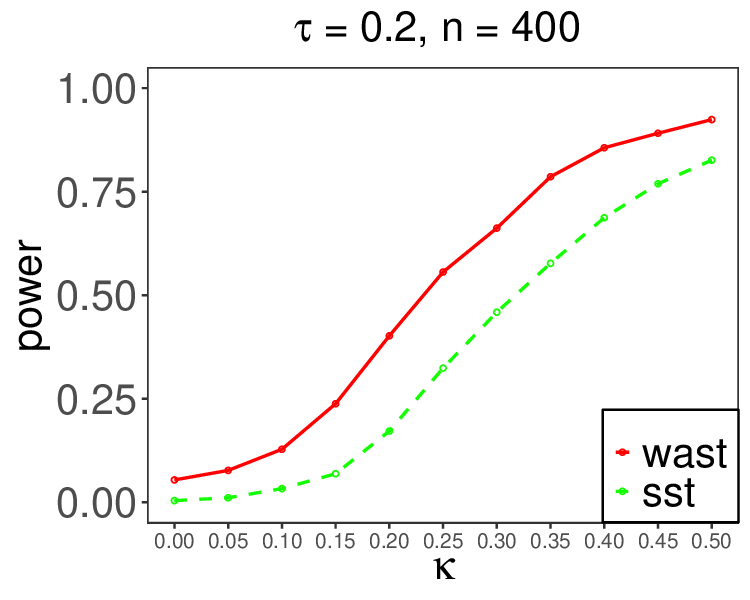}
		\includegraphics[scale=0.3]{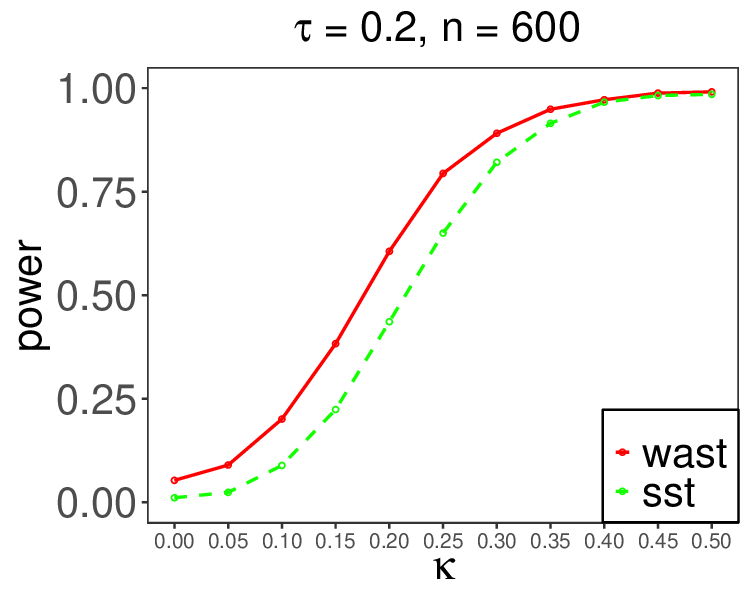} \\
		\includegraphics[scale=0.3]{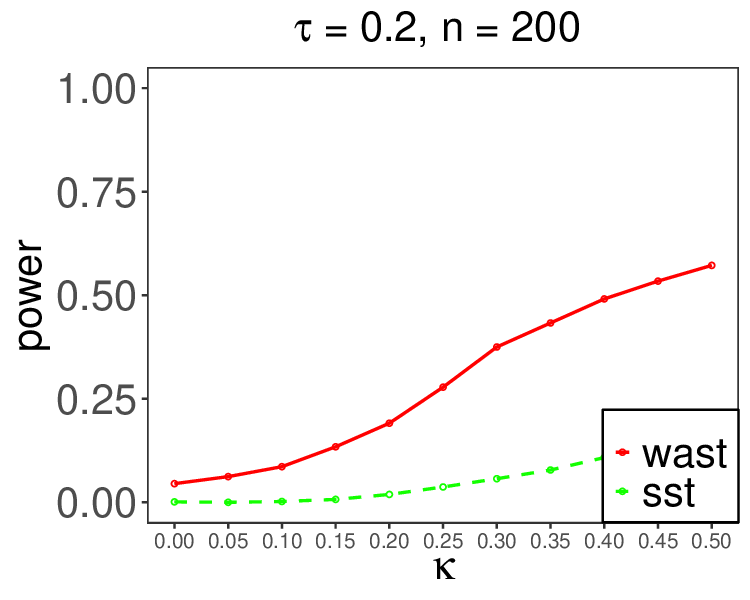}
		\includegraphics[scale=0.3]{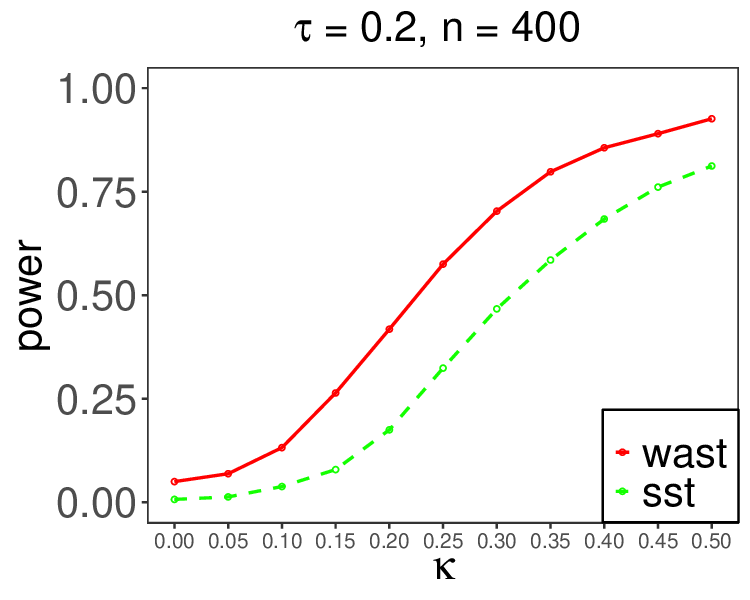}
		\includegraphics[scale=0.3]{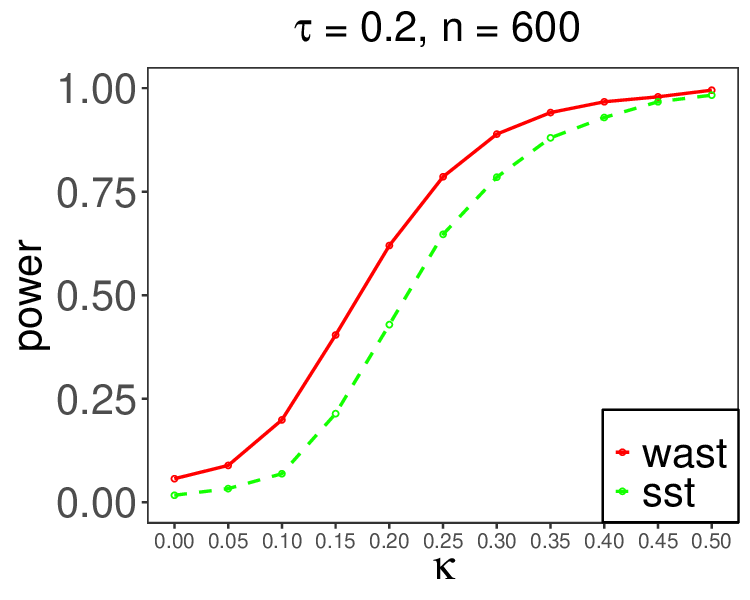}
		\caption{\it Powers of test statistic for quantile regression with $\tau=0.2$ and with large numbers of dense $\bZ$ by the proposed WAST (red solid line) and SST (green dashed line). From top to bottom, each row depicts the powers for $(p,q)=(2,100)$, $(p,q)=(2,500)$, $(p,q)=(6,100)$, $(p,q)=(6,500)$, $(p,q)=(11,100)$, and $(p,q)=(11,500)$.}
		\label{fig_qr20_dense}
	\end{center}
\end{figure}

\begin{figure}[!ht]
	\begin{center}
		\includegraphics[scale=0.3]{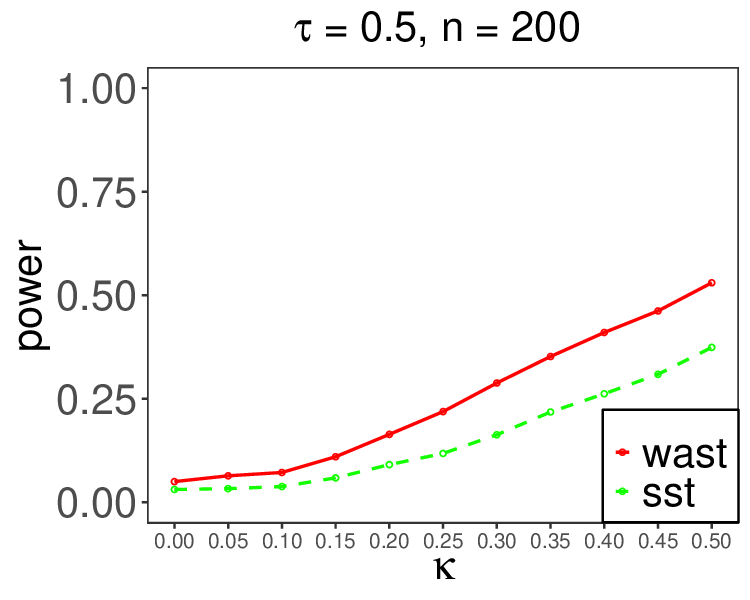}
		\includegraphics[scale=0.3]{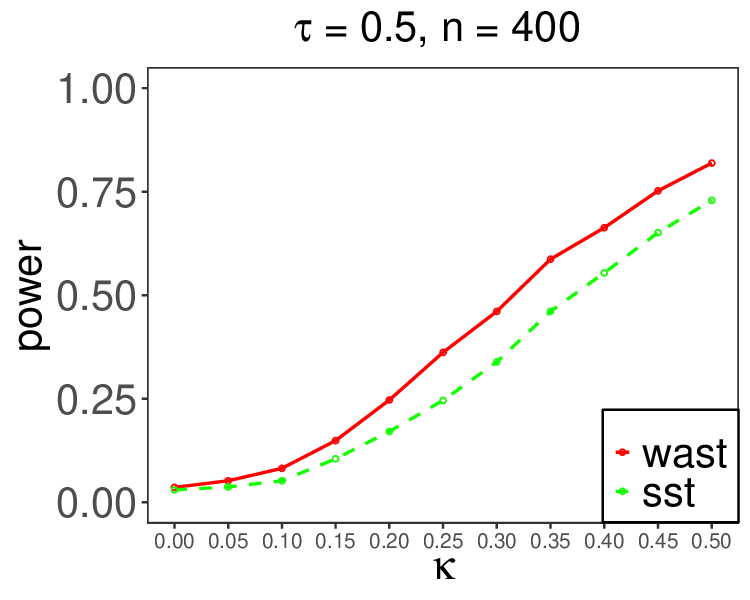}
		\includegraphics[scale=0.3]{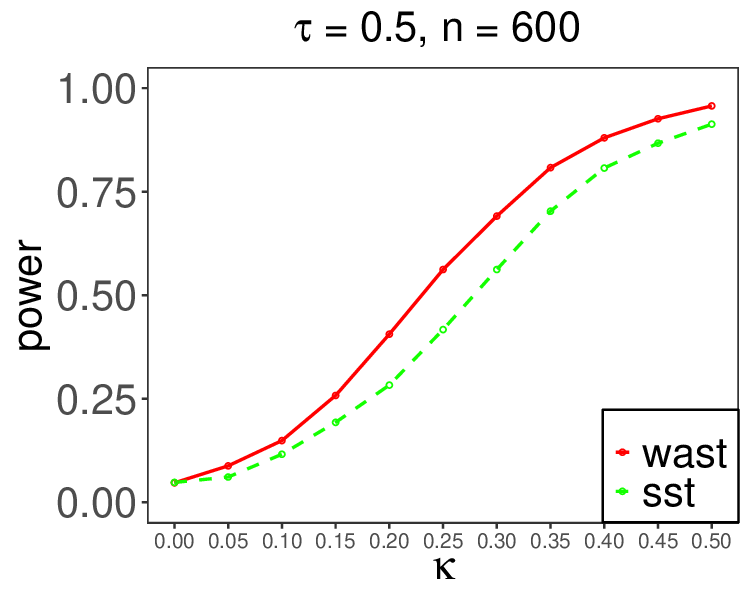}     \\
		\includegraphics[scale=0.3]{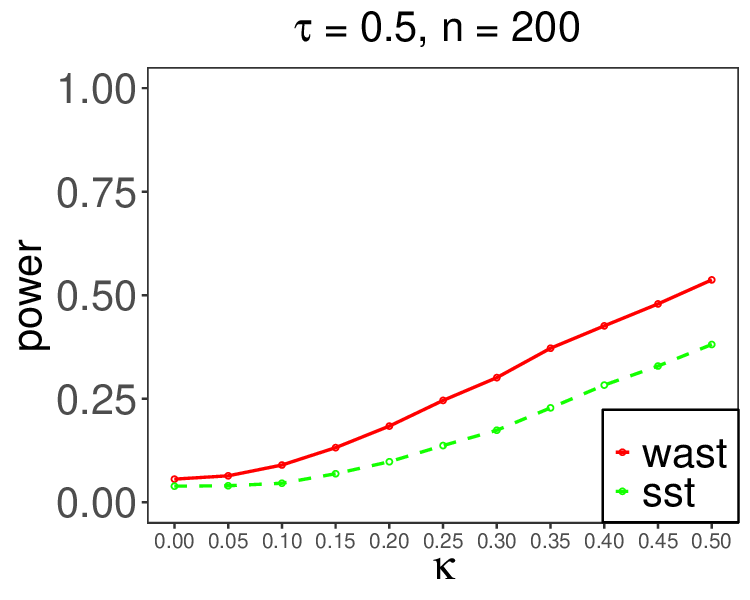}
		\includegraphics[scale=0.3]{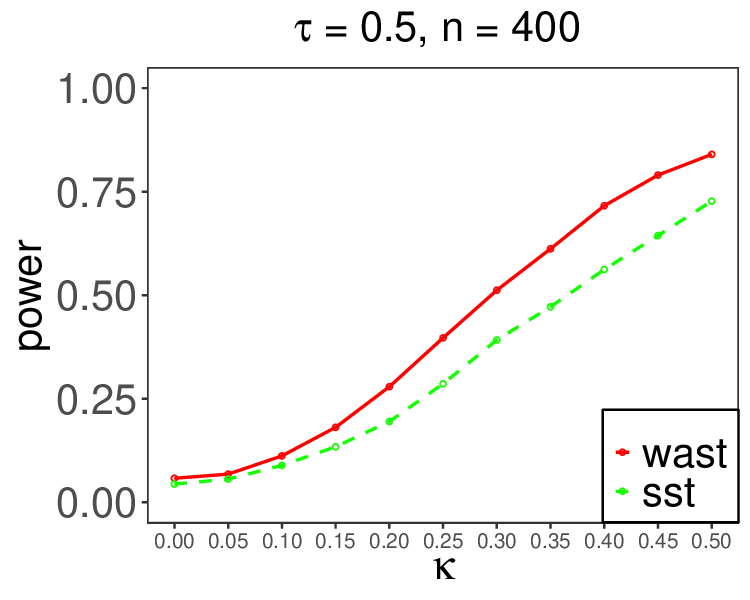}
		\includegraphics[scale=0.3]{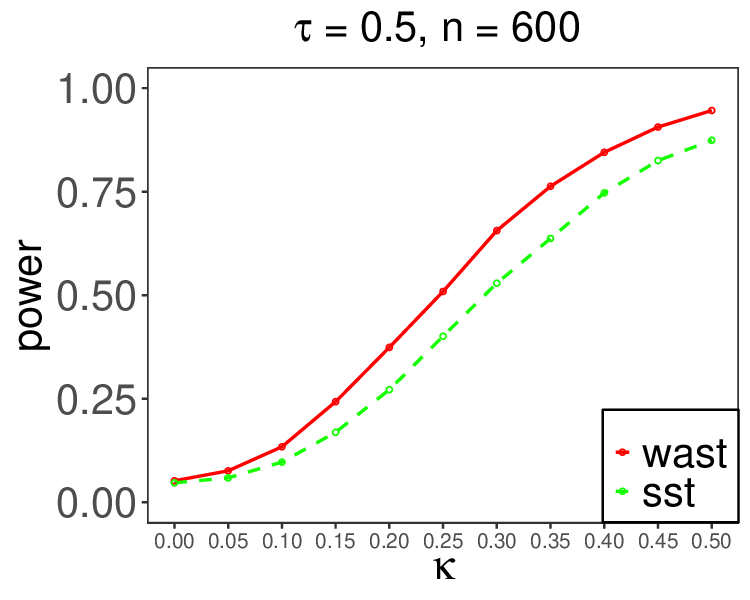}     \\
		\includegraphics[scale=0.3]{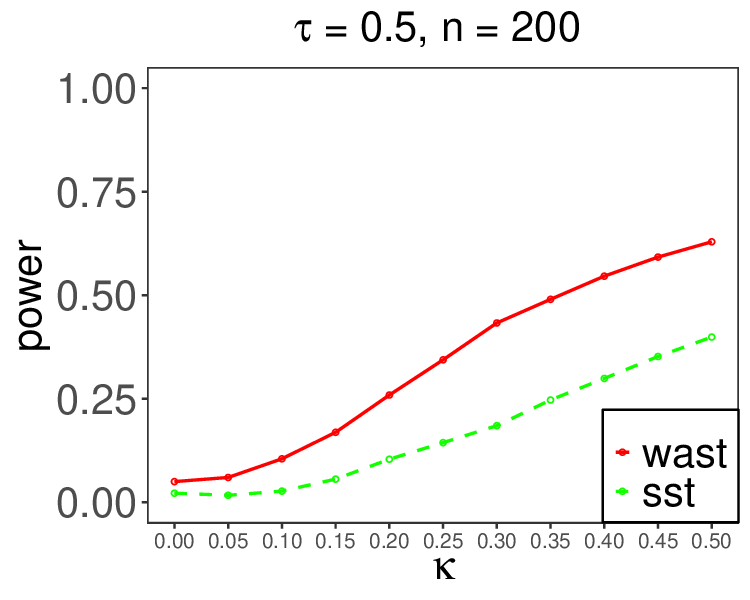}
		\includegraphics[scale=0.3]{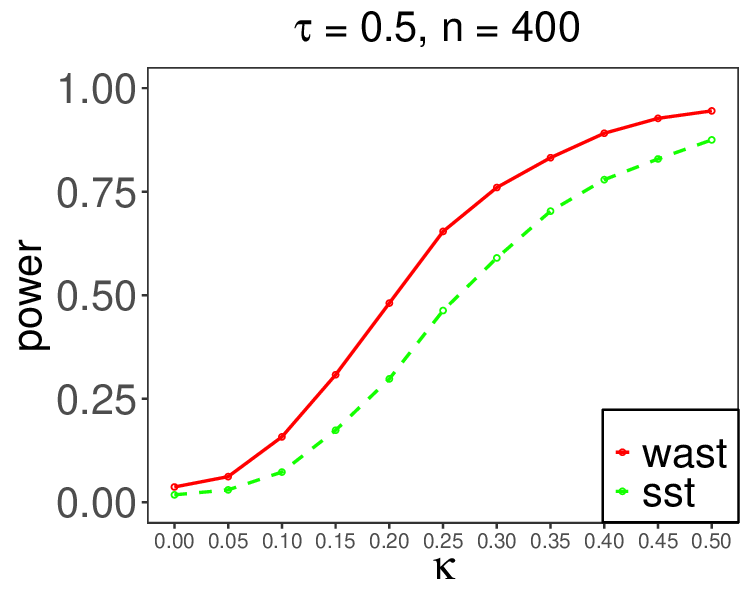}
		\includegraphics[scale=0.3]{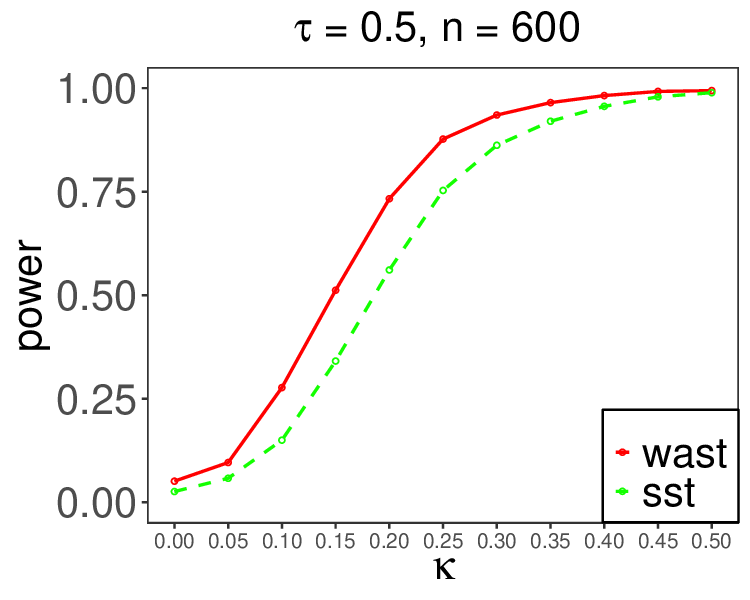}     \\
		\includegraphics[scale=0.3]{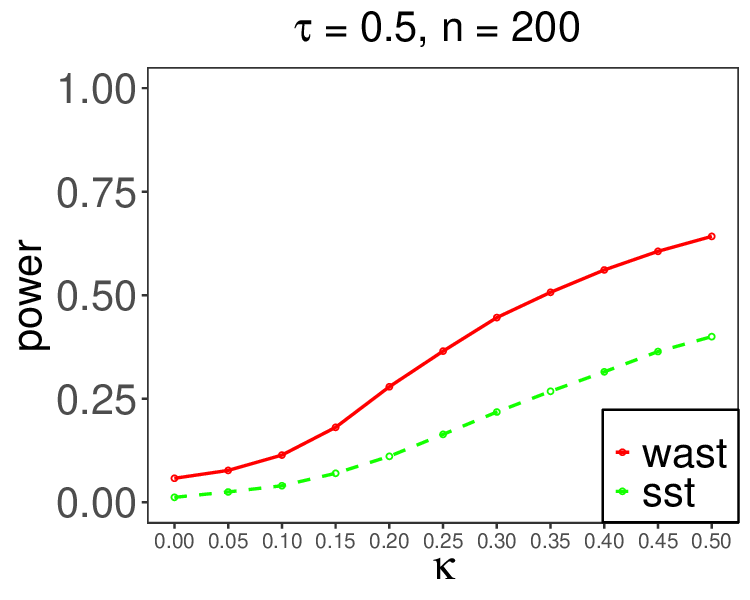}
		\includegraphics[scale=0.3]{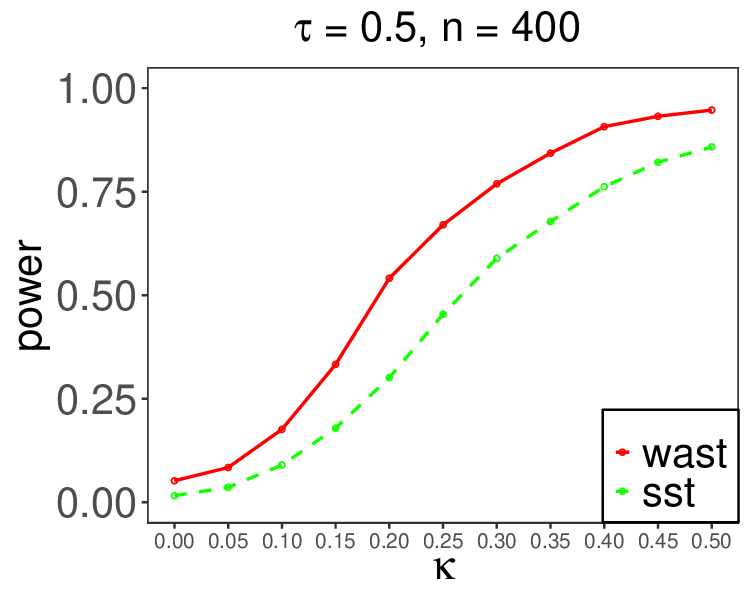}
		\includegraphics[scale=0.3]{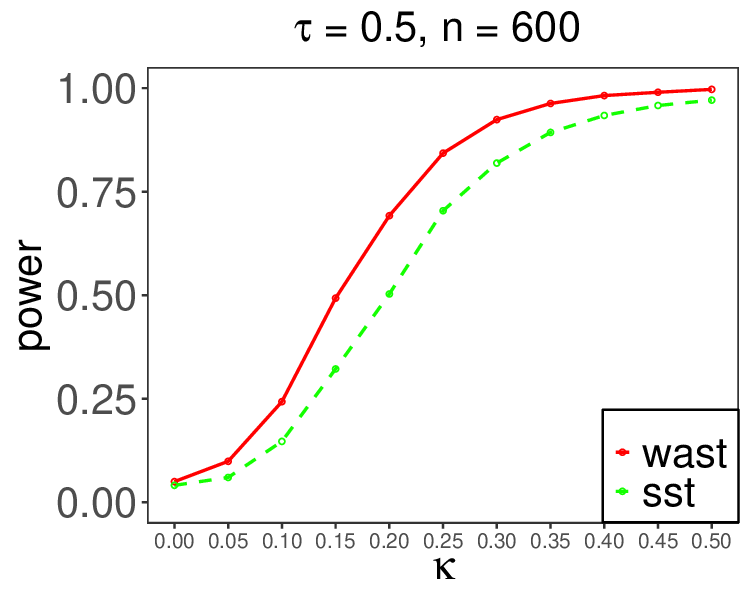}     \\
		\includegraphics[scale=0.3]{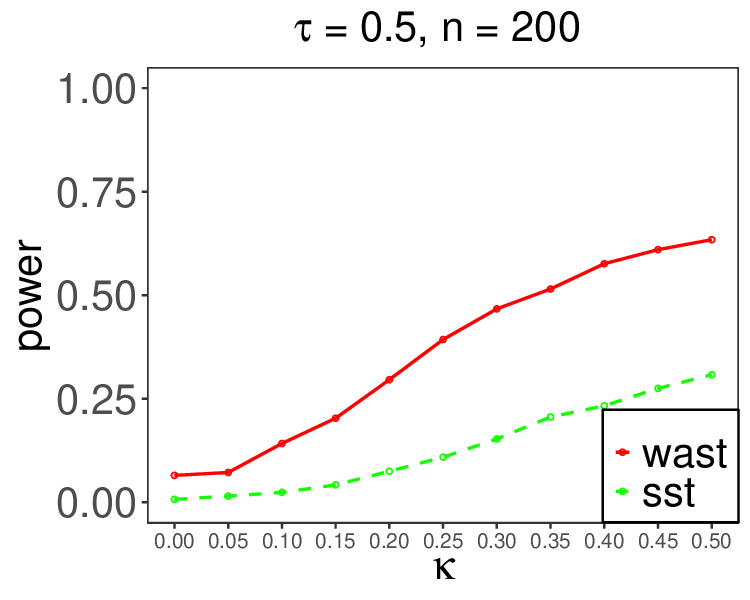}
		\includegraphics[scale=0.3]{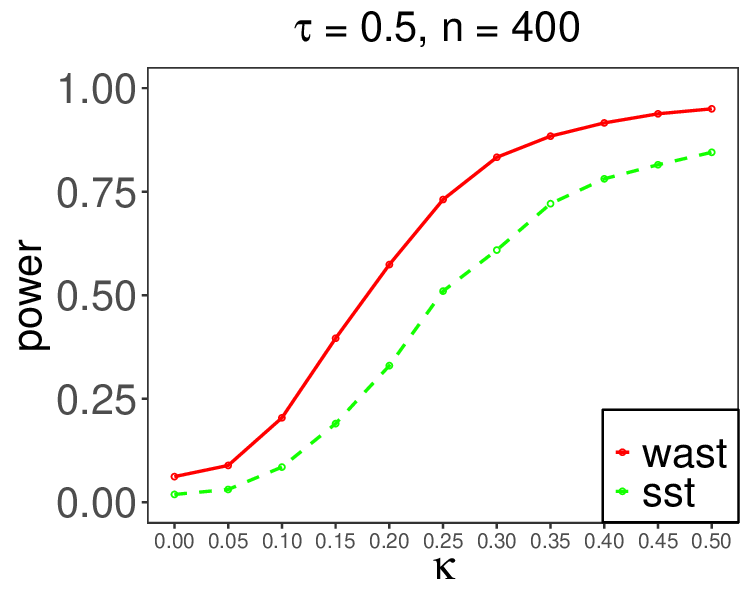}
		\includegraphics[scale=0.3]{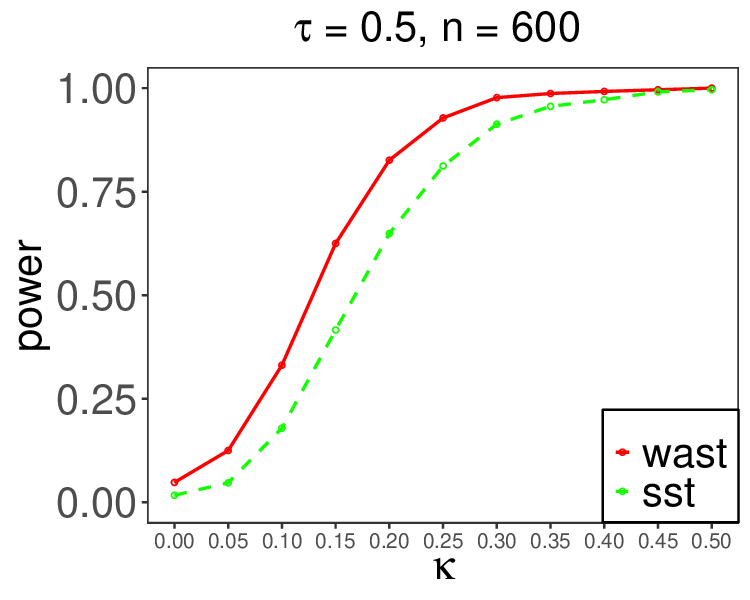}    \\
		\includegraphics[scale=0.3]{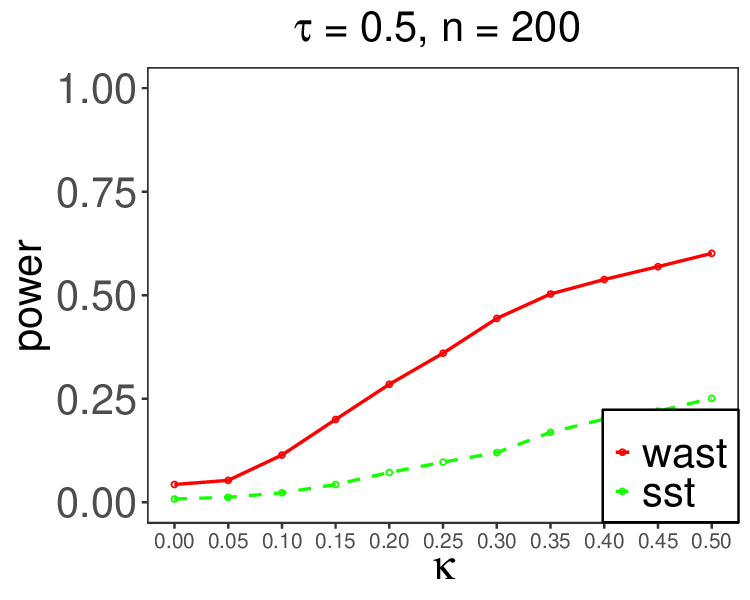}
		\includegraphics[scale=0.3]{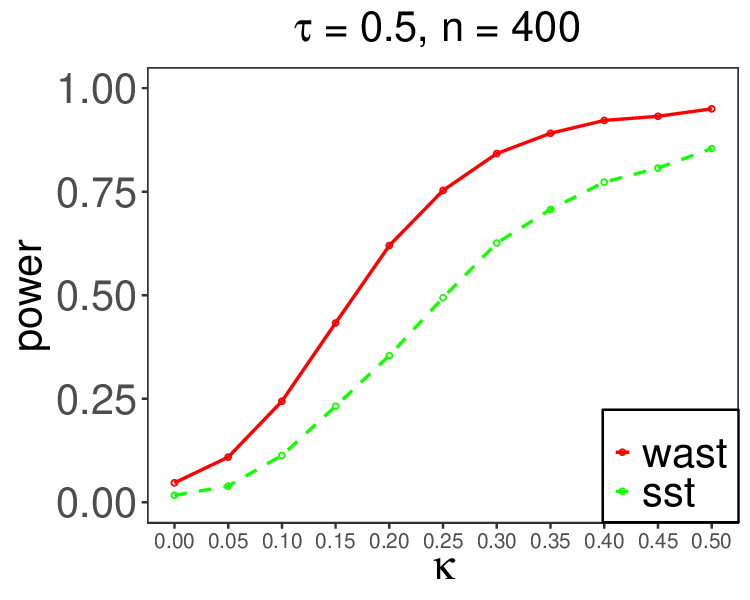}
		\includegraphics[scale=0.3]{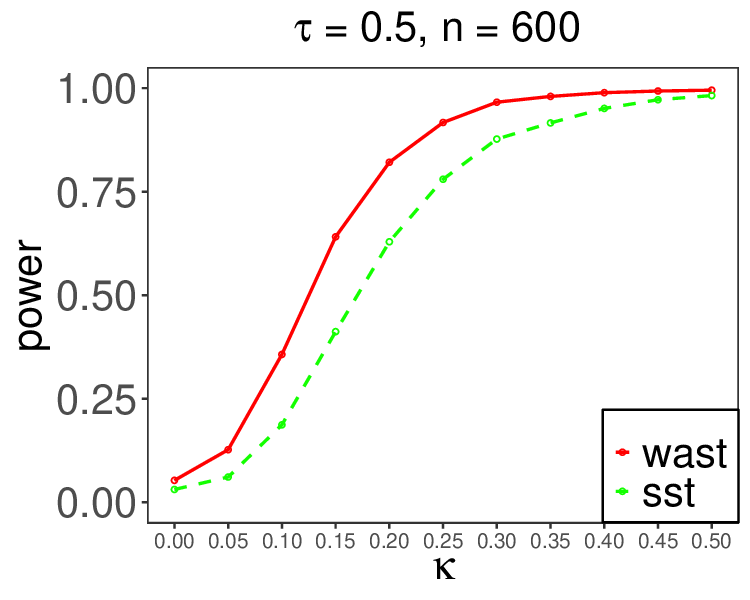}
		\caption{\it Powers of test statistic for quantile regression with $\tau=0.5$ by the proposed WAST (red solid line) and SST (green dashed line). From top to bottom, each row depicts the powers for $(p,q)=(2,100)$, $(p,q)=(2,500)$, $(p,q)=(6,100)$, $(p,q)=(6,500)$, $(p,q)=(11,100)$, and $(p,q)=(11,500)$.}
		\label{fig_qr50_dense}
	\end{center}
\end{figure}

\begin{figure}[!ht]
	\begin{center}
		\includegraphics[scale=0.3]{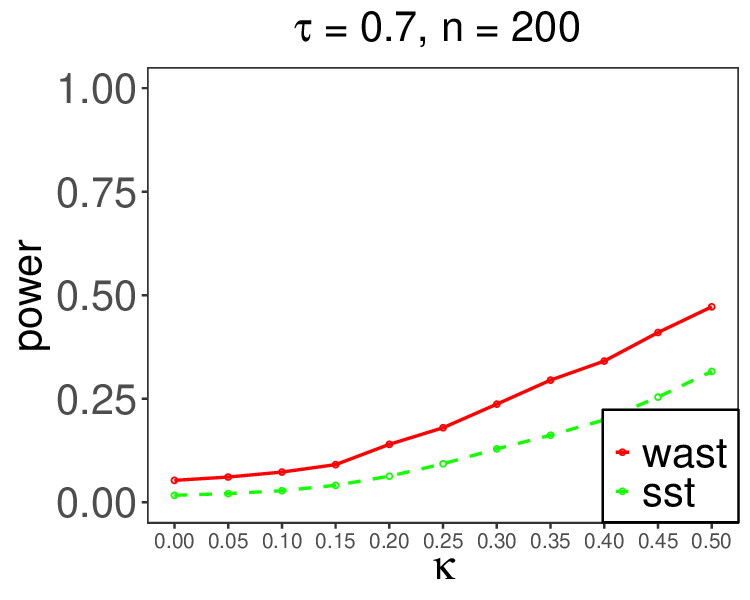}
		\includegraphics[scale=0.3]{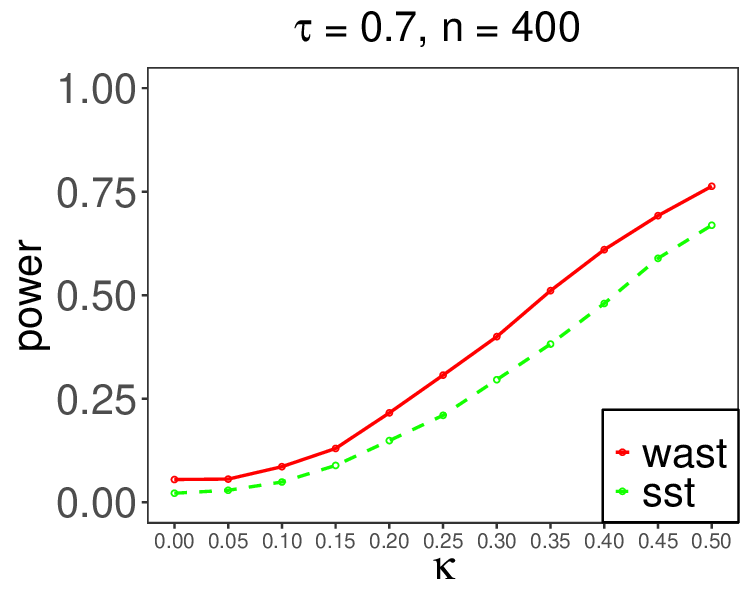}
		\includegraphics[scale=0.3]{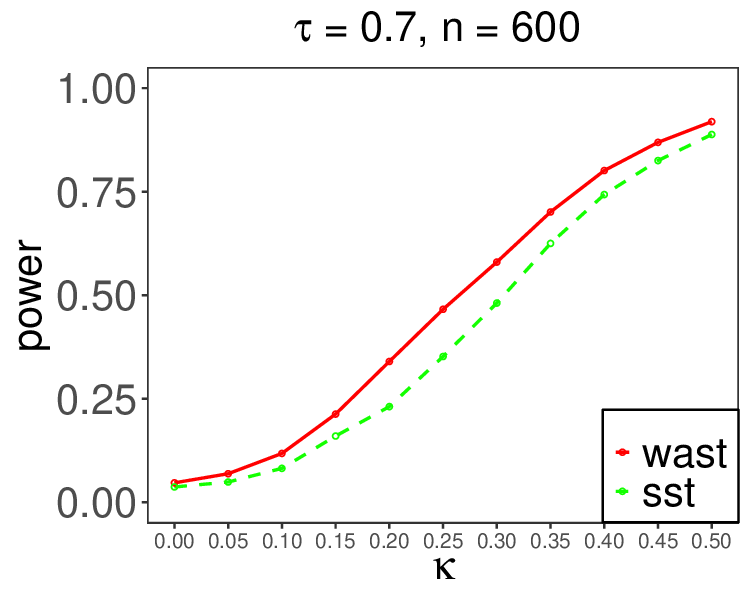}     \\
		\includegraphics[scale=0.3]{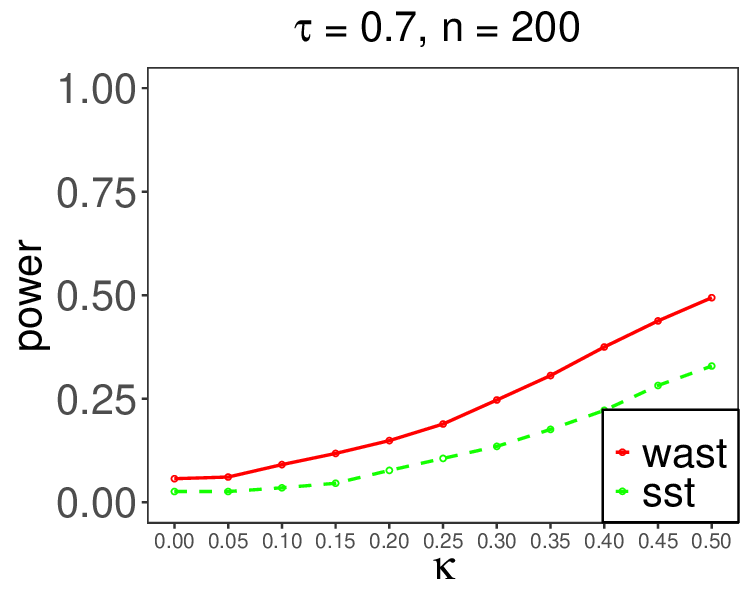}
		\includegraphics[scale=0.3]{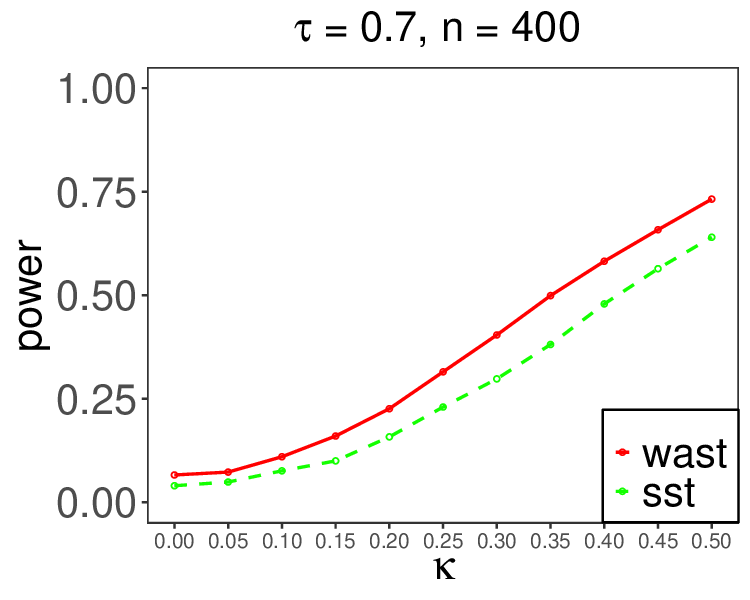}
		\includegraphics[scale=0.3]{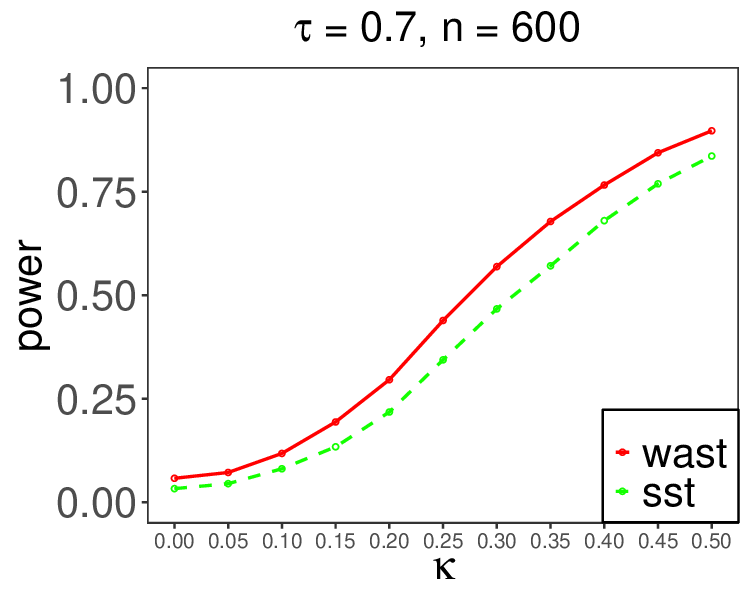}     \\
		\includegraphics[scale=0.3]{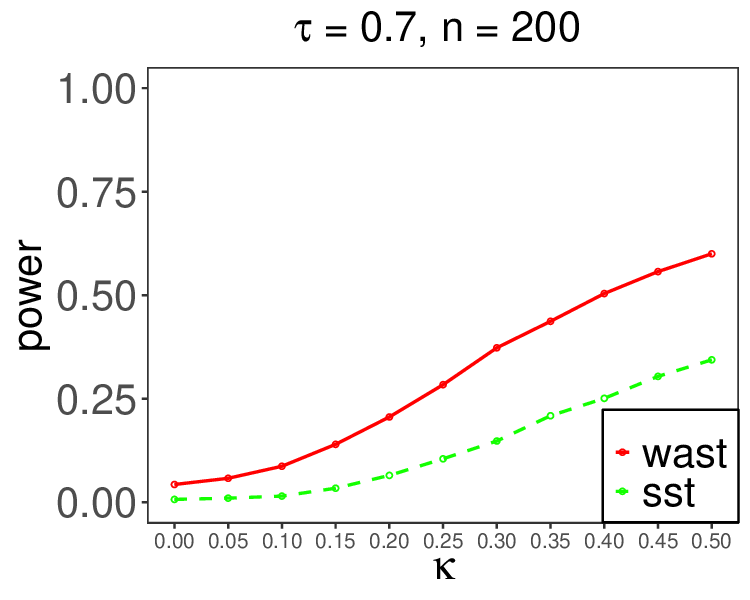}
		\includegraphics[scale=0.3]{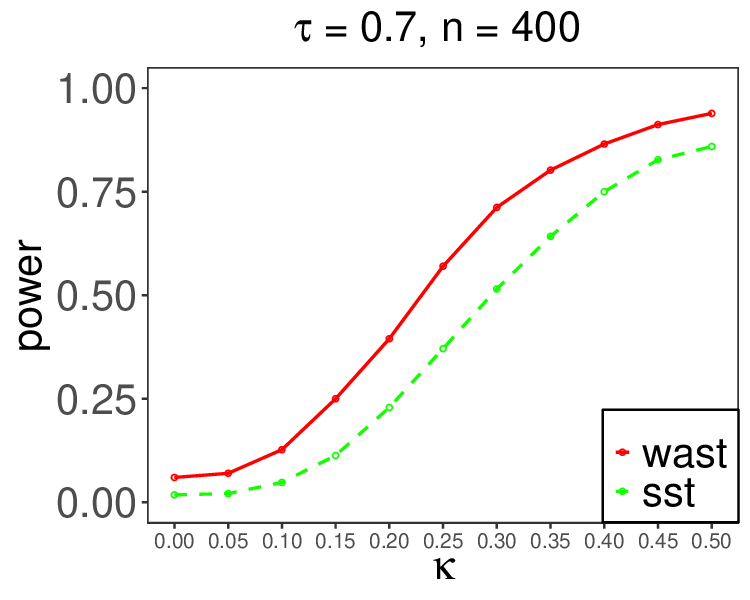}
		\includegraphics[scale=0.3]{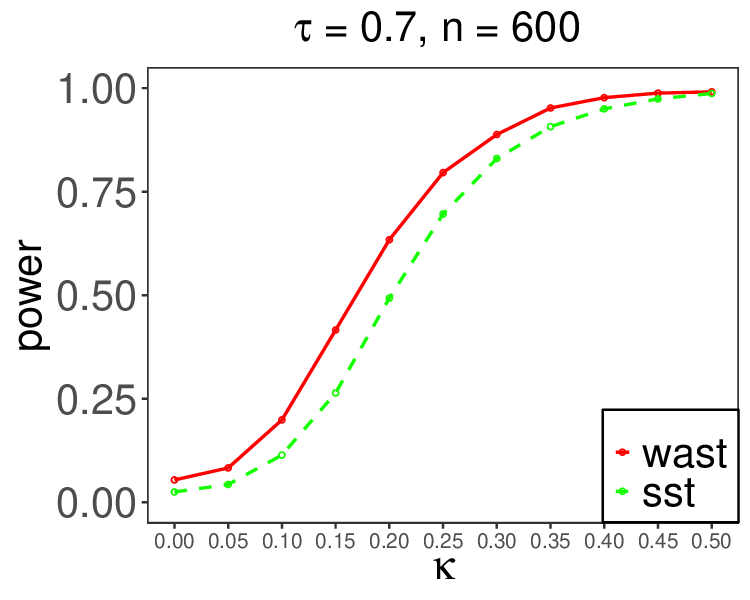}     \\
		\includegraphics[scale=0.3]{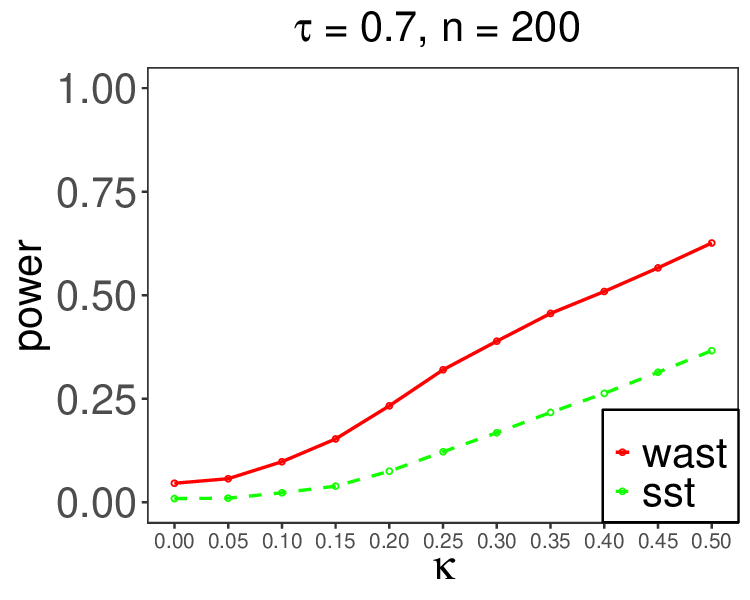}
		\includegraphics[scale=0.3]{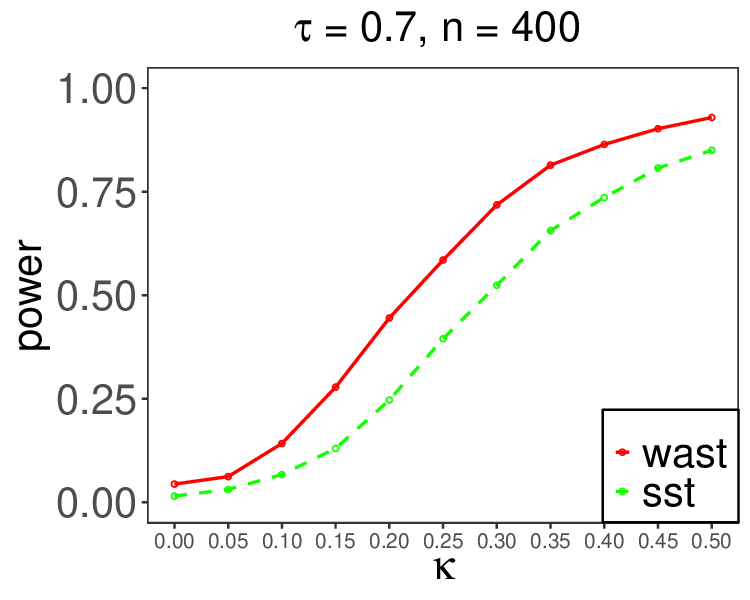}
		\includegraphics[scale=0.3]{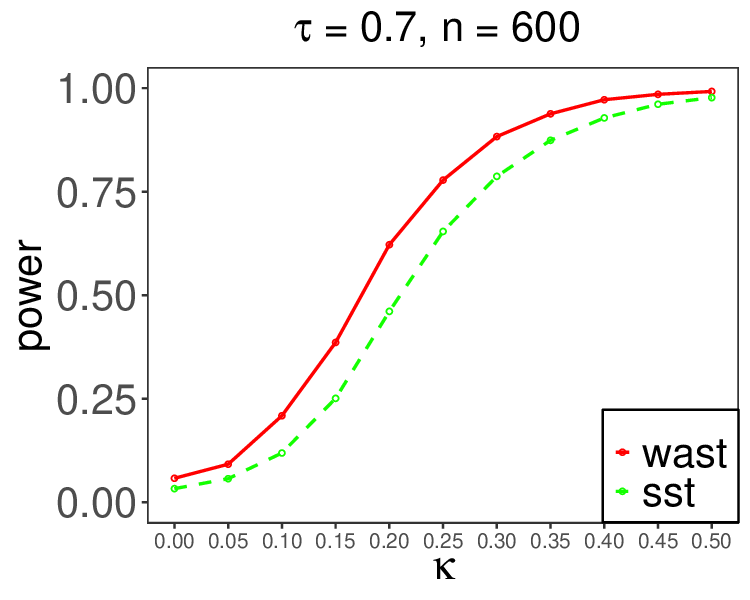}     \\
		\includegraphics[scale=0.3]{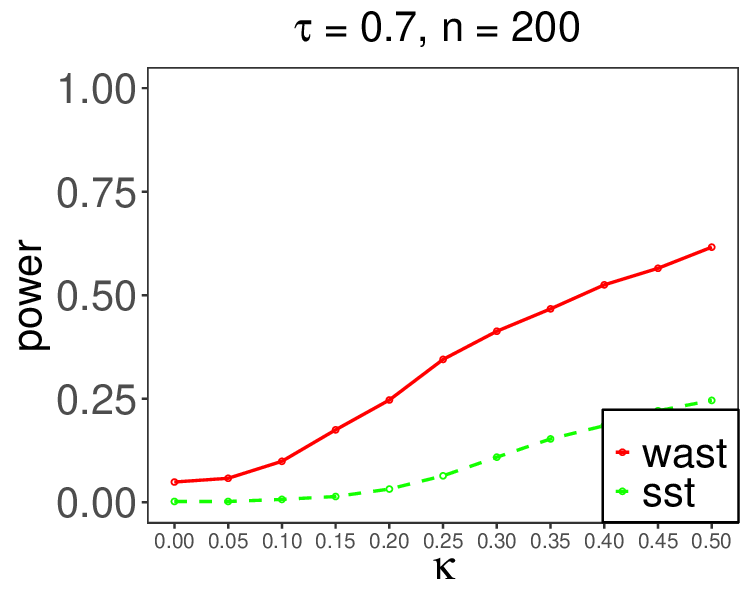}
		\includegraphics[scale=0.3]{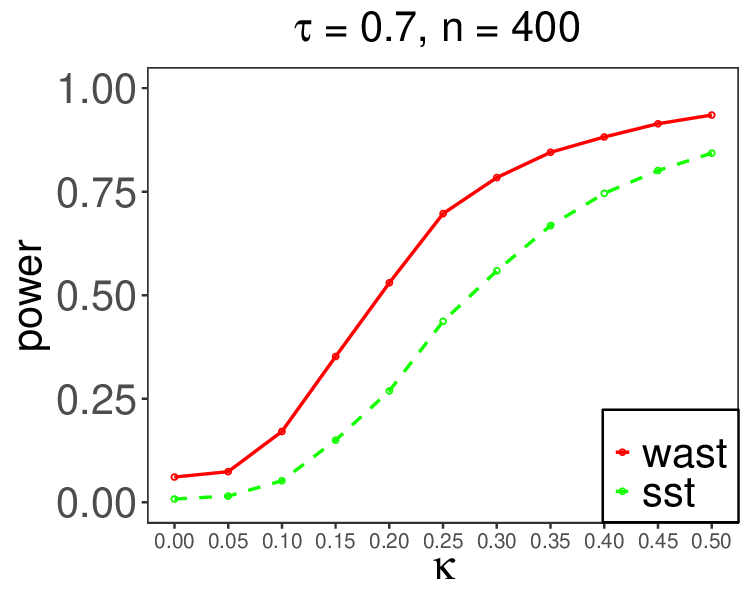}
		\includegraphics[scale=0.3]{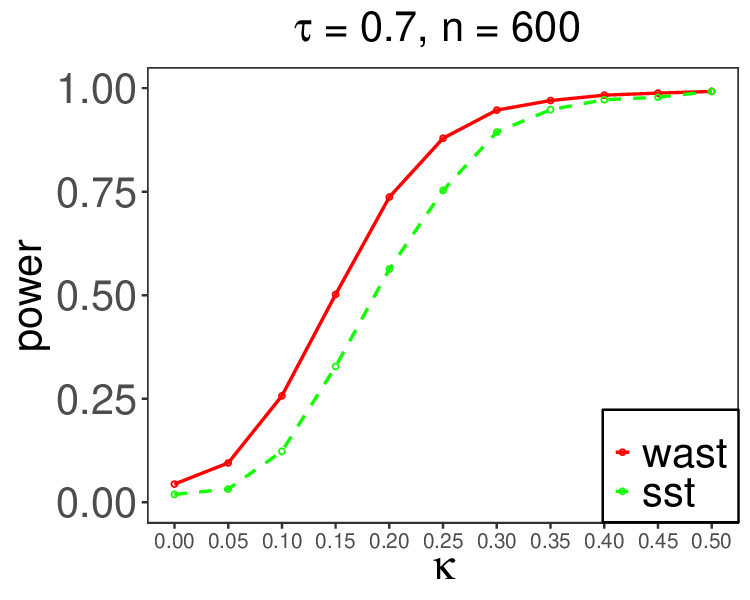}    \\
		\includegraphics[scale=0.3]{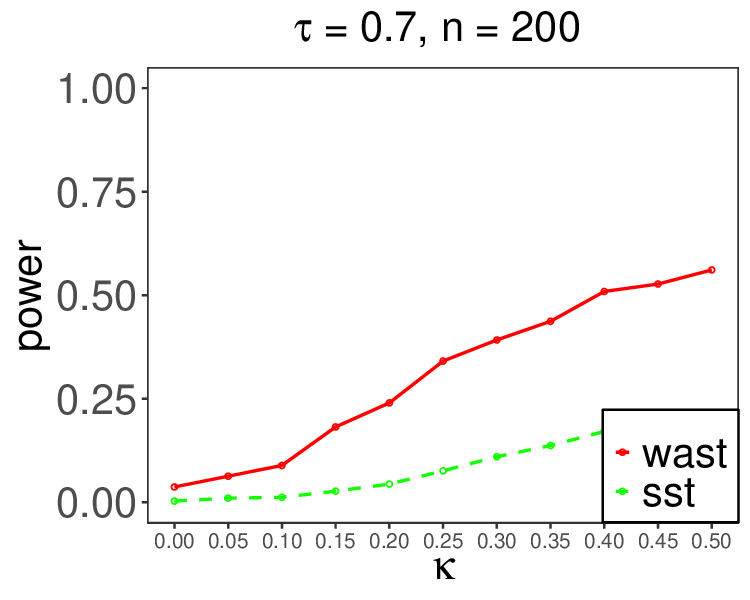}
		\includegraphics[scale=0.3]{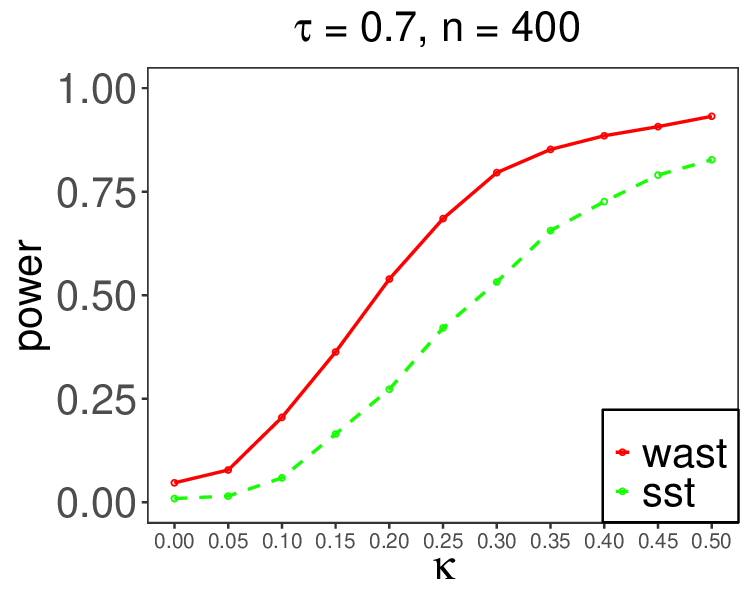}
		\includegraphics[scale=0.3]{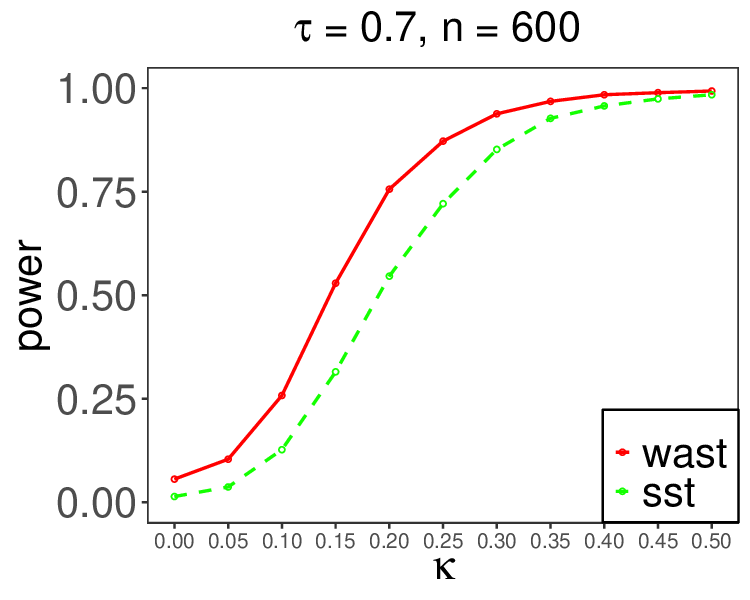}
		\caption{\it Powers of test statistic for quantile regression with $\tau=0.7$ and with large numbers of dense $\bZ$ by the proposed WAST (red solid line) and SST (green dashed line). From top to bottom, each row depicts the powers for $(p,q)=(2,100)$, $(p,q)=(2,500)$, $(p,q)=(6,100)$, $(p,q)=(6,500)$, $(p,q)=(11,100)$, and $(p,q)=(11,500)$.}
		\label{fig_qr70_dense}
	\end{center}
\end{figure}

\begin{figure}[!ht]
	\begin{center}
		\includegraphics[scale=0.3]{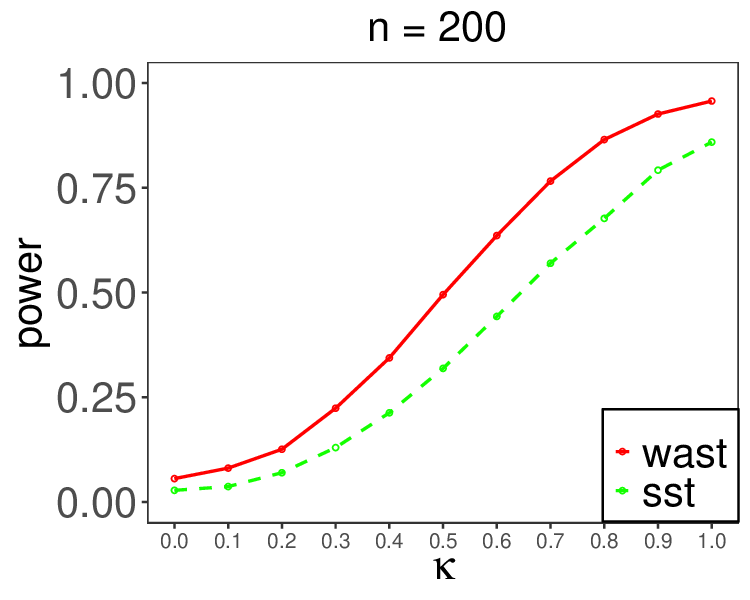}
		\includegraphics[scale=0.3]{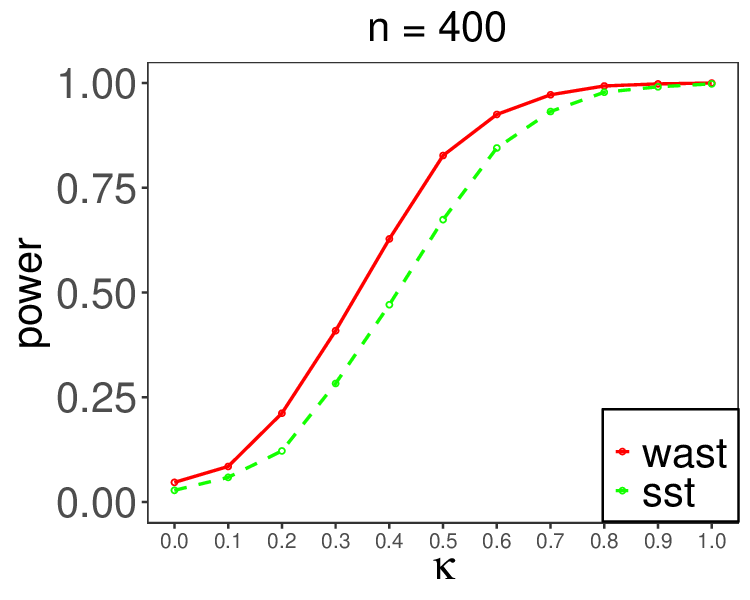}
		\includegraphics[scale=0.3]{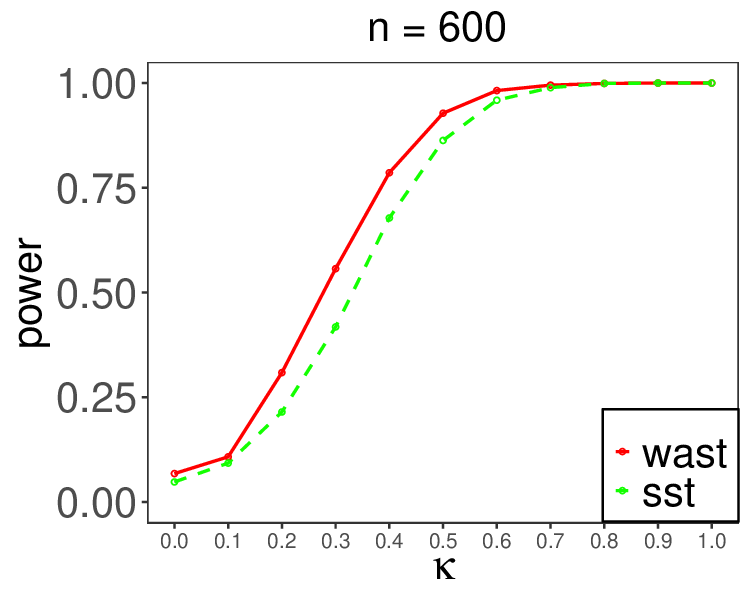}      \\
		\includegraphics[scale=0.3]{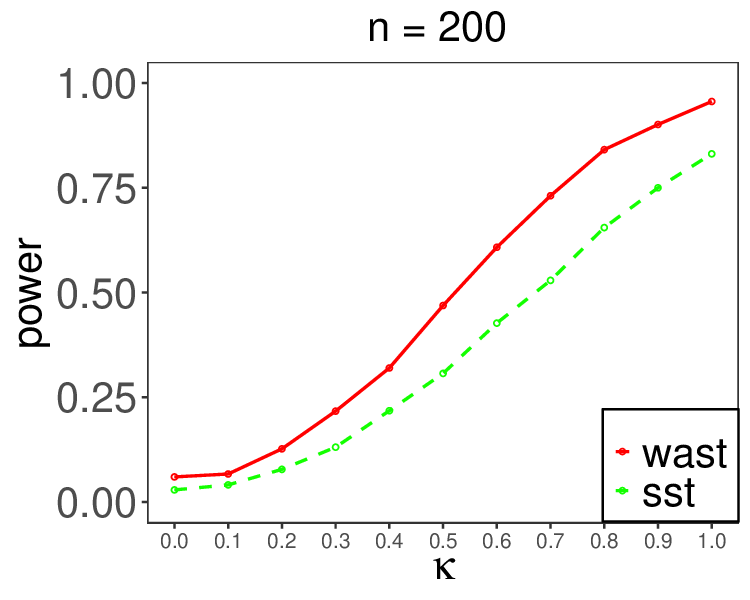}
		\includegraphics[scale=0.3]{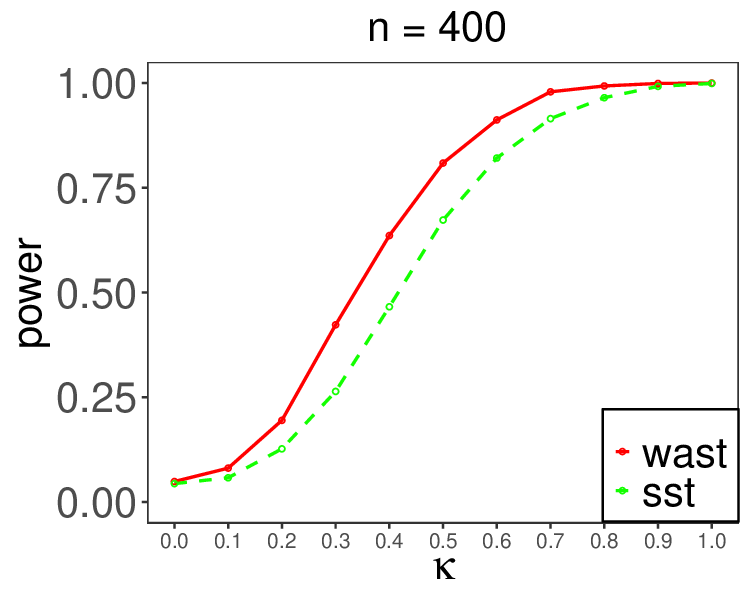}
		\includegraphics[scale=0.3]{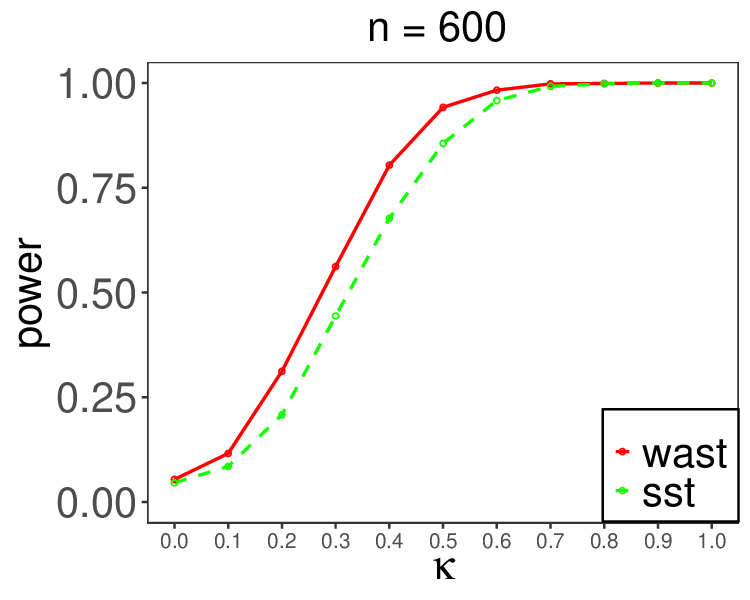}      \\
		\includegraphics[scale=0.3]{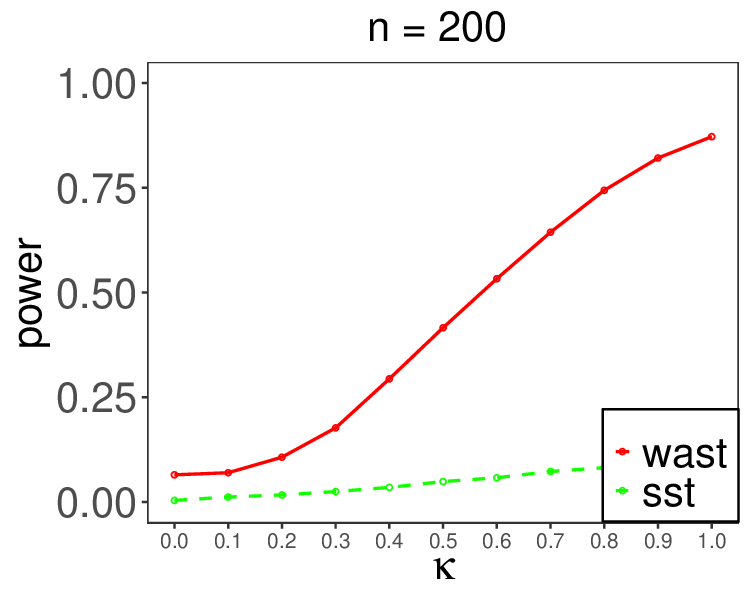}
		\includegraphics[scale=0.3]{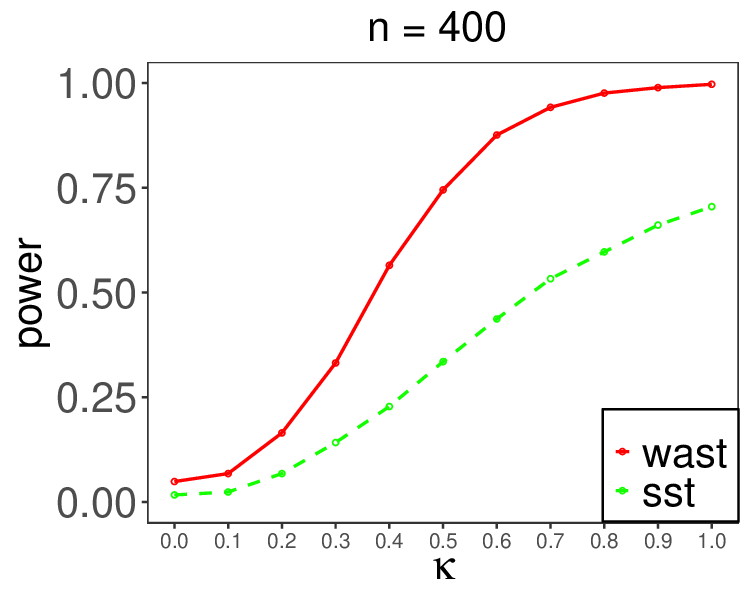}
		\includegraphics[scale=0.3]{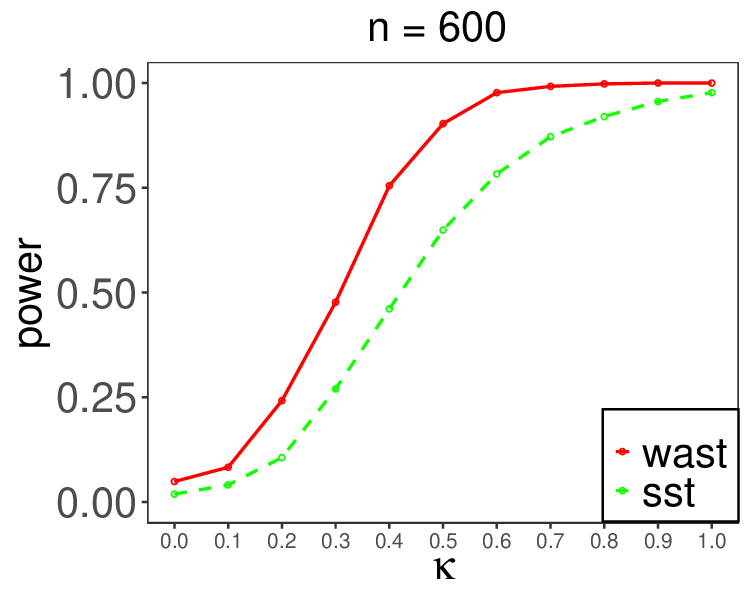}      \\
		\includegraphics[scale=0.3]{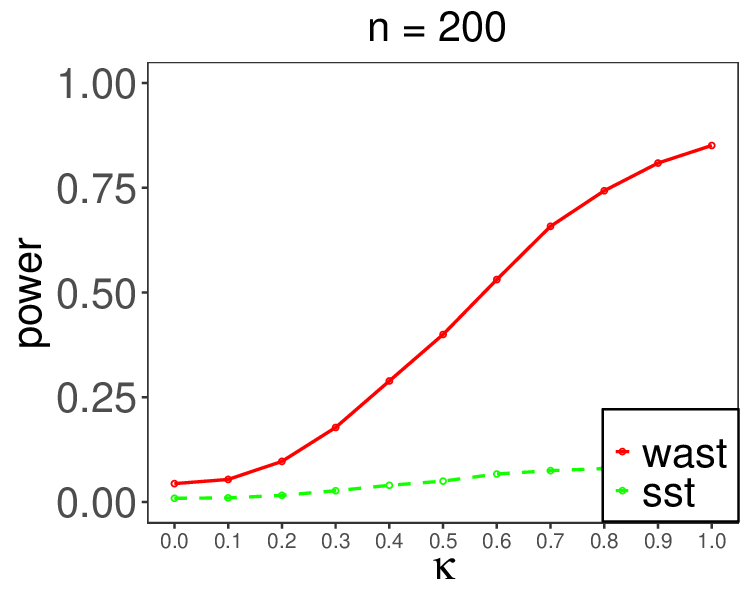}
		\includegraphics[scale=0.3]{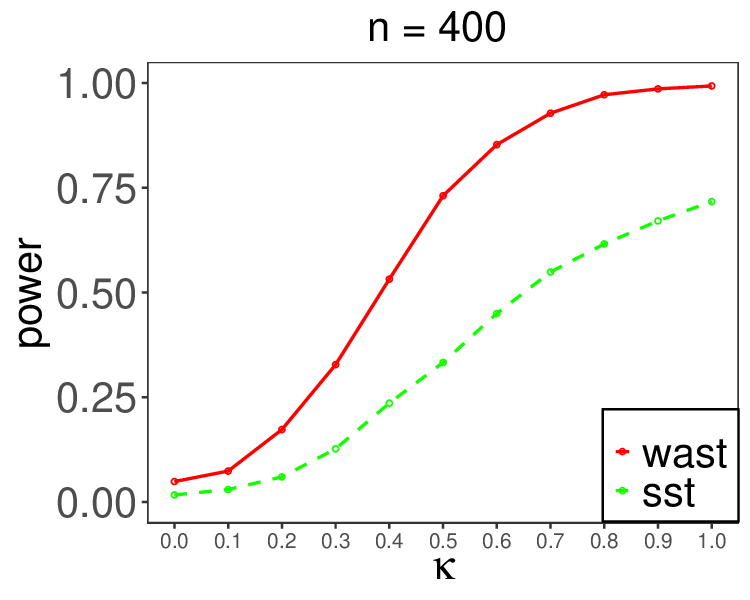}
		\includegraphics[scale=0.3]{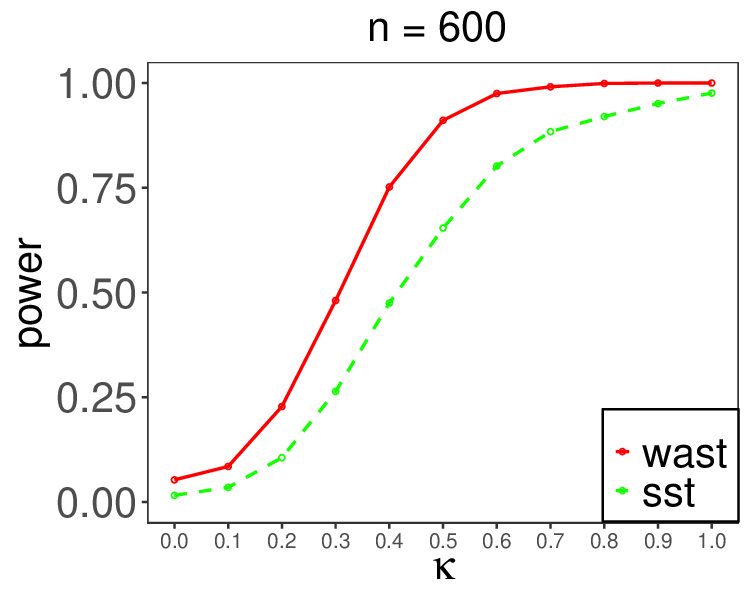}      \\
		\includegraphics[scale=0.3]{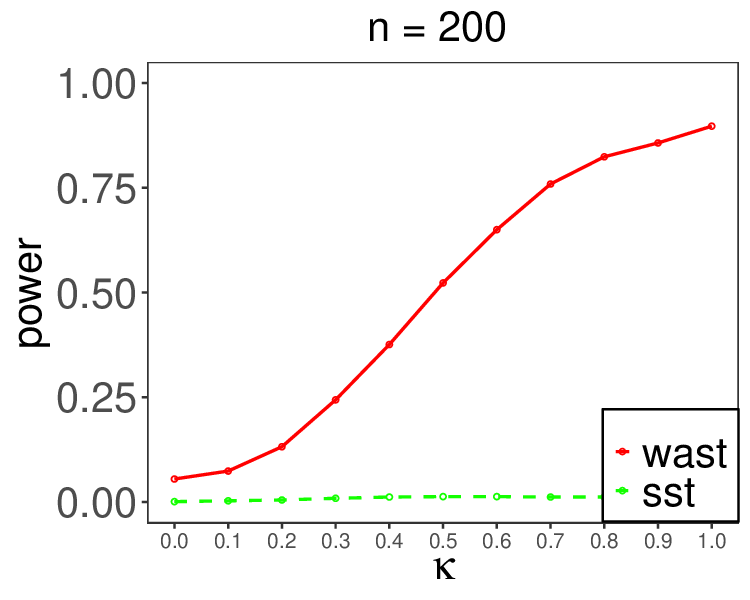}
		\includegraphics[scale=0.3]{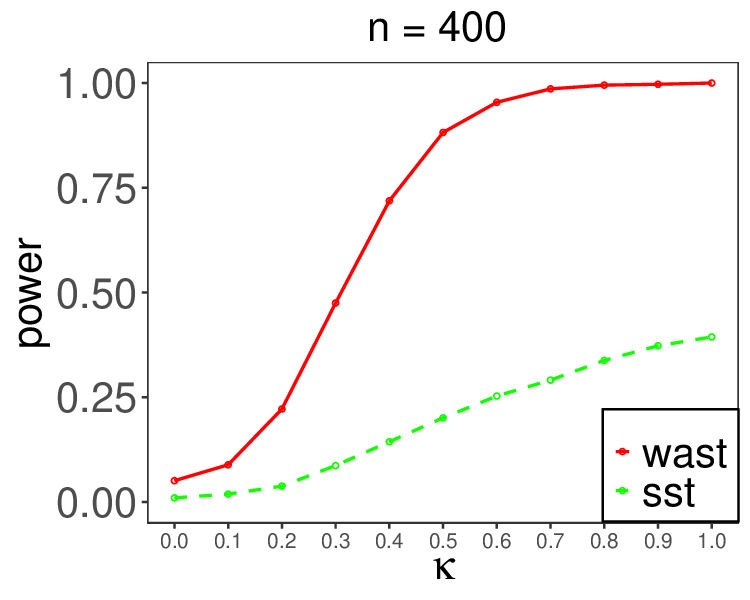}
		\includegraphics[scale=0.3]{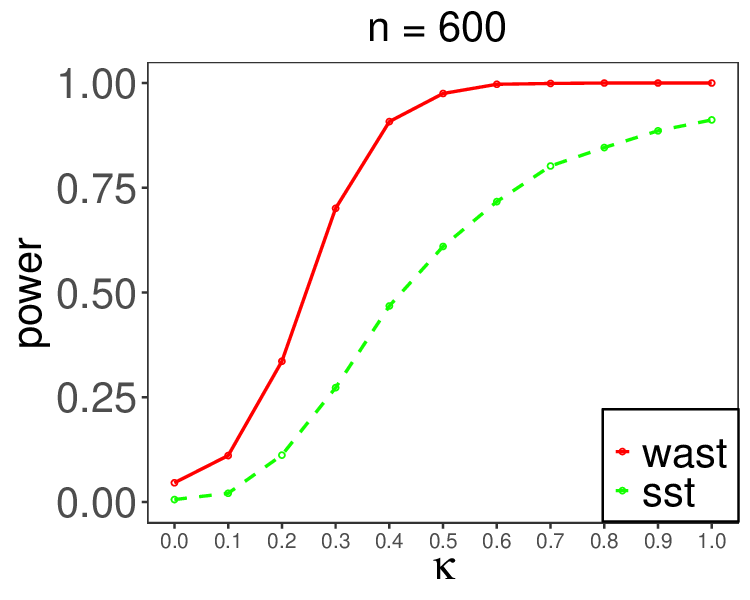}     \\
		\includegraphics[scale=0.3]{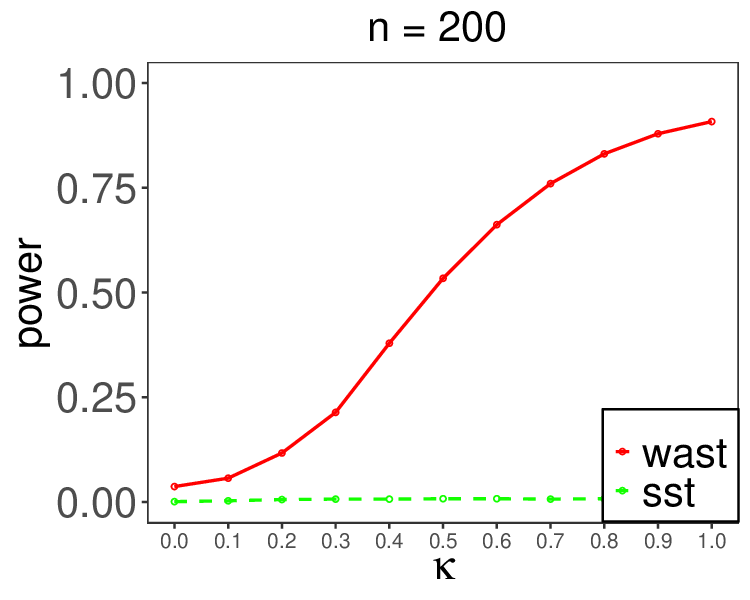}
		\includegraphics[scale=0.3]{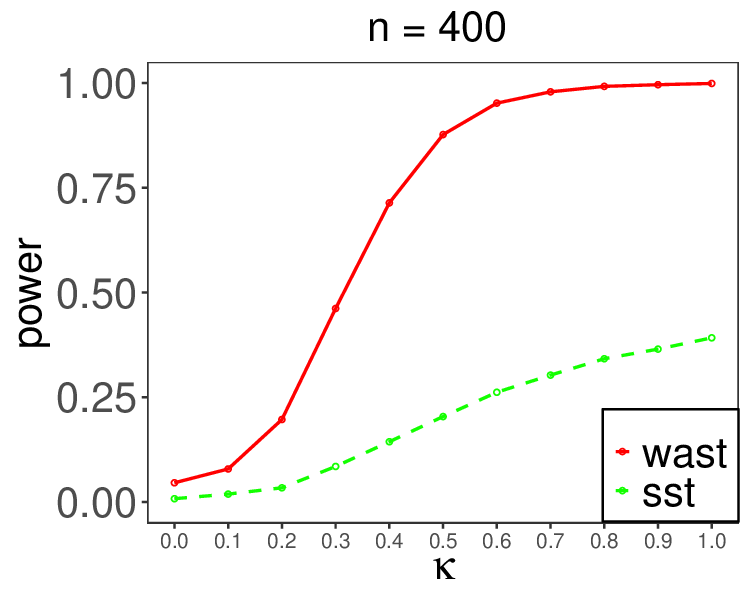}
		\includegraphics[scale=0.3]{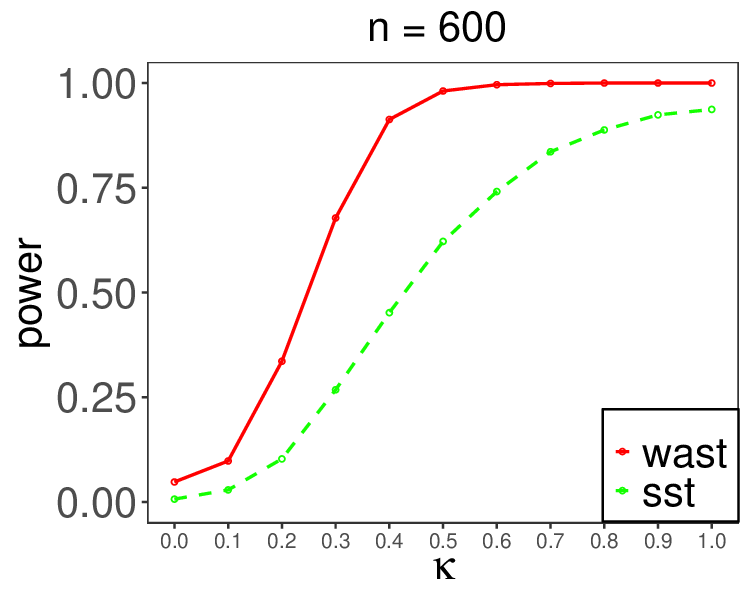}
		\caption{\it Powers of test statistic for semiparametric model (B1+P1) with large numbers of dense $\bZ$ by the proposed WAST (red solid line) and SST (green dashed line). From top to bottom, each row depicts the powers for $(p,q)=(2,100)$, $(p,q)=(2,500)$, $(p,q)=(6,100)$, $(p,q)=(6,500)$, $(p,q)=(11,100)$, and $(p,q)=(11,500)$.}
		\label{fig_semiparam1_dense}
	\end{center}
\end{figure}

\begin{figure}[!ht]
	\begin{center}
		\includegraphics[scale=0.3]{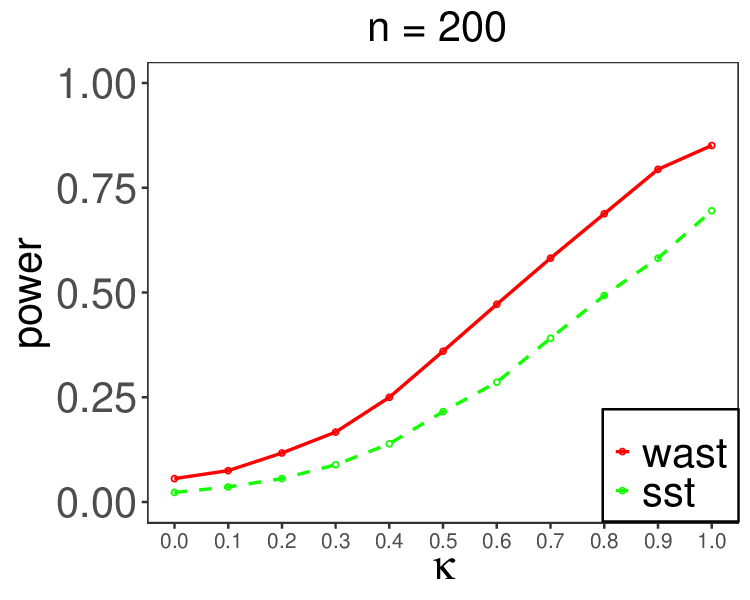}
		\includegraphics[scale=0.3]{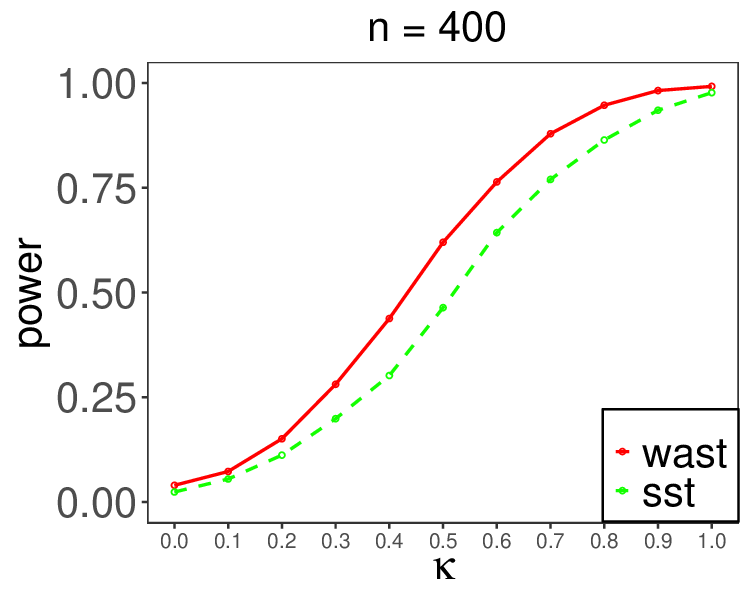}
		\includegraphics[scale=0.3]{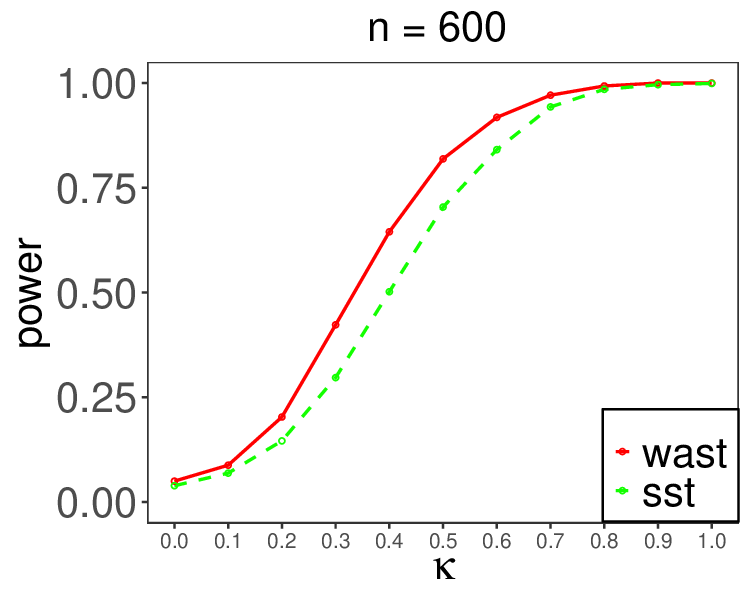}     \\
		\includegraphics[scale=0.3]{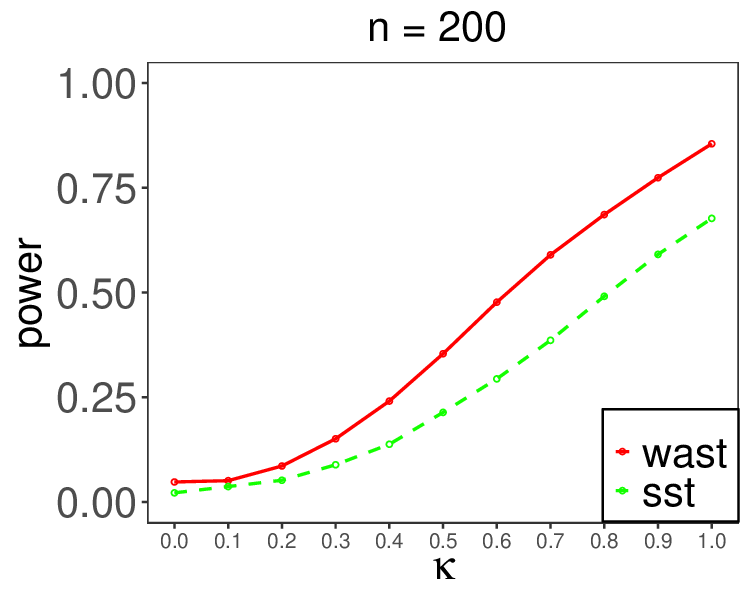}
		\includegraphics[scale=0.3]{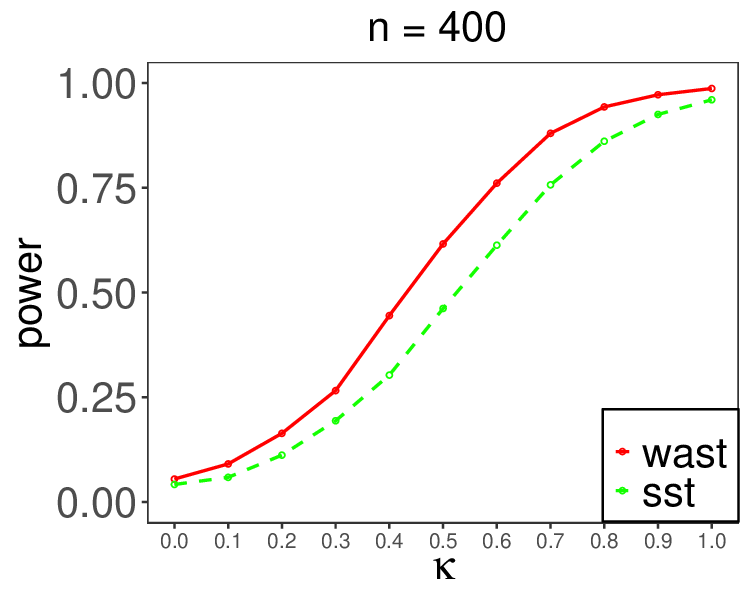}
		\includegraphics[scale=0.3]{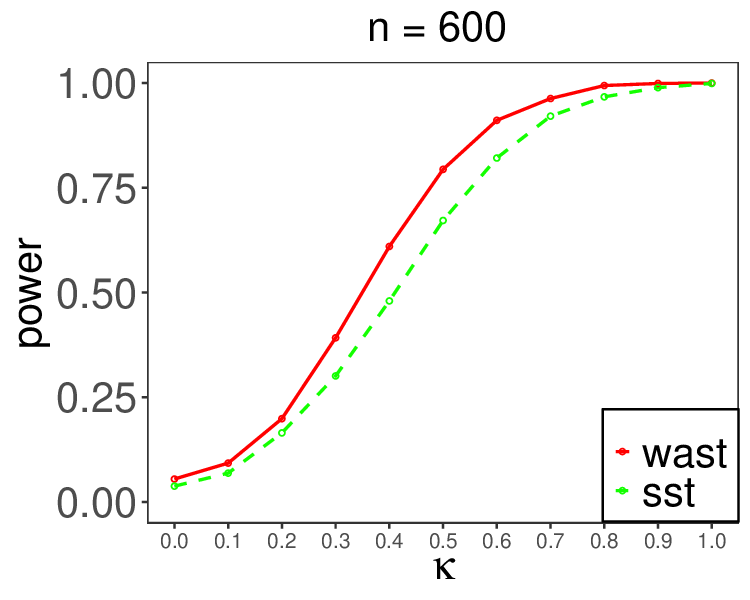}     \\
		\includegraphics[scale=0.3]{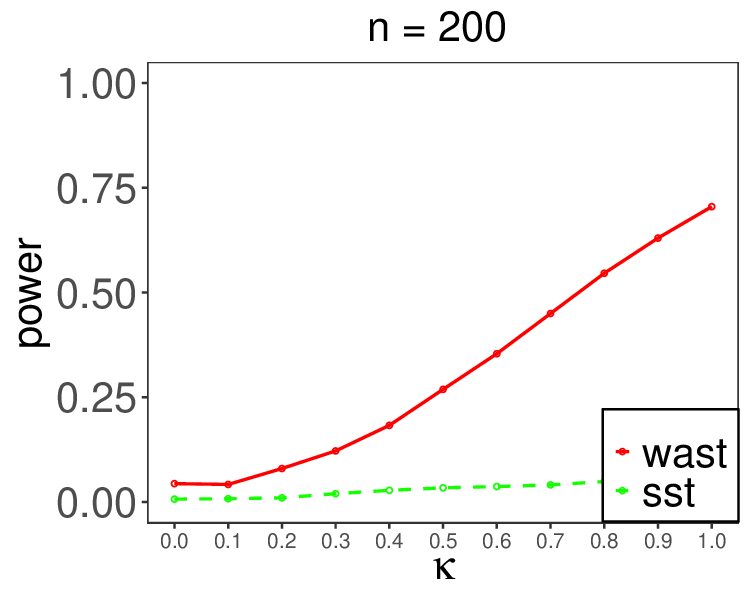}
		\includegraphics[scale=0.3]{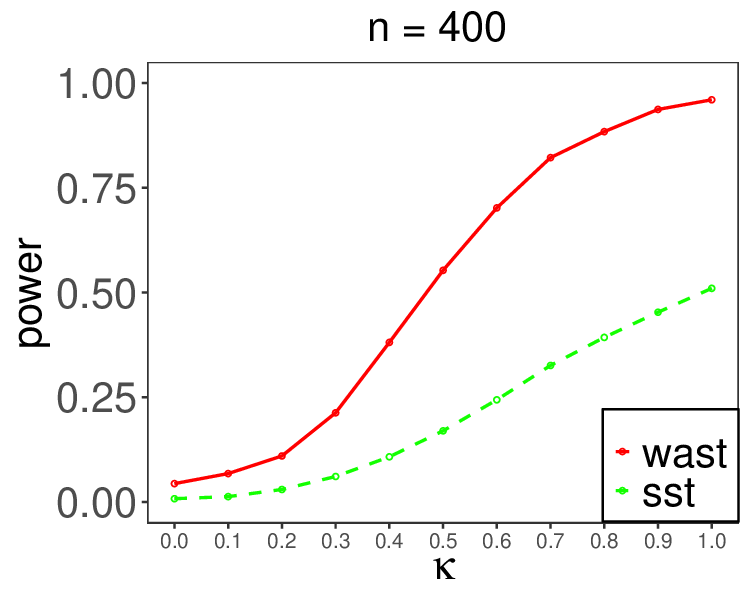}
		\includegraphics[scale=0.3]{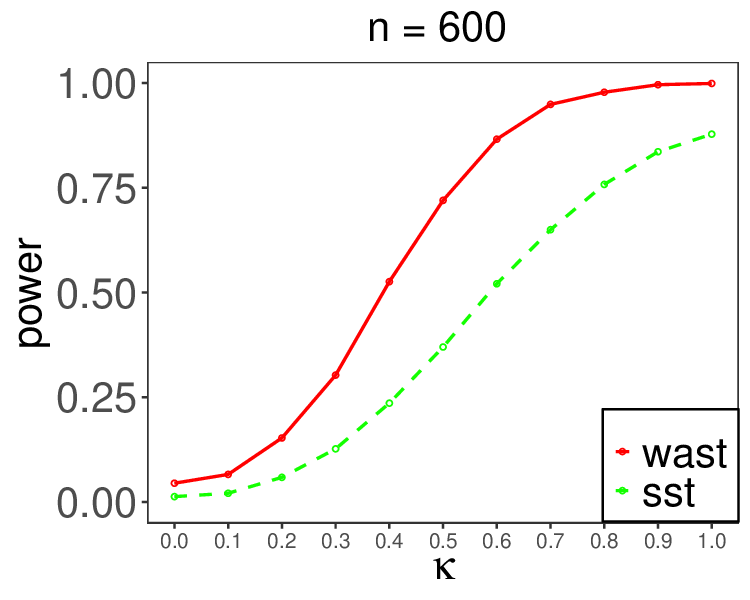}     \\
		\includegraphics[scale=0.3]{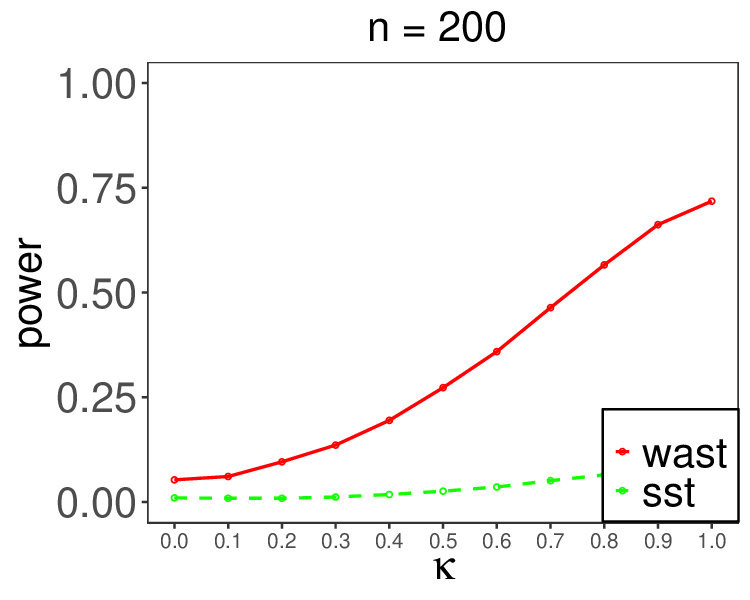}
		\includegraphics[scale=0.3]{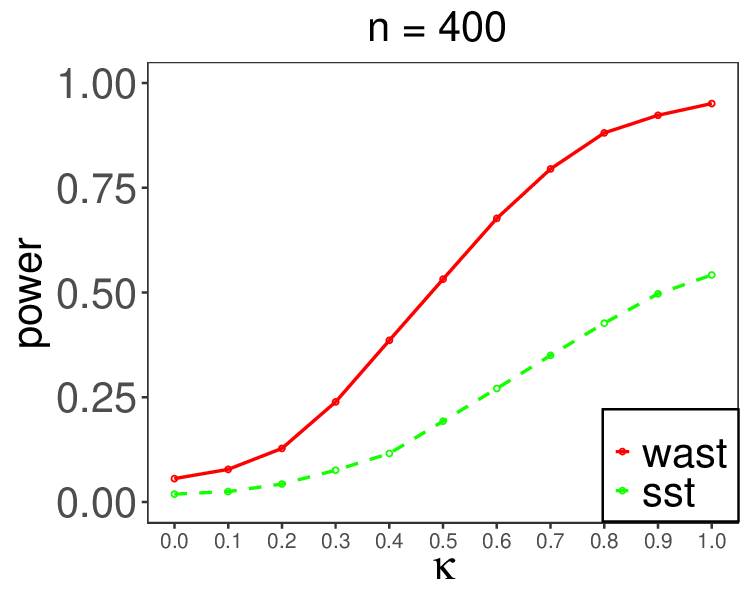}
		\includegraphics[scale=0.3]{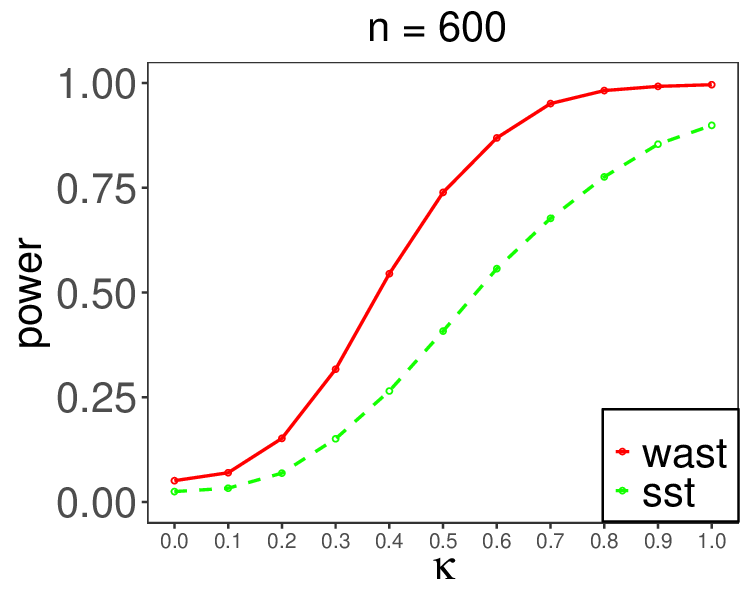}     \\
		\includegraphics[scale=0.3]{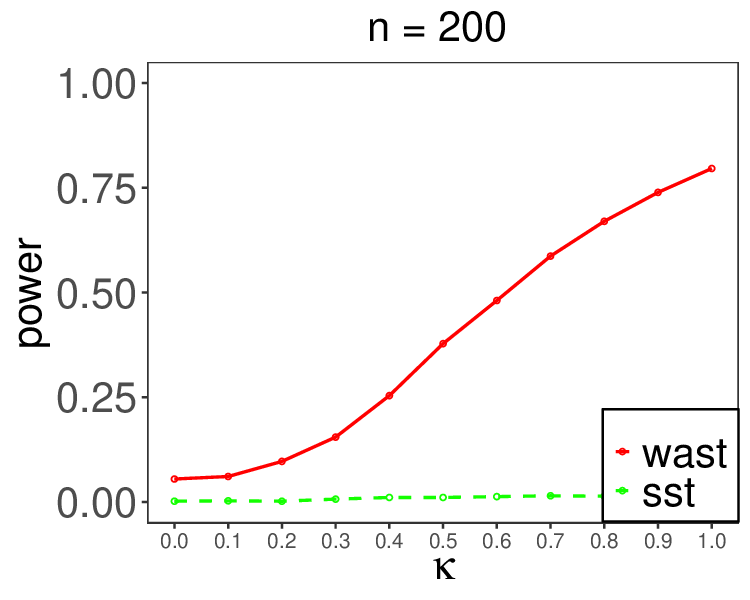}
		\includegraphics[scale=0.3]{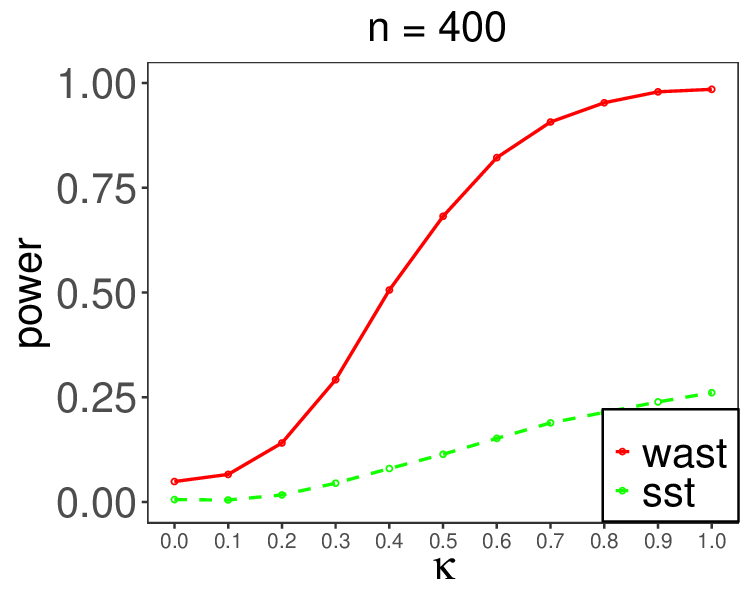}
		\includegraphics[scale=0.3]{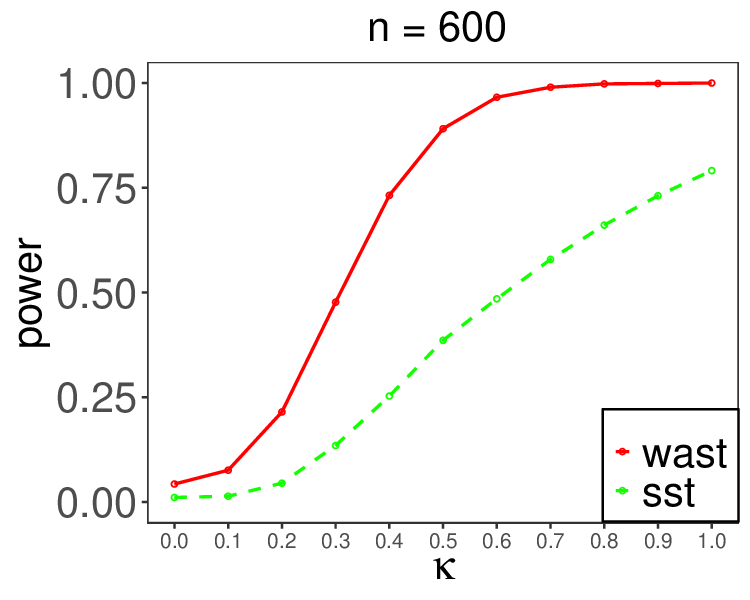}    \\
		\includegraphics[scale=0.3]{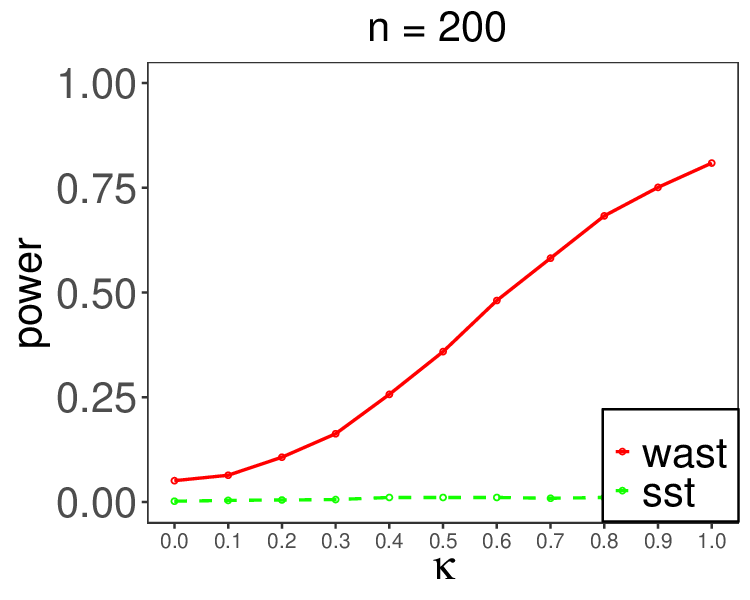}
		\includegraphics[scale=0.3]{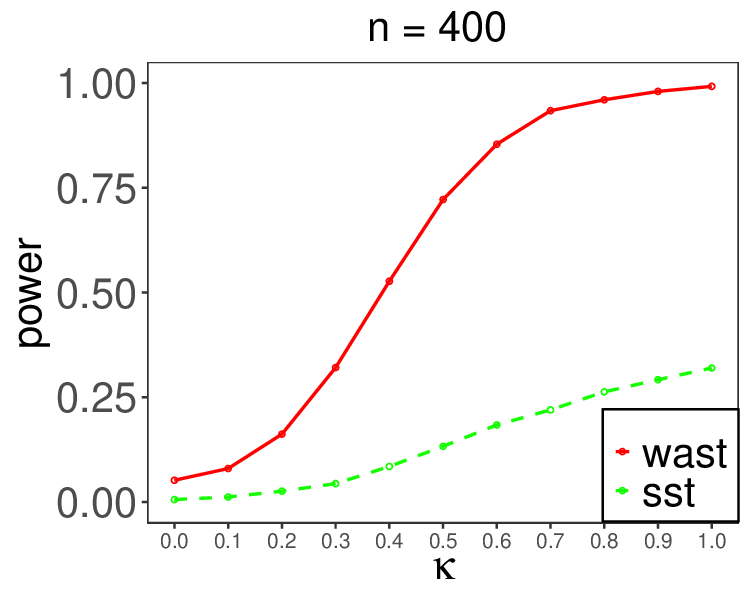}
		\includegraphics[scale=0.3]{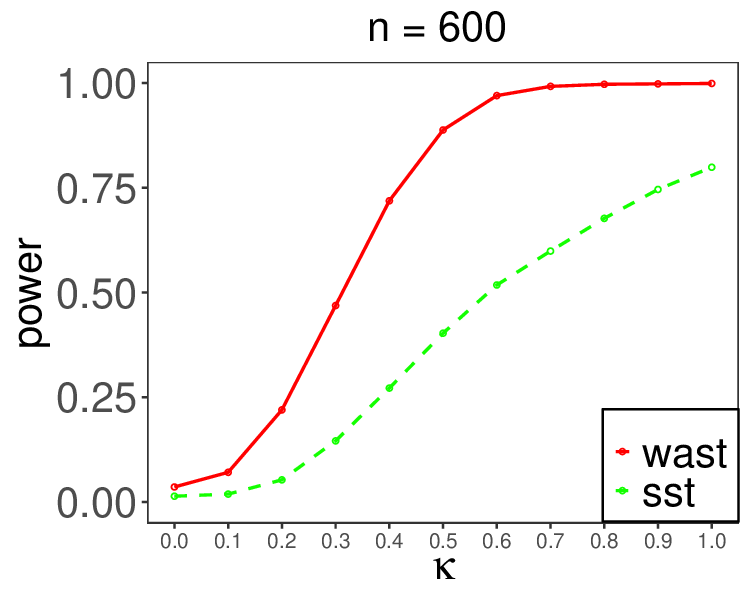}
		\caption{\it Powers of test statistic semiparametric model (B1+P1) with large numbers of dense $\bZ$ by the proposed WAST (red solid line) and SST (green dashed line). From top to bottom, each row depicts the powers for $(p,q)=(2,100)$, $(p,q)=(2,500)$, $(p,q)=(6,100)$, $(p,q)=(6,500)$, $(p,q)=(11,100)$, and $(p,q)=(11,500)$.}
		\label{fig_semiparam2_dense}
	\end{center}
\end{figure}

Type \uppercase\expandafter{\romannumeral1} errors ($\kappa=0$) for probit, quantile and semiparametric models are listed in Table \ref{table_size_ee_dense}. We can see from Table \ref{table_size_ee_dense} that the size of the proposed WAST are close to the nominal significance level $0.05$, but for most scenarios the size of the SST are much smaller than 0.05. From Figure \ref{fig_probit_dense}-\ref{fig_semiparam2_dense} we have same conclusion as these in Section \ref{simulation_ee_sparse}. We omit the detailed analysis here.

\begin{table}[htp!]
	\def~{\hphantom{0}}
    \tiny
	\caption{Type \uppercase\expandafter{\romannumeral1} errors of the proposed WAST and SST based on resampling for probit regression model (ProbitRE), quantile regression (QuantRE) and semiparametric model (SPMoldel) with large numbers of dense $\bZ$. The nominal significant level is 0.05.
}
	\resizebox{\textwidth}{!}{
    \begin{threeparttable}
		\begin{tabular}{llcccccccc}
			\hline
			\multirow{2}{*}{Model}&\multirow{2}{*}{$(p,q)$}
			&\multicolumn{2}{c}{ $n=200$} && \multicolumn{2}{c}{ $n=400$} && \multicolumn{2}{c}{ $n=600$} \\
			\cline{3-4} \cline{6-7} \cline{9-10}
			&&WAST& SST && WAST& SST && WAST& SST \\
			\cline{3-10}
			ProbitRE &$(2,100)$         & 0.035 & 0.000 && 0.058 & 0.001 && 0.040 & 0.013 \\
			&$(2,500)$                  & 0.035 & 0.002 && 0.047 & 0.007 && 0.056 & 0.013 \\
			&$(6,100)$                  & 0.052 & 0.000 && 0.041 & 0.000 && 0.056 & 0.000 \\
			&$(6,500)$                  & 0.056 & 0.000 && 0.046 & 0.000 && 0.056 & 0.003 \\
			&$(11,100)$                 & 0.056 & 0.000 && 0.048 & 0.000 && 0.049 & 0.000 \\
			&$(11,500)$                 & 0.061 & 0.000 && 0.040 & 0.000 && 0.053 & 0.000 \\
			[1 ex]
			SPModel &$(2,100)$          & 0.056 & 0.028 && 0.047 & 0.028 && 0.068 & 0.048 \\
			(B1+P1)&$(2,500)$           & 0.060 & 0.029 && 0.049 & 0.044 && 0.054 & 0.046 \\
			&$(6,100)$                  & 0.065 & 0.004 && 0.049 & 0.017 && 0.049 & 0.019 \\
			&$(6,500)$                  & 0.044 & 0.009 && 0.049 & 0.017 && 0.053 & 0.016 \\
			&$(11,100)$                 & 0.055 & 0.001 && 0.051 & 0.010 && 0.046 & 0.006 \\
			&$(11,500)$                 & 0.037 & 0.001 && 0.046 & 0.008 && 0.048 & 0.007 \\
			[1 ex]
			SPModel &$(2,100)$          & 0.056 & 0.023 && 0.040 & 0.024 && 0.050 & 0.039 \\
			(B2+P2)&$(2,500)$           & 0.048 & 0.022 && 0.055 & 0.042 && 0.055 & 0.038 \\
			&$(6,100)$                  & 0.044 & 0.007 && 0.044 & 0.008 && 0.045 & 0.013 \\
			&$(6,500)$                  & 0.053 & 0.010 && 0.056 & 0.019 && 0.051 & 0.025 \\
			&$(11,100)$                 & 0.055 & 0.002 && 0.049 & 0.006 && 0.043 & 0.011 \\
			&$(11,500)$                 & 0.051 & 0.002 && 0.052 & 0.006 && 0.036 & 0.014 \\
			[1 ex]
			QuantRE &$(2,100)$          & 0.054 & 0.016 && 0.059 & 0.025 && 0.058 & 0.032 \\
			($\tau=0.2$)&$(2,500)$      & 0.054 & 0.018 && 0.050 & 0.017 && 0.050 & 0.032 \\
			&$(6,100)$                  & 0.059 & 0.005 && 0.054 & 0.010 && 0.057 & 0.017 \\
			&$(6,500)$                  & 0.054 & 0.003 && 0.048 & 0.005 && 0.062 & 0.018 \\
			&$(11,100)$                 & 0.050 & 0.000 && 0.054 & 0.004 && 0.053 & 0.011 \\
			&$(11,500)$                 & 0.045 & 0.001 && 0.050 & 0.007 && 0.057 & 0.017 \\
			[1 ex]
			QuantRE &$(2,100)$          & 0.050 & 0.031 && 0.036 & 0.030 && 0.047 & 0.048 \\
			($\tau=0.5$)&$(2,500)$      & 0.056 & 0.039 && 0.058 & 0.044 && 0.052 & 0.047 \\
			&$(6,100)$                  & 0.050 & 0.022 && 0.037 & 0.018 && 0.051 & 0.026 \\
			&$(6,500)$                  & 0.058 & 0.012 && 0.052 & 0.016 && 0.050 & 0.041 \\
			&$(11,100)$                 & 0.065 & 0.007 && 0.062 & 0.019 && 0.048 & 0.017 \\
			&$(11,500)$                 & 0.043 & 0.008 && 0.047 & 0.017 && 0.053 & 0.031 \\
			[1 ex]
			QuantRE &$(2,100)$          & 0.053 & 0.017 && 0.055 & 0.022 && 0.047 & 0.037 \\
			($\tau=0.7$)&$(2,500)$      & 0.057 & 0.026 && 0.066 & 0.040 && 0.058 & 0.033 \\
			&$(6,100)$                  & 0.043 & 0.007 && 0.060 & 0.018 && 0.054 & 0.025 \\
			&$(6,500)$                  & 0.046 & 0.009 && 0.044 & 0.015 && 0.058 & 0.033 \\
			&$(11,100)$                 & 0.049 & 0.002 && 0.061 & 0.008 && 0.044 & 0.019 \\
			&$(11,500)$                 & 0.037 & 0.003 && 0.047 & 0.009 && 0.056 & 0.014 \\
			\hline
		\end{tabular}
\end{threeparttable}
	}
	\label{table_size_ee_dense}
\end{table}

\subsection{Change plane analysis with sparse \texorpdfstring{$\bZ$}{} for GLMs}\label{simulation_glm_sparse}

Consider the GLM in the main paper with the canonical parameter
\begin{align*}
\mu = \tbX\trans\balpha + \bX\trans\bbeta\bone(\bZ\trans\btheta\geq0).
\end{align*}
For Gaussian and Poisson families,
we set $\alpha_1 = 0.5$,  and for binomial family $\alpha_1$ is chosen so that the proportion of cases in the data set is ${1}/{3}$ on average under the null hypothesis. $\theta_2,\cdots,\theta_q$ are generated randomly from uniform distribution $U(1, 3)$, and $\theta_1$ is chosen as the negative of the 0.65 percentile of $Z_2\theta_2+\cdots+Z_q\theta_q$, which means that $\bZ\trans\btheta$ divides the population into two groups with 0.35 and 0.65 observations, respectively. For the predictor, we generate $(Z_1, Z_2,\cdots,Z_q)\trans$ from a multivariate normal distribution with mean $\bzero$ and covariance $\bI_q$. We set $\tX_j=\bone(v_j>0)$ with $j=2,\cdots,r$, $X_k = \bone(v_k>0)$ with $k=2,\cdots,p$, and $\tX_1=1$, $X_1=1$, and $Z_1=1$.

For all models with change plane analysis, we evaluated the power under a sequence of alternative models indexed by $\kappa$, that is $H_1^{\kappa}: \bbeta^{\kappa}=\kappa\bbeta^*$ with $\kappa = i/10$ for semiparametric model and $\kappa=i/20$ for others, $i=1,\cdots,10$, where $\bbeta^*=(1,\cdots,1)\trans$. We set sample size $n=(300, 600)$ for GLMs, 1000 repetitions and 1000 bootstrap samples, and report in Figure \ref{fig_gaussian_dense}-\ref{fig_poisson_dense} the performance for the WAST, SST and SLRT. We calculate the SST and SLRT over $\{\btheta^{(k)}=(\theta^{(k)}_1,\cdots,\theta^{(k)}_q)\trans: k=1,\cdots,K\}$ with the number of threshold values $K=1000$. Let $\btheta^{(k)}_{-1}=\tilde{\btheta}^{(k)}_{-1}/\|\tilde{\btheta}^{(k)}_{-1}\|$, where $\btheta^{(k)}_{-1}=(\theta^{(k)}_2,\cdots,\theta^{(k)}_q)\trans$ and $\tilde{\btheta}^{(k)}_{-1}=(\tilde{\theta}^{(k)}_2,\cdots,\tilde{\theta}^{(k)}_q)\trans$, and $\tilde{\btheta}^{(k)}_{-1}$ is drawn independently form $(r-1)$-variate standard normal distribution.  For each $\theta_1^{(k)}$, $k=1,\cdots,K$, we selected it by equal grid search in the range from the lower 10th percentile to upper 10th percentile of the data points of $\{\theta^{(k)}_2Z_{2i}+\cdots+\theta^{(k)}_qZ_{qi}\}_{i=1}^n$, which is same as that in \cite{2020Threshold}. Here we consider four combinations of $(p,q)=(2, 100), (2, 500), (6, 100), (6, 500), (11, 100), (11, 500)$.

\begin{figure}[!ht]
	\begin{center}
		\includegraphics[scale=0.26]{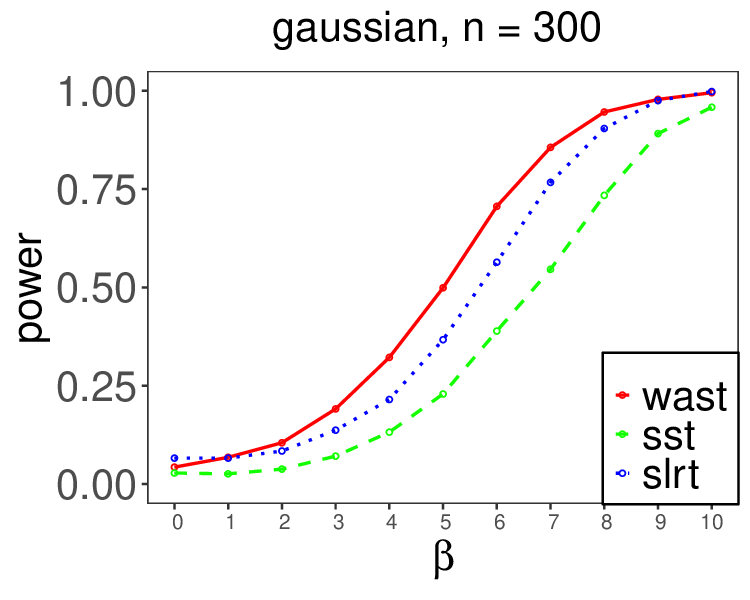}
		\includegraphics[scale=0.26]{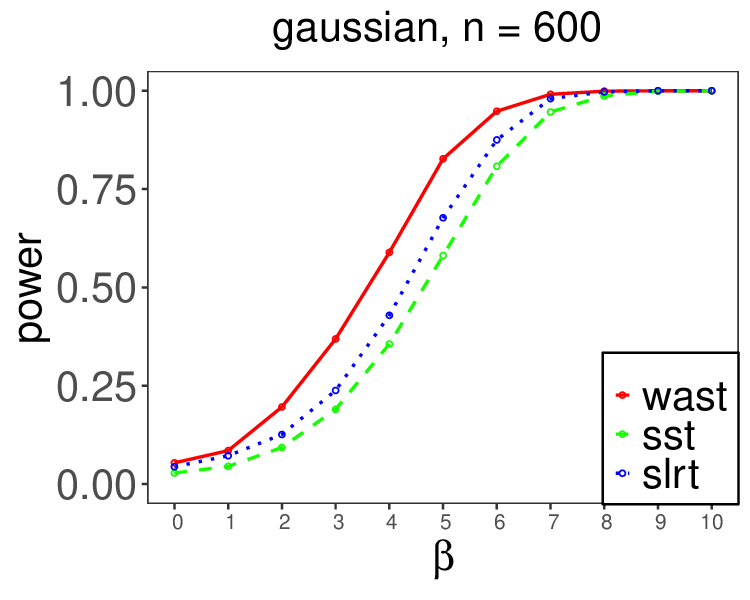}
		\includegraphics[scale=0.26]{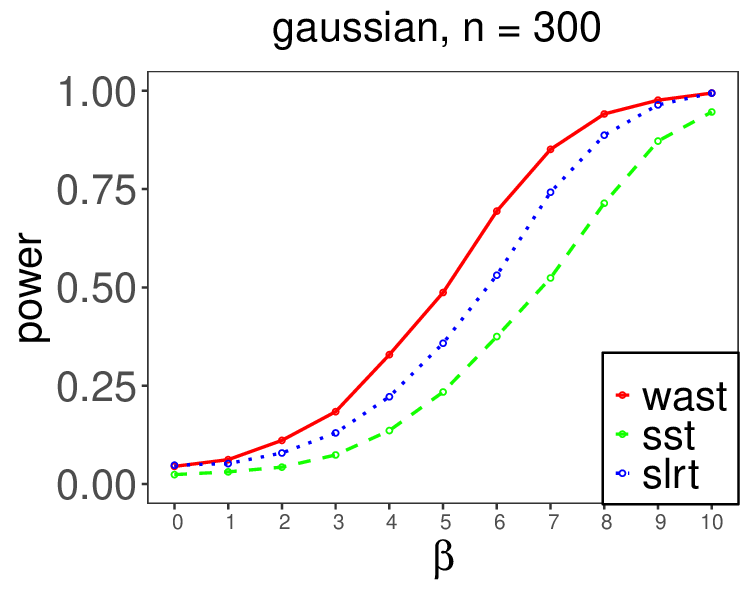}
		\includegraphics[scale=0.26]{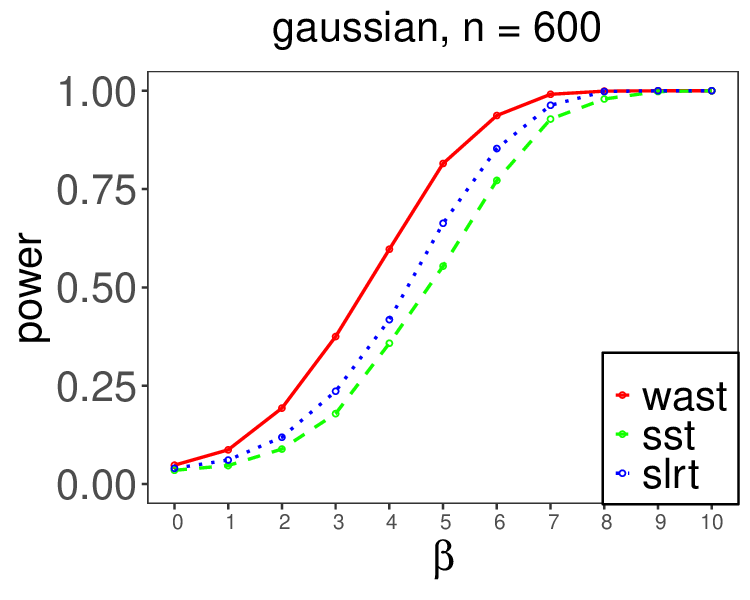}  \\
		\includegraphics[scale=0.26]{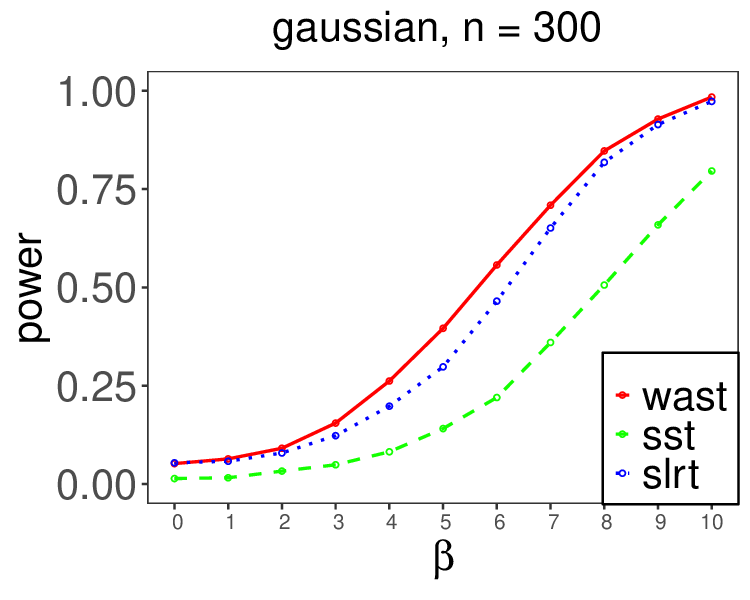}
		\includegraphics[scale=0.26]{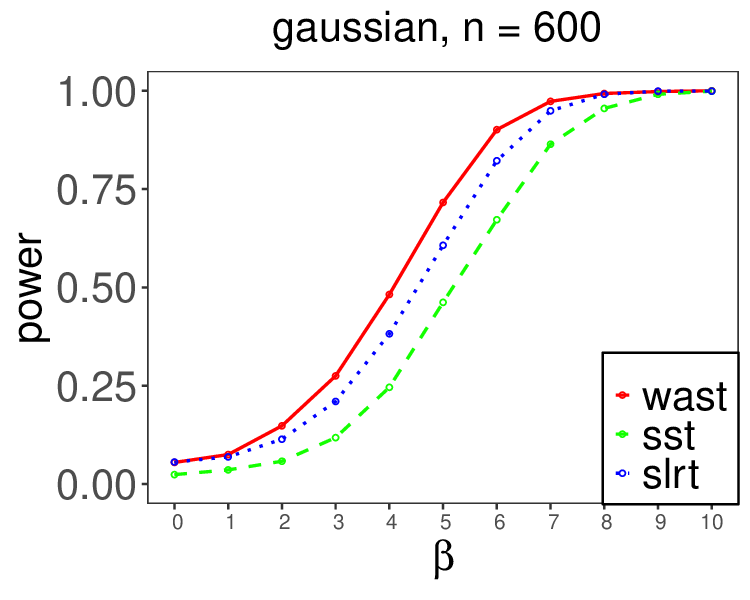}
		\includegraphics[scale=0.26]{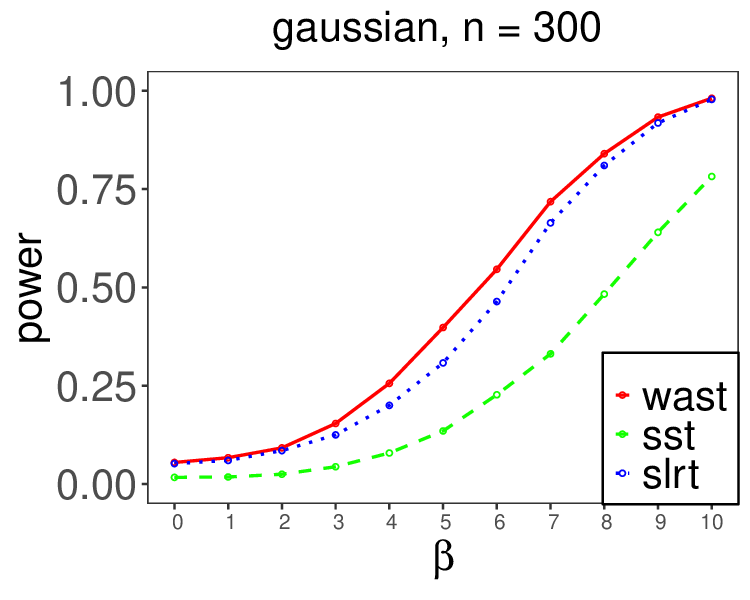}
		\includegraphics[scale=0.26]{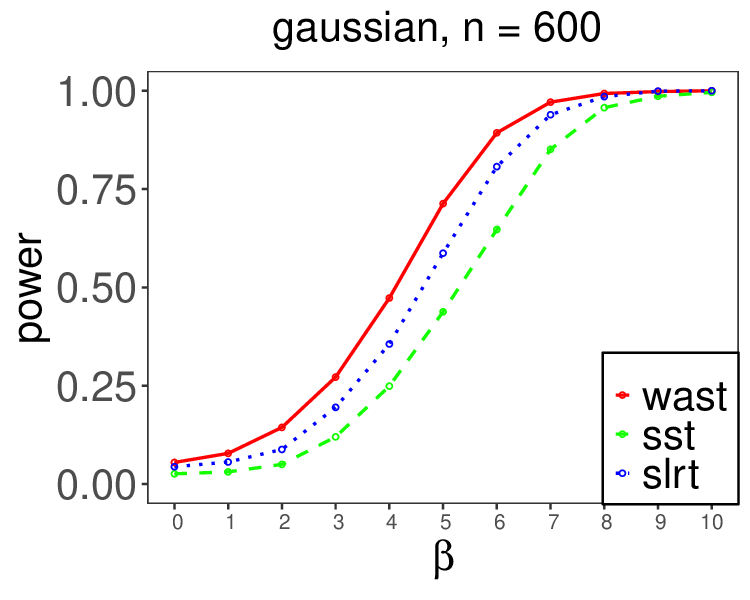}  \\
		\includegraphics[scale=0.26]{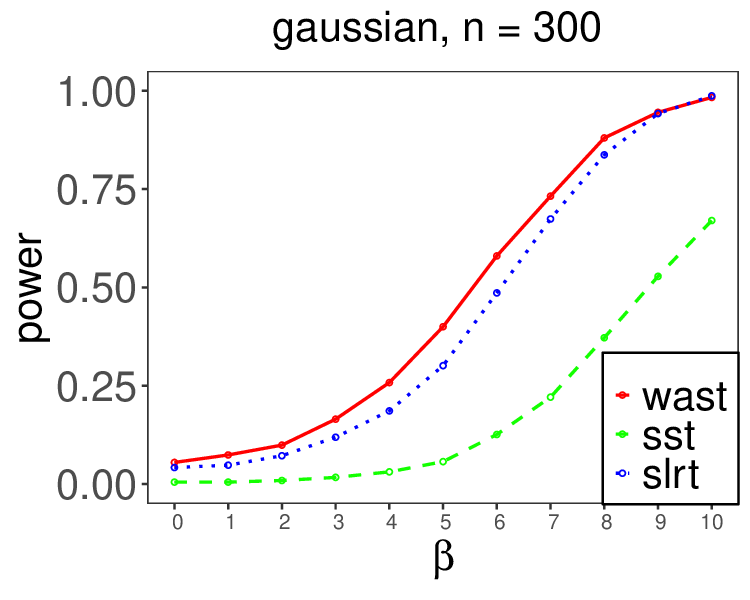}
		\includegraphics[scale=0.26]{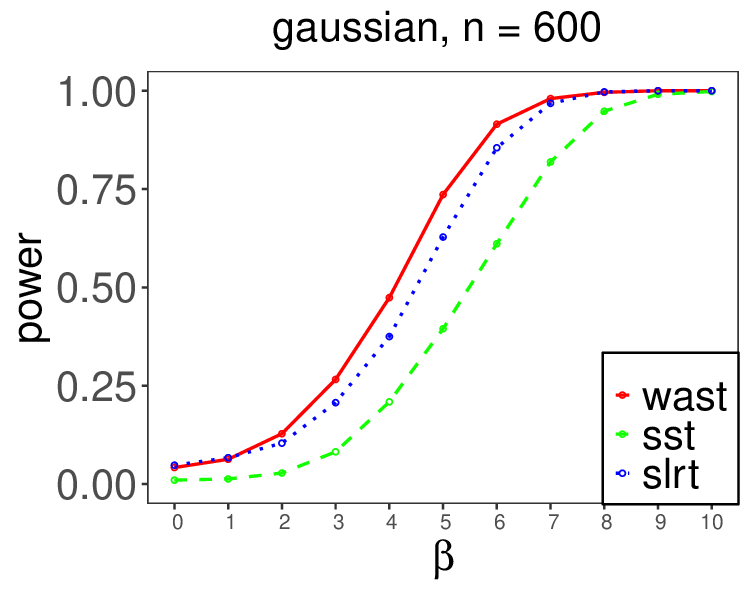}
		\includegraphics[scale=0.26]{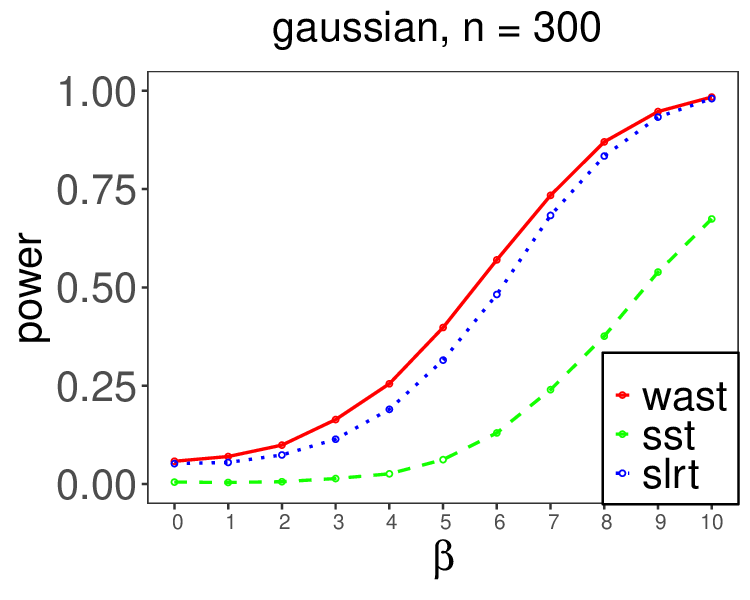}
		\includegraphics[scale=0.26]{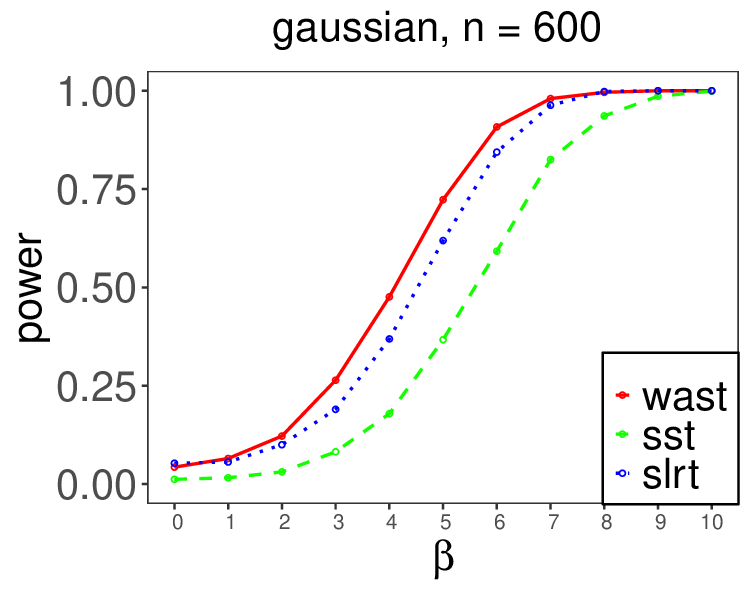}
		\caption{\it Powers of test statistic for GLM with Gaussian and with large numbers of sparse $\bZ$ by the proposed WAST (red solid line), SST (green dashed line) and SLRT (blue dotted line). From top to bottom, each row depicts the powers for $p=2, 6, 11$. From left to right, each column depicts the powers for $q=100, 100, 500, 500$.
		}
		\label{fig_gaussian_sparse}
	\end{center}
\end{figure}

\begin{figure}[!ht]
	\begin{center}
		\includegraphics[scale=0.26]{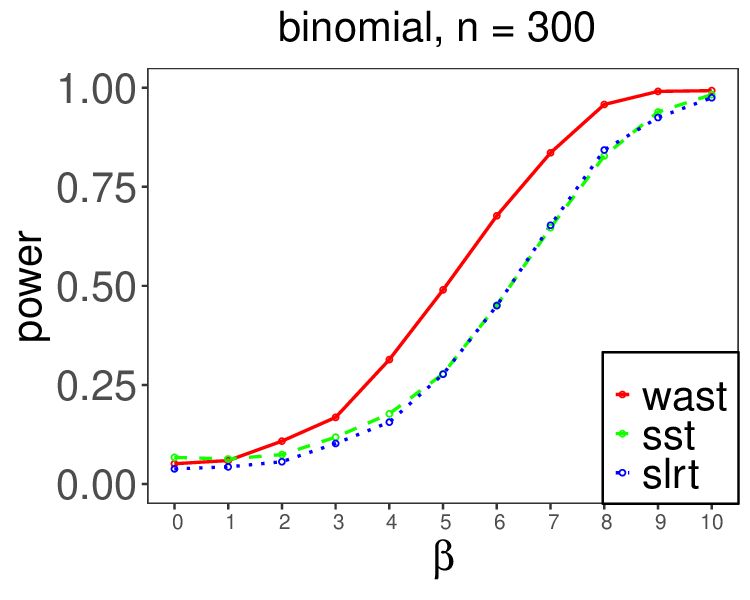}
		\includegraphics[scale=0.26]{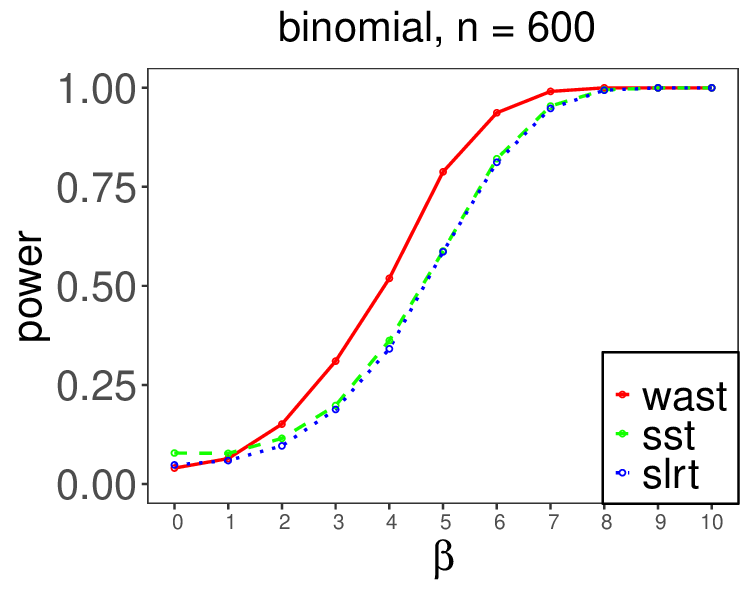}
		\includegraphics[scale=0.26]{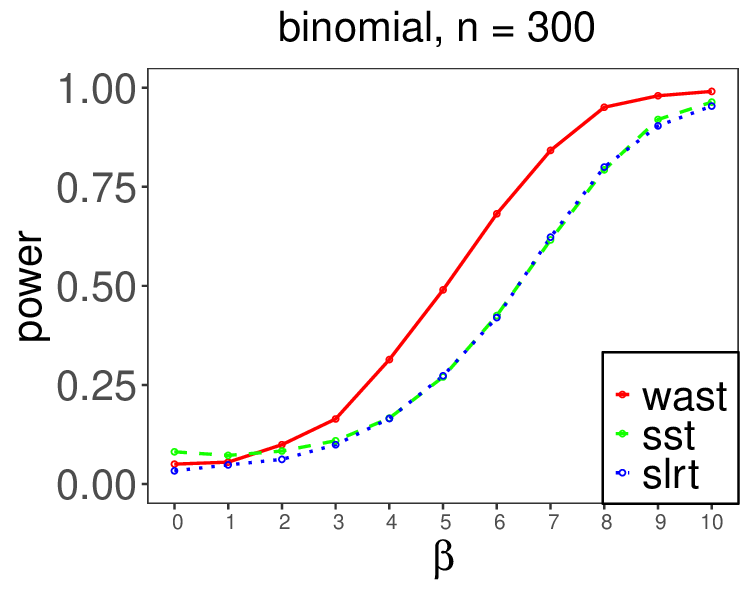}
		\includegraphics[scale=0.26]{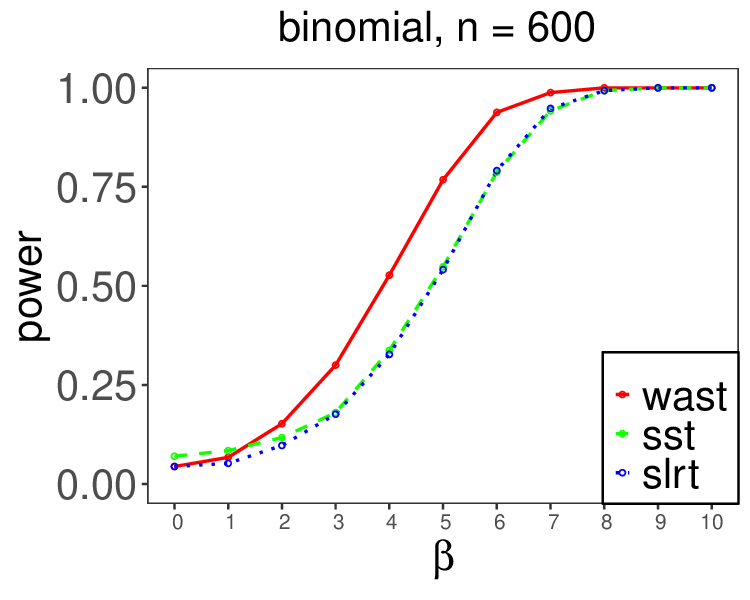}  \\
		\includegraphics[scale=0.26]{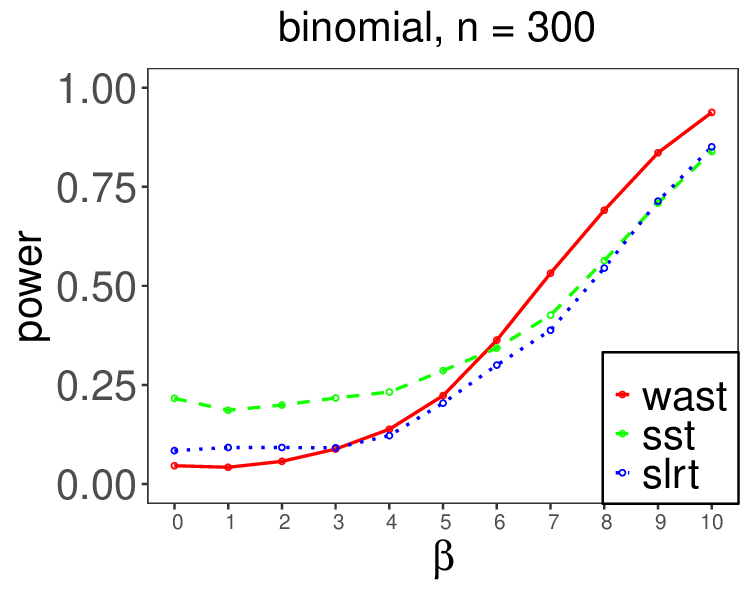}
		\includegraphics[scale=0.26]{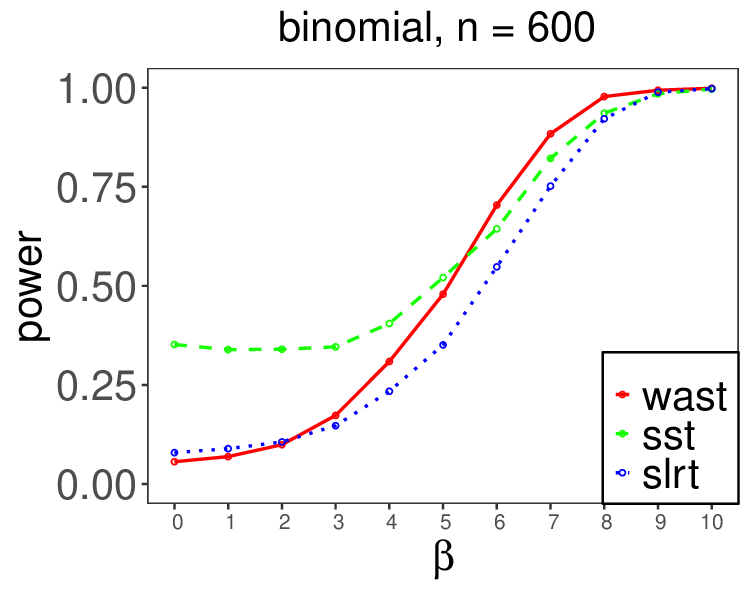}
		\includegraphics[scale=0.26]{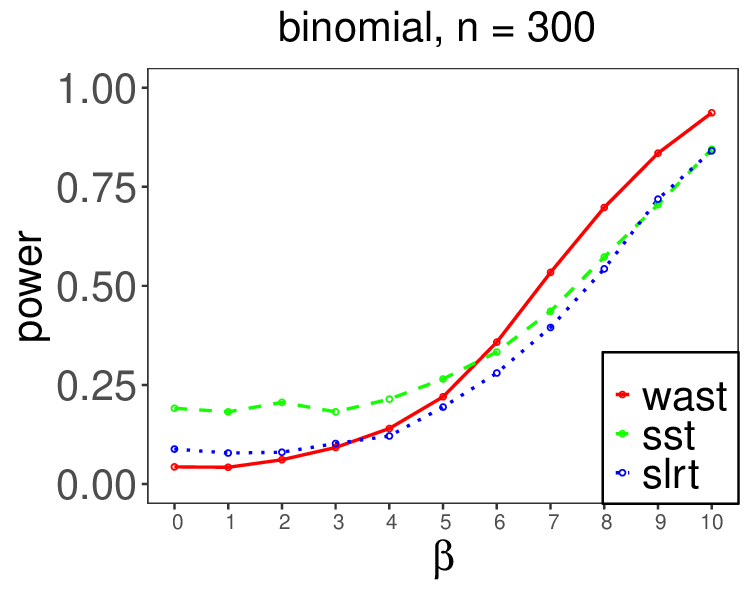}
		\includegraphics[scale=0.26]{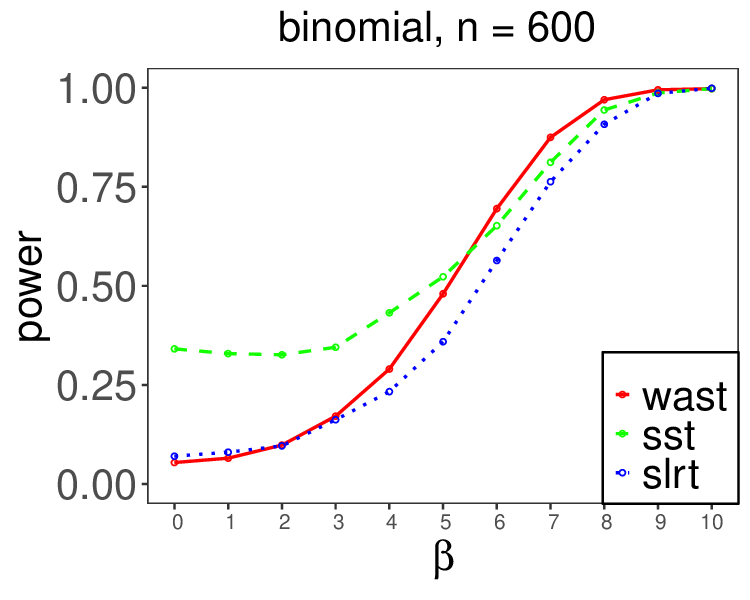}  \\
		\includegraphics[scale=0.26]{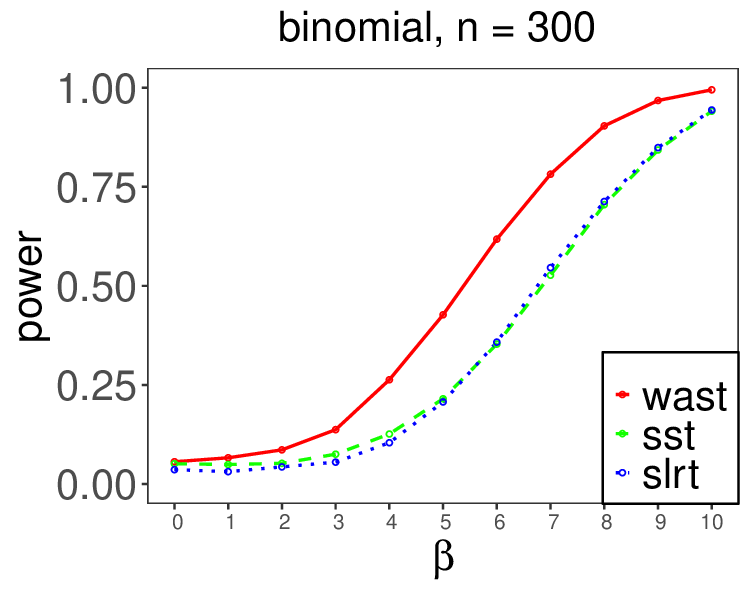}
		\includegraphics[scale=0.26]{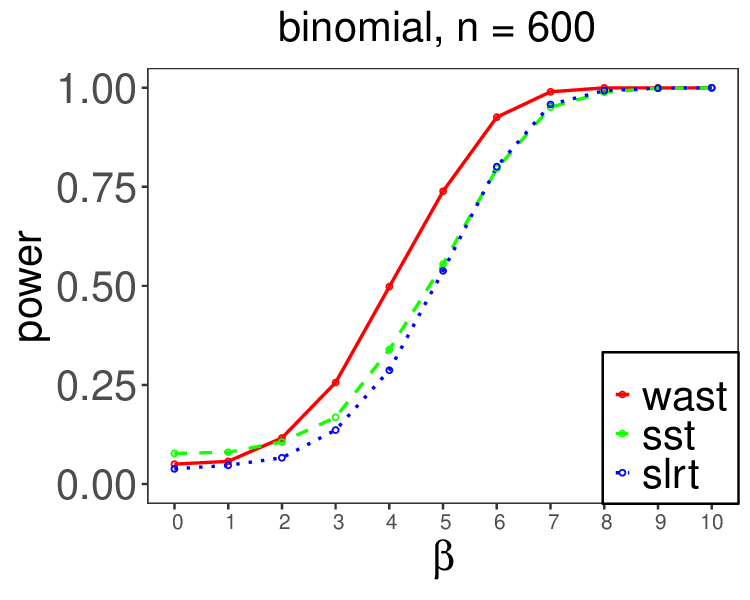}
		\includegraphics[scale=0.26]{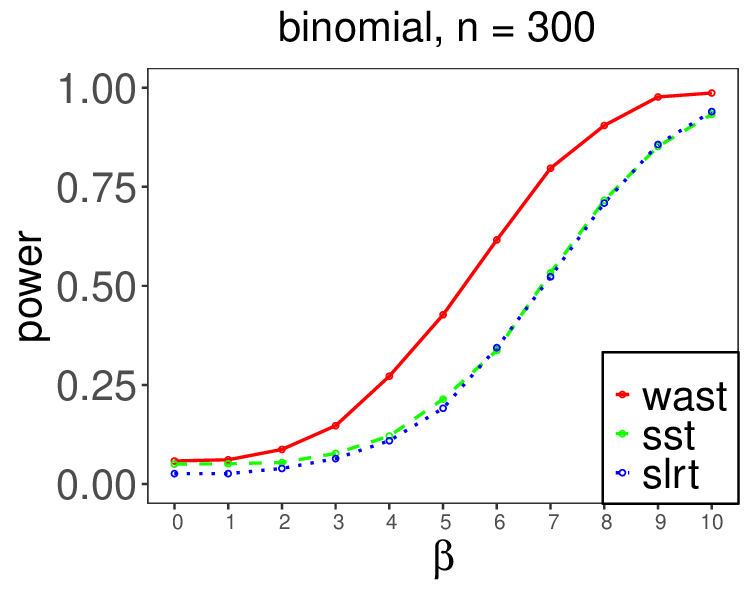}
		\includegraphics[scale=0.26]{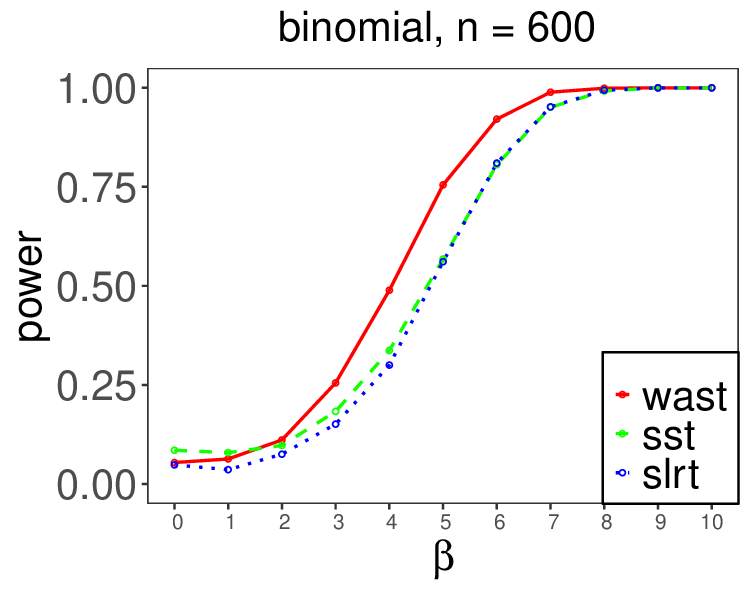}
		\caption{\it Powers of test statistic for GLM with binomial and with large numbers of sparse $\bZ$ by the proposed WAST (red solid line), SST (green dashed line) and SLRT (blue dotted line). From top to bottom, each row depicts the powers for $p=2, 6, 11$. From left to right, each column depicts the powers for $q=100, 100, 500, 500$.
		}
		\label{fig_binomial_sparse}
	\end{center}
\end{figure}

\begin{figure}[!ht]
	\begin{center}
		\includegraphics[scale=0.26]{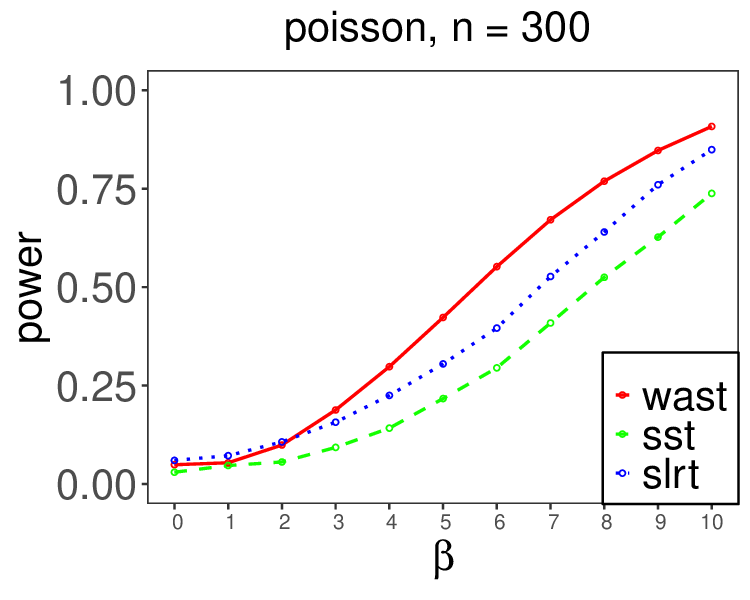}
		\includegraphics[scale=0.26]{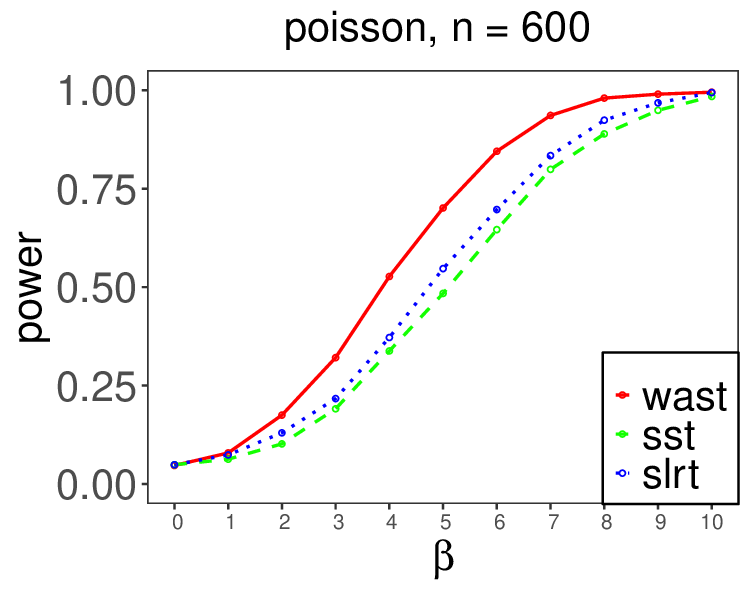}
		\includegraphics[scale=0.26]{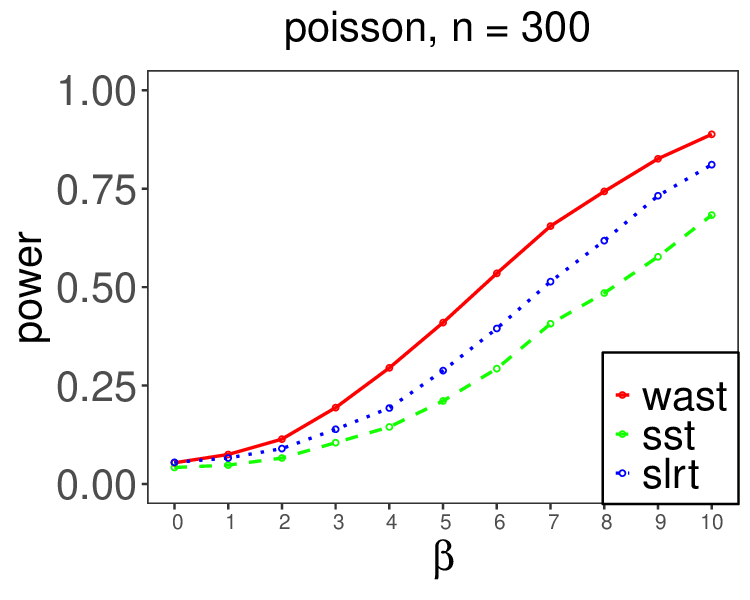}
		\includegraphics[scale=0.26]{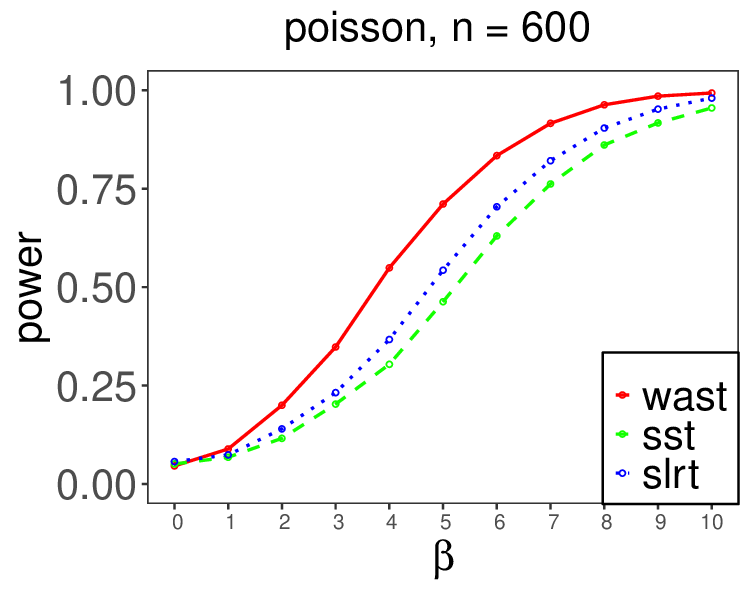}  \\
		\includegraphics[scale=0.26]{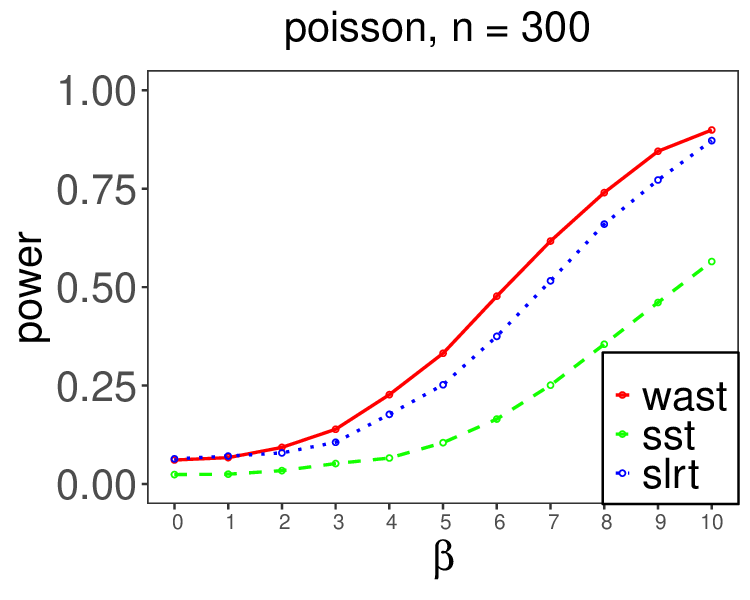}
		\includegraphics[scale=0.26]{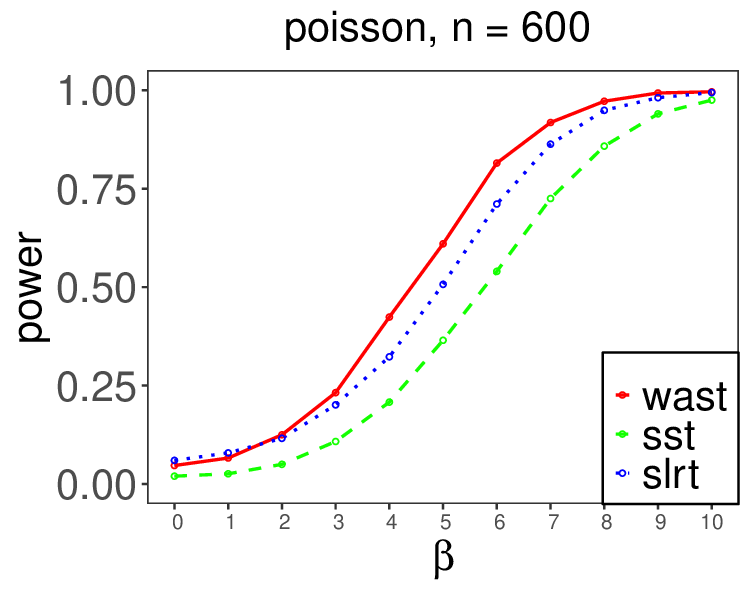}
		\includegraphics[scale=0.26]{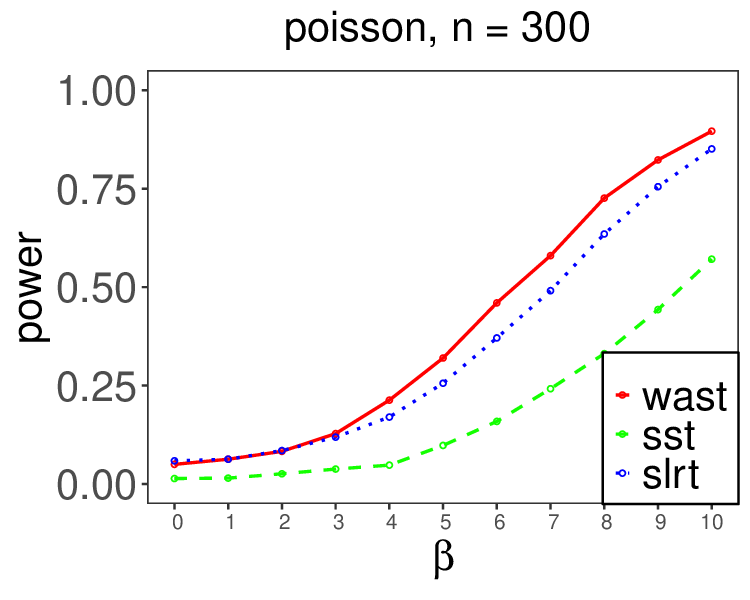}
		\includegraphics[scale=0.26]{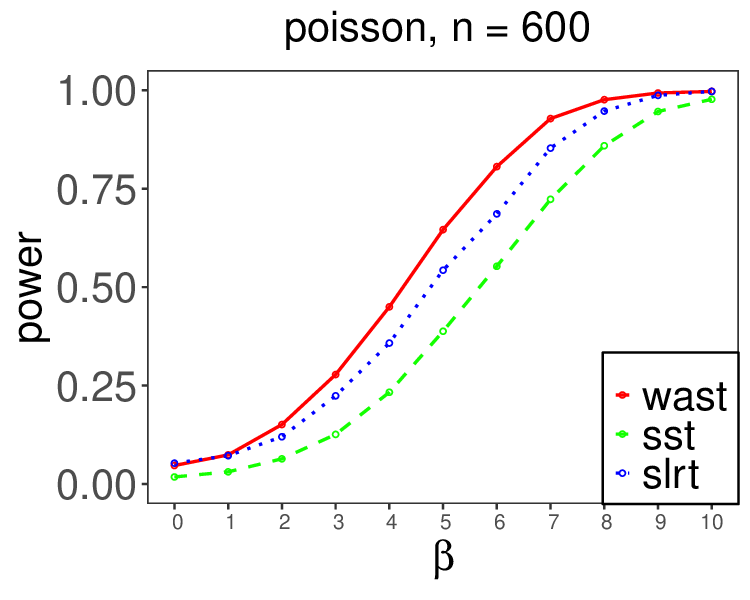}  \\
		\includegraphics[scale=0.26]{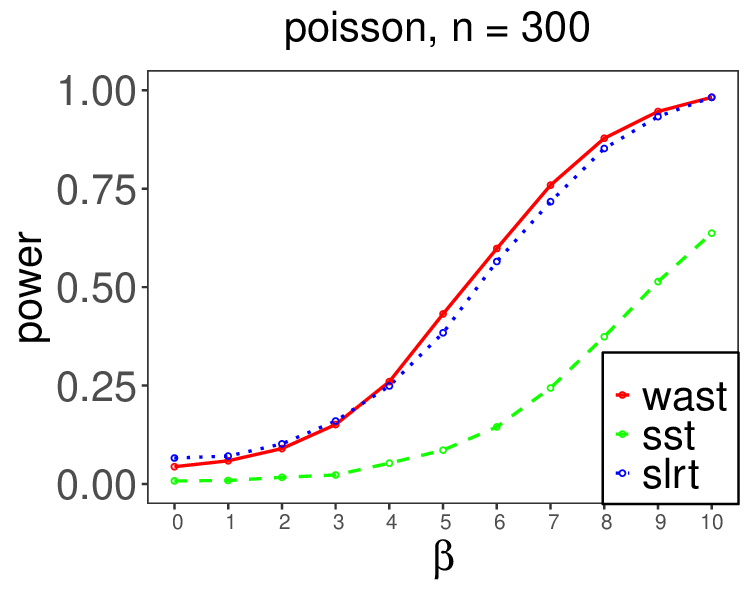}
		\includegraphics[scale=0.26]{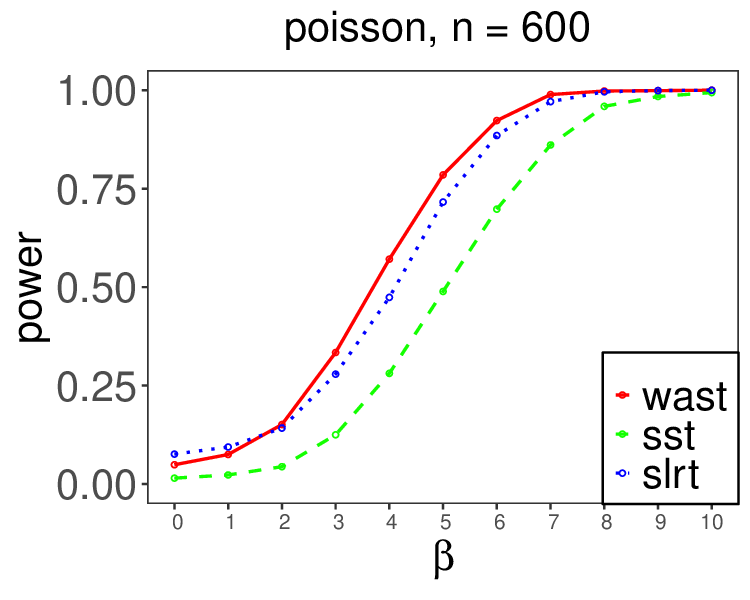}
		\includegraphics[scale=0.26]{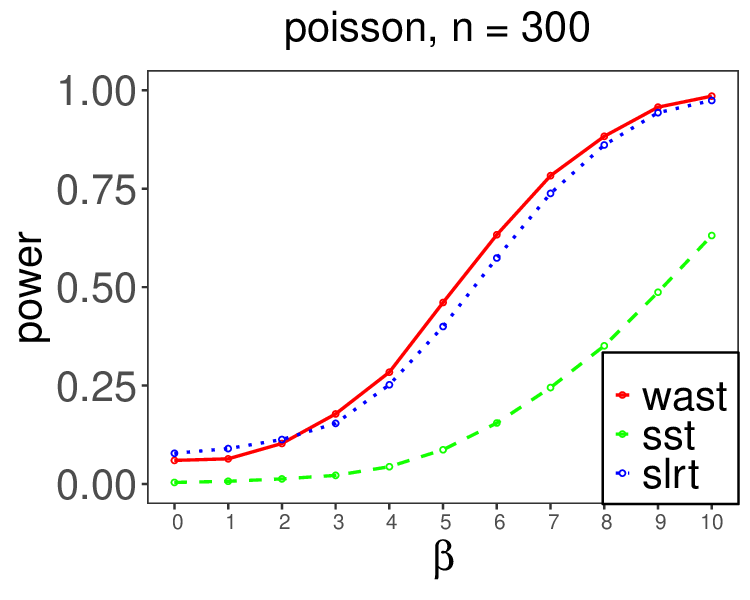}
		\includegraphics[scale=0.26]{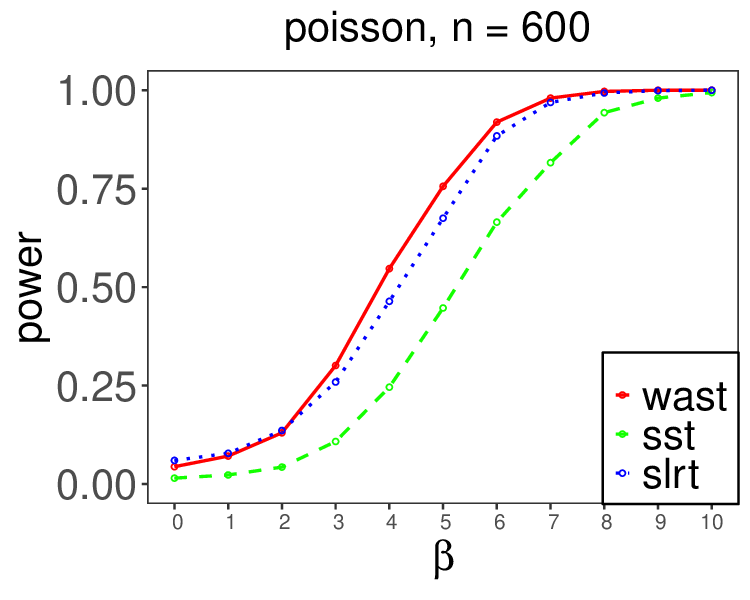}
		\caption{\it Powers of test statistic for GLM with Poisson and with large numbers of sparse $\bZ$ by the proposed WAST (red solid line), SST (green dashed line) and SLRT (blue dotted line). From top to bottom, each row depicts the powers for $p=2, 6, 11$. From left to right, each column depicts the powers for $q=100, 100, 500, 500$.
		}
		\label{fig_poisson_sparse}
	\end{center}
\end{figure}

Type \uppercase\expandafter{\romannumeral1} errors ($\kappa=0$) for GLMs are list in Table \ref{table_size_glm_sparse}. We can see from Table \ref{table_size_glm_sparse} that the sizes of the proposed WAST and the SLRT are close to the nominal significance level $0.05$, but for most scenarios the size of the SST are much smaller than 0.05.
Figure \ref{fig_gaussian_sparse}-\ref{fig_poisson_sparse} indicate that powers become greater as sample size $n$ increases, which are verified the asymptotic theory. The proposed WAST has comparable power with the SST for the semiparametric model, but the size of the SST is much less than the nominal level 0.05 in 
Figure \ref{fig_gaussian_sparse}-\ref{fig_poisson_sparse} show that the powers of the proposed WAST are uniformly greater than the competitors SST and SLRT.

\begin{table}[htp!]
	\def~{\hphantom{0}}
	\tiny
	\caption{Type \uppercase\expandafter{\romannumeral1} errors of the proposed WAST, SST and SLRT based on resampling for GLMs with large numbers of sparse $\bZ$. The nominal significant level is 0.05.
	}
	\resizebox{\textwidth}{!}{
		\begin{threeparttable}
			\begin{tabular}{llcccccccc}
				\hline
				\multirow{2}{*}{Model}&\multirow{2}{*}{$(p,q)$}
				&\multicolumn{3}{c}{ $n=300$} && \multicolumn{3}{c}{ $n=600$} \\
				\cline{3-5} \cline{7-9}
				&& WAST& SST & SLRT&& WAST& SST & SLRT\\
				\cline{3-9}
				Gaussian &$(2,100)$         & 0.043 & 0.028 & 0.066 && 0.054 & 0.028 & 0.044 \\
				&$(2,500)$                  & 0.045 & 0.024 & 0.048 && 0.048 & 0.035 & 0.040 \\
				&$(6,100)$                  & 0.052 & 0.014 & 0.054 && 0.055 & 0.024 & 0.056 \\
				&$(6,500)$                  & 0.055 & 0.017 & 0.052 && 0.055 & 0.026 & 0.044 \\
				&$(11,100)$                 & 0.055 & 0.005 & 0.042 && 0.042 & 0.010 & 0.048 \\
				&$(11,500)$                 & 0.058 & 0.005 & 0.052 && 0.043 & 0.012 & 0.053 \\
				[1 ex]
				binomial &$(2,100)$          &0.051 & 0.067 & 0.038 && 0.040 & 0.078 & 0.048 \\\
				&$(2,500)$                 & 0.050 & 0.081 & 0.033 && 0.044 & 0.070 & 0.044 \\
				&$(6,100)$                  & 0.046 & 0.216 & 0.084 && 0.056 & 0.352 & 0.079 \\
				&$(6,500)$                  & 0.043 & 0.191 & 0.088 && 0.054 & 0.341 & 0.070 \\
				&$(11,100)$                 & 0.056 & 0.051 & 0.036 && 0.050 & 0.077 & 0.038 \\
				&$(11,500)$                 & 0.058 & 0.050 & 0.026 && 0.054 & 0.085 & 0.048 \\
				[1 ex]
				Poisson &$(2,100)$          & 0.049 & 0.030 & 0.060 && 0.047 & 0.049 & 0.049 \\
				&$(2,500)$                  & 0.054 & 0.042 & 0.055 && 0.046 & 0.051 & 0.057 \\
				&$(6,100)$                  & 0.061 & 0.024 & 0.064 && 0.047 & 0.020 & 0.060 \\
				&$(6,500)$                  & 0.050 & 0.014 & 0.059 && 0.047 & 0.018 & 0.053 \\
				&$(11,100)$                 & 0.044 & 0.008 & 0.066 && 0.049 & 0.015 & 0.076 \\
				&$(11,500)$                 & 0.060 & 0.004 & 0.078 && 0.044 & 0.015 & 0.060 \\
				\hline
			\end{tabular}
		\end{threeparttable}
	}
	\label{table_size_glm_sparse}
\end{table}

\subsection{Change plane analysis with dense \texorpdfstring{$\bZ$}{} for GLMs}\label{simulation_glm_dense}

For the sparse $\bZ$, we generate $\theta_2,\cdots,\theta_q$ from the uniform distribution $U(0,3)$ and set the rest being zero. Other settings are same as these in Section \ref{simulation_glm_sparse}.

\begin{figure}[!ht]
	\begin{center}
		\includegraphics[scale=0.26]{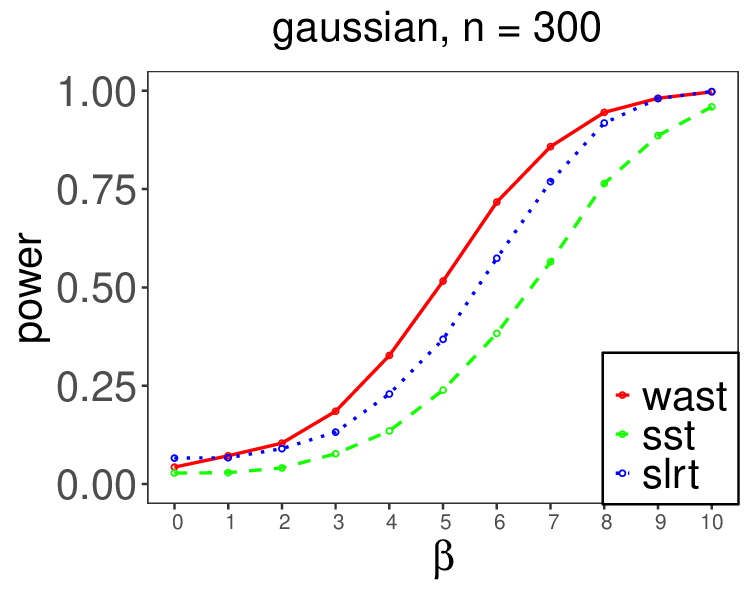}
		\includegraphics[scale=0.26]{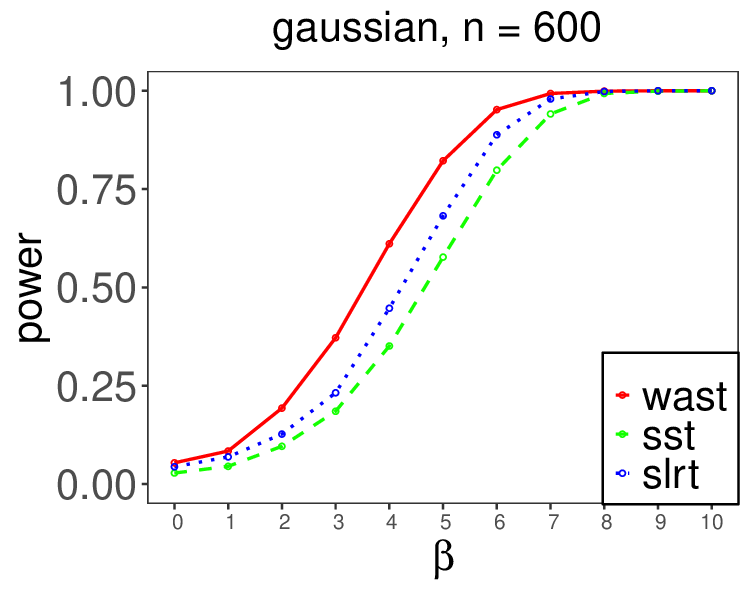}
		\includegraphics[scale=0.26]{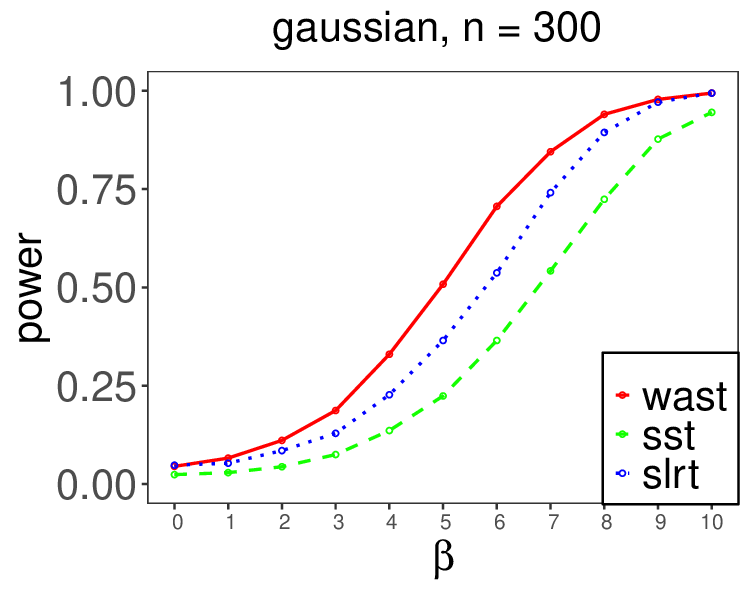}
		\includegraphics[scale=0.26]{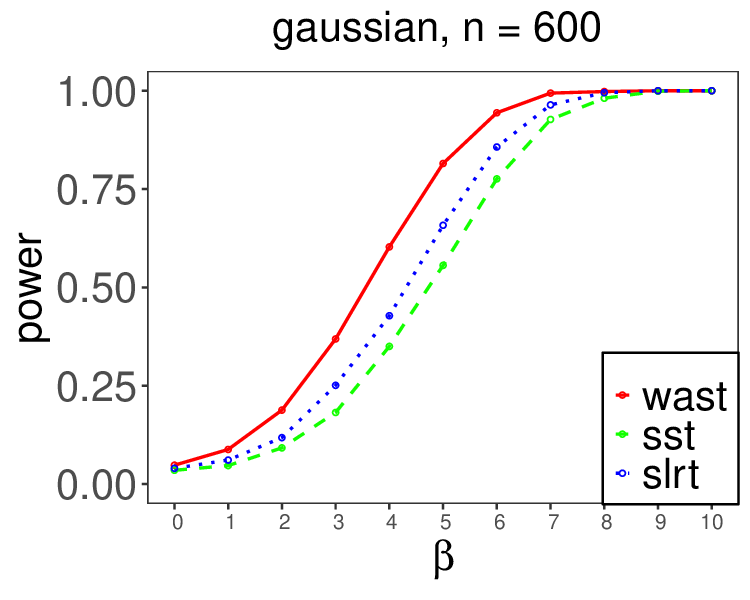}  \\
		\includegraphics[scale=0.26]{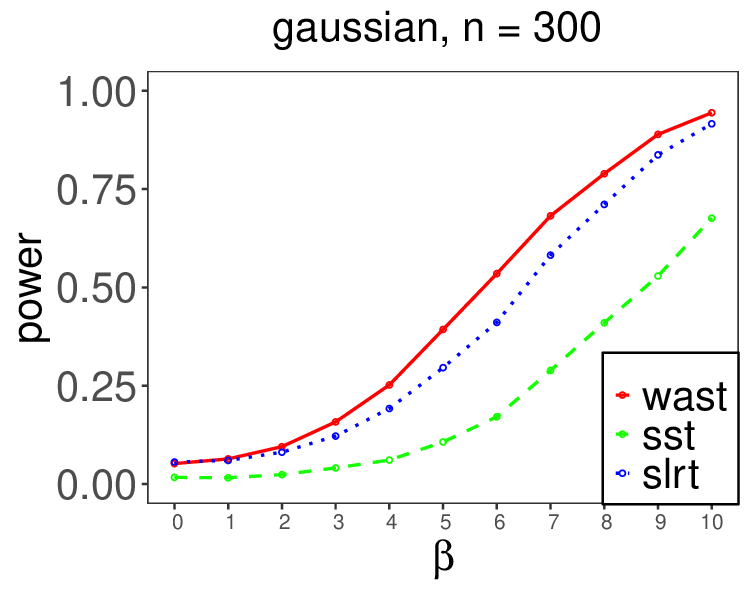}
		\includegraphics[scale=0.26]{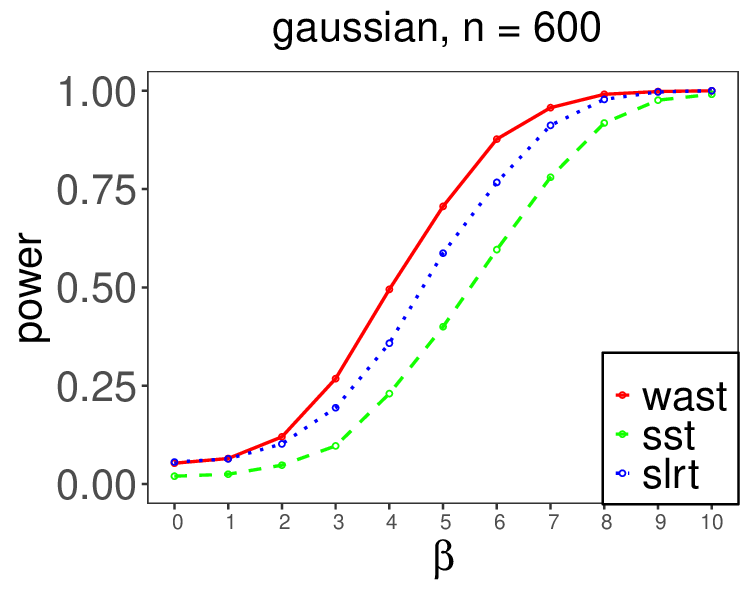}
		\includegraphics[scale=0.26]{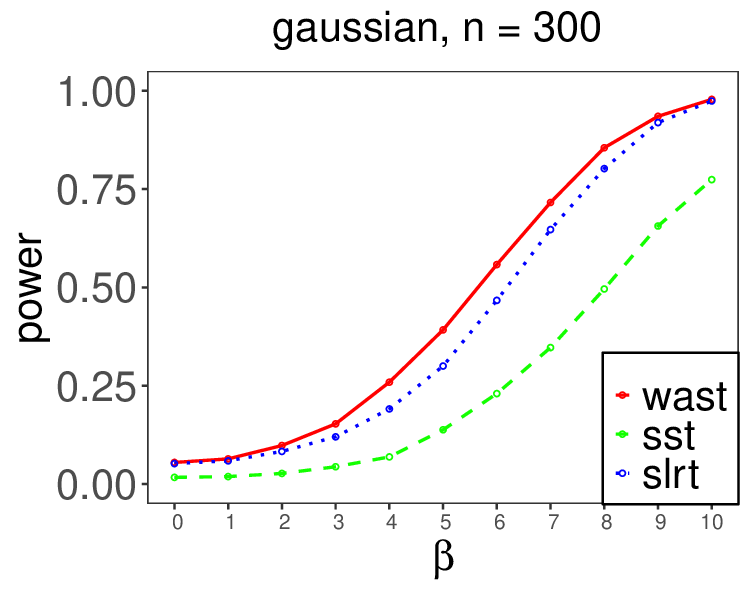}
		\includegraphics[scale=0.26]{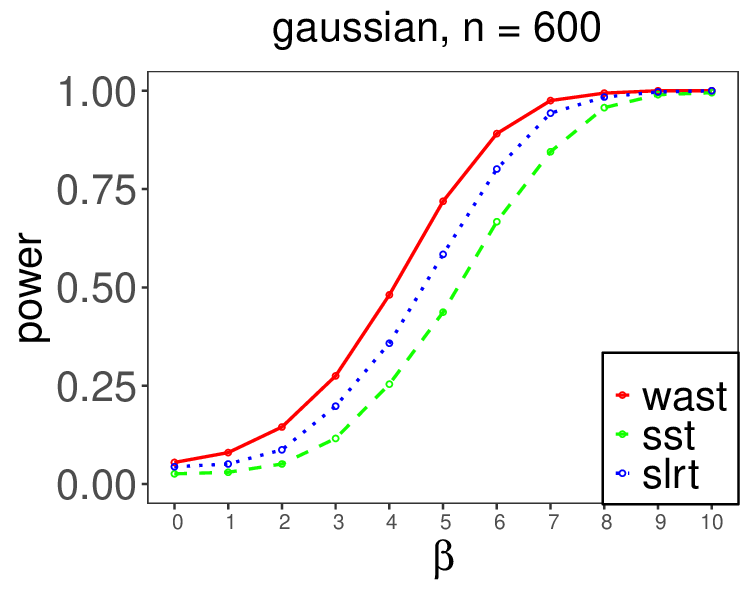}  \\
		\includegraphics[scale=0.26]{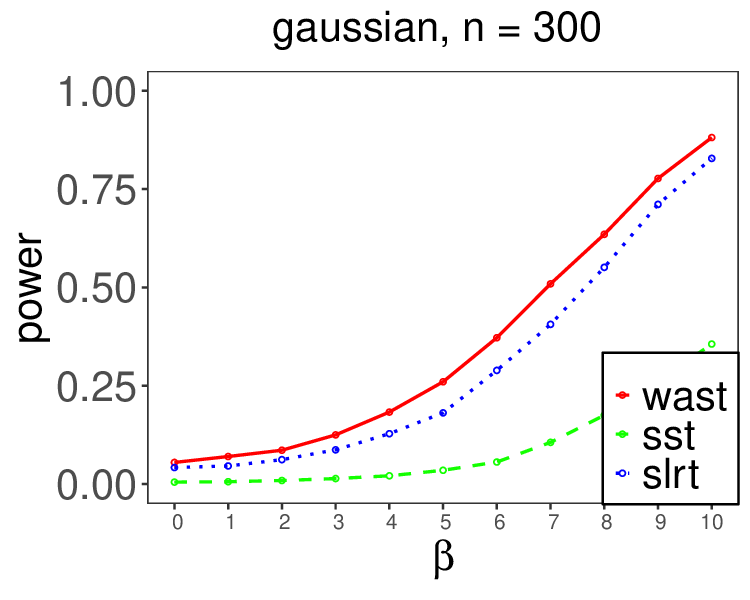}
		\includegraphics[scale=0.26]{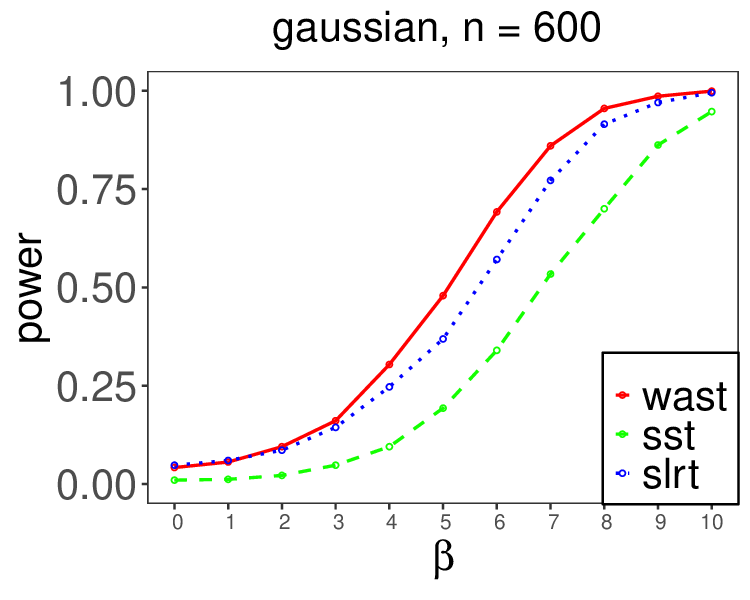}
		\includegraphics[scale=0.26]{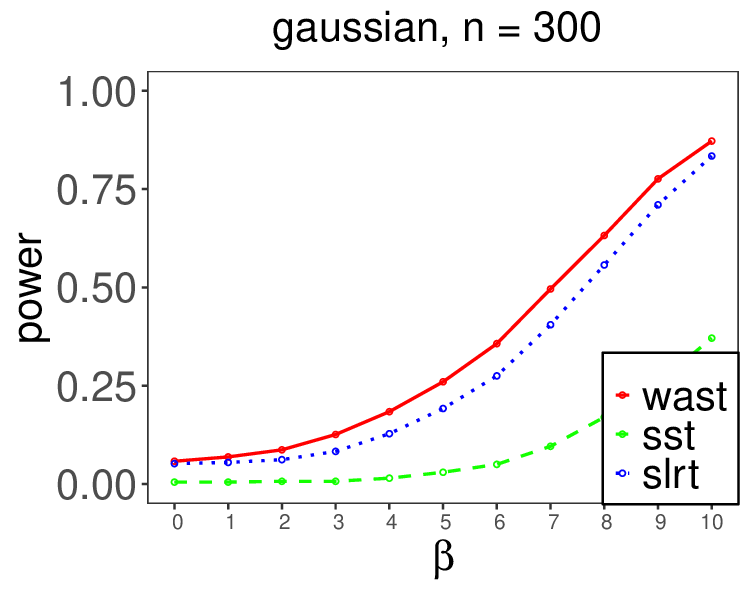}
		\includegraphics[scale=0.26]{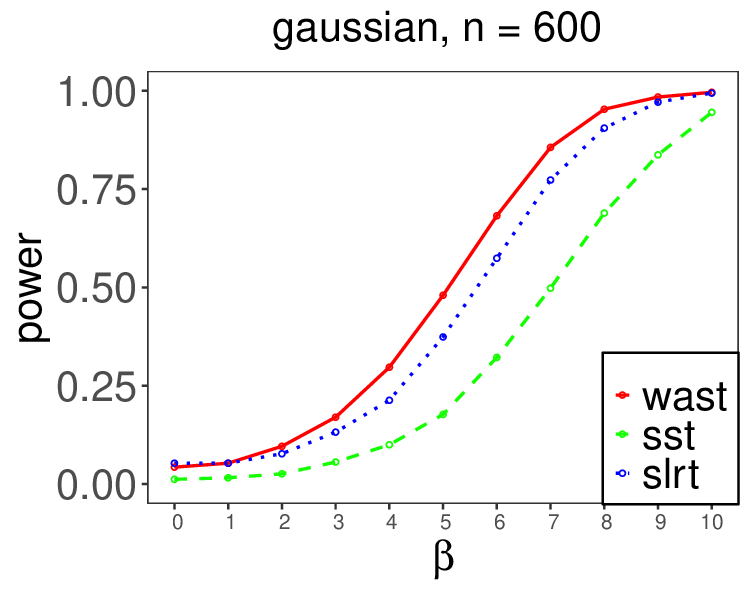}
		\caption{\it Powers of test statistic for GLM with Gaussian and with large numbers of dense $\bZ$ by the proposed WAST (red solid line), SST (green dashed line) and SLRT (blue dotted line). From top to bottom, each row depicts the powers for $p=2, 6, 11$. From left to right, each column depicts the powers for $q=100, 100, 500, 500$.
		}
		\label{fig_gaussian_dense}
	\end{center}
\end{figure}

\begin{figure}[!ht]
	\begin{center}
		\includegraphics[scale=0.26]{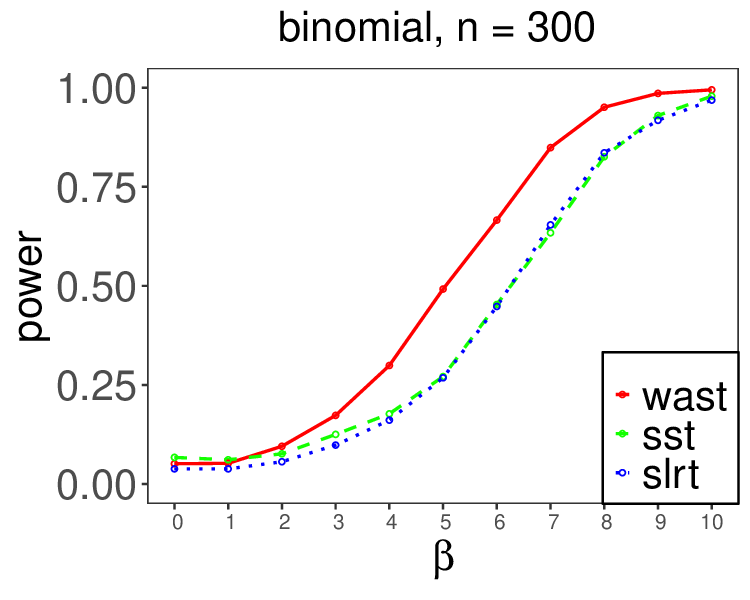}
		\includegraphics[scale=0.26]{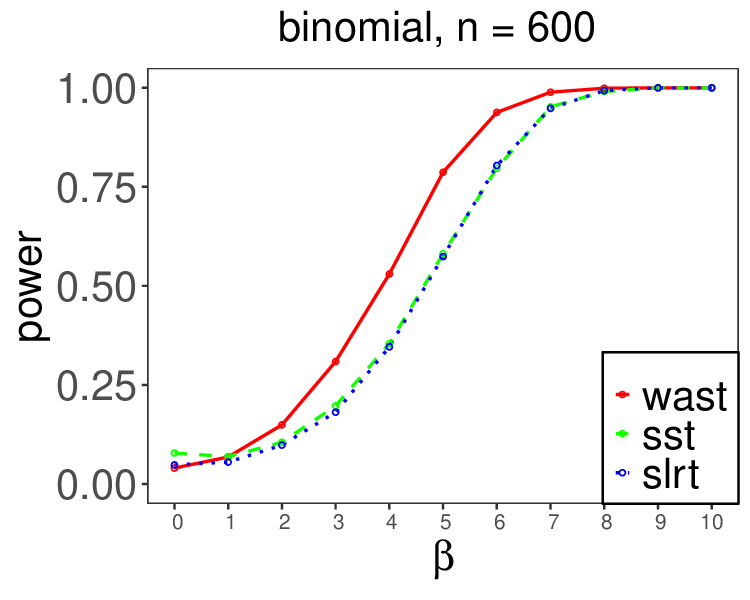}
		\includegraphics[scale=0.26]{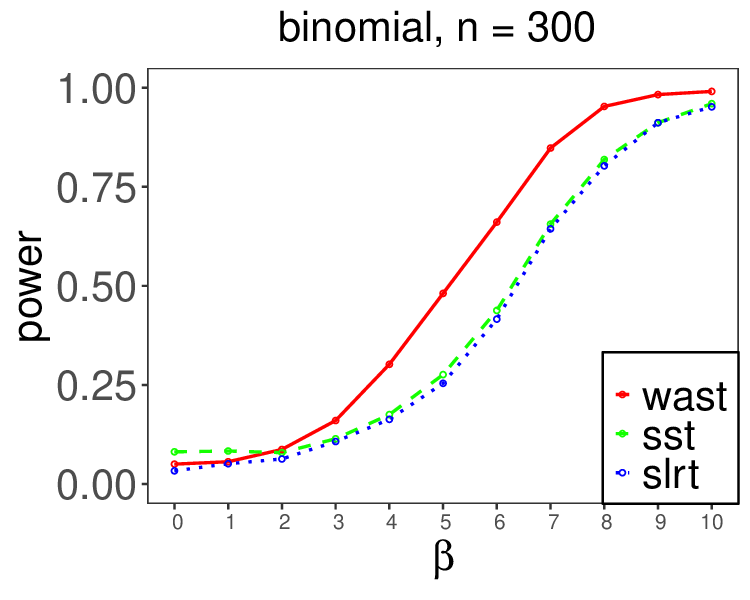}
		\includegraphics[scale=0.26]{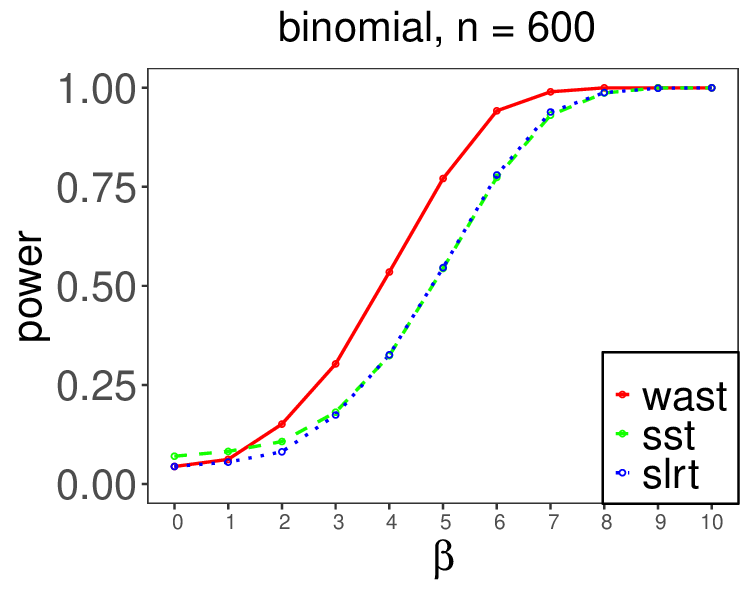}  \\
		\includegraphics[scale=0.26]{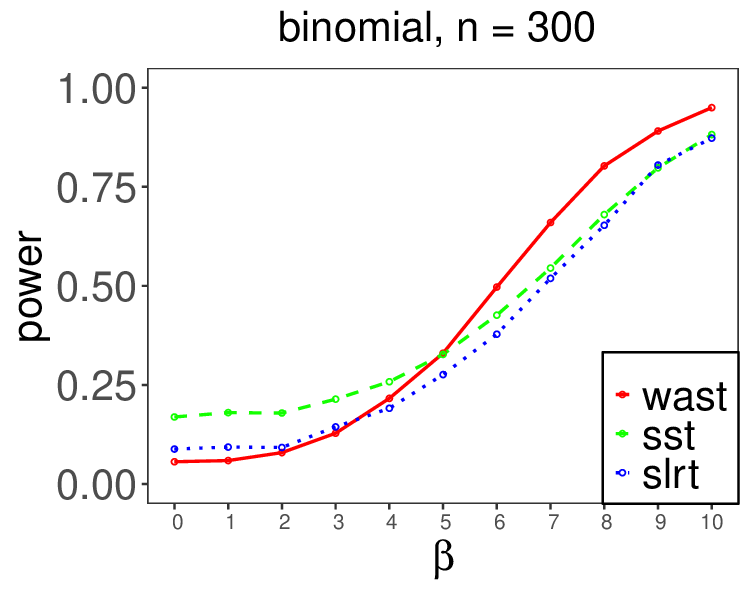}
		\includegraphics[scale=0.26]{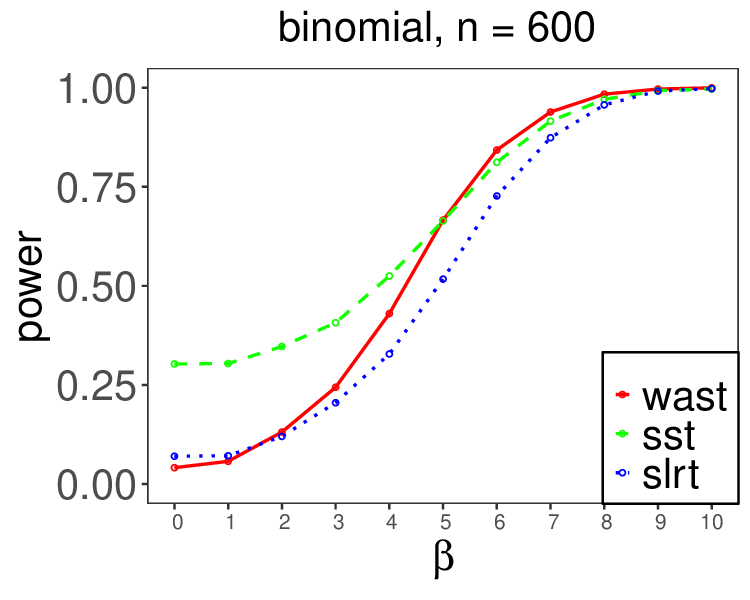}
		\includegraphics[scale=0.26]{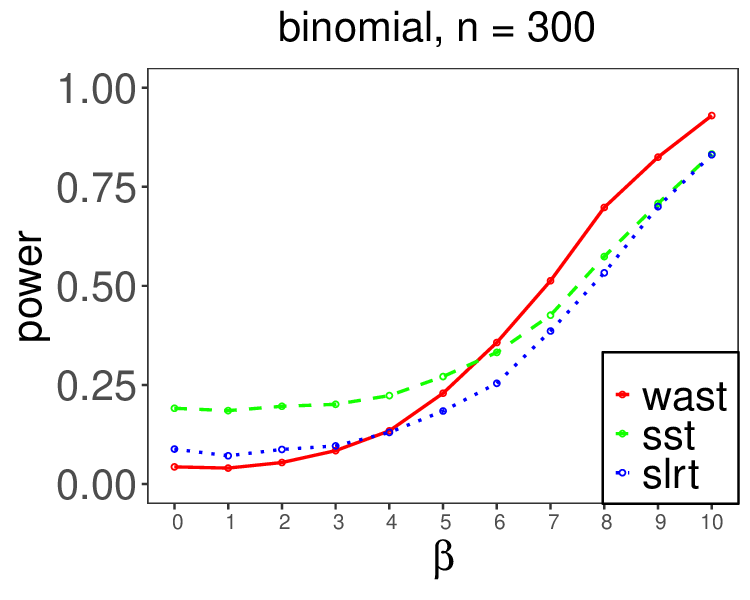}
		\includegraphics[scale=0.26]{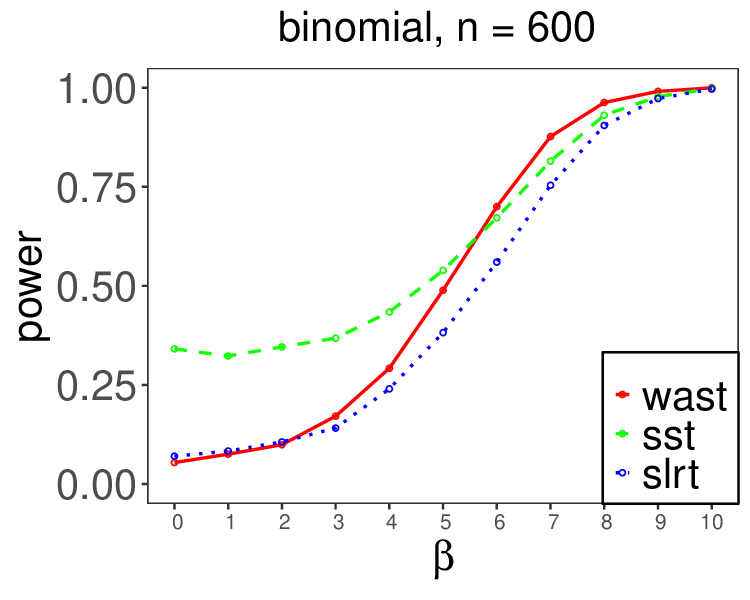}  \\
		\includegraphics[scale=0.26]{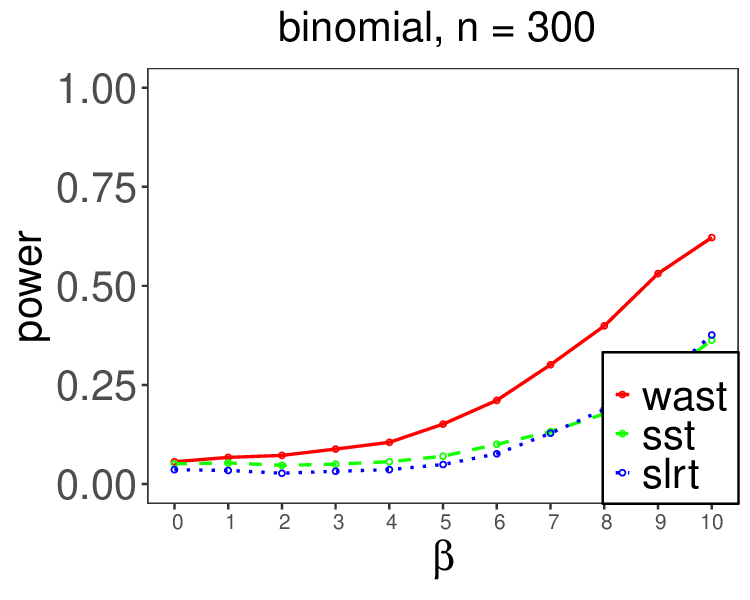}
		\includegraphics[scale=0.26]{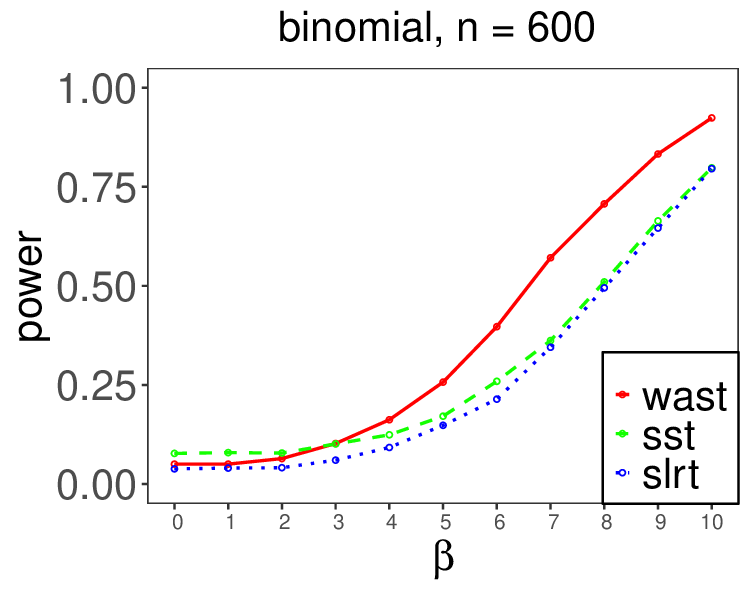}
		\includegraphics[scale=0.26]{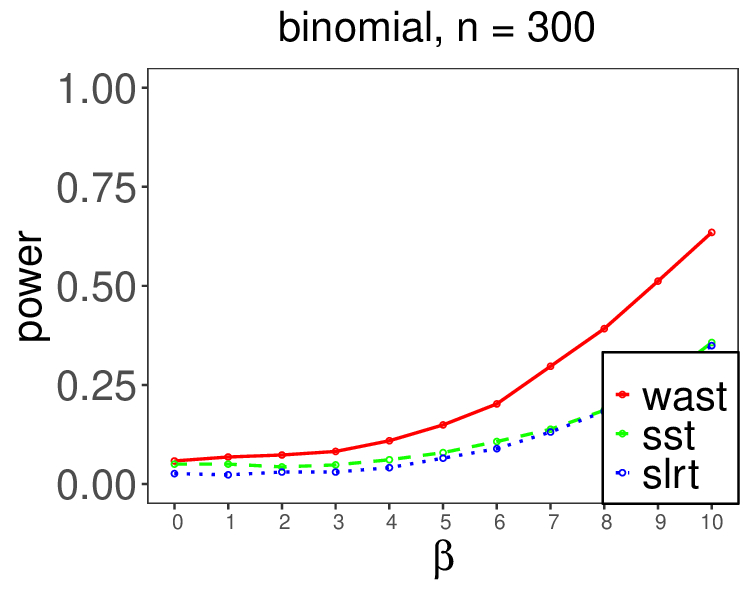}
		\includegraphics[scale=0.26]{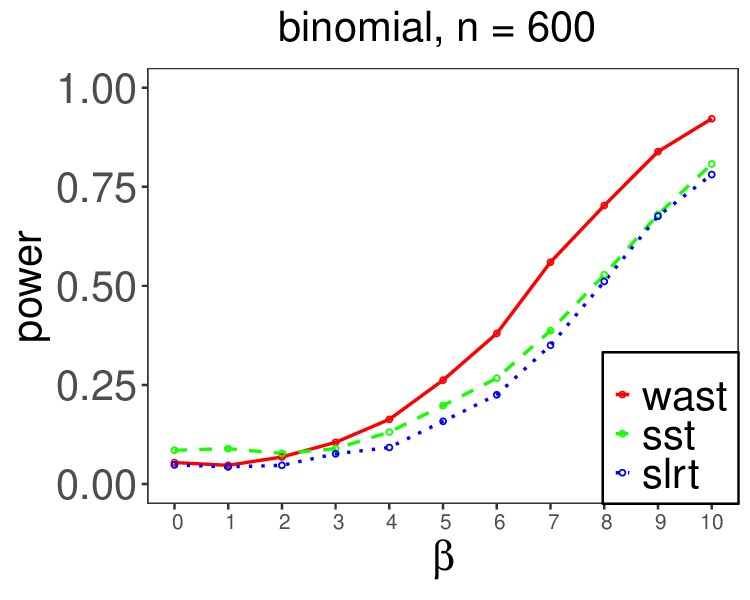}
		\caption{\it Powers of test statistic for GLM with binomial and with large numbers of dense $\bZ$ by the proposed WAST (red solid line), SST (green dashed line) and SLRT (blue dotted line). From top to bottom, each row depicts the powers for $p=2, 6, 11$. From left to right, each column depicts the powers for $q=100, 100, 500, 500$.
		}
		\label{fig_binomial_dense}
	\end{center}
\end{figure}

\begin{figure}[!ht]
	\begin{center}
		\includegraphics[scale=0.26]{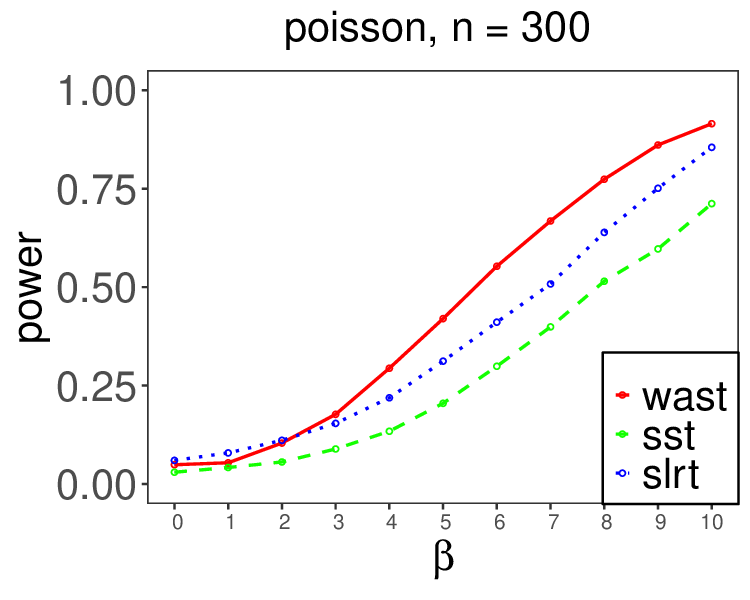}
		\includegraphics[scale=0.26]{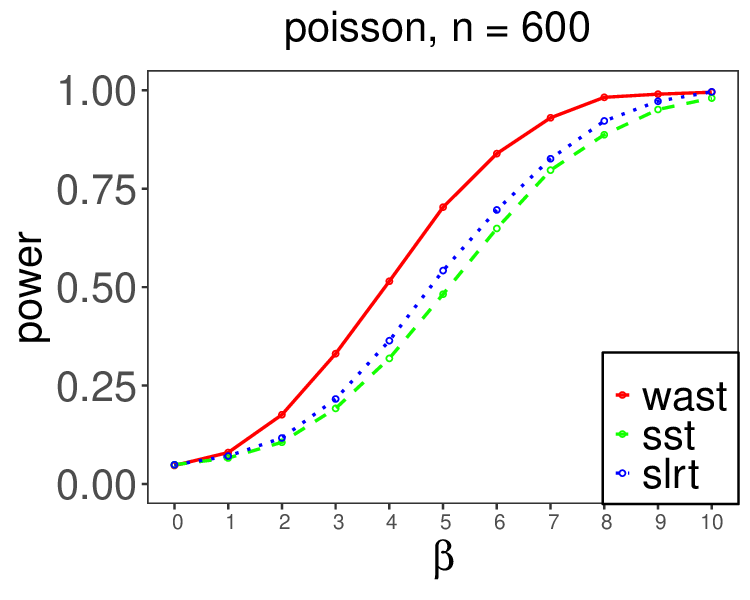}
		\includegraphics[scale=0.26]{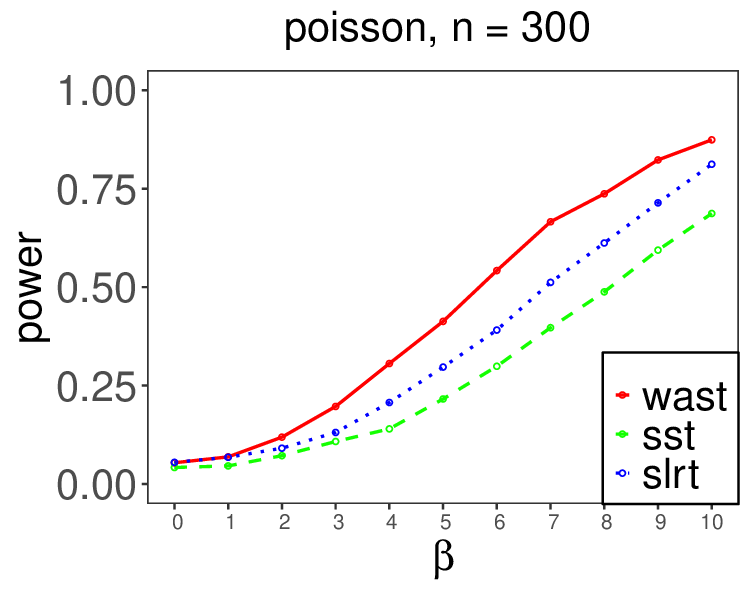}
		\includegraphics[scale=0.26]{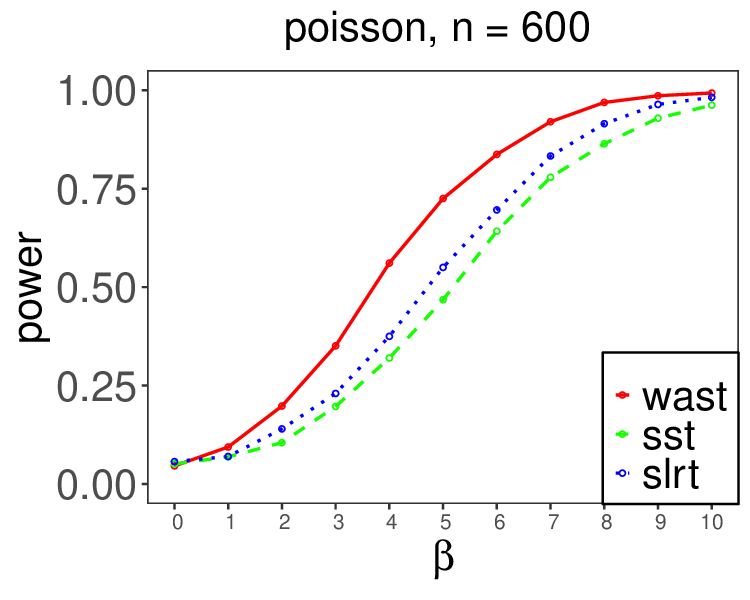}  \\
		\includegraphics[scale=0.26]{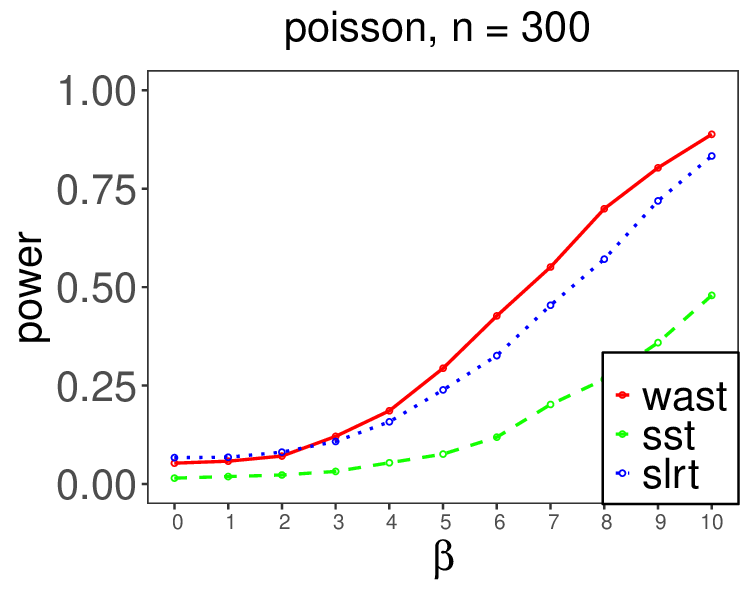}
		\includegraphics[scale=0.26]{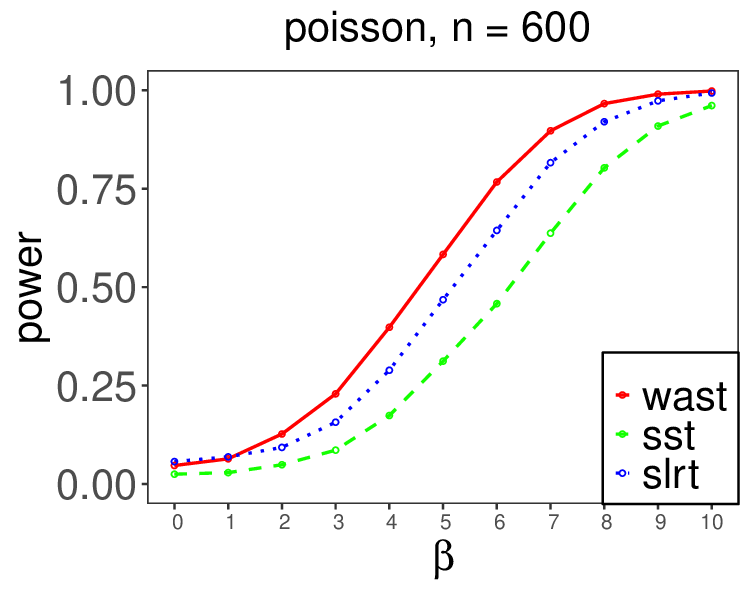}
		\includegraphics[scale=0.26]{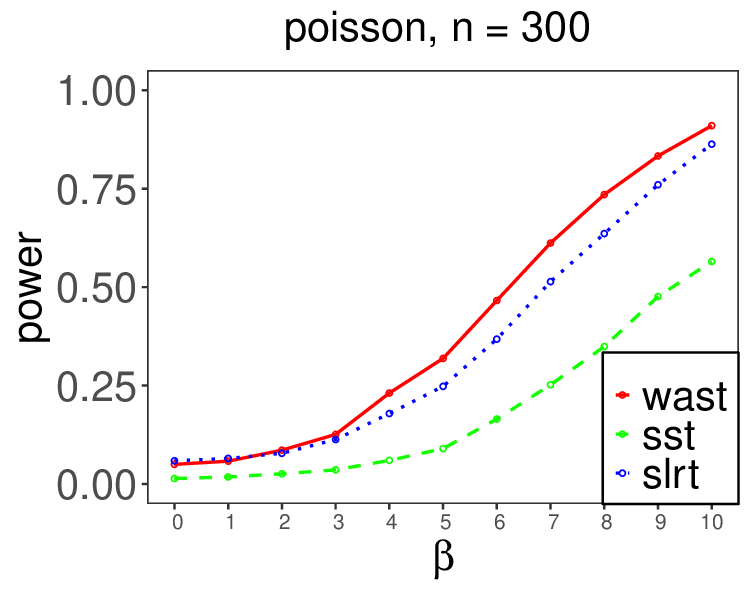}
		\includegraphics[scale=0.26]{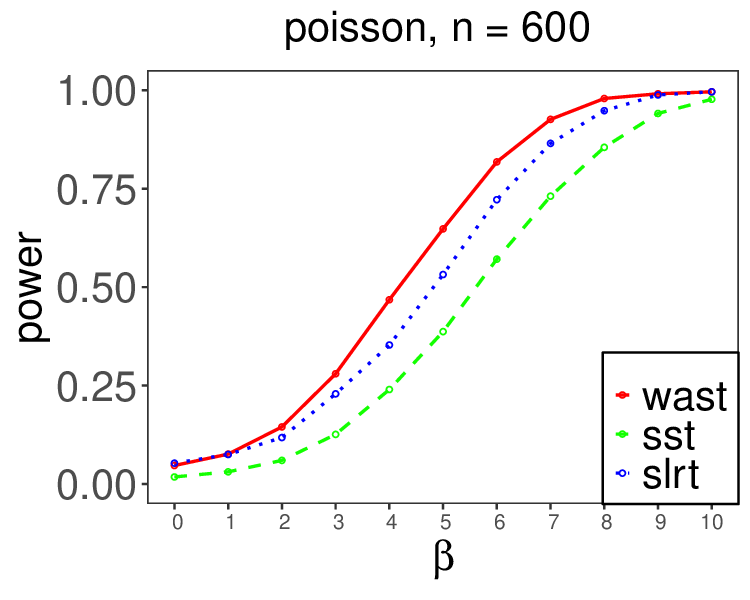}  \\
		\includegraphics[scale=0.26]{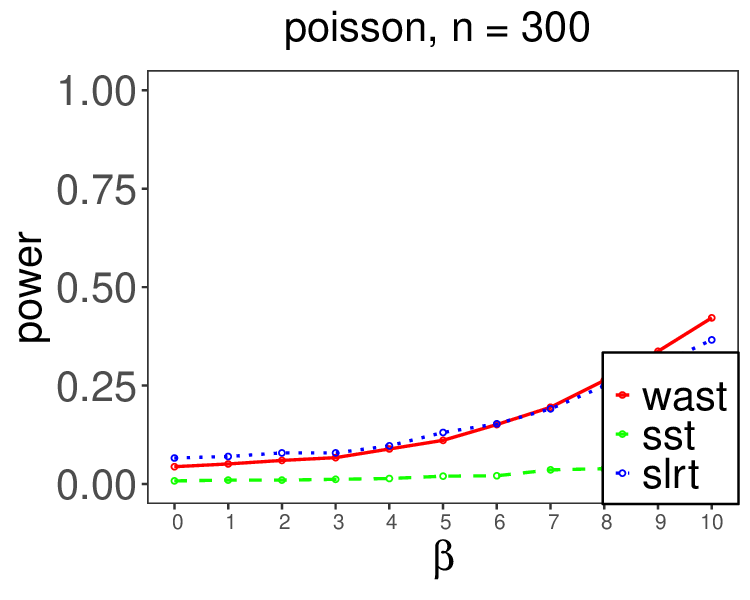}
		\includegraphics[scale=0.26]{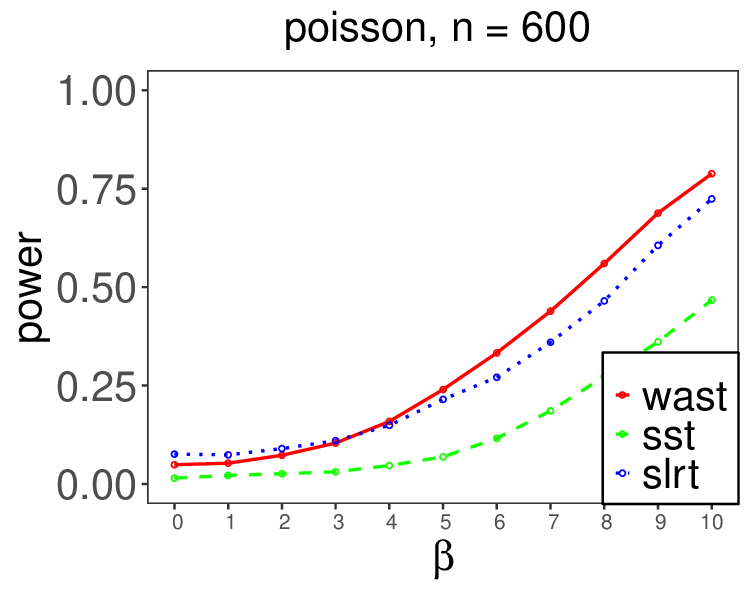}
		\includegraphics[scale=0.26]{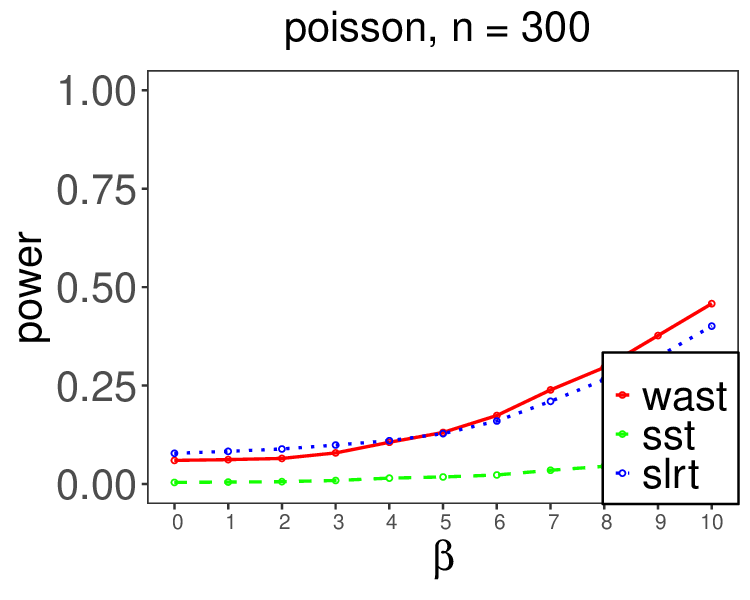}
		\includegraphics[scale=0.26]{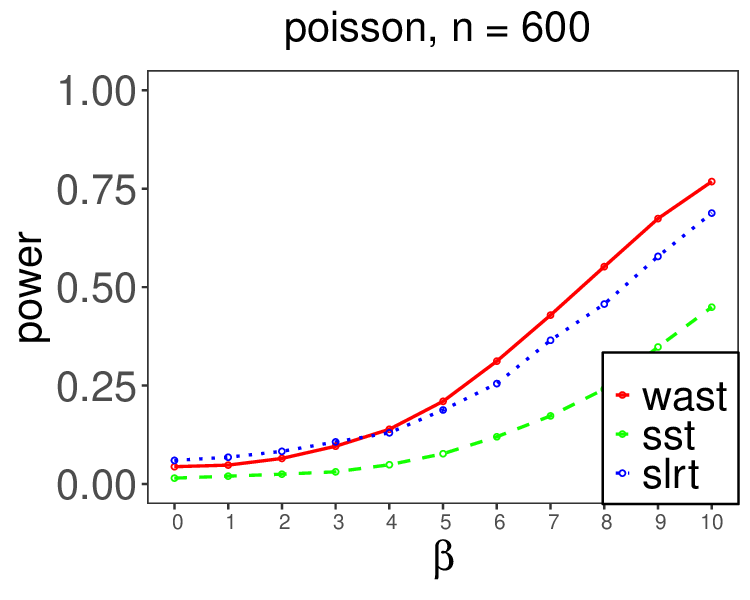}
		\caption{\it Powers of test statistic for GLM with Poisson and with large numbers of dense $\bZ$ by the proposed WAST (red solid line), SST (green dashed line) and SLRT (blue dotted line). From top to bottom, each row depicts the powers for $p=2, 6, 11$. From left to right, each column depicts the powers for $q=100, 100, 500, 500$.
		}
		\label{fig_poisson_dense}
	\end{center}
\end{figure}

Type \uppercase\expandafter{\romannumeral1} errors ($\kappa=0$) for probit, quantile and semiparametric models are listed in Table \ref{table_size_glm_dense}. We can see from Table \ref{table_size_glm_dense} that the size of the proposed WAST are close to the nominal significance level $0.05$, but for most scenarios the size of the SST are much smaller than 0.05. From Figure \ref{fig_gaussian_dense}-\ref{fig_poisson_dense} we have same conclusion as these in Section \ref{simulation_glm_sparse}. We omit the detailed analysis here.

\begin{table}[htp!]
	\def~{\hphantom{0}}
    \tiny
	\caption{Type \uppercase\expandafter{\romannumeral1} errors of the proposed WAST, SST and SLRT  based on resampling for GLMs with large numbers of dense $\bZ$. The nominal significant level is 0.05.
}
	\resizebox{\textwidth}{!}{
    \begin{threeparttable}
		\begin{tabular}{llcccccccc}
			\hline
			\multirow{2}{*}{Model}&\multirow{2}{*}{$(p,q)$}
			&\multicolumn{3}{c}{ $n=300$} && \multicolumn{3}{c}{ $n=600$} \\
			\cline{3-5} \cline{7-9}
			&& WAST& SST & SLRT&& WAST& SST & SLRT\\
			\cline{3-9}
			Gaussian &$(2,100)$         & 0.043 & 0.028 & 0.066 && 0.054 & 0.028 & 0.044 \\
			&$(2,500)$                  & 0.045 & 0.024 & 0.048 && 0.048 & 0.035 & 0.040 \\
			&$(6,100)$                  & 0.052 & 0.017 & 0.056 && 0.053 & 0.020 & 0.056 \\
			&$(6,500)$                  & 0.055 & 0.017 & 0.052 && 0.055 & 0.026 & 0.044 \\
			&$(11,100)$                 & 0.055 & 0.005 & 0.042 && 0.042 & 0.010 & 0.048 \\
			&$(11,500)$                 & 0.058 & 0.005 & 0.052 && 0.043 & 0.012 & 0.053 \\
			[1 ex]
			binomial &$(2,100)$          &0.051 & 0.067 & 0.038 && 0.040 & 0.078 & 0.048 \\\
			&$(2,500)$                  & 0.050 & 0.081 & 0.033 && 0.044 & 0.070 & 0.044 \\
			&$(6,100)$                  & 0.056 & 0.169 & 0.088 && 0.041 & 0.303 & 0.070 \\
			&$(6,500)$                  & 0.043 & 0.191 & 0.088 && 0.054 & 0.341 & 0.070 \\
			&$(11,100)$                 & 0.056 & 0.051 & 0.036 && 0.050 & 0.077 & 0.038 \\
			&$(11,500)$                 & 0.058 & 0.050 & 0.026 && 0.054 & 0.085 & 0.048 \\
			[1 ex]
			Poisson &$(2,100)$          & 0.049 & 0.030 & 0.060 && 0.047 & 0.049 & 0.049 \\
			&$(2,500)$                  & 0.054 & 0.042 & 0.055 && 0.046 & 0.051 & 0.057 \\
			&$(6,100)$                  & 0.053 & 0.015 & 0.067 && 0.047 & 0.025 & 0.057 \\
			&$(6,500)$                  & 0.050 & 0.014 & 0.059 && 0.047 & 0.018 & 0.053 \\
			&$(11,100)$                 & 0.044 & 0.008 & 0.066 && 0.049 & 0.015 & 0.076 \\
			&$(11,500)$                 & 0.060 & 0.004 & 0.078 && 0.044 & 0.015 & 0.060 \\
			\hline
		\end{tabular}
\end{threeparttable}
	}
	\label{table_size_glm_dense}
\end{table}

\subsection{Comparison of the WAST }
\begin{figure}[!ht]
	\begin{center}
		\includegraphics[scale=0.4]{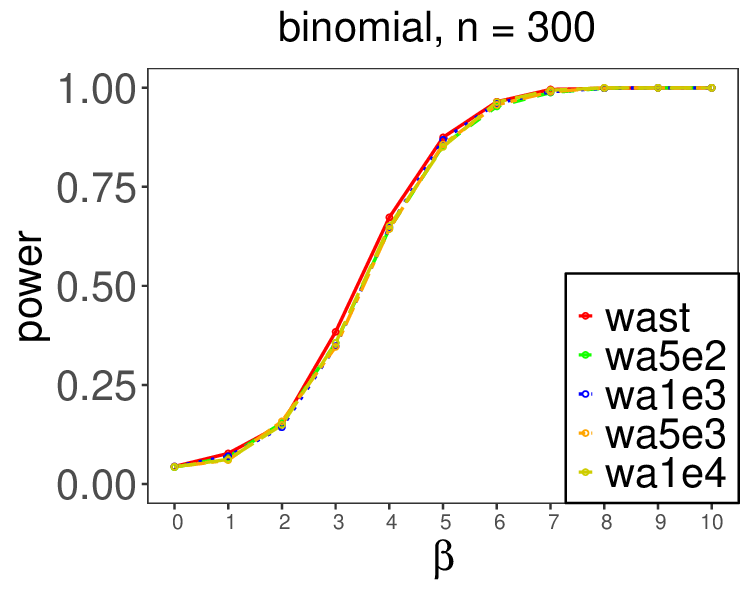}
		\includegraphics[scale=0.4]{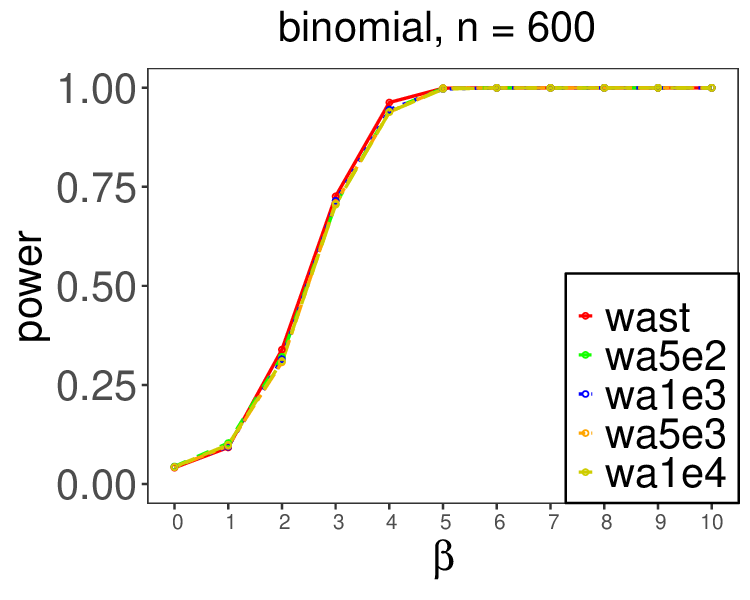}
		\includegraphics[scale=0.4]{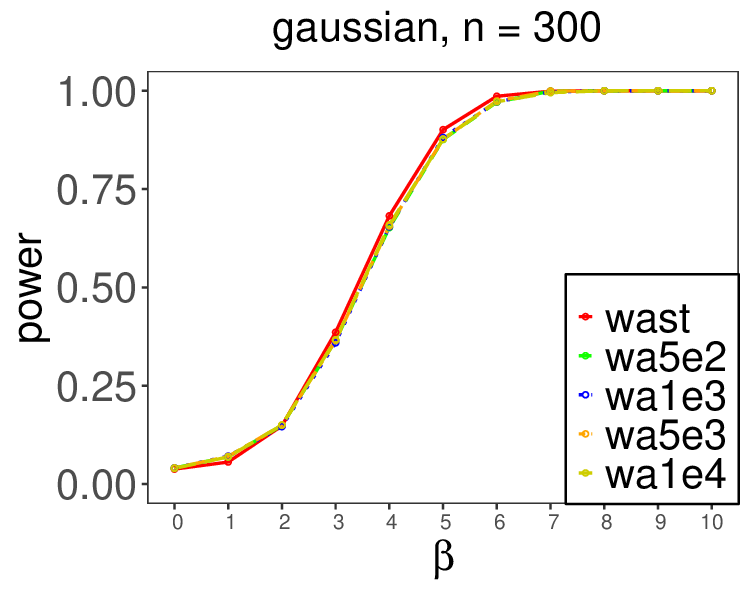}
		\includegraphics[scale=0.4]{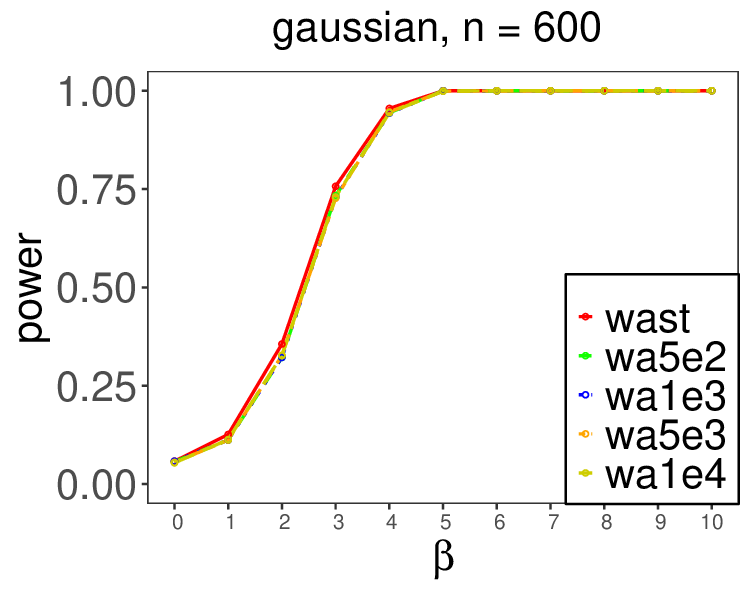}
		\includegraphics[scale=0.4]{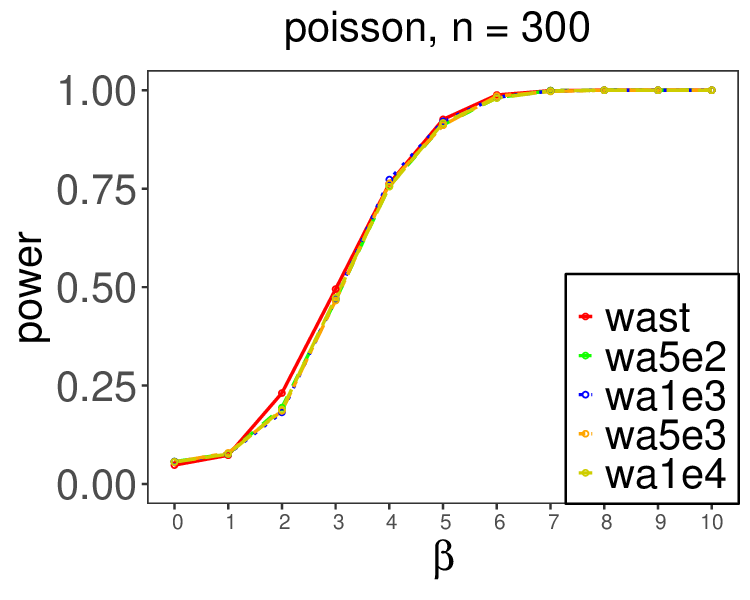}
		\includegraphics[scale=0.4]{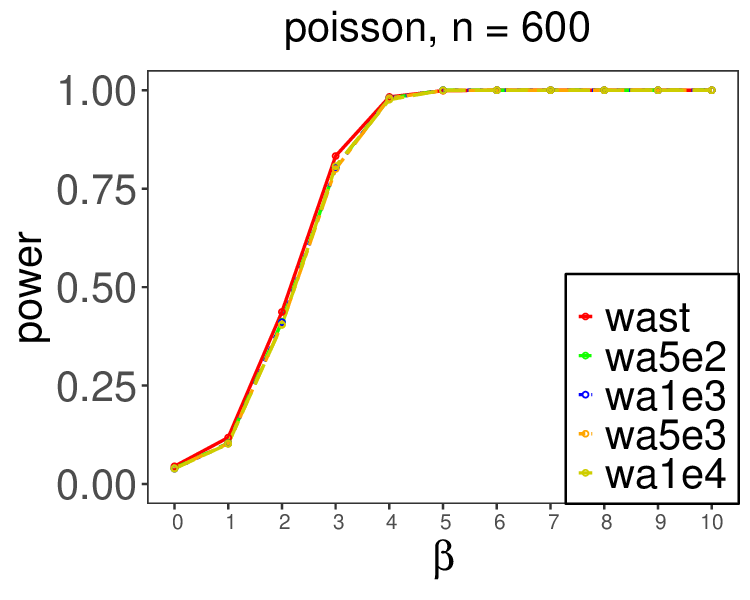}
		\includegraphics[scale=0.4]{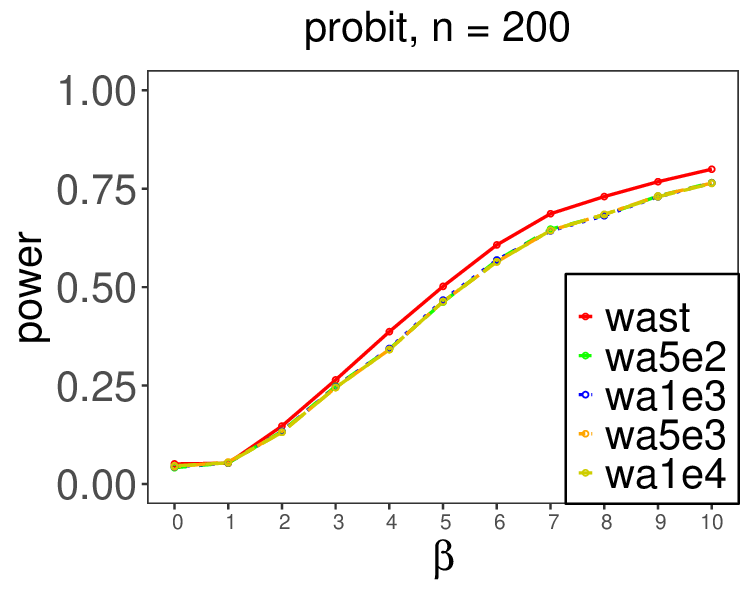}
		\includegraphics[scale=0.4]{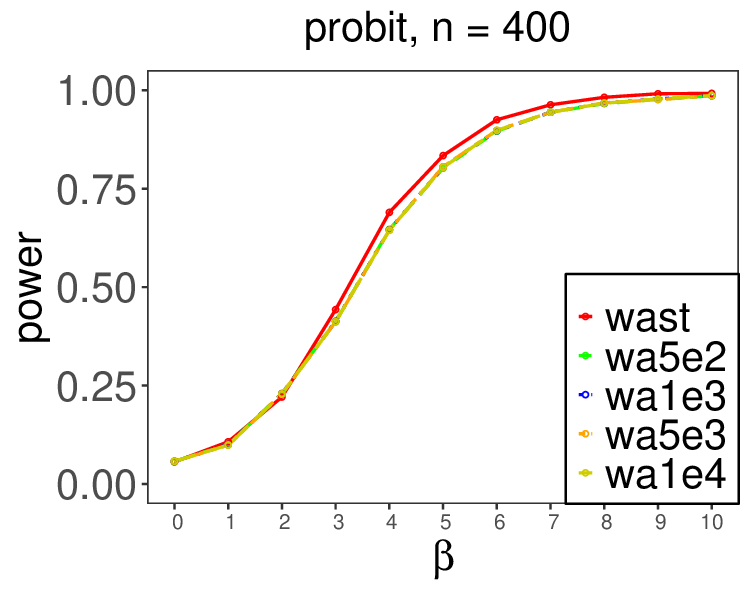}
		\includegraphics[scale=0.4]{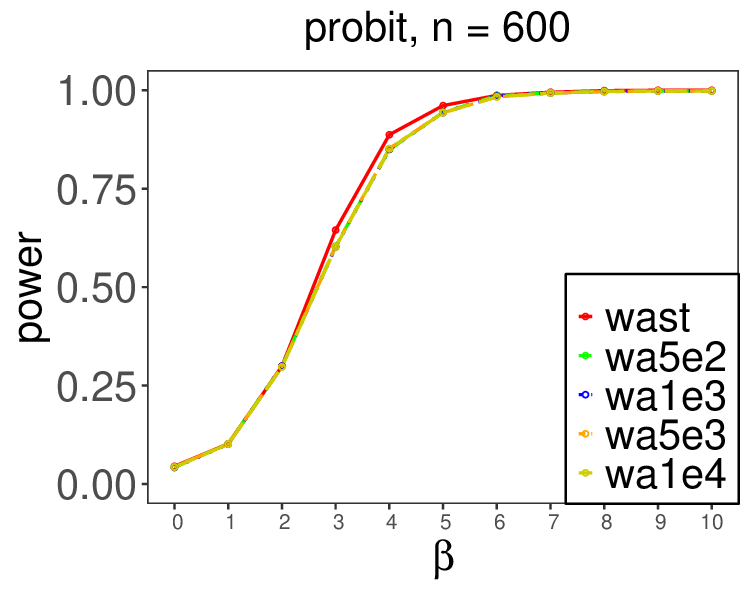}
		\includegraphics[scale=0.4]{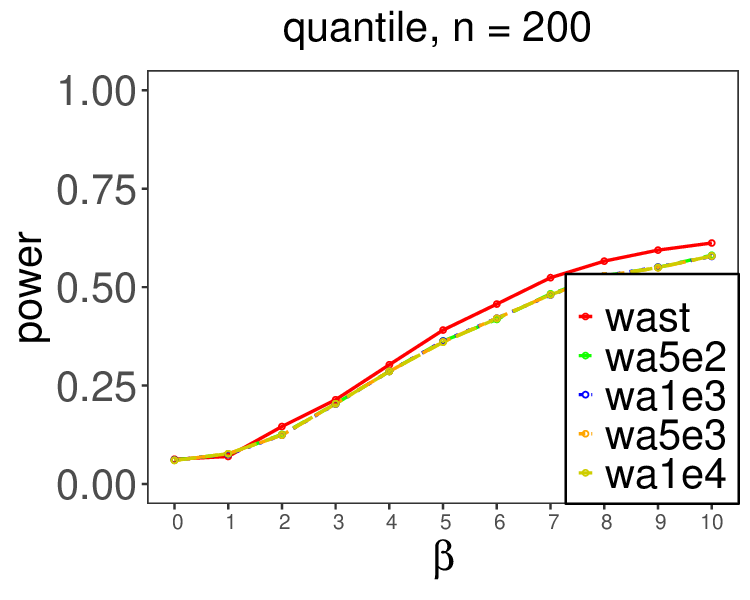}
		\includegraphics[scale=0.4]{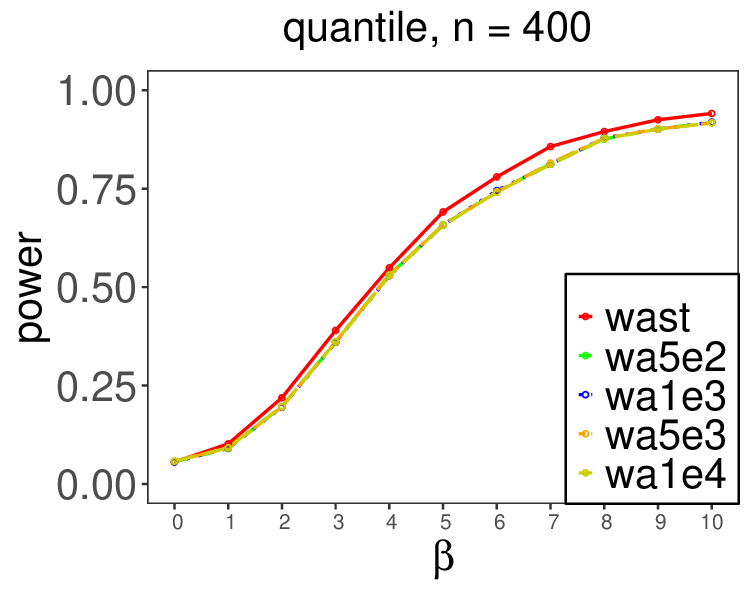}
		\includegraphics[scale=0.4]{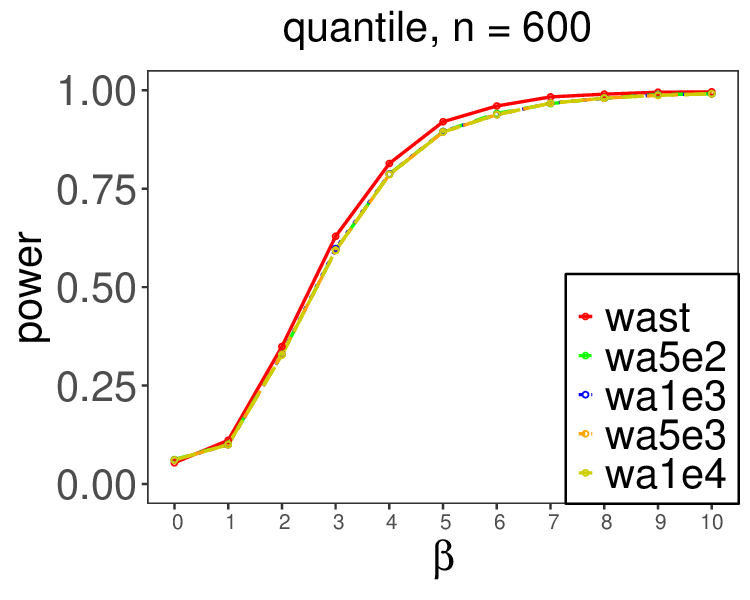}
		\includegraphics[scale=0.4]{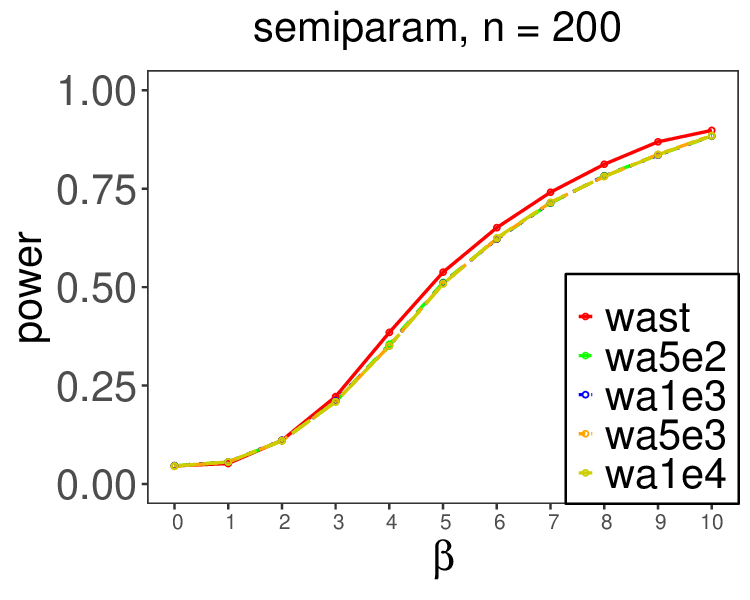}
		\includegraphics[scale=0.4]{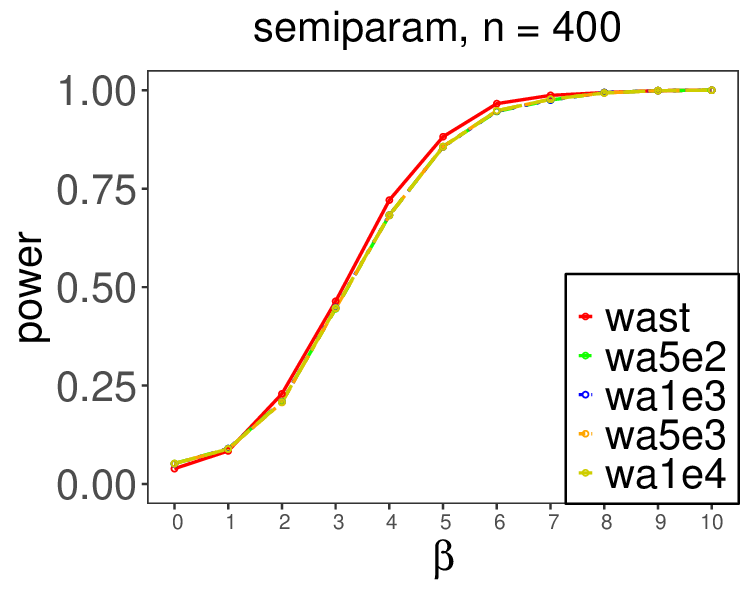}
		\includegraphics[scale=0.4]{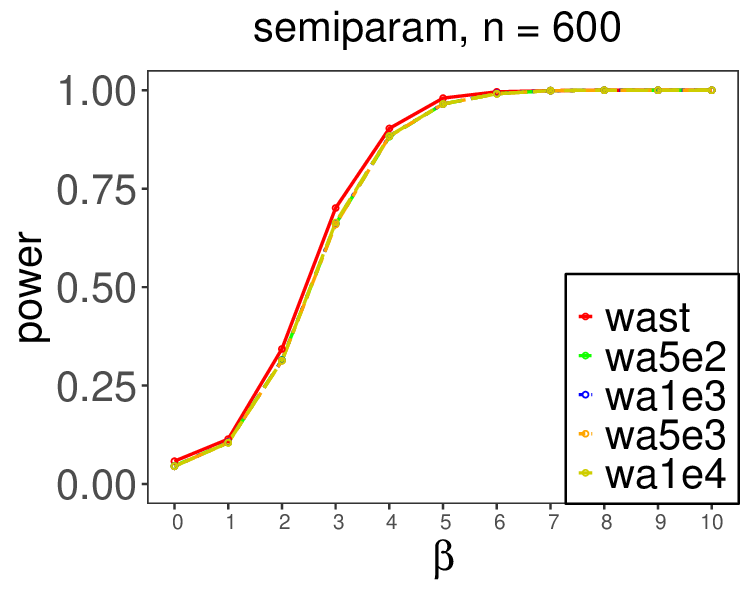}
		\caption{Powers of test statistics by WAST (red solid line) and WAST with approximation (``wa5e2" denotes $N=500$, ``wa1e3" denotes $N=1000$, ``wa5e3" denotes $N=5000$, ``wa1e4" denotes $N=10000$). }
		\label{fig_comparison}
	\end{center}
\end{figure}

In this section, we compare the WAST with close form (15) and with approximation (14). The simulation settings are same as that in Section 5 of main paper, in which we set $(r, p, q)=(2, 11, 100)$ and apply sparse scenario with the number of nonzero parameters being 10. We set $N=500, 1000, 5000, 10000$ in the formula (14) of the main paper. Figure \ref{fig_comparison} indicate that the WAST has slightly greater power than the WAST with approximation for each distribution and sample size. The number $N$ of the generated $z$s is not sensitive for the WAST with approximation when $N\geq500$.

\clearpage

\section{Simulation Studies for small \texorpdfstring{$q$}{}}\label{simulations_smallq}
\setcounter{figure}{0}
\setcounter{table}{0}
\def\thetable{C.\arabic{table}}
\def\thefigure{C.\arabic{figure}}
\def\thethm{C.\arabic{thm}}

\subsection{Single change plane for GLMs}
For fair comparison, we use the same settings as considered by \cite{2021Threshold}.
Consider the GLMs with the canonical parameter
\begin{align*}
\mu = \tbX\trans\balpha + \bX\trans\bbeta\bone(\bZ\trans\btheta\geq0),
\end{align*}
where $\bbeta = \kappa\bone_p$ and $\alpha_2=\cdots=\alpha_r=\log(1.4)$. For Gaussian and Poisson families,
we set $\alpha_1 = 0.5$,  and for binomial family $\alpha_1$ is chosen so that the proportion of cases in the data set is ${1}/{3}$ on average under the null hypothesis. $\theta_2,\cdots,\theta_q$ are equally spaced numbers from $-1$ to $1$, and $\theta_1$ is chosen as the negative of the 0.65 percentile of $Z_2\theta_2+\cdots+Z_q\theta_q$, which means that $\bZ\trans\btheta$ divides the population into two groups with 0.35 and 0.65 observations, respectively. For the predictor, we generate $(v_2, \cdots,v_{\max\{r,p\}}, Z_2,\cdots,Z_q)\trans$ from a multivariate normal distribution with mean $\bzero$ and covariance $\Sigma=(\tilde{\rho}_{ij})$, where $\tilde{\rho}_{ij}=1$ if $i=j$ and $\tilde{\rho}_{ij}=\rho$ otherwise. Here, we consider both $\rho=0$ and $\rho=0.3$. We set $\tX_j=\bone(v_j>0)$ with $j=2,\cdots,r$, $X_k = \bone(v_k>0)$ with $k=2,\cdots,p$, and $\tX_1=1$, $X_1=1$, and $Z_1=1$.

We evaluate the power under a sequence of alternative models indexed by $\kappa$, i.e., $H_1^{\kappa}: \bbeta = \kappa\bone_p$ with $\kappa$ equally spaced in the range $(0,\kappa_{\max}]$. We set the sample size as $n=(300, 600)$, and we set 1000 repetitions and 1000 bootstrap samples. For comparison, we consider three methods here: (i) the proposed WAST, (ii) the SST, and (iii) the approximated supremum of likelihood-ratio test statistic (SLRT) given by \cite{2021Threshold}. We calculate both SLRT and SST over $\{\btheta^{(k)}=(\theta^{(k)}_1,\cdots,\theta^{(k)}_q)\trans: k=1,\cdots,K\}$. Let $\btheta^{(k)}_{-1}=\tilde{\btheta}^{(k)}_{-1}/\|\tilde{\btheta}^{(k)}_{-1}\|$, where $\btheta^{(k)}_{-1}=(\theta^{(k)}_2,\cdots,\theta^{(k)}_q)\trans$ and $\tilde{\btheta}^{(k)}_{-1}=(\tilde{\theta}^{(k)}_2,\cdots,\tilde{\theta}^{(k)}_q)\trans$, and $\tilde{\btheta}^{(k)}_{-1}$ is drawn independently from the multivariate normal distribution $\mathcal{N}(\bzero_{q-1},\bone_{q-1})$. For each $\theta_1^{(k)}$, $k=1,\cdots,K$, we select it by an equal grid search in the range from the lower tenth percentile to the upper tenth percentile of the data points of $\{\theta^{(k)}_2Z_{i2}+\cdots+\theta^{(k)}_qZ_{iq}\}_{i=1}^n$, which is the same as that in \cite{2021Threshold}. Here, we set $K=1000$.

Table~\ref{table_size} lists the type-\uppercase\expandafter{\romannumeral1} errors ($\kappa=0$) with the nominal significance level of 0.05. As can be seen, the sizes of the proposed method are close to 0.05. For the SST, its sizes are far from 0.05 in most of the scenarios, and it performs the worst of these three methods. For the SLRT, its sizes are close to 0.05 except for $(r,p,q)=(2, 51,11)$ and $(r,p,q)=(6, 51,11)$.

\begin{figure}[!ht]
\begin{center}
\includegraphics[scale=0.3]{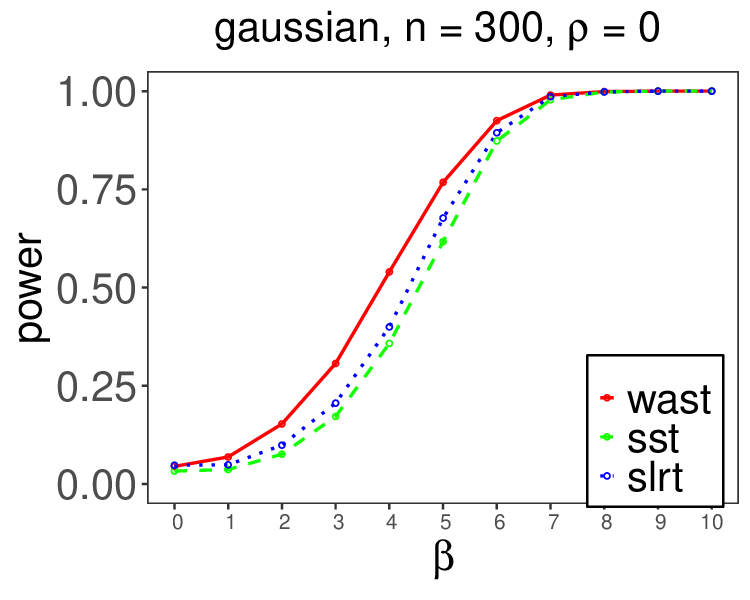}
\includegraphics[scale=0.3]{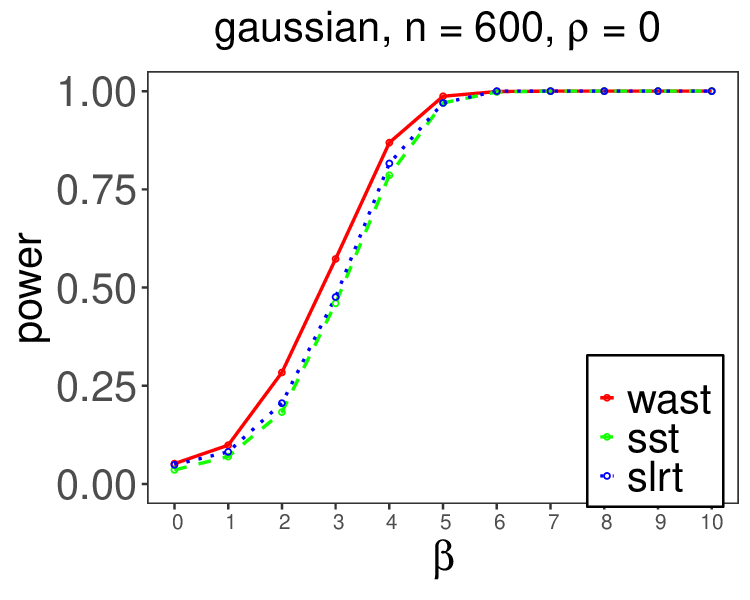}
\includegraphics[scale=0.3]{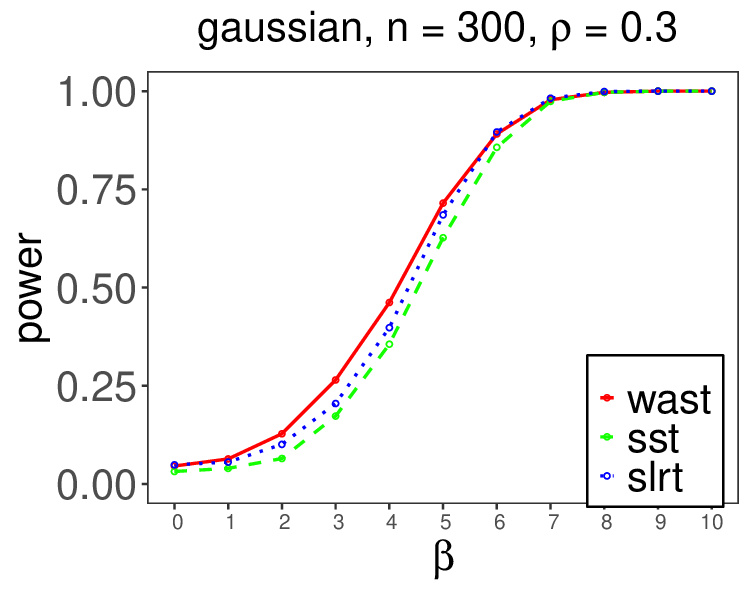}
\includegraphics[scale=0.3]{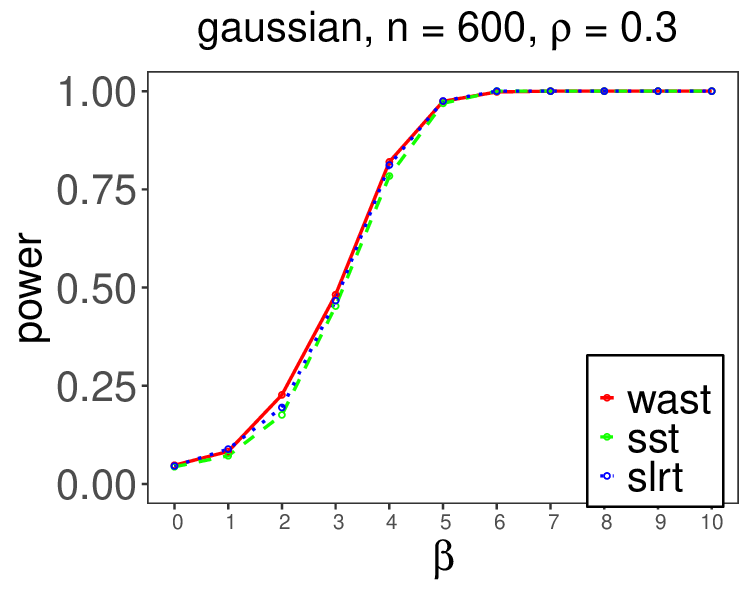}
\includegraphics[scale=0.3]{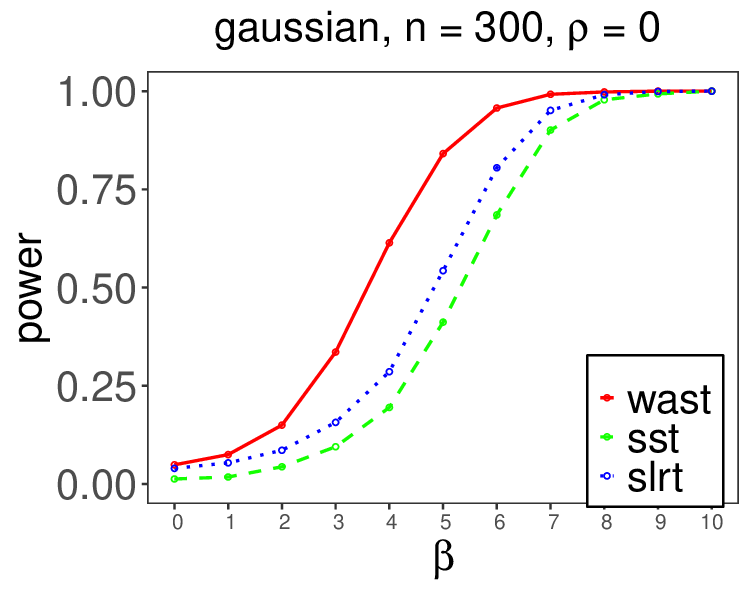}
\includegraphics[scale=0.3]{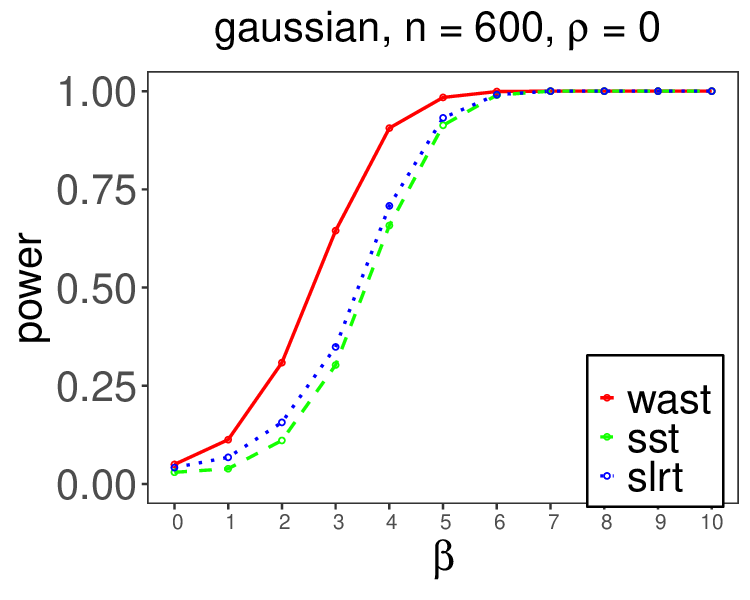}
\includegraphics[scale=0.3]{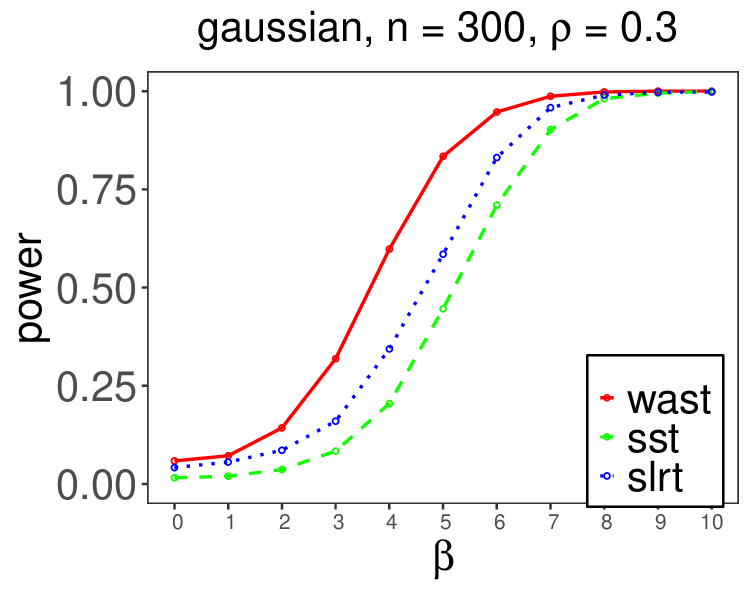}
\includegraphics[scale=0.3]{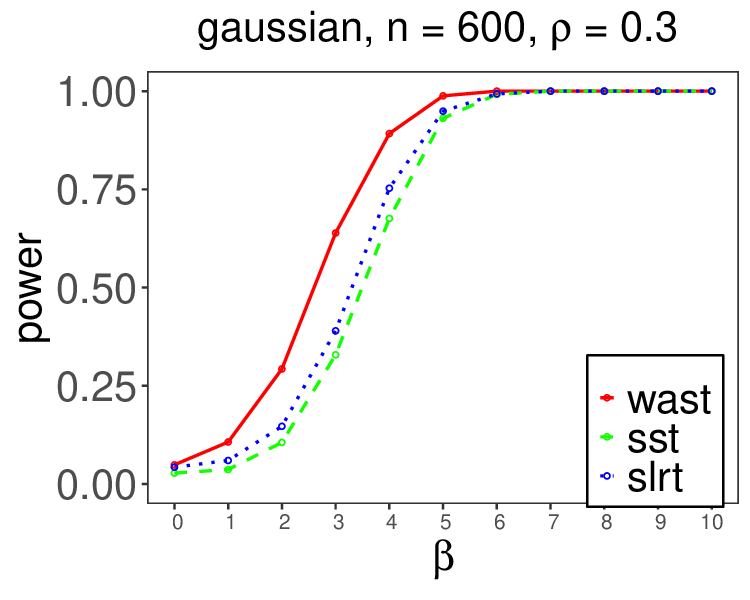}
\includegraphics[scale=0.3]{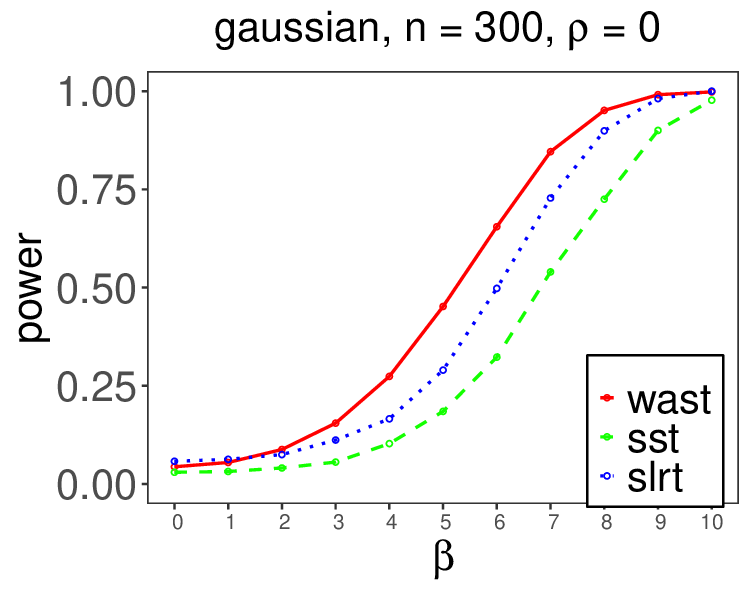}
\includegraphics[scale=0.3]{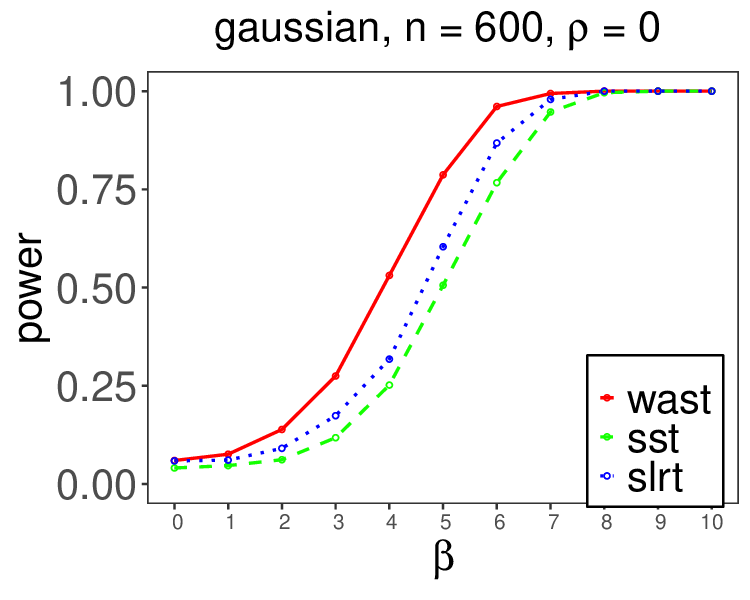}
\includegraphics[scale=0.3]{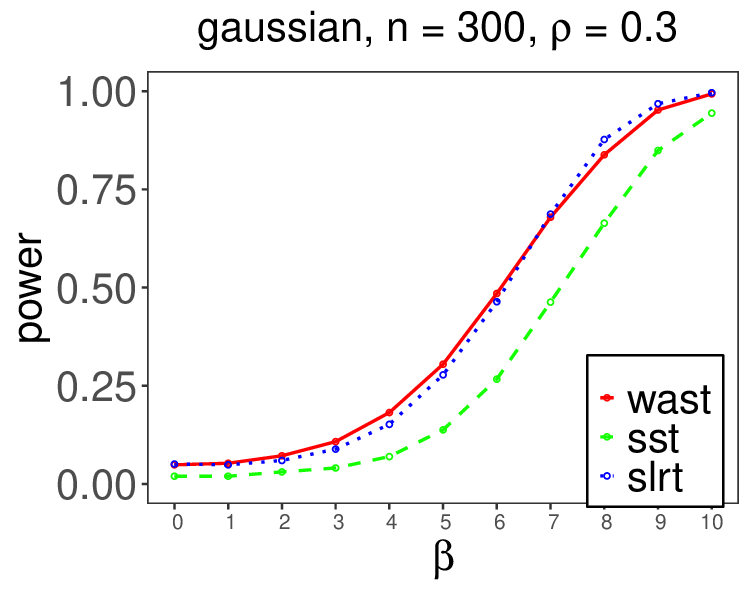}
\includegraphics[scale=0.3]{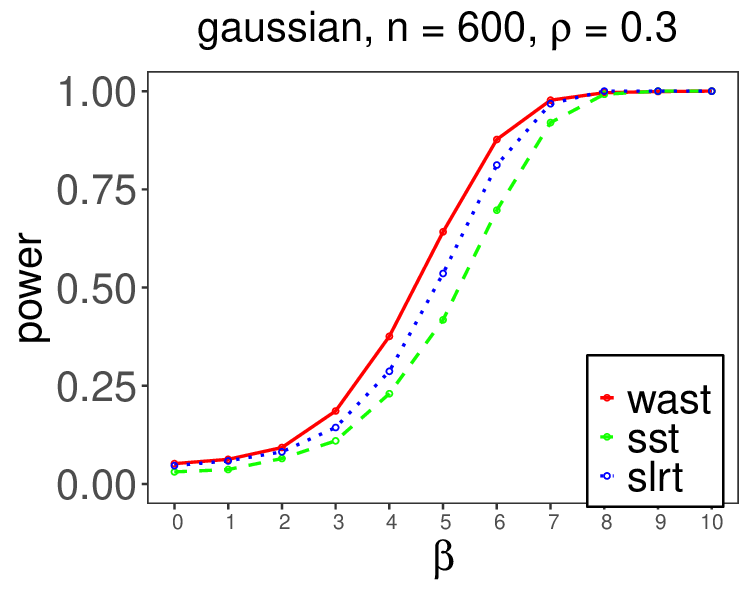}
\includegraphics[scale=0.3]{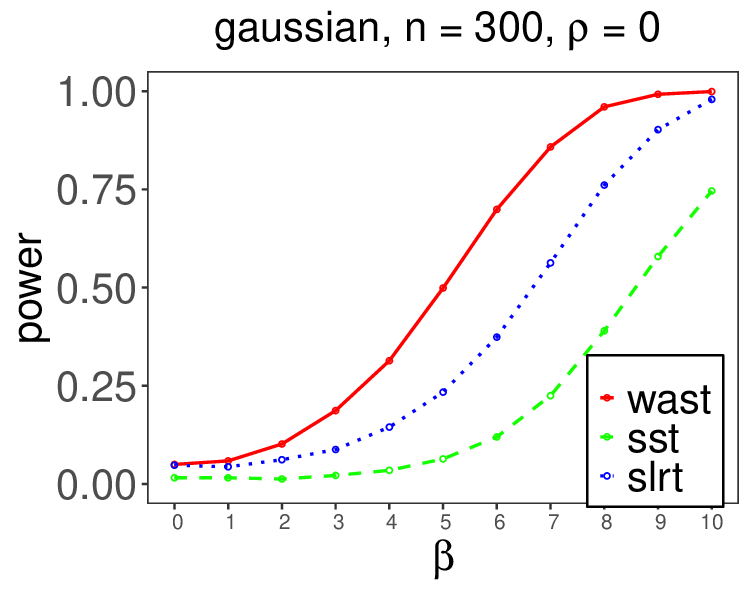}
\includegraphics[scale=0.3]{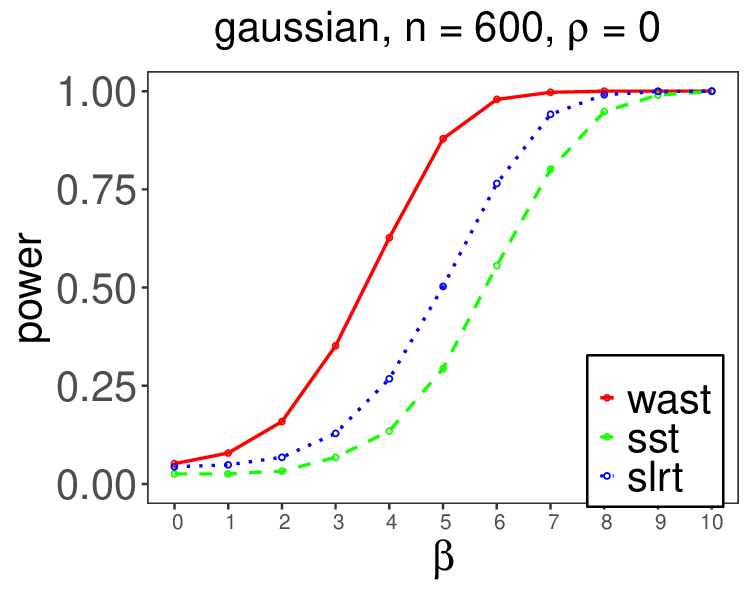}
\includegraphics[scale=0.3]{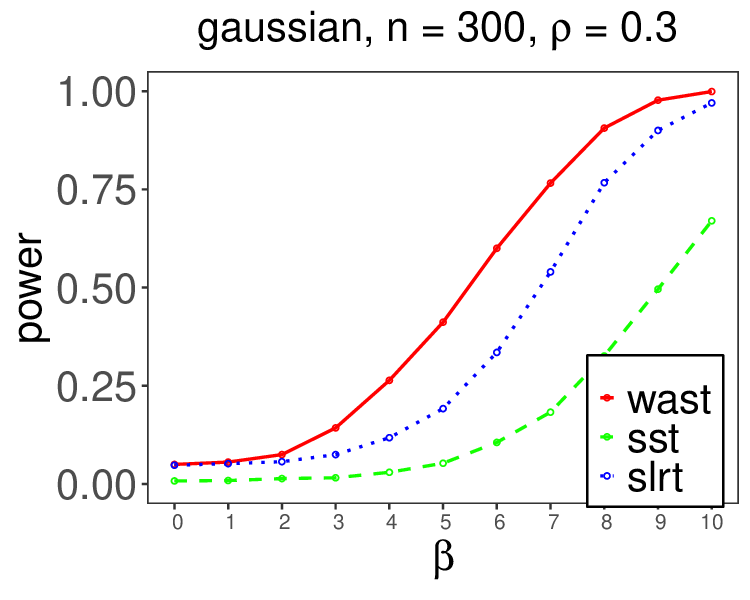}
\includegraphics[scale=0.3]{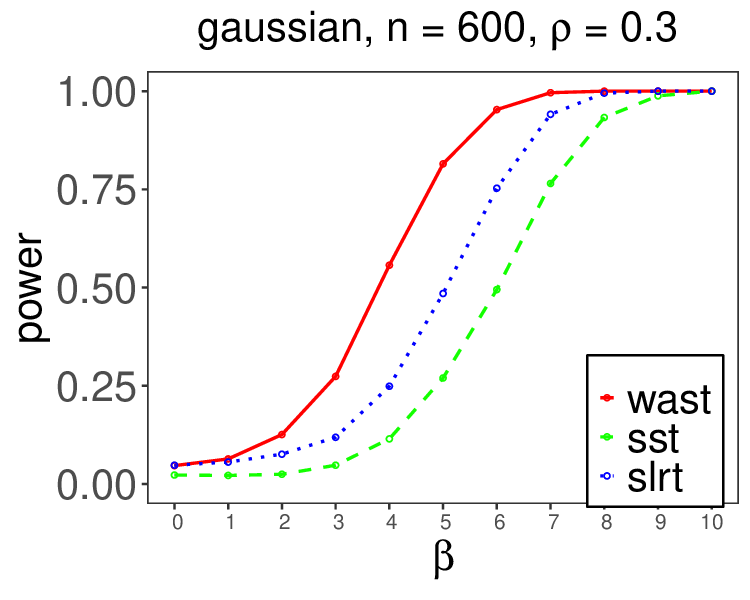}
\includegraphics[scale=0.3]{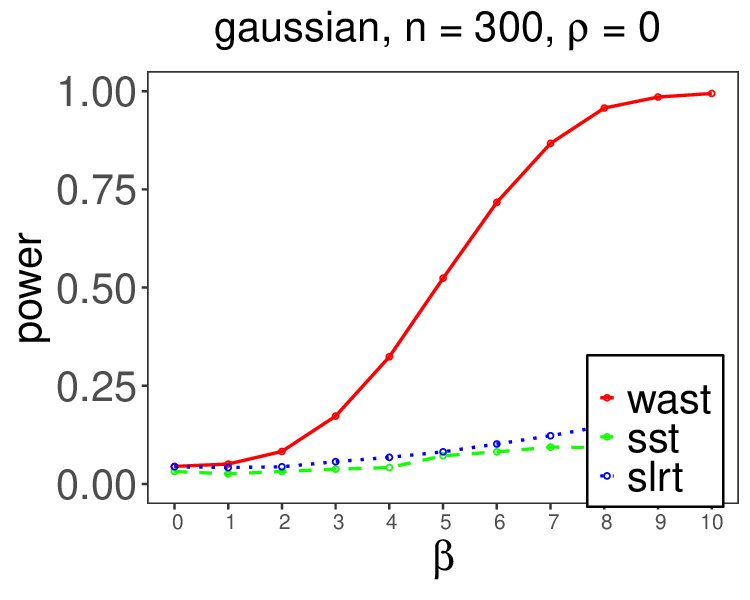}
\includegraphics[scale=0.3]{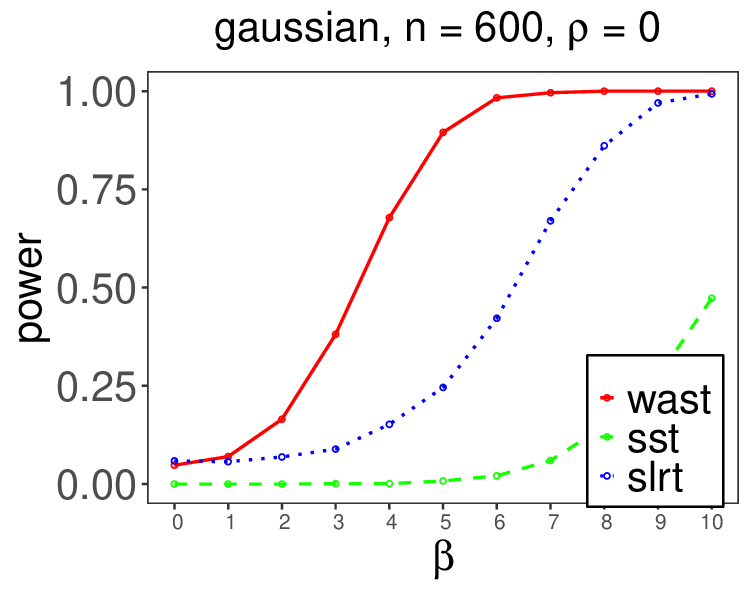}
\includegraphics[scale=0.3]{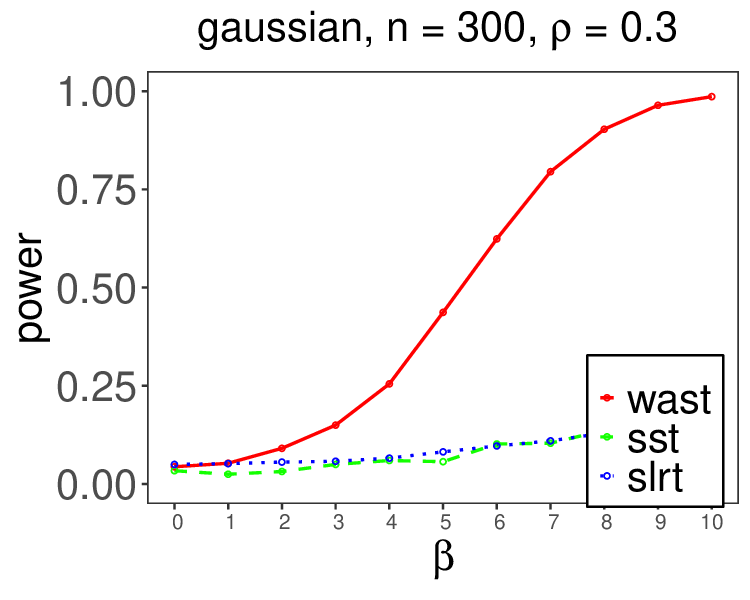}
\includegraphics[scale=0.3]{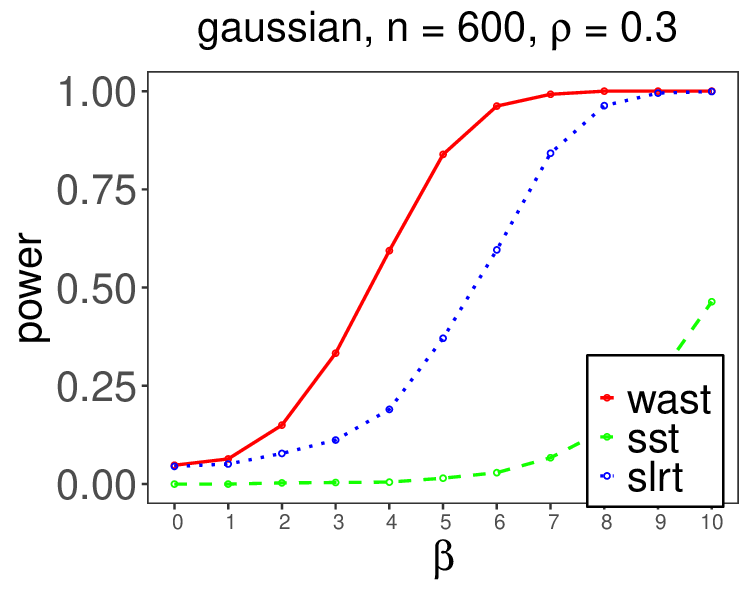}
\includegraphics[scale=0.3]{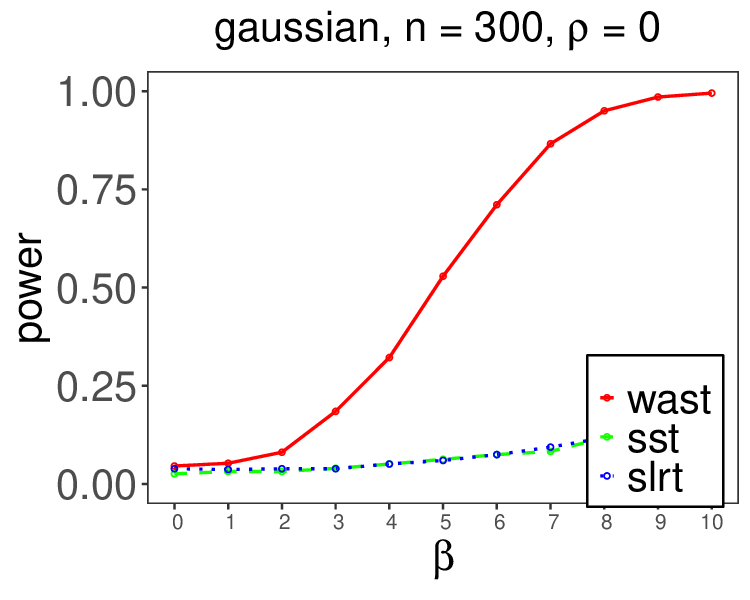}
\includegraphics[scale=0.3]{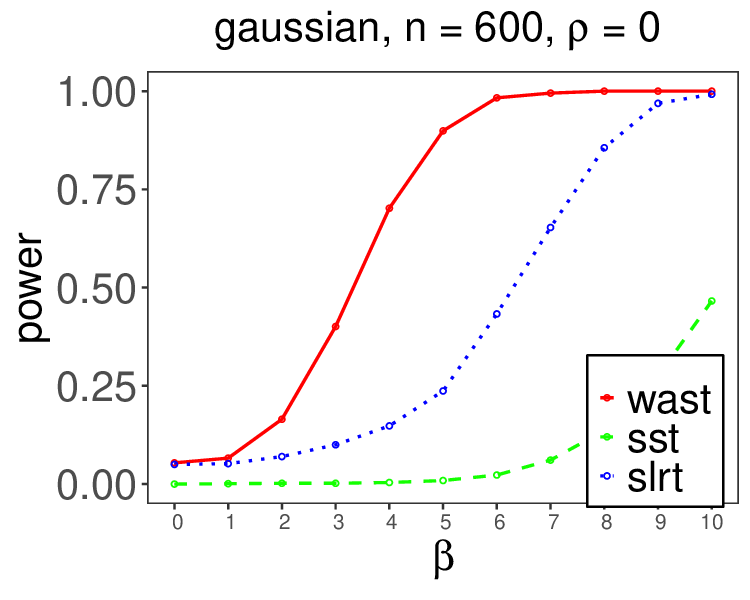}
\includegraphics[scale=0.3]{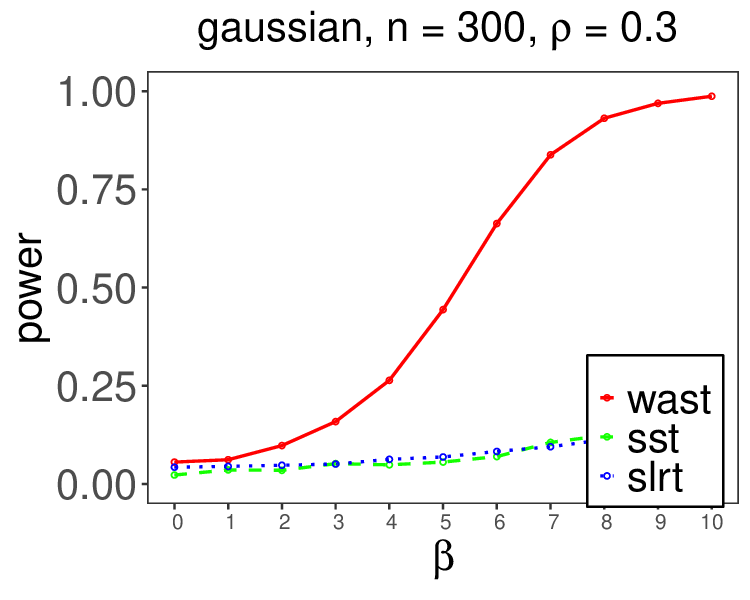}
\includegraphics[scale=0.3]{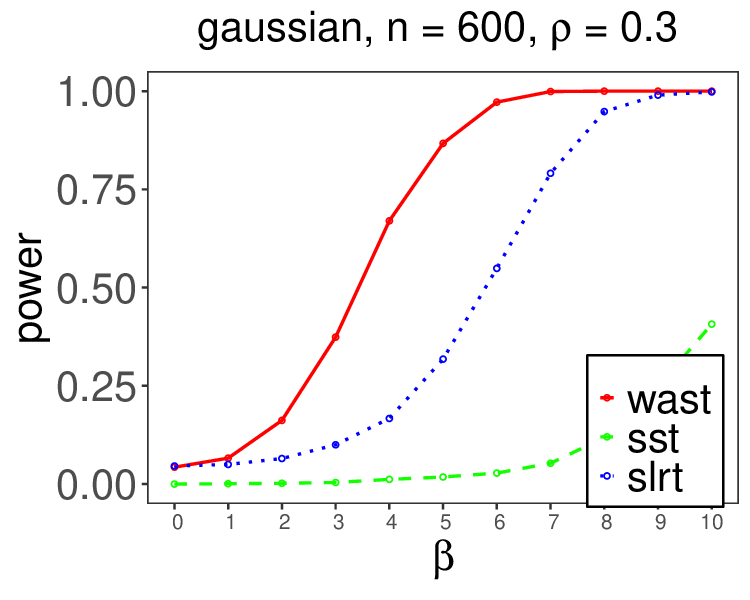}
\caption{Powers of testing linear model with Gaussian error by WAST (red solid line), SST (green dashed line), and SLRT (blue dotted line) for $n=(300,600)$. From top to bottom, each row depicts the powers for the cases $(r,p,q)=(2,2,3)$, $(6,6,3)$, $(2,2,11)$, $(6,6,11)$, $(2,51,11)$, and $(6,51,11)$.}
\label{fig_gaussian}
\end{center}
\end{figure}

\begin{figure}[!ht]
	\begin{center}
		\includegraphics[scale=0.285]{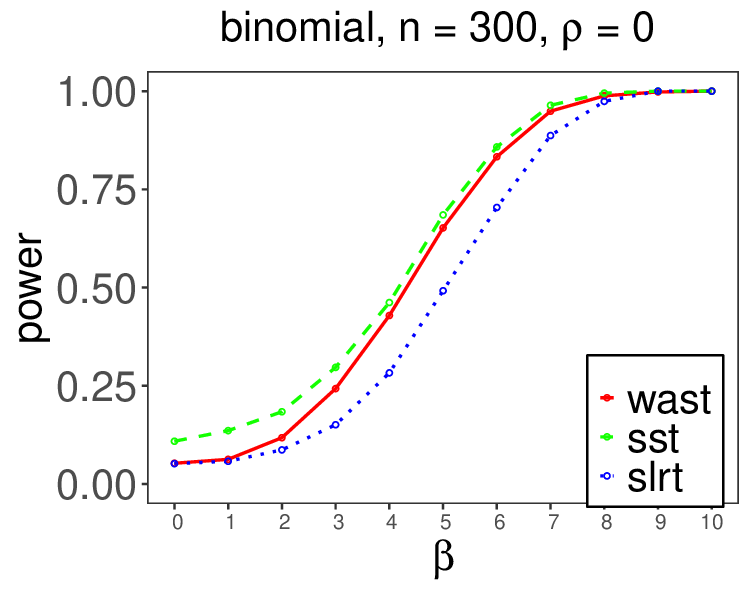}
		\includegraphics[scale=0.285]{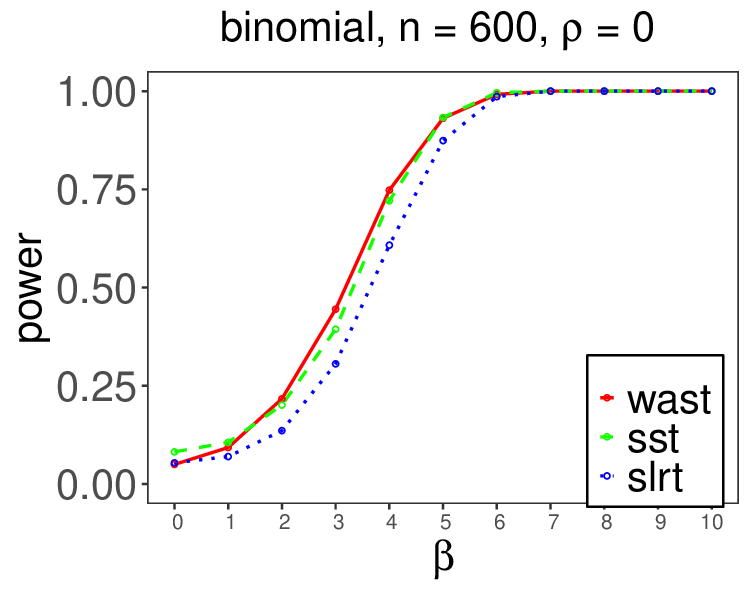}
		\includegraphics[scale=0.285]{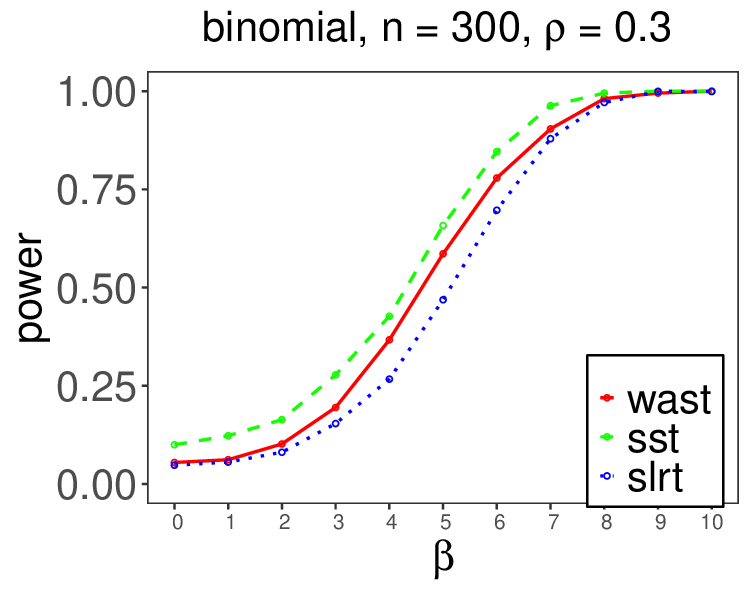}
		\includegraphics[scale=0.285]{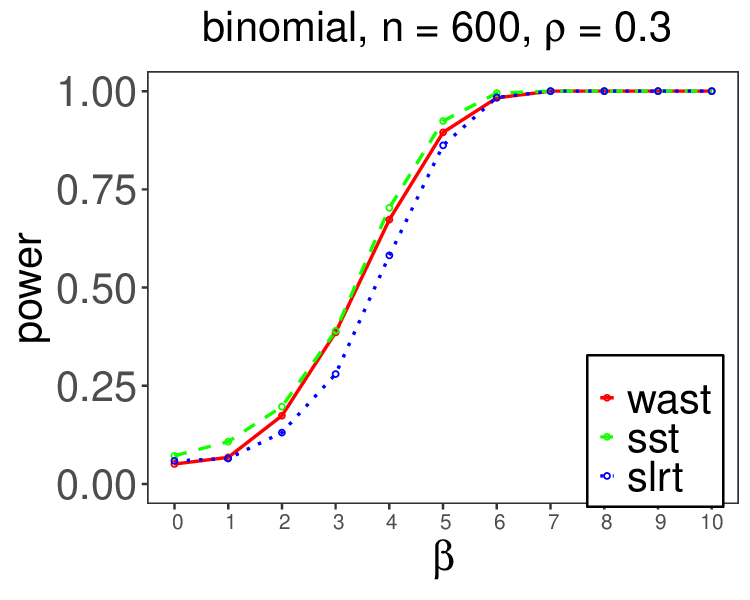}
		\includegraphics[scale=0.285]{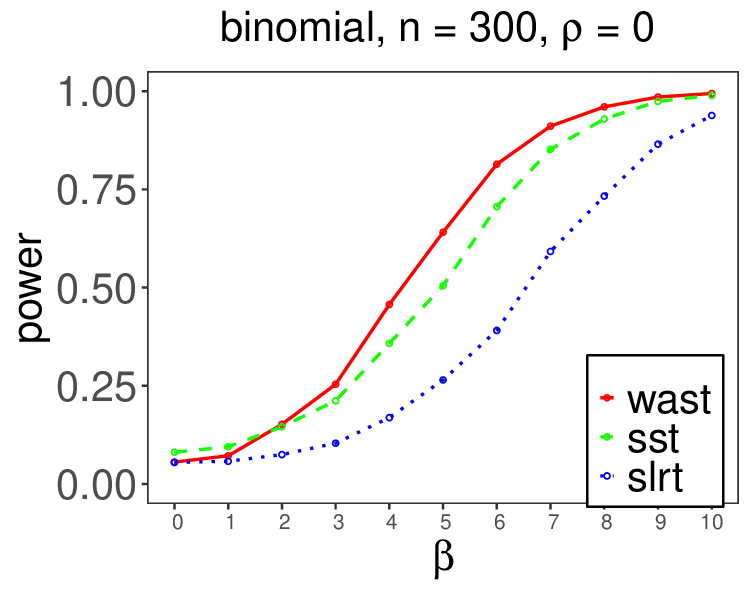}
		\includegraphics[scale=0.285]{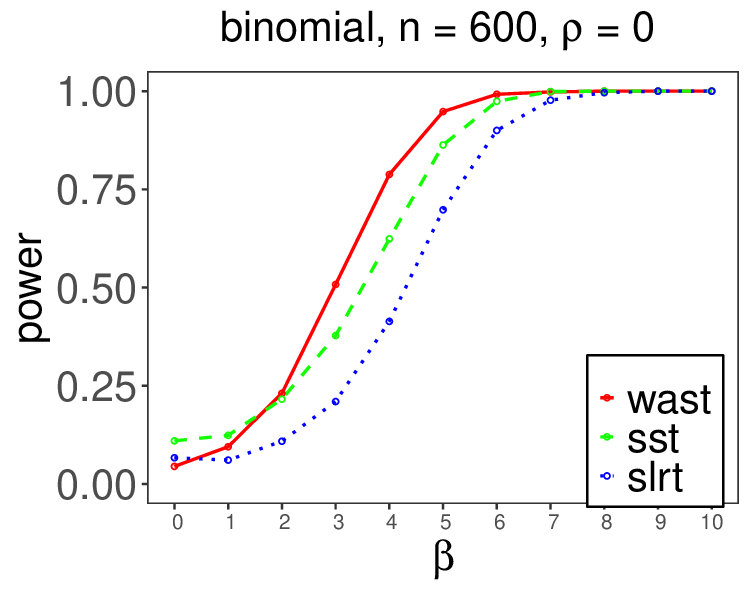}
		\includegraphics[scale=0.285]{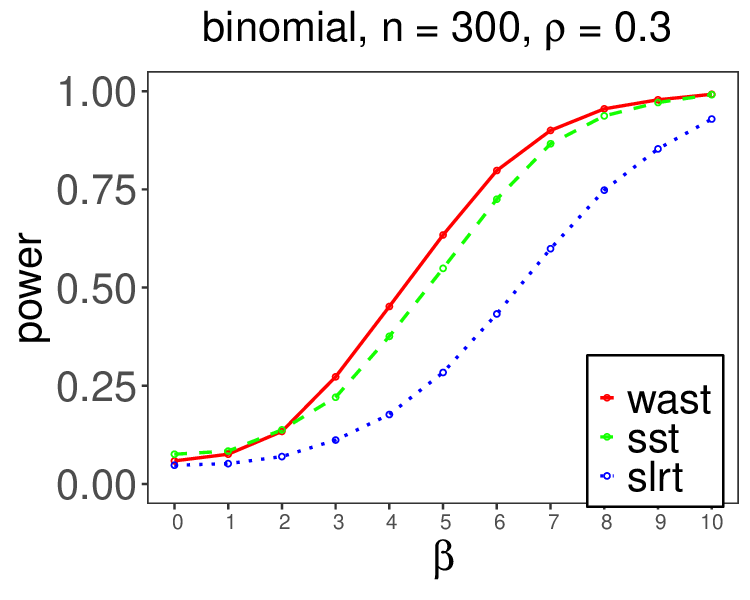}
		\includegraphics[scale=0.285]{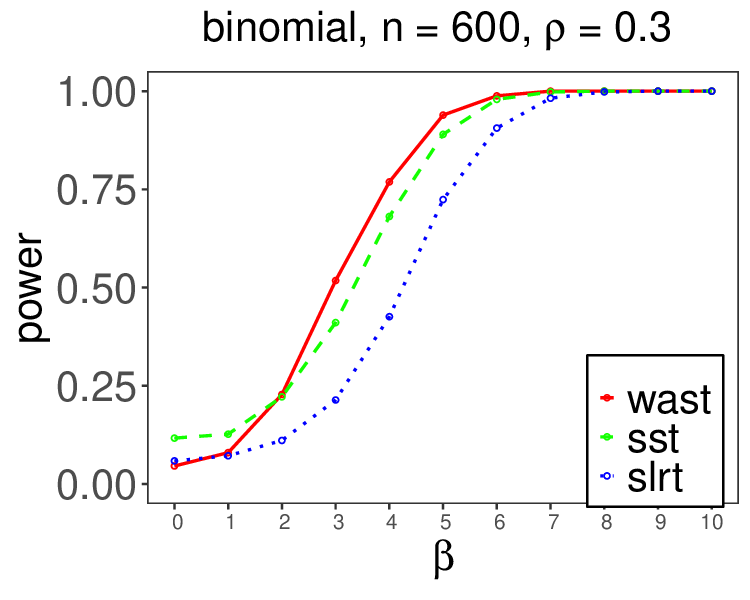}
		\includegraphics[scale=0.285]{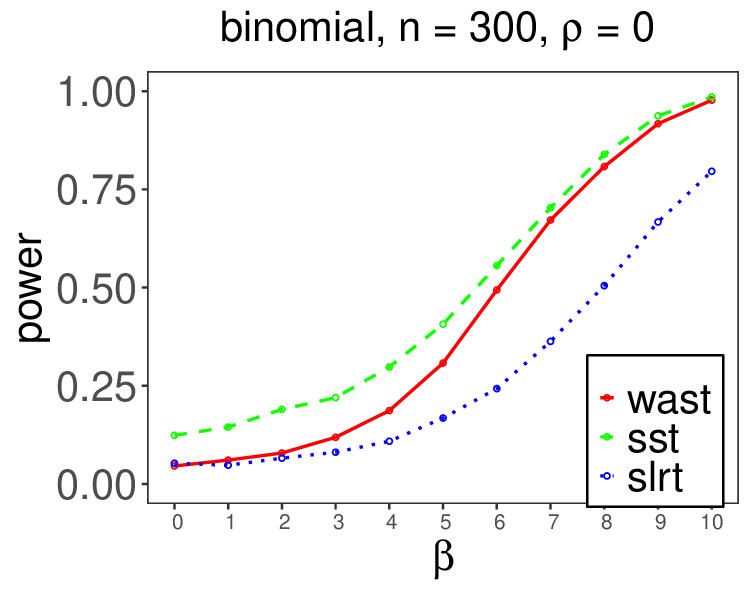}
		\includegraphics[scale=0.285]{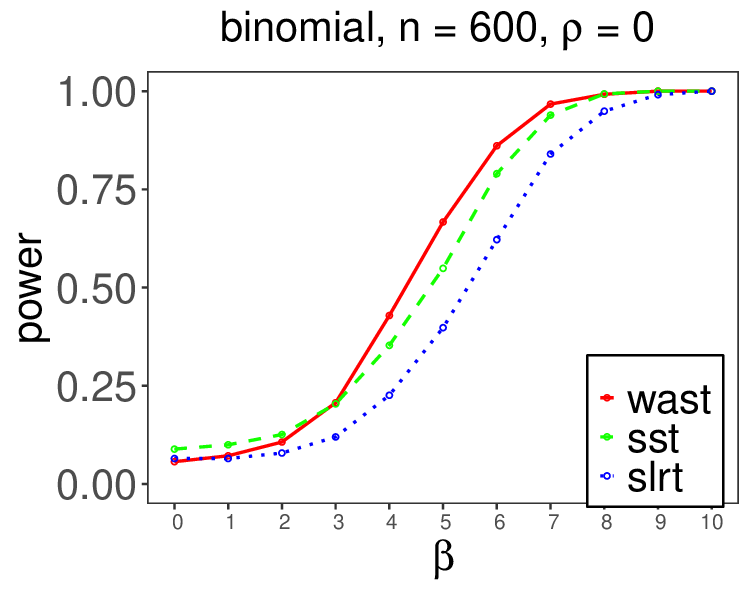}
		\includegraphics[scale=0.285]{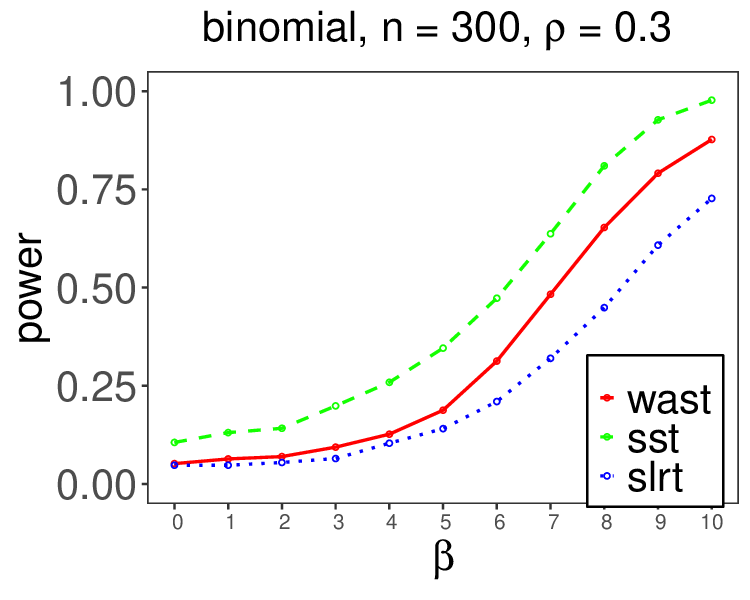}
		\includegraphics[scale=0.285]{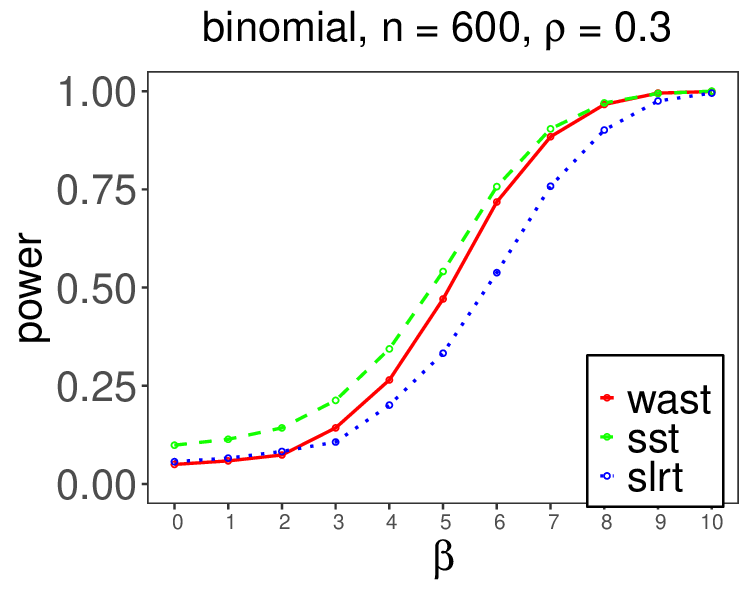}
		\includegraphics[scale=0.285]{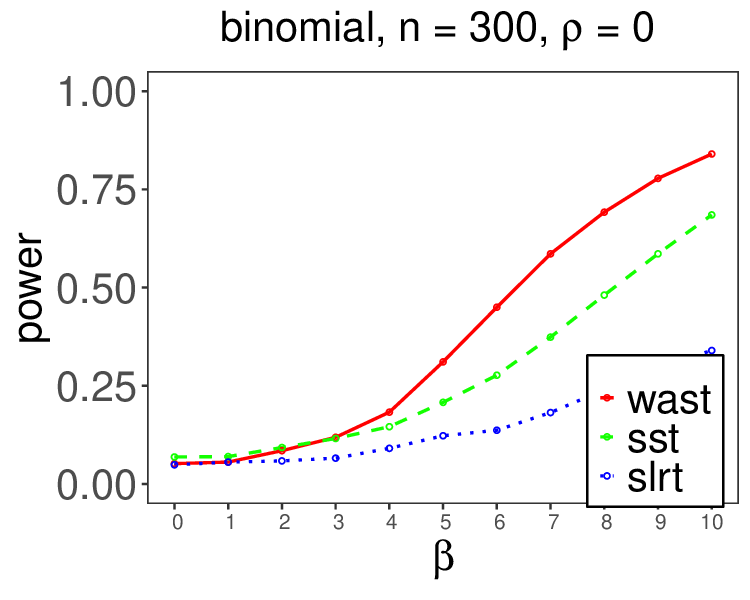}
		\includegraphics[scale=0.285]{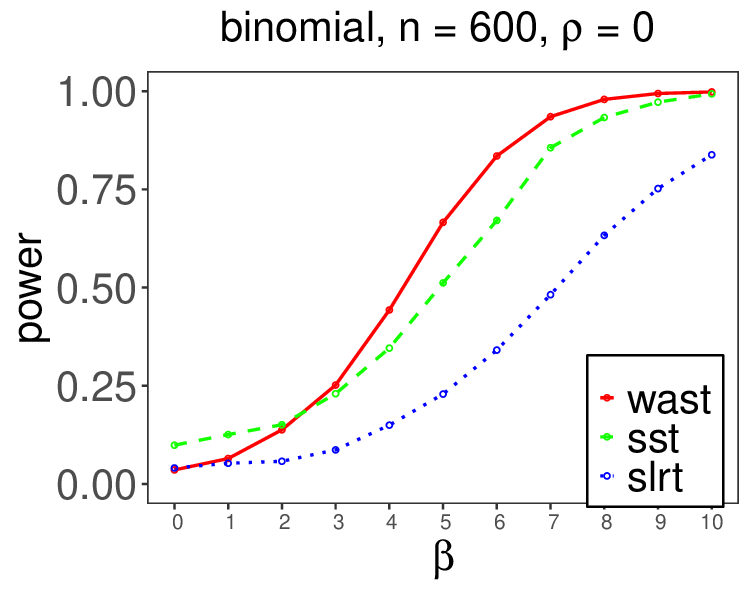}
		\includegraphics[scale=0.285]{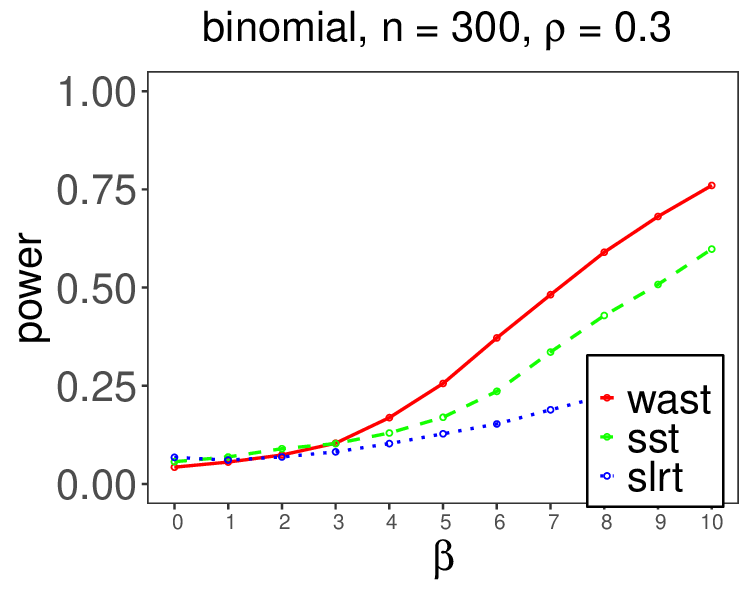}
		\includegraphics[scale=0.285]{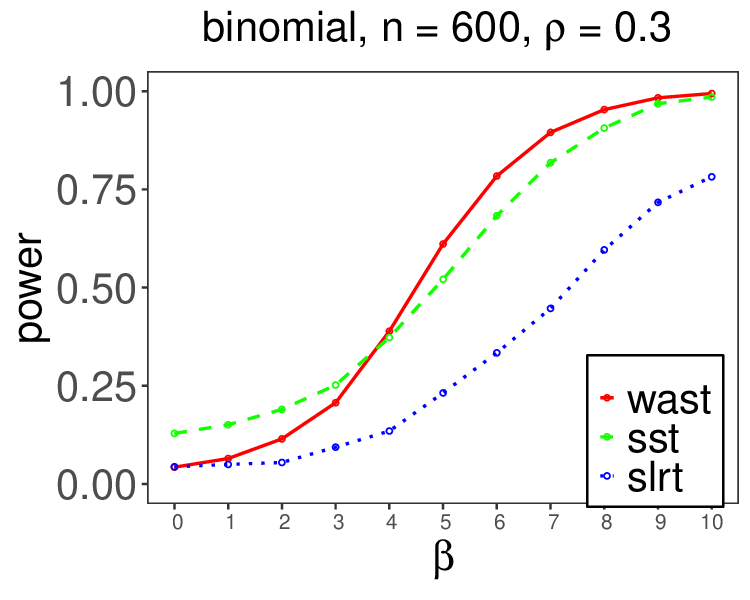}
		\includegraphics[scale=0.285]{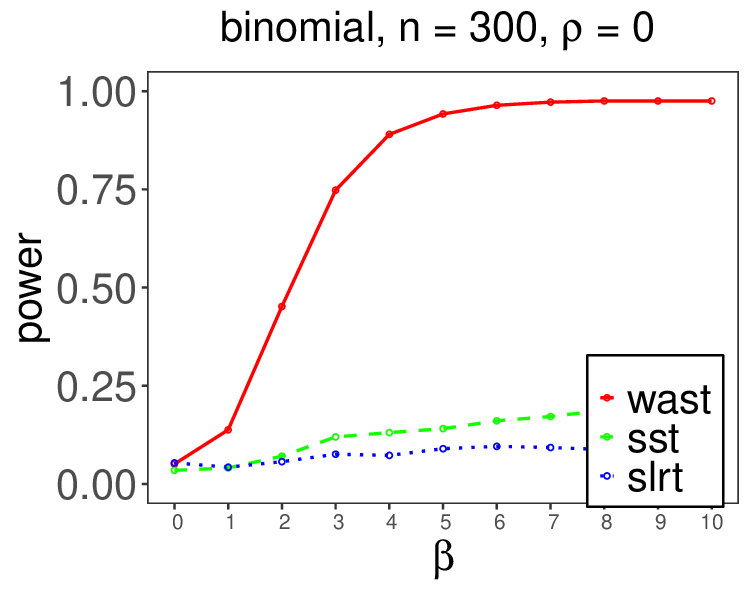}
		\includegraphics[scale=0.285]{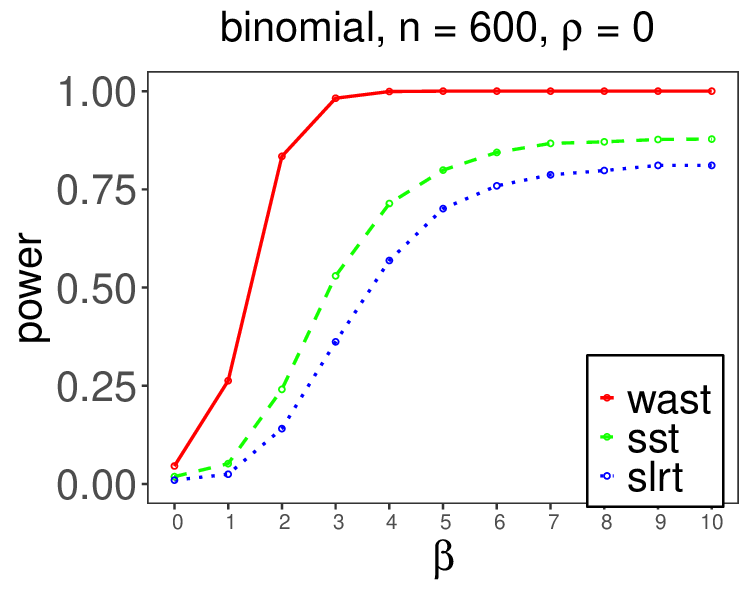}
		\includegraphics[scale=0.285]{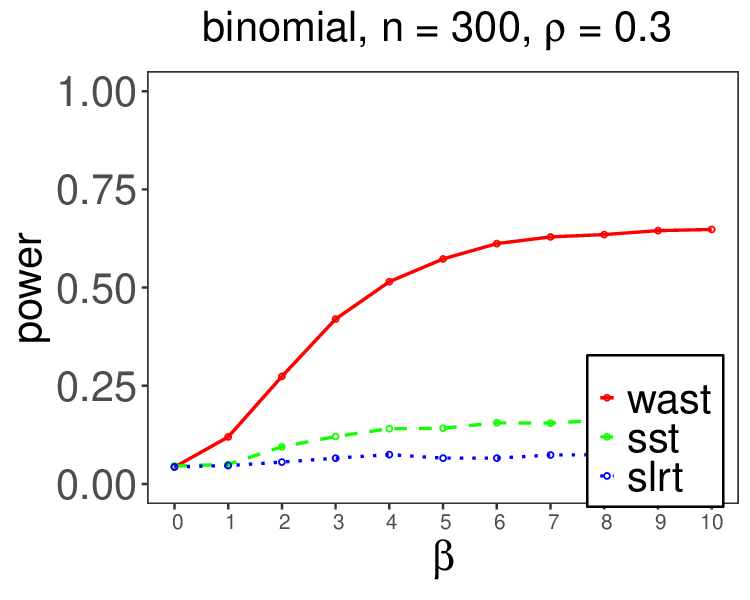}
		\includegraphics[scale=0.285]{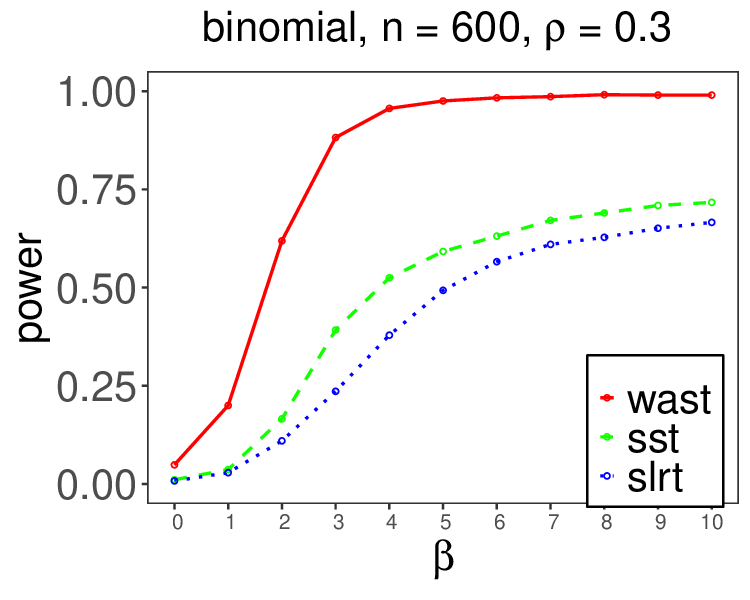}
		\includegraphics[scale=0.285]{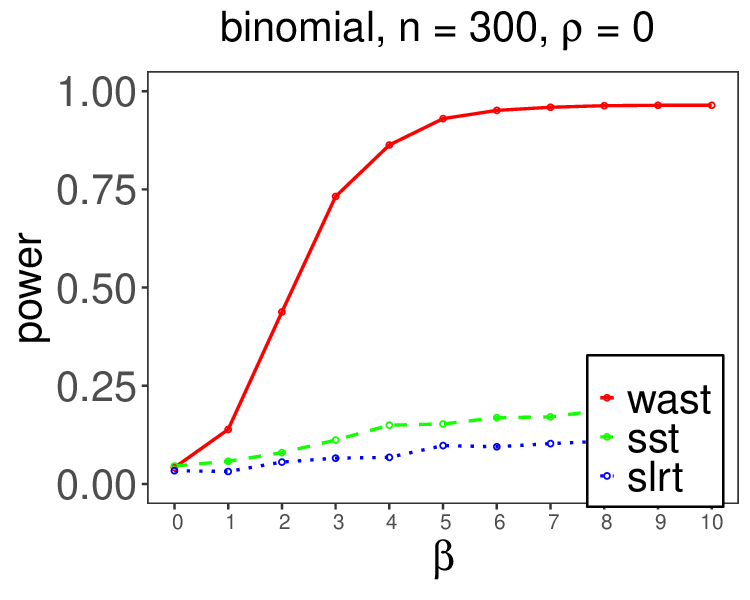}
		\includegraphics[scale=0.285]{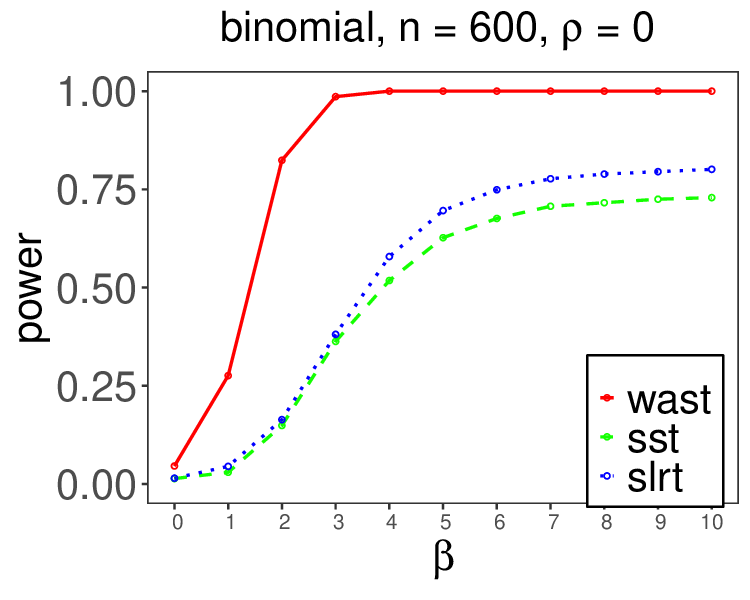}
		\includegraphics[scale=0.285]{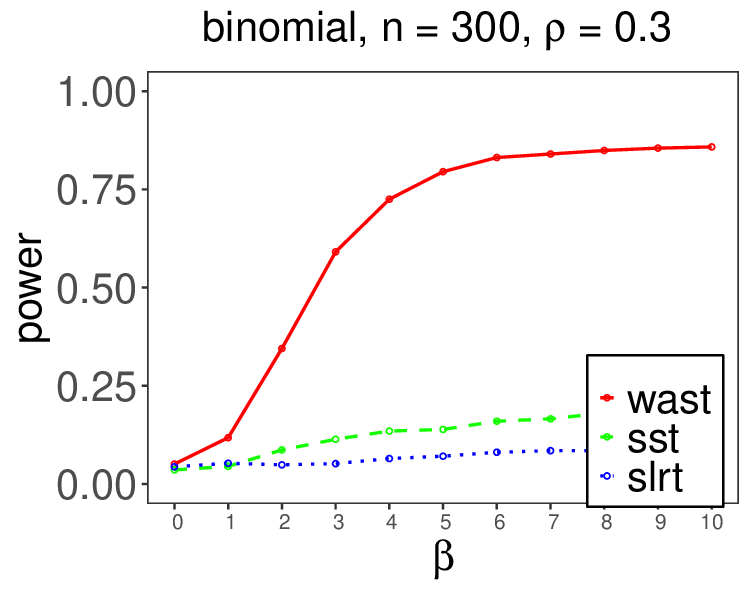}
		\includegraphics[scale=0.285]{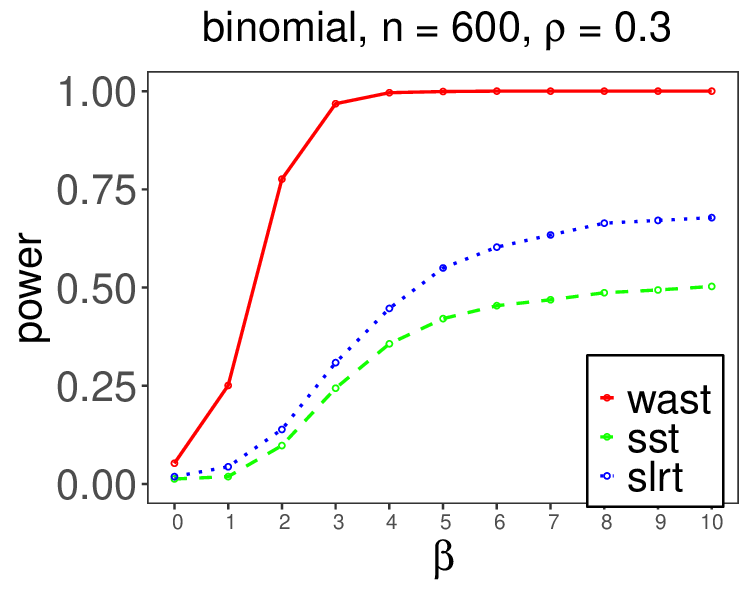}
		\caption{Powers of testing logistic regression by WAST (red solid line), SST (green dashed line), and SLRT (blue dotted line) for $n=(300,600)$. From top to bottom, each row depicts the powers for the cases $(r,p,q)=(2,2,3)$, $(6,6,3)$, $(2,2,11)$, $(6,6,11)$, $(2,51,11)$, and $(6,51,11)$.}
		\label{fig_binomial}
	\end{center}
\end{figure}

\begin{figure}[!ht]
	\begin{center}
		\includegraphics[scale=0.3]{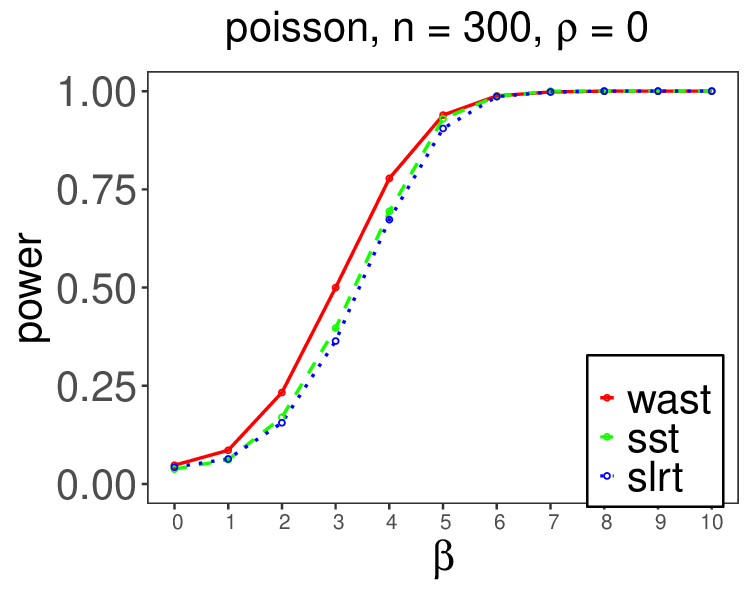}
		\includegraphics[scale=0.3]{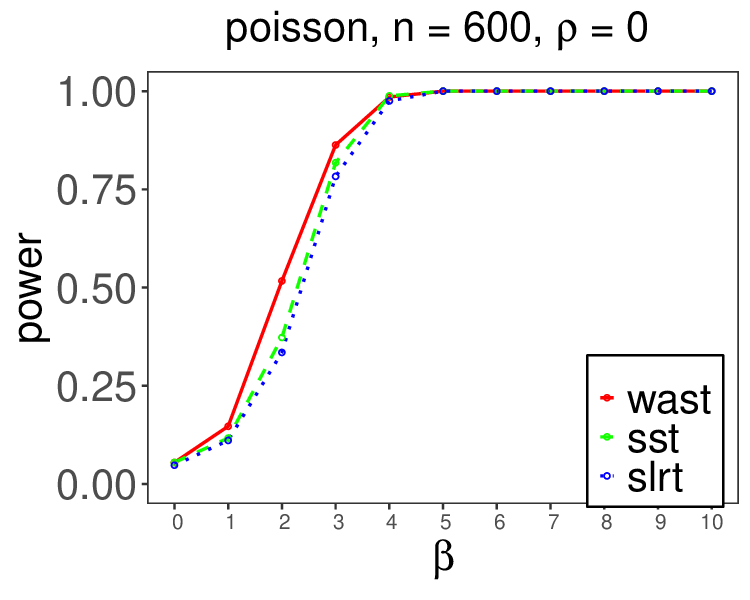}
		\includegraphics[scale=0.3]{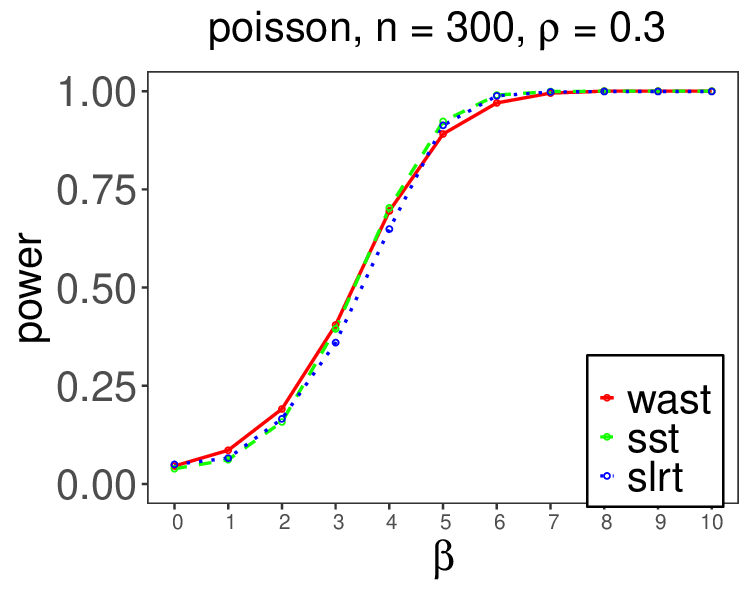}
		\includegraphics[scale=0.3]{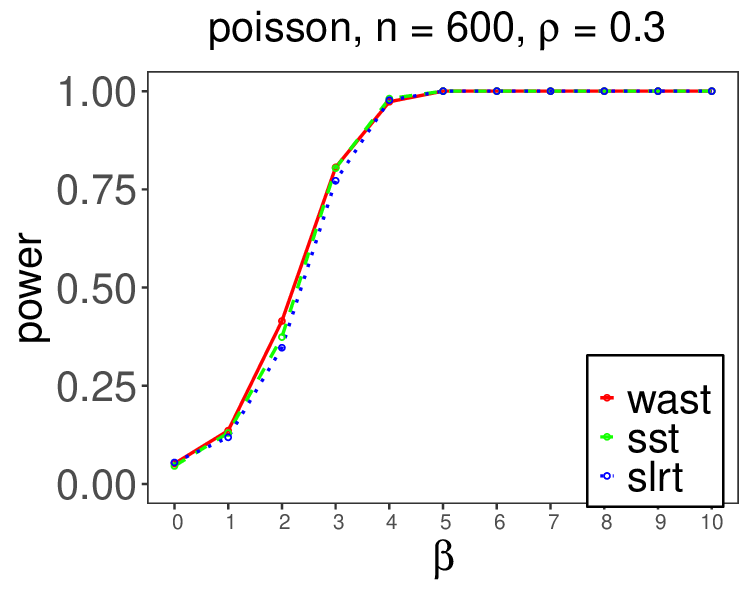}
		\includegraphics[scale=0.3]{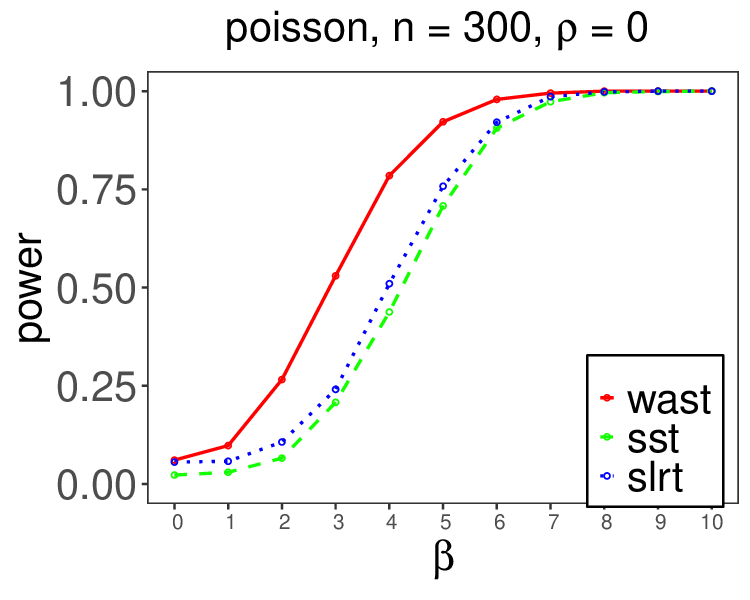}
		\includegraphics[scale=0.3]{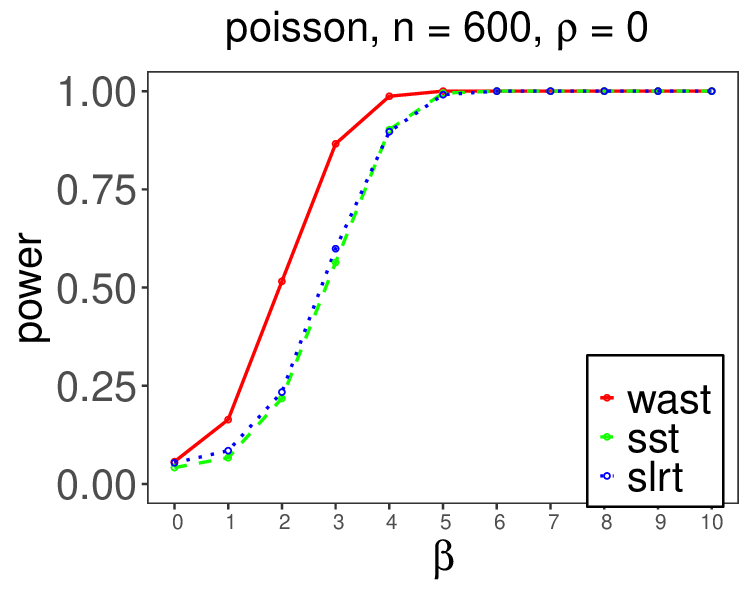}
		\includegraphics[scale=0.3]{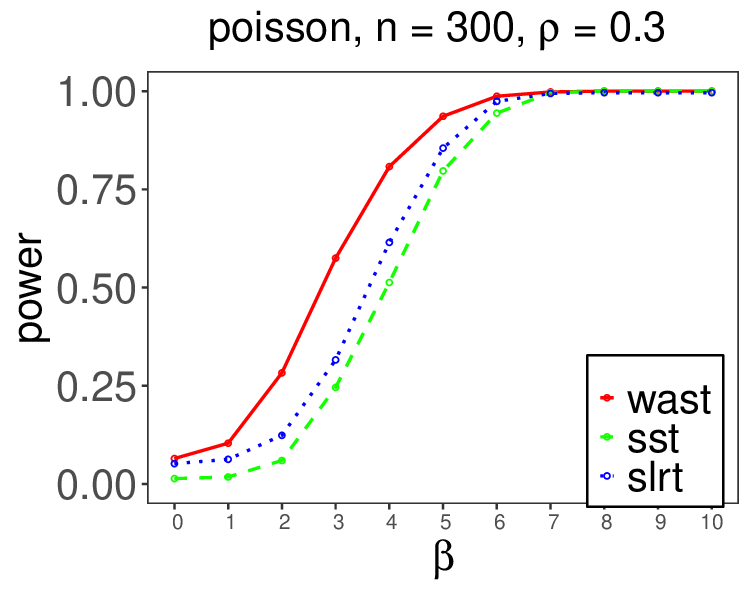}
		\includegraphics[scale=0.3]{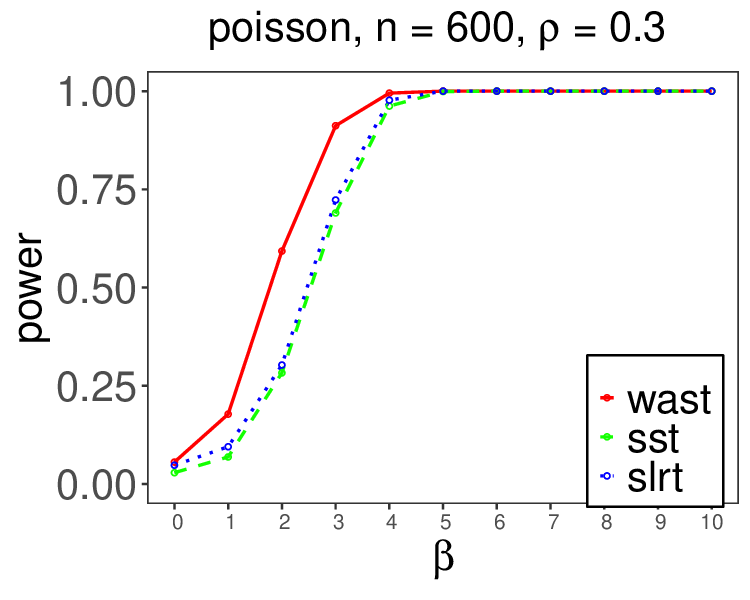}	
		\includegraphics[scale=0.3]{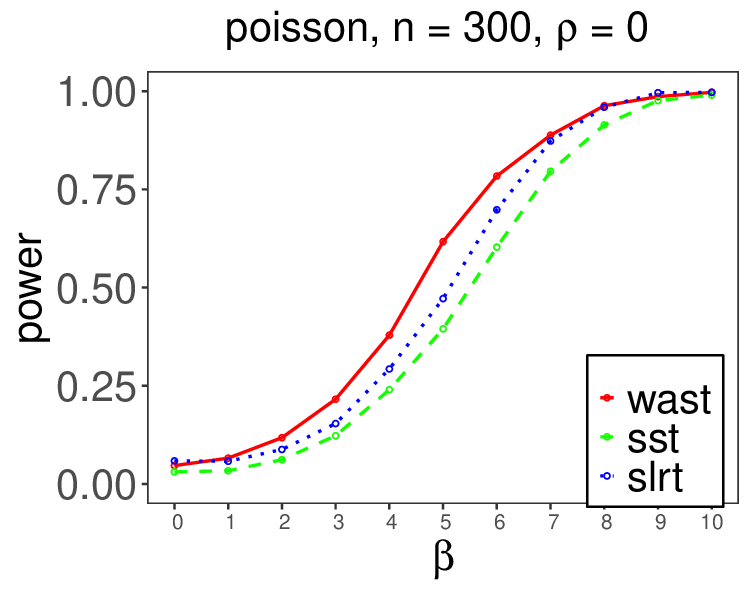}
		\includegraphics[scale=0.3]{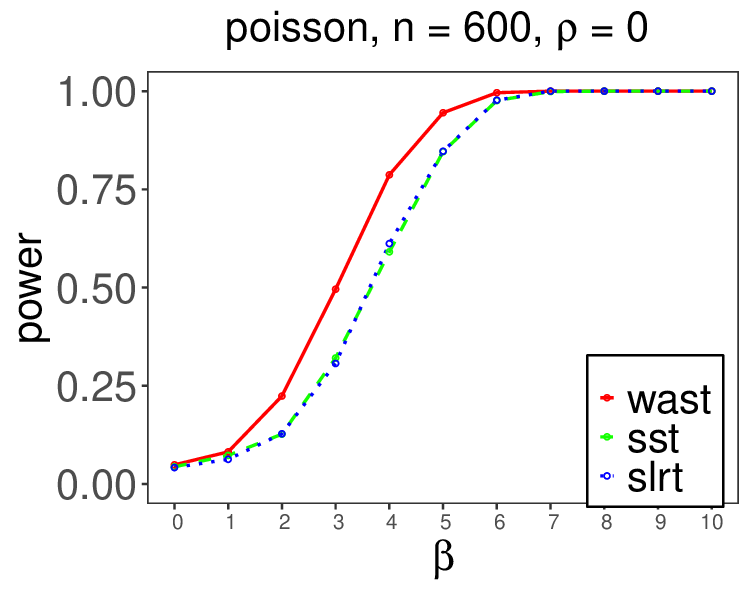}
		\includegraphics[scale=0.3]{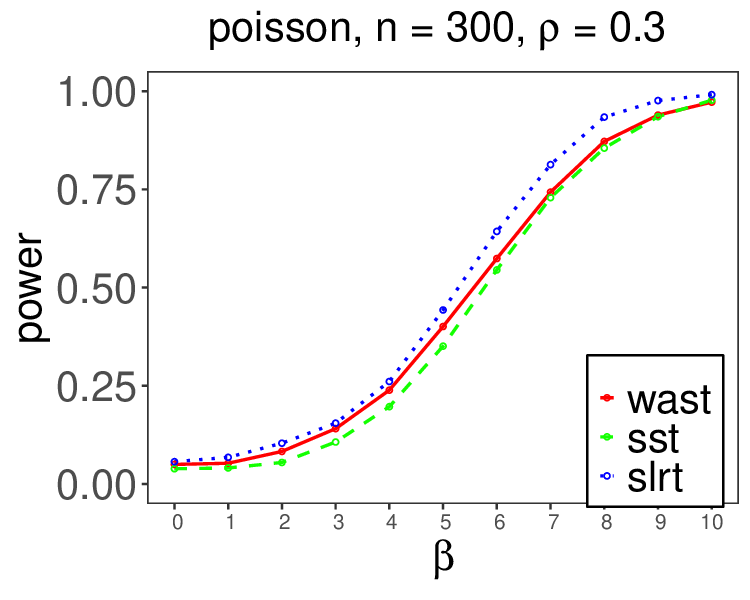}
		\includegraphics[scale=0.3]{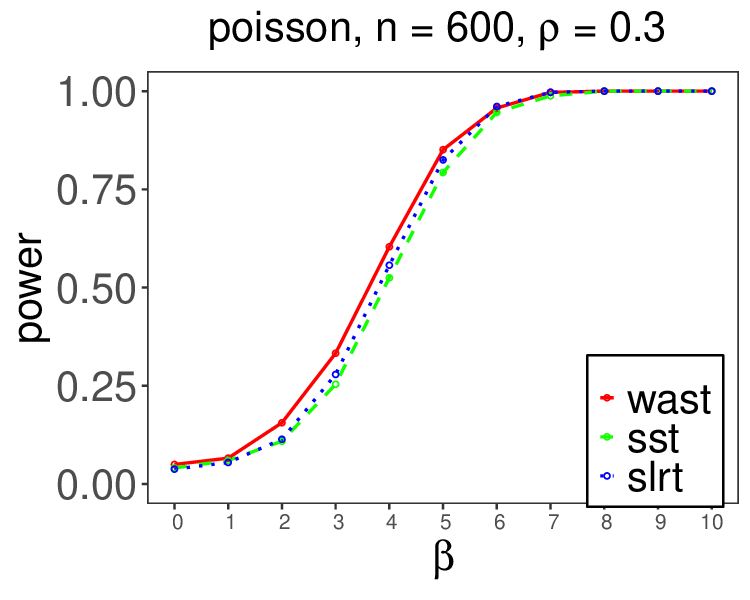}
		\includegraphics[scale=0.3]{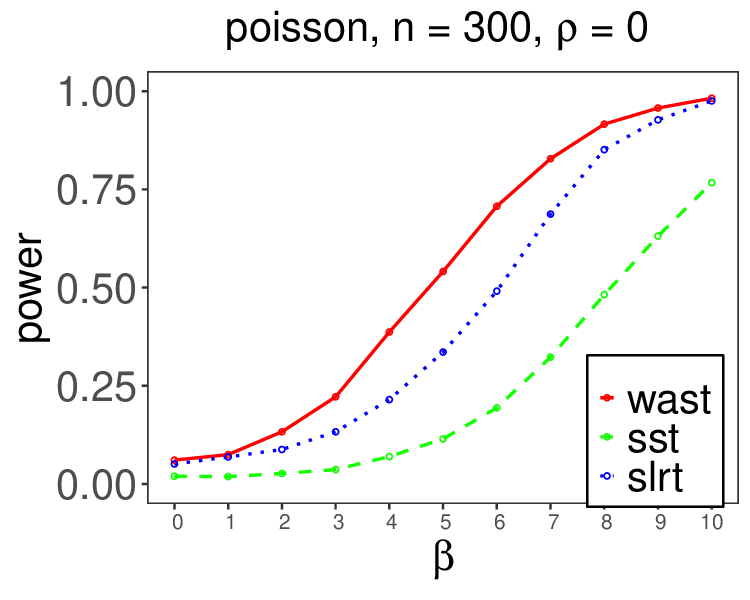}
		\includegraphics[scale=0.3]{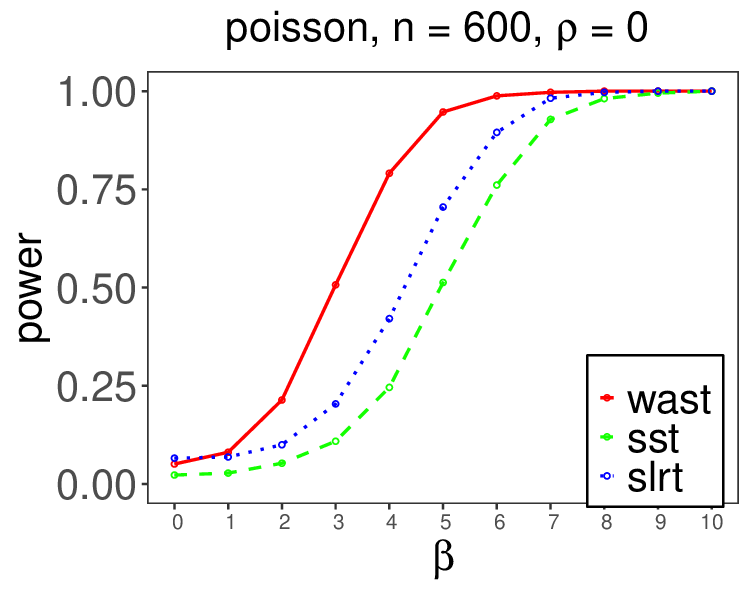}
		\includegraphics[scale=0.3]{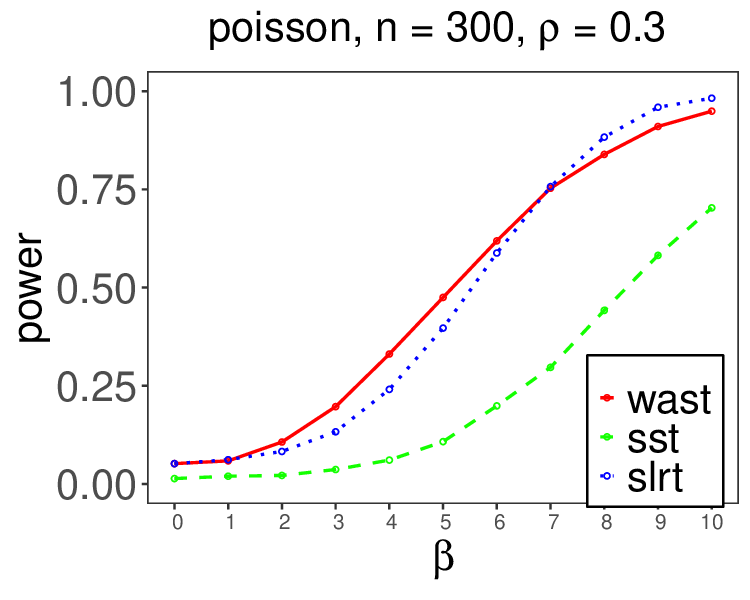}
		\includegraphics[scale=0.3]{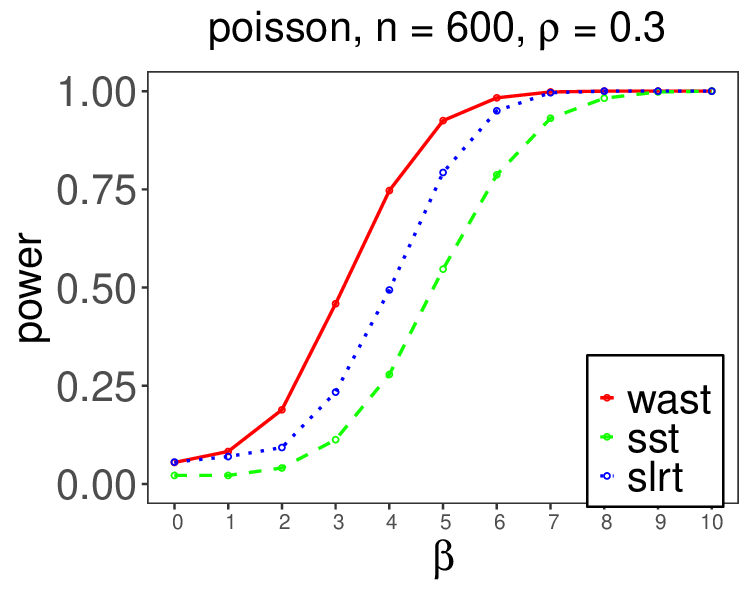}
		\includegraphics[scale=0.3]{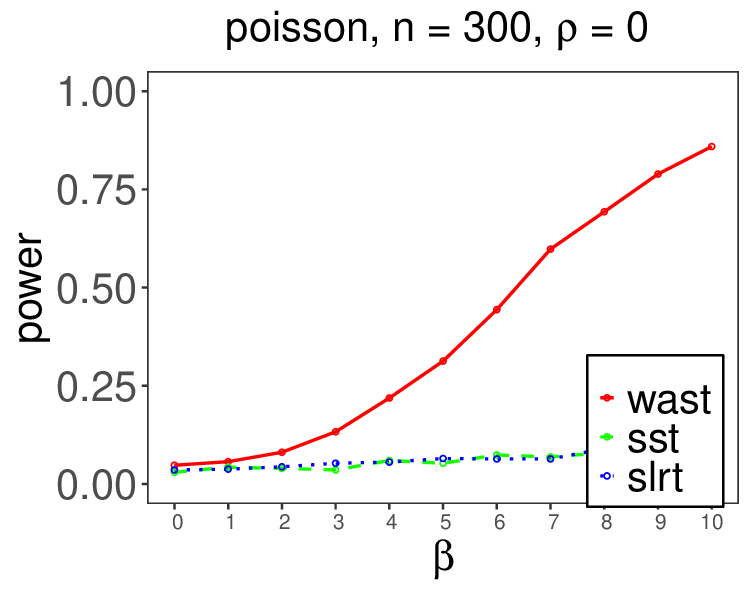}
		\includegraphics[scale=0.3]{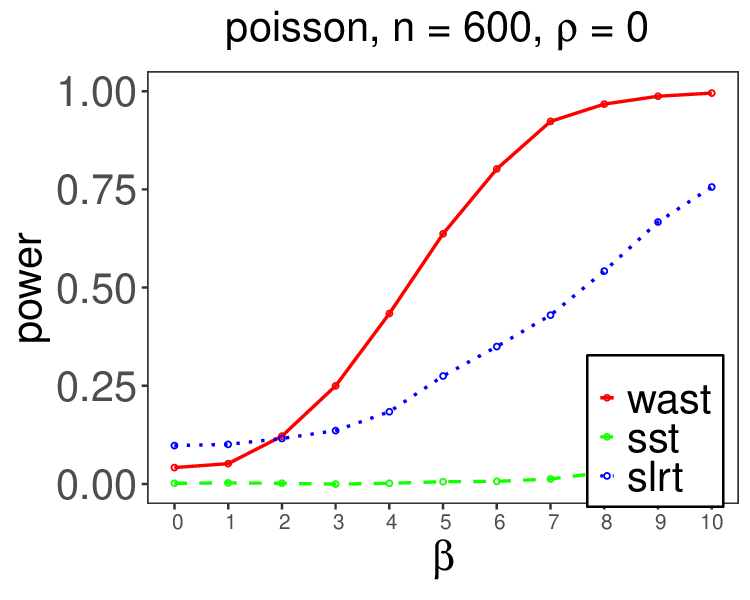}
		\includegraphics[scale=0.3]{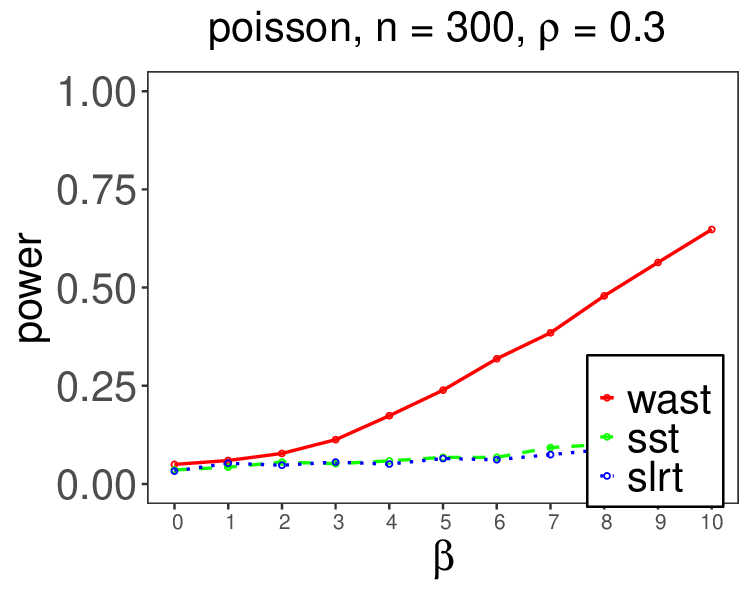}
		\includegraphics[scale=0.3]{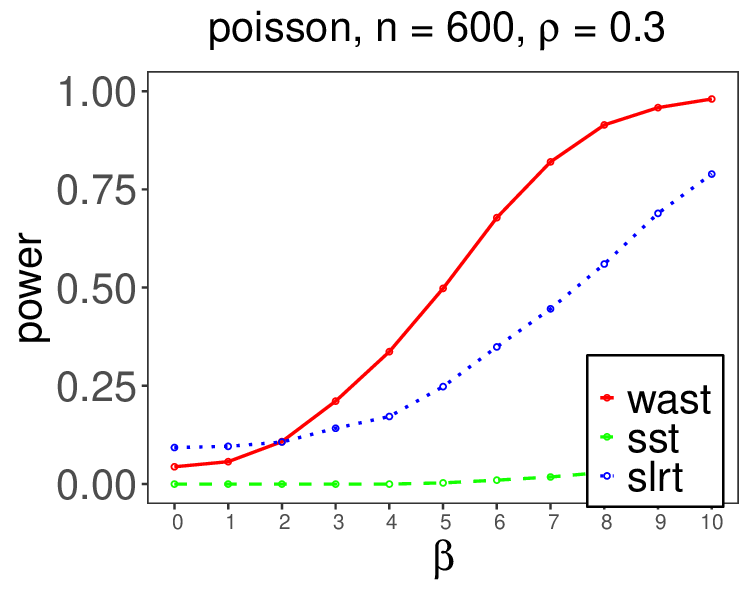}
		\includegraphics[scale=0.3]{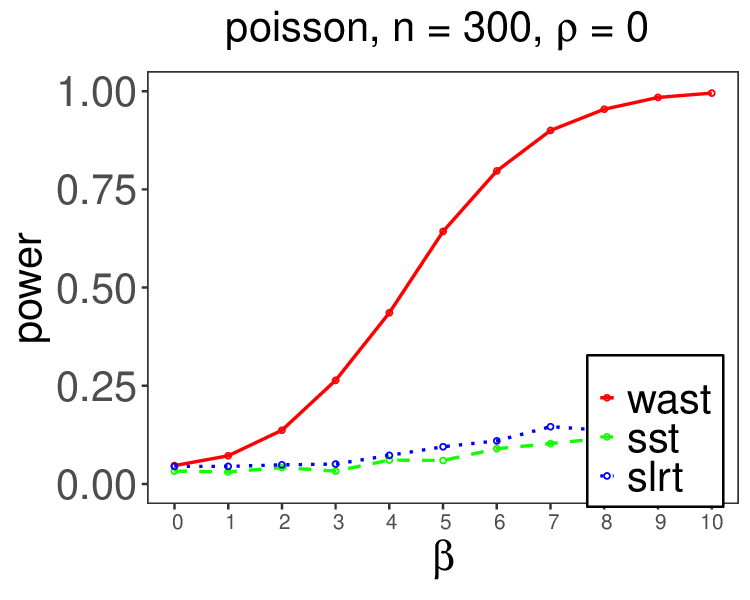}
		\includegraphics[scale=0.3]{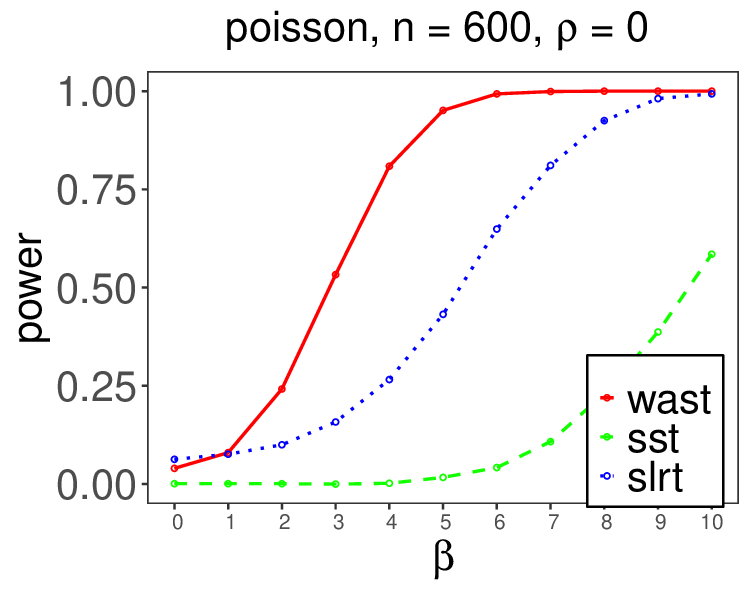}
		\includegraphics[scale=0.3]{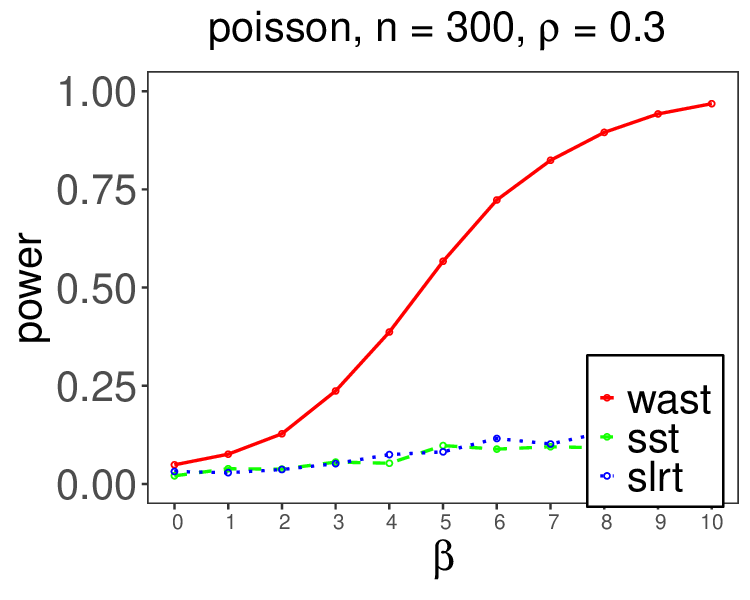}
		\includegraphics[scale=0.3]{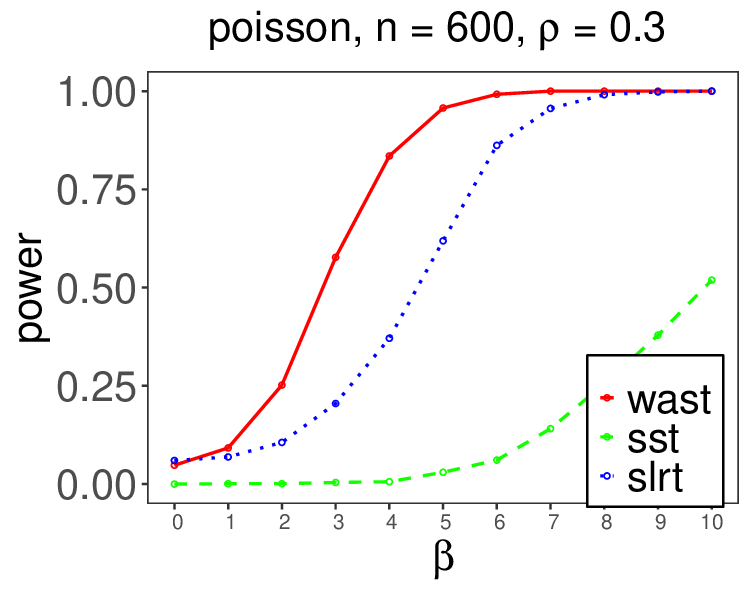}
		\caption{\it Powers of testing Poisson regression by the proposed WAST (red solid line), SST (green dashed line), and SLRT (blue dotted line) for $n=(300,600)$. From top to bottom, each row panel depicts the powers for the case $(r,p,q)=(2,2,3)$, $(6,6,3)$, $(2,2,11)$, $(6,6,11)$, $(2,51,11)$, and $(6,51,11)$.
}
		\label{fig_poisson}
	\end{center}
\end{figure}

Figure~\ref{fig_gaussian} - \ref{fig_poisson} show the powers for all the scenarios considered in Table~\ref{table_size}. As can be seen, the powers increase as the sample size $n$ increases, which verifies the asymptotic theory as expected. For the Gaussian, binomial and Poisson families, the power of the proposed WAST increases faster than those of the SST and SLRT, and this superiority is particularly evident when $n$ is small or $p$ and $q$ increase. In summary, the proposed test statistic shows very competitive performance.

\begin{table}
	\caption{\label{table_size} Type-\uppercase\expandafter{\romannumeral1} errors of the WAST, SST, and SLRT for GLMs.}
	\resizebox{\textwidth}{!}{
		\begin{threeparttable}
			\begin{tabular}{ll*{15}{c}}\\
				\hline
				\multirow{3}{*}{Family}&\multirow{3}{*}{$(r,p,q)$}
				&\multicolumn{7}{c}{ $\rho=0$} && \multicolumn{7}{c}{$\rho=0.3$}\\
				\cline{3-9} \cline{11-17}
				& &\multicolumn{3}{c}{ $n=300$} && \multicolumn{3}{c}{ $n=600$} && \multicolumn{3}{c}{ $n=300$} && \multicolumn{3}{c}{ $n=600$}\\
				\cline{3-5} \cline{7-9} \cline{11-13} \cline{15-17}
				&
				& WAST& SST & SLRT&& WAST& SST & SLRT&& WAST& SST & SLRT&& WAST& SST & SLRT\\
				\cline{3-17}
				Gaussian &$(2,2,3)$ & 0.045 & 0.033 & 0.048 && 0.052 & 0.036 & 0.049 && 0.046 & 0.032 & 0.049 && 0.048 & 0.044 & 0.046 \\
				&$(6,6,3)$ & 0.049 & 0.013 & 0.040 && 0.050 & 0.030 & 0.042 && 0.059 & 0.016 & 0.042 && 0.049 & 0.028 & 0.043 \\
				& $(2,2,11)$ & 0.044 & 0.030 & 0.058 && 0.060 & 0.041 & 0.059 && 0.049 & 0.020 & 0.051 && 0.052 & 0.031 & 0.047 \\
				&$(6,6,11)$ & 0.050 & 0.016 & 0.048 && 0.052 & 0.026 & 0.044 && 0.050 & 0.008 & 0.048 && 0.047 & 0.023 & 0.048 \\
				& $(2,51,11)$ & 0.045 & 0.032 & 0.044 && 0.048 & 0.000 & 0.059 && 0.044 & 0.034 & 0.050 && 0.048 & 0.000 & 0.045 \\
				&$(6,51,11)$ & 0.046 & 0.026 & 0.039 && 0.054 & 0.000 & 0.050 && 0.056 & 0.023 & 0.043 && 0.043 & 0.000 & 0.046 \\
				[1 ex]
				Binomial &$(2,2,3)$ & 0.053 & 0.109 & 0.052 && 0.050 & 0.082 & 0.054 && 0.055 & 0.100 & 0.048 && 0.051 & 0.072 & 0.059 \\
				&$(6,6,3)$ & 0.056 & 0.081 & 0.055 && 0.045 & 0.110 & 0.067 && 0.059 & 0.076 & 0.048 && 0.046 & 0.117 & 0.059 \\
				& $(2,2,11)$ & 0.046 & 0.124 & 0.053 && 0.057 & 0.089 & 0.065 && 0.052 & 0.106 & 0.048 && 0.050 & 0.099 & 0.057 \\
				&$(6,6,11)$ & 0.052 & 0.069 & 0.049 && 0.036 & 0.099 & 0.041 && 0.043 & 0.057 & 0.068 && 0.043 & 0.129 & 0.044 \\
				& $(2,51,11)$ & 0.051 & 0.035 & 0.054 && 0.046 & 0.019 & 0.010 && 0.043 & 0.046 & 0.044 && 0.049 & 0.011 & 0.008 \\
				&$(6,51,11)$ & 0.042 & 0.046 & 0.034 && 0.046 & 0.014 & 0.015 && 0.051 & 0.036 & 0.044 && 0.053 & 0.013 & 0.019 \\
				[1 ex]
				Poisson &$(2,2,3)$ & 0.048 & 0.038 & 0.043 && 0.055 & 0.055 & 0.048 && 0.046 & 0.039 & 0.050 && 0.052 & 0.046 & 0.055 \\
				&$(6,6,3)$ & 0.061 & 0.023 & 0.056 && 0.057 & 0.042 & 0.054 && 0.065 & 0.014 & 0.052 && 0.056 & 0.029 & 0.048 \\
				& $(2,2,11)$ & 0.047 & 0.031 & 0.059 && 0.049 & 0.044 & 0.042 && 0.050 & 0.039 & 0.057 && 0.050 & 0.041 & 0.038 \\
				&$(6,6,11)$ & 0.061 & 0.020 & 0.051 && 0.051 & 0.023 & 0.066 && 0.052 & 0.014 & 0.052 && 0.055 & 0.022 & 0.056 \\
				& $(2,51,11)$ & 0.048 & 0.030 & 0.036 && 0.042 & 0.002 & 0.098 && 0.050 & 0.036 & 0.033 && 0.044 & 0.000 & 0.093 \\
				&$(6,51,11)$ & 0.047 & 0.033 & 0.045 && 0.040 & 0.001 & 0.063 && 0.049 & 0.021 & 0.032 && 0.048 & 0.000 & 0.060 \\
				\hline
			\end{tabular}
			\begin{tablenotes}
				\item The $\Gv$ $\bZ$ generated from a normal distribution with zero mean and unit variance. The nominal significance level is 0.05.
			\end{tablenotes}
		\end{threeparttable}
	}
\end{table}

\subsection{Single change plane for GLMs with \texorpdfstring{${\bf Z}$}{} from \texorpdfstring{$t_3$}{}}\label{scp1}
Consider the GLM in model (1.1) with the same setting as Section 5 but the predictors. Let $\bvarsigma=(\varsigma_1, \cdots,\varsigma_d)\trans$, where $d=\max(r,p)+q-2$. We generate $\varsigma_1, \cdots,\varsigma_d$ independently from $t_3$ distribution with degrees of freedom 3. We set $(v_2, \cdots,v_{\max\{r,p\}}, Z_2,\cdots,Z_r)\trans=\Sigma^{1/2}\bvarsigma$, where $\Sigma=(\tilde{\rho}_{ij})$ with $\tilde{\rho}_{ij}=1$ if $i=j$ and $\rho$ otherwise. We set $\tX_j=\bone(v_j>0)$ with $j=2,\cdots,r$,  $X_k = \bone(v_k>0)$ with $k=2,\cdots,p$, and $\tX_1=1$, $X_1=1$ and $Z_1=1$. In this case, the range of $\bZ\trans\btheta$ is larger than that when $\bZ$ from multivariate standard normal distribution as used in Section 5.

The type \uppercase\expandafter{\romannumeral1} errors of the proposed WAST and the SST and SLRT are listed in Table \ref{table_size_t3}. From Table \ref{table_size_t3}, it can be found that the type \uppercase\expandafter{\romannumeral1} errors for the proposed WAST and the SLRT are close to the nominal significant level 0.05 when $p$ is small, such as $p=2,6$. When $p$ is large ($p=51$), the WAST still achieves the type \uppercase\expandafter{\romannumeral1} error close to 0.05, while the type \uppercase\expandafter{\romannumeral1} errors of the SLRT are much larger than 0.05. As Table 1 in the main paper, the type \uppercase\expandafter{\romannumeral1} errors are far from 0.05 in most cases.

\begin{figure}[!ht]
	\begin{center}
		\includegraphics[scale=0.3]{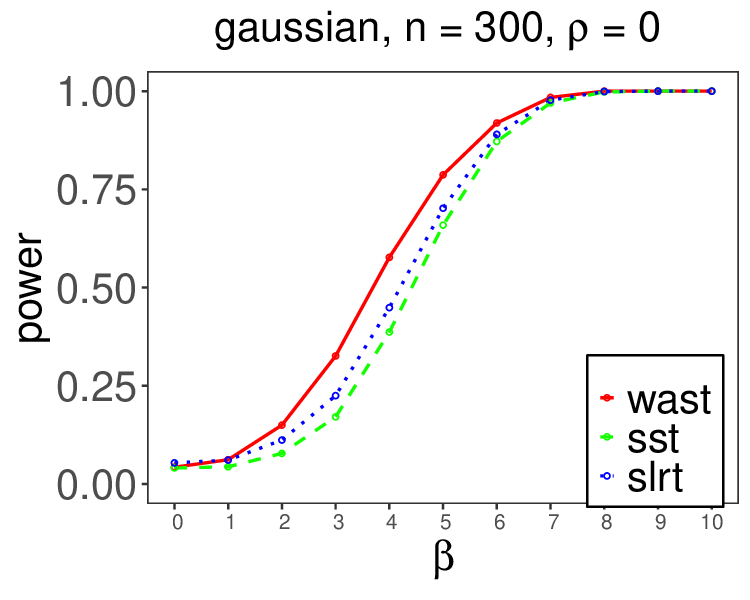}
		\includegraphics[scale=0.3]{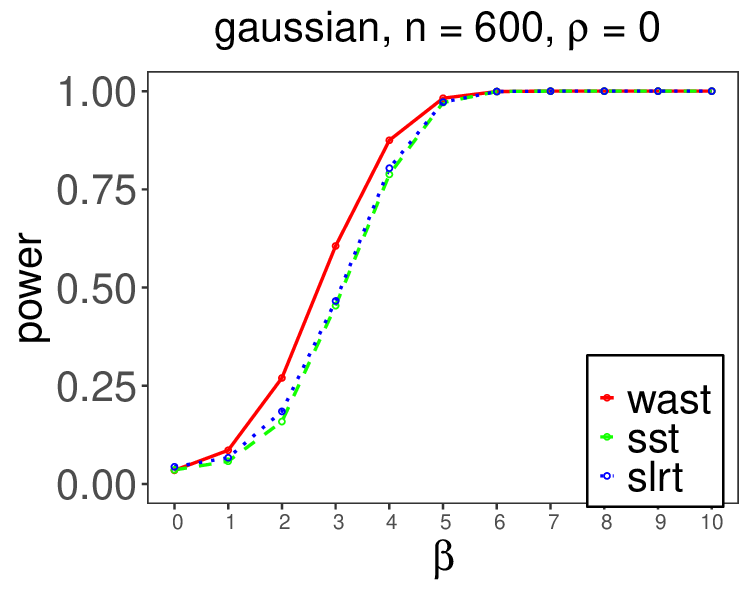}
		\includegraphics[scale=0.3]{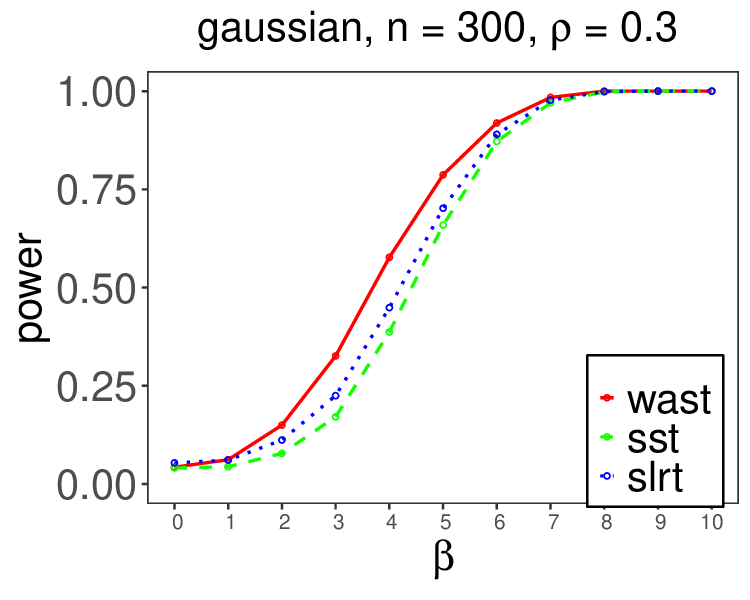}
		\includegraphics[scale=0.3]{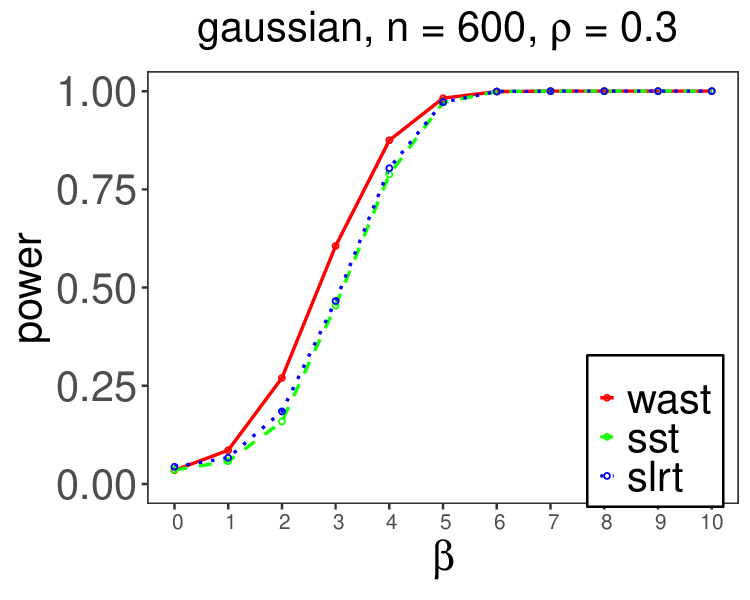}
		\includegraphics[scale=0.3]{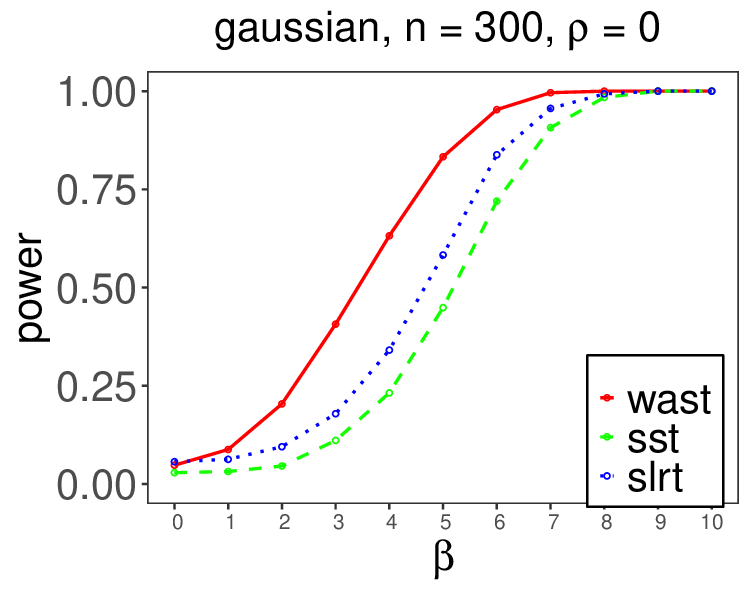}
		\includegraphics[scale=0.3]{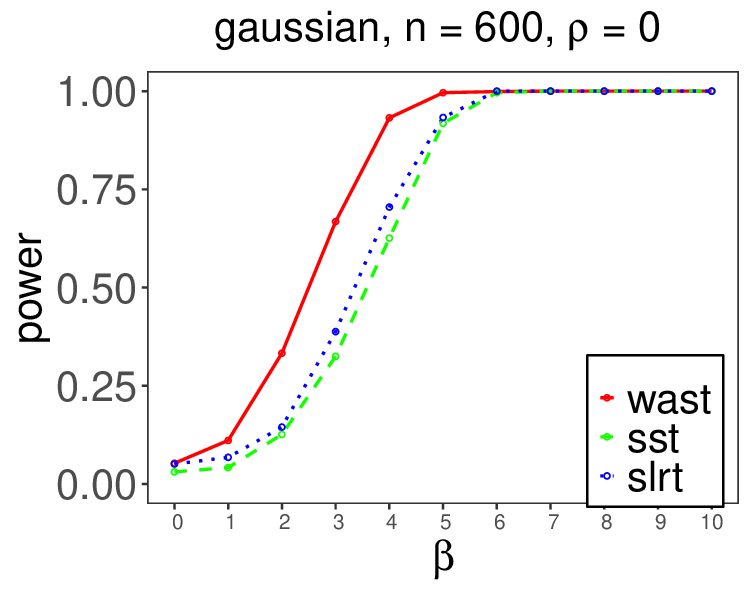}
		\includegraphics[scale=0.3]{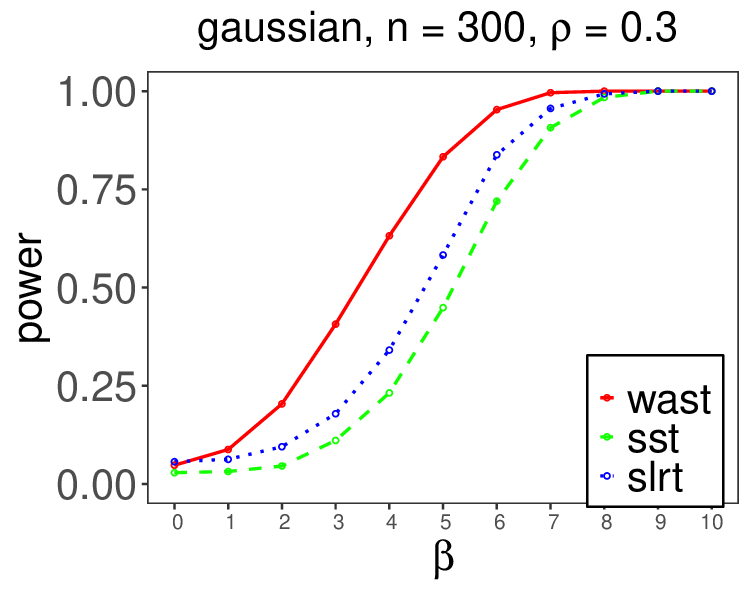}
		\includegraphics[scale=0.3]{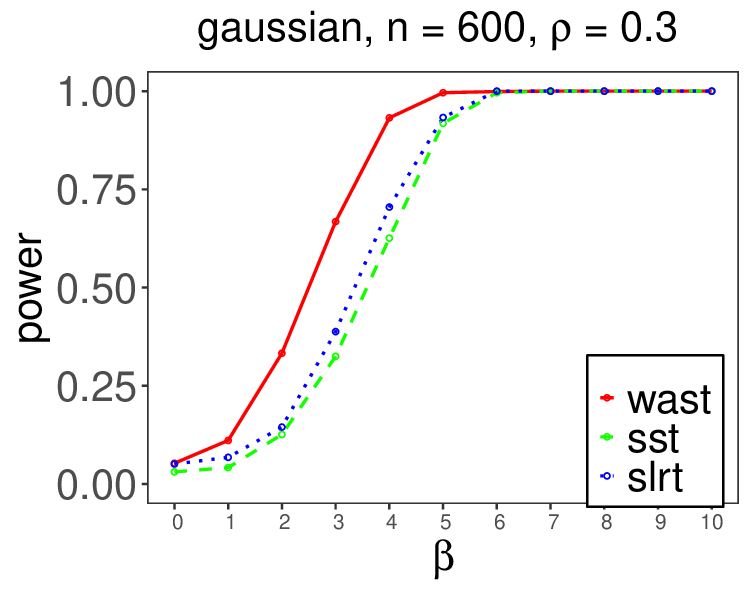}	
		\includegraphics[scale=0.3]{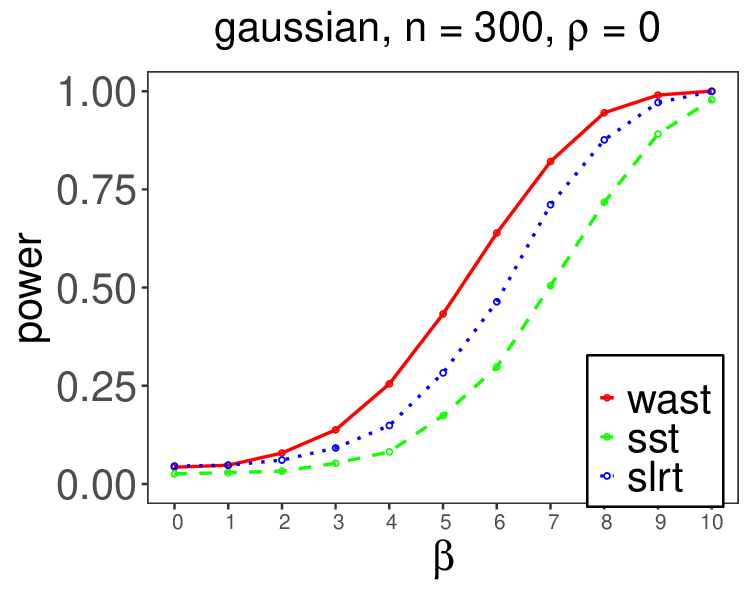}
		\includegraphics[scale=0.3]{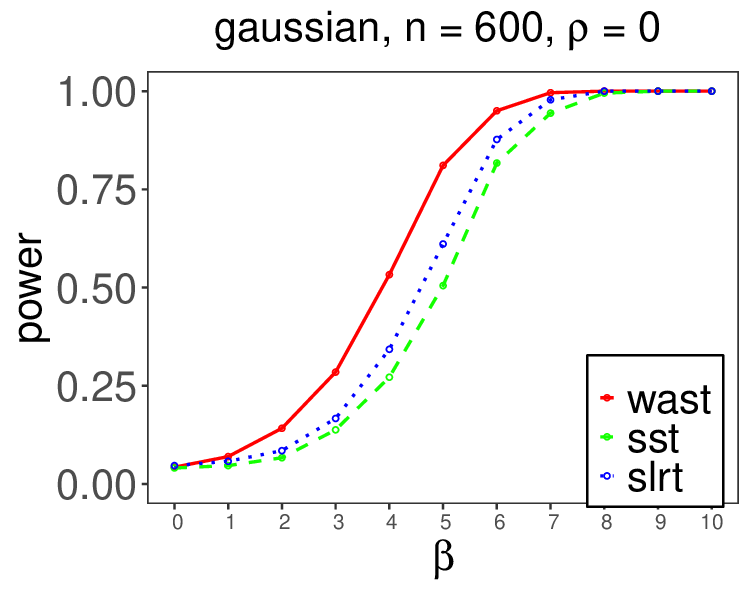}
		\includegraphics[scale=0.3]{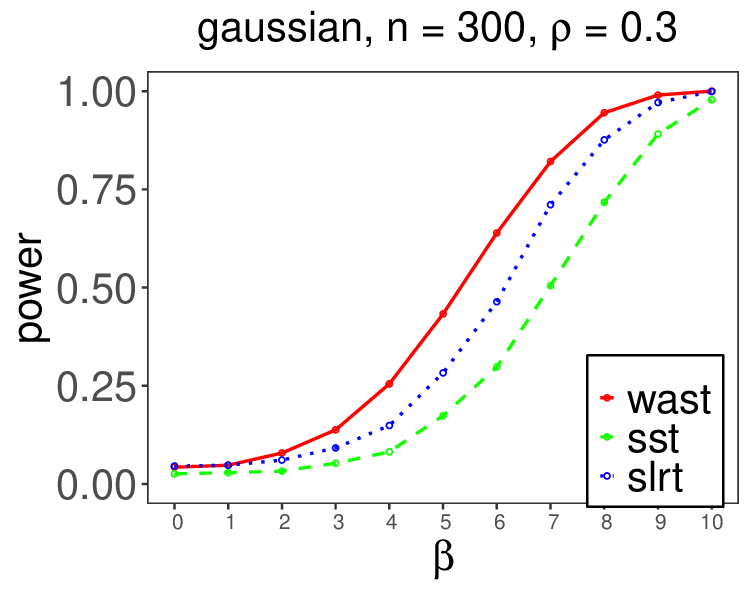}
		\includegraphics[scale=0.3]{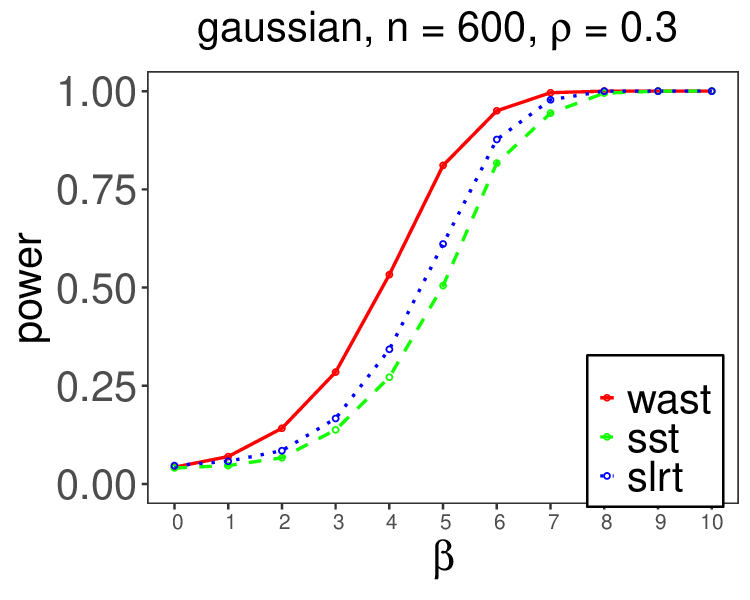}
		\includegraphics[scale=0.3]{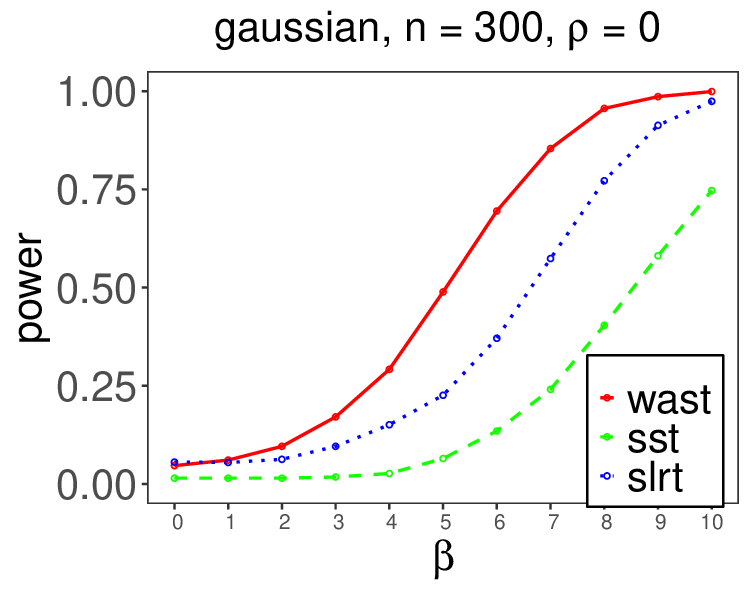}
		\includegraphics[scale=0.3]{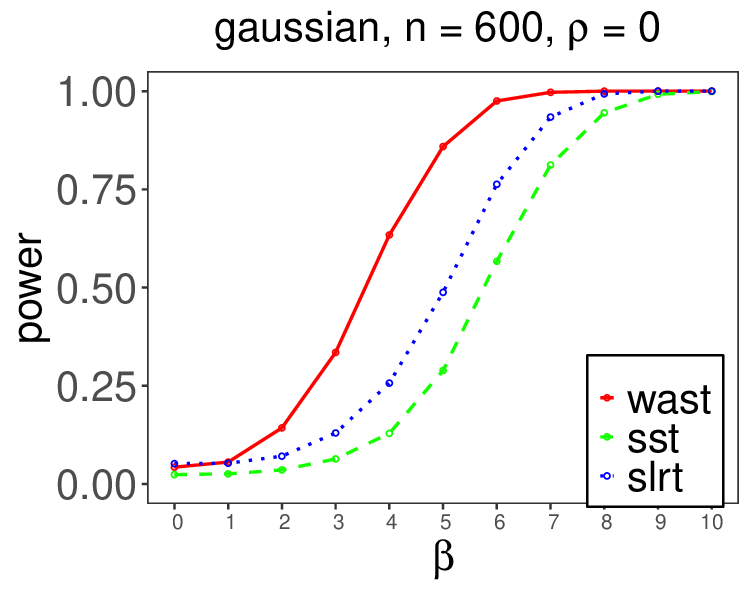}
		\includegraphics[scale=0.3]{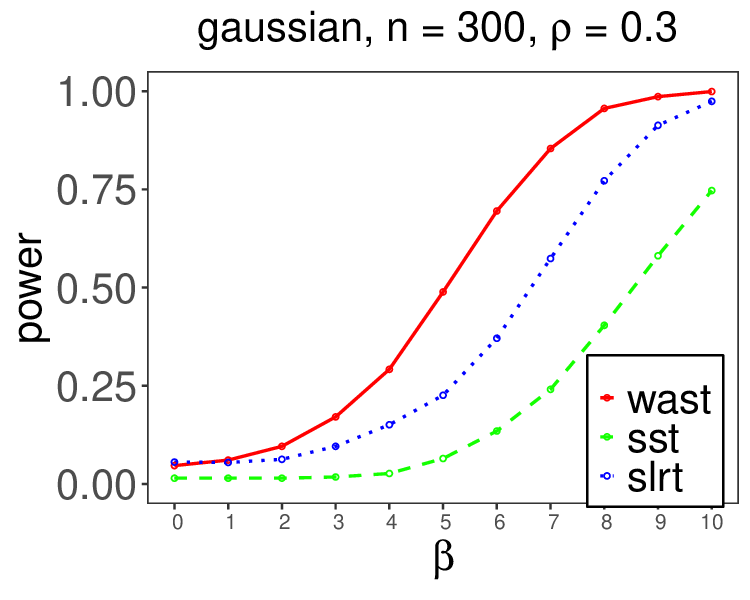}
		\includegraphics[scale=0.3]{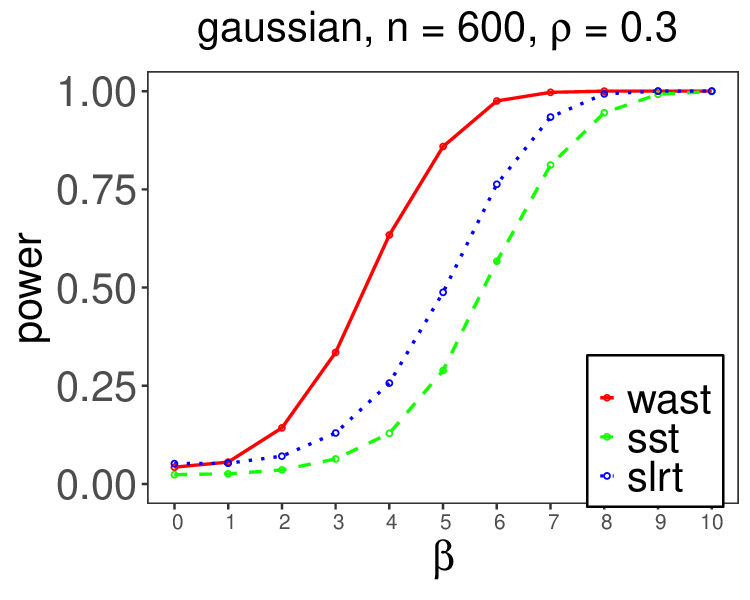}
		\includegraphics[scale=0.3]{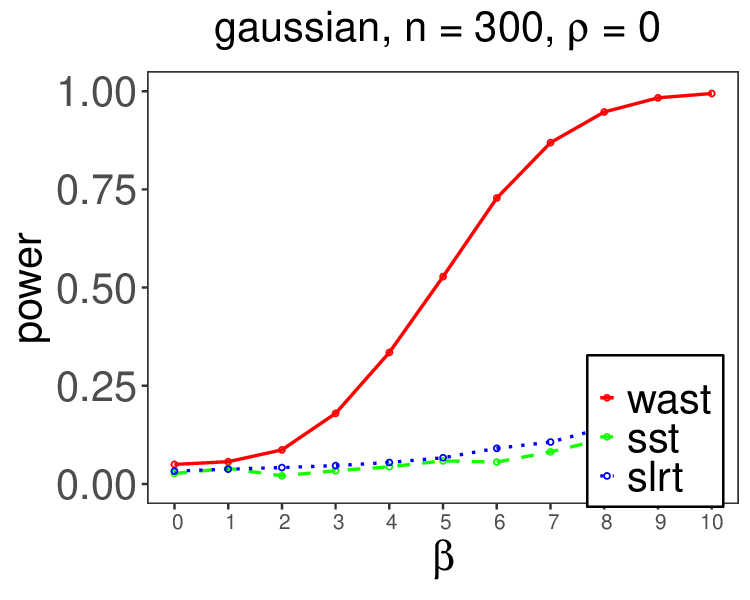}
		\includegraphics[scale=0.3]{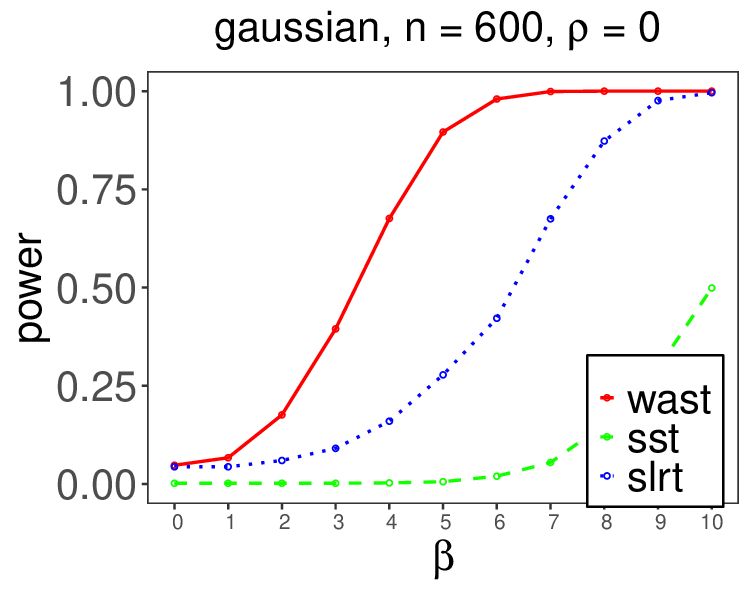}
		\includegraphics[scale=0.3]{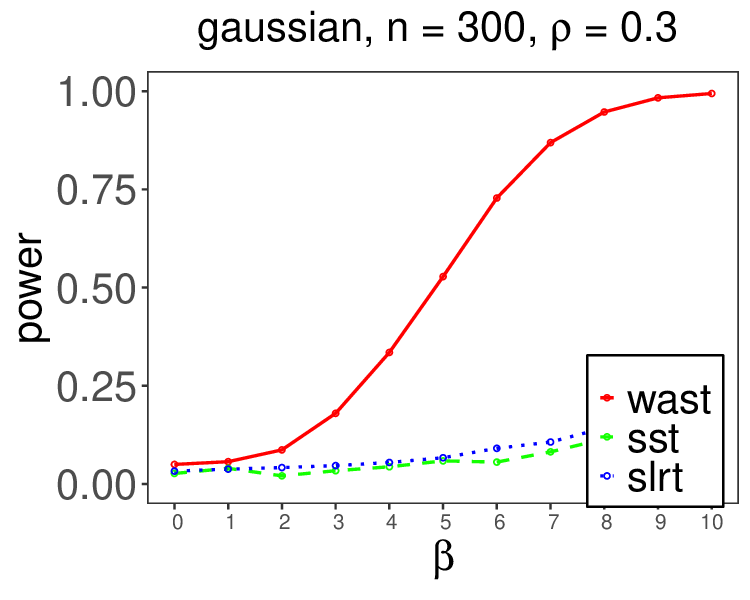}
		\includegraphics[scale=0.3]{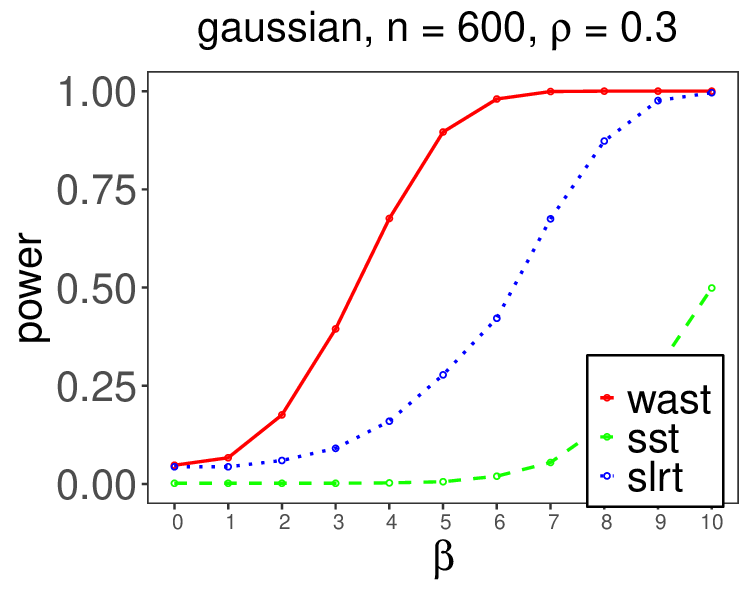}
		\includegraphics[scale=0.3]{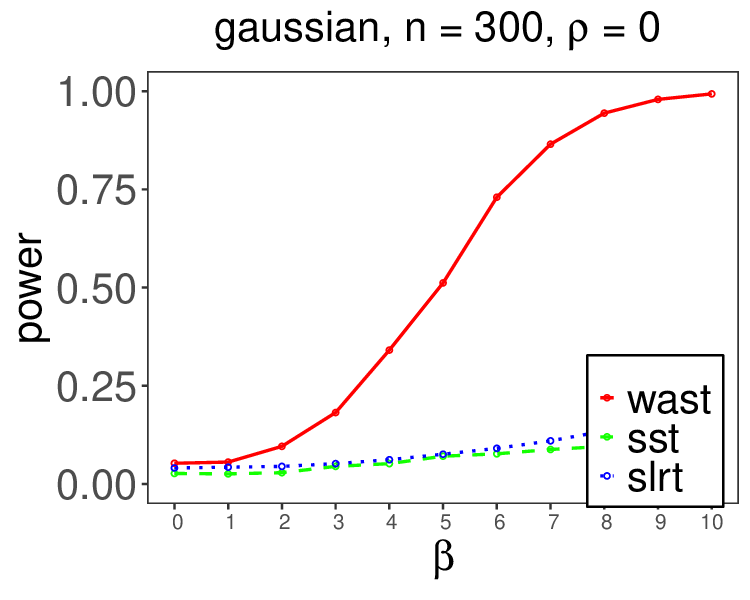}
		\includegraphics[scale=0.3]{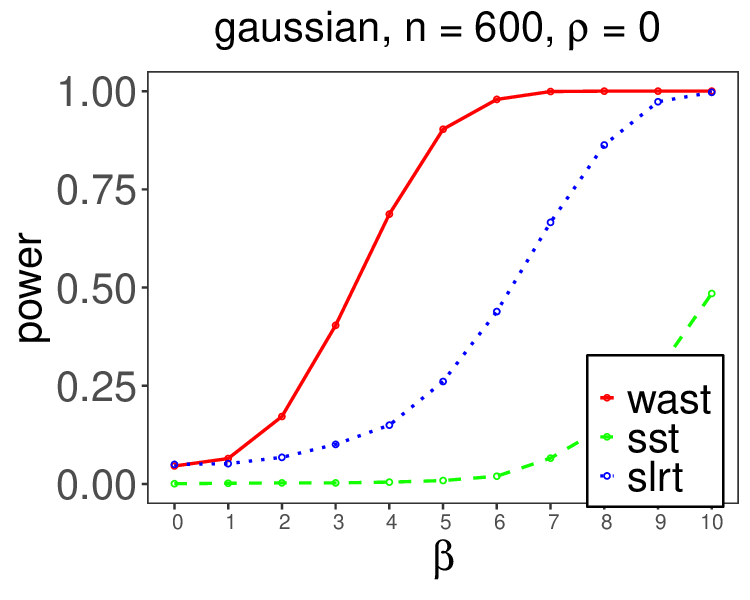}
		\includegraphics[scale=0.3]{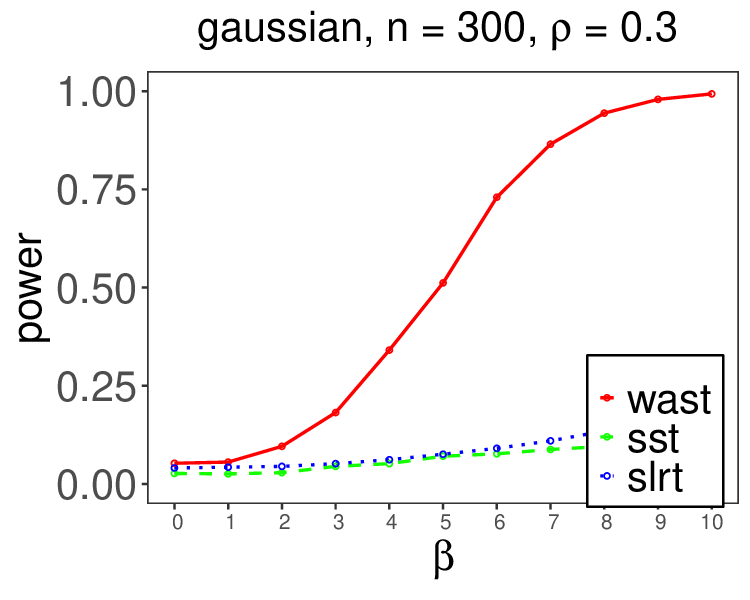}
		\includegraphics[scale=0.3]{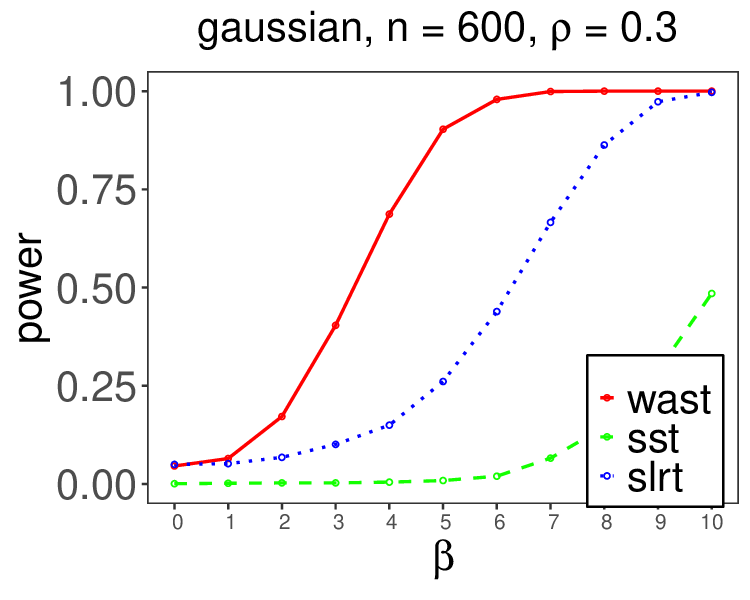}

		\caption{\it Powers of testing linear model with Gaussian error by the proposed WAST (red solid line), SST (green dashed line), and SLRT (blue dotted line) for $n=(300,600)$. From top to bottom, each row panel depicts the powers for the case $(r,p,q)=(2,2,3)$, $(6,6,3)$, $(2,2,11)$, $(6,6,11)$, $(2,51,11)$, and $(6,51,11)$.
}
		\label{fig_gaussian_zt3}
	\end{center}
\end{figure}

\begin{figure}[!ht]
	\begin{center}
		\includegraphics[scale=0.3]{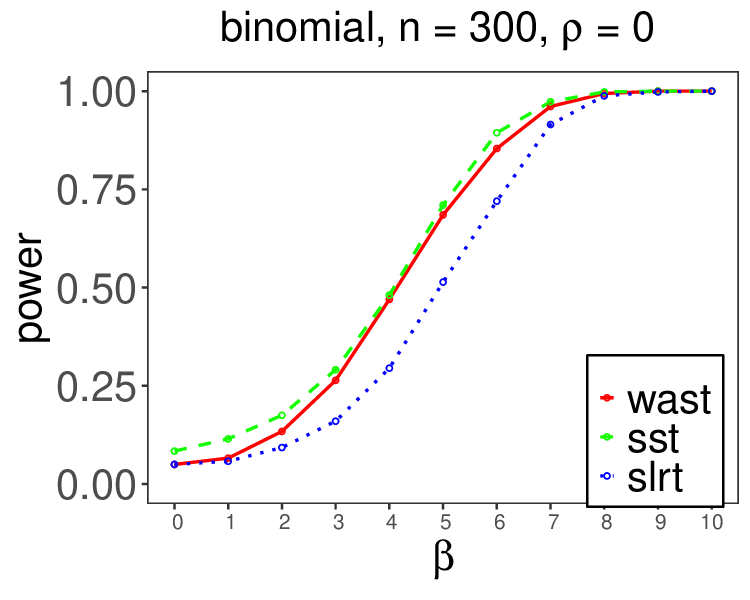}
		\includegraphics[scale=0.3]{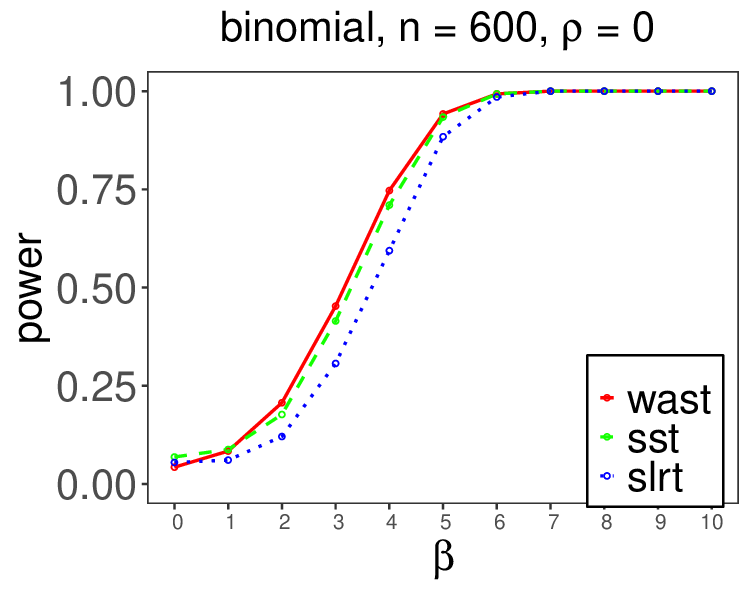}
		\includegraphics[scale=0.3]{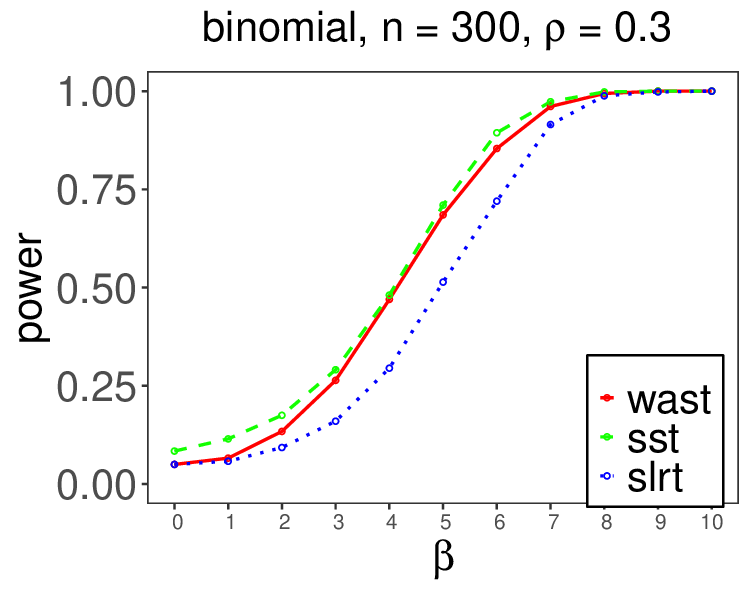}
		\includegraphics[scale=0.3]{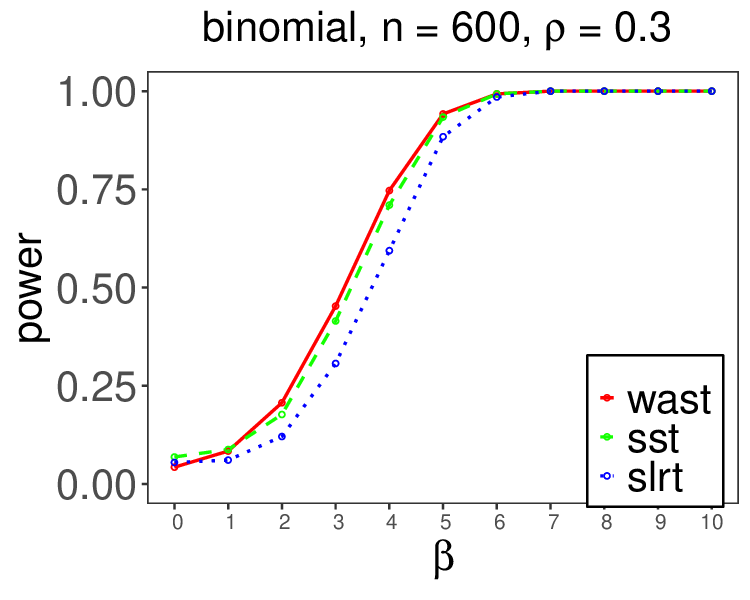}
		\includegraphics[scale=0.3]{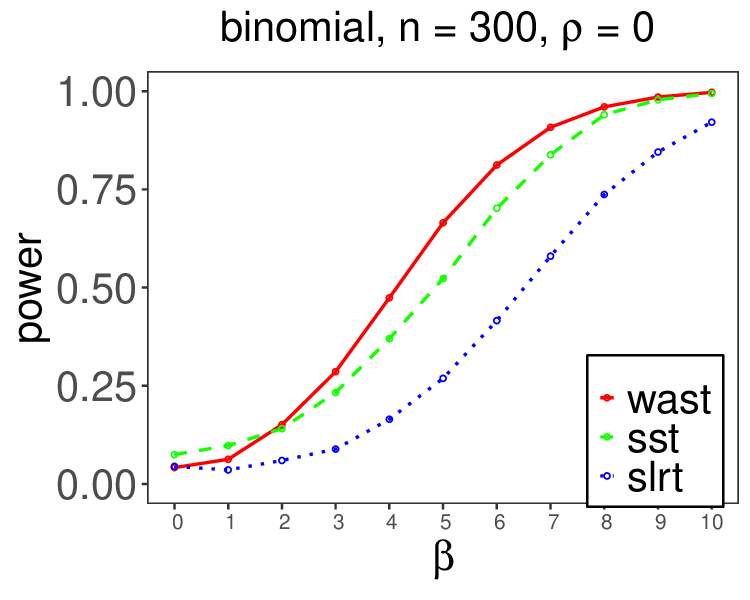}
		\includegraphics[scale=0.3]{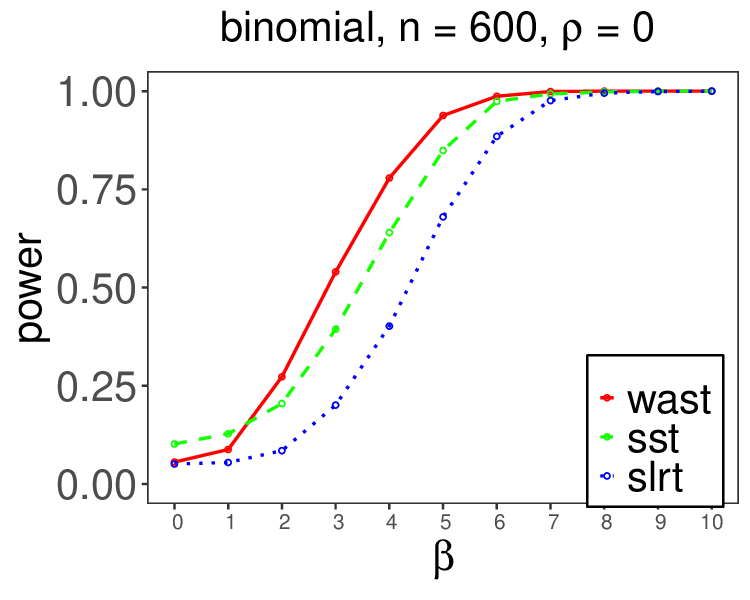}
		\includegraphics[scale=0.3]{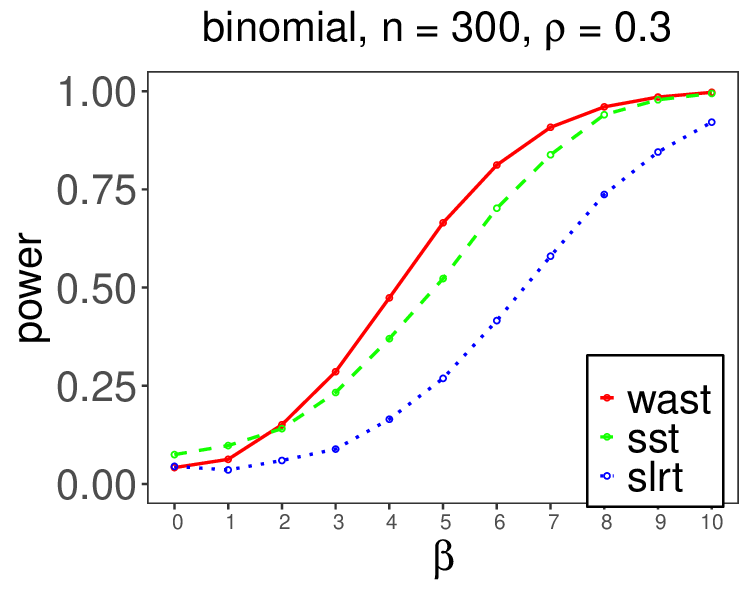}
		\includegraphics[scale=0.3]{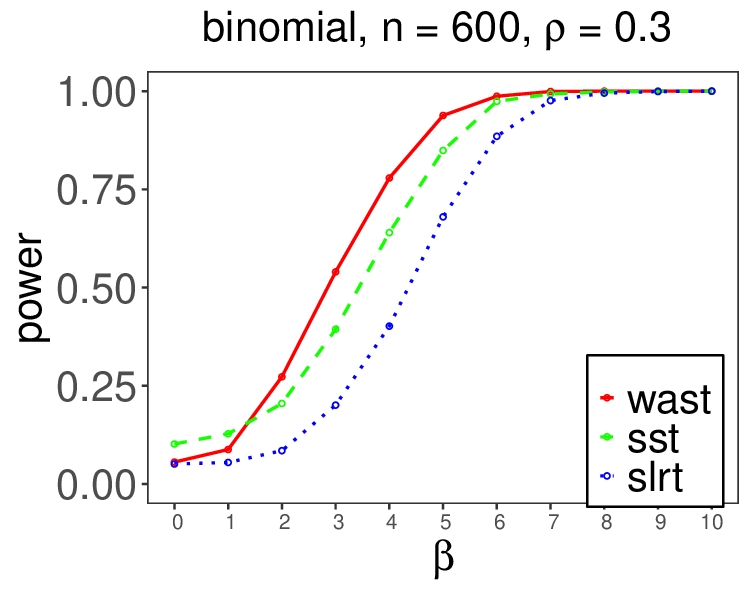}	
		\includegraphics[scale=0.3]{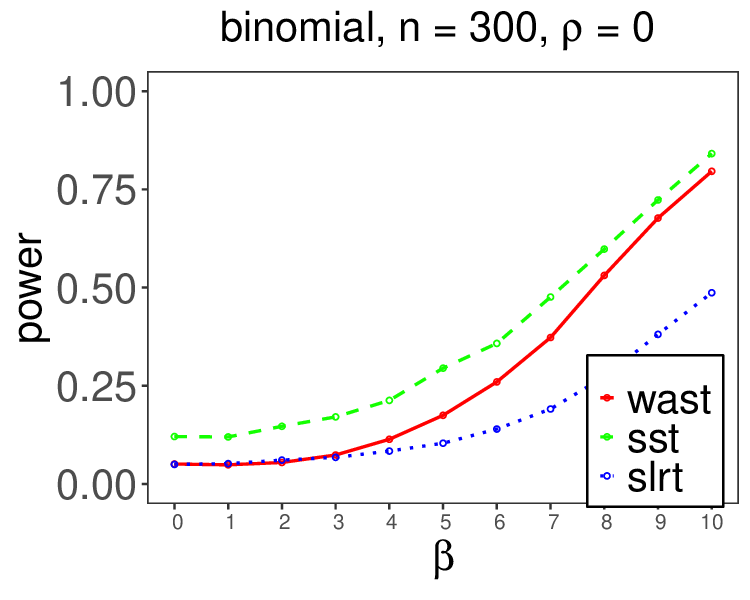}
		\includegraphics[scale=0.3]{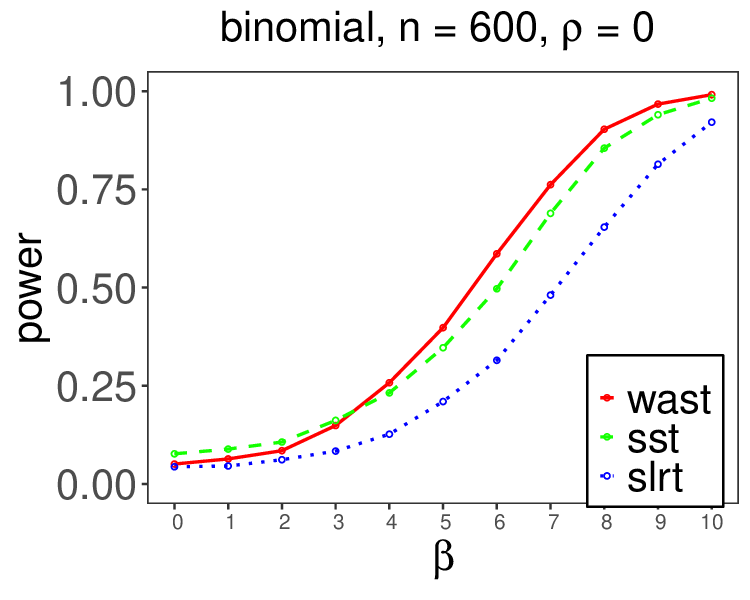}
		\includegraphics[scale=0.3]{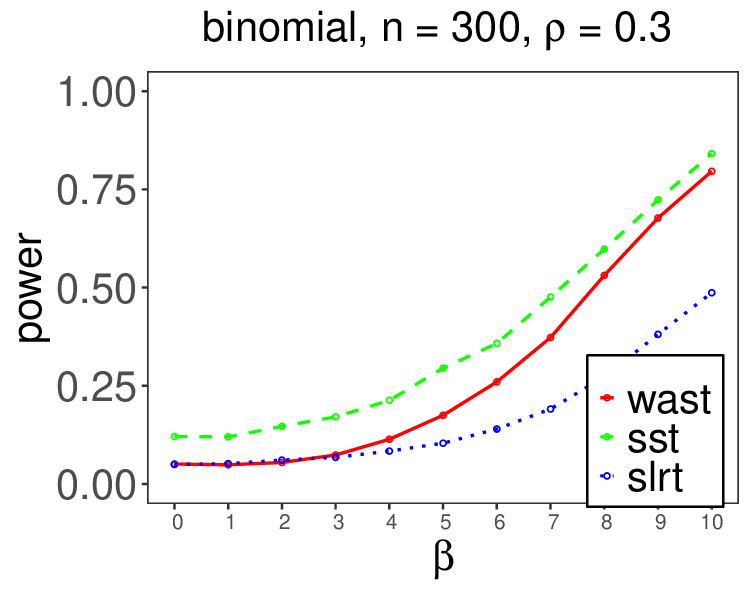}
		\includegraphics[scale=0.3]{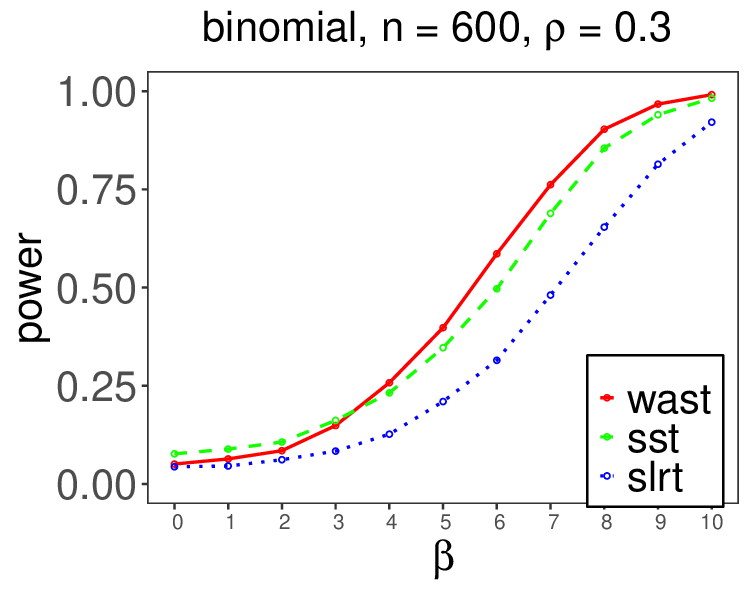}
		\includegraphics[scale=0.3]{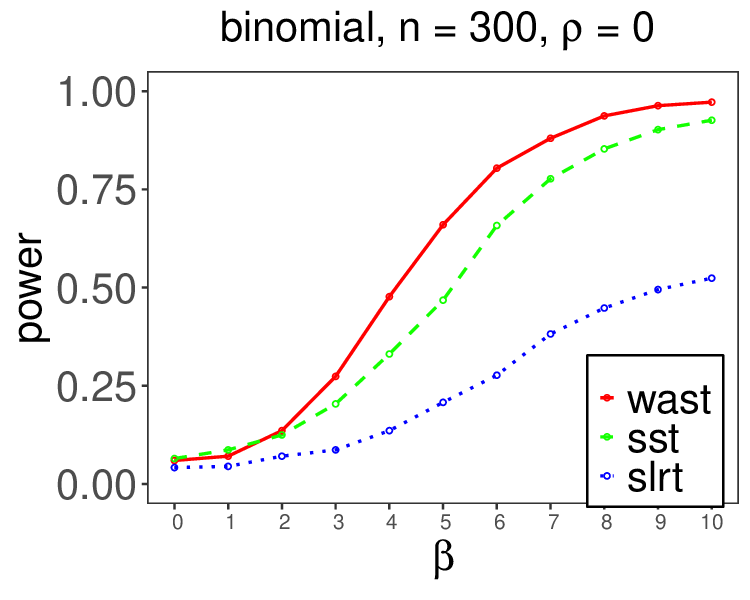}
		\includegraphics[scale=0.3]{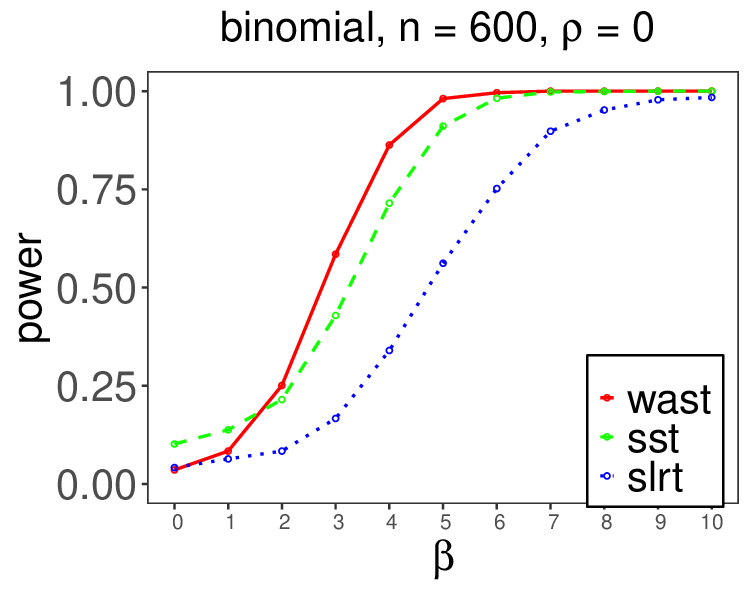}
		\includegraphics[scale=0.3]{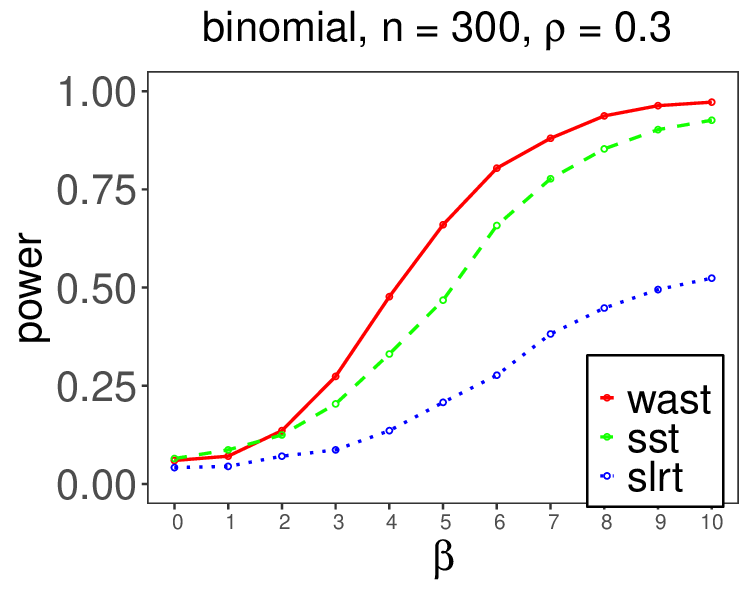}
		\includegraphics[scale=0.3]{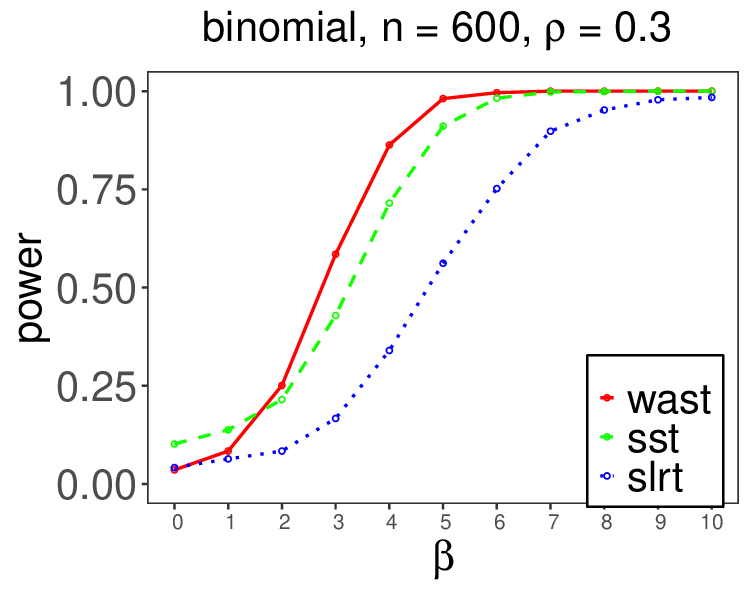}
		\includegraphics[scale=0.3]{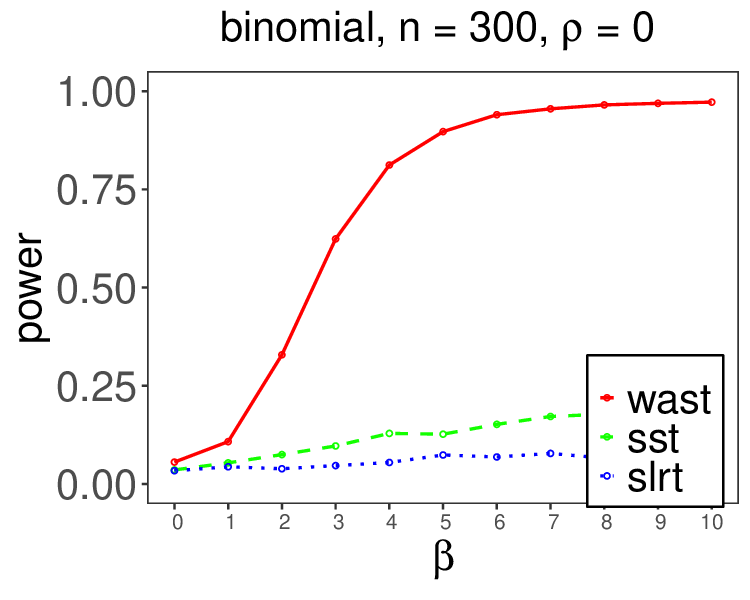}
		\includegraphics[scale=0.3]{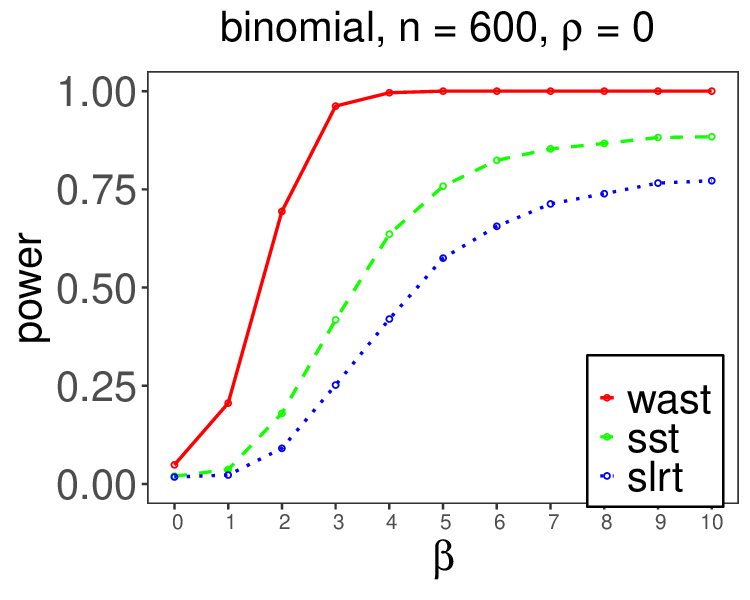}
		\includegraphics[scale=0.3]{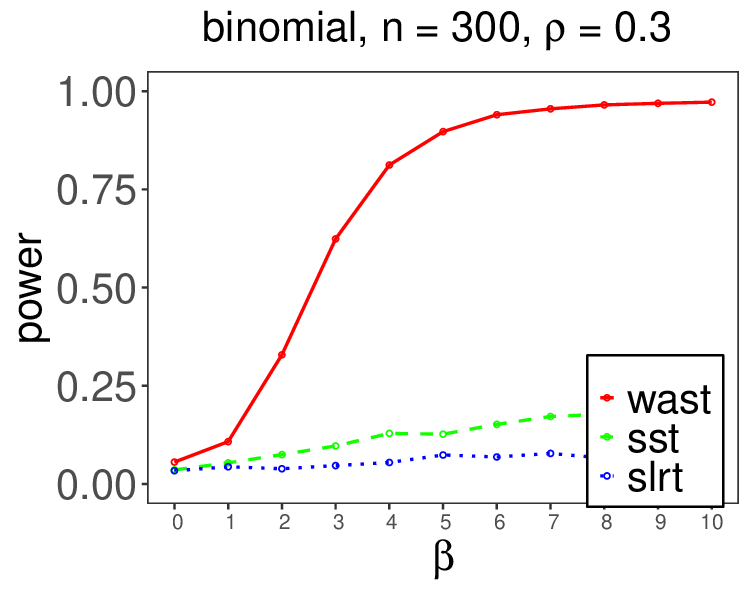}
		\includegraphics[scale=0.3]{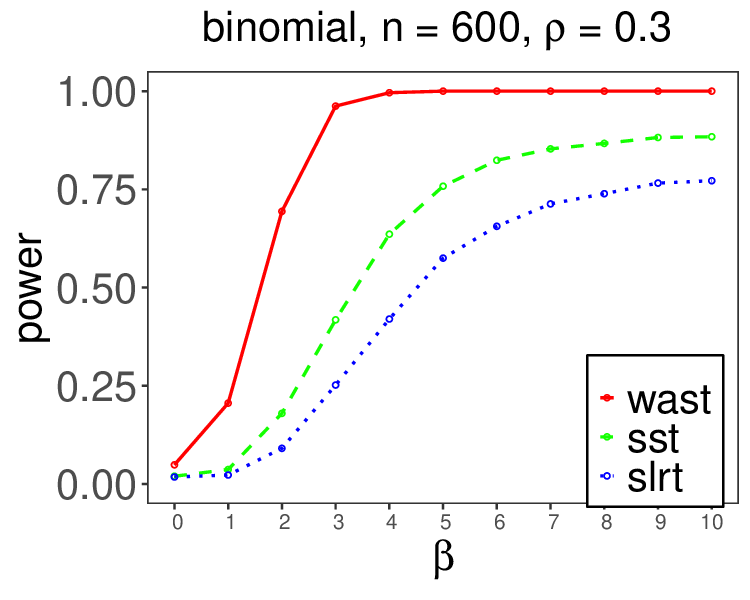}
		\includegraphics[scale=0.3]{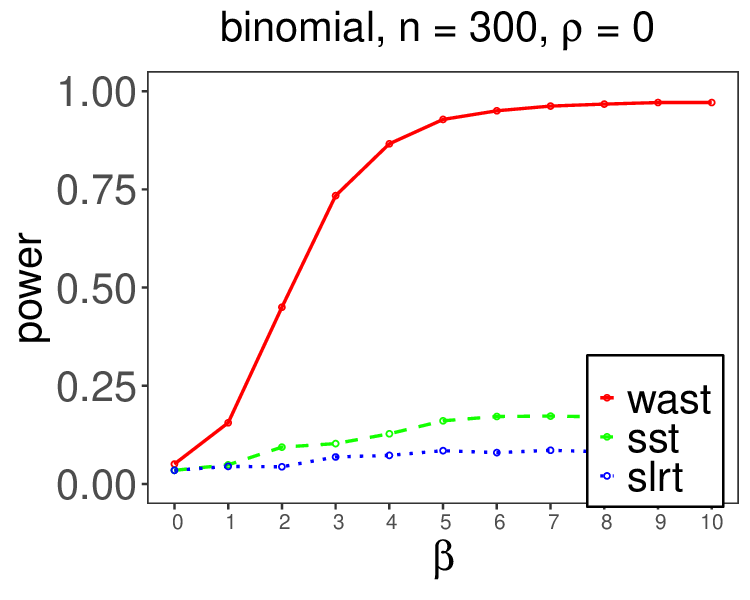}
		\includegraphics[scale=0.3]{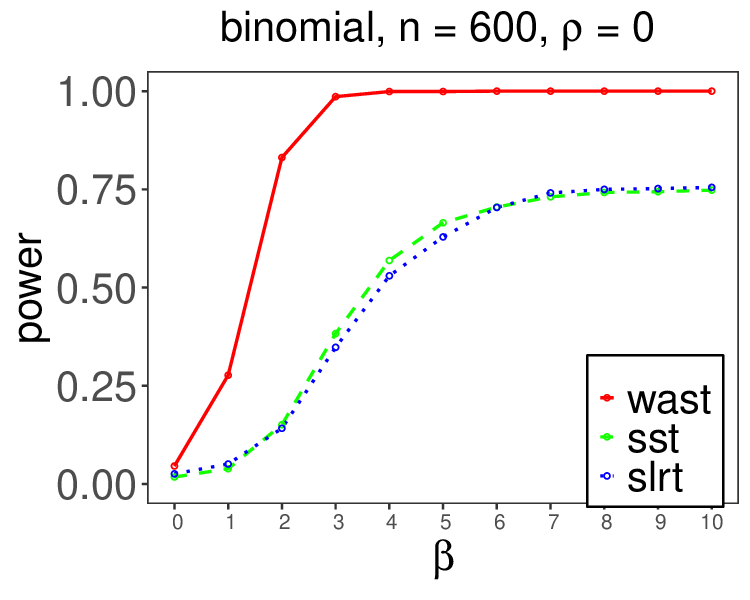}
		\includegraphics[scale=0.3]{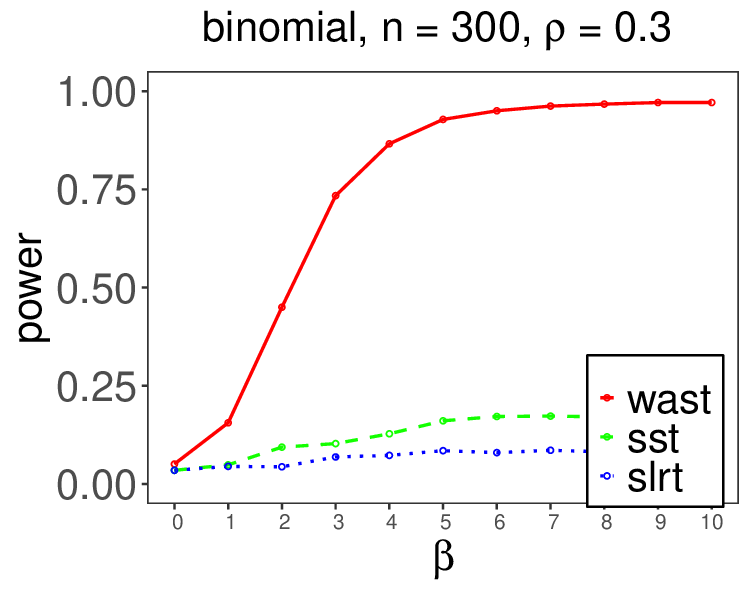}
		\includegraphics[scale=0.3]{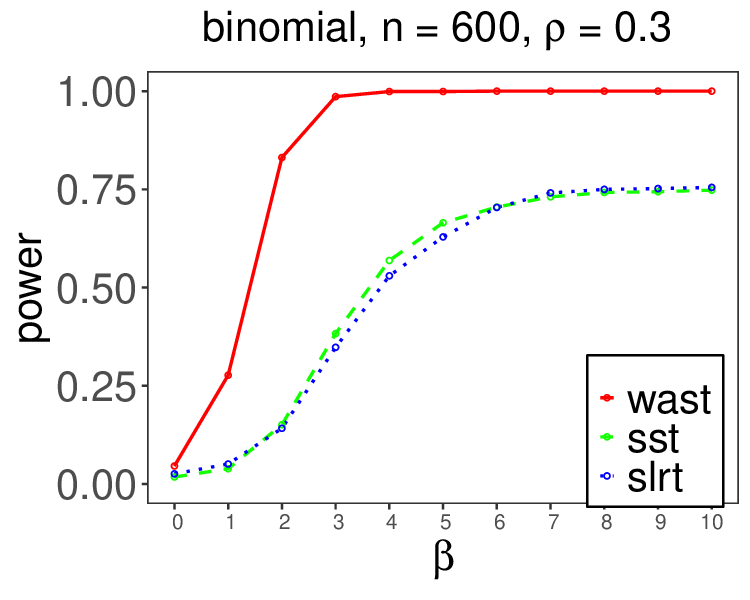}
		\caption{\it Powers of testing logistic regression by the proposed WAST (red solid line), SST (green dashed line), and SLRT (blue dotted line) for $n=(300,600)$. From top to bottom, each row panel depicts the powers for the case $(r,p,q)=(2,2,3)$, $(6,6,3)$, $(2,2,11)$, $(6,6,11)$, $(2,51,11)$, and $(6,51,11)$.
}
		\label{fig_binomial_zt3}
	\end{center}
\end{figure}

\begin{figure}[!ht]
	\begin{center}
		\includegraphics[scale=0.3]{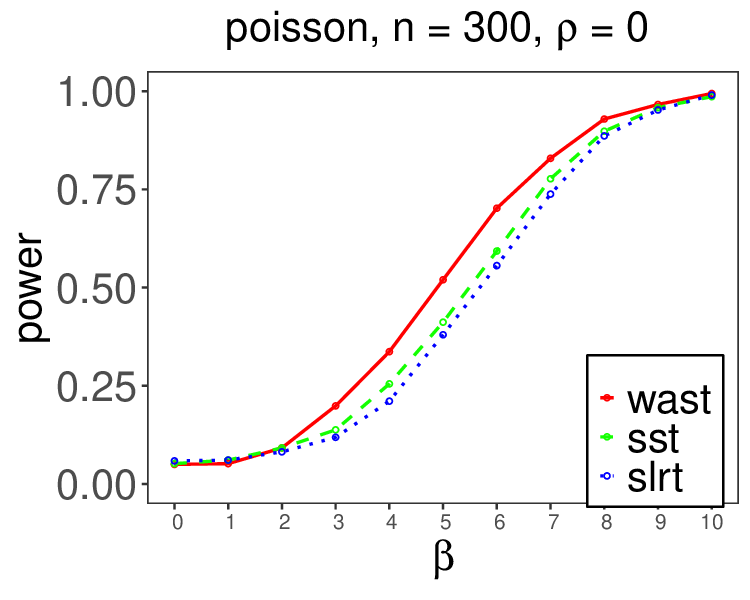}
		\includegraphics[scale=0.3]{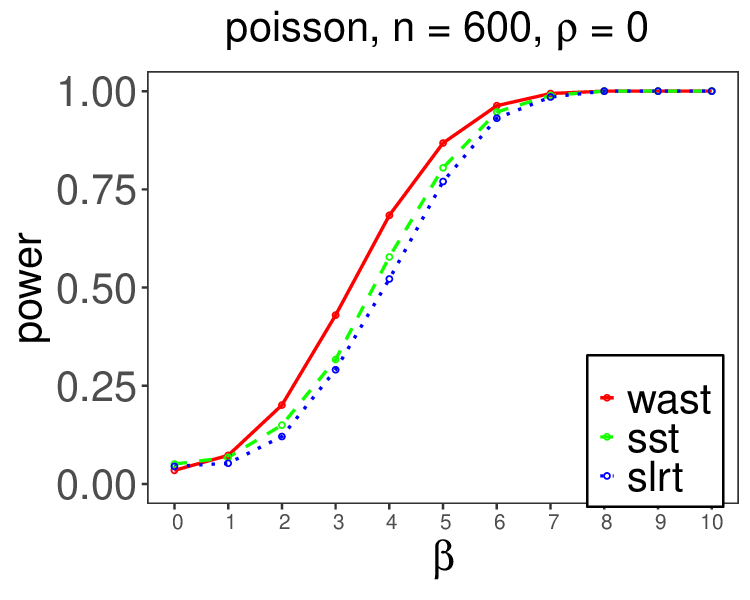}
		\includegraphics[scale=0.3]{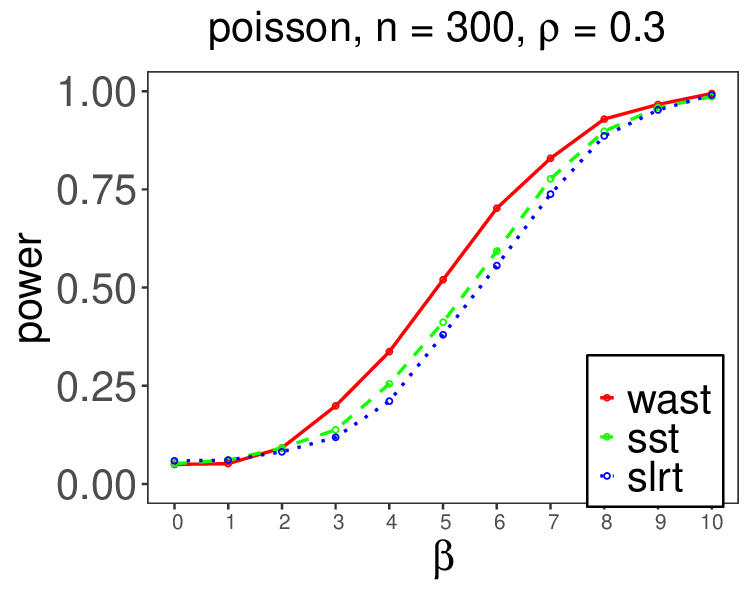}
		\includegraphics[scale=0.3]{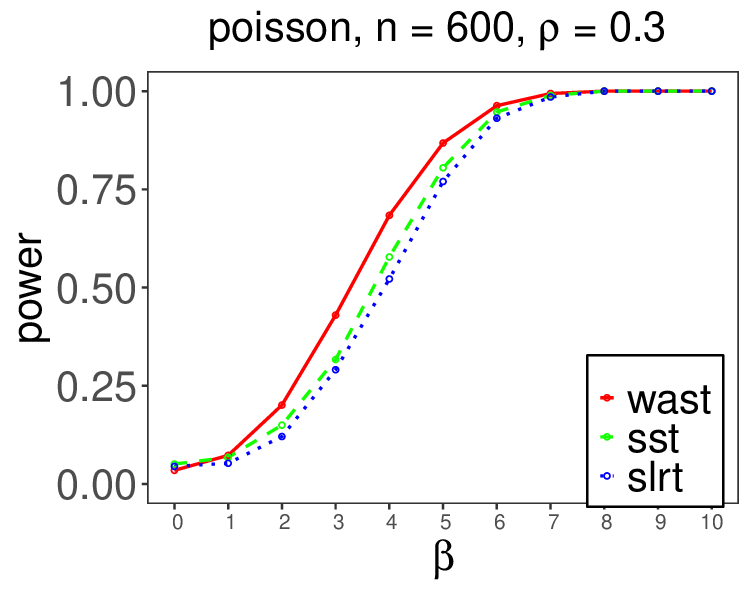}
		\includegraphics[scale=0.3]{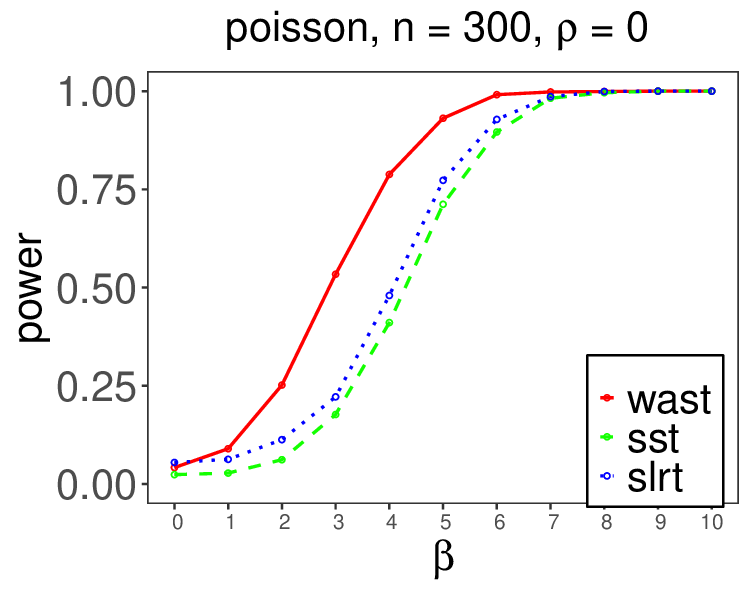}
		\includegraphics[scale=0.3]{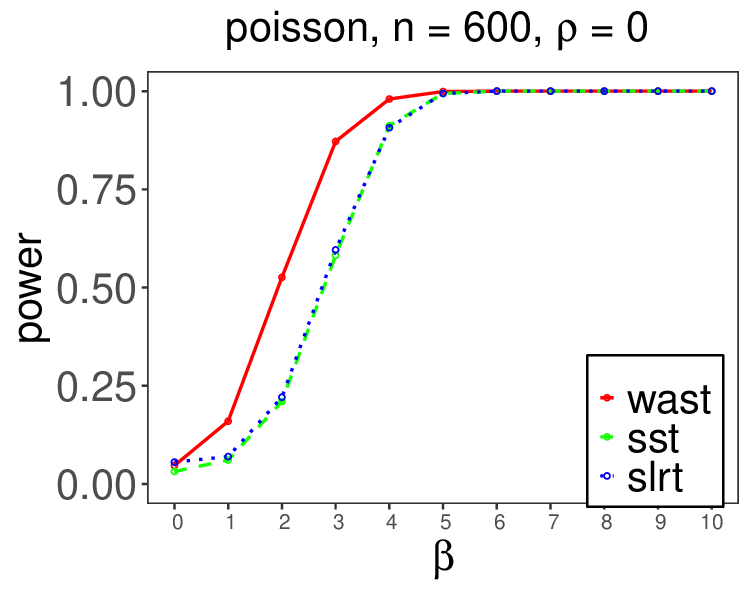}
		\includegraphics[scale=0.3]{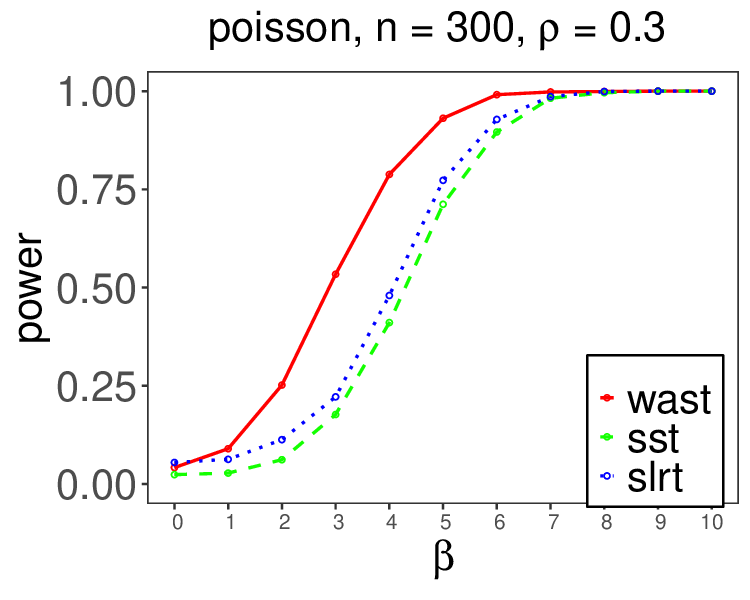}
		\includegraphics[scale=0.3]{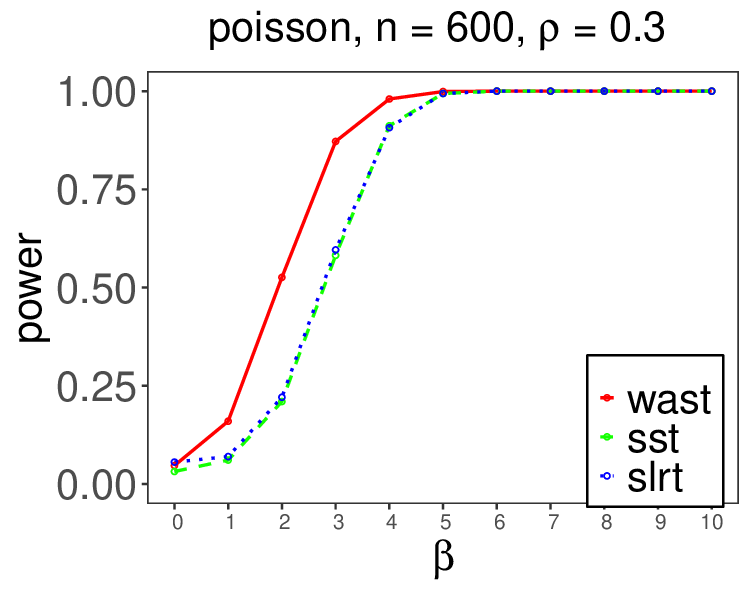}
		\includegraphics[scale=0.3]{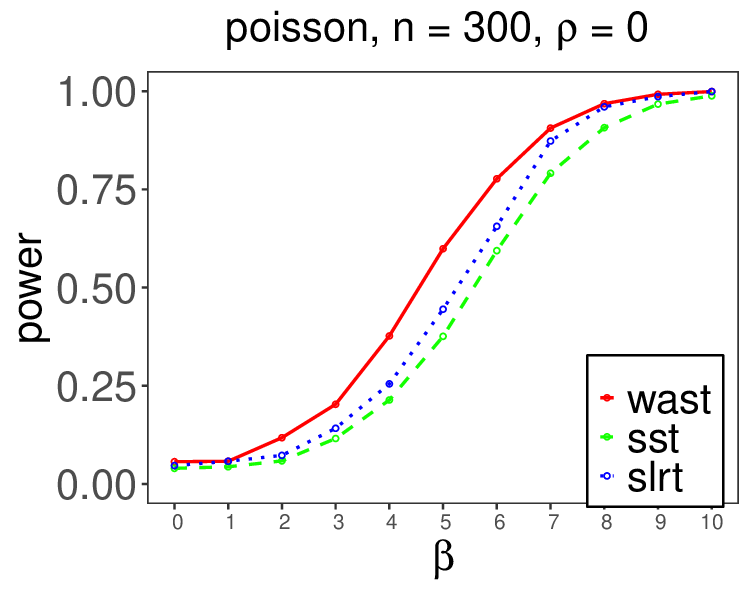}
		\includegraphics[scale=0.3]{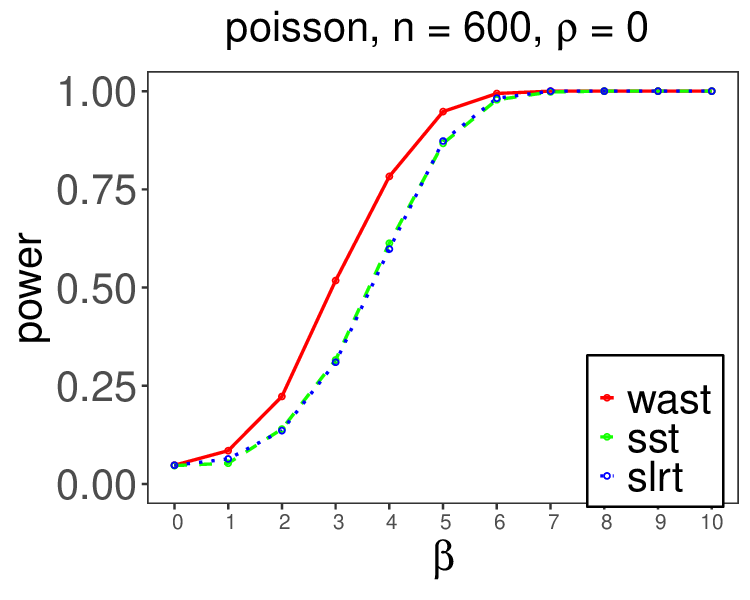}
		\includegraphics[scale=0.3]{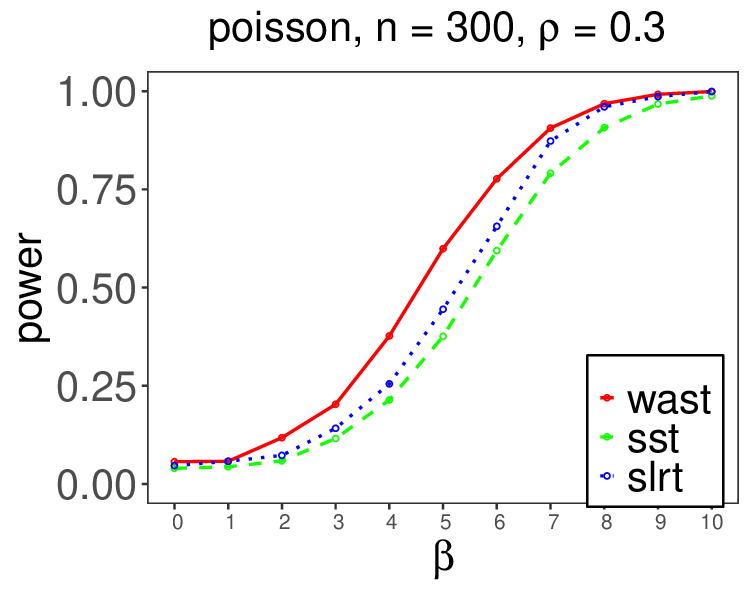}
		\includegraphics[scale=0.3]{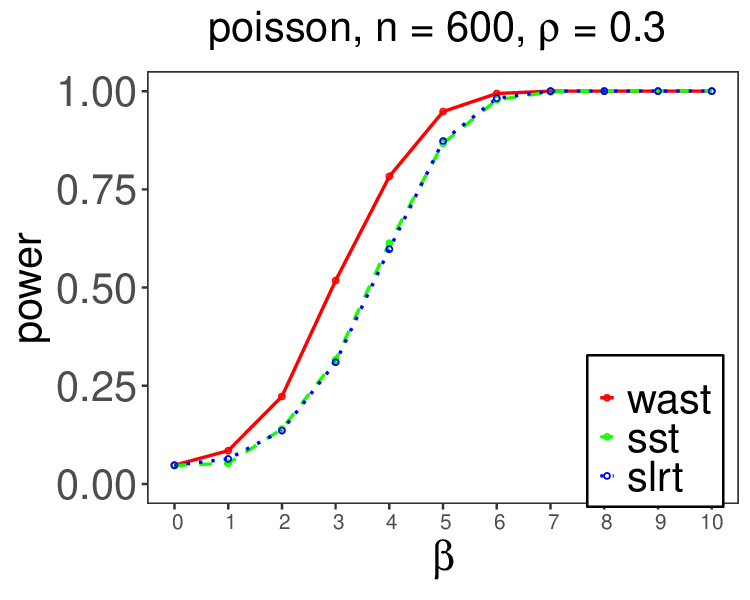}
		\includegraphics[scale=0.3]{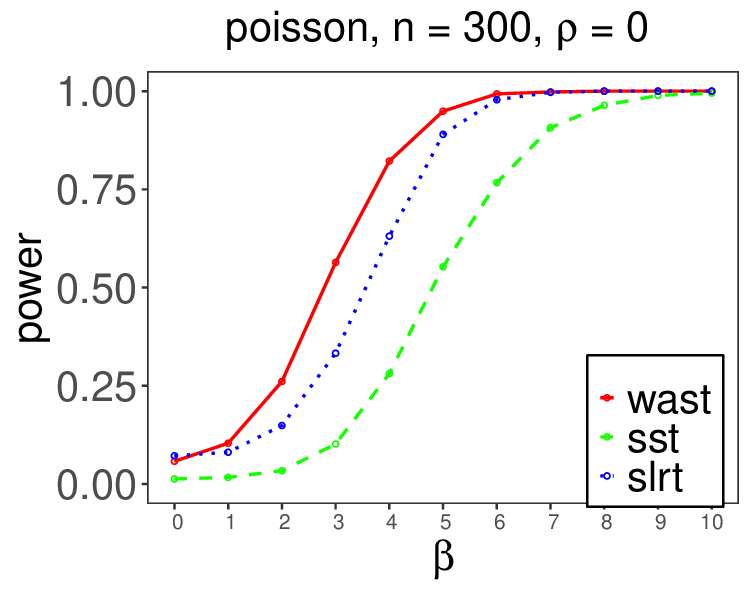}
		\includegraphics[scale=0.3]{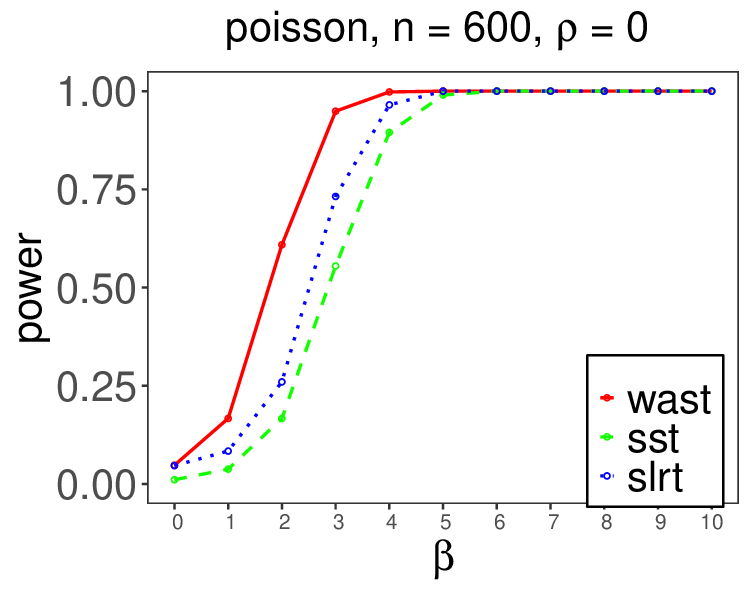}
		\includegraphics[scale=0.3]{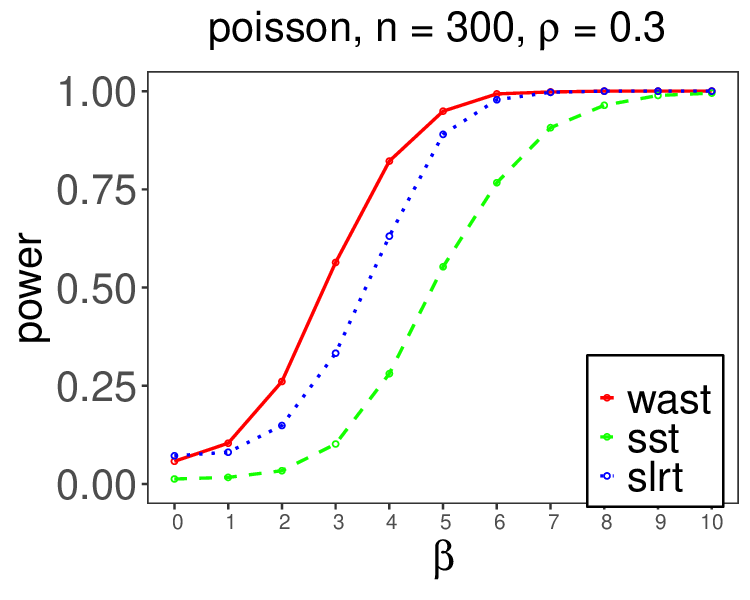}
		\includegraphics[scale=0.3]{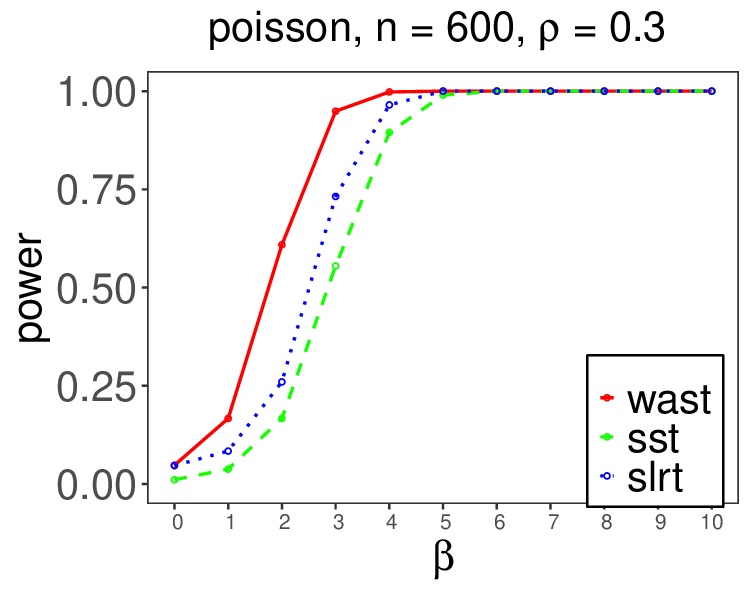}
		\includegraphics[scale=0.3]{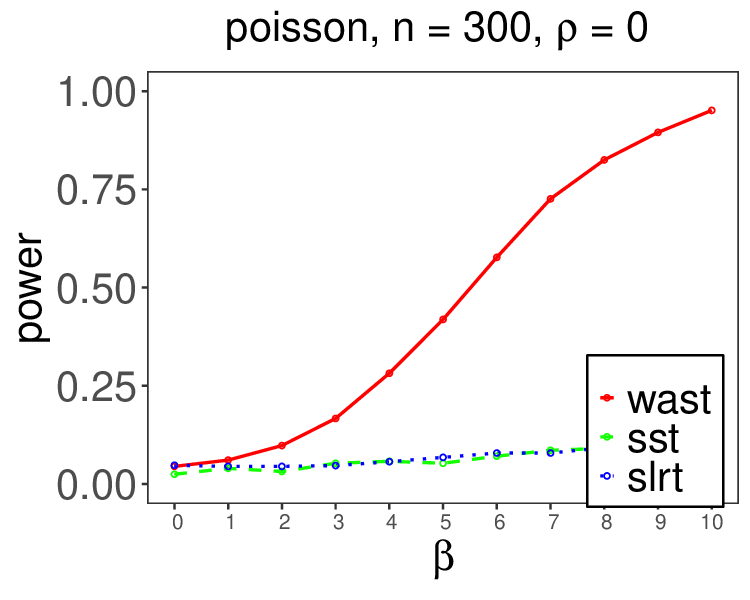}
		\includegraphics[scale=0.3]{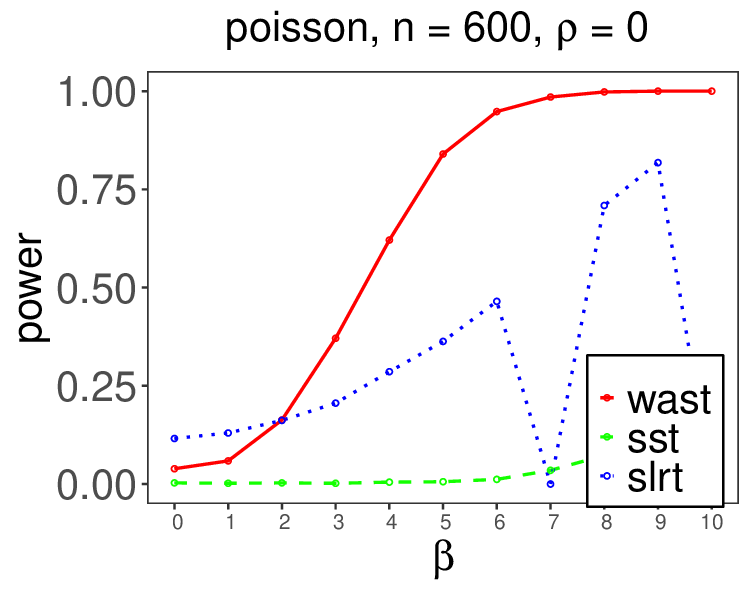}
		\includegraphics[scale=0.3]{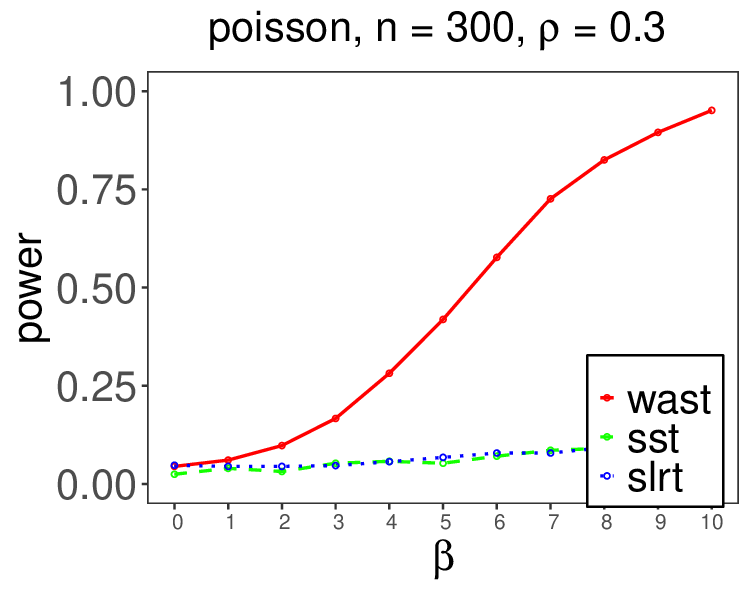}
		\includegraphics[scale=0.3]{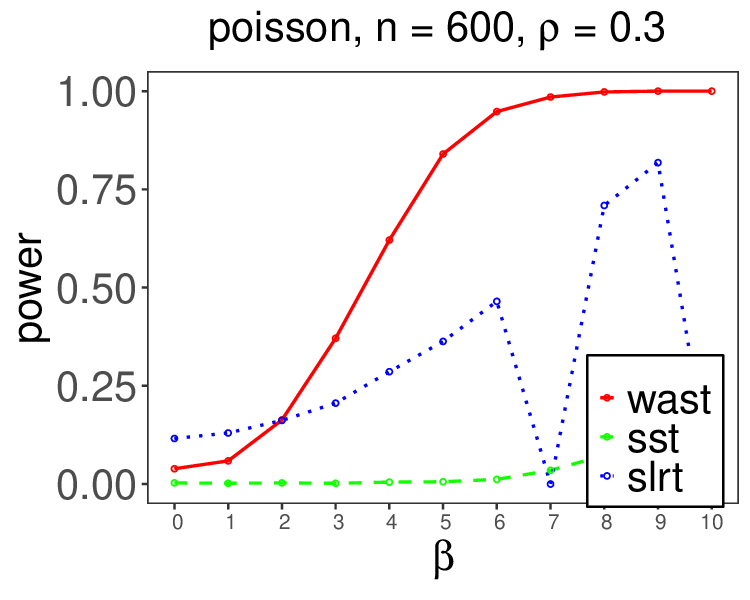}
		\includegraphics[scale=0.3]{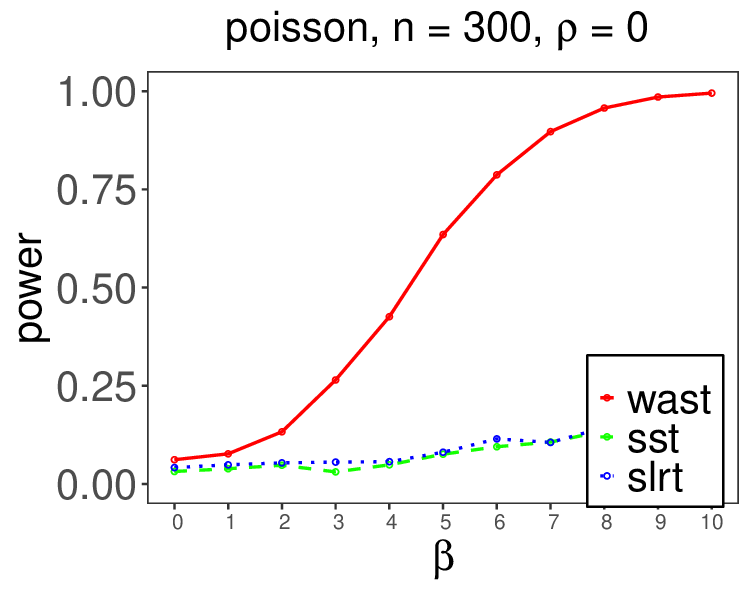}
		\includegraphics[scale=0.3]{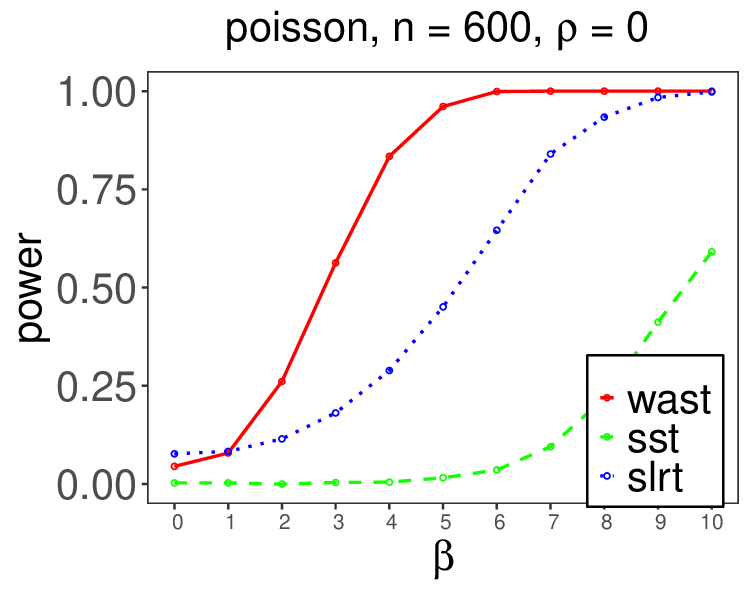}
		\includegraphics[scale=0.3]{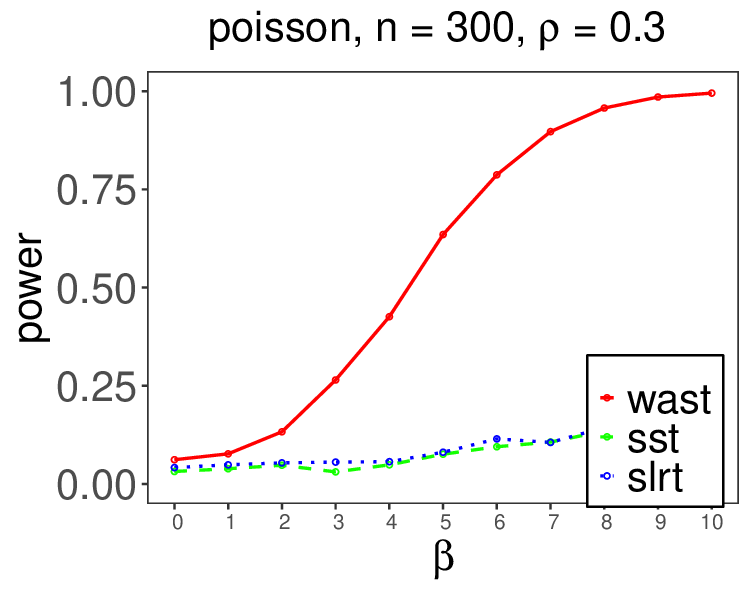}
		\includegraphics[scale=0.3]{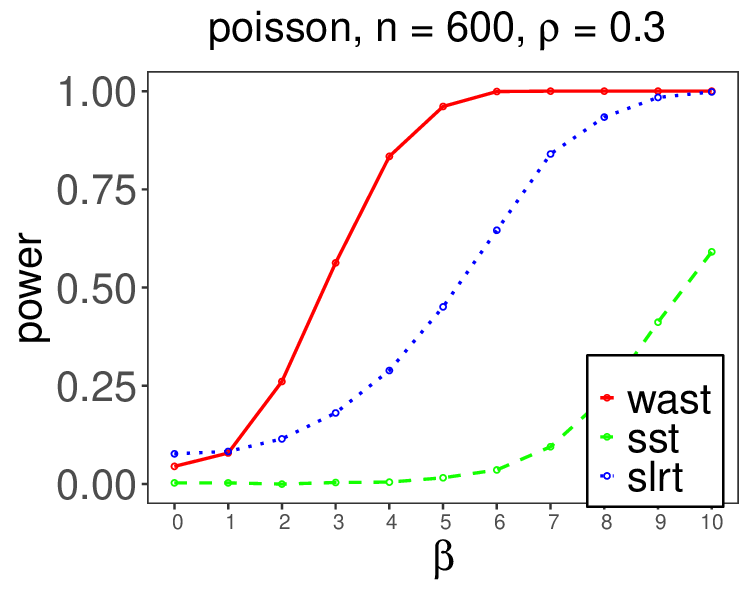}
		\caption{\it Powers of testing Poisson regression by the proposed WAST (red solid line), SST (green dashed line), and SLRT (blue dotted line) for $n=(300,600)$. From top to bottom, each row panel depicts the powers for the case $(r,p,q)=(2,2,3)$, $(6,6,3)$, $(2,2,11)$, $(6,6,11)$, $(2,51,11)$, and $(6,51,11)$.
}
		\label{fig_poisson_zt3}
	\end{center}
\end{figure}

The power curves of the WAST, SST and SLRT are depicted in Figures \ref{fig_gaussian_zt3}-\ref{fig_poisson_zt3}. From Figures \ref{fig_gaussian_zt3}-\ref{fig_poisson_zt3}, we see that the proposed method have higher power uniformly than those of the SLRT and SST, which illustrates the robustness of the proposed method and motivates researchers to apply the proposed method in practice.

\begin{figure}[!ht]
	\begin{center}
		\includegraphics[scale=0.3]{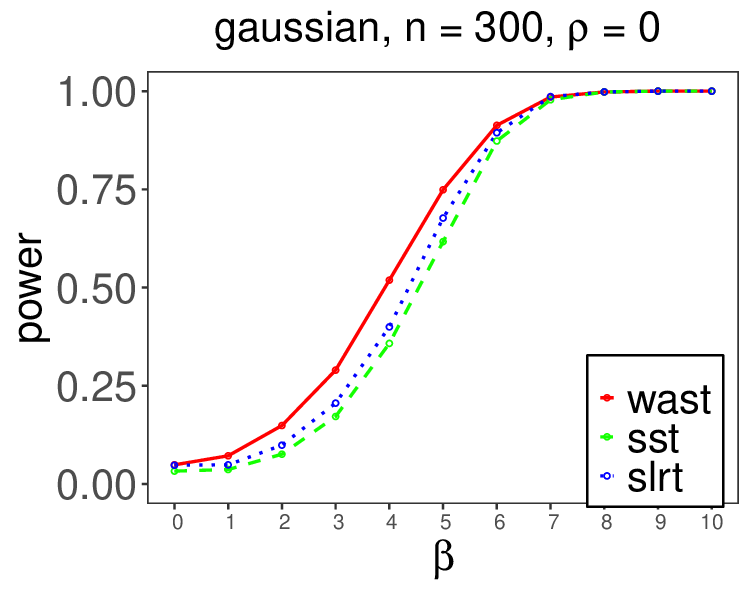}
		\includegraphics[scale=0.3]{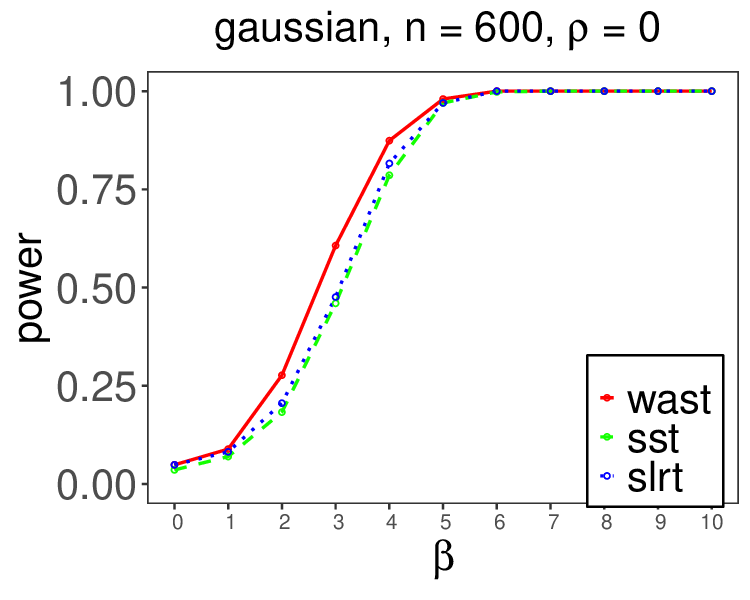}
		\includegraphics[scale=0.3]{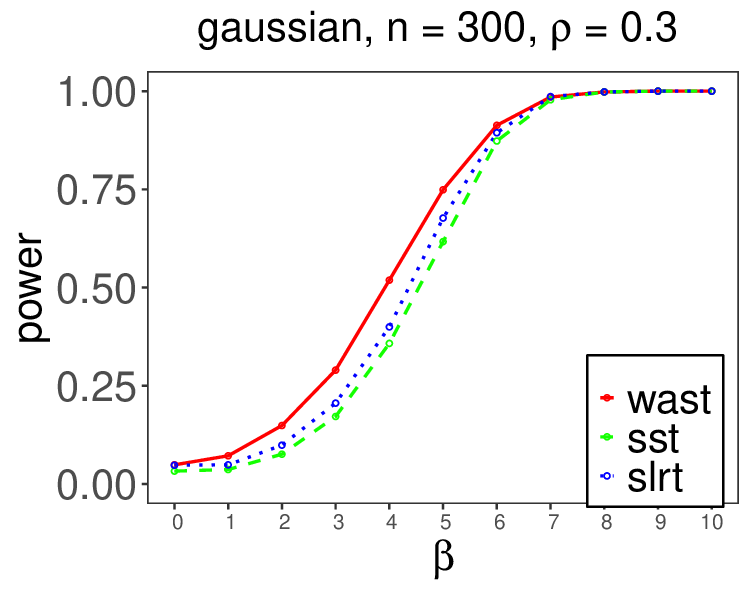}
		\includegraphics[scale=0.3]{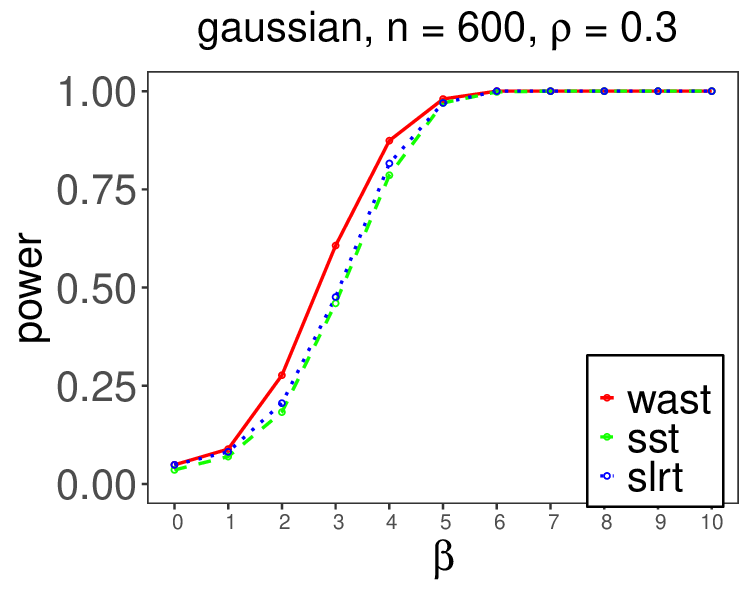}
		\includegraphics[scale=0.3]{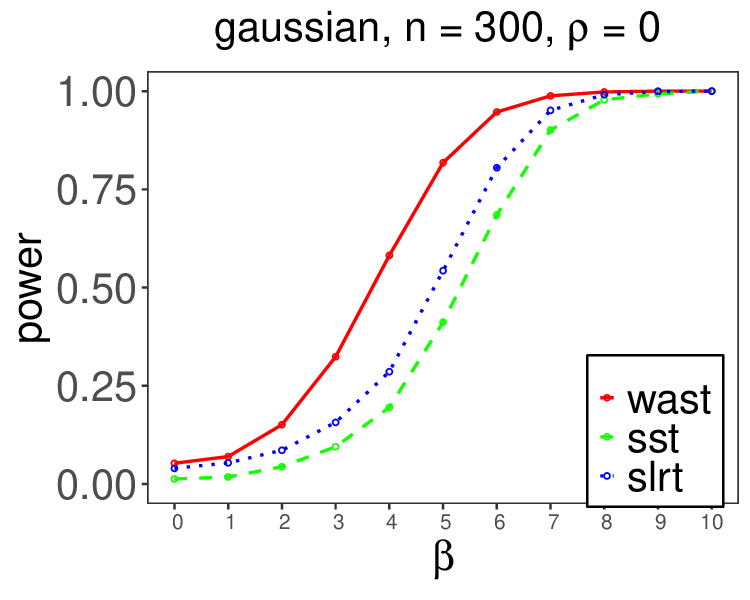}
		\includegraphics[scale=0.3]{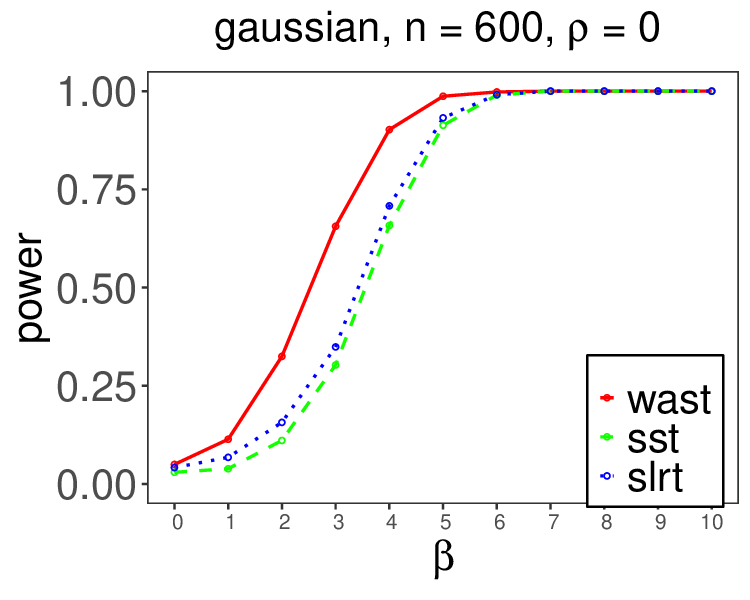}
		\includegraphics[scale=0.3]{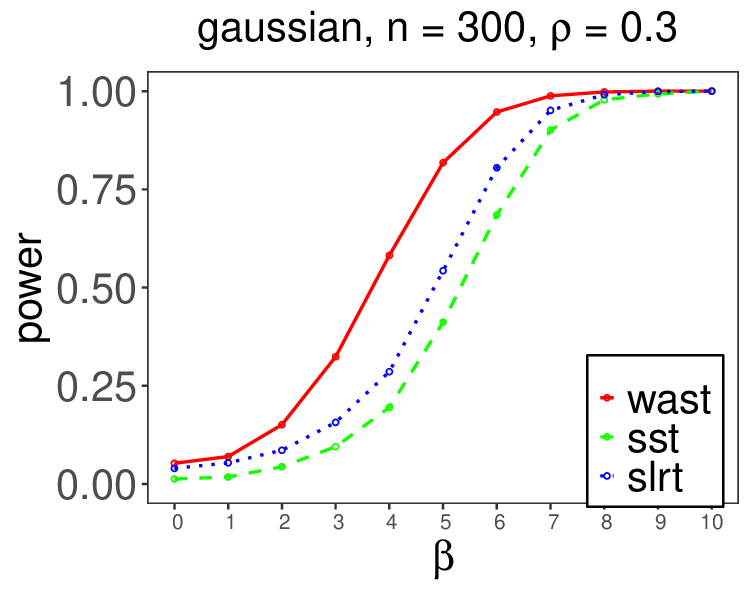}
		\includegraphics[scale=0.3]{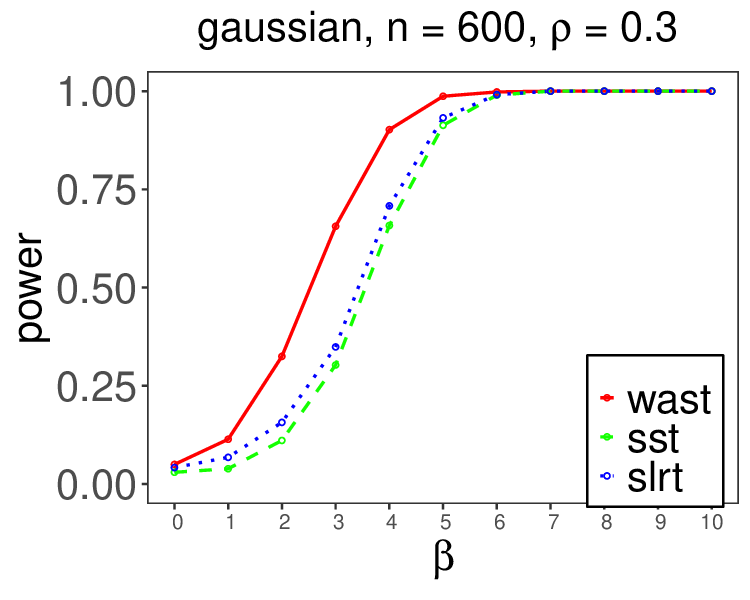} 		
		\includegraphics[scale=0.3]{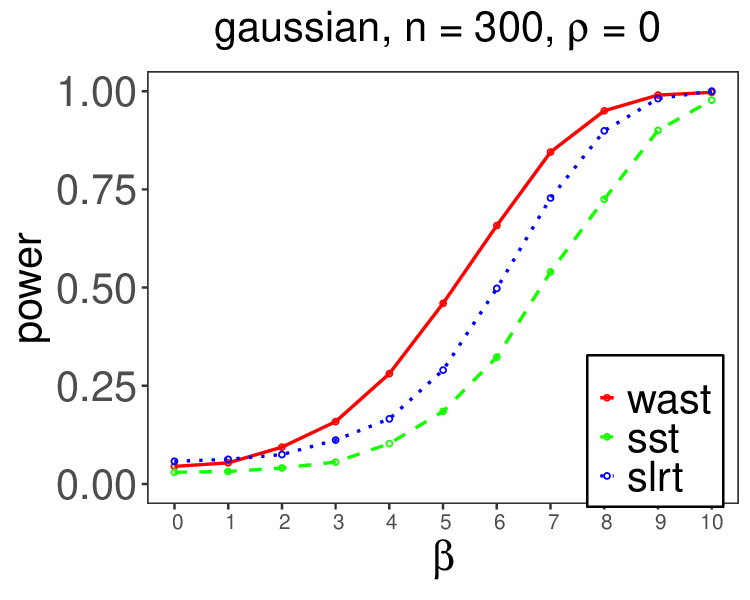}
		\includegraphics[scale=0.3]{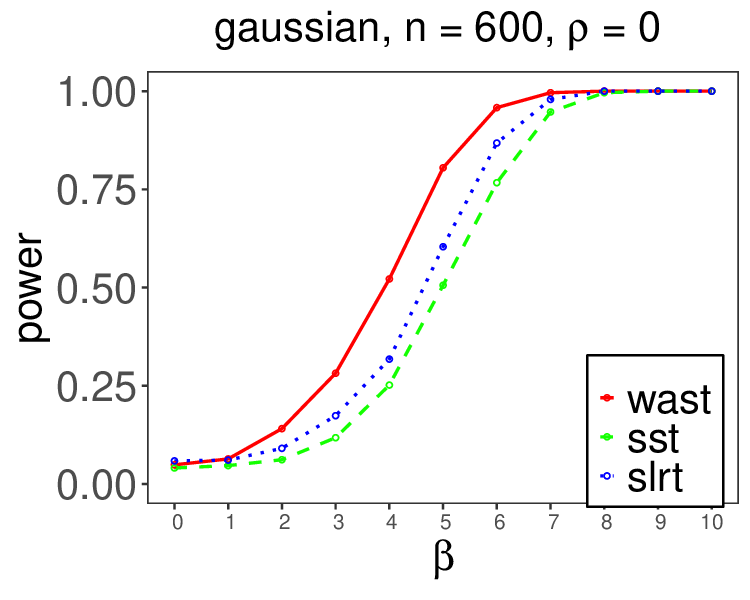}
		\includegraphics[scale=0.3]{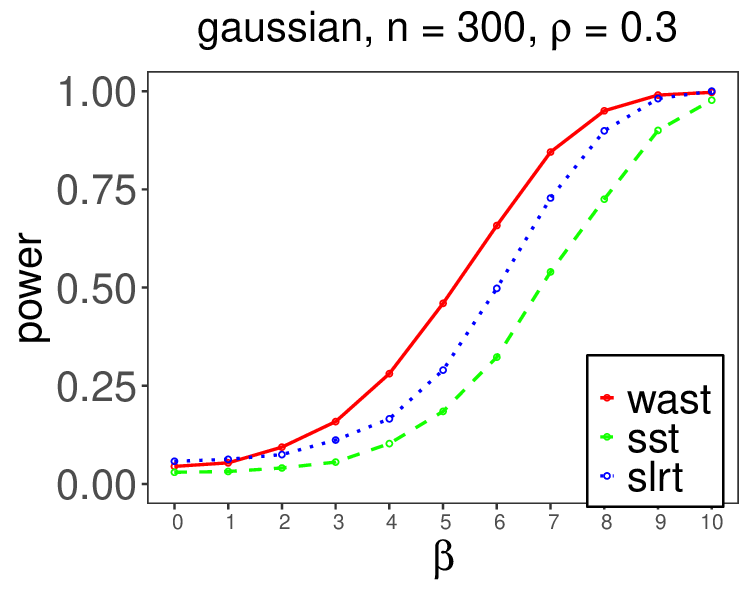}
		\includegraphics[scale=0.3]{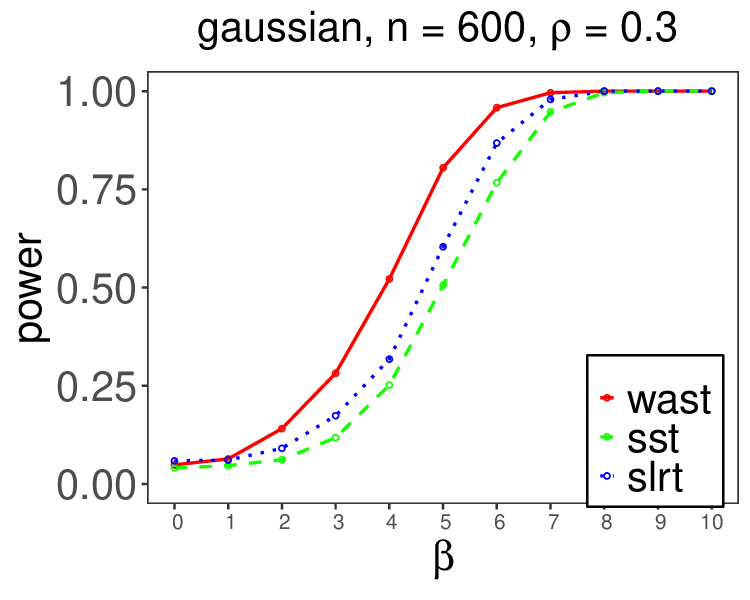}
		\includegraphics[scale=0.3]{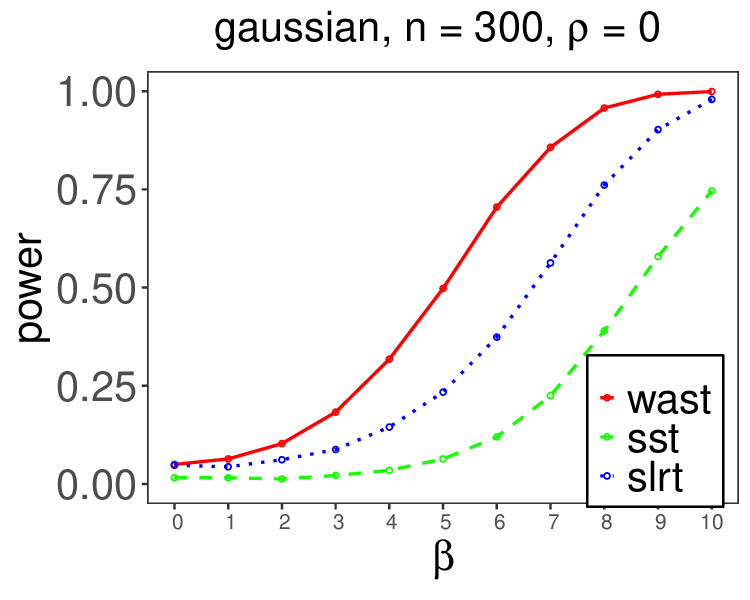}
		\includegraphics[scale=0.3]{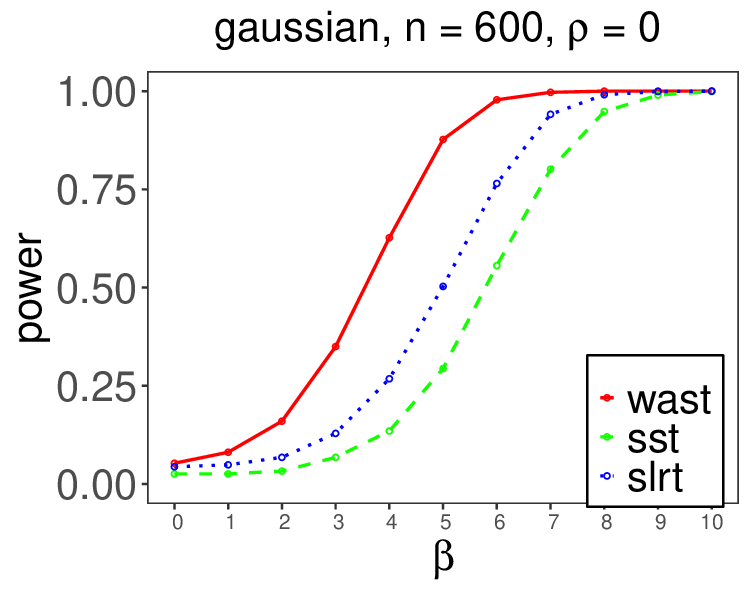}
		\includegraphics[scale=0.3]{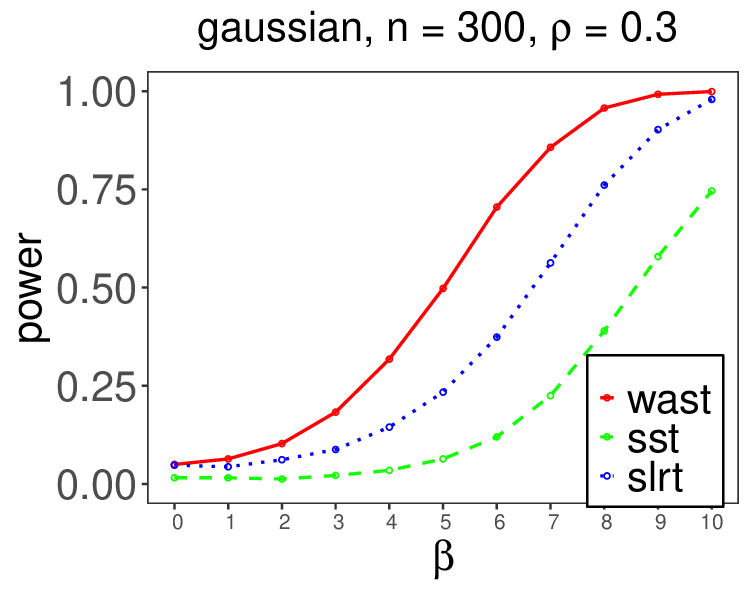}
		\includegraphics[scale=0.3]{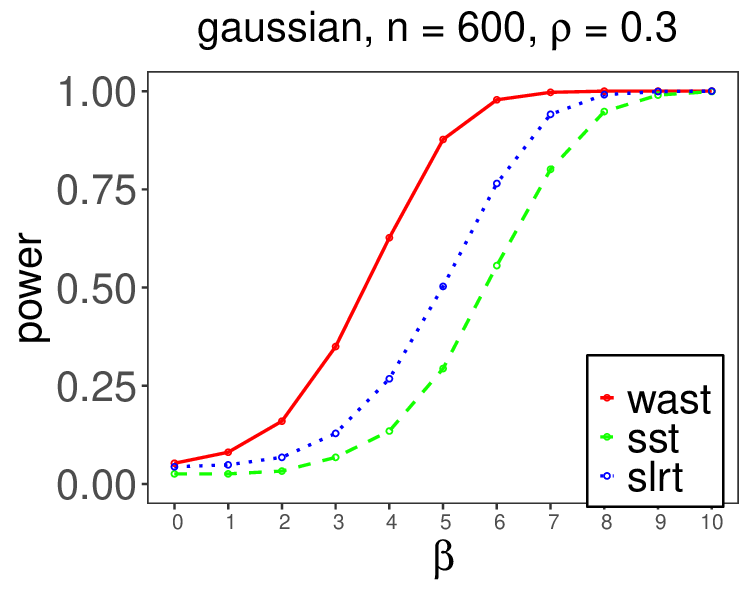}
		\includegraphics[scale=0.3]{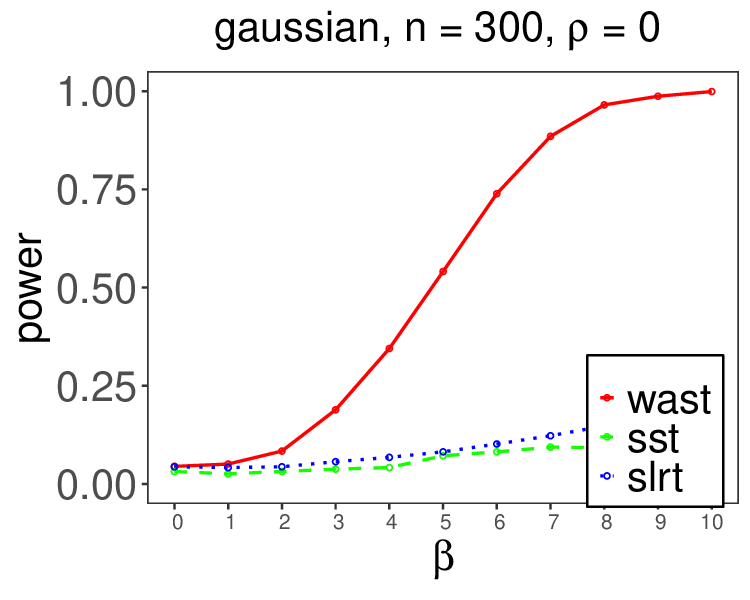}
		\includegraphics[scale=0.3]{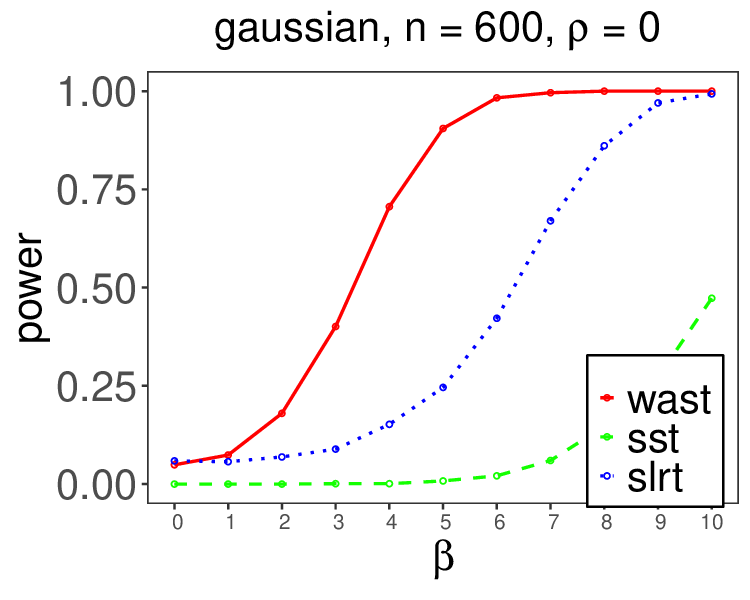}
		\includegraphics[scale=0.3]{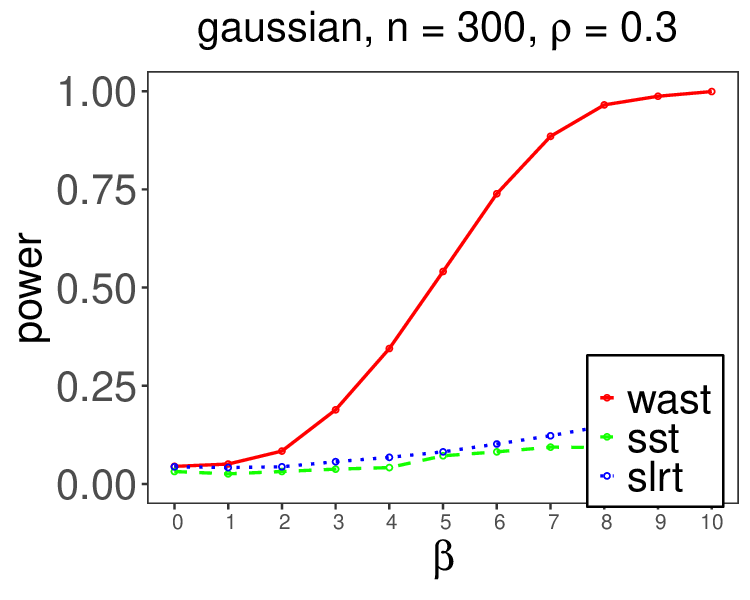}
		\includegraphics[scale=0.3]{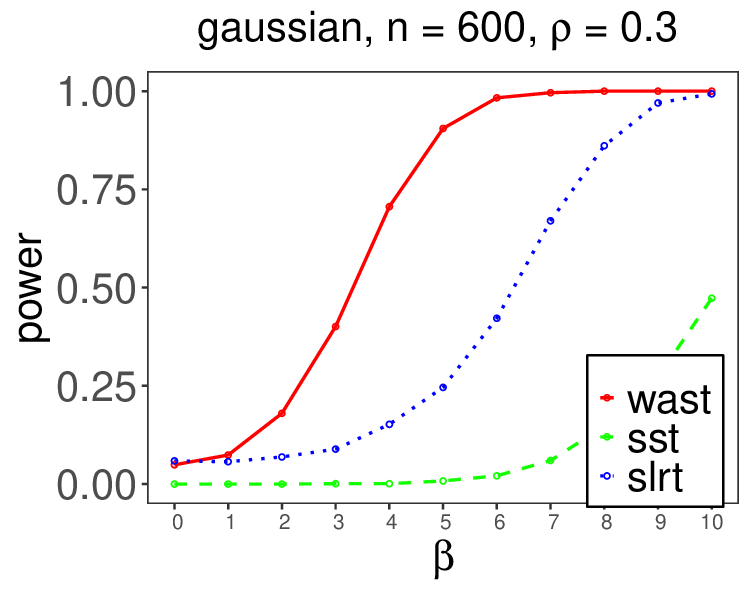}
		\includegraphics[scale=0.3]{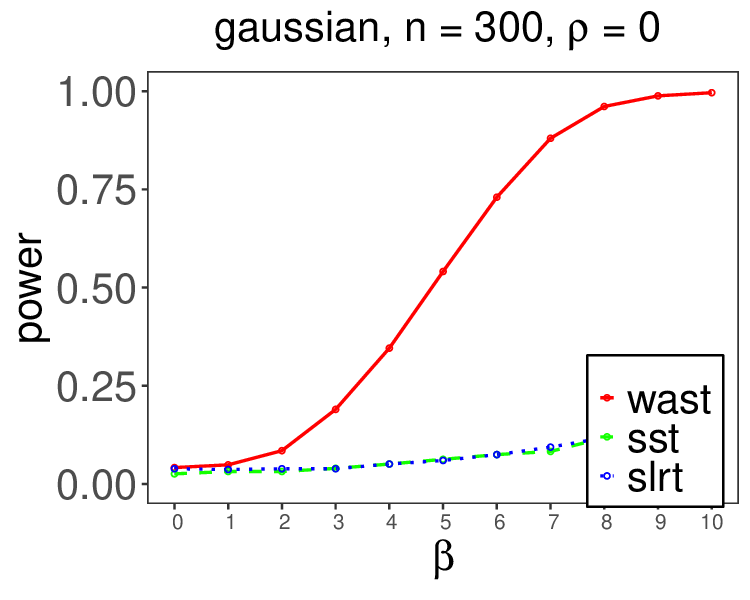}
		\includegraphics[scale=0.3]{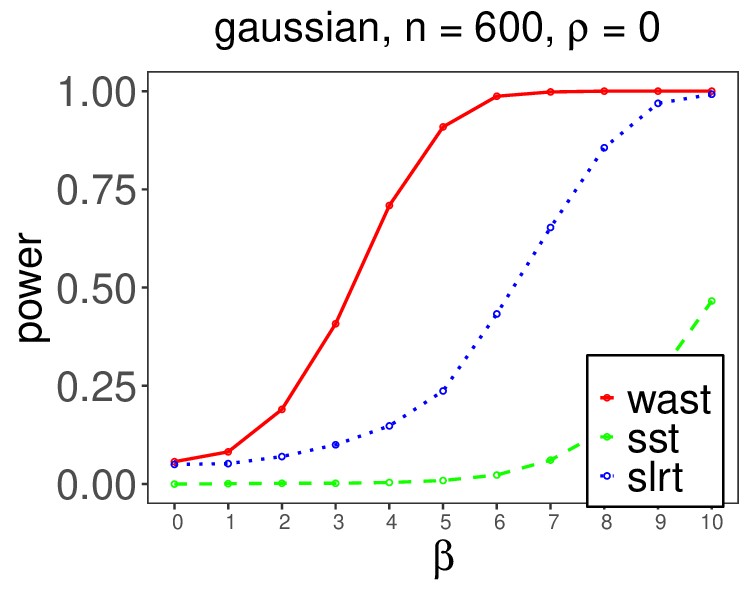}
		\includegraphics[scale=0.3]{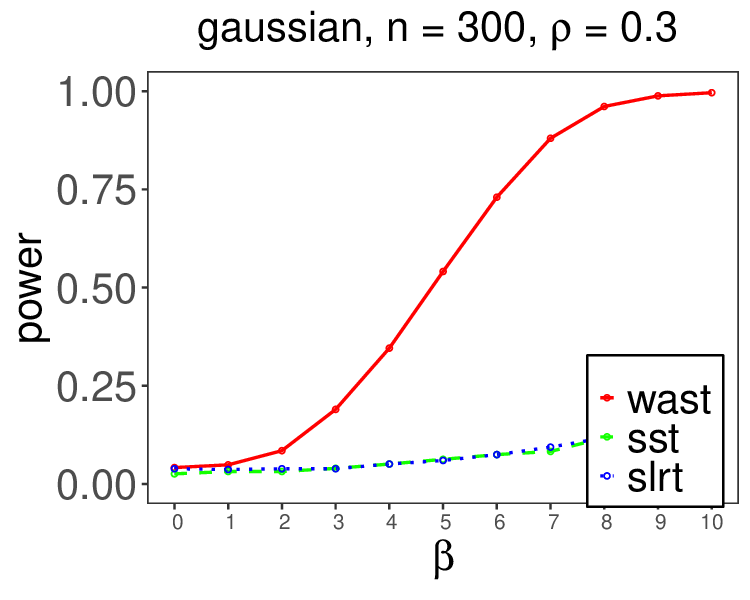}
		\includegraphics[scale=0.3]{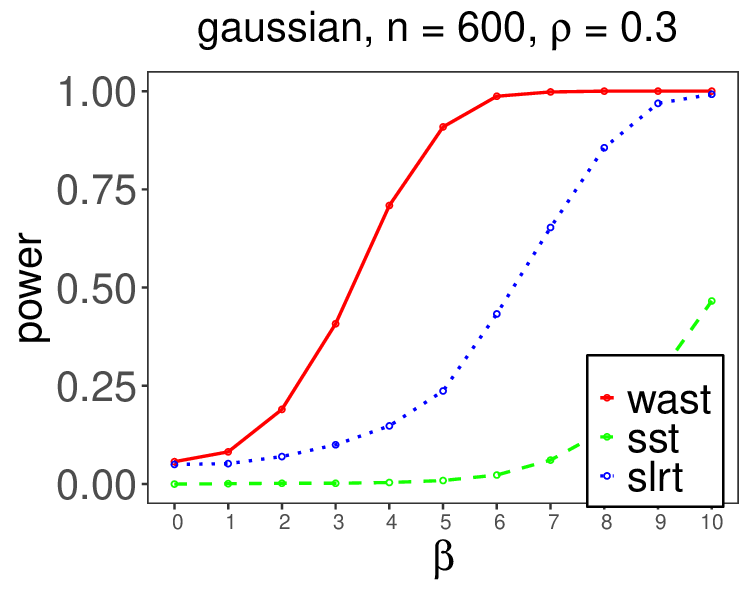}

		\caption{\it Powers of testing linear model with Gaussian error by the proposed WAST (red solid line), SST (green dashed line), and SLRT (blue dotted line) for $n=(300,600)$. From top to bottom, each row panel depicts the powers for the case $(r,p,q)=(2,2,3)$, $(6,6,3)$, $(2,2,11)$, $(6,6,11)$, $(2,51,11)$, and $(6,51,11)$.
}
		\label{fig_gaussian_zs4}
	\end{center}
\end{figure}

\begin{figure}[!ht]
	\begin{center}
		\includegraphics[scale=0.3]{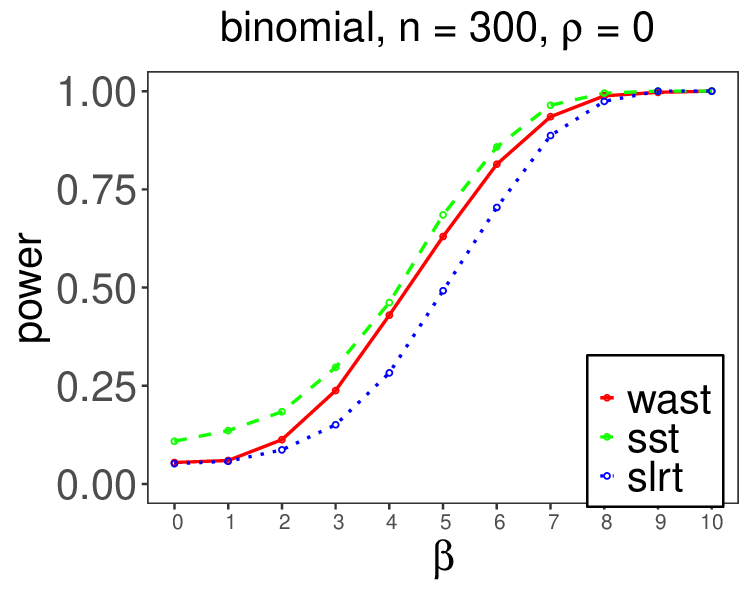}
		\includegraphics[scale=0.3]{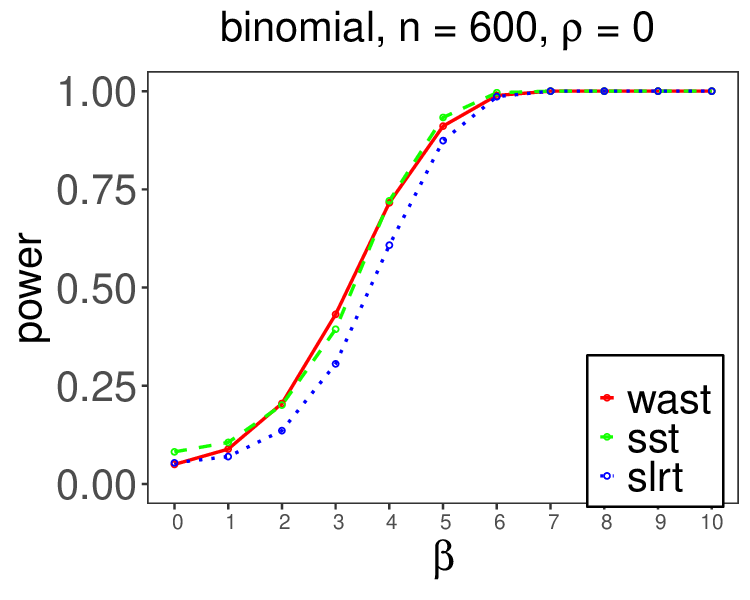}
		\includegraphics[scale=0.3]{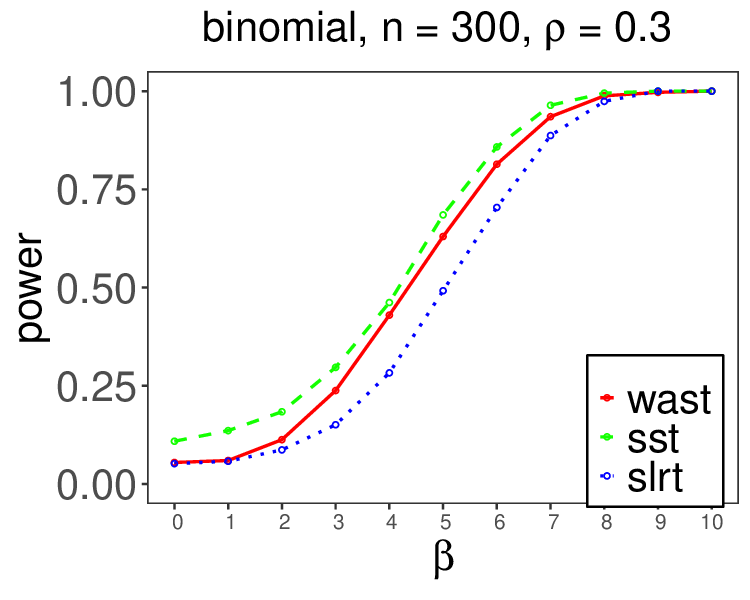}
		\includegraphics[scale=0.3]{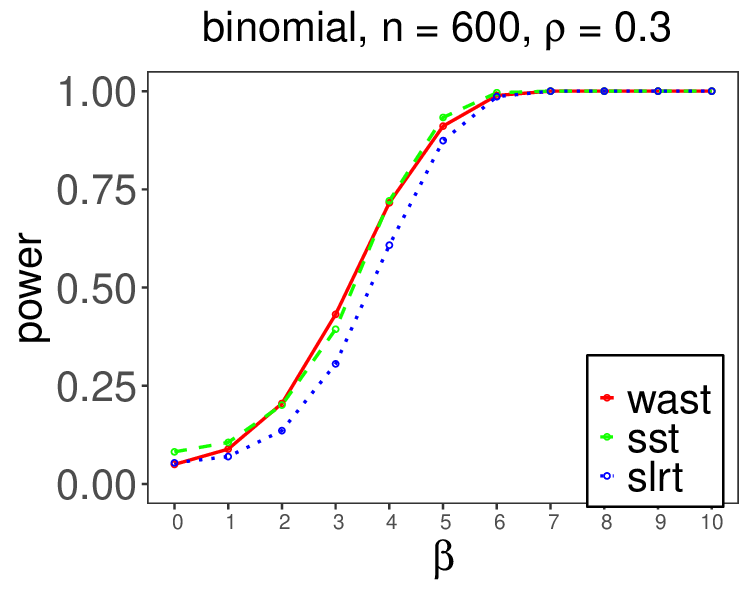}
		\includegraphics[scale=0.3]{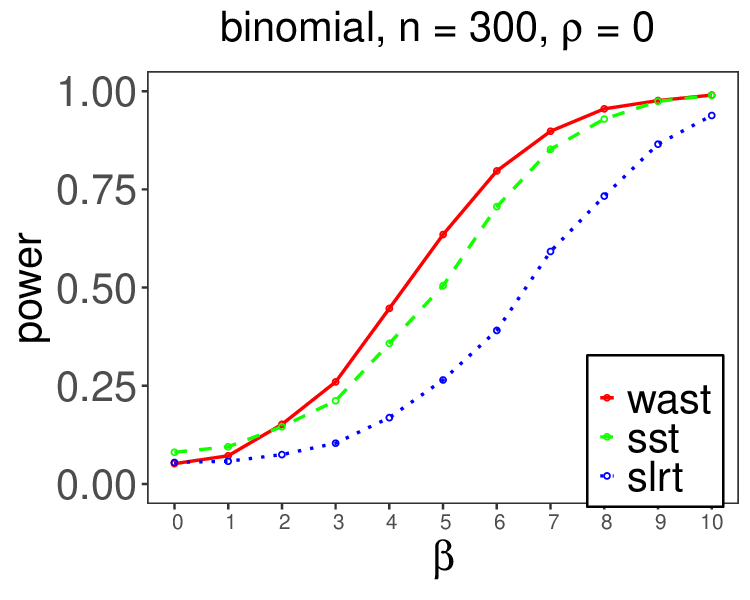}
		\includegraphics[scale=0.3]{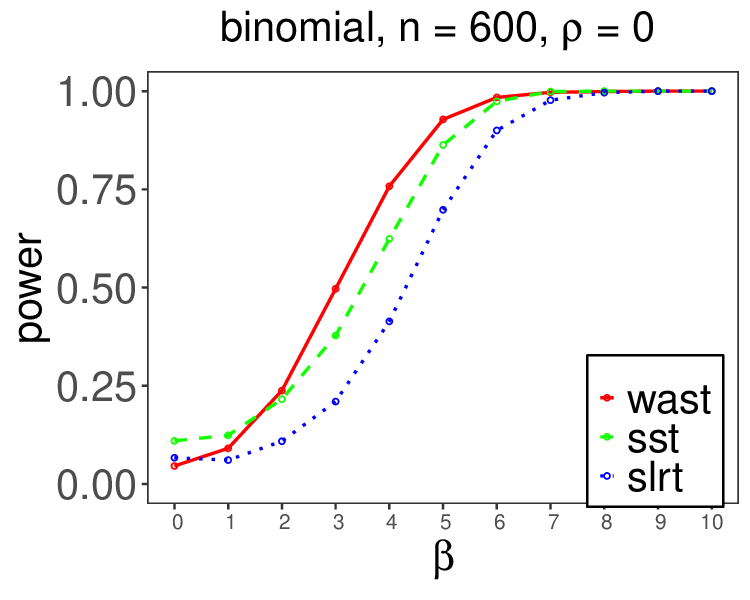}
		\includegraphics[scale=0.3]{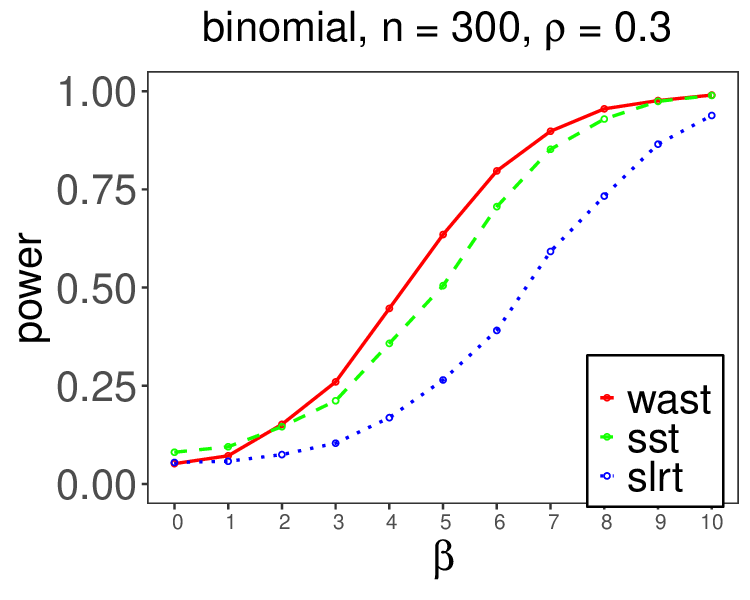}
		\includegraphics[scale=0.3]{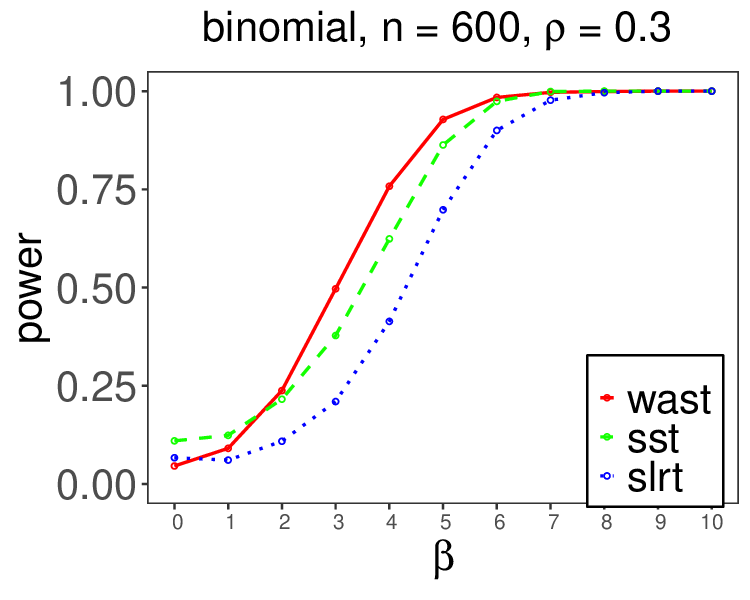}
		\includegraphics[scale=0.3]{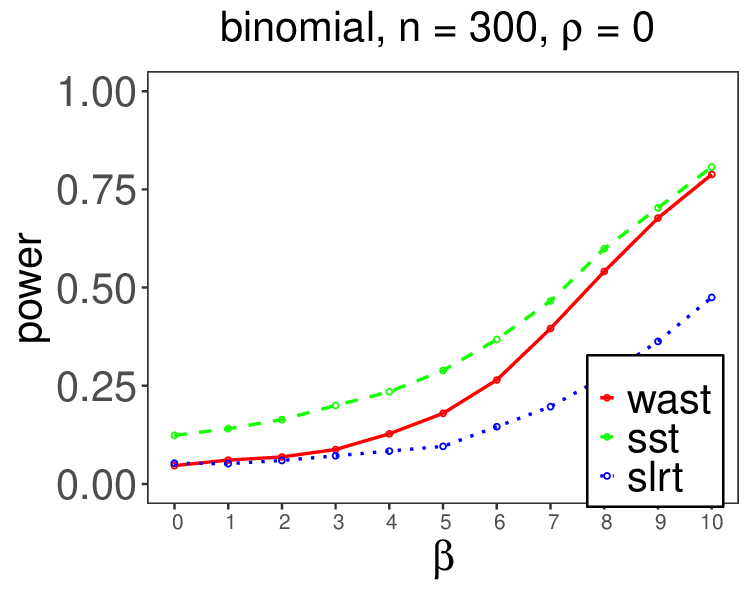}
		\includegraphics[scale=0.3]{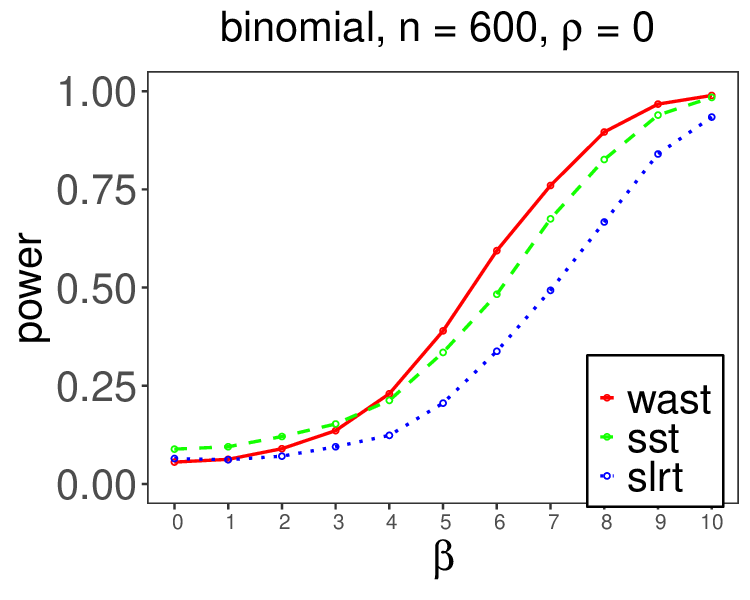}
		\includegraphics[scale=0.3]{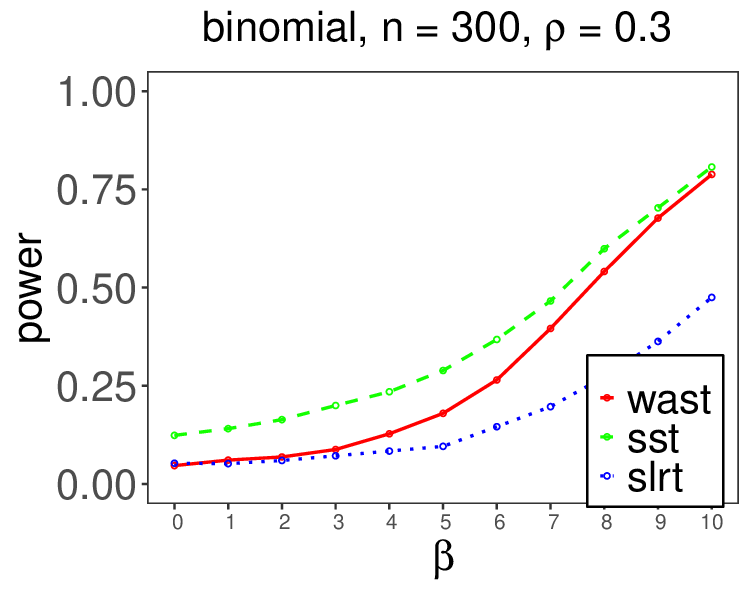}
		\includegraphics[scale=0.3]{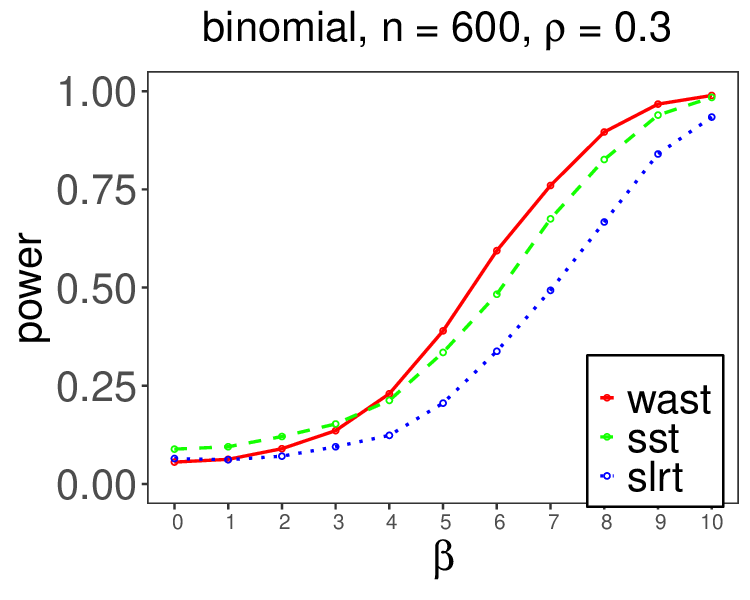}
		\includegraphics[scale=0.3]{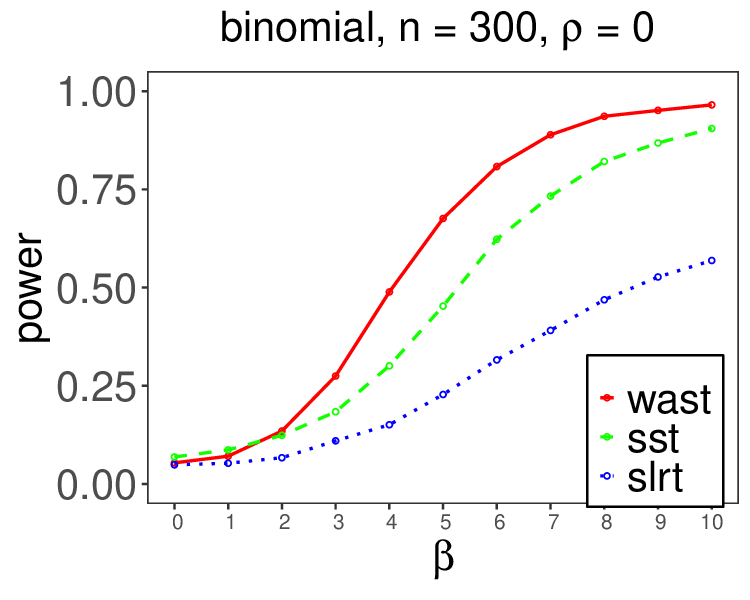}
		\includegraphics[scale=0.3]{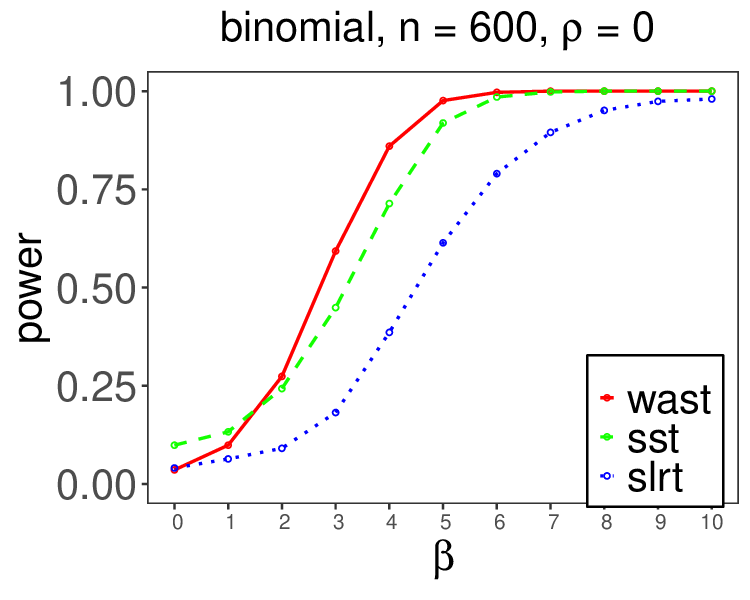}
		\includegraphics[scale=0.3]{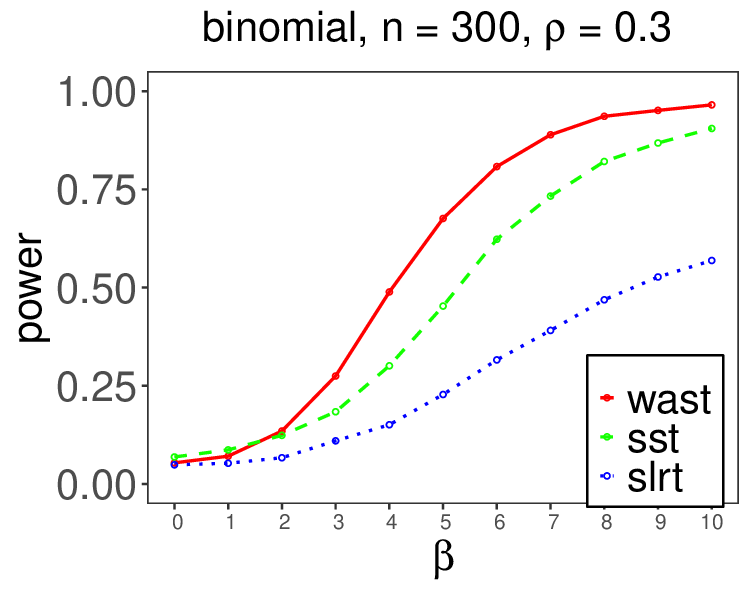}
		\includegraphics[scale=0.3]{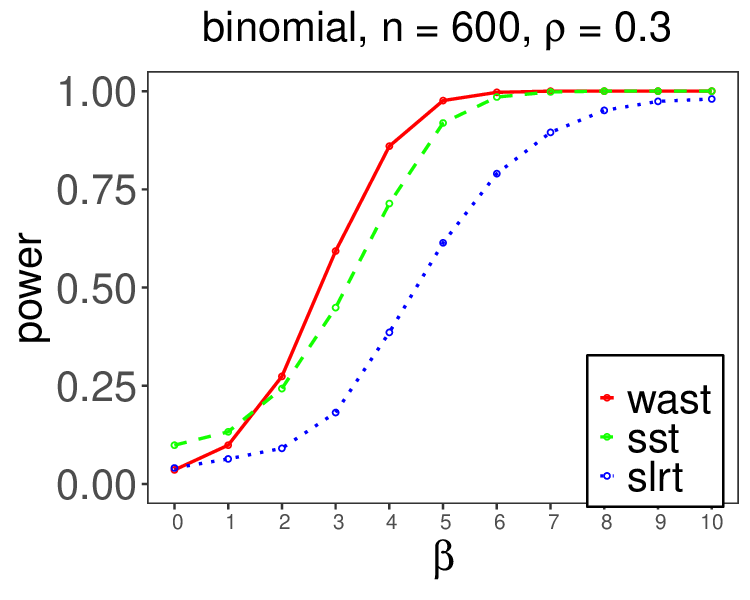} 		
		\includegraphics[scale=0.3]{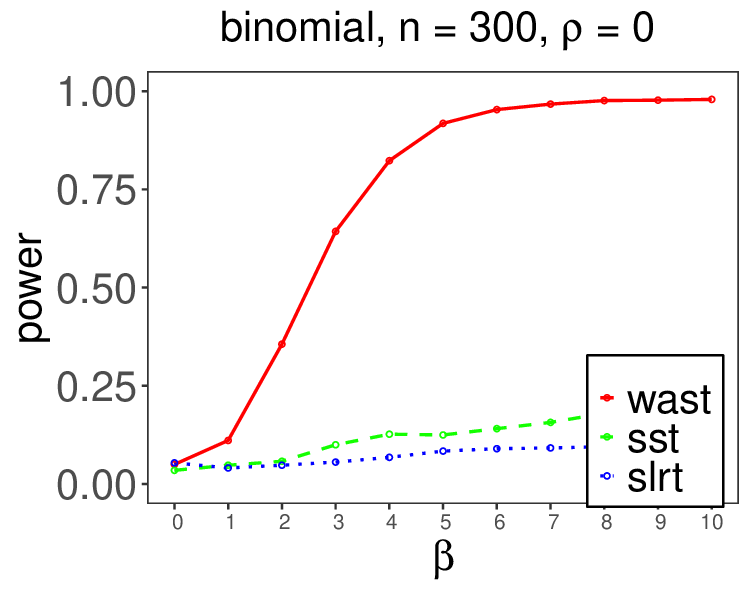}
		\includegraphics[scale=0.3]{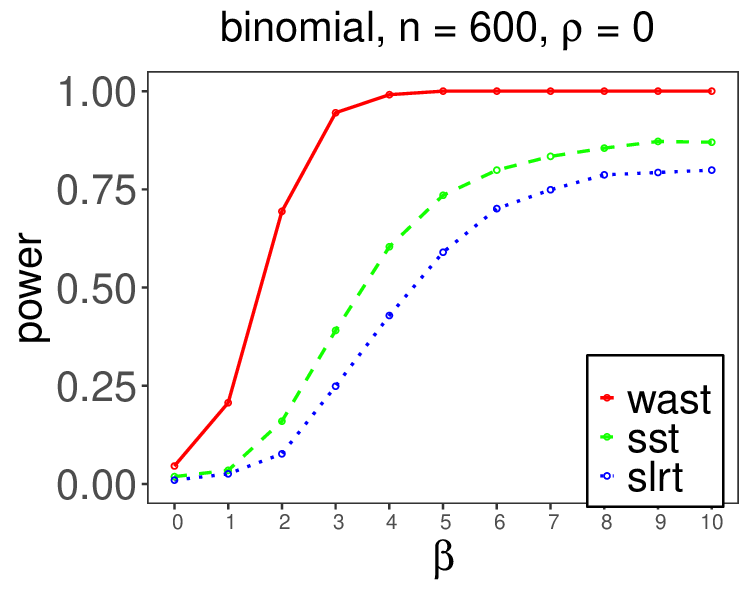}
		\includegraphics[scale=0.3]{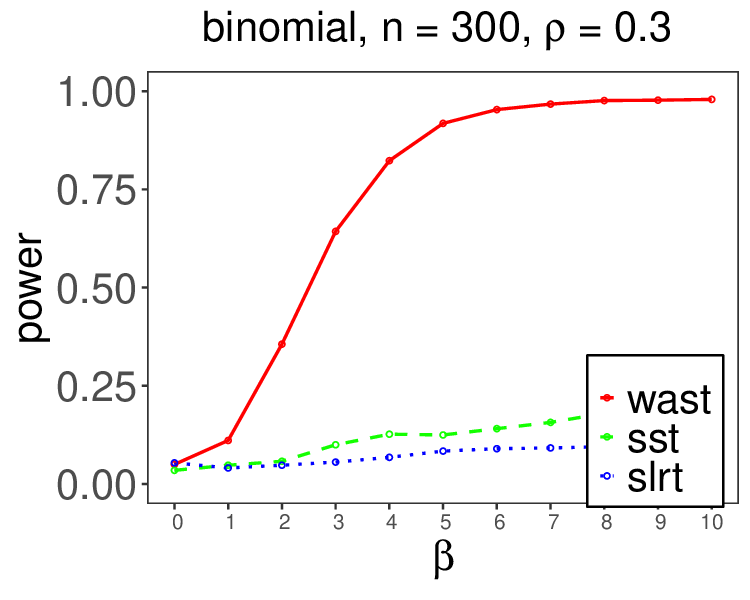}
		\includegraphics[scale=0.3]{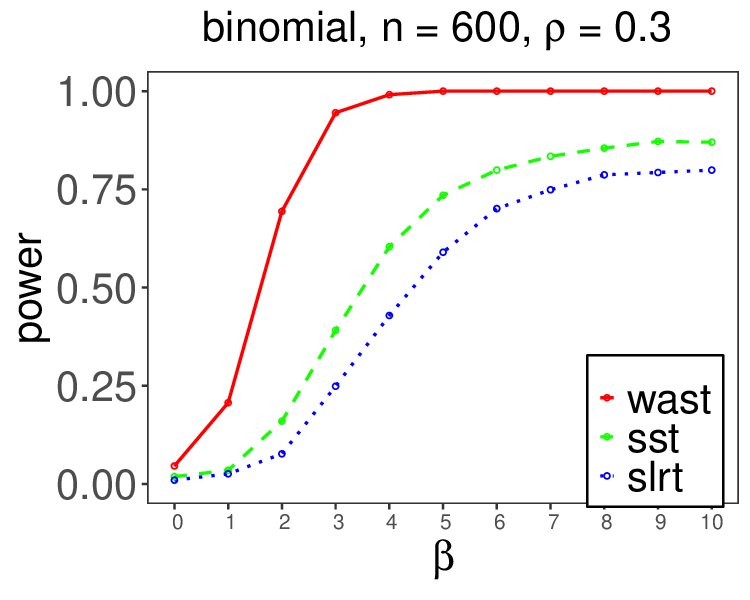}
		\includegraphics[scale=0.3]{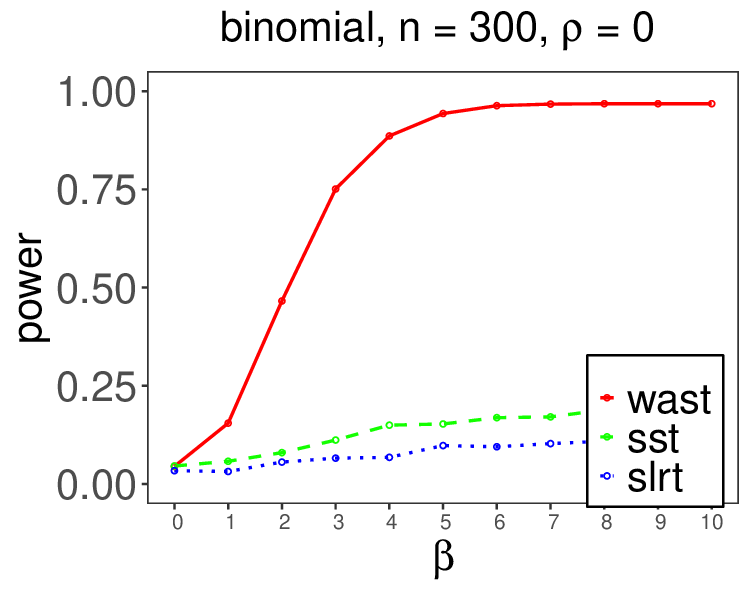}
		\includegraphics[scale=0.3]{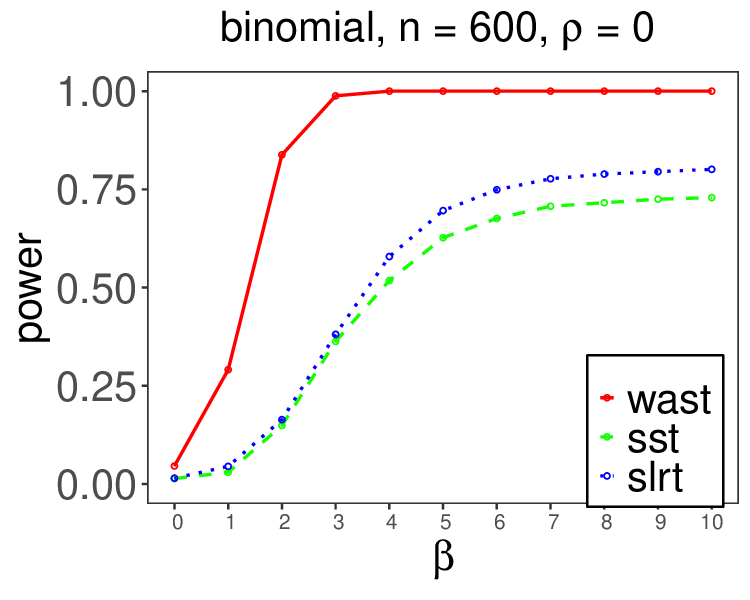}
		\includegraphics[scale=0.3]{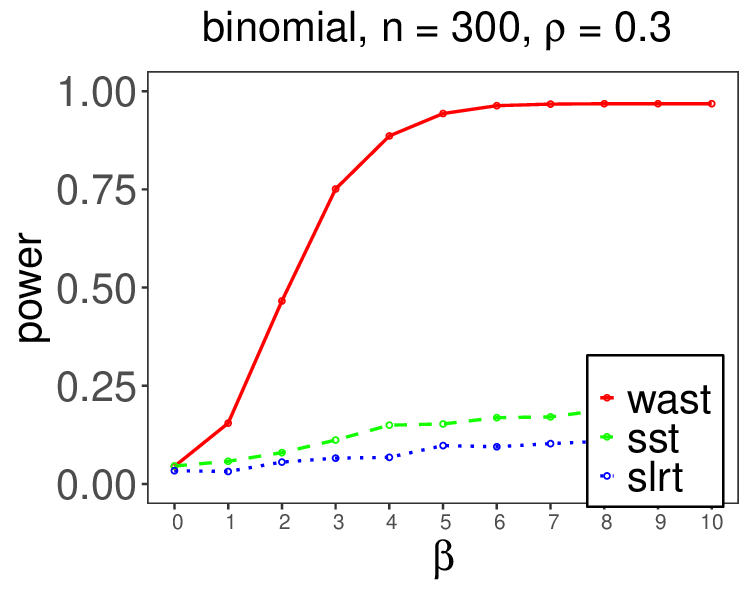}
		\includegraphics[scale=0.3]{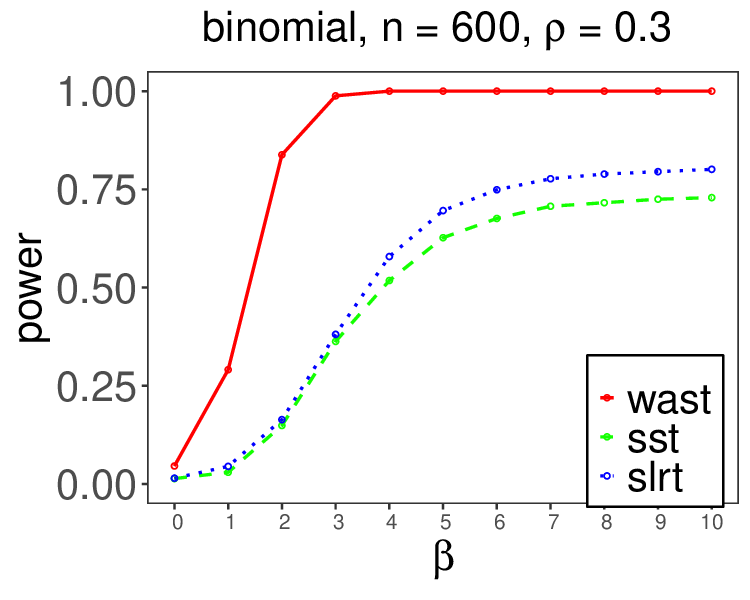}
		\caption{\it Powers of testing logistic regression by the proposed WAST (red solid line), SST (green dashed line), and SLRT (blue dotted line) for $n=(300,600)$. From top to bottom, each row panel depicts the powers for the case $(r,p,q)=(2,2,3)$, $(6,6,3)$, $(2,2,11)$, $(6,6,11)$, $(2,51,11)$, and $(6,51,11)$.
}
		\label{fig_binomial_zs4}
	\end{center}
\end{figure}

\begin{figure}[!ht]
	\begin{center}
		\includegraphics[scale=0.3]{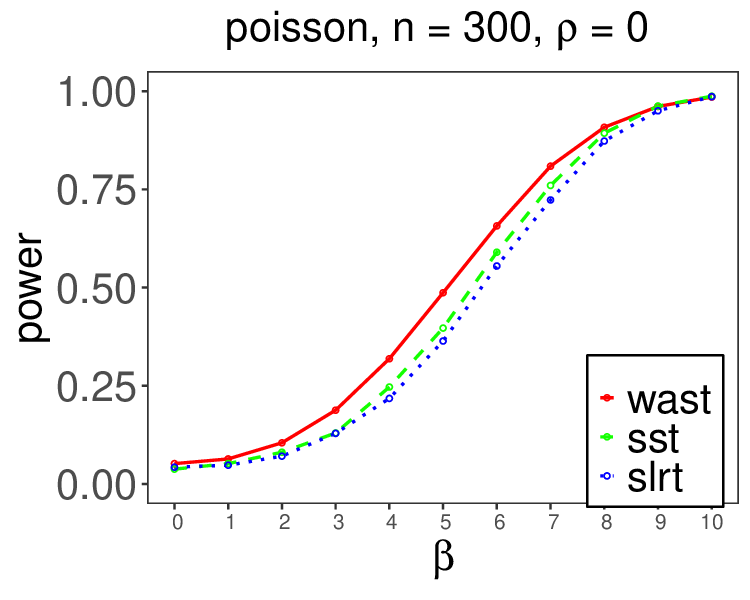}
		\includegraphics[scale=0.3]{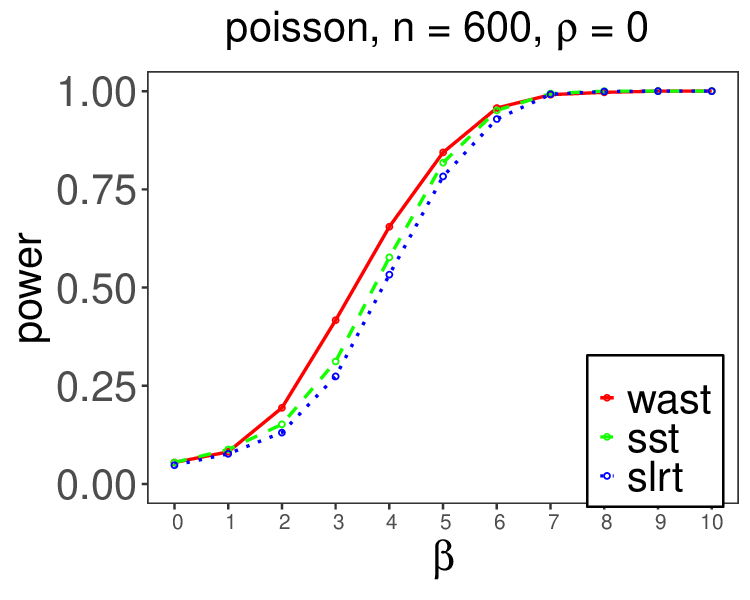}
		\includegraphics[scale=0.3]{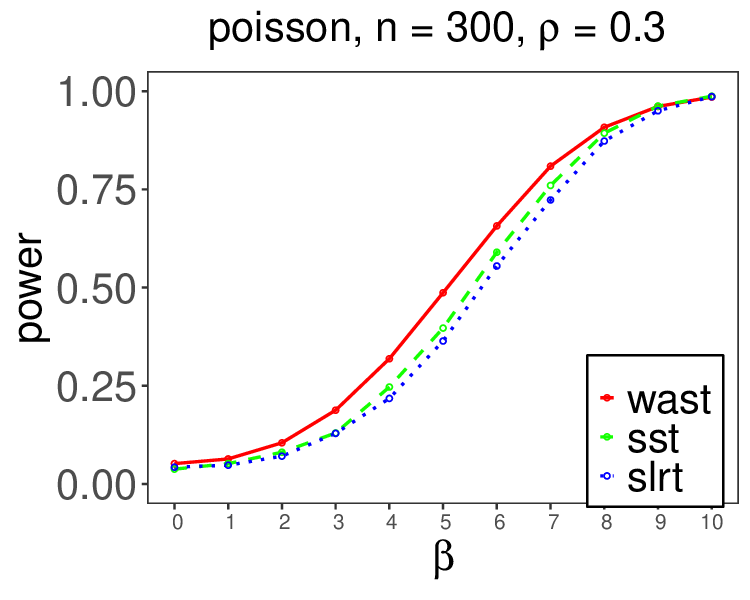}
		\includegraphics[scale=0.3]{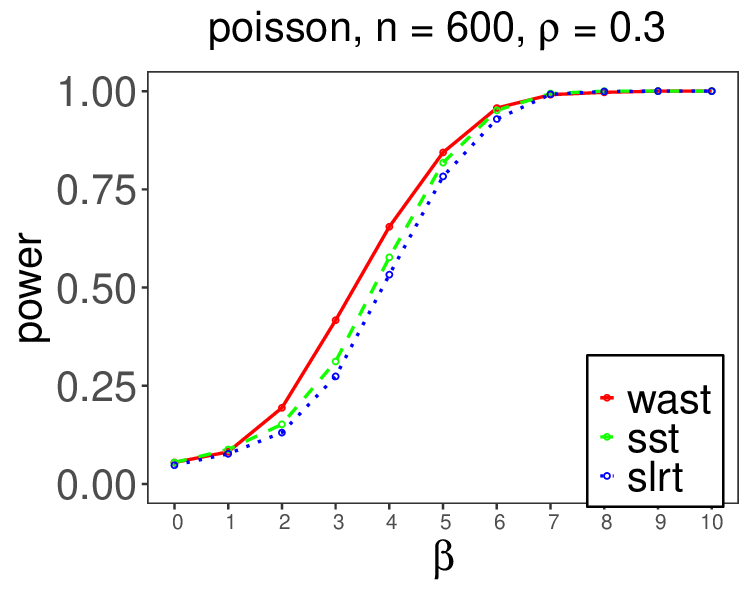}
		\includegraphics[scale=0.3]{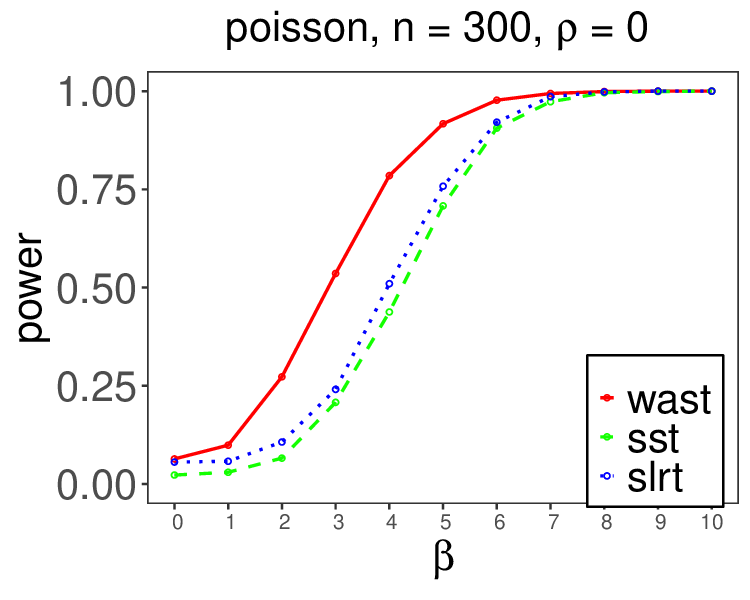}
		\includegraphics[scale=0.3]{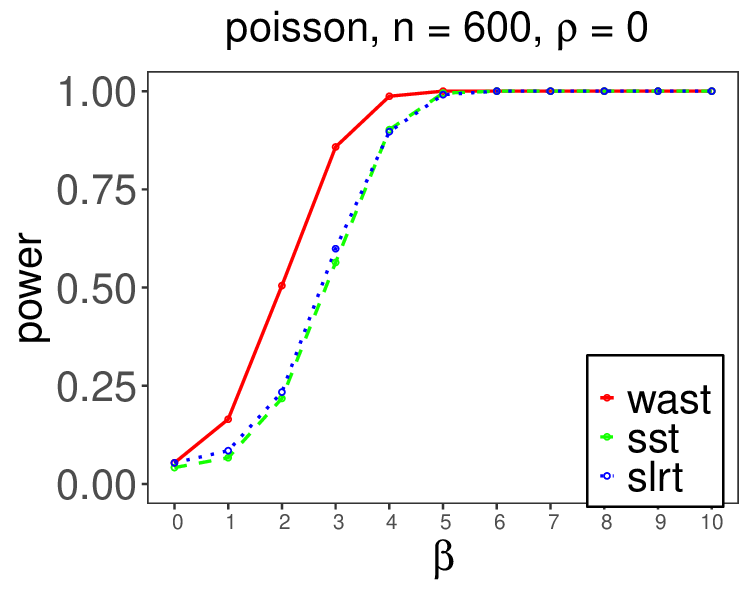}
		\includegraphics[scale=0.3]{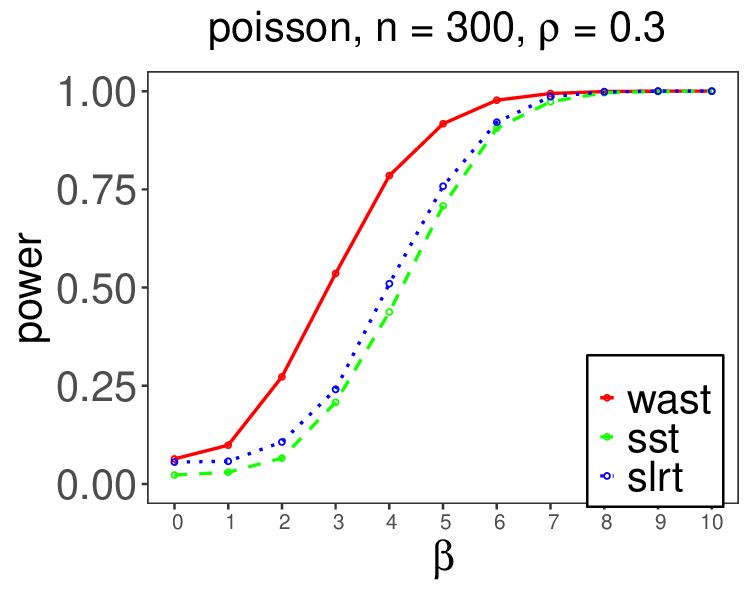}
		\includegraphics[scale=0.3]{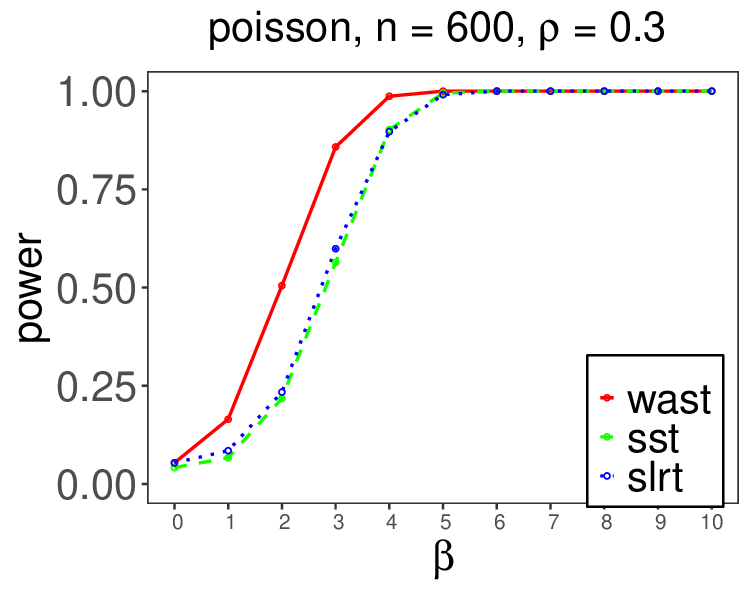}		
		\includegraphics[scale=0.3]{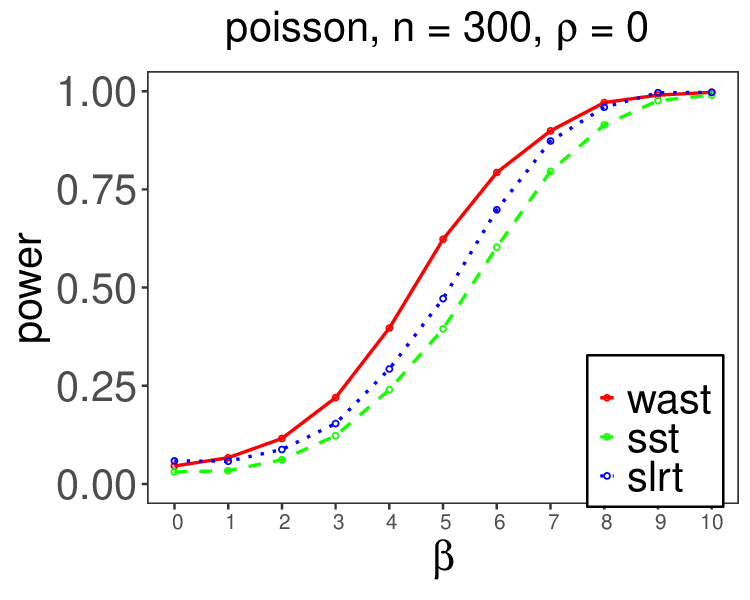}
		\includegraphics[scale=0.3]{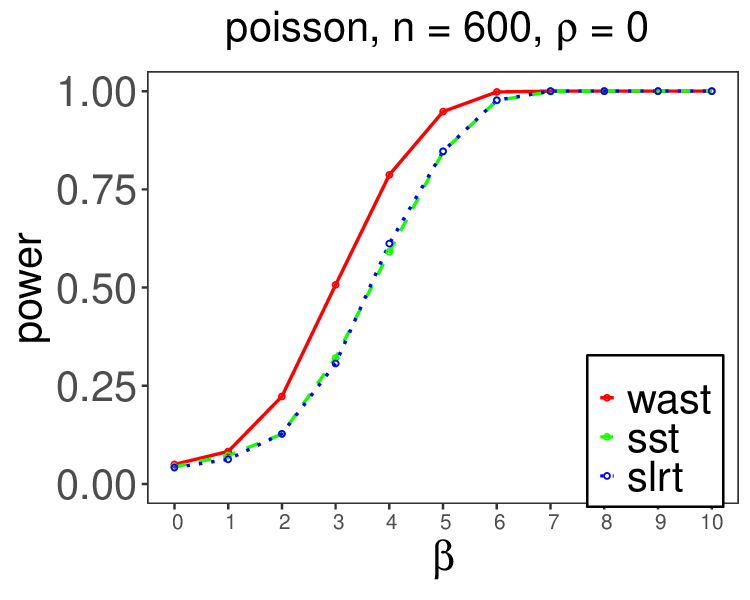}
		\includegraphics[scale=0.3]{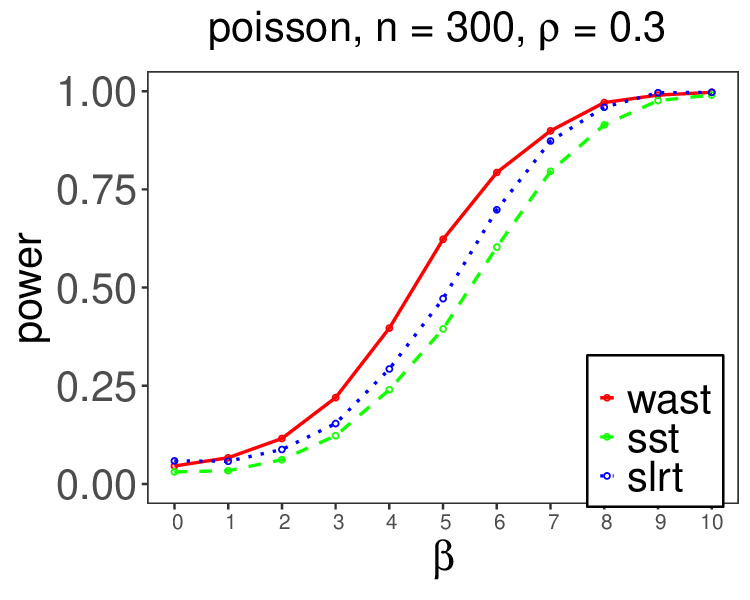}
		\includegraphics[scale=0.3]{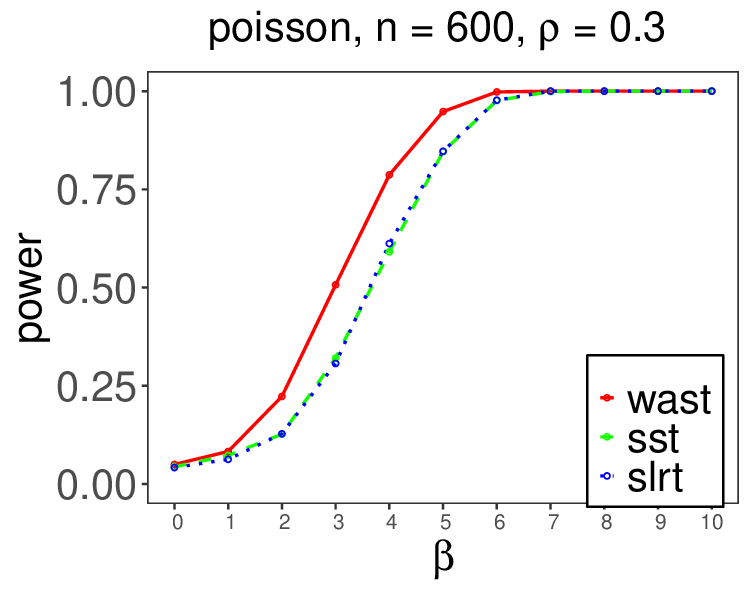}
		\includegraphics[scale=0.3]{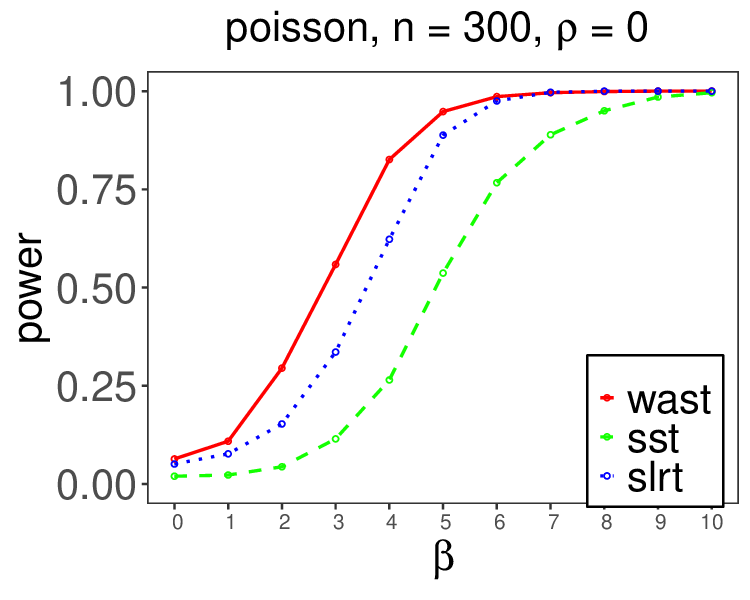}
		\includegraphics[scale=0.3]{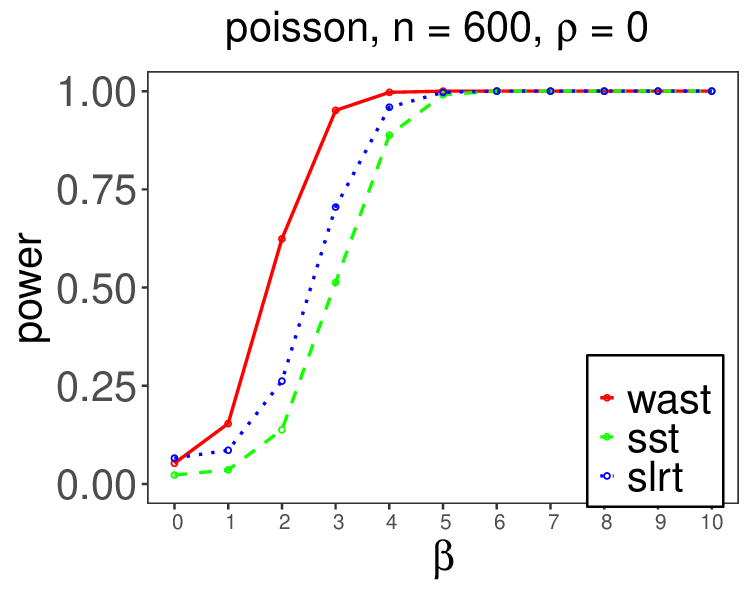}
		\includegraphics[scale=0.3]{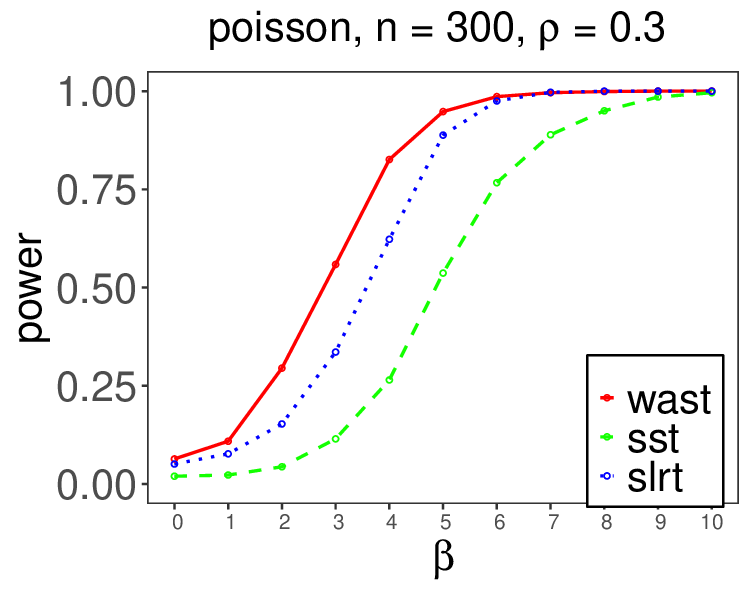}
		\includegraphics[scale=0.3]{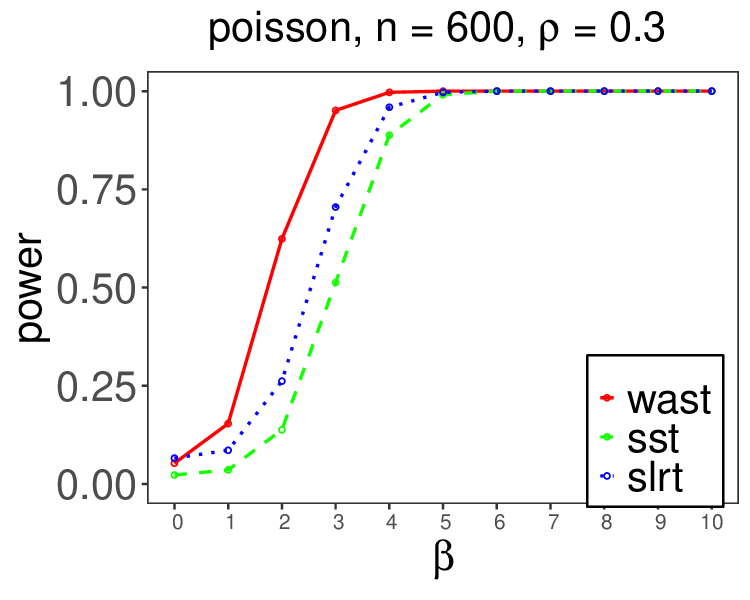} 		
		\includegraphics[scale=0.3]{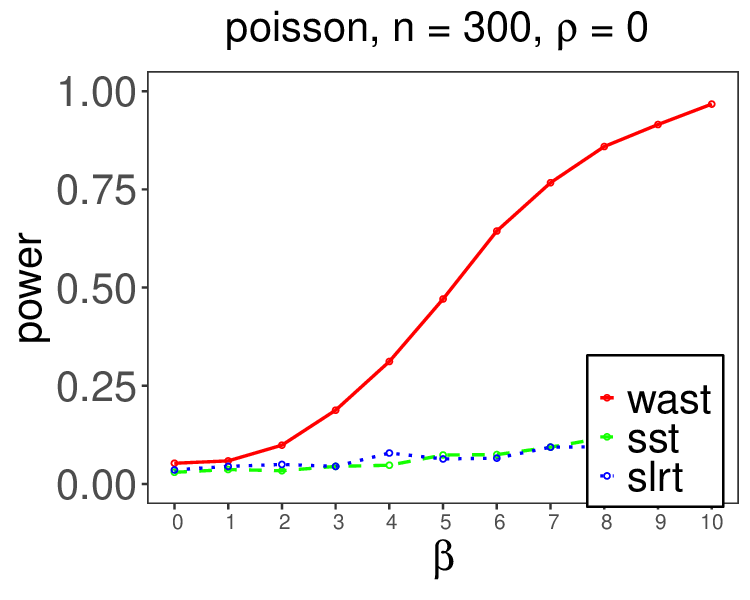}
		\includegraphics[scale=0.3]{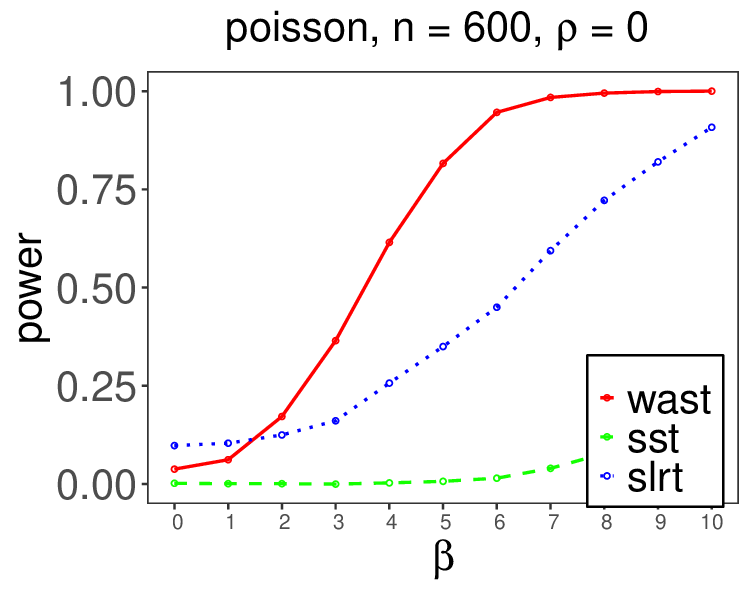}
		\includegraphics[scale=0.3]{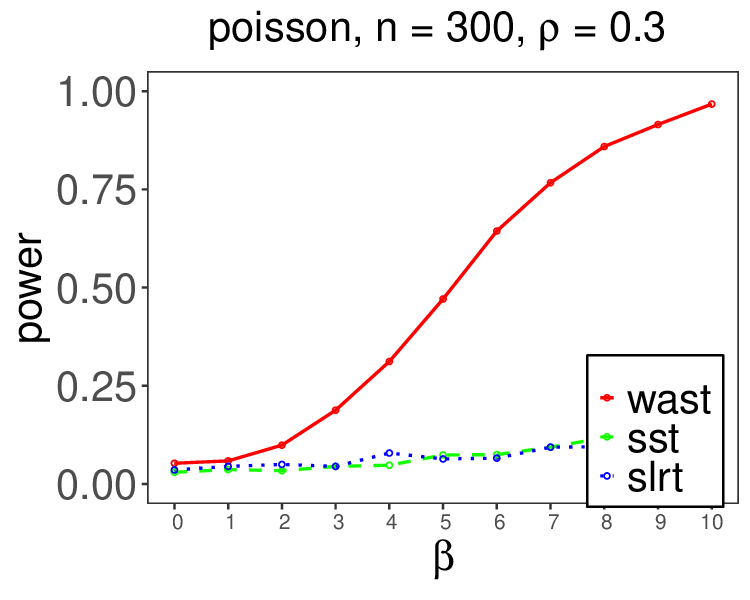}
		\includegraphics[scale=0.3]{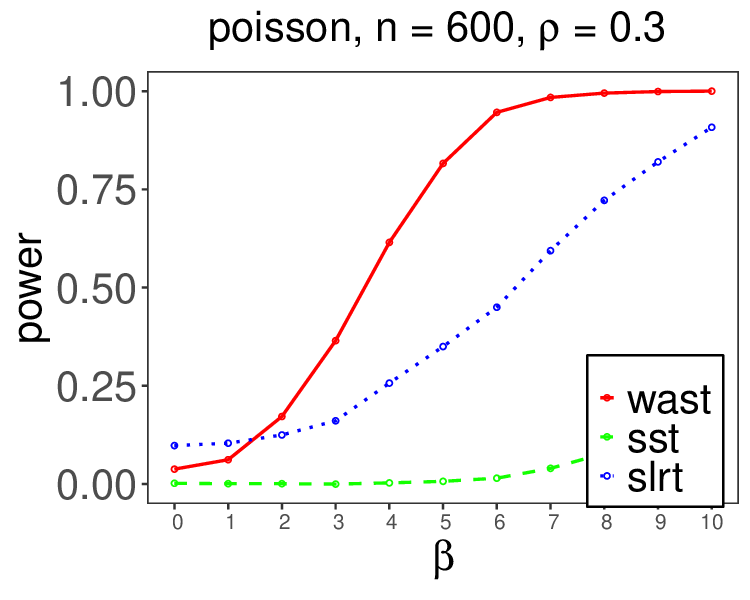}
		\includegraphics[scale=0.3]{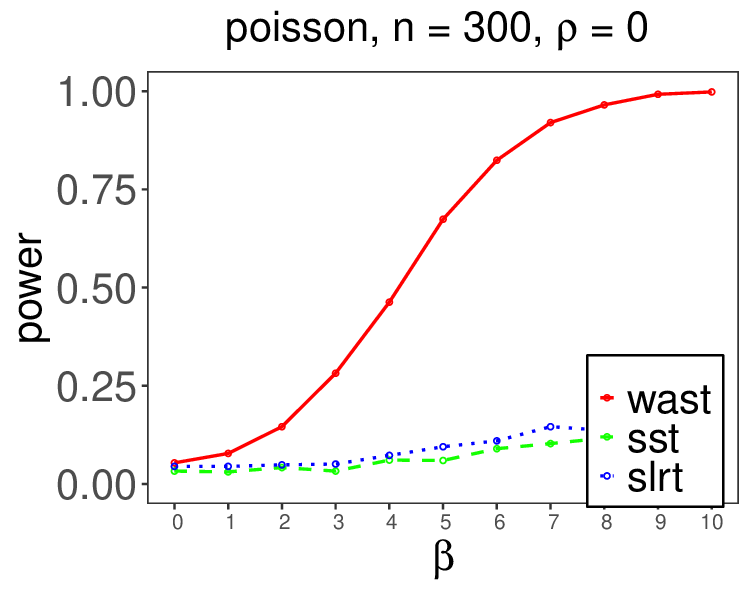}
		\includegraphics[scale=0.3]{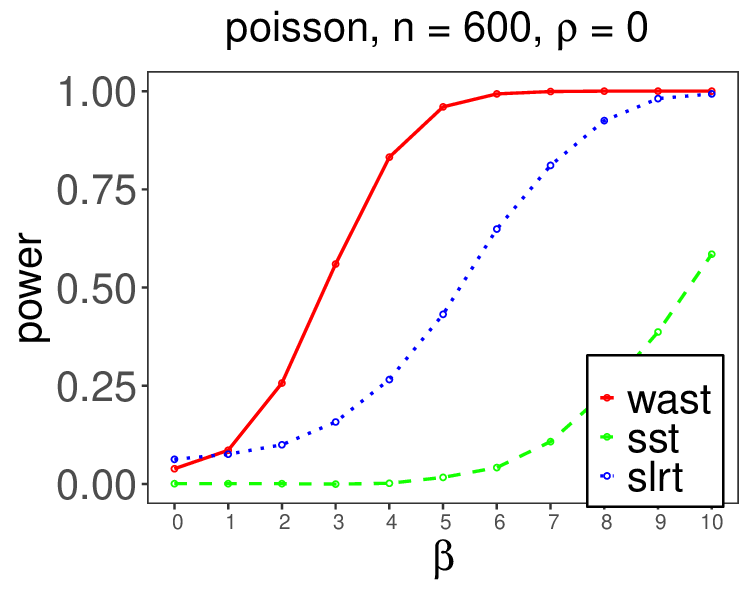}
		\includegraphics[scale=0.3]{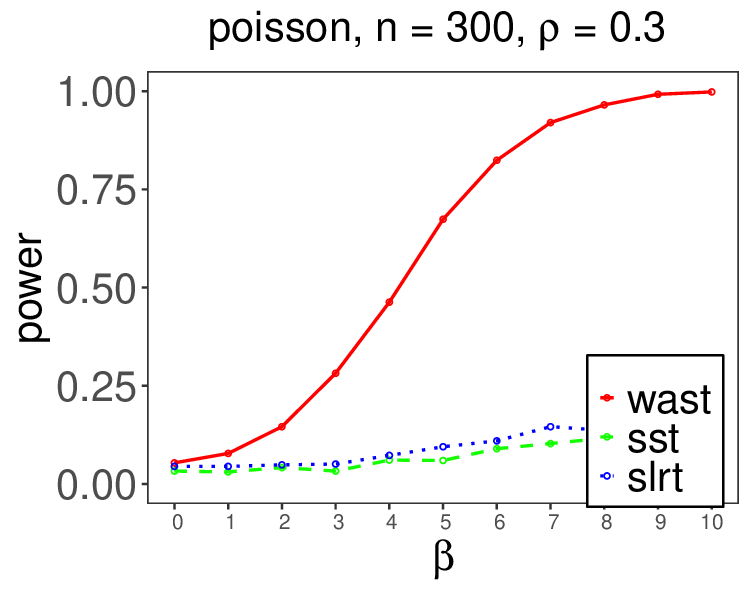}
		\includegraphics[scale=0.3]{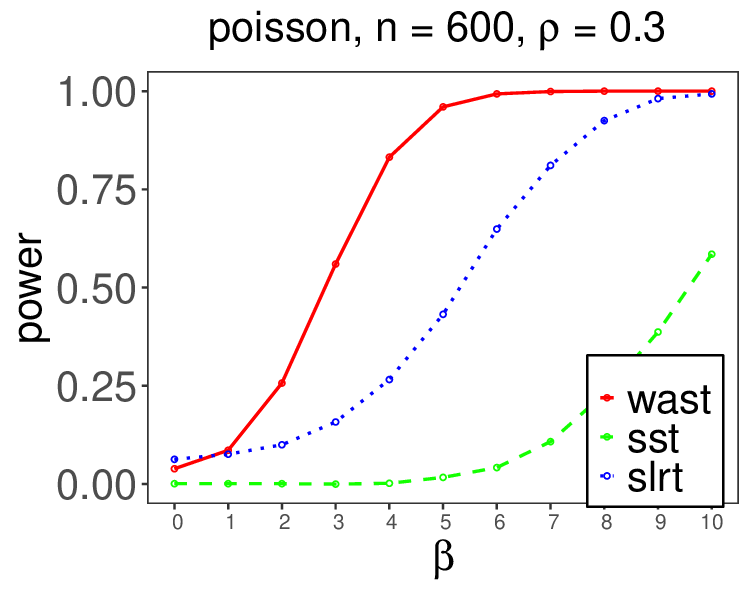}
		\caption{\it Powers of testing Poisson regression by the proposed WAST (red solid line), SST (green dashed line), and SLRT (blue dotted line) for $n=(300,600)$. From top to bottom, each row panel depicts the powers for the case $(r,p,q)=(2,2,3)$, $(6,6,3)$, $(2,2,11)$, $(6,6,11)$, $(2,51,11)$, and $(6,51,11)$.
}
	\label{fig_poisson_zs4}
	\end{center}
\end{figure}

We also generate $\varsigma_1, \cdots,\varsigma_d$ independently from normal distribution with mean zero and standard deviation $5$. Since the range of $\bZ\trans\btheta$ is large, the performance of the proposed test statistic are similar to that when $\varsigma_1, \cdots,\varsigma_d$ are independently from $t_3$ distribution with degrees of freedom 3. We list the sizes in Table \ref{table_size_s4} and depict power curves in Figures \ref{fig_gaussian_zs4}-\ref{fig_poisson_zs4} of the proposed test with two competitors SST and SLRT.
The same performance the proposed test statistic can be found from Table \ref{table_size_s4} and Figures \ref{fig_gaussian_zs4}-\ref{fig_poisson_zs4}.

\begin{table}
	\def~{\hphantom{0}}
	\caption{\label{table_size_t3} Type \uppercase\expandafter{\romannumeral1} errors of the proposed WAST, SST and SLRT.}
	\resizebox{\textwidth}{!}{
	\begin{threeparttable}
		\begin{tabular}{ll*{15}{c}}\\
			\hline
			\multirow{3}{*}{Family}&\multirow{3}{*}{$(r,p,q)$}
			&\multicolumn{7}{c}{$\rho=0$} && \multicolumn{7}{c}{$\rho=0.3$}\\
			\cline{3-9} \cline{11-17}
			& &\multicolumn{3}{c}{ $n=300$} && \multicolumn{3}{c}{ $n=600$} && \multicolumn{3}{c}{ $n=300$} && \multicolumn{3}{c}{ $n=600$}\\
			\cline{3-5} \cline{7-9} \cline{11-13} \cline{15-17}
			&
			&   WAST& SST & SLRT&& WAST& SST & SLRT&& WAST& SST & SLRT&& WAST& SST & SLRT\\
			\cline{3-17}
			Gaussian &$(2,2,3)$         & 0.043 & 0.041 & 0.054 && 0.035 & 0.036 & 0.044 && 0.045 & 0.038 & 0.062 && 0.045 & 0.036 & 0.046 \\
			&$(6,6,3)$                  & 0.048 & 0.029 & 0.057 && 0.053 & 0.031 & 0.051 && 0.049 & 0.021 & 0.043 && 0.051 & 0.029 & 0.050 \\
			& $(2,2,11)$                & 0.043 & 0.026 & 0.046 && 0.043 & 0.041 & 0.047 && 0.041 & 0.020 & 0.043 && 0.043 & 0.029 & 0.042 \\
			&$(6,6,11)$                 & 0.047 & 0.015 & 0.056 && 0.043 & 0.024 & 0.052 && 0.045 & 0.007 & 0.045 && 0.047 & 0.020 & 0.045 \\
			& $(2,51,11)$               & 0.050 & 0.027 & 0.033 && 0.048 & 0.002 & 0.044 && 0.064 & 0.039 & 0.040 && 0.040 & 0.002 & 0.040 \\
			&$(6,51,11)$                & 0.053 & 0.027 & 0.041 && 0.046 & 0.001 & 0.050 && 0.056 & 0.030 & 0.036 && 0.043 & 0.002 & 0.035 \\
			[1 ex]
			Binomial &$(2,2,3)$         & 0.050 & 0.084 & 0.050 && 0.043 & 0.069 & 0.055 && 0.044 & 0.083 & 0.040 && 0.040 & 0.071 & 0.060 \\
			&$(6,6,3)$                  & 0.042 & 0.075 & 0.045 && 0.056 & 0.102 & 0.051 && 0.036 & 0.064 & 0.047 && 0.054 & 0.111 & 0.052 \\
			& $(2,2,11)$                & 0.051 & 0.121 & 0.050 && 0.051 & 0.077 & 0.044 && 0.056 & 0.104 & 0.048 && 0.054 & 0.086 & 0.048 \\
			&$(6,6,11)$                 & 0.060 & 0.065 & 0.042 && 0.036 & 0.102 & 0.042 && 0.049 & 0.064 & 0.050 && 0.043 & 0.125 & 0.046 \\
			& $(2,51,11)$               & 0.056 & 0.036 & 0.034 && 0.049 & 0.020 & 0.018 && 0.042 & 0.042 & 0.043 && 0.045 & 0.013 & 0.009 \\
			&$(6,51,11)$                & 0.051 & 0.035 & 0.035 && 0.046 & 0.018 & 0.026 && 0.044 & 0.041 & 0.034 && 0.053 & 0.011 & 0.023 \\
			[1 ex]
			Poisson &$(2,2,3)$          & 0.050 & 0.052 & 0.059 && 0.035 & 0.051 & 0.045 && 0.048 & 0.046 & 0.049 && 0.039 & 0.042 & 0.048 \\
			&$(6,6,3)$                  & 0.042 & 0.024 & 0.055 && 0.047 & 0.032 & 0.056 && 0.036 & 0.016 & 0.054 && 0.056 & 0.028 & 0.059 \\
			& $(2,2,11)$                & 0.057 & 0.040 & 0.047 && 0.048 & 0.047 & 0.048 && 0.055 & 0.035 & 0.066 && 0.054 & 0.049 & 0.054 \\
			&$(6,6,11)$                 & 0.058 & 0.013 & 0.072 && 0.048 & 0.011 & 0.047 && 0.057 & 0.013 & 0.052 && 0.052 & 0.030 & 0.059 \\
			& $(2,51,11)$               & 0.045 & 0.025 & 0.048 && 0.039 & 0.003 & 0.116 && 0.044 & 0.039 & 0.047 && 0.059 & 0.000 & 0.113 \\
			&$(6,51,11)$                & 0.062 & 0.032 & 0.042 && 0.045 & 0.003 & 0.077 && 0.048 & 0.028 & 0.036 && 0.054 & 0.001 & 0.086 \\
			\hline
		\end{tabular}
		\begin{tablenotes}
		\item The setting for single change plane is from Section \ref{scp1}  with $\Gv$ $\bZ$ generated from $t_3$ distribution. The nominal significant level is 0.05.	
		\end{tablenotes}
	\end{threeparttable}
	}
\end{table}

\begin{table}
	\def~{\hphantom{0}}
	\caption{\label{table_size_s4} Type \uppercase\expandafter{\romannumeral1} errors of the proposed WAST, SST and SLRT.}
	\resizebox{\textwidth}{!}{
	\begin{threeparttable}
		\begin{tabular}{ll*{15}{c}}\\
			\hline
			\multirow{3}{*}{Family}&\multirow{3}{*}{$(r,p,q)$}
			&\multicolumn{7}{c}{ $\rho=0$} && \multicolumn{7}{c}{$\rho=0.3$}\\
			\cline{3-9} \cline{11-17}
			& &\multicolumn{3}{c}{ $n=300$} && \multicolumn{3}{c}{ $n=600$} && \multicolumn{3}{c}{ $n=300$} && \multicolumn{3}{c}{ $n=600$}\\
			\cline{3-5} \cline{7-9} \cline{11-13} \cline{15-17}
			&
			&   WAST& SST & SLRT&& WAST& SST & SLRT&& WAST& SST & SLRT&& WAST& SST & SLRT\\
			\cline{3-17}
			Gaussian &$(2,2,3)$         & 0.049 & 0.033 & 0.048 && 0.049 & 0.036 & 0.049 && 0.045 & 0.032 & 0.049 && 0.053 & 0.044 & 0.046 \\
			&$(6,6,3)$                  & 0.053 & 0.013 & 0.040 && 0.050 & 0.030 & 0.042 && 0.062 & 0.016 & 0.042 && 0.051 & 0.028 & 0.043 \\
			& $(2,2,11)$                & 0.045 & 0.030 & 0.058 && 0.049 & 0.041 & 0.059 && 0.048 & 0.020 & 0.051 && 0.046 & 0.031 & 0.047 \\
			&$(6,6,11)$                 & 0.050 & 0.016 & 0.048 && 0.053 & 0.026 & 0.044 && 0.050 & 0.008 & 0.048 && 0.047 & 0.023 & 0.048 \\
			& $(2,51,11)$               & 0.045 & 0.032 & 0.044 && 0.049 & 0.000 & 0.059 && 0.046 & 0.034 & 0.050 && 0.046 & 0.000 & 0.045 \\
			&$(6,51,11)$                & 0.042 & 0.026 & 0.039 && 0.057 & 0.000 & 0.050 && 0.057 & 0.023 & 0.043 && 0.048 & 0.000 & 0.046 \\
			[1 ex]
			Binomial &$(2,2,3)$         & 0.055 & 0.109 & 0.052 && 0.050 & 0.082 & 0.054 && 0.055 & 0.100 & 0.048 && 0.050 & 0.072 & 0.059 \\
			&$(6,6,3)$                  & 0.052 & 0.081 & 0.055 && 0.046 & 0.110 & 0.067 && 0.056 & 0.076 & 0.048 && 0.051 & 0.117 & 0.059 \\
			& $(2,2,11)$                & 0.047 & 0.124 & 0.053 && 0.056 & 0.089 & 0.065 && 0.051 & 0.106 & 0.048 && 0.052 & 0.099 & 0.057 \\
			&$(6,6,11)$                 & 0.054 & 0.069 & 0.049 && 0.036 & 0.099 & 0.041 && 0.044 & 0.057 & 0.068 && 0.046 & 0.129 & 0.044 \\
			& $(2,51,11)$               & 0.050 & 0.035 & 0.054 && 0.046 & 0.019 & 0.010 && 0.041 & 0.046 & 0.044 && 0.054 & 0.011 & 0.008 \\
			&$(6,51,11)$                & 0.045 & 0.046 & 0.034 && 0.046 & 0.014 & 0.015 && 0.052 & 0.036 & 0.044 && 0.055 & 0.013 & 0.019 \\
			[1 ex]
			Poisson &$(2,2,3)$          & 0.052 & 0.038 & 0.043 && 0.055 & 0.055 & 0.048 && 0.044 & 0.039 & 0.050 && 0.052 & 0.046 & 0.055 \\
			&$(6,6,3)$                  & 0.064 & 0.023 & 0.056 && 0.054 & 0.042 & 0.054 && 0.065 & 0.014 & 0.052 && 0.051 & 0.029 & 0.048 \\
			& $(2,2,11)$                & 0.046 & 0.031 & 0.059 && 0.050 & 0.044 & 0.042 && 0.051 & 0.039 & 0.057 && 0.050 & 0.041 & 0.038 \\
			&$(6,6,11)$                 & 0.064 & 0.020 & 0.051 && 0.053 & 0.023 & 0.066 && 0.051 & 0.014 & 0.052 && 0.053 & 0.022 & 0.056 \\
			& $(2,51,11)$               & 0.053 & 0.030 & 0.036 && 0.038 & 0.002 & 0.098 && 0.051 & 0.036 & 0.033 && 0.043 & 0.000 & 0.093 \\
			&$(6,51,11)$                & 0.054 & 0.033 & 0.045 && 0.039 & 0.001 & 0.063 && 0.046 & 0.021 & 0.032 && 0.045 & 0.000 & 0.060 \\
			\hline
		\end{tabular}
		\begin{tablenotes}
		\item The setting for single change plane is from Section \ref{scp1} with $\Gv$ $\bZ$ generated from normal distribution with mean 0 and standard deviation 5. The nominal significant level is 0.05.	
		\end{tablenotes}
	\end{threeparttable}
	}
\end{table}

\subsection{Multiple change-plane analysis}\label{scp2}

It is straightforward to extend the proposed test procedure for models with a change plane to the case of multiple change planes, as considered by \citet{2021Multithreshold}. For instance, we consider the quantile regression model with multiple change planes:
\begin{align*}
	Q_{Y_i}(\tau|\tbX_i, \bX_i, \bZ_i) = \tbX_i\trans\balpha(\tau)+\sum_{t=1}^{T}\bX_{ti}\trans\bbeta_t(\tau)\bone(\bZ_{ti}\trans\btheta_t(\tau)\geq0),
\end{align*}
where $\balpha(\tau)=(\alpha_1(\tau), \cdots, \alpha_r(\tau))\trans\in\Theta_{\alpha}\subseteq\mathbb{R}^{r}$, $\bbeta_t(\tau)=(\beta_{t1}(\tau), \cdots, \beta_{tp_t}(\tau))\trans\in\Theta_{\beta_t}\subseteq\mathbb{R}^{p_t}$, and $\btheta_t(\tau)=(\theta_{t1}(\tau), \cdots, \theta_{tq_t}(\tau))\trans\in\Theta_{\btheta_t}\subseteq\mathbb{R}^{q_t}$ are unknown parameters, and the number of change planes $T$ is a preset integer, $\bX_i=(\bX_{1i}\trans, \cdots, \bX_{Ti}\trans)\trans$ and $\bZ_i=(\bZ_{1i}\trans, \cdots, \bZ_{Ti}\trans)\trans$. We still denote $\bbeta(\tau)=(\bbeta_1(\tau)\trans, \cdots, \bbeta_T(\tau)\trans)\trans$ and $\btheta(\tau)=(\btheta_1(\tau)\trans, \cdots, \btheta_T(\tau)\trans)\trans$,
and $\Theta_{\beta}=\Theta_{\beta_1}\times\cdots\times\Theta_{\beta_T}$ and $\Theta_{\theta}=\Theta_{\theta_1}\times\cdots\times\Theta_{\theta_T}$.

Testing whether $\bbeta_t(\tau)=\bzero$ for all $t\in\{1,\cdots,T\}$ is of interest, i.e.,
\begin{align}\label{test_multi}
	H_0: \bbeta(\tau)=\bzero \quad \mathrm{versus} \quad H_1: \bbeta(\tau)\neq \bzero.
\end{align}
Similar to the WAST in the main paper, we propose the integrated statistic
\begin{align*}
	T_n^{(m)} = \frac{1}{n(n-1)}\sum_{i\neq j}\tilde{\omega}_{ij}\left[\bone(Y_i-\tbX_i\trans\hbalpha(\tau)\leq0)-\tau\right]\left[\bone(Y_j-\tbX_j\trans\hbalpha(\tau)\leq0)-\tau\right]
\end{align*}
with
\begin{align}\label{rho_m}
	\tilde{\omega}_{ij} = \sum_{t=1}^{T}\bX_{ti}\trans\bX_{tj}\int_{\btheta_t\in\Theta_{\theta_t}}\bone(\bZ_{ti}\trans\btheta_t\geq0) \bone(\bZ_{tj}\trans\btheta_t\geq0)w_t(\btheta_t)d\btheta_t.
\end{align}
Taking the weight $w(\btheta)$ to be the standard multivariate normal density leads to, for $i\neq j$,
\begin{align*} \tilde{\omega}_{ij}=\frac{1}{4}\sum_{t=1}^{T}\bX_{ti}\trans\bX_{tj}+\frac{1}{2\pi}\sum_{t=1}^{T}\bX_{ti}\trans\bX_{tj}\arctan\left(\frac{\varrho_{t, ij}}{\sqrt{1-\varrho_{t, ij}^2}}\right),
\end{align*}
where $\varrho_{t, ij} = \bZ_{ti}\trans\bZ_{tj}(\|\bZ_{ti}\|\|\bZ_{tj}\|)^{-1}$.

As in model (7) in the main paper, we consider the empirical estimating functions
\begin{align}\label{model1b}
	\begin{split}
		\Psi_{n}(\balpha, \bbeta, \btheta)=&\sum_{i=1}^{n} \bW_i(\btheta)\psi_0^{(m)}(\bV_i, \balpha, \bbeta, \btheta), \\
		\Psi_{1n}(\balpha)=&\sum_{i=1}^{n} \psi_1(\bV_i, \balpha),
	\end{split}
\end{align}
where $\bW_i(\btheta)=(\bX_{1i}\trans\bone(\bZ_{1i}\trans\btheta_1\geq 0), \cdots, \bX_{Ti}\trans\bone(\bZ_{Ti}\trans\btheta_T\geq 0))\trans$. For ease of expression, we denote $\psi^{(m)}(\bV, \balpha, \bbeta, \btheta)=\bW(\btheta)\psi_0^{(m)}(\bV, \balpha, \bbeta, \btheta)$, and $\psi_0^{(m)}(\bV_i, \balpha)=\psi_0^{(m)}(\bV_i, \balpha, \bzero, \btheta)$, where $\psi_0^{(m)}(\bV_i, \balpha)$ is free of $\btheta$, and $\psi^{(m)}(\bV, \balpha, \bbeta, \btheta)$ satisfies $\bE\psi^{(m)}(\bV, \balpha, \bbeta, \btheta)=0$ with true parameters $\balpha_0$, $\bbeta_0$, and $\btheta_0$. We propose the integrated statistic with $\tilde{\omega}_{ij}$ defined in \eqref{rho_m}, i.e.,
\begin{align*}
	T_n^{(m)} = \frac{1}{n(n-1)}\sum_{i\neq j}\tilde{\omega}_{ij}\psi_0^{(m)}(\bV_i, \hbalpha)\psi_0^{(m)}(\bV_j, \hbalpha).
\end{align*}

Let $\tK^t(\btheta)=\partial\bE[\bX_{ti}\bone(\bZ_{ti}\trans\btheta_1\geq0)\psi_0^{(m)}(\bV, \balpha_0)]/\partial\balpha\trans$, 
$\tK(\btheta)=(\tK^1(\btheta)\trans, \cdots, \tK^T(\btheta)\trans)\trans$,
and the kernel of a U-statistic under the null hypothesis be
\begin{align}\label{kernel0c}
	\begin{split}
		h^{(m)}(\bV_i, \bV_j)
		=& \tilde{\omega}_{ij}\psi_0^{(m)}(\bV_i, \balpha_0)\psi_0^{(m)}(\bV_j, \balpha_0)- \psi_1(\bV_i, \balpha_0)\trans \tK_{j}\\
		&- \tK_{i}\trans \psi_1(\bV_j, \balpha_0) + \psi_1(\bV_i, \balpha_0)\trans \tH \psi_1(\bV_j, \balpha_0),
	\end{split}
\end{align}
where
\begin{align*}
	\tK_i
	=\sum_{t=1}^{T}\int_{\btheta_t\in\Theta_{\theta_t}}
	J\trans \tK^t(\btheta)\trans
	\bX_{ti}\bone(\bZ_{ti}\trans\btheta_t\geq0)\psi_0^{(m)}(\bV_i, \balpha_0) w_t(\btheta_t)d\btheta_t
\end{align*}
and
\begin{align*}
	\tH
	=\sum_{t=1}^{T}\int_{\btheta_t\in\Theta_{\theta_t}}
	J\trans\tK^t(\btheta)\trans
	\tK^t(\btheta)J w_t(\btheta_t)d\btheta_t.
\end{align*}

We are now ready to derive the asymptotic distributions of $T_n^{(m)}$ under the null and local alternative hypotheses.
\begin{thm}\label{coroll31}
	If Assumptions~{(A1)--(A5)} hold, then under the null hypothesis, we have
	\begin{align*}
		nT_n^{(m)} -\mu_0 \lkonv \zeta,
	\end{align*}
	where $\mu_0=-2\bE[\psi_1(\bV_1, \balpha_0)\trans \tK_{1}]
	+\bE[\psi_1(\bV_1, \balpha_0)\trans \tH \psi_1(\bV_1, \balpha_0)]$, $\zeta$ is a random variable of the form $\zeta=\sum_{j=1}^{\infty}\lambda_{j}(\chi^2_{1j}-1)$, and $\chi^2_{11}, \chi^2_{12}, ...$ are independent $\chi^2_{1}$ variables.
	Here, $\{\lambda_{j}\}$ are the eigenvalues of the kernel $h^{(m)}(\bv_{1}, \bv_{2})$ under $f(\bv, \balpha_0, \bzero, \btheta_0)$, i.e., they are the solutions of $\lambda_{j}g_{j}(\bv_2)=\int_{0}^{\infty}h^{(m)}(\bv_1, \bv_2)g_{j}(\bv_1)f(\bv_1, \balpha_0, \bzero, \btheta_0)d\bv_1$ for nonzero $g_{j}$.
\end{thm}

\noindent{\bf Proof.} 
	We show this Theorem as the same as the proof of Theorem 1.
	Define notations as follows.
	\begin{align*}
		\tK =& \left(
		\begin{array}{c}
			\int_{\btheta_1\in\Theta_{\theta_1}} \rP[\bX_{1i}\trans\bone(\bZ_{1i}\trans\btheta_1\geq0)\partial\psi_0^{(m)}(\bV_i, \balpha_0)/\partial\balpha\trans]w(\btheta_1) d\btheta_1\\
			\vdots\\
			\int_{\btheta_T\in\Theta_{\theta_T}} \rP[\bX_{Ti}\trans\bone(\bZ_{Ti}\trans\btheta_T\geq0)\partial\psi_0^{(m)}(\bV_i, \balpha_0)/\partial\balpha\trans]w(\btheta_T) d\btheta_T
		\end{array}
		\right),\\
		\tpsi_{i} =& \left(
		\begin{array}{c}
			\int_{\btheta_1\in\Theta_{\theta_1}} \bX_{1i}\trans\bone(\bZ_{1i}\trans\btheta_1\geq0)w(\btheta_1) d\btheta_1\\
			\vdots\\
			\int_{\btheta_T\in\Theta_{\theta_T}} \bX_{Ti}\trans\bone(\bZ_{Ti}\trans\btheta_T\geq0)w(\btheta_T) d\btheta_T
		\end{array}
		\right),\\
		\tpsi_{ij}(\balpha) =& \tilde{\omega}_{ij}\psi_0^{(m)}(\bV_i, \balpha)\psi_0^{(m)}(\bV_j, \balpha).
	\end{align*}
	By Assumption (A1), we have
	\begin{align*}
		\Pn[\bW(\btheta)\psi_0^{(m)}(\bV, \hbalpha)]
		=&\Pn[\bW(\btheta)\psi_0^{(m)}(\bV, \hbalpha_0)]	- \tK(\btheta)J\Psi_{1n}(\balpha_0)+o_p(n^{-1/2}),
	\end{align*}
	and
	\begin{align*}
		\frac{1}{n}\sum_{i\neq j}\tpsi_{ij}(\hbalpha)
		=&\frac{1}{n}\sum_{i\neq j}\tilde{\omega}_{ij}\psi_0^{(m)}(\bV_i, \hbalpha)\psi_0^{(m)}(\bV_j, \hbalpha)\\
		=&\frac{1}{n}\sum_{i\neq j}\left[\tpsi_{ij}(\balpha_0) - \frac{1}{n}\Psi_{1n}(\balpha_0)\trans \tK_{i}- \frac{1}{n}\tK_{j}\trans \Psi_{1n}(\balpha_0) + \frac{1}{n^2}\Psi_{1n}(\balpha_0)\trans \tH\Psi_{1n}(\balpha_0)\right]\\
		&+\frac{1}{n}\sum_{i\neq j}\left[\tpsi_{i}	- \frac{1}{n}\tK J\Psi_{1n}(\balpha_0)\right] o_p(n^{-1/2})\\
		&+\frac{1}{n}\sum_{i\neq j}\left[\tpsi_{j}-\frac{1}{n}\tK J\Psi_{1n}(\balpha_0)\right] o_p(n^{-1/2})+o_p(n^{-1}).
	\end{align*}
	By notations, we have, for any $i\neq j$,
	\begin{align*}
		\bE[\tK_{i}] = \bzero\quad\mbox{and}\quad \bE[\tpsi_{ij}(\balpha_0)] = 0.
	\end{align*}

	It can be shown by the central limit theorem that
	\begin{align*}
		\frac{1}{n}\sum_{i\neq j}\left[\tpsi_{i}	- \frac{1}{n}\tK J\Psi_{1n}(\balpha_0)\right]
		=\frac{n-1}{n}\sum_{i=1}^n\left[\tpsi_{i}	- \frac{1}{n}\tK J\Psi_{1n}(\balpha_0)\right]
		=O_p(n^{1/2}).
	\end{align*}
	which implies that
	\begin{align*}
		T_n =&\frac{1}{n(n-1)}\sum_{i\neq j}\tpsi_{ij}(\hbalpha)\\
		=&\frac{1}{n(n-1)}\sum_{i\neq j}\left[\tpsi_{ij}(\balpha_0) - \frac{1}{n}\Psi_{1n}(\balpha_0)\trans \tK_{i}- \frac{1}{n}\tK_{j}\trans \Psi_{1n}(\balpha_0) + \frac{1}{n^2}\Psi_{1n}(\balpha_0)\trans \tH\Psi_{1n}(\balpha_0)\right]\\
		&+o_p(n^{-1})\\
		\equiv&T_{n1}+T_{n2}+T_{n3}+o_p(n^{-1}),
	\end{align*}
	where
	\begin{align*}
		T_{n1} =& \frac{1}{n(n-1)}\sum_{i\neq j}\tpsi_{ij}(\balpha_0),\\
		T_{n2} =& -\frac{2}{n^2}\sum_{i=1}^n \Psi_{1n}(\balpha_0)\trans K_{i},\\
		T_{n3} =& \frac{1}{n^2}\Psi_{1n}(\balpha_0)\trans H\Psi_{1n}(\balpha_0).
	\end{align*}
	As the same arguments as the proof of Theorem 2, we have
	\begin{align*}
		T_n - \frac{\mu_0}{n} = \tT_n + R_n +o_p(n^{-1}),
	\end{align*}
	where
	\begin{align*}
		\tT_n = \frac{1}{n(n-1)}\sum_{i\neq j}h(\bV_i,\bV_j)
	\end{align*}
	with the kernel given in (22) in the main paper
	\begin{align*}
		\begin{split}
			h^{(m)}(\bV_i,\bV_j)
			=& \tilde{\omega}_{ij}\psi_0^{(m)}(\bV_i, \balpha_0)\psi_0^{(m)}(\bV_j, \balpha_0)- \psi_1(\bV_i, \balpha_0)\trans \tK_{j}\\
			&- \tK_{i}\trans \psi_1(\bV_j, \balpha_0) + \psi_1(\bV_i, \balpha_0)\trans \tH \psi_1(\bV_j, \balpha_0),
		\end{split}
	\end{align*}
	and $R_n=o_p(n^{-1})$.

	According to Theorem in 5.5.2 of \cite{Serfling1980} by checking it conditions, we can show that
	\begin{align*}
		n\tT_n \lkonv  \zeta,
	\end{align*}
	where $\zeta$ is a random variable of the form $\zeta=\sum_{j=1}^{\infty}\lambda_{j}(\chi^2_{1j}-1)$, and $\chi^2_{11},\chi^2_{12},...$ are independent $\chi^2_{1}$ variables, that is, $\zeta$ has characteristic function
	\begin{align*}
		E\left[e^{it\zeta}\right]=\prod_{j=1}^{\infty}(1-2it\lambda_{j})^{-\frac{1}{2}}e^{-it\lambda_{j}}.
	\end{align*}

\begin{thm}\label{coroll32}
	If Assumptions~{(A1)--(A6)} hold, then under the local alternative hypothesis that is $\bbeta=n^{-1/2}\bxi$ with a deterministic vector $\bxi\in\Theta_{\beta}$, we have
	\begin{align*}
		nT_n^{(m)}-\mu_0 \lkonv \zeta,
	\end{align*}
	where $\mu_0$ is defined in Theorem~\ref{coroll31}, $\zeta$ is a random variable of the form $\zeta=\sum_{j=1}^{\infty}\lambda_{j}(\chi^2_{1j}(\mu_{aj})-1)$, and $\chi^2_{11}(\mu_{a1}), \chi^2_{12}(\mu_{a2}), \cdots$ are independent noncentral $\chi^2_{1}$ variables.
	Here, $\{\lambda_{j}\}$ are the eigenvalues of the kernel $h^{(m)}(\bv_1, \bv_2)$ defined in \eqref{kernel0c} under $f(\bv, \balpha_0, \bzero, \btheta_0)$, i.e., they are the solutions of $\lambda_{j}g_{j}(\bv_2)=\int_{0}^{\infty}h^{(m)}(\bv_1, \bv_2)g_{j}(\bv_1)f(\bv_1, \balpha_0, \bzero, \btheta_0)d\bv_1$ for nonzero $g_{j}$, and each noncentrality parameter of $\chi^2_{1j}(\mu_{aj})$ is
	\begin{align*}
		\mu_{aj} = \bE\left[\phi_j(\bV_{0})\bxi\trans\partial \log(f(\bV_{0}, \balpha_0, \bzero, \btheta_0))/\partial\bbeta\right], \quad j=1, 2, \cdots,
	\end{align*}
	where $\{\phi_j(\bv)\}$ denote orthonormal eigenfunctions corresponding to the eigenvalues $\{\lambda_j\}$, and $\bV_0$ is generated from the null distribution $f(\bv, \balpha_0, \bzero, \btheta_0)$.
\end{thm}

The proof of this theorem can be completed along the same lines of Theorem 2 and Theorem \ref{coroll31}.

Next, we demonstrate the finite-sample performance of the proposed WAST in multiple change-plane models. 
Consider the GLM in example~2 with the canonical parameter
\begin{align*}
\mu = \tbX\trans\balpha + \bX\trans\bbeta_1\bone(\bZ_1\trans\btheta_1\geq0)+\bX\trans\bbeta_2\bone(\bZ_2\trans\btheta_2\geq0),
\end{align*}
where $\balpha$ is the same as that in Section~5 in the main paper, $\bbeta_1=\bbeta_2 = \tau\bone_q$, and $\tbX$, $\bX$, and $\bZ_1=\bZ_2$ are generated as in Section~5 in the main paper. Let $\bZ_{-1,s}$ be $\bZ_s$ excluding the first component, and let $\btheta_{-1,s}$ be $\btheta_s$ excluding the first component, $s=1,2$. $\theta_{11}$ and $\theta_{21}$ are chosen as the negatives of the 0.3 and 0.6 percentiles of $\bZ_{-1,1}\trans\btheta_{-1,1}$ and $\bZ_{-1,1}\trans\btheta_{-1,1}$, respectively, which means that $\bZ_{1}\trans\btheta_{1}$ divides the population into two groups with 0.7 and 0.3 observations, and $\bZ_{1}\trans\btheta_{1}$ divides the population into two groups with 0.4 and 0.6 observations. When $(r,p,q)=(2,2,3)$ or $(r,p,q)=(6,6,3)$, the maximal value of the range $\kappa_{\max}$ is $\log(0.7)$, $\log(0.35)$ and $\log(0.7)$ for Gaussian, binomial and Poisson regression, respectively. While $(r,p,q)=(2,2,11)$ or $(r,p,q)=(6,6,11)$, the maximal value of the range $\kappa_{\max}$ is $\log(0.9)$, $\log(0.75)$, and $\log(0.9)$ for Gaussian, binomial, and Poisson regression, respectively. The other settings are the same as in Section~5 in the main paper.

\begin{figure}[!ht]
\begin{center}
\includegraphics[scale=0.3]{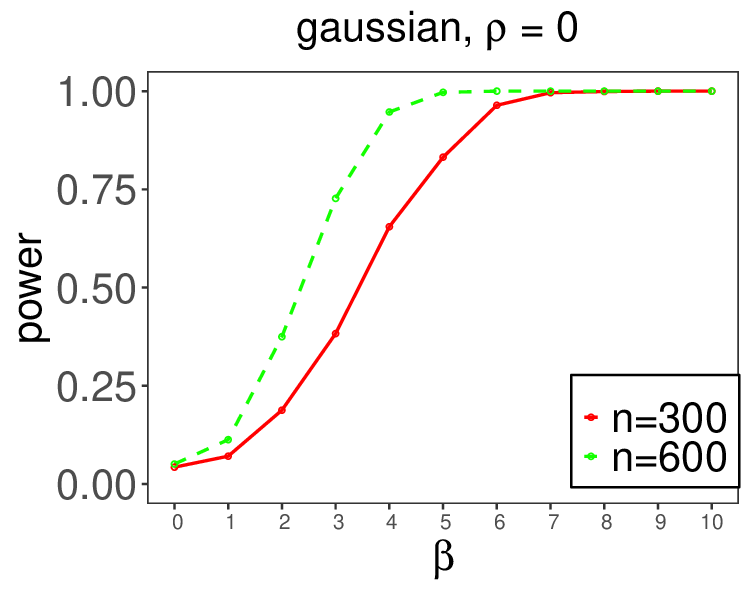}
\includegraphics[scale=0.3]{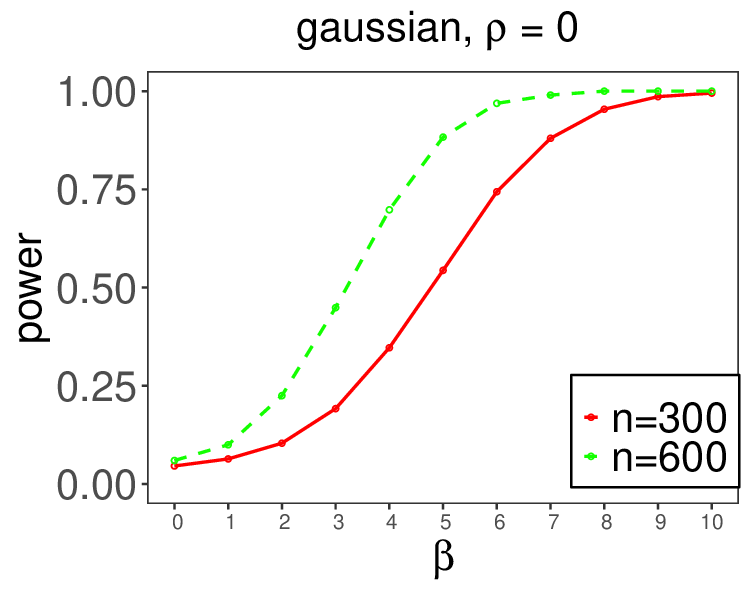}
\includegraphics[scale=0.3]{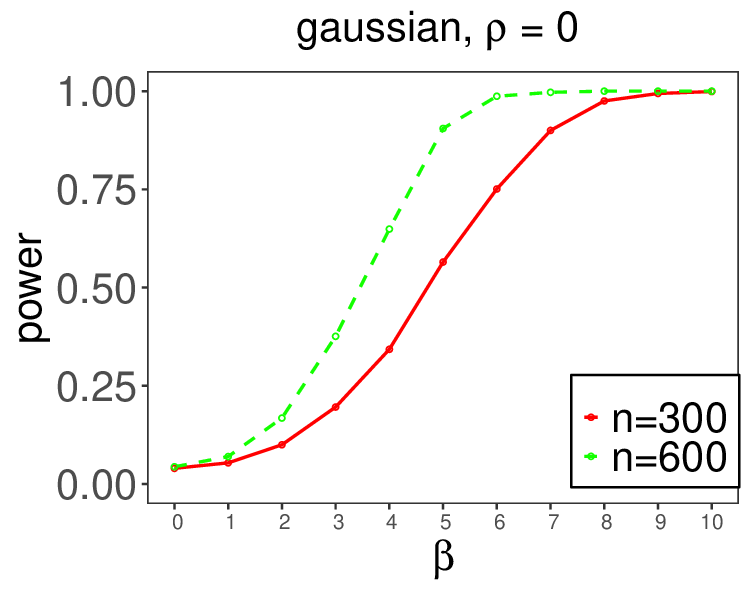}
\includegraphics[scale=0.3]{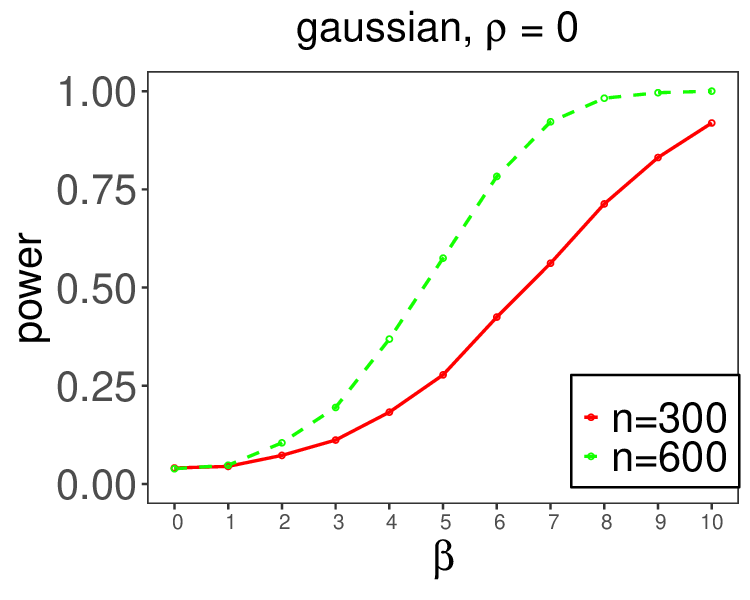}
\includegraphics[scale=0.3]{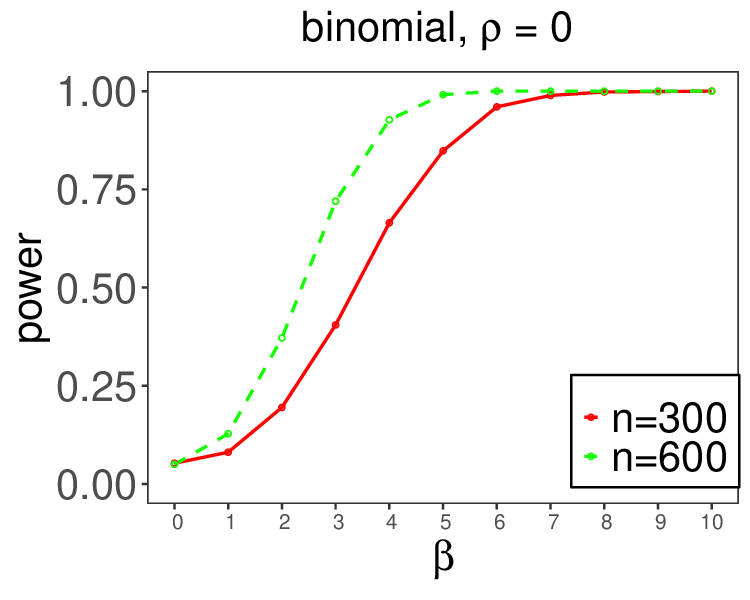}
\includegraphics[scale=0.3]{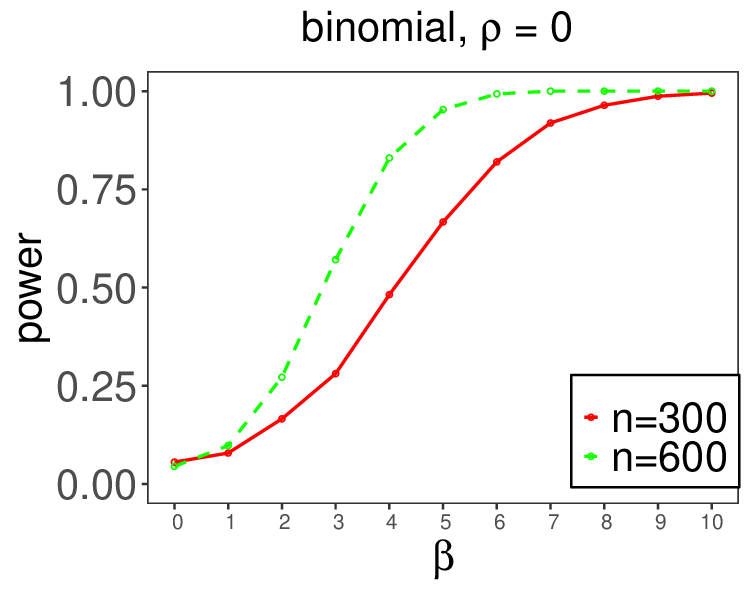}
\includegraphics[scale=0.3]{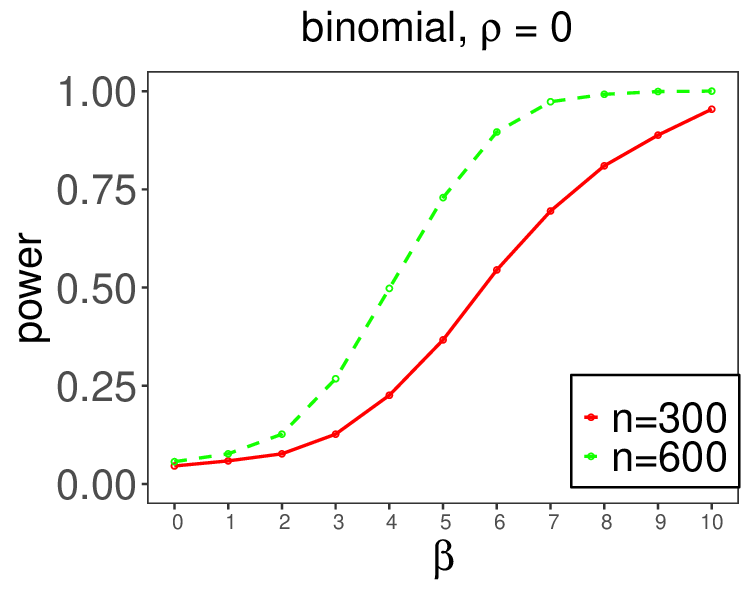}
\includegraphics[scale=0.3]{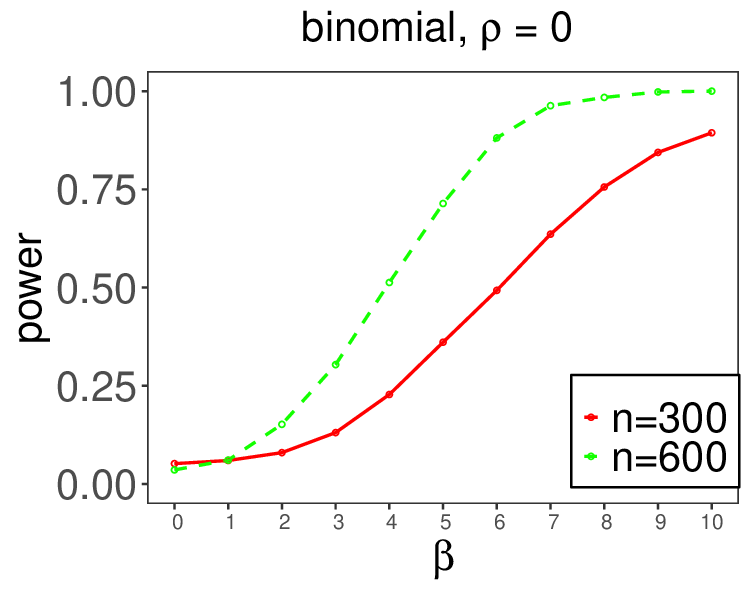}
\includegraphics[scale=0.3]{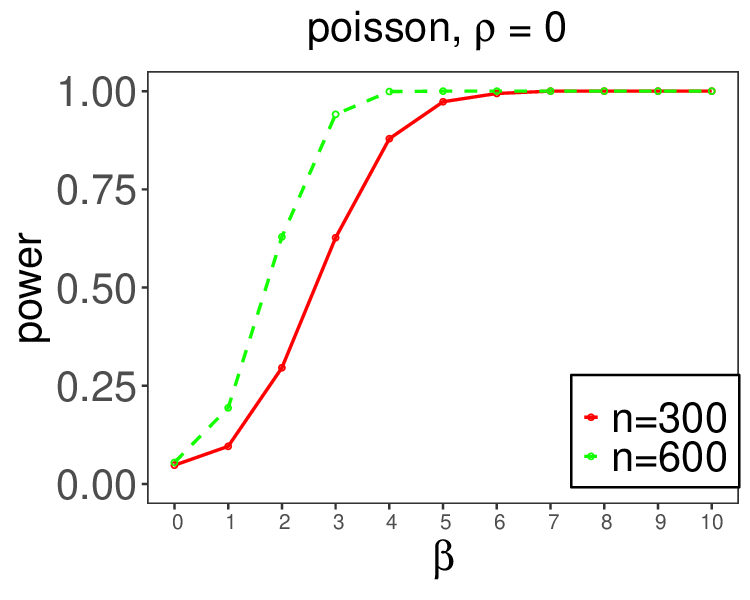}
\includegraphics[scale=0.3]{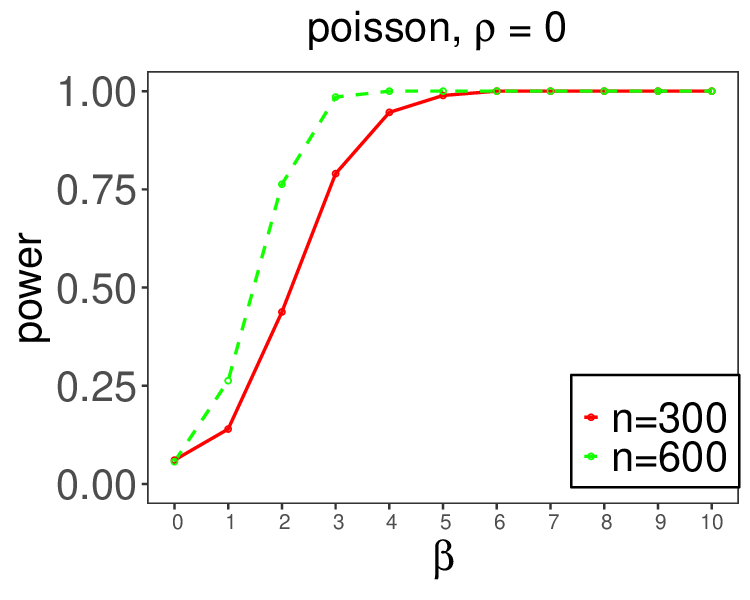}
\includegraphics[scale=0.3]{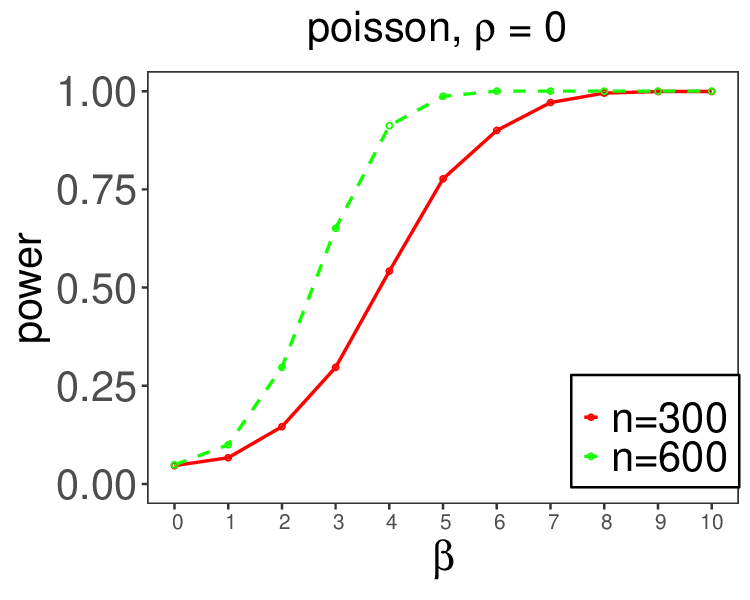}
\includegraphics[scale=0.3]{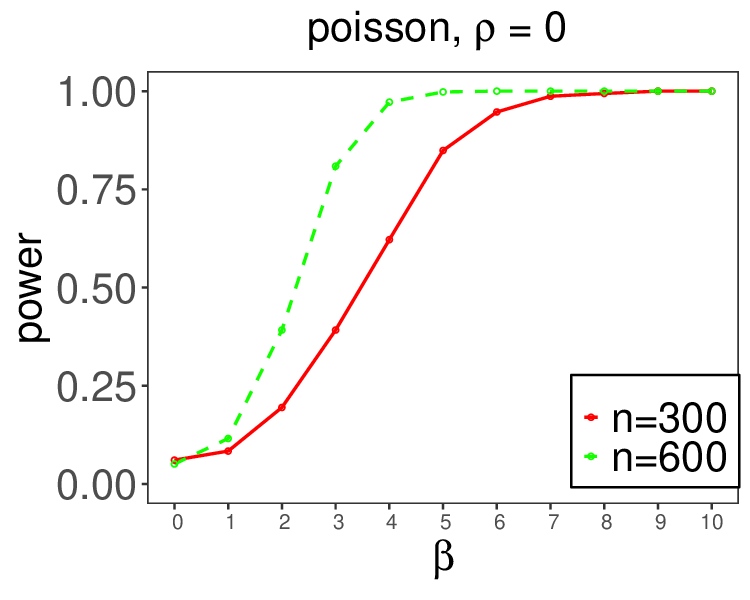}
\includegraphics[scale=0.3]{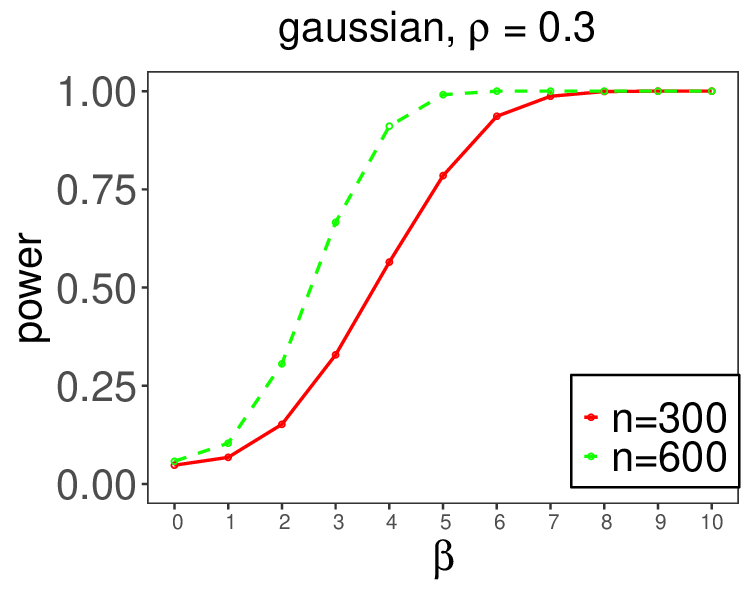}
\includegraphics[scale=0.3]{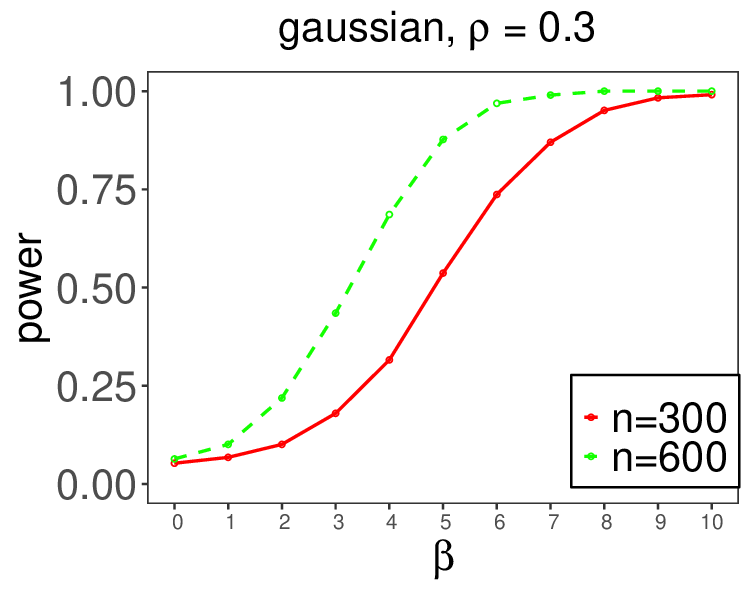}
\includegraphics[scale=0.3]{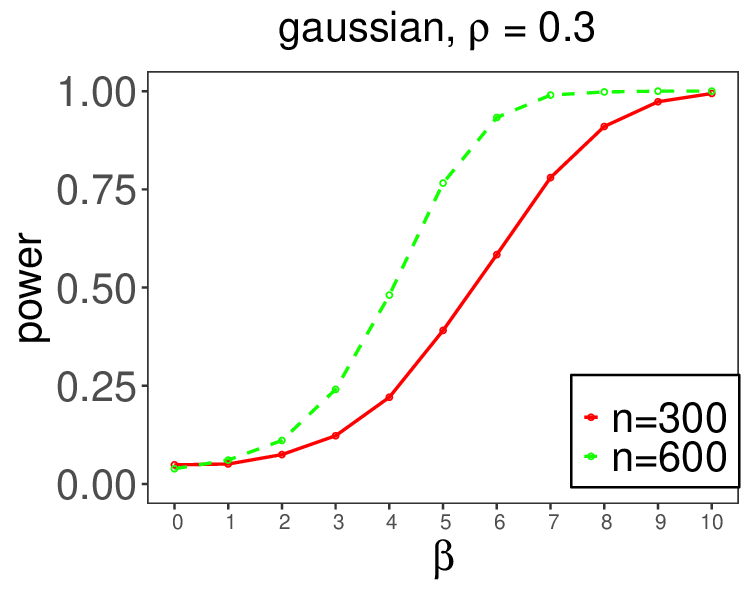}
\includegraphics[scale=0.3]{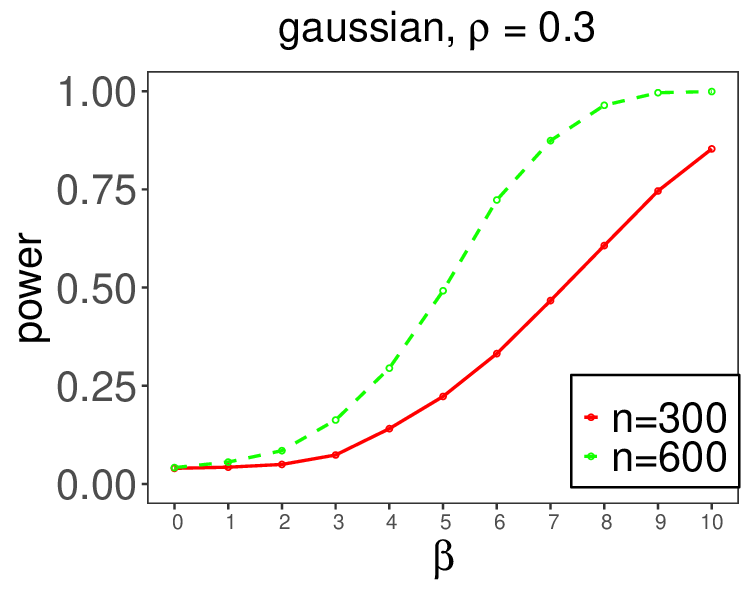}
\includegraphics[scale=0.3]{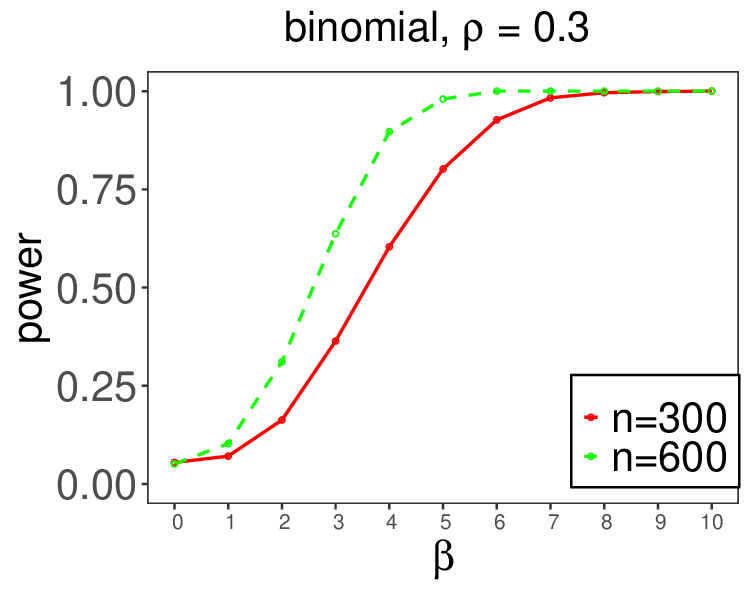}
\includegraphics[scale=0.3]{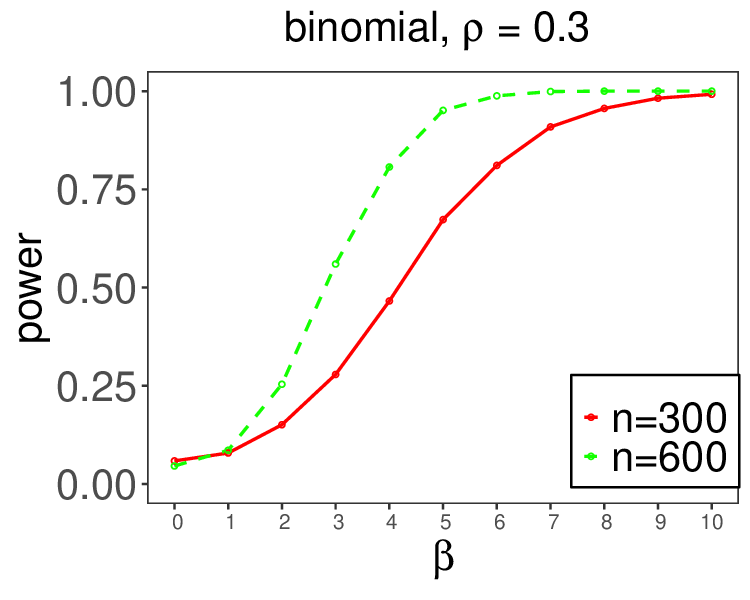}
\includegraphics[scale=0.3]{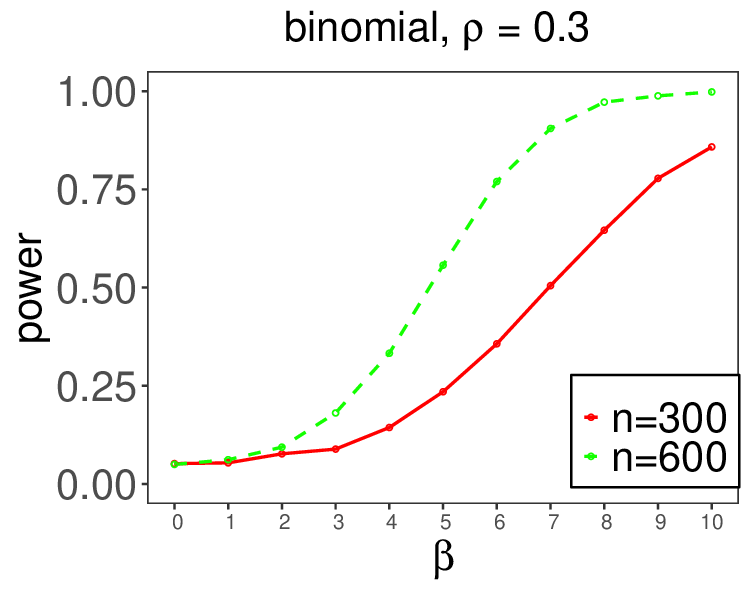}
\includegraphics[scale=0.3]{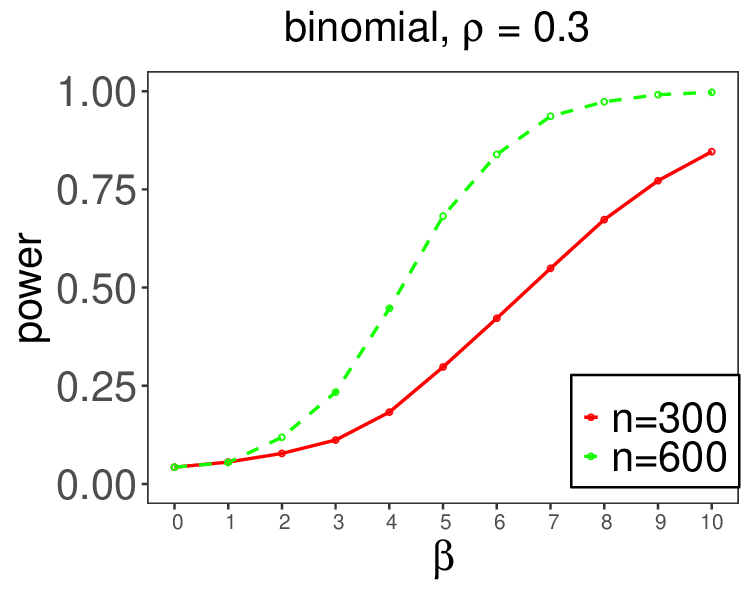}
\includegraphics[scale=0.3]{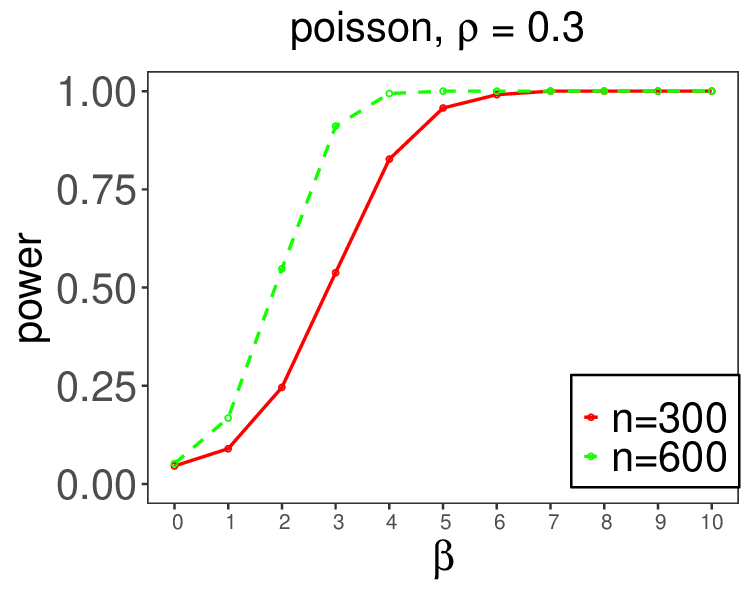}
\includegraphics[scale=0.3]{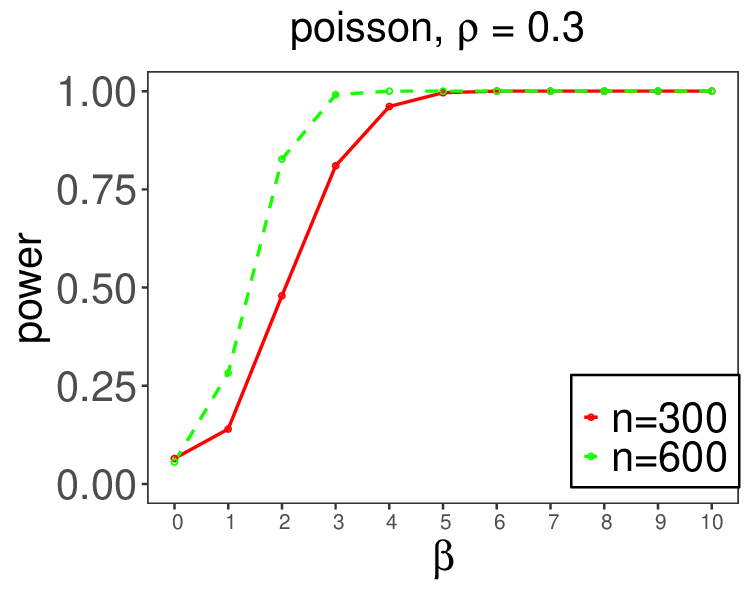}
\includegraphics[scale=0.3]{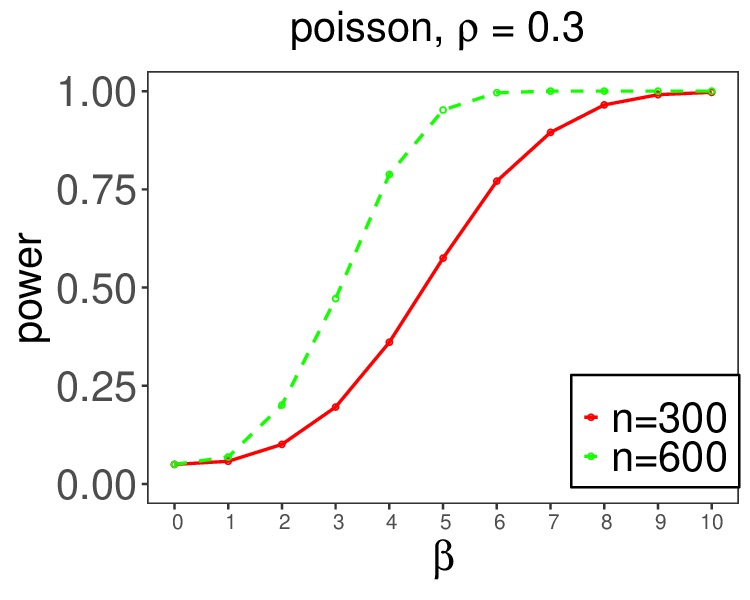}
\includegraphics[scale=0.3]{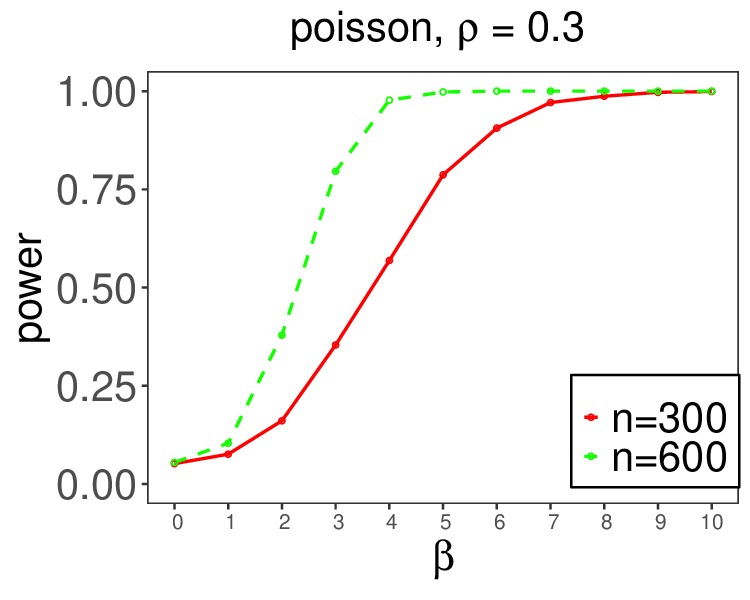}
\caption{Powers of proposed WAST with $n=300$ (red solid line) and $n=600$ (green dashed line) for multiple change planes. From left to right, each column depicts the powers for the cases $(r,p,q)=(2,2,3)$, $(6,6,3)$, $(2,2,11)$, and $(6,6,11)$.}
\label{fig_multisubg}
\end{center}
\end{figure}

Table~\ref{table_multisubg} lists the sizes of the proposed WAST with multiple change planes; as can be seen, the sizes are close to the nominal level of 0.05, as expected. Figure~\ref{fig_multisubg} shows the power of the proposed WAST with multiple change planes; as can be seen, all the powers (i) increase as the sample size increases and (ii) increase rapidly as $\bbeta$ departs from the null hypothesis.

\begin{table}
\caption{\label{table_multisubg} Type-\uppercase\expandafter{\romannumeral1} errors of proposed WAST for multiple change planes.}
\centering
\begin{threeparttable}
\begin{tabular}{llccccccc}\\
\hline
\multirow{2}{*}{$(r,p,q)$}& \multirow{2}{*}{$n$}
& \multicolumn{3}{c}{ $\rho=0$} && \multicolumn{3}{c}{ $\rho=0.3$} \\
\cline{3-5} \cline{7-9}
&& Gaussian& Binomial & Poisson&& Gaussian& Binomial & Poisson\\
\cline{3-9}
$(2,2,3)$& $300$ & 0.043 & 0.048 & 0.053 && 0.055 & 0.048 & 0.046 \\
& $600$ & 0.051 & 0.058 & 0.050 && 0.051 & 0.055 & 0.052 \\
[1 ex]
$(6,6,3)$& $300$ & 0.046 & 0.053 & 0.056 && 0.059 & 0.061 & 0.065 \\
& $600$ & 0.060 & 0.064 & 0.045 && 0.046 & 0.057 & 0.056 \\
[1 ex]
$(2,2,11)$& $300$ & 0.040 & 0.049 & 0.046 && 0.052 & 0.047 & 0.050 \\
& $600$ & 0.044 & 0.039 & 0.057 && 0.050 & 0.049 & 0.050 \\
[1 ex]
$(6,6,11)$& $300$ & 0.041 & 0.040 & 0.052 && 0.043 & 0.061 & 0.052 \\
& $600$ & 0.039 & 0.042 & 0.036 && 0.043 & 0.051 & 0.055 \\
\hline
\end{tabular}
\begin{tablenotes}
\item The settings are from Section~\ref{scp2} with $\Gv$ $\bZ$ generated from a normal distribution with zero mean and unit variance. The nominal significance level is 0.05.
\end{tablenotes}
\end{threeparttable}
\end{table}

\subsection{Change plane analysis for quantile, probit and semiparametric models}\label{simulation_cp}
For fair comparison, in \ref{simulation_cp}--\ref{simulation_cpz_cauchy}, we consider the settings as similar to the existing literature \citet{2011Testing,2017Change}, which are different from the GLMs in Section 5 in the main paper.

We consider the probit regression model in the main paper
\begin{align*}
	f(\bV_i, \balpha, \bbeta, \btheta)=\Phi(0.5+\tX_{1i}+\bX_i\trans\bbeta\bone(\bZ_i\trans\btheta\geq0))^{Y_i}\Phi(-0.5-\tX_{1i}-\bX_i\trans\bbeta\bone(\bZ_i\trans\btheta\geq0))^{1-Y_i}.
\end{align*}
We generate $X_{1i}$ independently from binomial distribution $b(1,0.5)$, $\bX_i$ from $p$-variate normal distribution with mean $\bzero$ and variance $\sqrt{2}I$, and $Z_{1i}=1$ and $(Z_{2i},\cdots,Z_{qi})\trans$ from $(q-1)$-variate standard normal distribution.
We set $\bbeta=(1,\cdots,1)\trans$ under $H_0$, and $(\theta_2, \cdots,\theta_q)\trans = (1,2,\cdots,2)\trans$ under $H_1$, where $\theta_1$ is chosen as the negative of the 0.65 percentile of $Z_2\theta_2+\cdots+Z_q\theta_q$, which means that $\bZ\trans\btheta$ divides the population into two groups with 0.35 and 0.65 observations, respectively.

Consider the quantile regression model in the main paper,
\begin{align*}
	Y_i = 0.5+\tX_i+\bX_i\trans\bbeta(\tau)\bone(\bZ_i\trans\btheta(\tau)\geq0) + \eps_i,
\end{align*}
where $\tX_i$ is generated from standard normal distribution, $\bX_i$ and $\bZ_i$ are generated as same as these in probit regression model. The error $\eps_i$ is generated from standard normal distribution.
The parameters $\bbeta$ under $H_0$ and $\btheta$ under $H_1$ are same as these in probit model.

Consider the semiparametric model in the main paper,
\begin{align*}
	Y_i = \gamma(\tbX_i,\balpha_2) +\bX_i\trans\bbeta A_i \bone(\bZ_i\trans \btheta\geq0)+\eps_i.
\end{align*}
As the same settings provided by \cite{2017Change}, we consider two settings for the propensity score model $\pi(\tbX_i,\balpha_1)$, denoted by P1 and P2:
\begin{enumerate}[({P}1)]
	\item $\pi(\tbX_i,\balpha_1)=0.5$;
	\item $\pi(\tbX_i,\balpha_1)=\frac{\exp(0.5\tX_{2i}+0.5\tX_{3i})}{1+\exp(0.5\tX_{2i}+0.5\tX_{3i})}$,
\end{enumerate}
and two baseline mean functions for $\gamma(\tbX_i,\balpha_2)$, denoted by B1 and B2:
\begin{enumerate}[({B}1)]
	\item $\gamma(\tbX_i,\balpha_2) = 1+0.5\tX_{2i}+\tX_{3i}^2$;
	\item $\gamma(\tbX_i,\balpha_2) = 1+\sin(\tX_{2i}+\tX_{3i})$.
\end{enumerate}
We calculate the WAST and SST in the main paper by fitting a linear model for baseline function $\gamma(\tbX_i,\balpha_2)$ and a logistic regression for  the propensity score model $\pi(\tbX_i,\balpha_1)$. Therefore, two baseline models $\gamma(\tbX_i,\balpha_2)$ are misspecified. We generate $v_{1i}$ independently from binomial distribution $b(1,0.5)$, $v_{2i}$ independently from uniform distribution $U(-1,1)$, and $X_{ji}$ independently from uniform distribution $U(-1,1)$ , where $i=1,\cdots,n$ and $j=1,\cdots,p$. We set $\tbX_i=(1, v_{1i}, v_{2i})\trans$ and $\bZ_i=(1, v_{1i}, v_{2i}, X_{1i},\cdots,X_{(p-3)i})\trans$.
The parameters $\bbeta$ under $H_0$ and $\btheta$ under $H_1$ are same as these in probit model.

\begin{figure}[!ht]
	\begin{center}
		\includegraphics[scale=0.33]{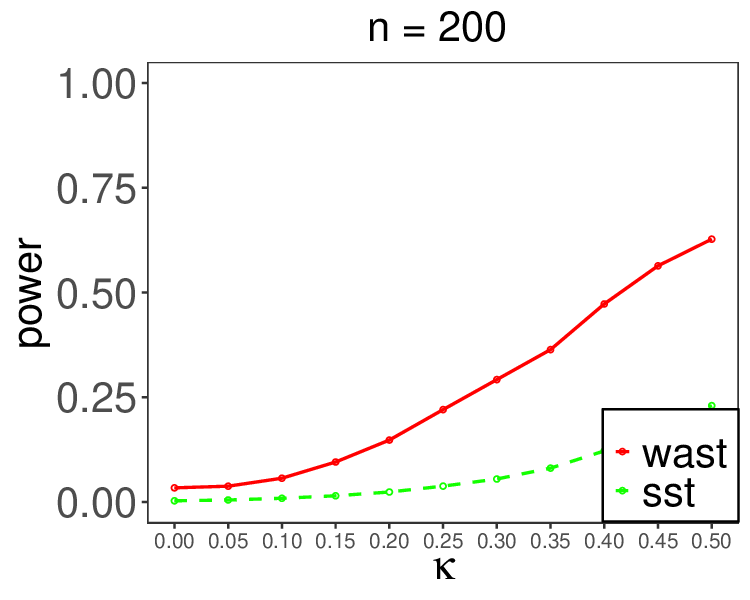}
		\includegraphics[scale=0.33]{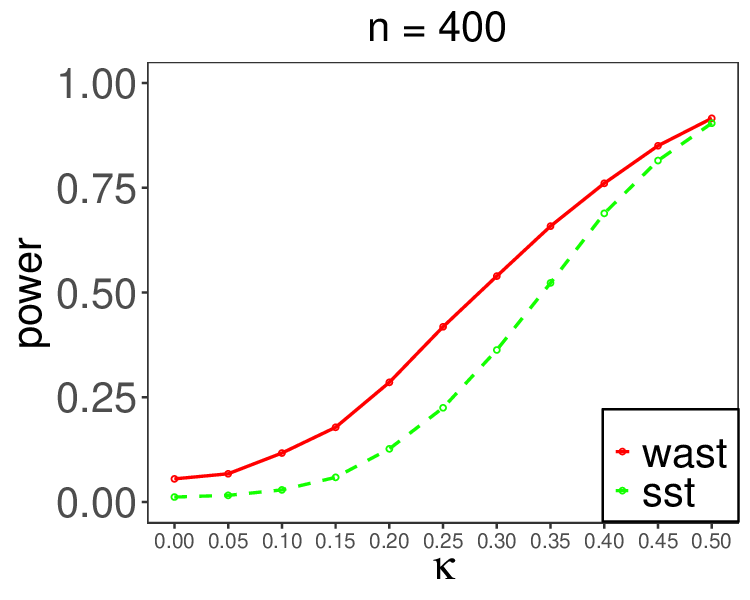}
		\includegraphics[scale=0.33]{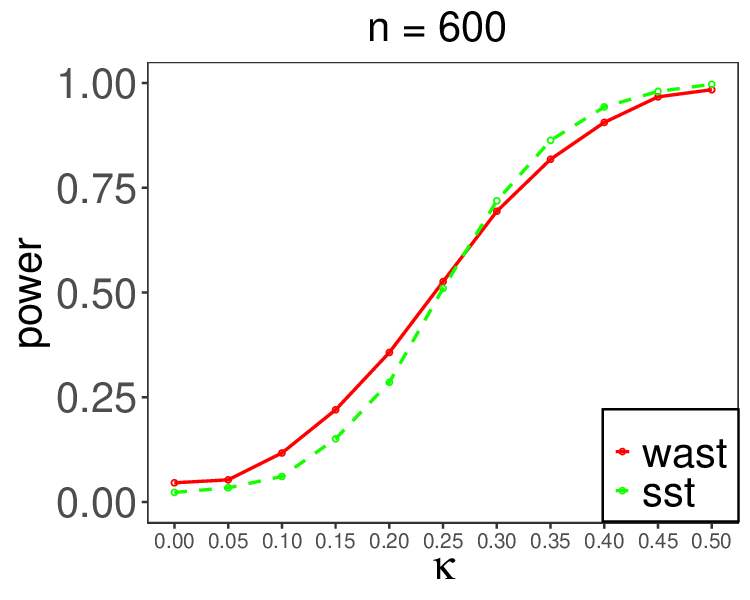}      \\
		\includegraphics[scale=0.33]{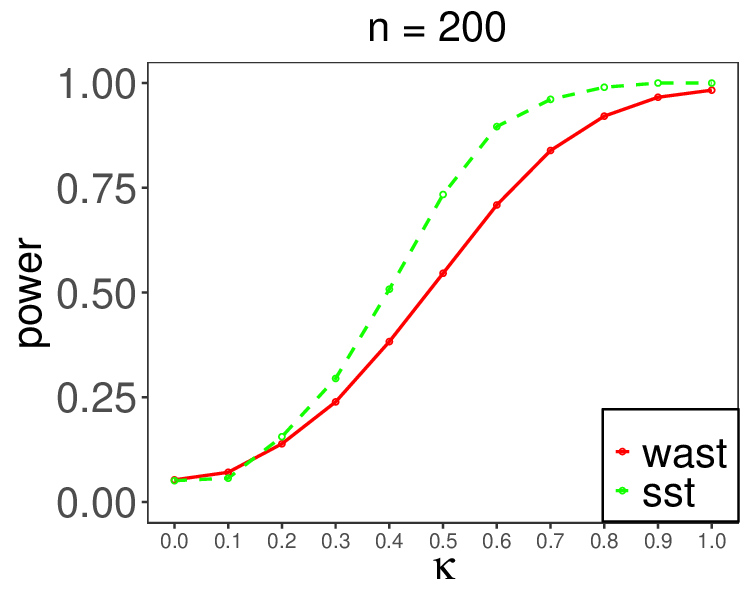}
		\includegraphics[scale=0.33]{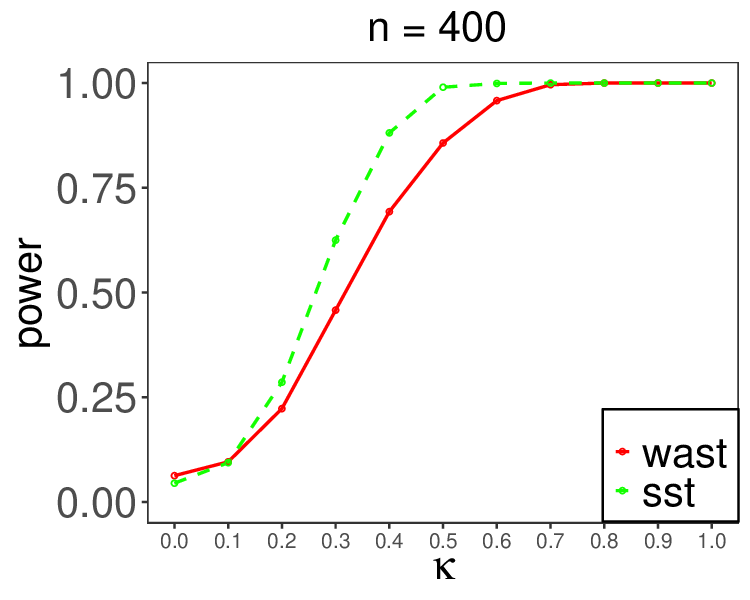}
		\includegraphics[scale=0.33]{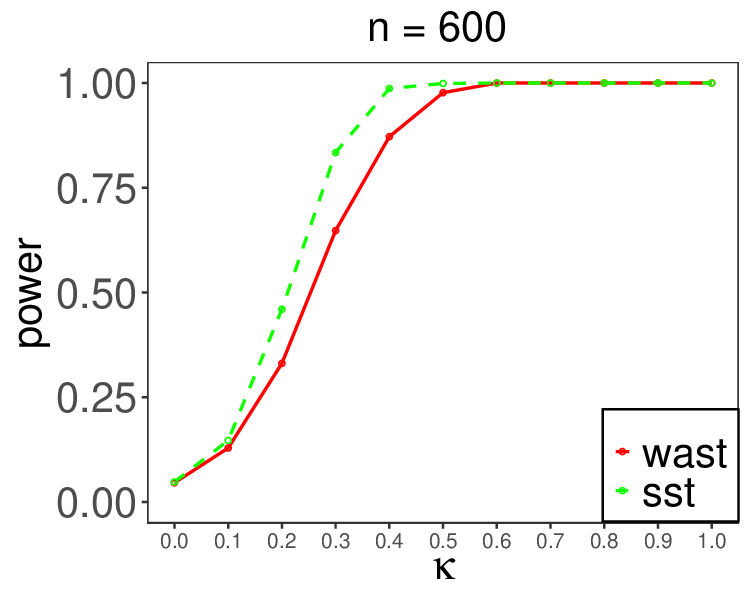}         \\
		\includegraphics[scale=0.33]{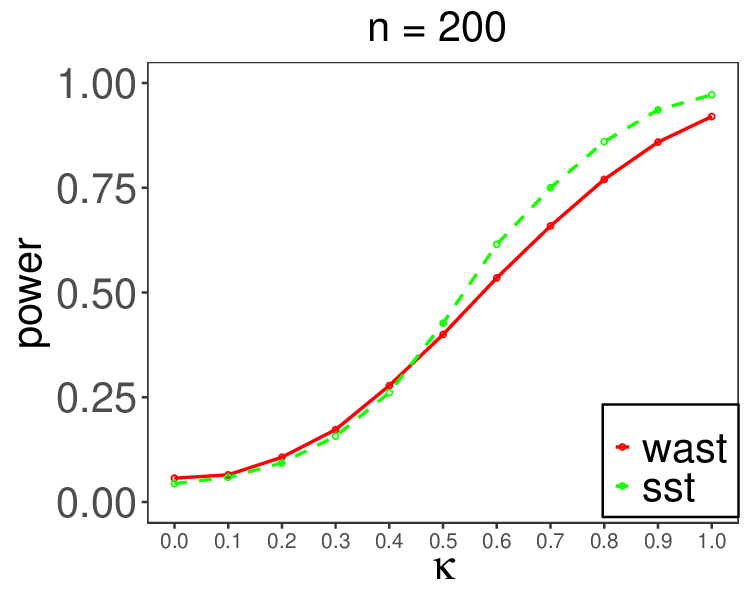}
		\includegraphics[scale=0.33]{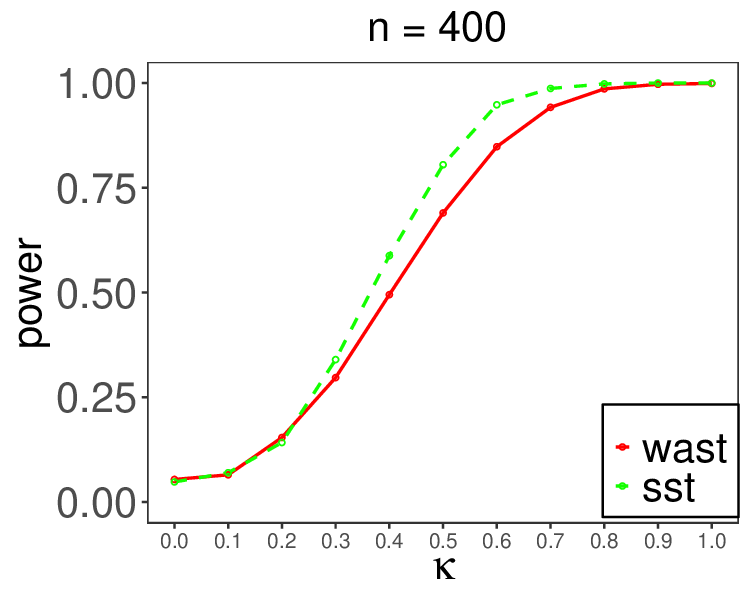}
		\includegraphics[scale=0.33]{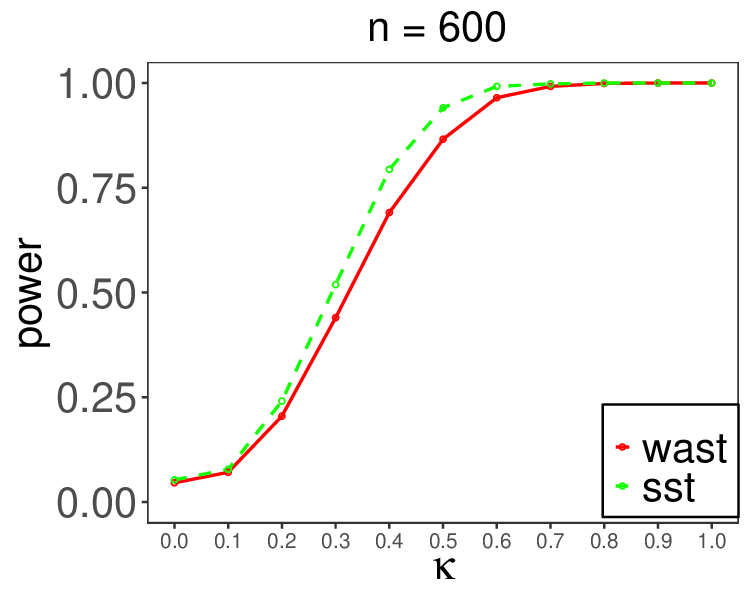} \\
		\includegraphics[scale=0.33]{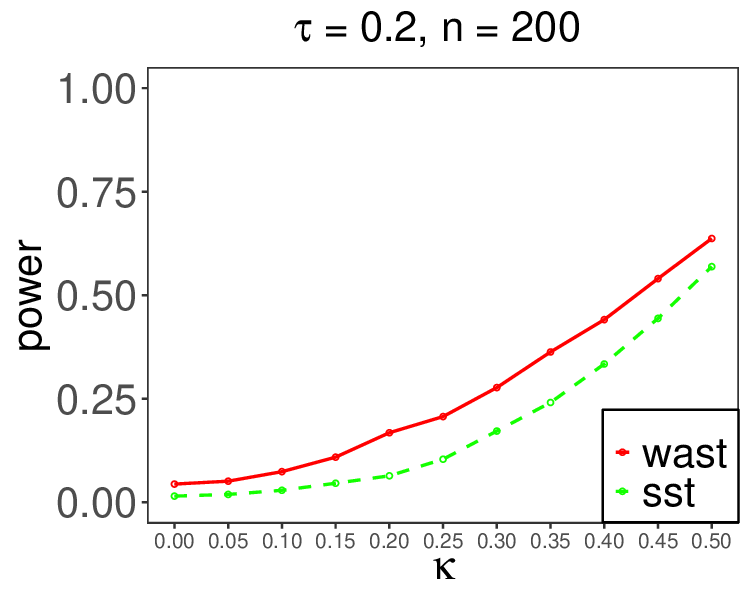}
		\includegraphics[scale=0.33]{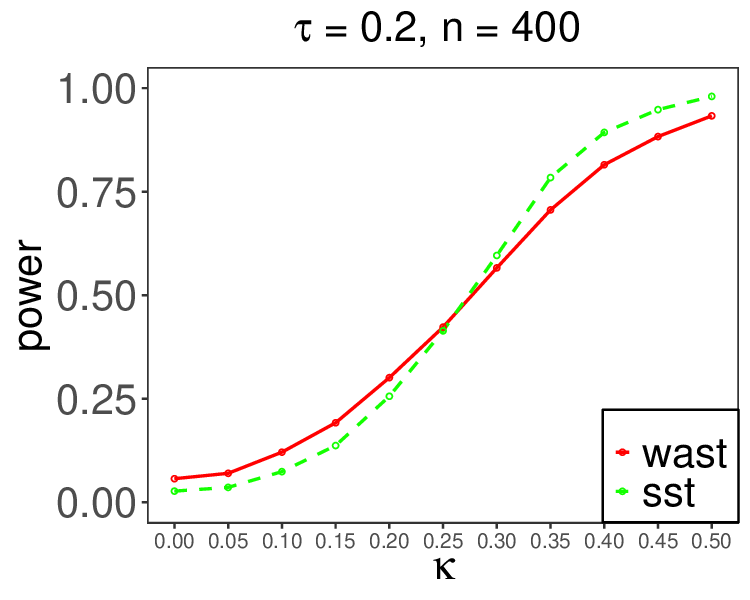}
		\includegraphics[scale=0.33]{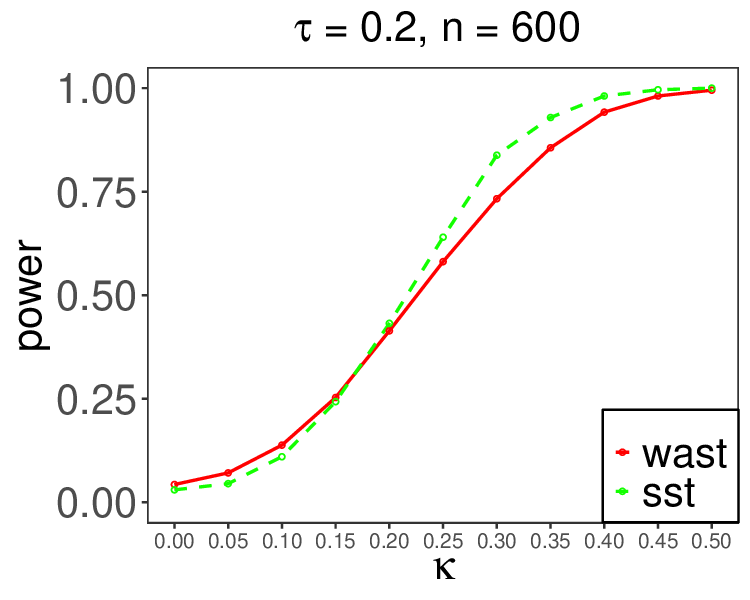}  \\
		\includegraphics[scale=0.33]{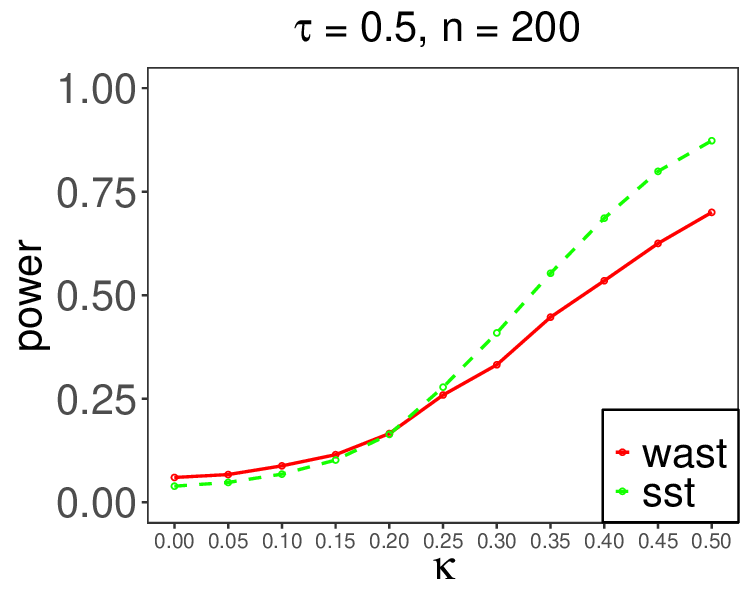}
		\includegraphics[scale=0.33]{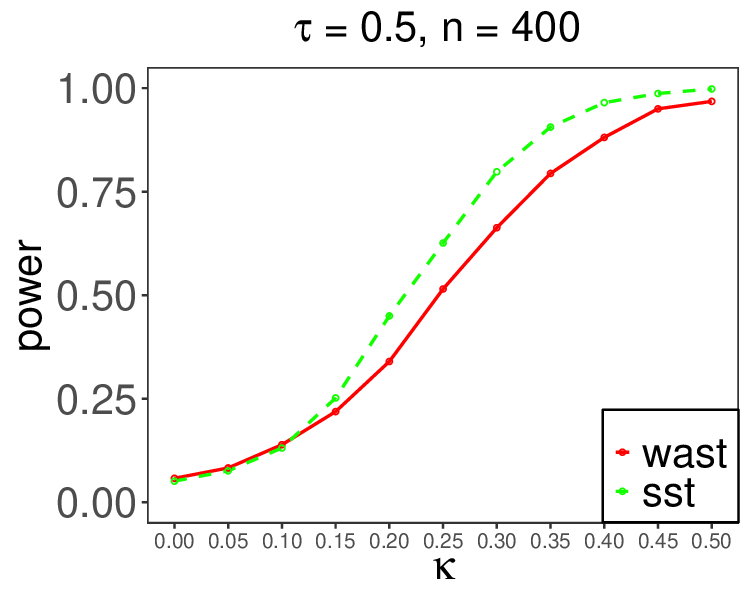}
		\includegraphics[scale=0.33]{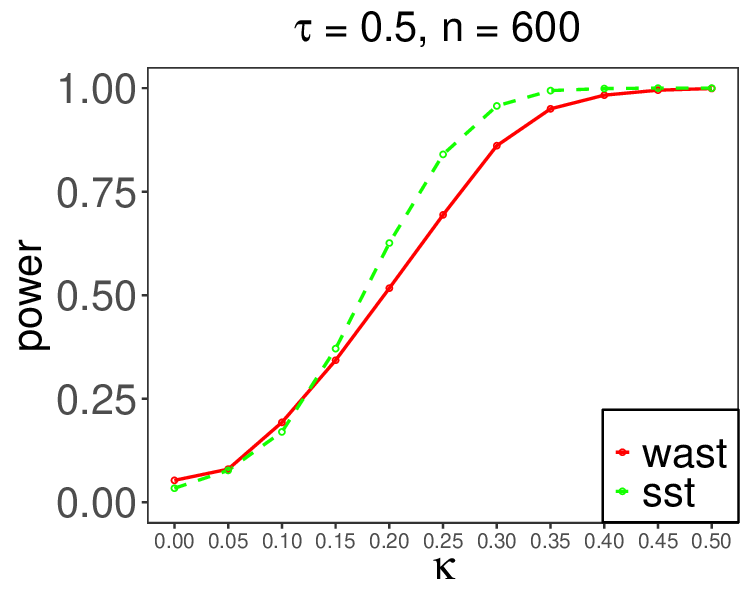}  \\
		\includegraphics[scale=0.33]{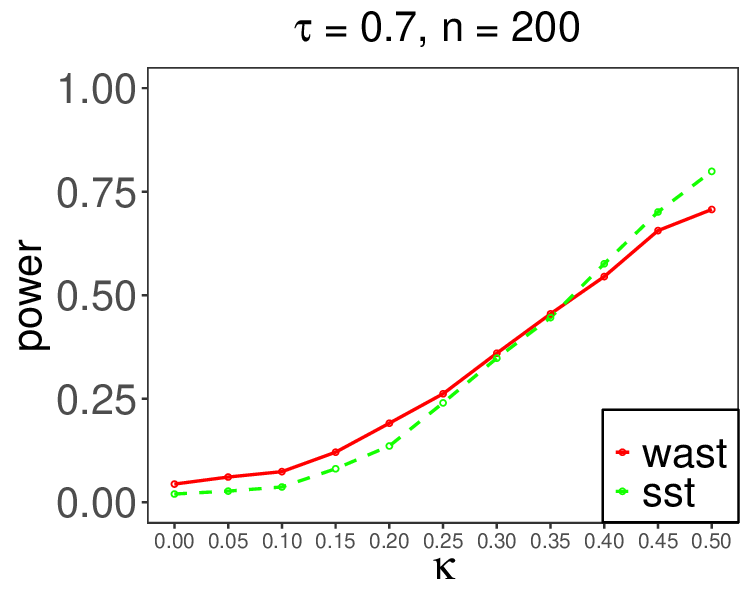}
		\includegraphics[scale=0.33]{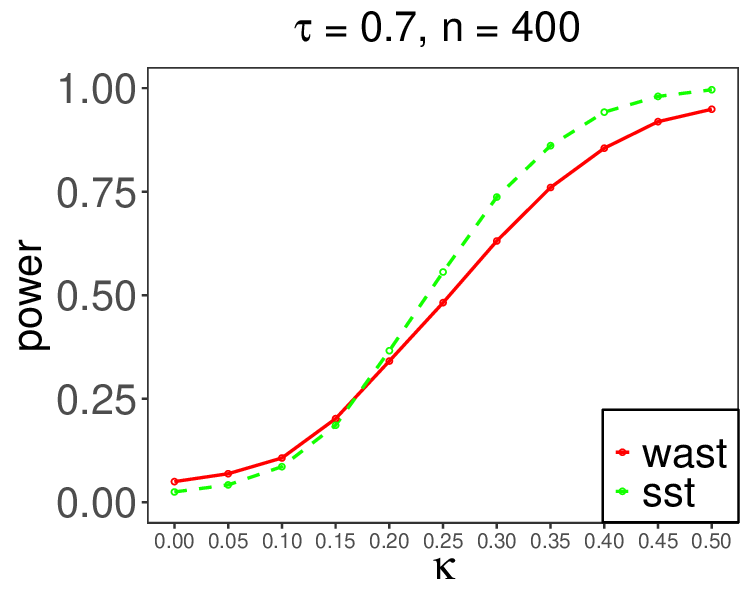}
		\includegraphics[scale=0.33]{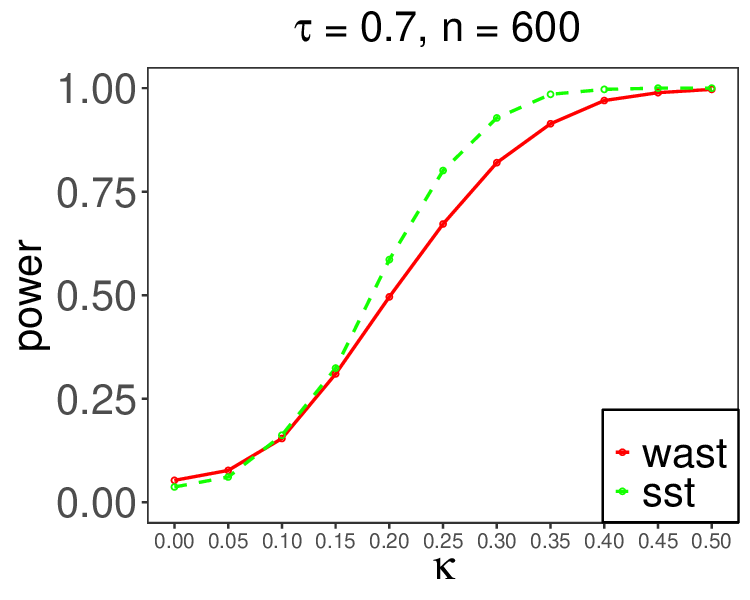}
		\caption{\it Powers of test statistic by the proposed WAST (red solid line) and SST (green dashed line) for $(p,q)=(1,3)$. From top to bottom, each row depicts the powers for probit model, semiparametric model, quantile regression with $\tau=0.2$, $\tau=0.5$ and $\tau=0.7$, respectively. Here the $\Gv$ $Z$ is generated from multivariate normal distribution with mean $\bzero$ and covariance $I$.}
		\label{fig_qr13}
	\end{center}
\end{figure}

\begin{figure}[!ht]
	\begin{center}
		\includegraphics[scale=0.33]{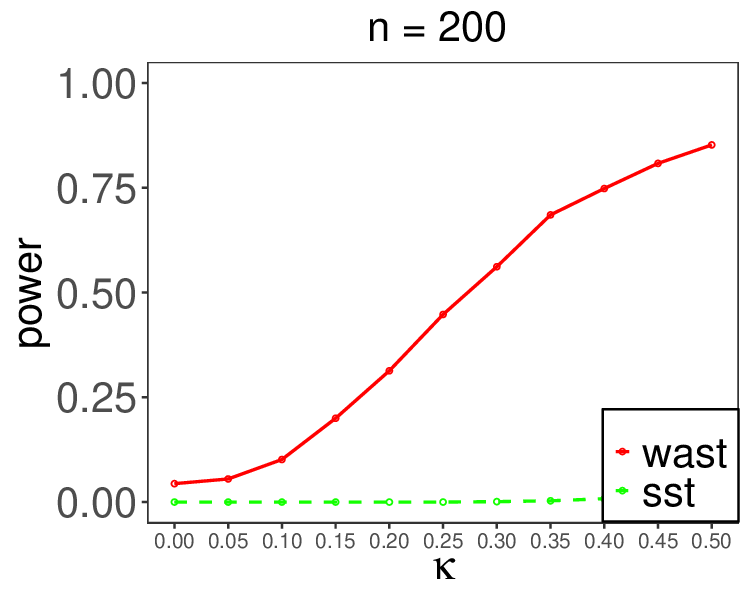}
		\includegraphics[scale=0.33]{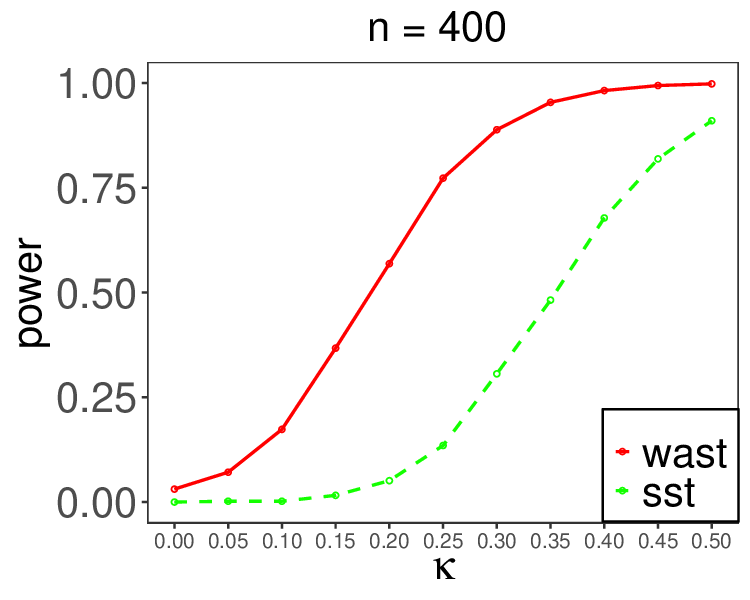}
		\includegraphics[scale=0.33]{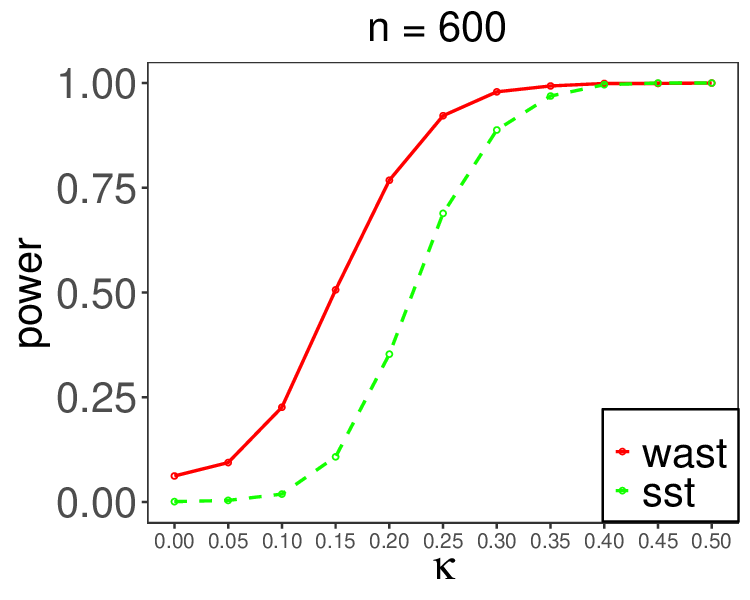}     \\
		\includegraphics[scale=0.33]{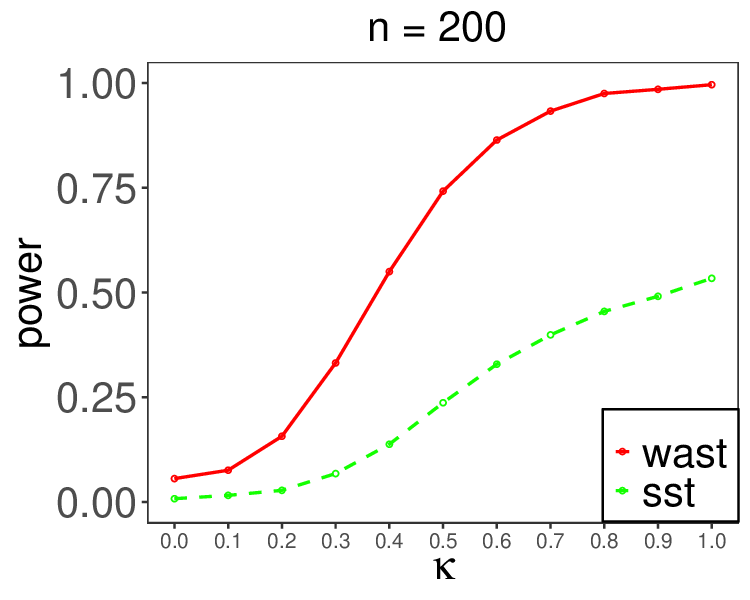}
		\includegraphics[scale=0.33]{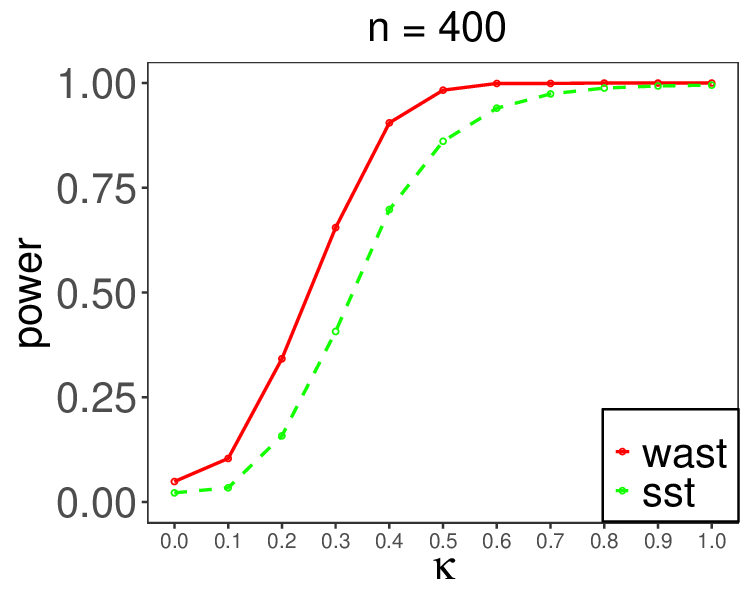}
		\includegraphics[scale=0.33]{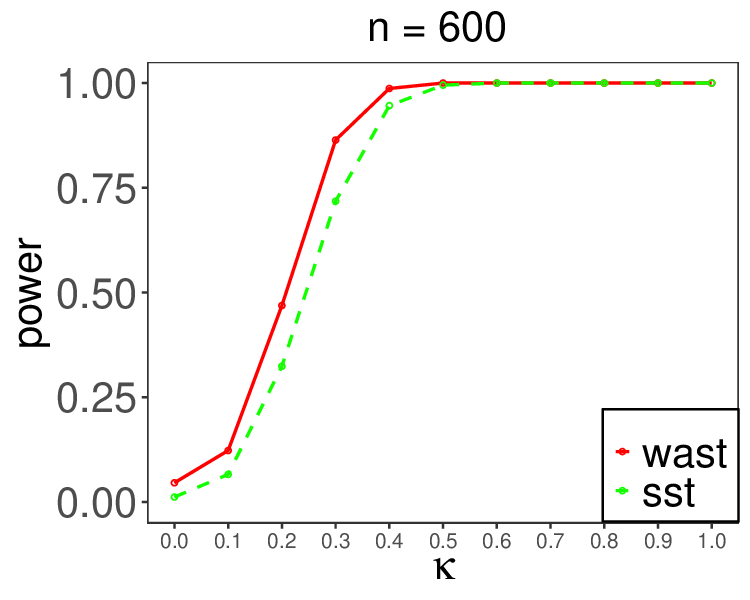}        \\
		\includegraphics[scale=0.33]{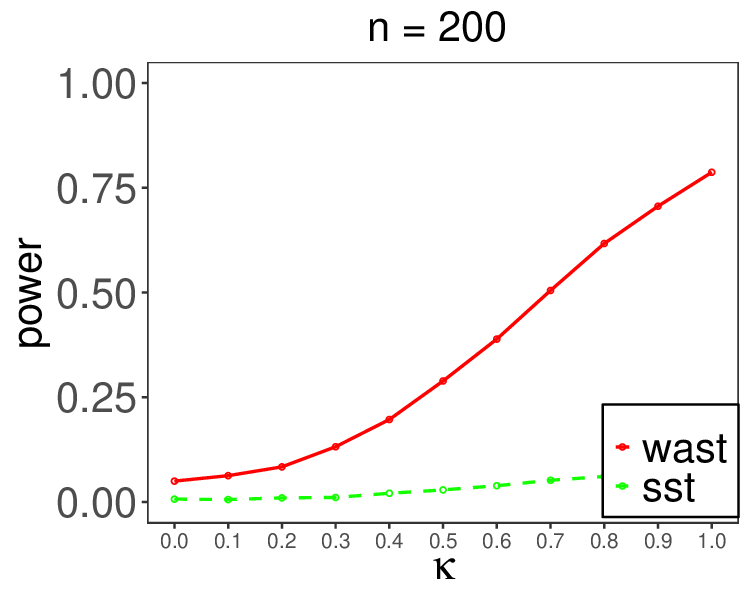}
		\includegraphics[scale=0.33]{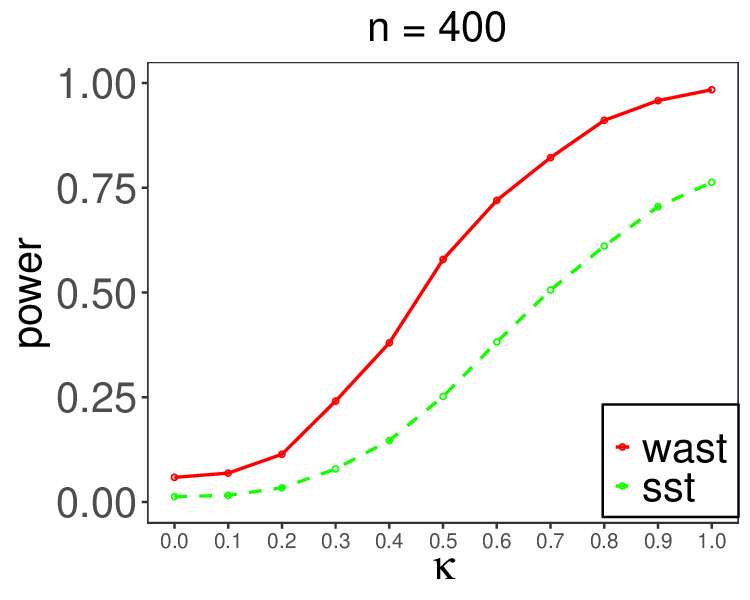}
		\includegraphics[scale=0.33]{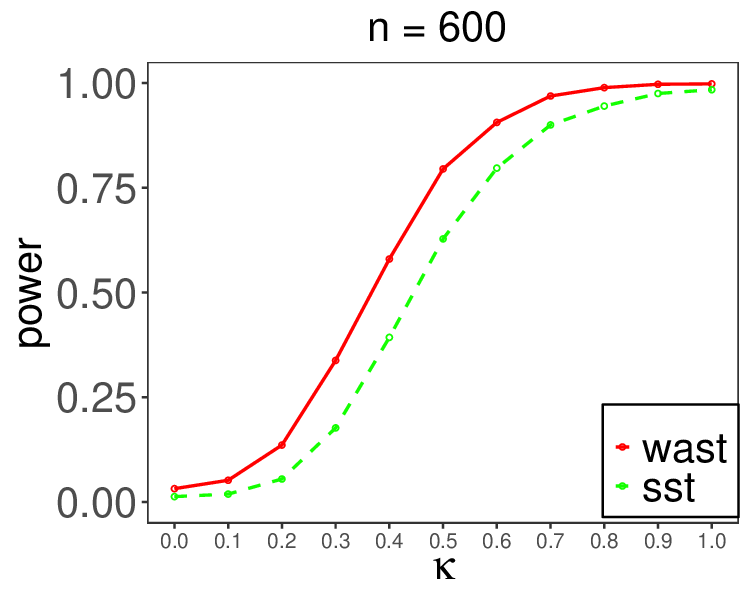} \\
		\includegraphics[scale=0.33]{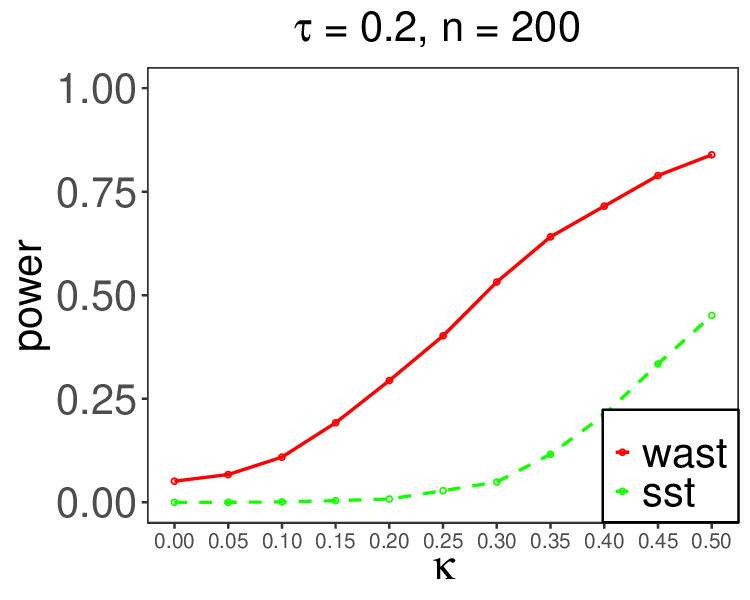}
		\includegraphics[scale=0.33]{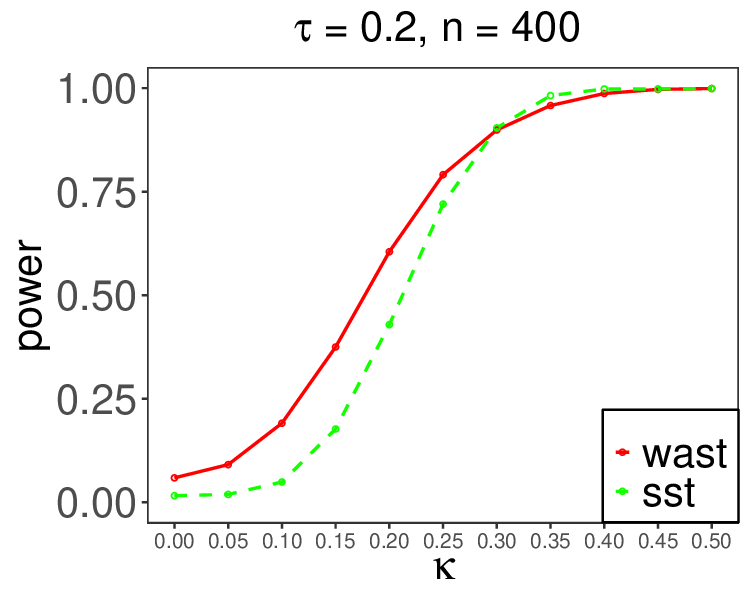}
		\includegraphics[scale=0.33]{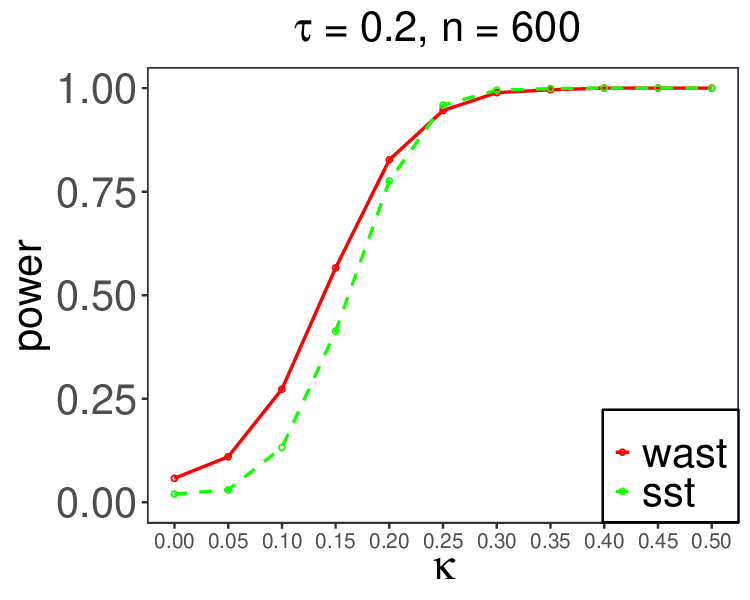} \\
		\includegraphics[scale=0.33]{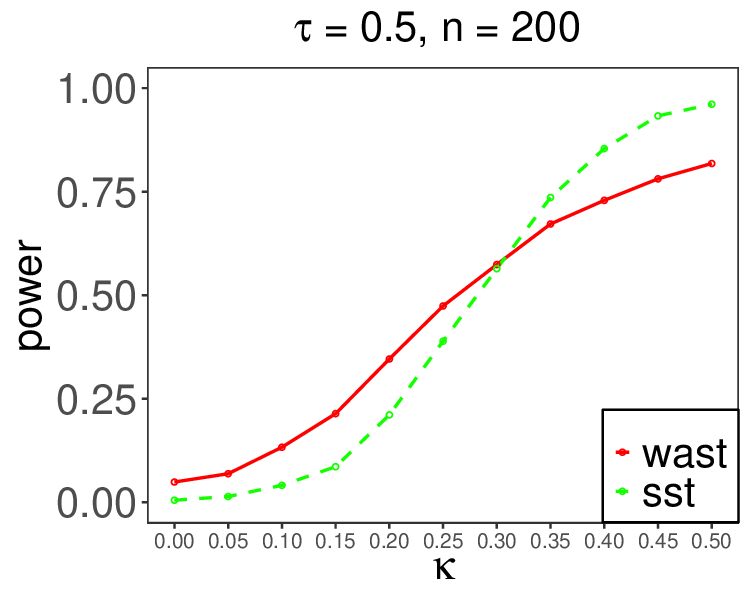}
		\includegraphics[scale=0.33]{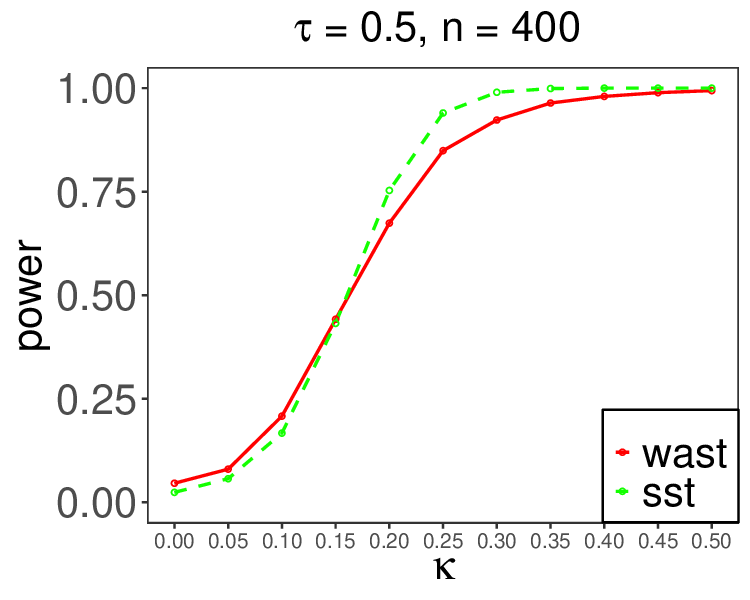}
		\includegraphics[scale=0.33]{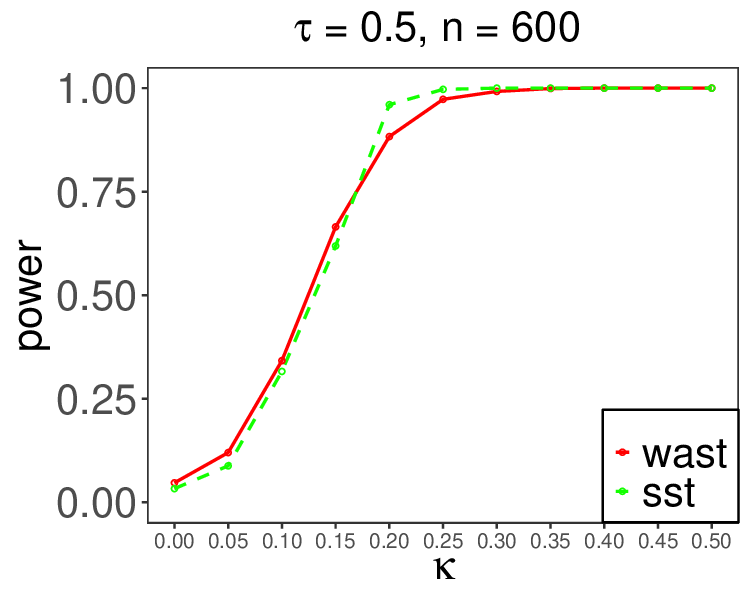} \\
		\includegraphics[scale=0.33]{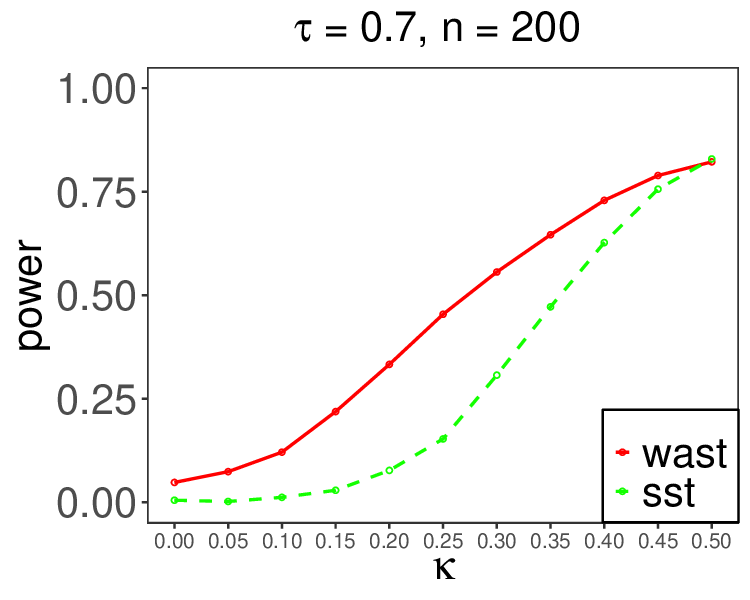}
		\includegraphics[scale=0.33]{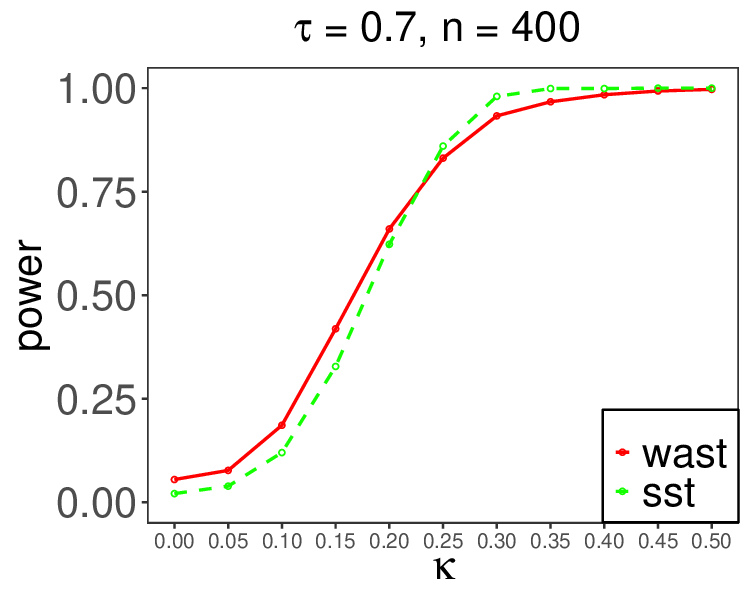}
		\includegraphics[scale=0.33]{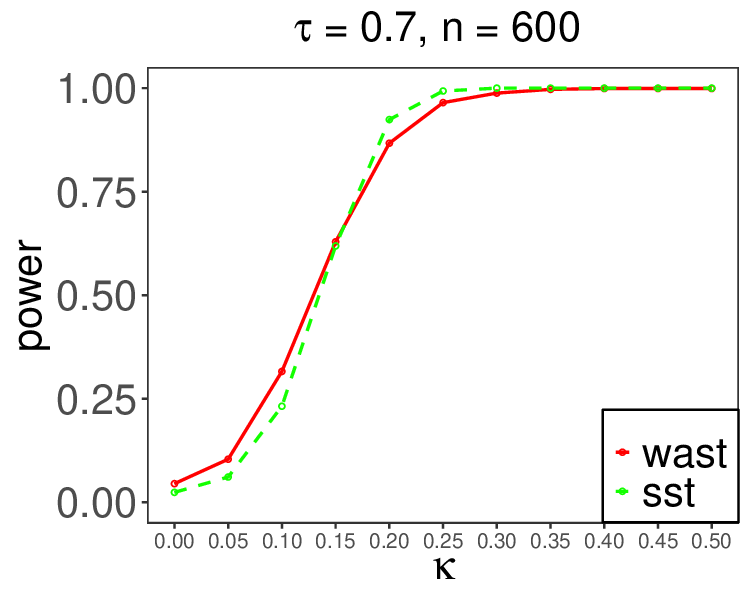}
		\caption{\it Powers of test statistic by the proposed WAST (red solid line) and SST (green dashed line) for $(p,q)=(5,5)$. From top to bottom, each row depicts the powers for probit model, semiparametric model, quantile regression with $\tau=0.2$, $\tau=0.5$ and $\tau=0.7$, respectively. Here the $\Gv$ $Z$ is generated from multivariate normal distribution with mean $\bzero$ and covariance $I$.}
		\label{fig_qr55}
	\end{center}
\end{figure}

\begin{figure}[!ht]
	\begin{center}
		\includegraphics[scale=0.33]{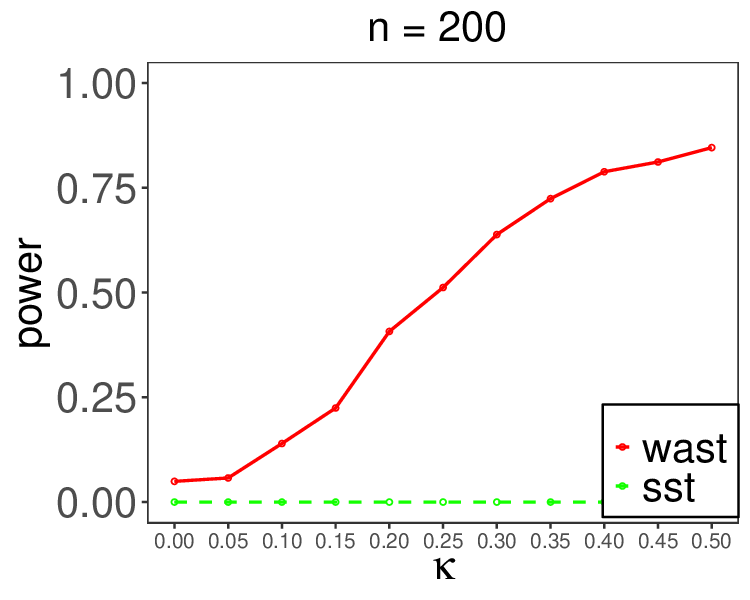}
		\includegraphics[scale=0.33]{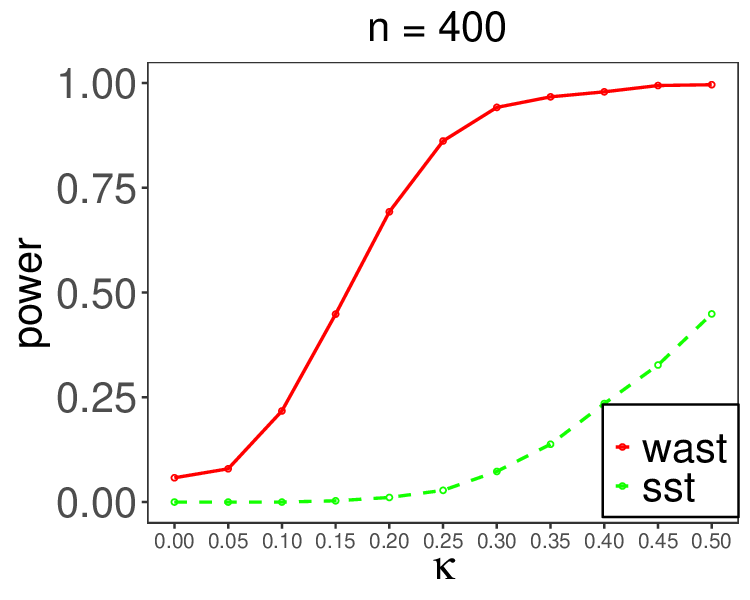}
		\includegraphics[scale=0.33]{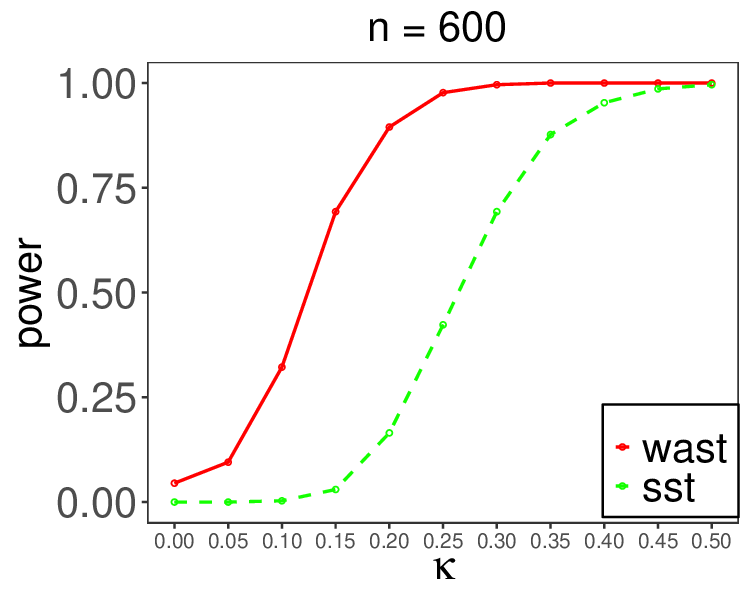}       \\
		\includegraphics[scale=0.33]{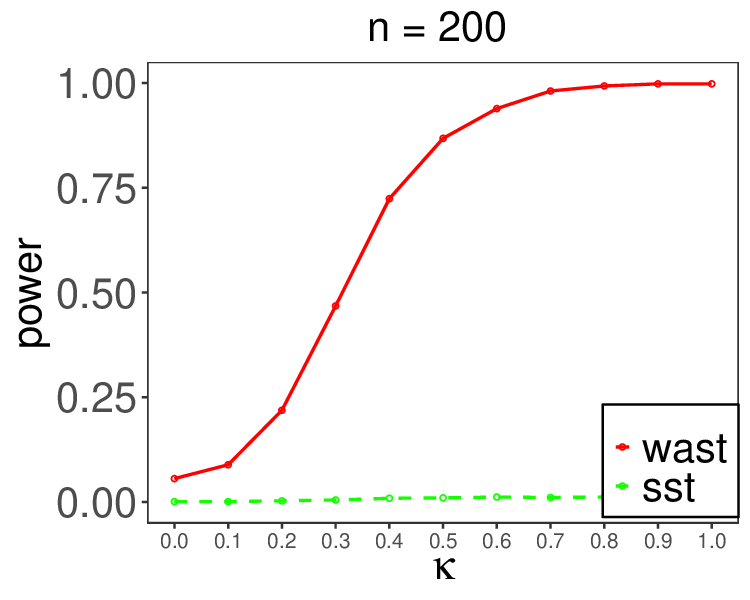}
		\includegraphics[scale=0.33]{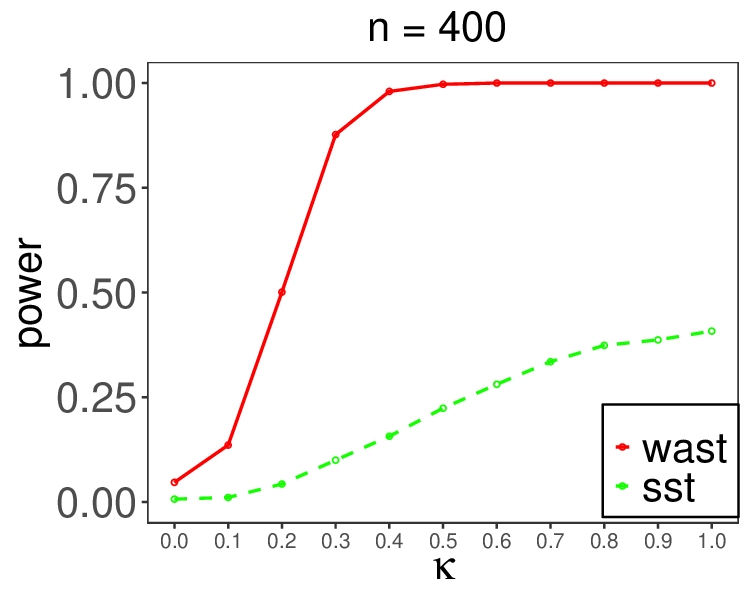}
		\includegraphics[scale=0.33]{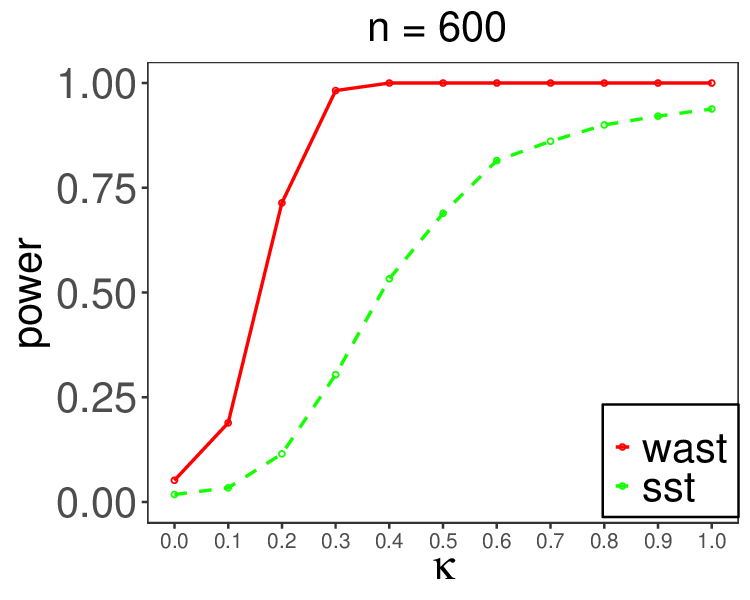}    \\
		\includegraphics[scale=0.33]{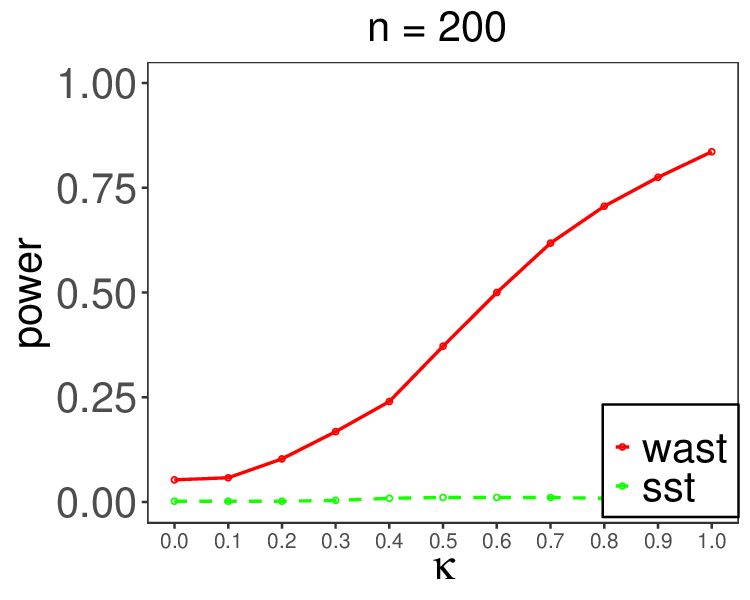}
		\includegraphics[scale=0.33]{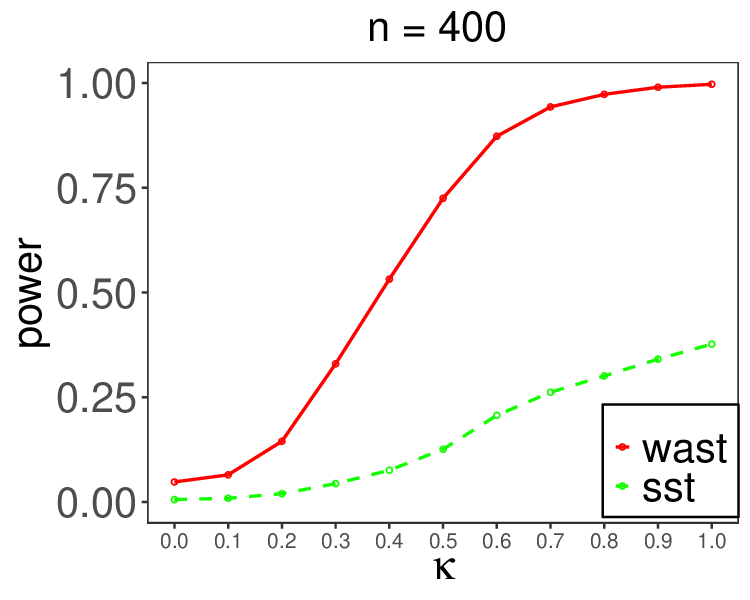}
		\includegraphics[scale=0.33]{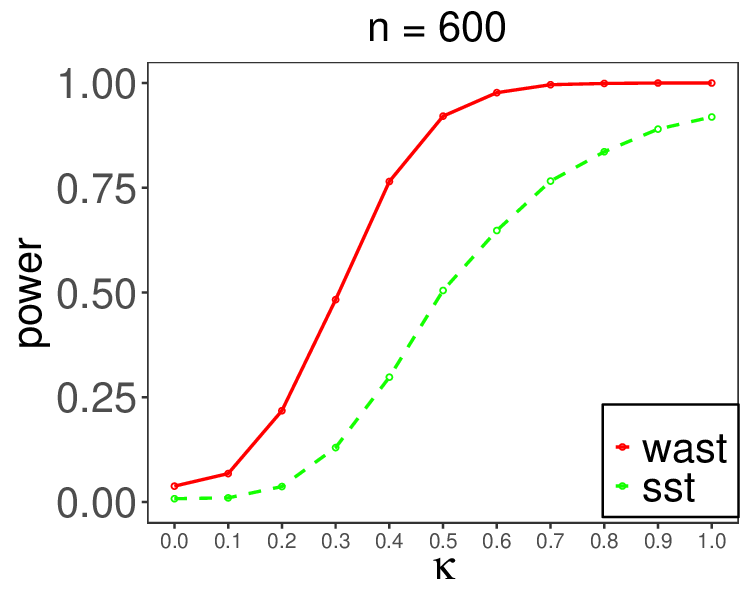}   \\
		\includegraphics[scale=0.33]{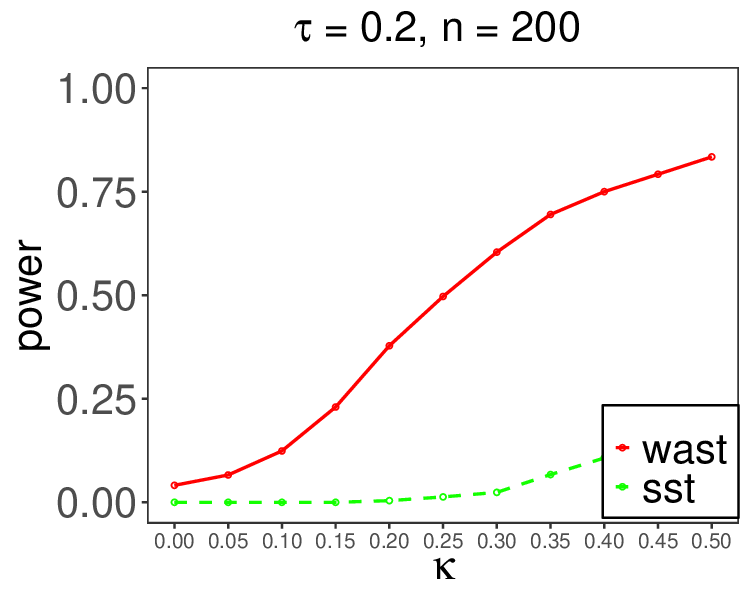}
		\includegraphics[scale=0.33]{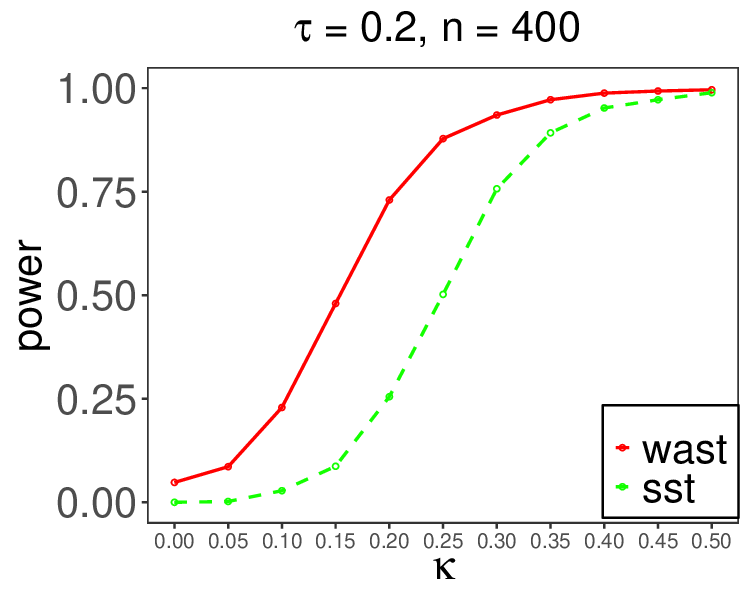}
		\includegraphics[scale=0.33]{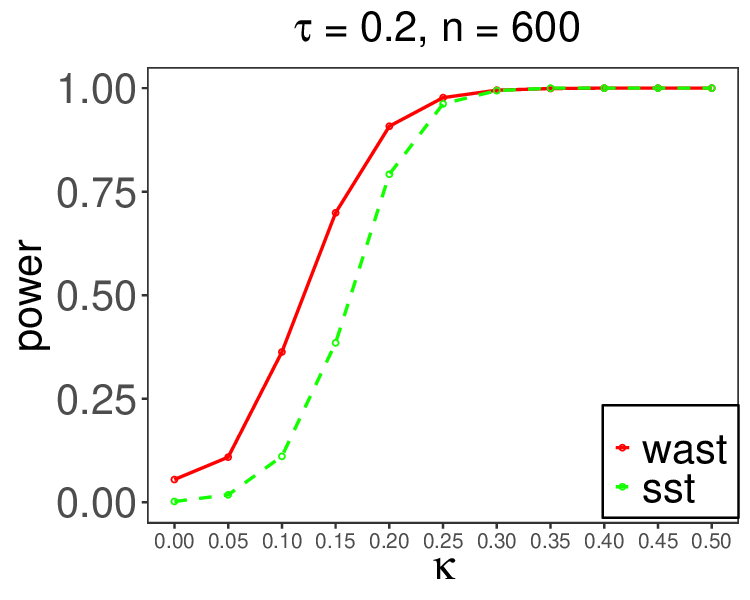}    \\
		\includegraphics[scale=0.33]{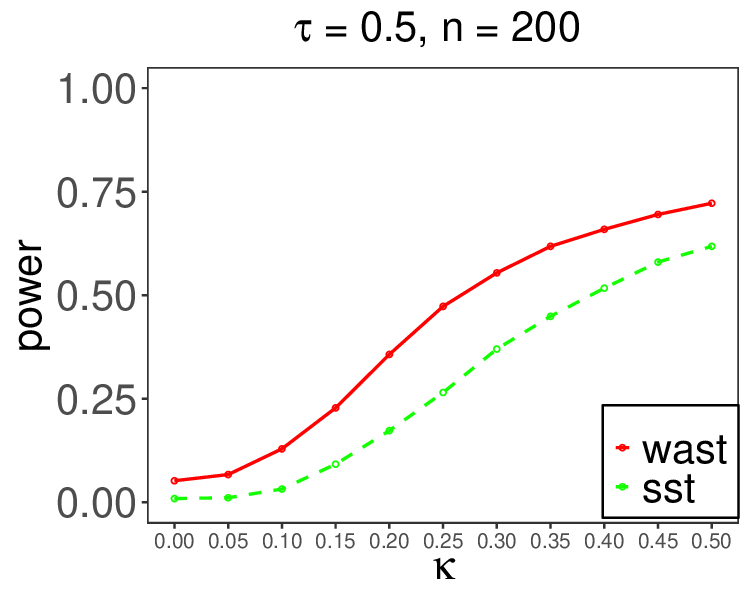}
		\includegraphics[scale=0.33]{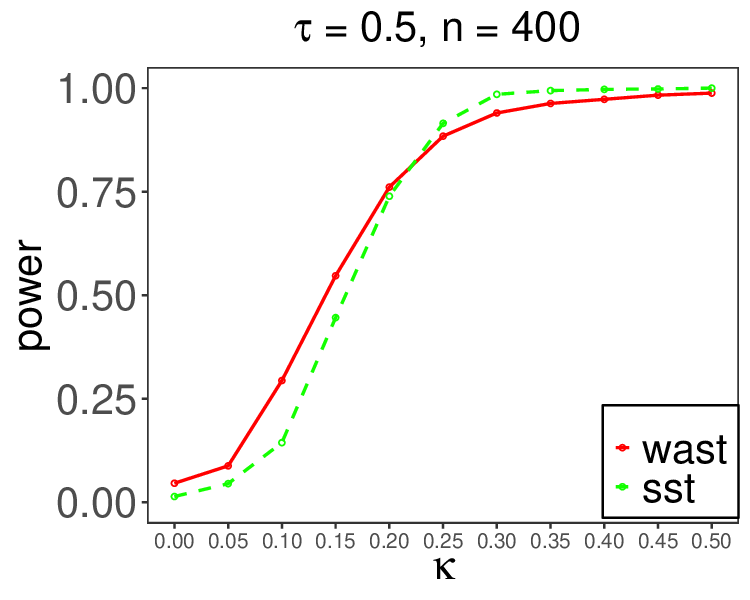}
		\includegraphics[scale=0.33]{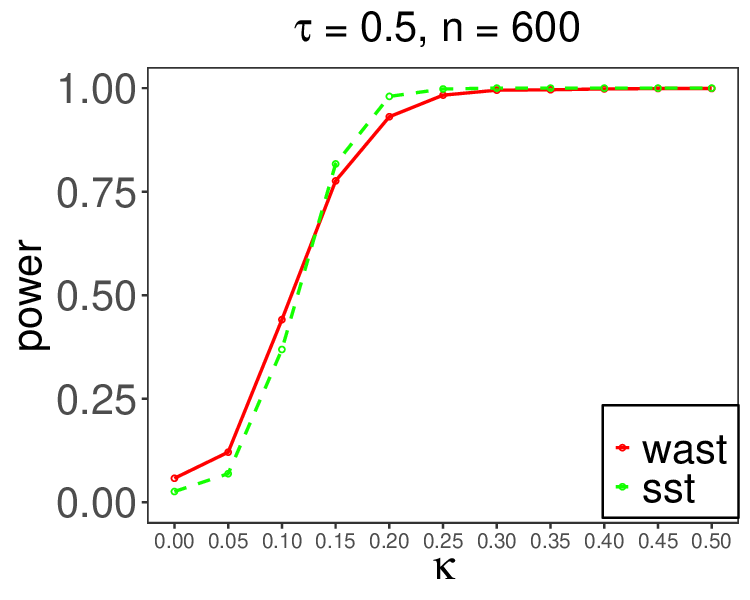}    \\
		\includegraphics[scale=0.33]{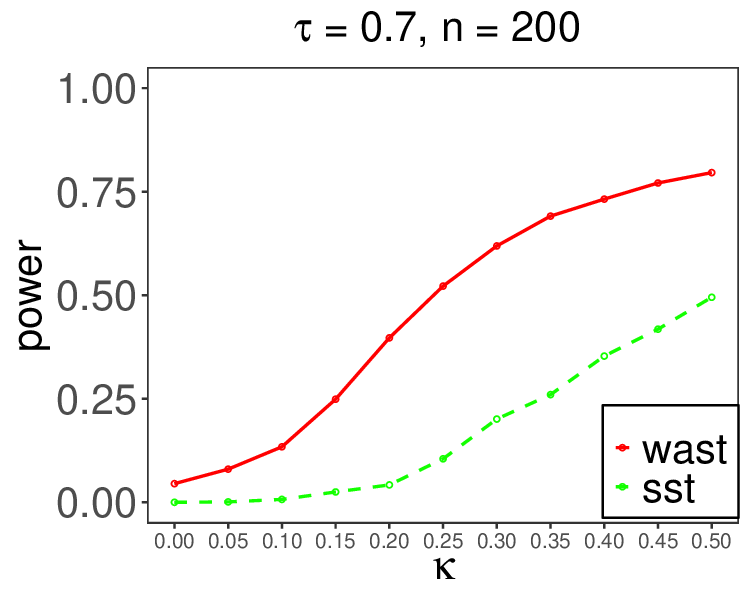}
		\includegraphics[scale=0.33]{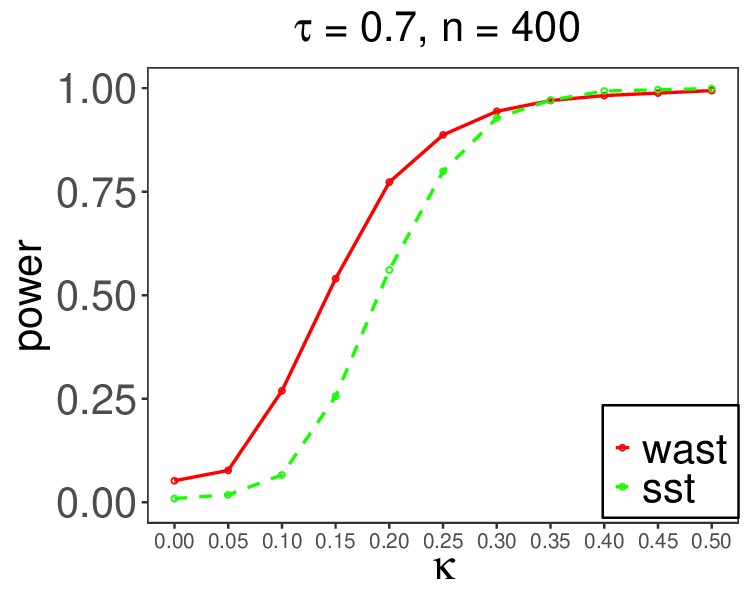}
		\includegraphics[scale=0.33]{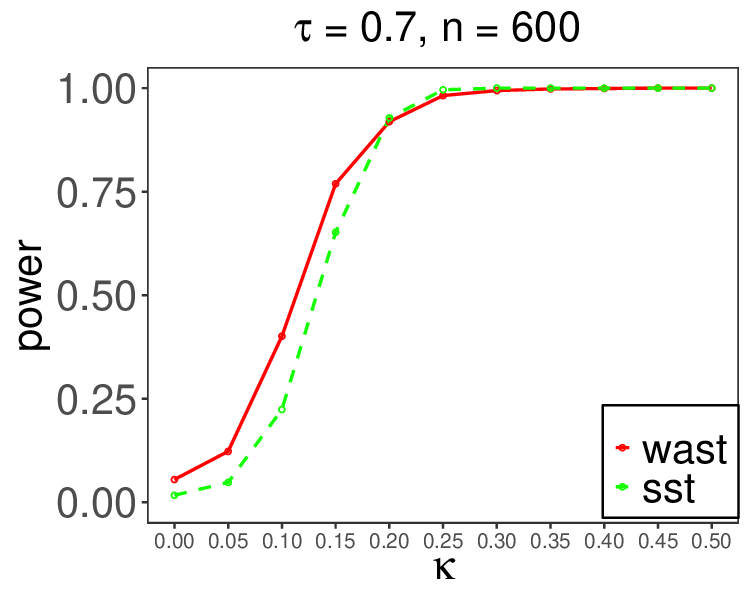}
		\caption{\it Powers of test statistic by the proposed WAST (red solid line) and SST (green dashed line) for $(p,q)=(10,10)$. From top to bottom, each row depicts the powers for probit model, semiparametric model, quantile regression with $\tau=0.2$, $\tau=0.5$ and $\tau=0.7$, respectively. Here the $\Gv$ $Z$ is generated from multivariate normal distribution with mean $\bzero$ and covariance $I$.}
		\label{fig_qr1010}
	\end{center}
\end{figure}

\begin{figure}[!ht]
	\begin{center}
		\includegraphics[scale=0.33]{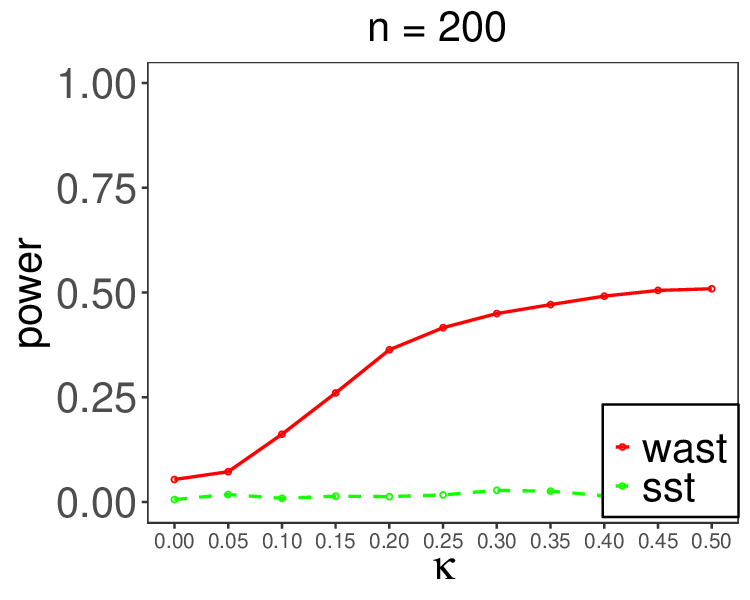}
		\includegraphics[scale=0.33]{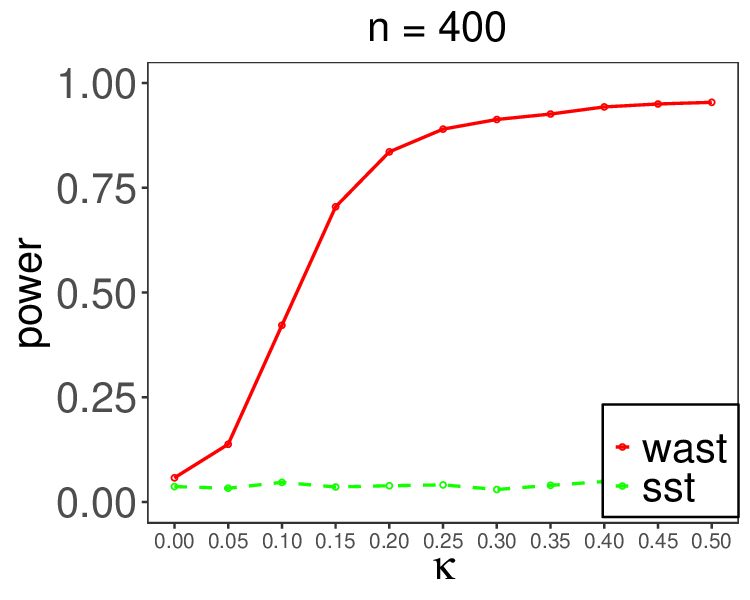}
		\includegraphics[scale=0.33]{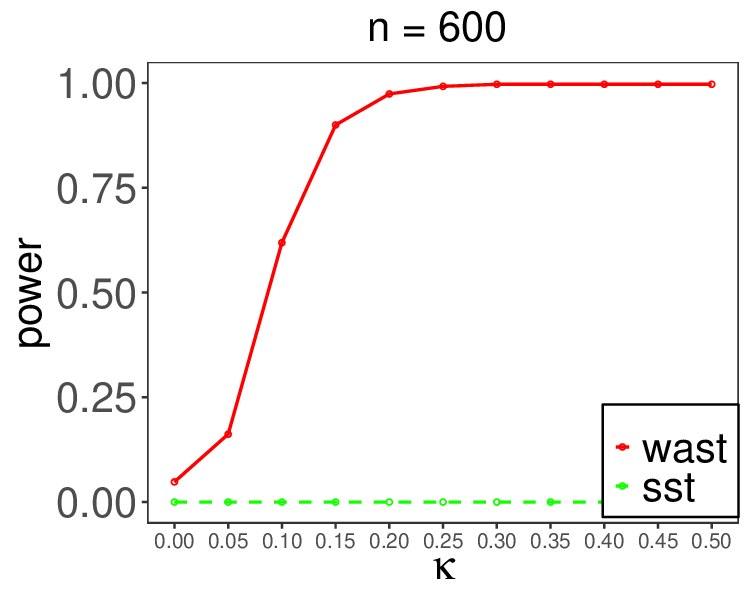}       \\
		\includegraphics[scale=0.33]{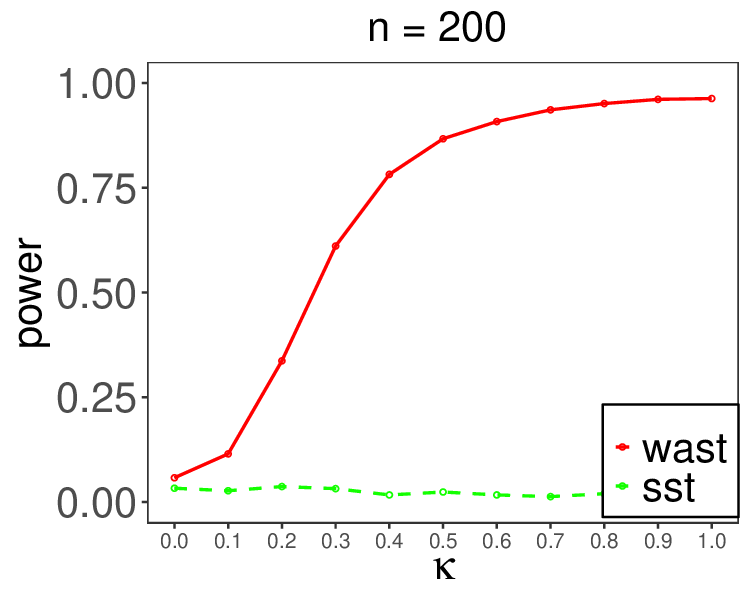}
		\includegraphics[scale=0.33]{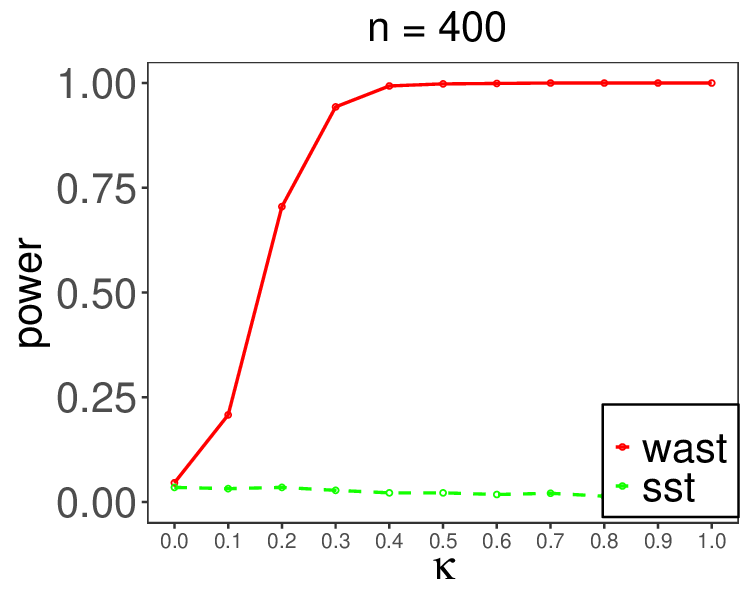}
		\includegraphics[scale=0.33]{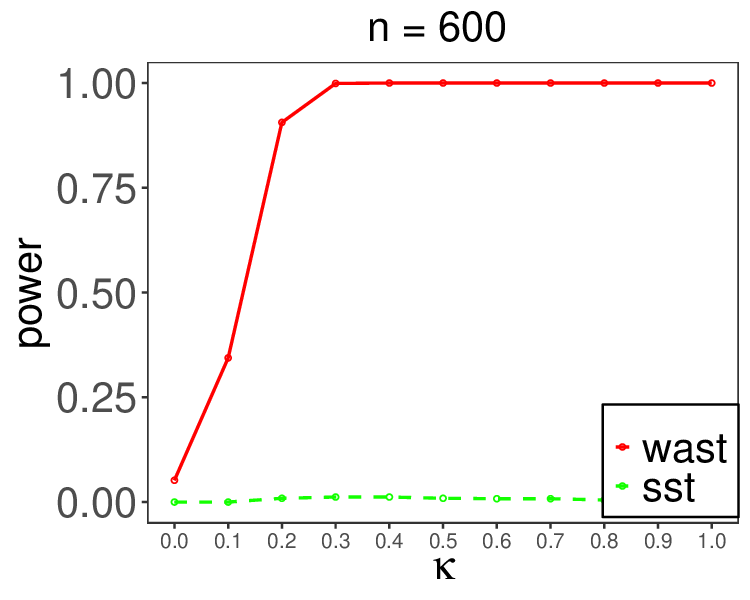}    \\
		\includegraphics[scale=0.33]{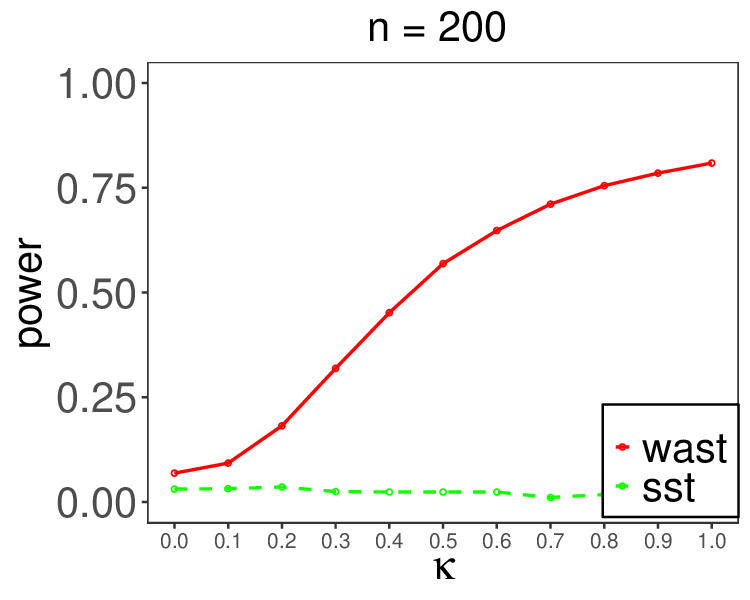}
		\includegraphics[scale=0.33]{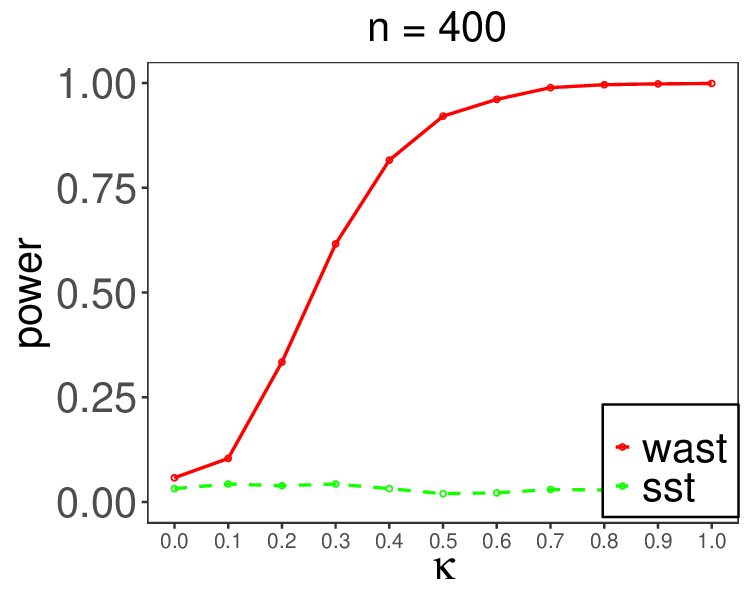}
		\includegraphics[scale=0.33]{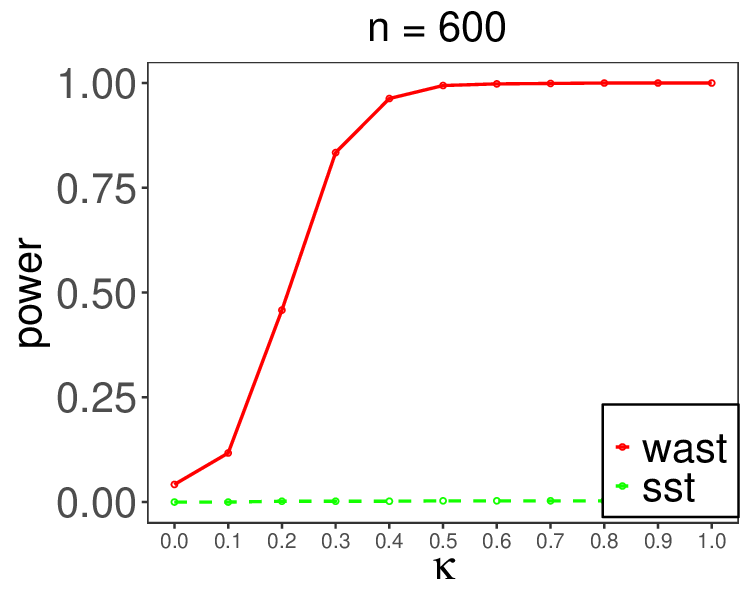}   \\
		\includegraphics[scale=0.33]{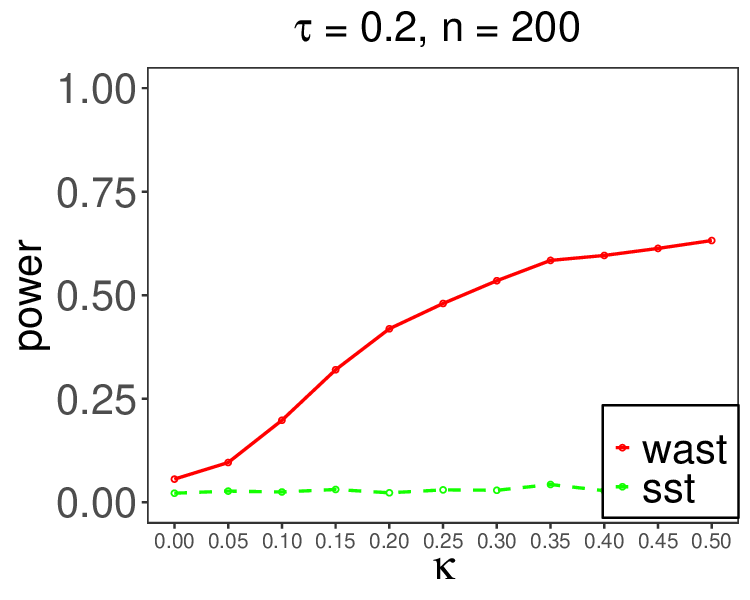}
		\includegraphics[scale=0.33]{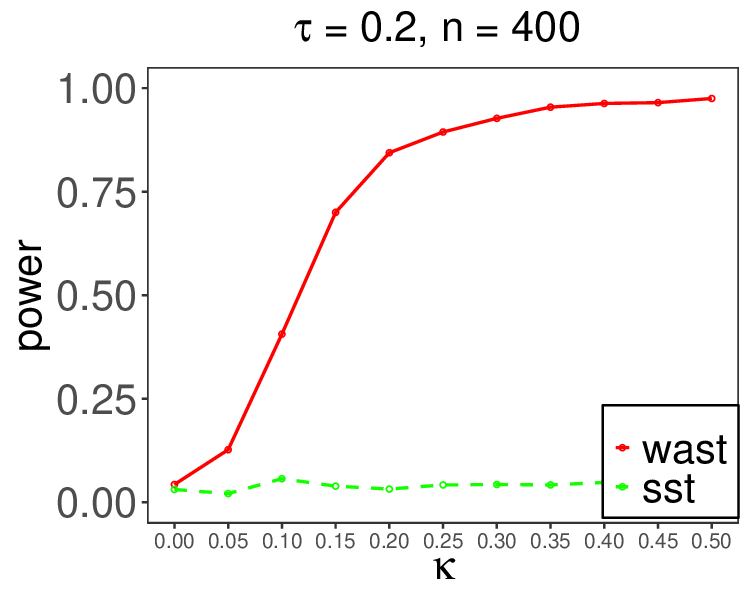}
		\includegraphics[scale=0.33]{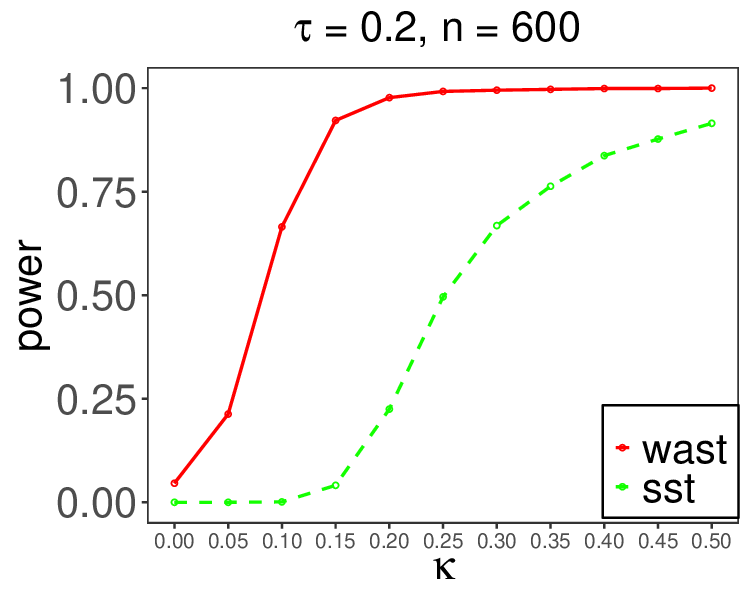}    \\
		\includegraphics[scale=0.33]{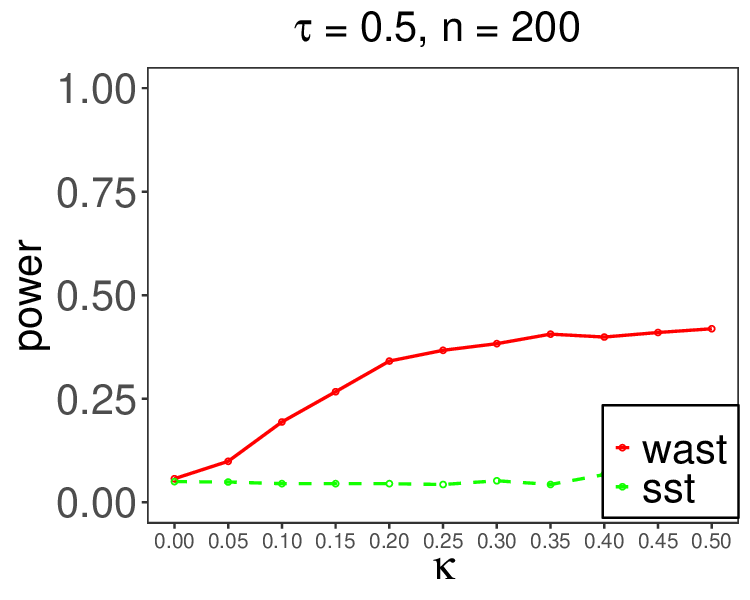}
		\includegraphics[scale=0.33]{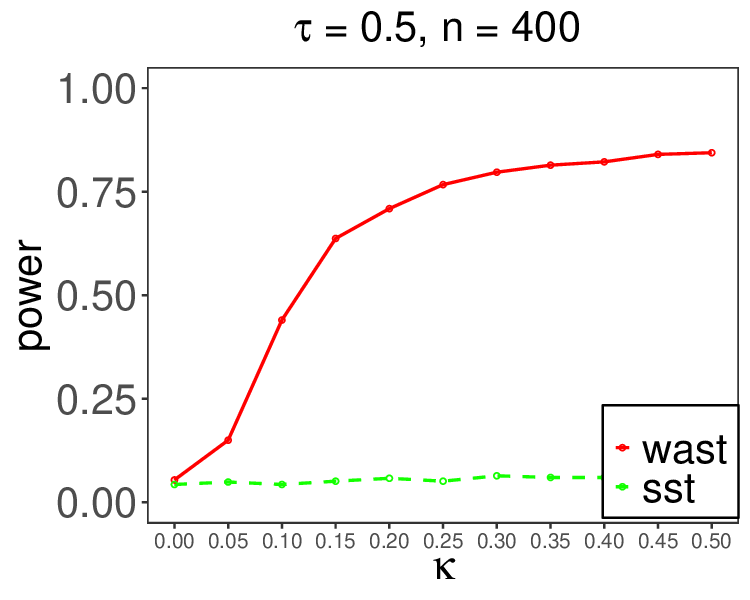}
		\includegraphics[scale=0.33]{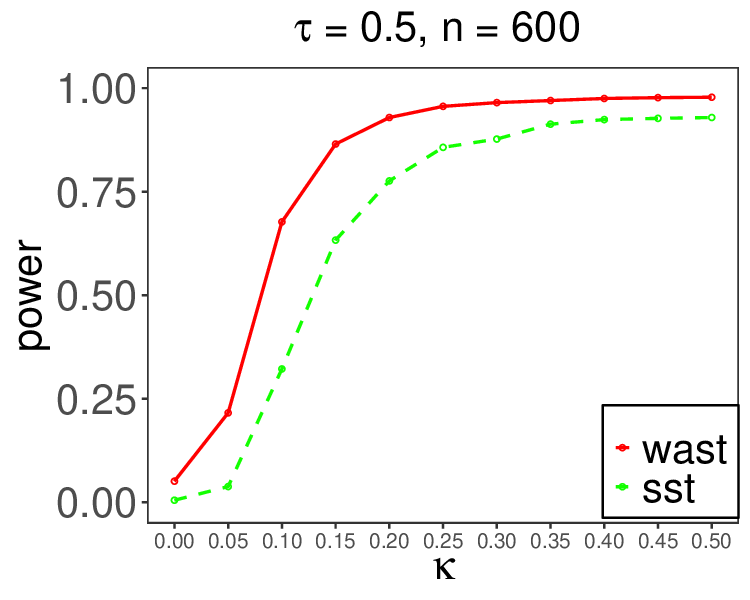}    \\
		\includegraphics[scale=0.33]{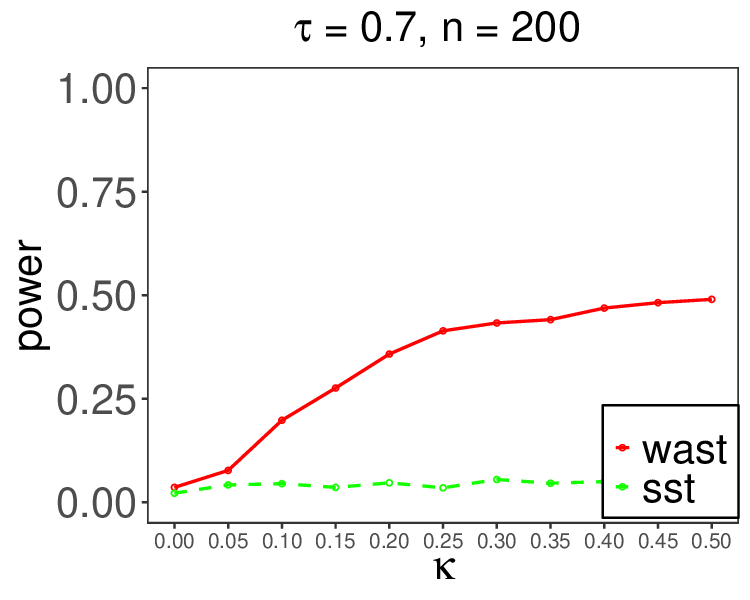}
		\includegraphics[scale=0.33]{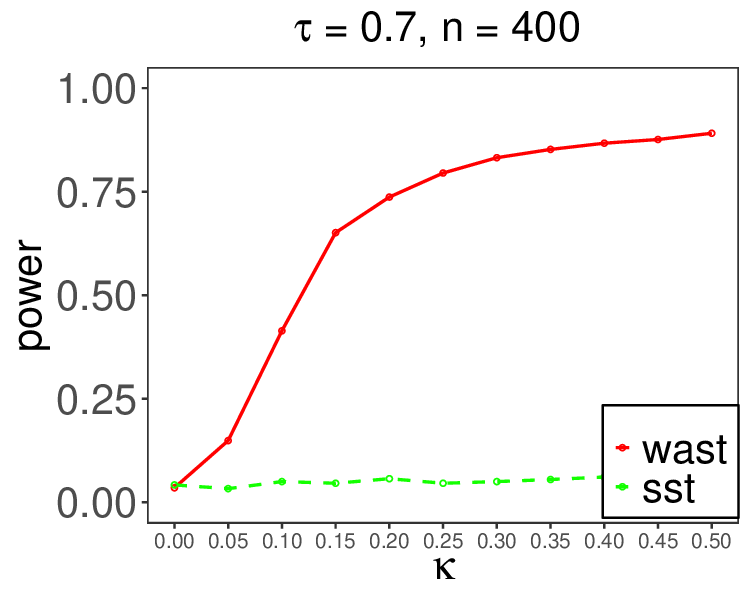}
		\includegraphics[scale=0.33]{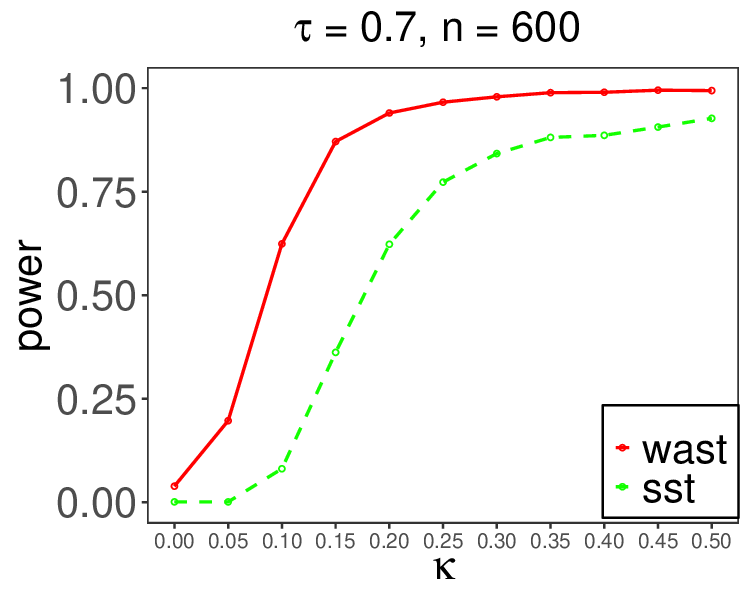}
		\caption{\it Powers of test statistic by the proposed WAST (red solid line) and SST (green dashed line) for $(p,q)=(50,20)$. From top to bottom, each row depicts the powers for probit model, semiparametric model, quantile regression with $\tau=0.2$, $\tau=0.5$ and $\tau=0.7$, respectively. Here the $\Gv$ $Z$ is generated from multivariate normal distribution with mean $\bzero$ and covariance $I$.}
		\label{fig_qr5020}
	\end{center}
\end{figure}

For all models with change plane analysis, we evaluated the power under a sequence of alternative models indexed by $\kappa$, that is $H_1^{\kappa}: \bbeta^{\kappa}=\kappa\bbeta^*$ with $\kappa = i/10$ for semiparametric model and $\kappa=i/20$ for others, $i=1,\cdots,10$, where $\bbeta^*=(1,\cdots,1)\trans$. We set sample size $n=(200, 400, 600)$, 1000 repetitions and 1000 bootstrap samples, and report in Figure \ref{fig_qr13}-\ref{fig_qr5020} the performance for both the WAST and SST. We calculate SST over $\{\btheta^{(k)}=(\theta^{(k)}_1,\cdots,\theta^{(k)}_q)\trans: k=1,\cdots,K\}$ with the number of threshold values $K=2000$. Let $\btheta^{(k)}_{-1}=\tilde{\btheta}^{(k)}_{-1}/\|\tilde{\btheta}^{(k)}_{-1}\|$, where $\btheta^{(k)}_{-1}=(\theta^{(k)}_2,\cdots,\theta^{(k)}_q)\trans$ and $\tilde{\btheta}^{(k)}_{-1}=(\tilde{\theta}^{(k)}_2,\cdots,\tilde{\theta}^{(k)}_q)\trans$, and $\tilde{\btheta}^{(k)}_{-1}$ is drawn independently form $(r-1)$-variate standard normal distribution.  For each $\theta_1^{(k)}$, $k=1,\cdots,K$, we selected it by equal grid search in the range from the lower 10th percentile to upper 10th percentile of the data points of $\{\theta^{(k)}_2Z_{2i}+\cdots+\theta^{(k)}_qZ_{qi}\}_{i=1}^n$, which is same as that in \cite{2020Threshold}. Here we consider four combinations of $(p,q)=(1, 3), (5, 5), (10, 10), (50, 20)$.

Type \uppercase\expandafter{\romannumeral1} errors ($\kappa=0$) with sample size $n=(200, 400, 600)$ are list in Table \ref{table_size1}. We can see from Table \ref{table_size1} that the size of the proposed WAST are close to the nominal significance level $0.05$, but for most scenarios the size of the SST are much smaller than 0.05.
Figure \ref{fig_qr13}-\ref{fig_qr55} indicate that powers become greater as sample size $n$ increases, which are verified the asymptotic theory. The proposed WAST has comparable power with the SST for the semiparametric model, but the size of the SST is much less than the nominal level 0.05 in Table \ref{table_size1}. SST has wrong size, and is less powerful than the WAST for probit and quantile regression models in four combinations of $(p,q)$. Since the size is much less than the nominal level 0.05 when dimension of $\bX$ not smaller than 5 in Table \ref{table_size1}, we only depict in Figure \ref{fig_qr1010} the power curves for the WAST. Figure \ref{fig_qr1010} shows that the power increases fast as the sample size becomes large as expected from the theory, and the power increases rapidly in all models and two combinations $(p,q)=(10, 10)$ and $(50, 20)$ as the alternative is farther away from the null.

\begin{table}[htp!]
	\def~{\hphantom{0}}
    \tiny
	\caption{Type \uppercase\expandafter{\romannumeral1} errors of the proposed WAST and SST based on resampling for probit regression model (ProbitRE), quantile regression (QuantRE) and semiparametric model (SPMoldel).
}
	\resizebox{\textwidth}{!}{
    \begin{threeparttable}
		\begin{tabular}{llcccccccc}
			\hline
			\multirow{2}{*}{Model}&\multirow{2}{*}{$(p,q)$}
			&\multicolumn{2}{c}{ $n=200$} && \multicolumn{2}{c}{ $n=400$} && \multicolumn{2}{c}{ $n=600$} \\
			\cline{3-4} \cline{6-7} \cline{9-10}
			&&WAST& SST && WAST& SST && WAST& SST \\
			\cline{3-10}
			ProbitRE &$(1,3)$         & 0.034 & 0.003 && 0.055 & 0.012 && 0.046 & 0.023   \\
			&$(5,5)$                  & 0.044 & 0.000 && 0.031 & 0.000 && 0.062 & 0.001   \\
			& $(10,10)$               & 0.049 & 0.000 && 0.058 & 0.000 && 0.045 & 0.000    \\
			&$(50,20)$                & 0.054 & 0.006 && 0.058 & 0.037 && 0.048 & 0.000    \\
			[1 ex]
			SPModel &$(1,3)$          & 0.053 & 0.051 && 0.063 & 0.045 && 0.046 & 0.048   \\
			(B1+P1)&$(5,5)$           & 0.056 & 0.008 && 0.049 & 0.022 && 0.046 & 0.012   \\
			& $(10,10)$               & 0.056 & 0.001 && 0.047 & 0.007 && 0.052 & 0.018 \\
			&$(50,20)$                & 0.058 & 0.033 && 0.046 & 0.035 && 0.052 & 0.000 \\
			[1 ex]
			SPModel &$(1,3)$          & 0.057 & 0.044 && 0.054 & 0.048 && 0.046 & 0.053    \\
			(B2+P2)&$(5,5)$           & 0.050 & 0.007 && 0.059 & 0.013 && 0.032 & 0.013    \\
			& $(10,10)$               & 0.053 & 0.002 && 0.048 & 0.006 && 0.038 & 0.008 \\
			&$(50,20)$                & 0.058 & 0.031 && 0.058 & 0.032 && 0.042 & 0.000 \\
			[1 ex]
			QuantRE &$(1,3)$          & 0.044 & 0.015 && 0.057 & 0.027 && 0.043 & 0.030   \\
			($\tau=0.2$)&$(5,5)$      & 0.051 & 0.000 && 0.059 & 0.016 && 0.058 & 0.020   \\
			& $(10,10)$               & 0.041 & 0.000 && 0.048 & 0.000 && 0.055 & 0.002 \\
			&$(50,20)$                & 0.056 & 0.022 && 0.043 & 0.031 && 0.046 & 0.000 \\
			[1 ex]
			QuantRE &$(1,3)$          & 0.060 & 0.039 && 0.058 & 0.051 && 0.053 & 0.034   \\
			($\tau=0.5$)&$(5,5)$      & 0.049 & 0.005 && 0.046 & 0.024 && 0.047 & 0.033   \\
			& $(10,10)$               & 0.052 & 0.009 && 0.046 & 0.014 && 0.058 & 0.026 \\
			&$(50,20)$                & 0.057 & 0.050 && 0.054 & 0.043 && 0.051 & 0.005 \\
			[1 ex]
			QuantRE &$(1,3)$          & 0.044 & 0.020 && 0.050 & 0.025 && 0.053 & 0.037   \\
			($\tau=0.7$)&$(5,5)$      & 0.048 & 0.005 && 0.055 & 0.021 && 0.045 & 0.024   \\
			& $(10,10)$               & 0.045 & 0.000 && 0.052 & 0.009 && 0.055 & 0.017 \\
			&$(50,20)$                & 0.036 & 0.022 && 0.035 & 0.042 && 0.039 & 0.001 \\
			\hline
		\end{tabular}
\begin{tablenotes}
\item The nominal significant level is 0.05. Here the $\Gv$ $\bZ$ is generated from multivariate normal distribution with mean $\bzero$ and covariance $\bI$.
\end{tablenotes}
\end{threeparttable}
	}
	\label{table_size1}
\end{table}

\subsection{Change plane analysis for quantile, probit and semiparametric models with \texorpdfstring{$Z$}{} from \texorpdfstring{$\mathcal{N}(0,25)$ }{}}\label{simulation_cpz}
We consider probit regression, quantile regression and semiparametric models with same settings as Section \ref{simulation_cp} in the main paper but different $\Gv$ $\bZ$ and error's distribution of quantile regression model. We generate $Z_{ji}$ independently from normal distribution with mean zero and standard deviation 5, and set $\bZ_i=(1,Z_{1i},\cdots,Z_{(q-1)i})\trans$, where $j=1,\cdots,q-1$. For quantile regression, we generate error $\eps_i$ independently from $t_3$ distribution with degrees freedom 3.

\begin{figure}[!ht]
	\begin{center}
		\includegraphics[scale=0.33]{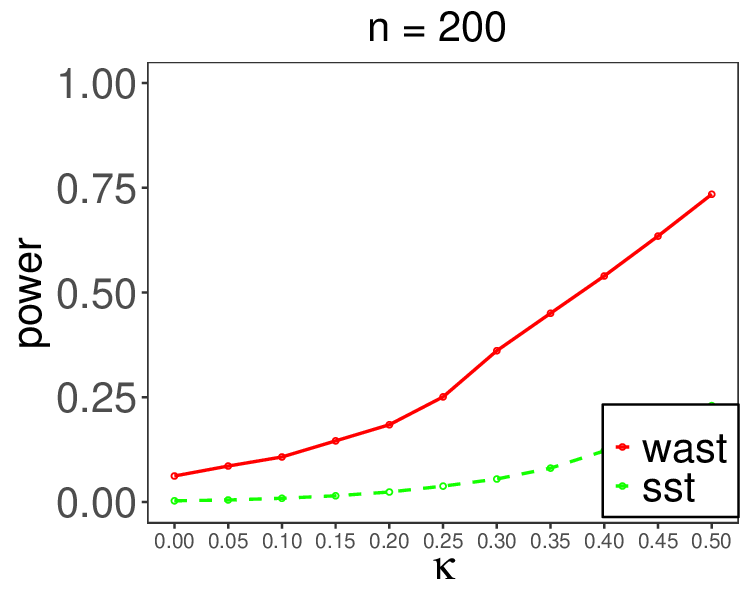}
		\includegraphics[scale=0.33]{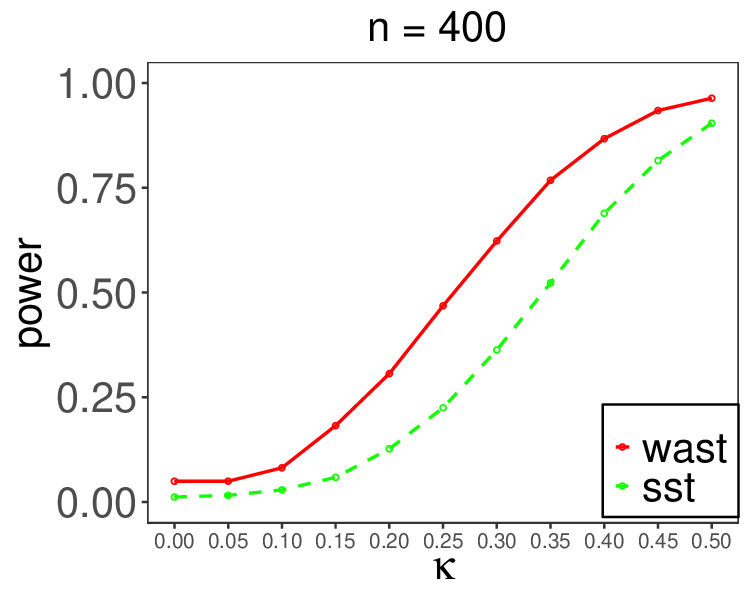}
		\includegraphics[scale=0.33]{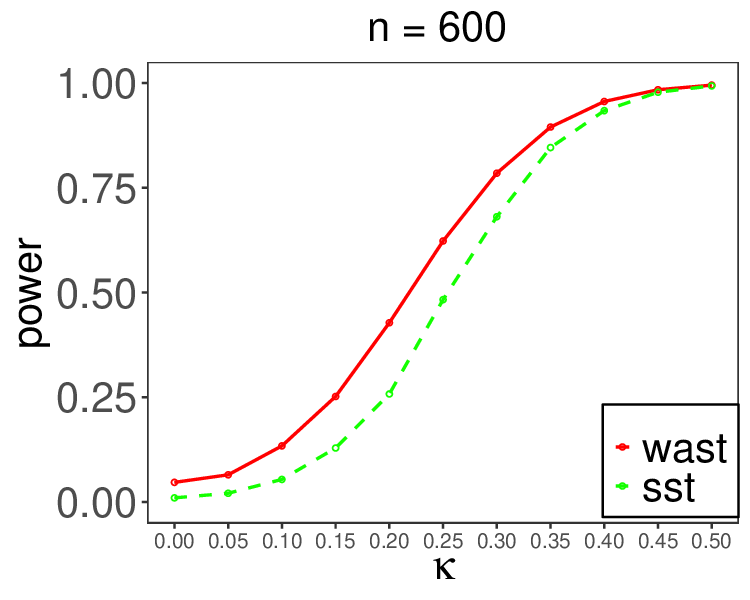}      \\
		\includegraphics[scale=0.33]{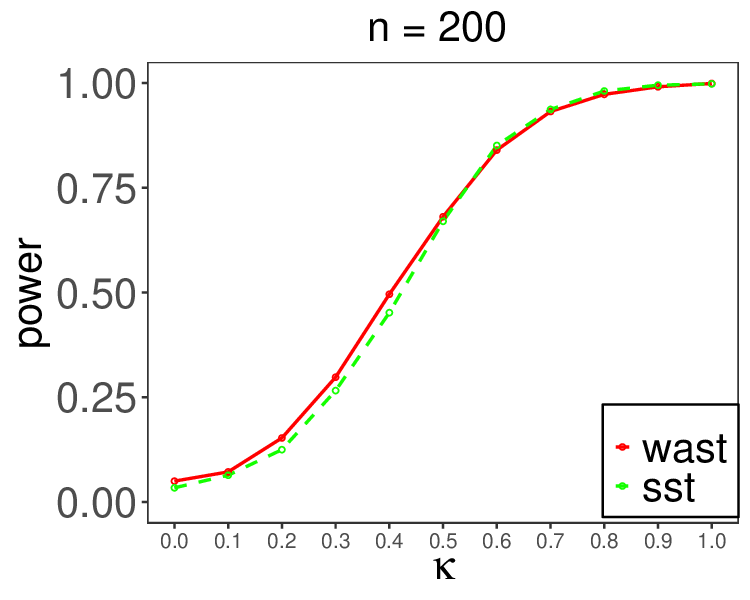}
		\includegraphics[scale=0.33]{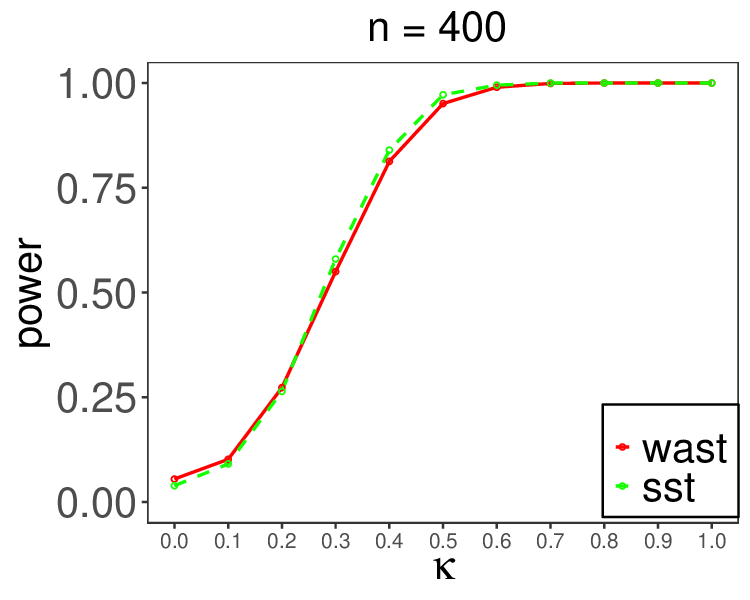}
		\includegraphics[scale=0.33]{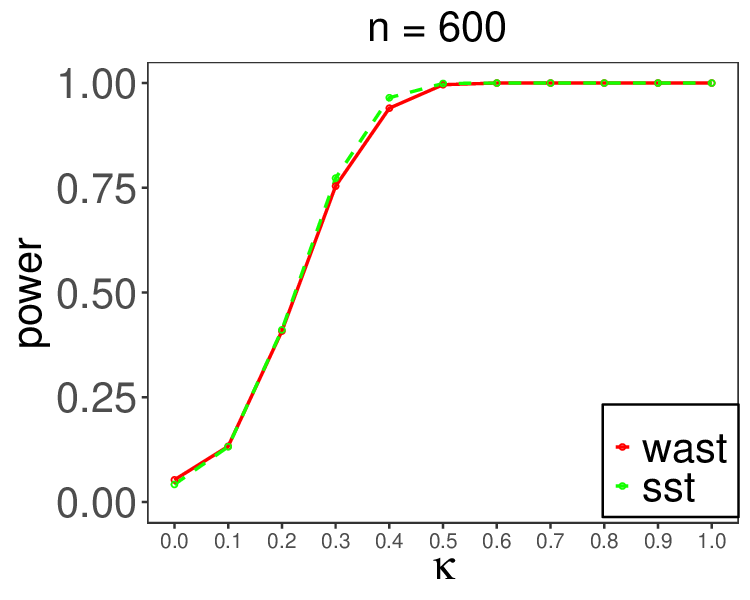}   \\
		\includegraphics[scale=0.33]{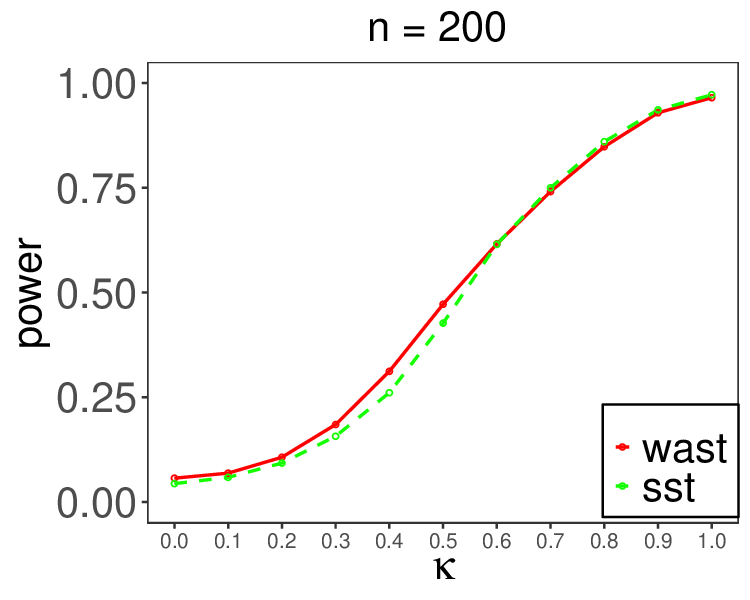}
		\includegraphics[scale=0.33]{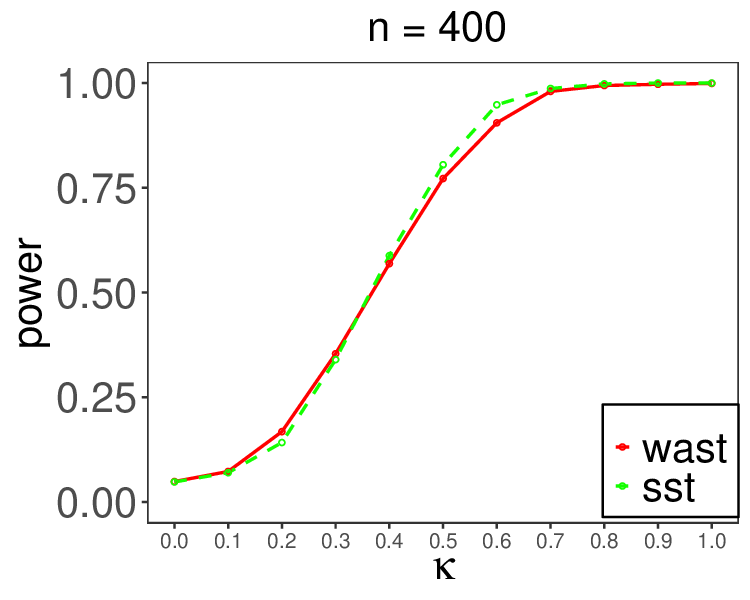}
		\includegraphics[scale=0.33]{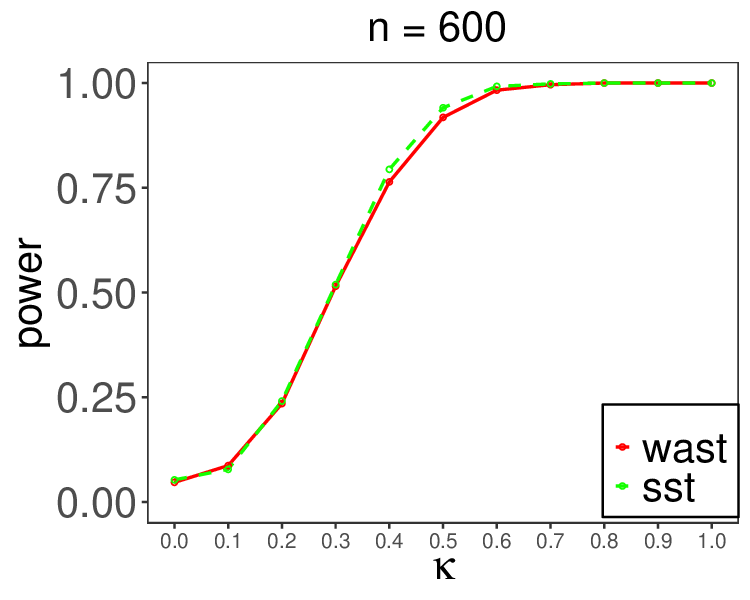}  \\
		\includegraphics[scale=0.33]{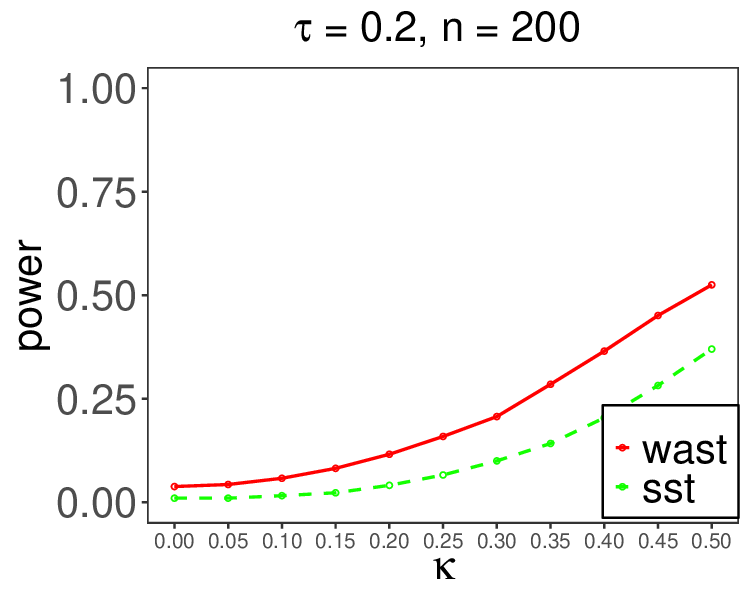}
		\includegraphics[scale=0.33]{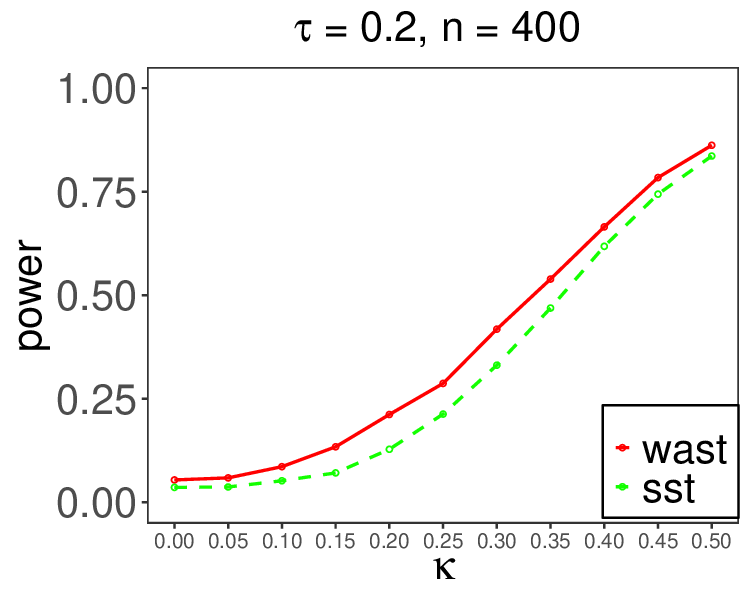}
		\includegraphics[scale=0.33]{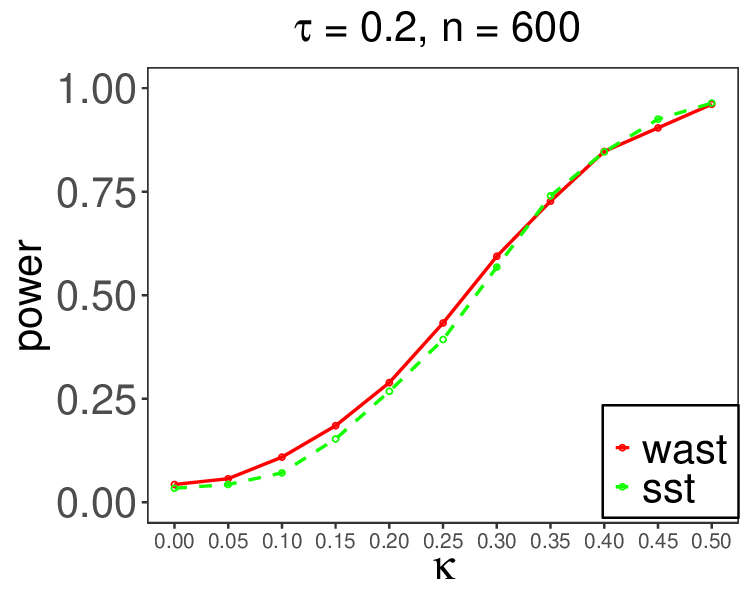}   \\
		\includegraphics[scale=0.33]{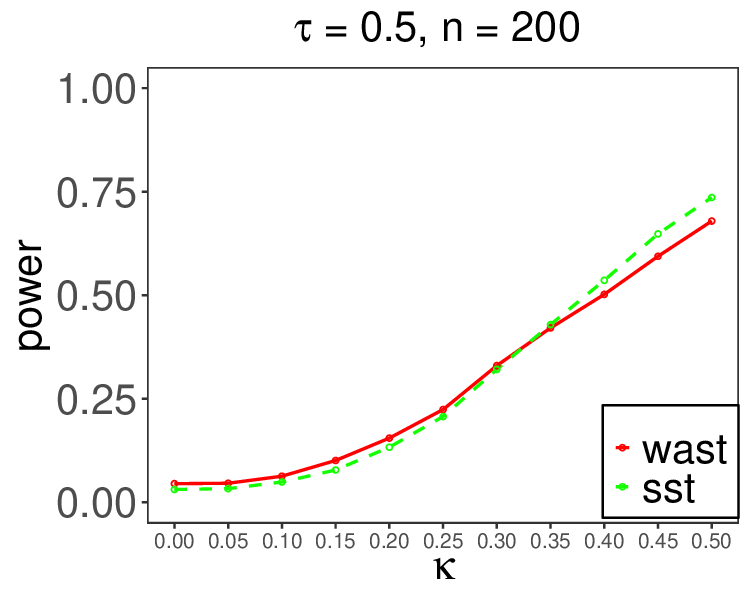}
		\includegraphics[scale=0.33]{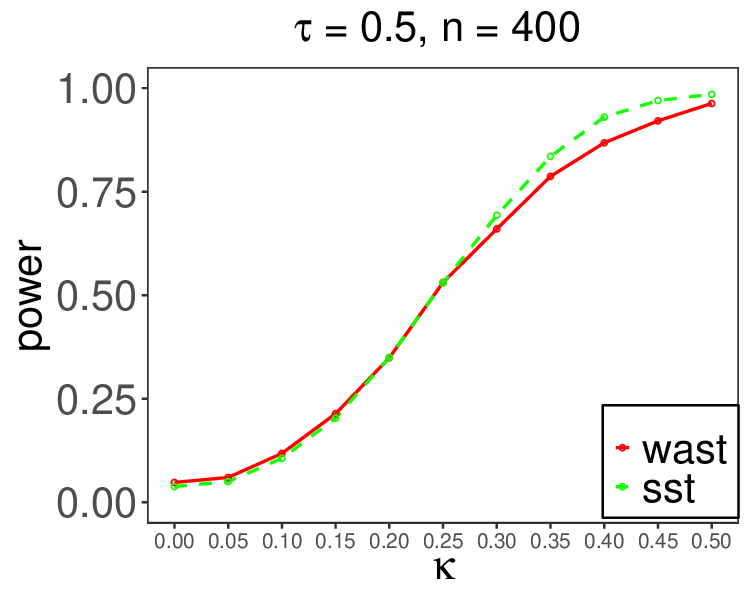}
		\includegraphics[scale=0.33]{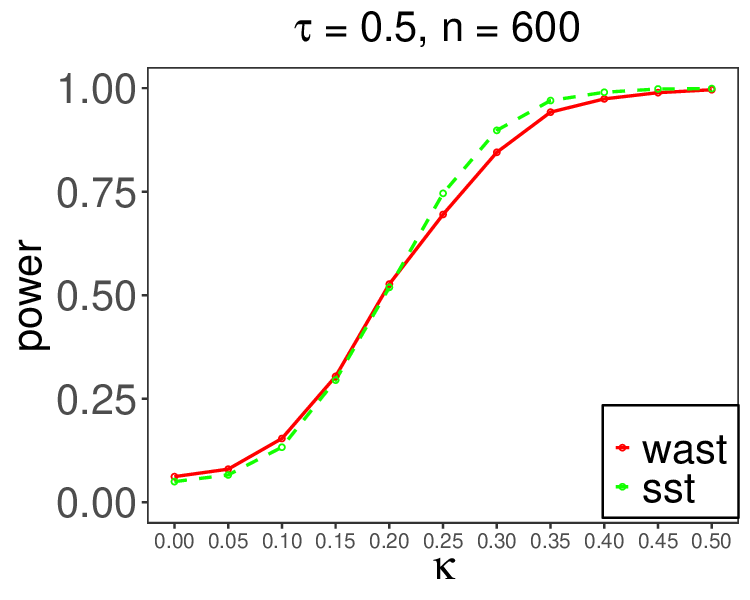}   \\
		\includegraphics[scale=0.33]{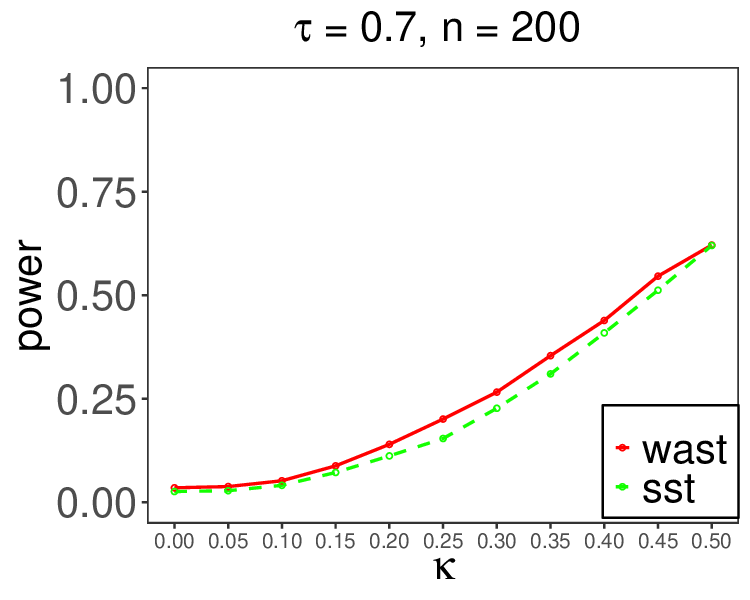}
		\includegraphics[scale=0.33]{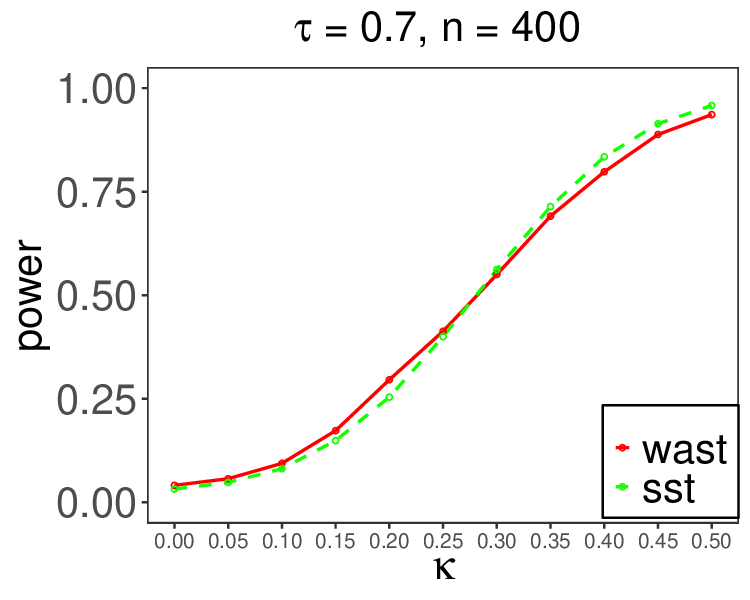}
		\includegraphics[scale=0.33]{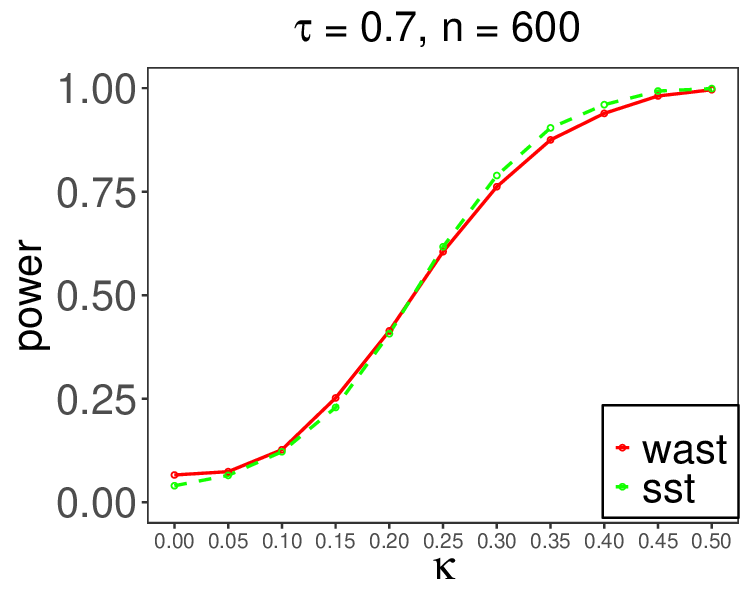}
		\caption{\it Powers of test statistic by the proposed WAST (red solid line) and SST (green dashed line) for $(p,q)=(1,3)$. From top to bottom, each row depicts the powers for probit model, semiparametric model, quantile regression with $\tau=0.2$, $\tau=0.5$ and $\tau=0.7$, respectively. Here the $\Gv$ $Z$ is generated from multivariate normal distribution with mean $\bzero$ and covariance $25I$.}
		\label{fig_qr13_1}
	\end{center}
\end{figure}

\begin{figure}[!ht]
	\begin{center}
		\includegraphics[scale=0.33]{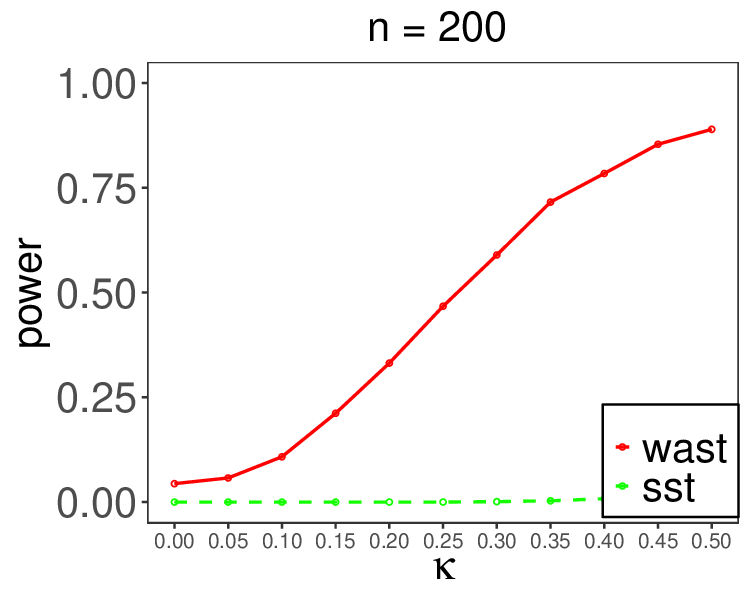}
		\includegraphics[scale=0.33]{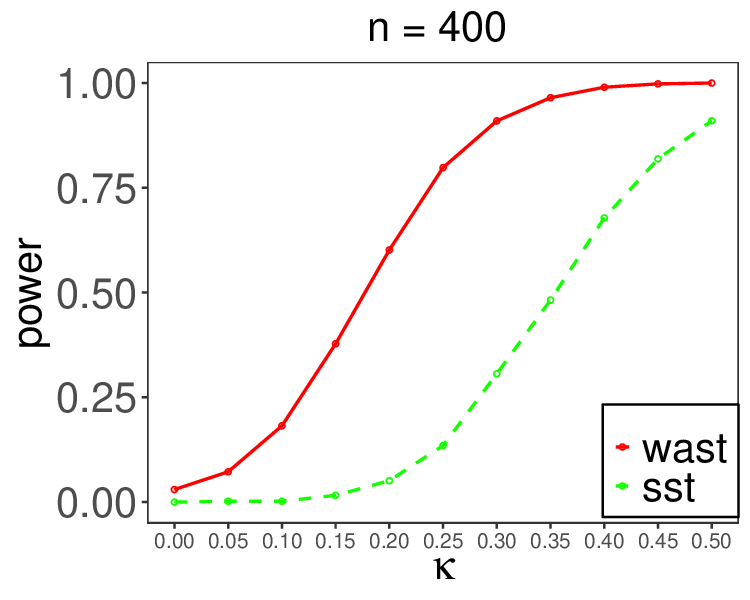}
		\includegraphics[scale=0.33]{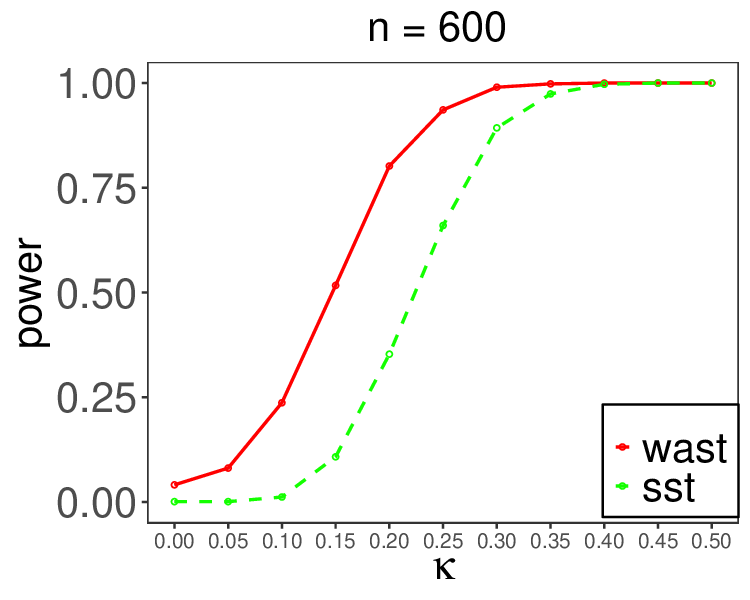}      \\
		\includegraphics[scale=0.33]{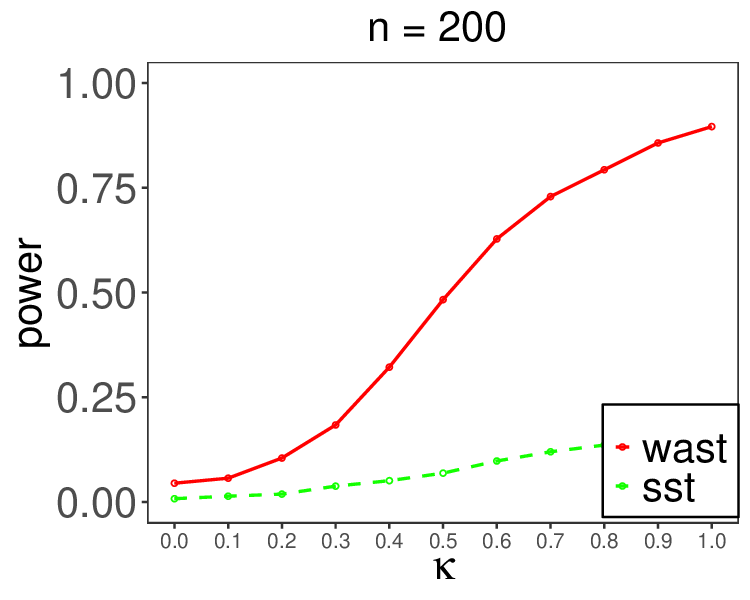}
		\includegraphics[scale=0.33]{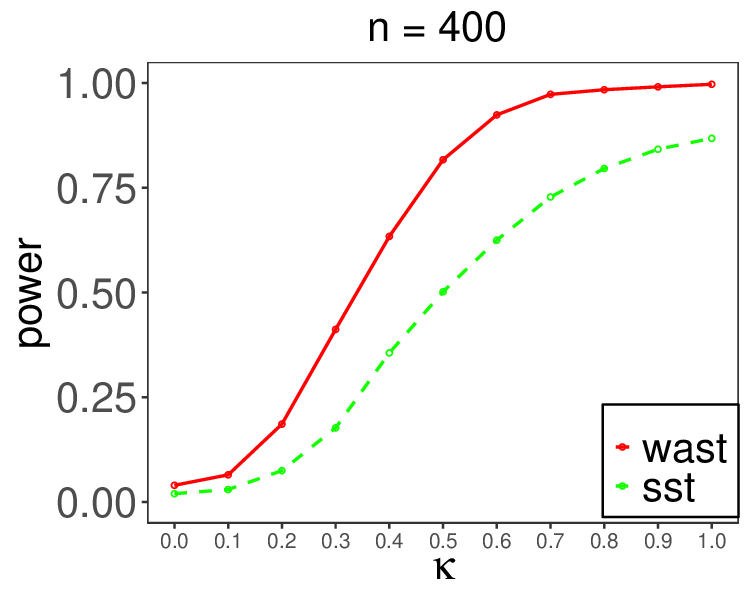}
		\includegraphics[scale=0.33]{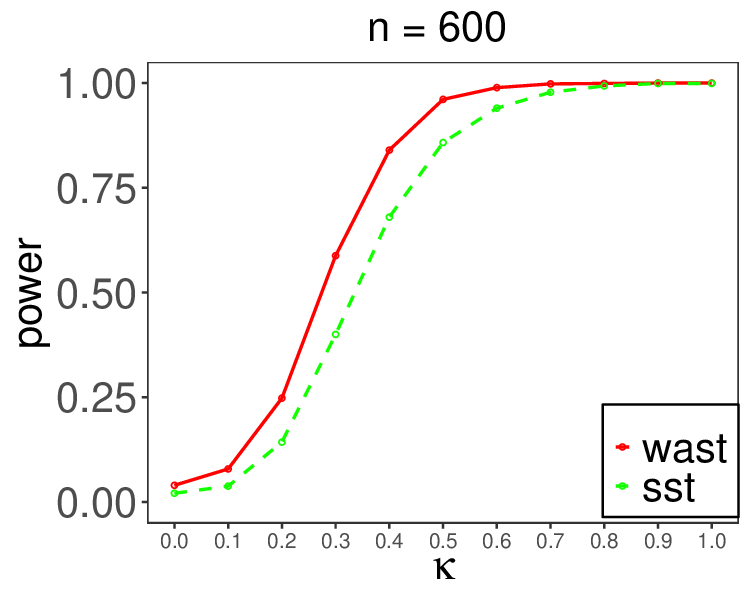}   \\
		\includegraphics[scale=0.33]{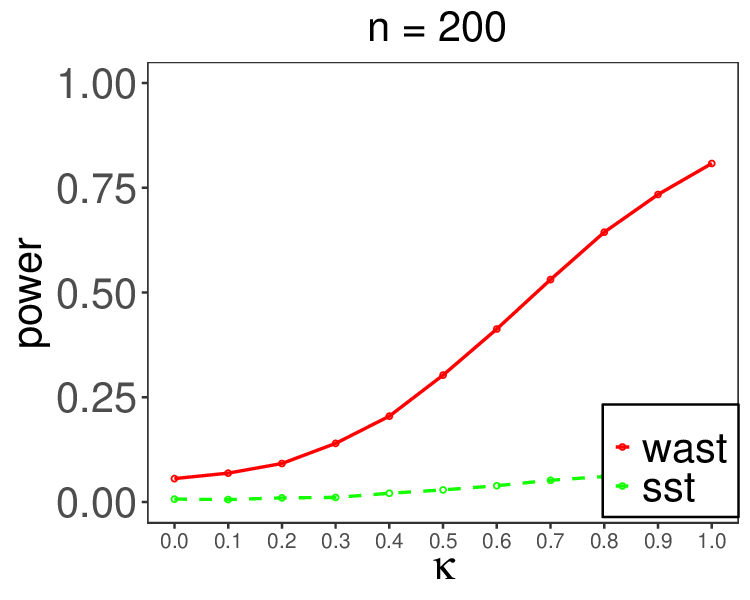}
		\includegraphics[scale=0.33]{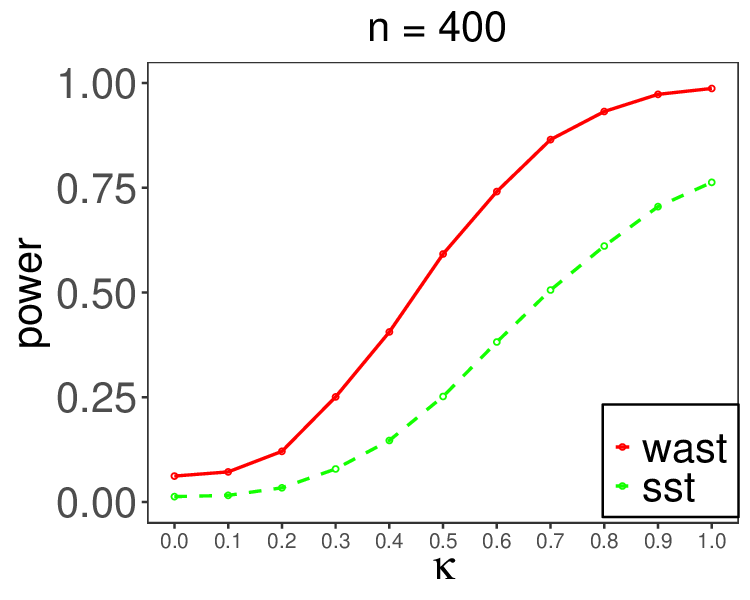}
		\includegraphics[scale=0.33]{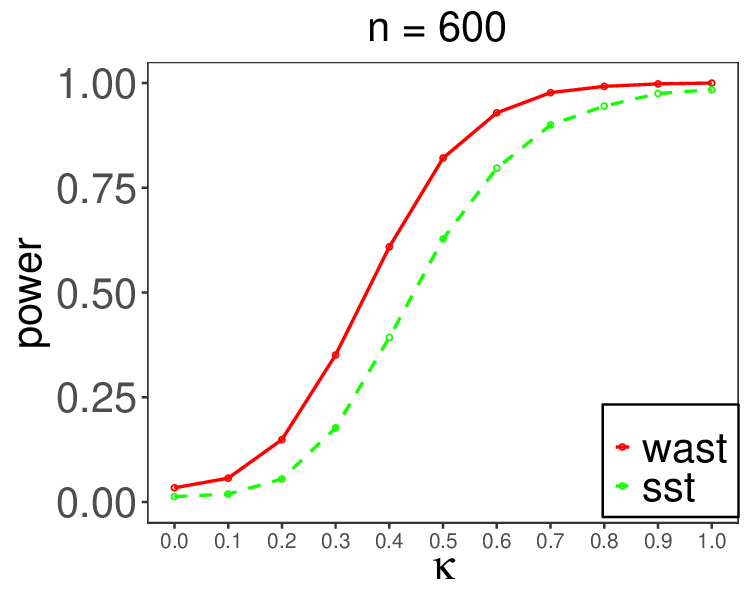}  \\
		\includegraphics[scale=0.33]{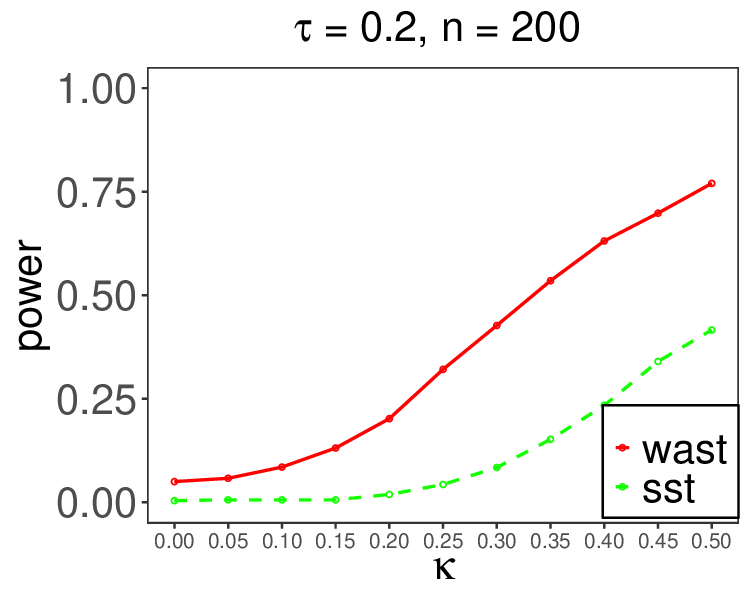}
		\includegraphics[scale=0.33]{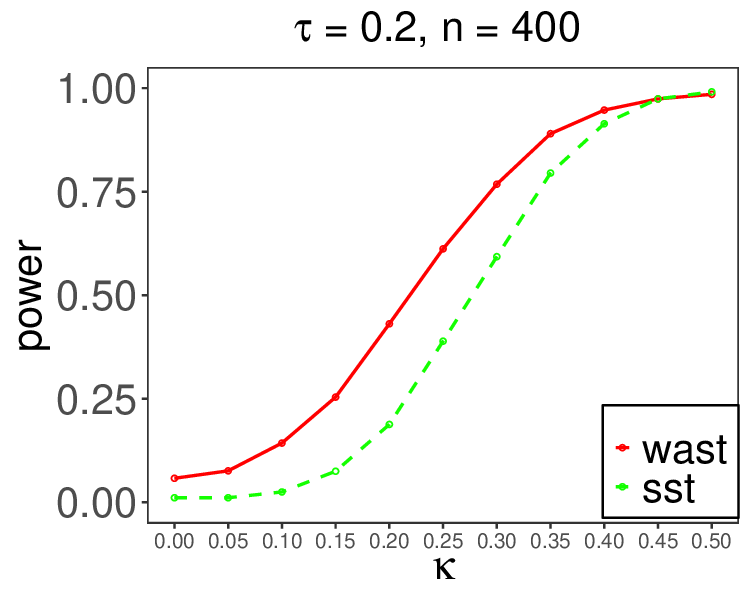}
		\includegraphics[scale=0.33]{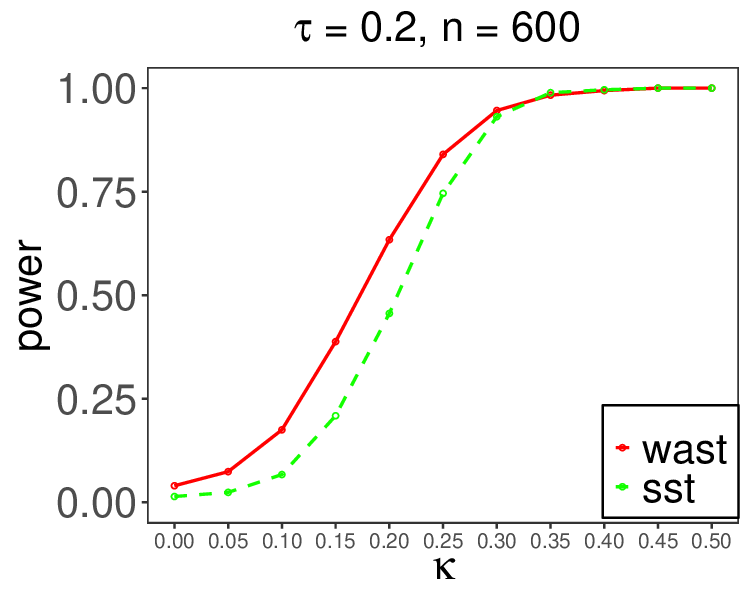}   \\
		\includegraphics[scale=0.33]{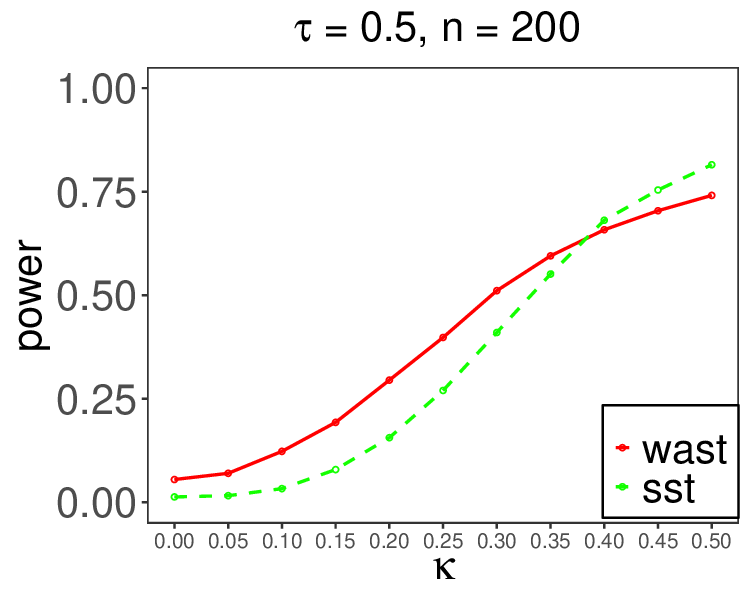}
		\includegraphics[scale=0.33]{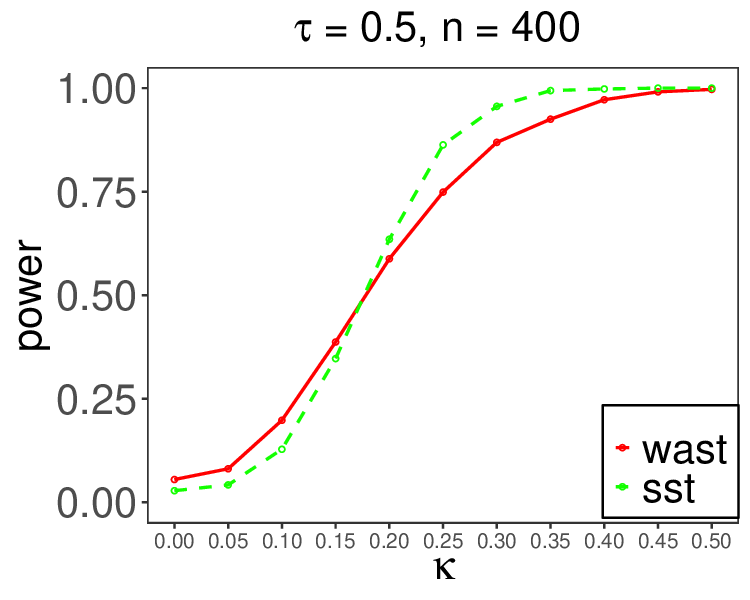}
		\includegraphics[scale=0.33]{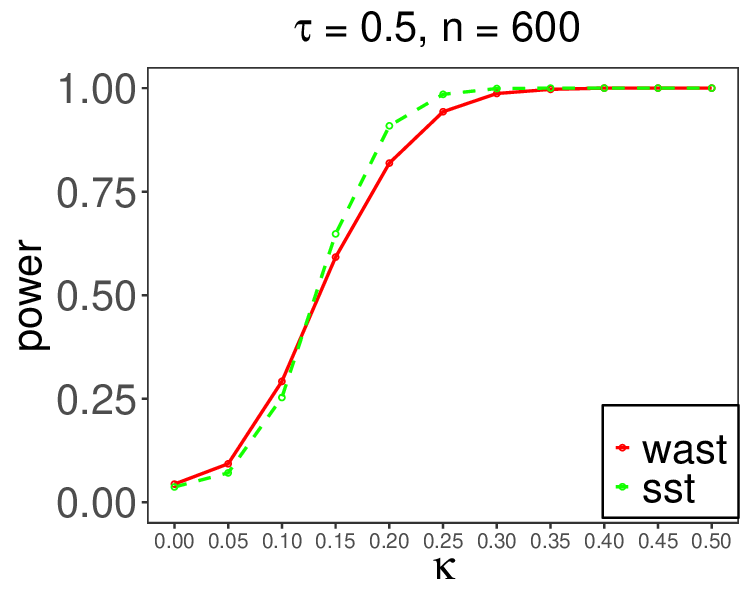}   \\
		\includegraphics[scale=0.33]{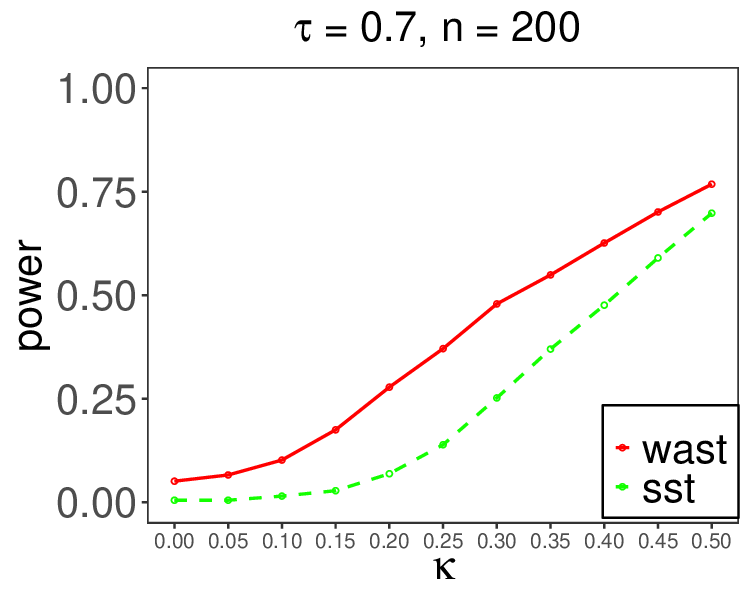}
		\includegraphics[scale=0.33]{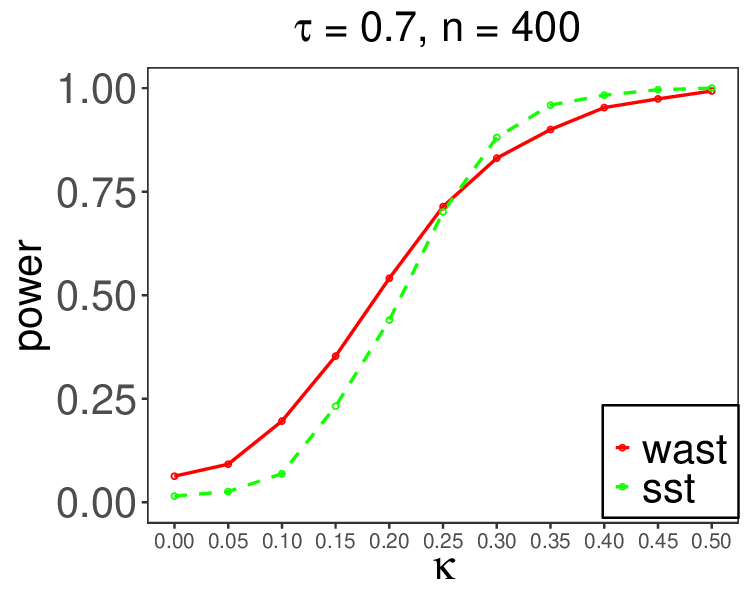}
		\includegraphics[scale=0.33]{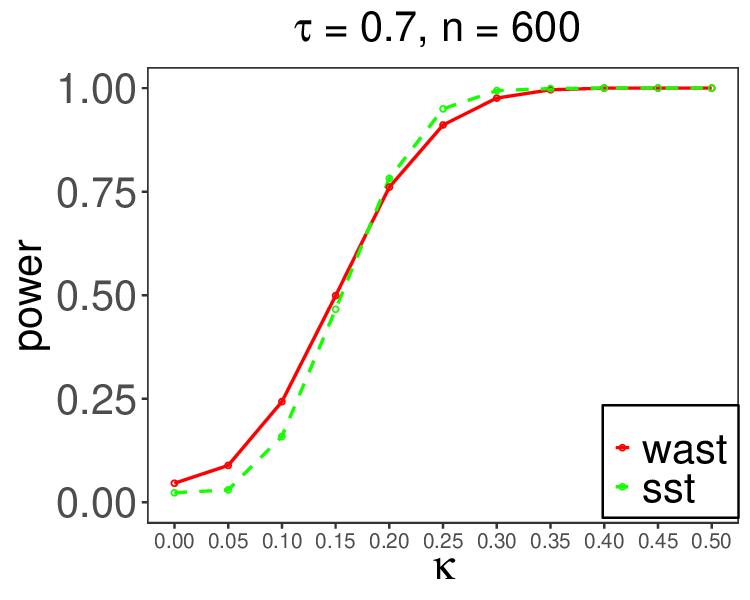}
		\caption{\it Powers of test statistic by the proposed WAST (red solid line) and SST (green dashed line) for $(p,q)=(5,5)$. From top to bottom, each row depicts the powers for probit model, semiparametric model, quantile regression with $\tau=0.2$, $\tau=0.5$ and $\tau=0.7$, respectively. Here the $\Gv$ $Z$ is generated from multivariate normal distribution with mean $\bzero$ and covariance $25I$.}
		\label{fig_qr55_1}
	\end{center}
\end{figure}

\begin{figure}[!ht]
	\begin{center}
		\includegraphics[scale=0.33]{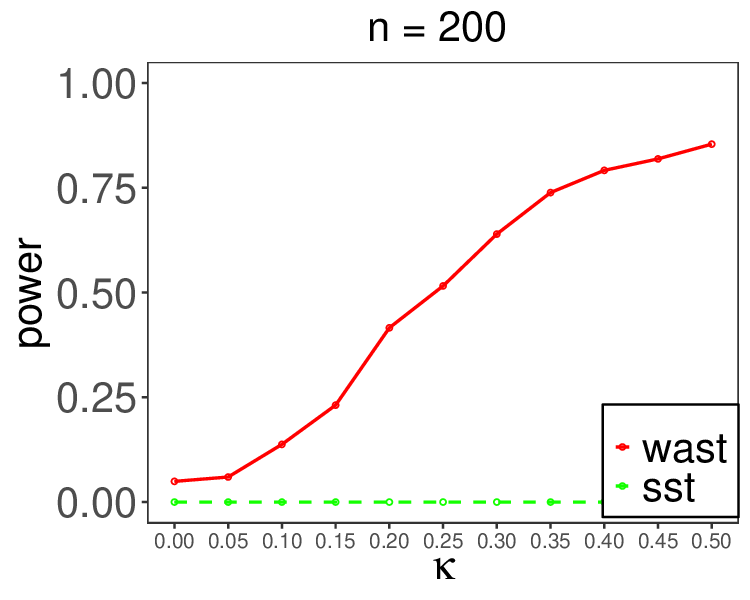}
		\includegraphics[scale=0.33]{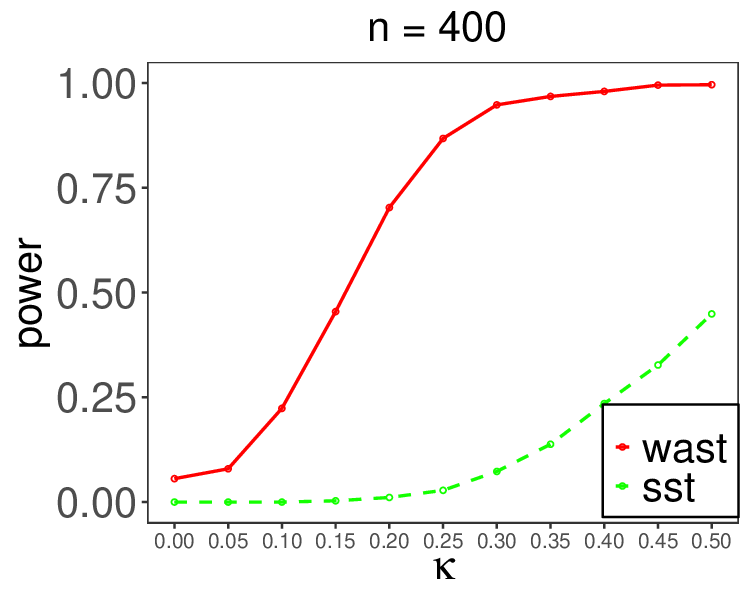}
		\includegraphics[scale=0.33]{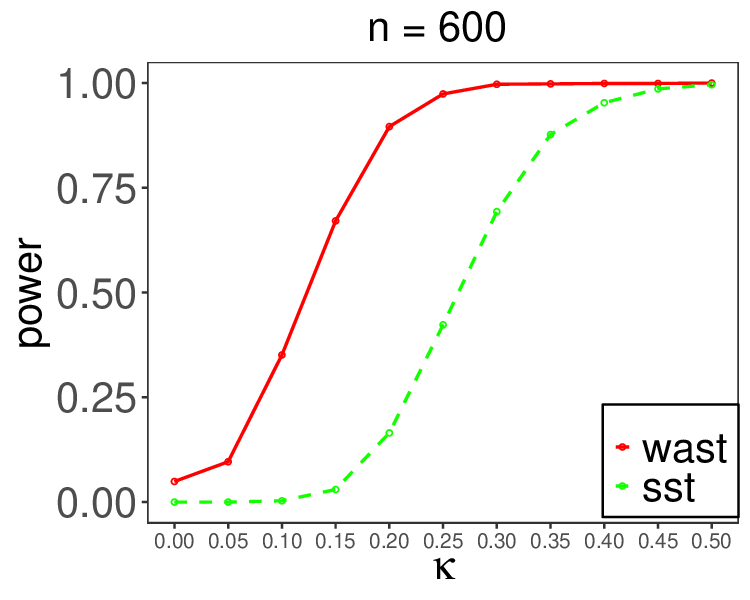}      \\
		\includegraphics[scale=0.33]{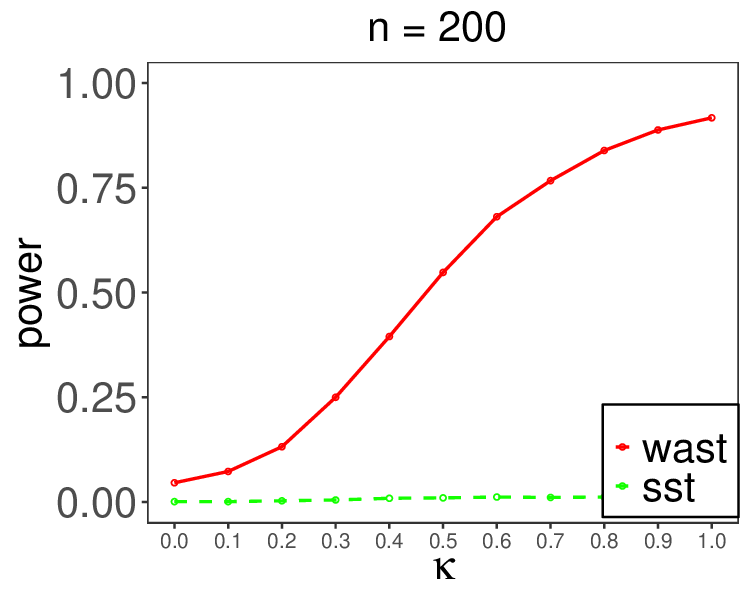}
		\includegraphics[scale=0.33]{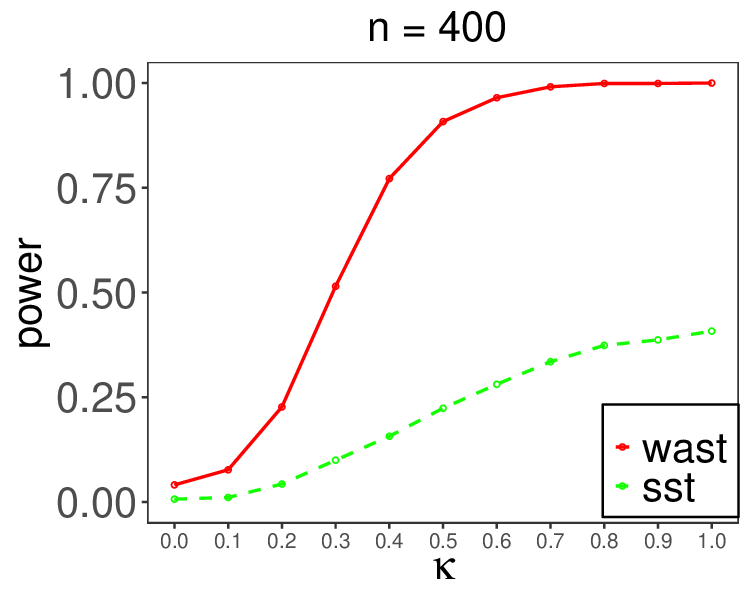}
		\includegraphics[scale=0.33]{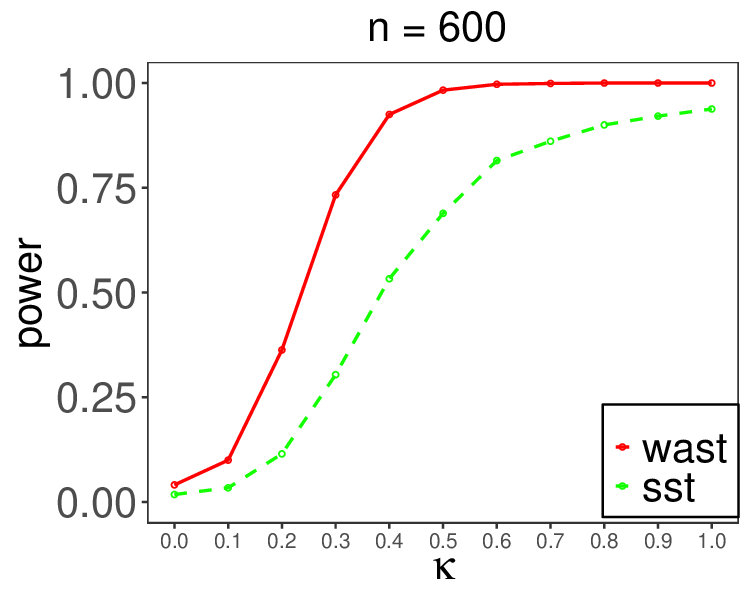}   \\
		\includegraphics[scale=0.33]{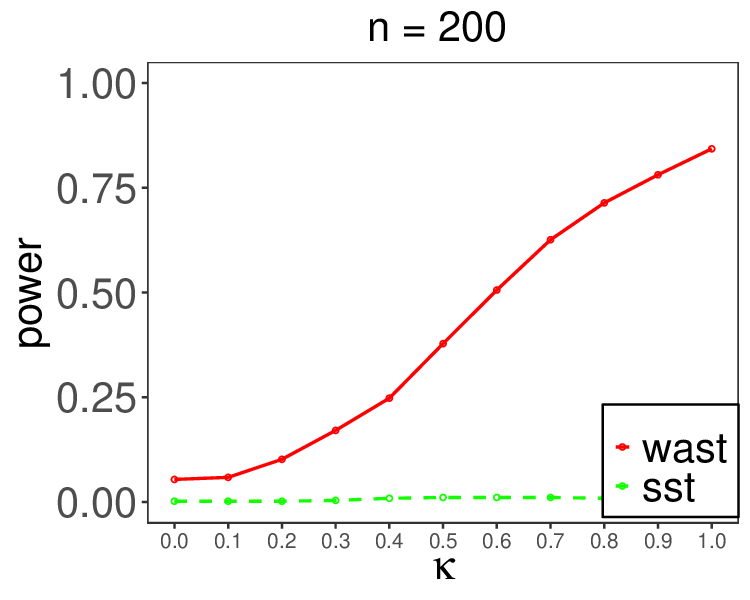}
		\includegraphics[scale=0.33]{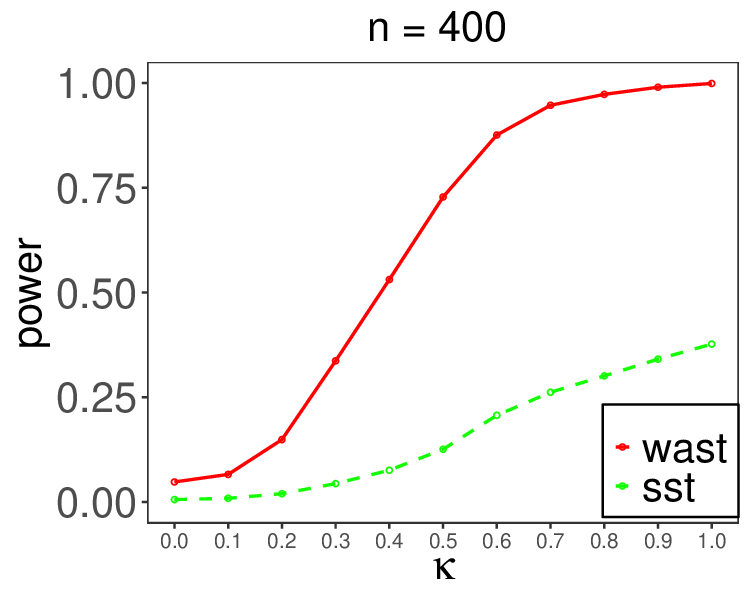}
		\includegraphics[scale=0.33]{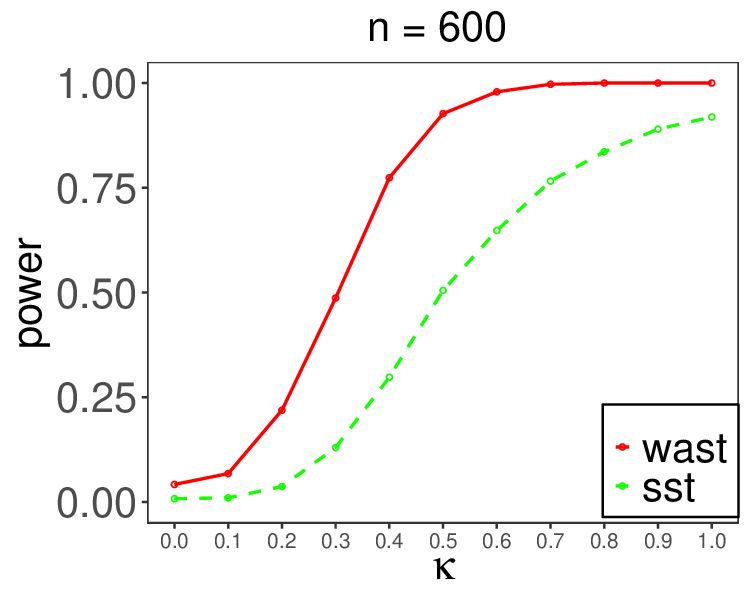}  \\
		\includegraphics[scale=0.33]{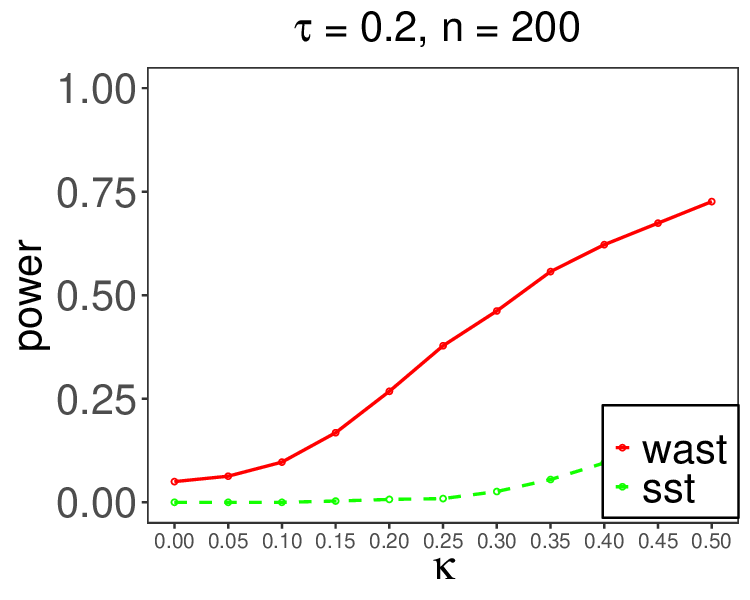}
		\includegraphics[scale=0.33]{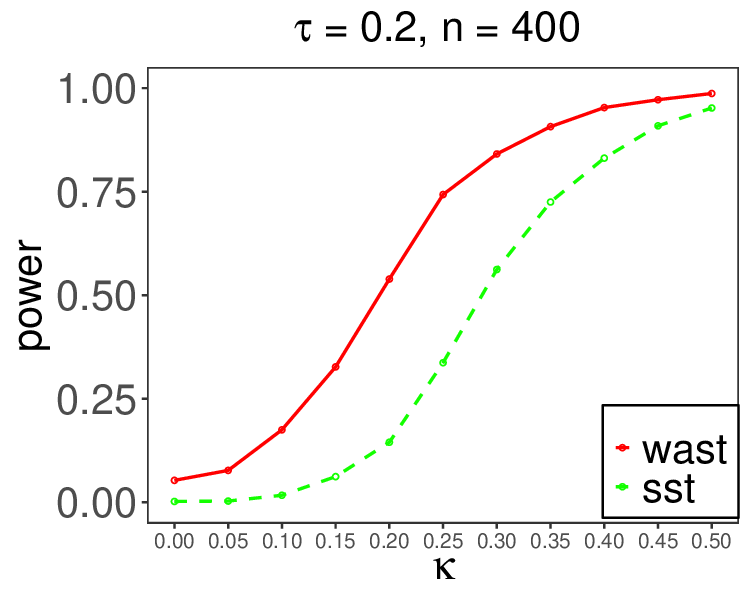}
		\includegraphics[scale=0.33]{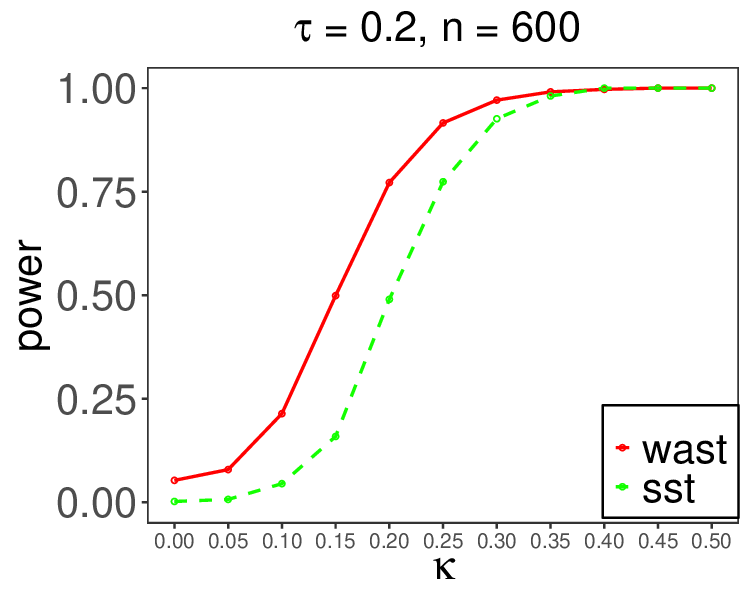}   \\
		\includegraphics[scale=0.33]{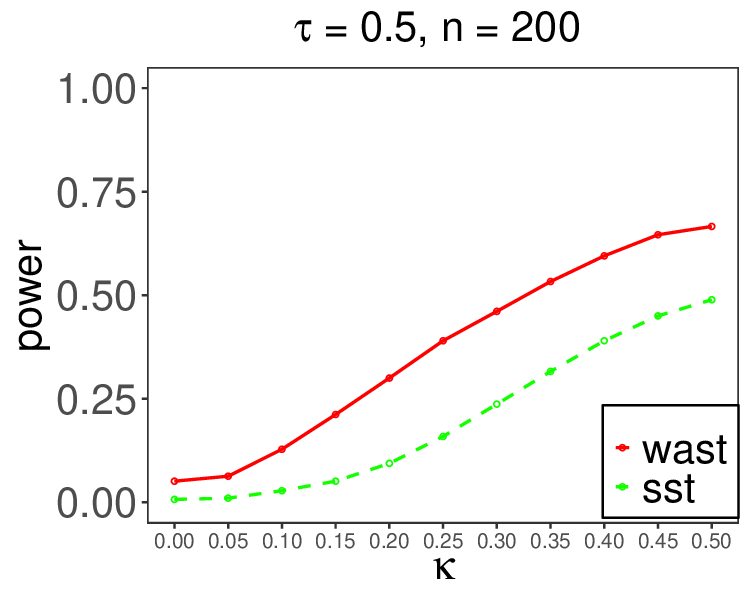}
		\includegraphics[scale=0.33]{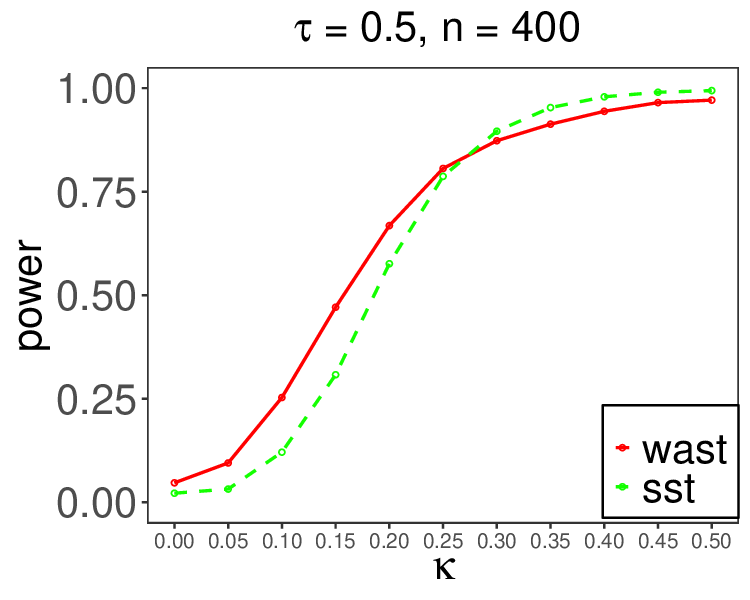}
		\includegraphics[scale=0.33]{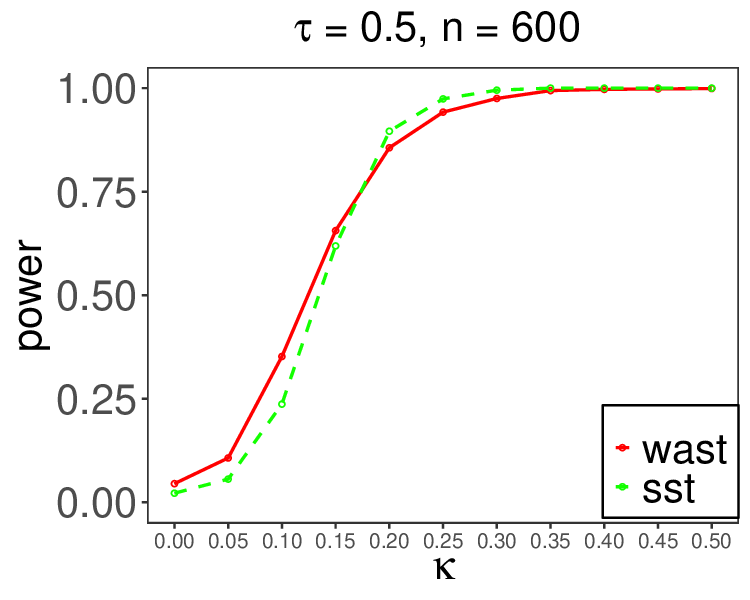}   \\
		\includegraphics[scale=0.33]{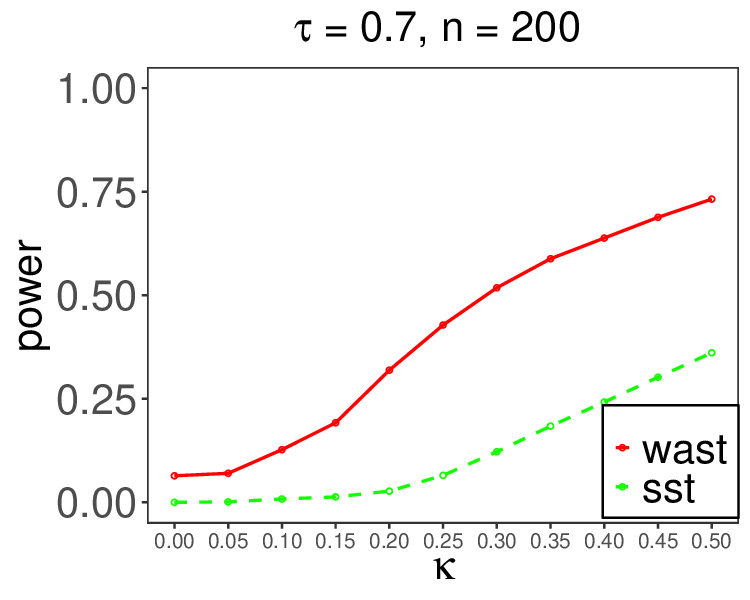}
		\includegraphics[scale=0.33]{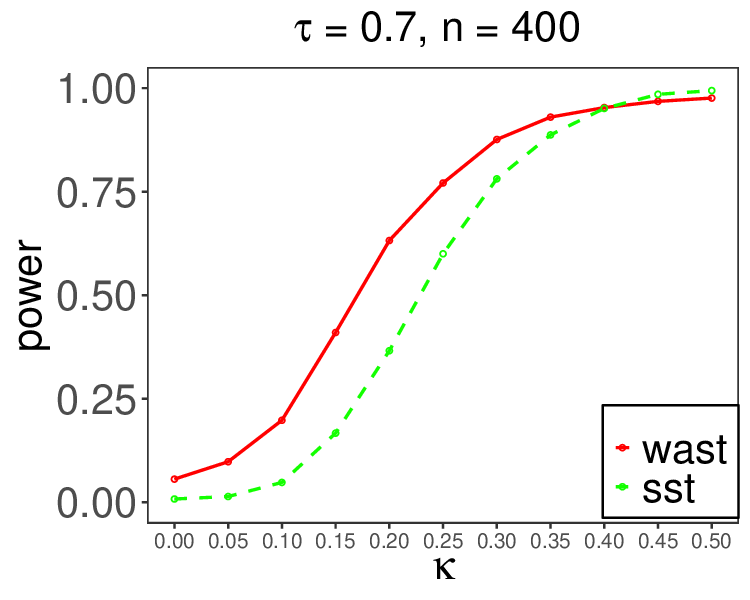}
		\includegraphics[scale=0.33]{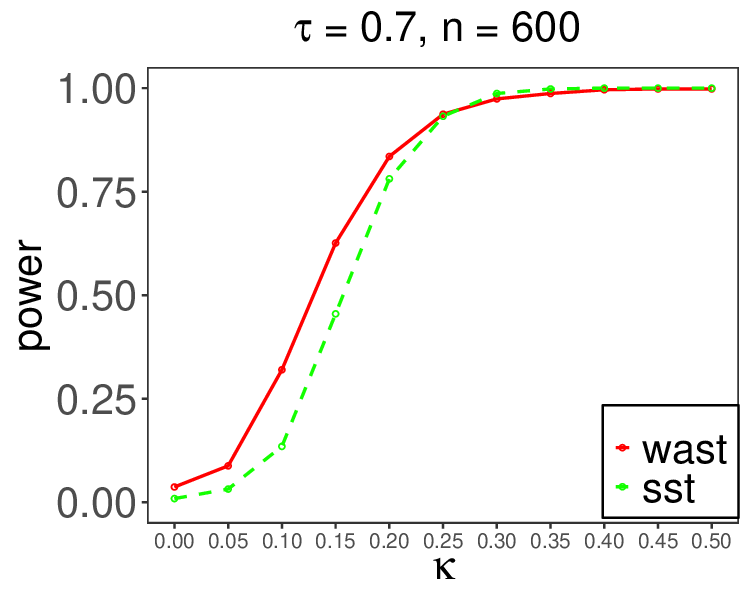}
		\caption{\it Powers of test statistic by the proposed WAST (red solid line) and SST (green dashed line) for $(p,q)=(10,10)$. From top to bottom, each row depicts the powers for probit model, semiparametric model, quantile regression with $\tau=0.2$, $\tau=0.5$ and $\tau=0.7$, respectively. Here the $\Gv$ $Z$ is generated from multivariate normal distribution with mean $\bzero$ and covariance $25I$.}
		\label{fig_qr1010_1}
	\end{center}
\end{figure}

\begin{figure}[!ht]
	\begin{center}
		\includegraphics[scale=0.33]{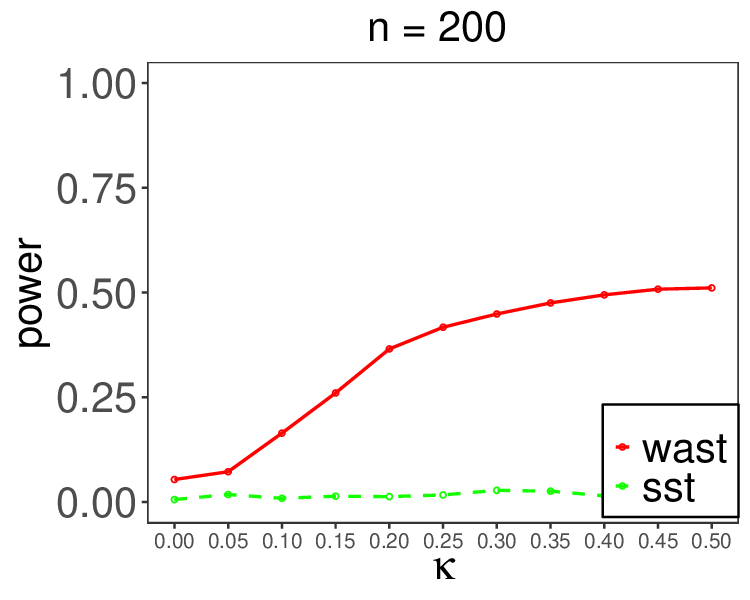}
		\includegraphics[scale=0.33]{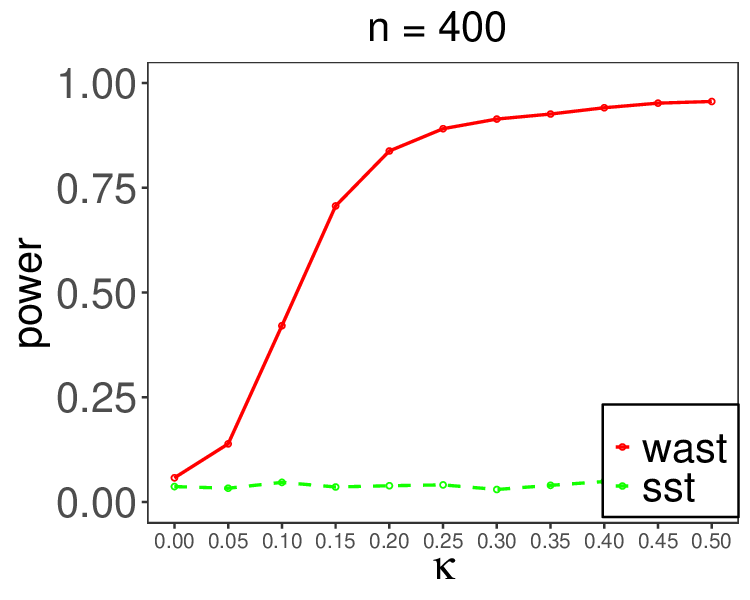}
		\includegraphics[scale=0.33]{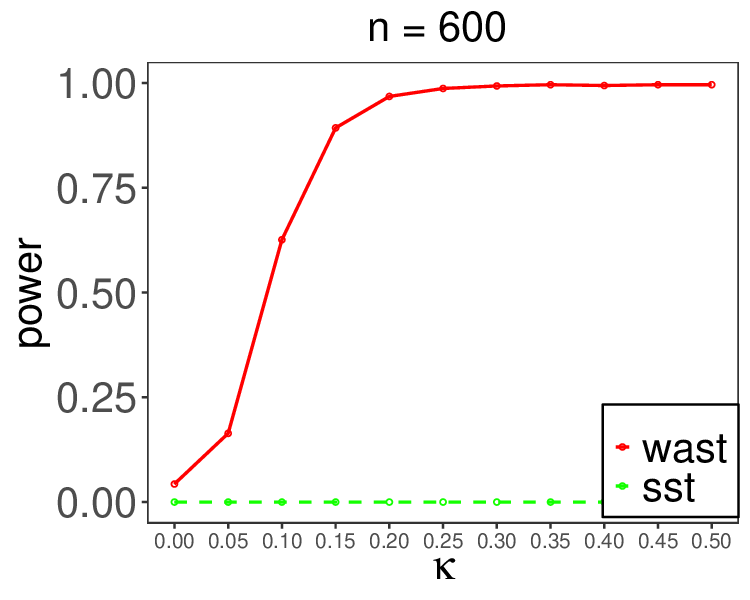}     \\
		\includegraphics[scale=0.33]{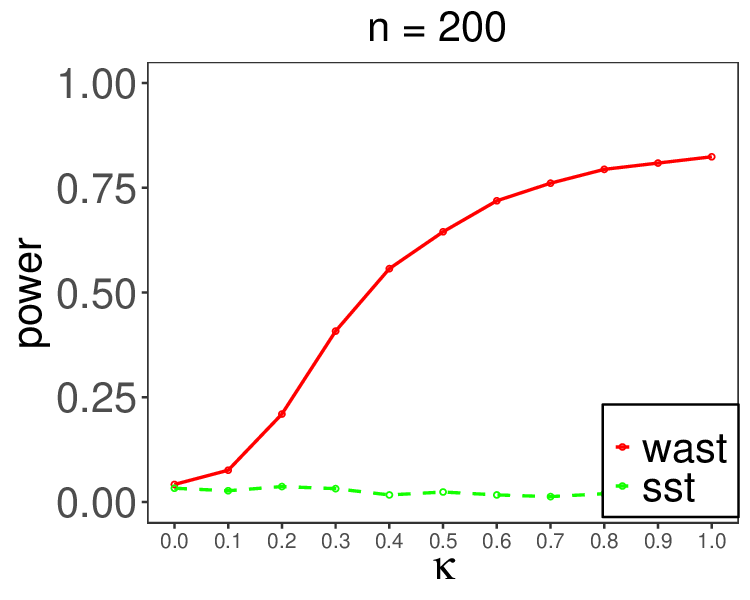}
		\includegraphics[scale=0.33]{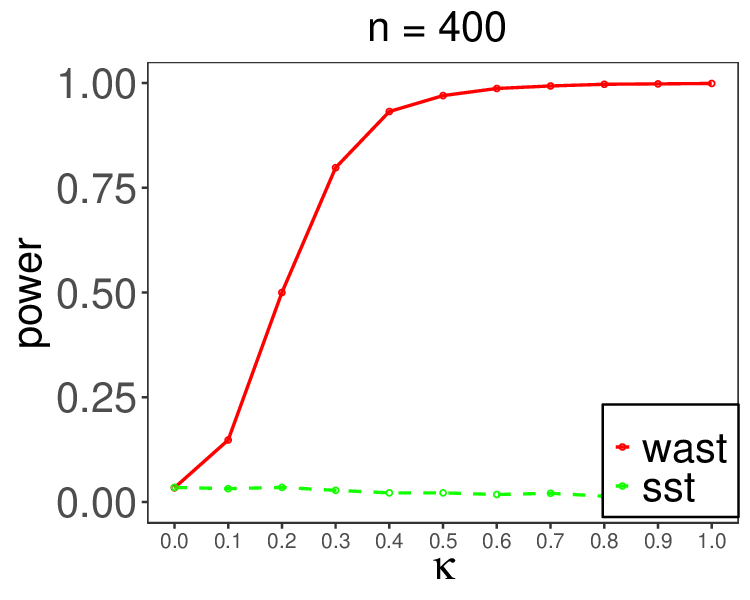}
		\includegraphics[scale=0.33]{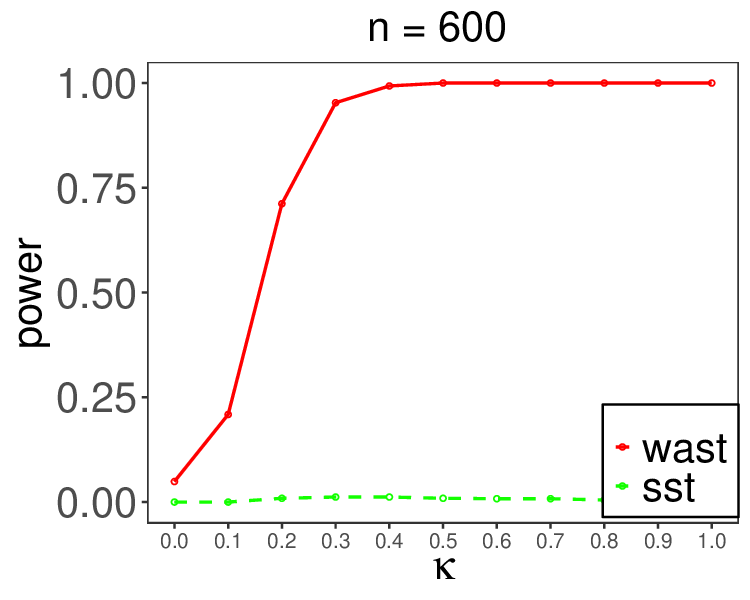}  \\
		\includegraphics[scale=0.33]{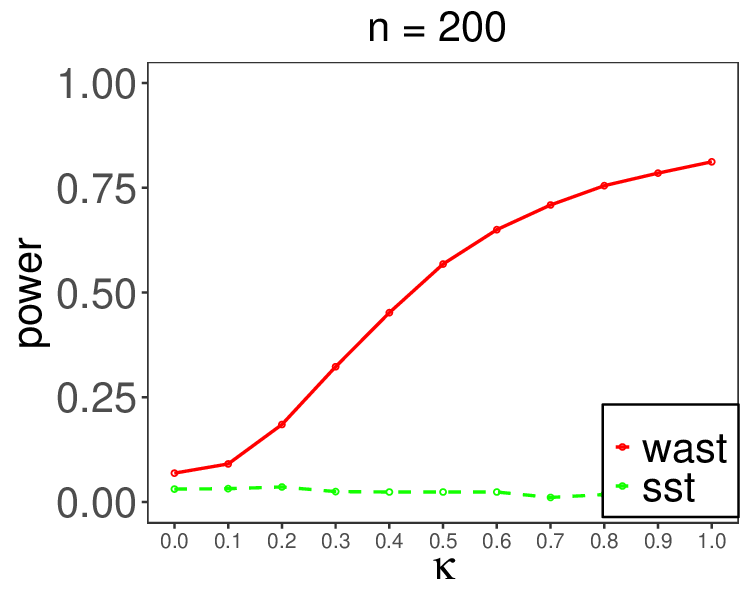}
		\includegraphics[scale=0.33]{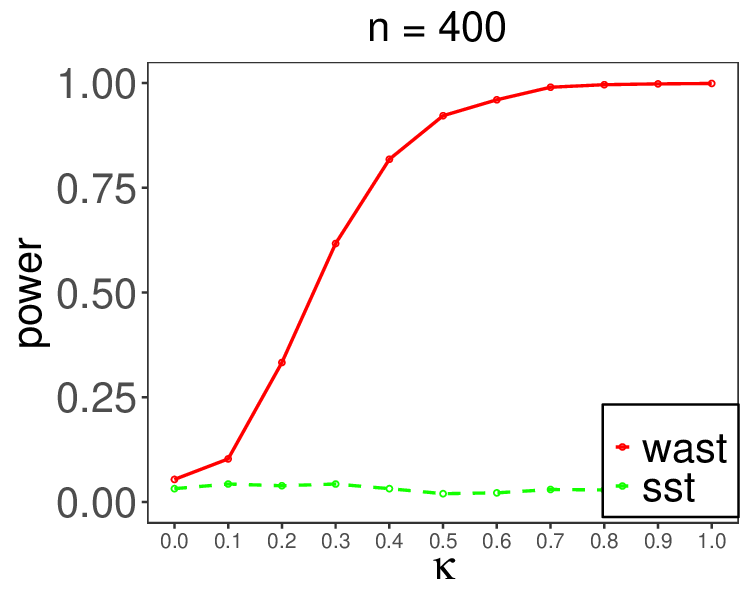}
		\includegraphics[scale=0.33]{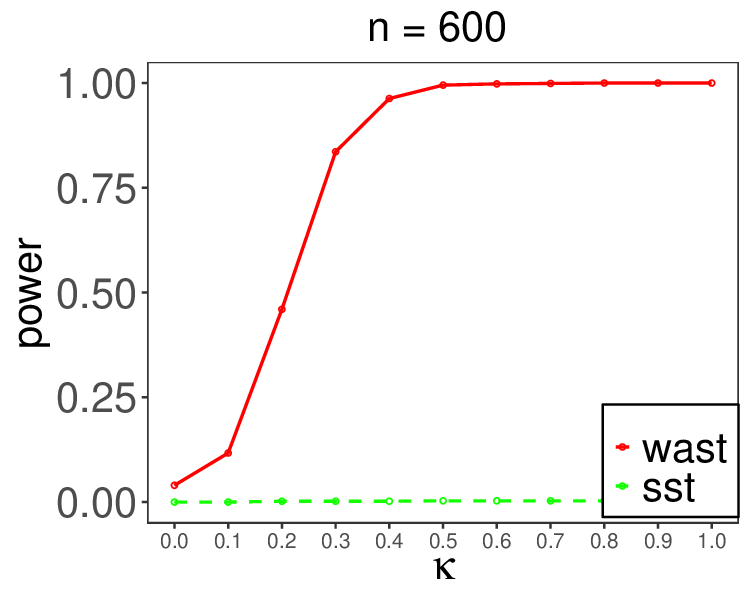} \\
		\includegraphics[scale=0.33]{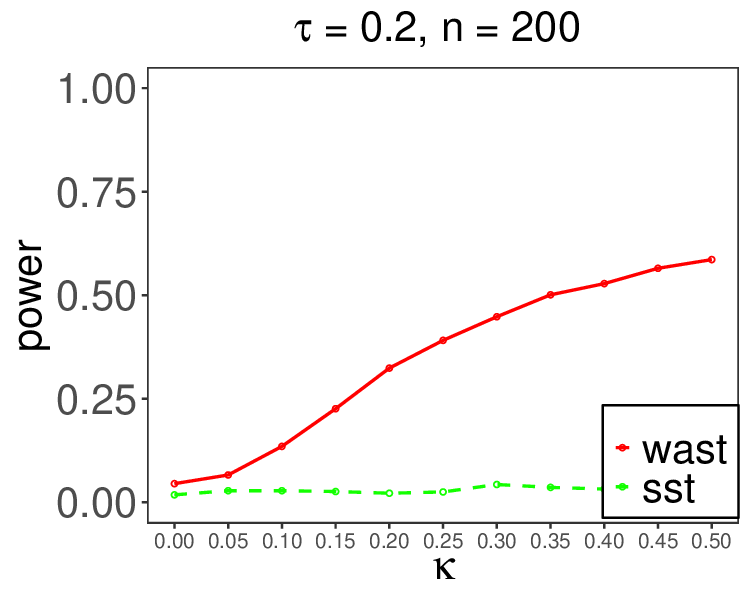}
		\includegraphics[scale=0.33]{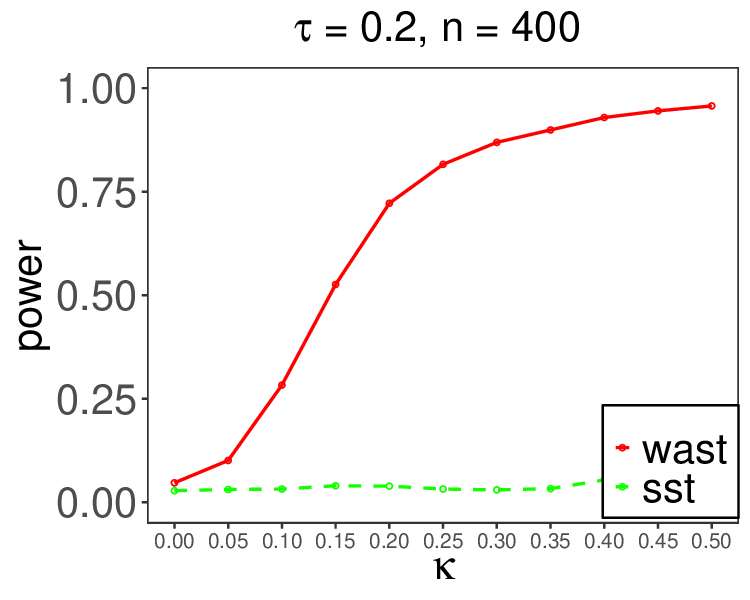}
		\includegraphics[scale=0.33]{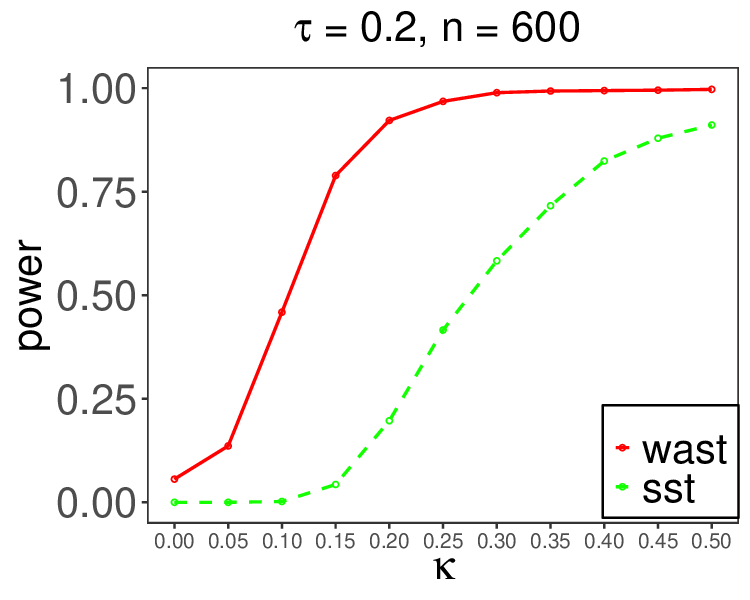}  \\
		\includegraphics[scale=0.33]{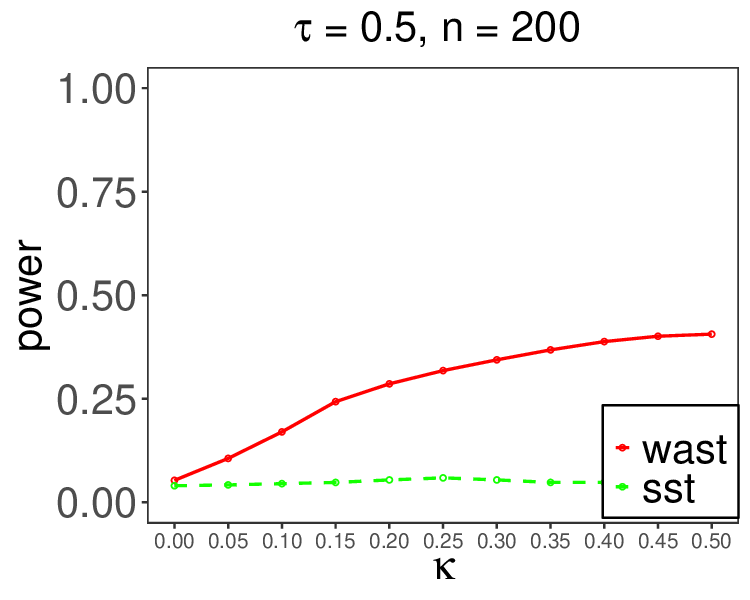}
		\includegraphics[scale=0.33]{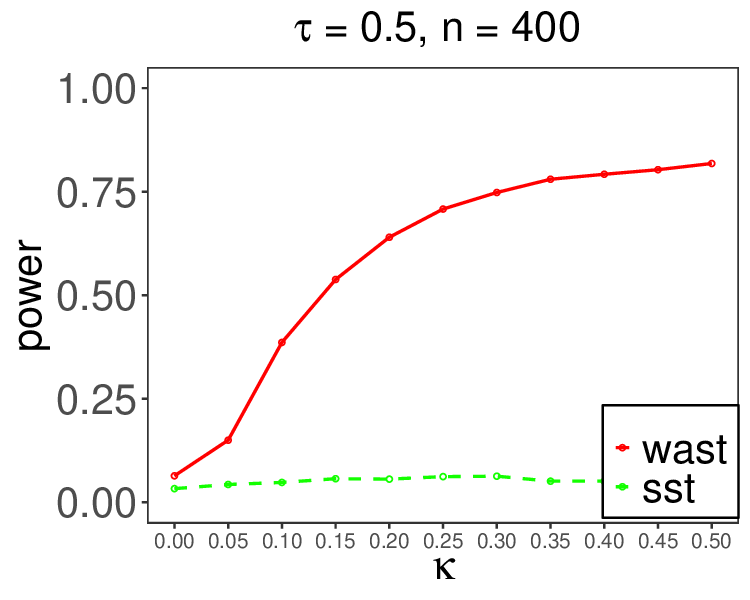}
		\includegraphics[scale=0.33]{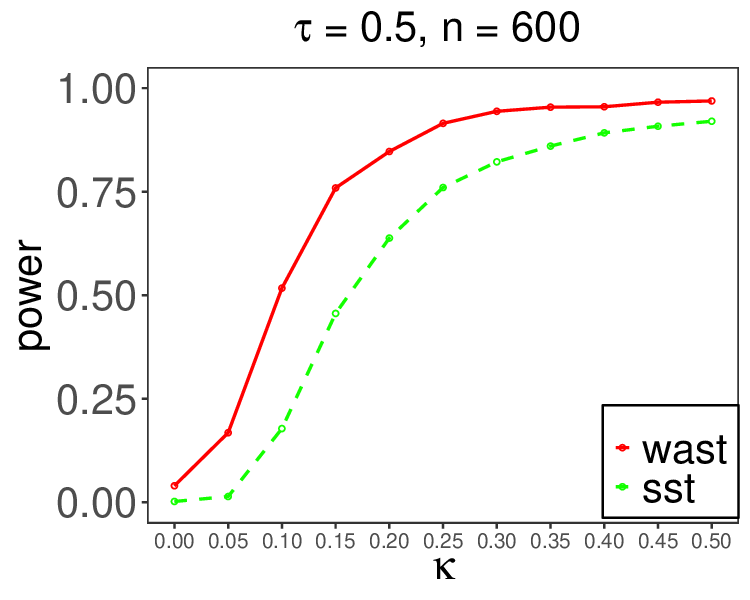}  \\
		\includegraphics[scale=0.33]{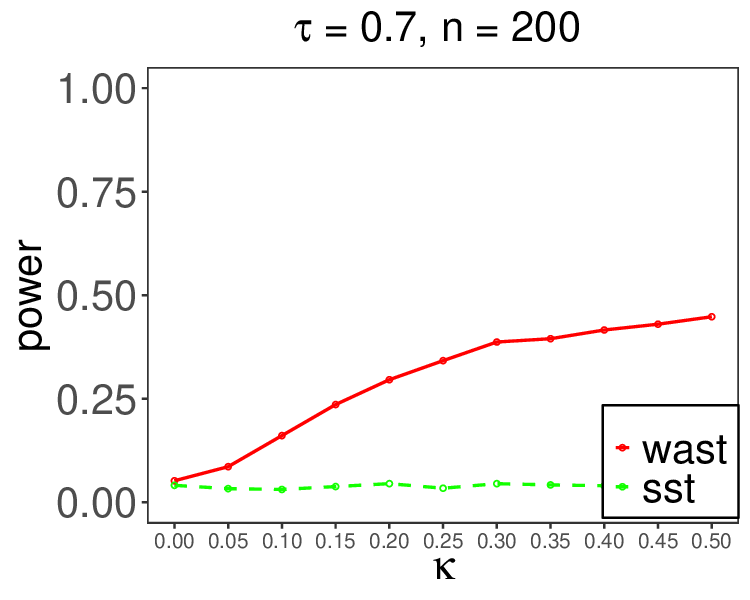}
		\includegraphics[scale=0.33]{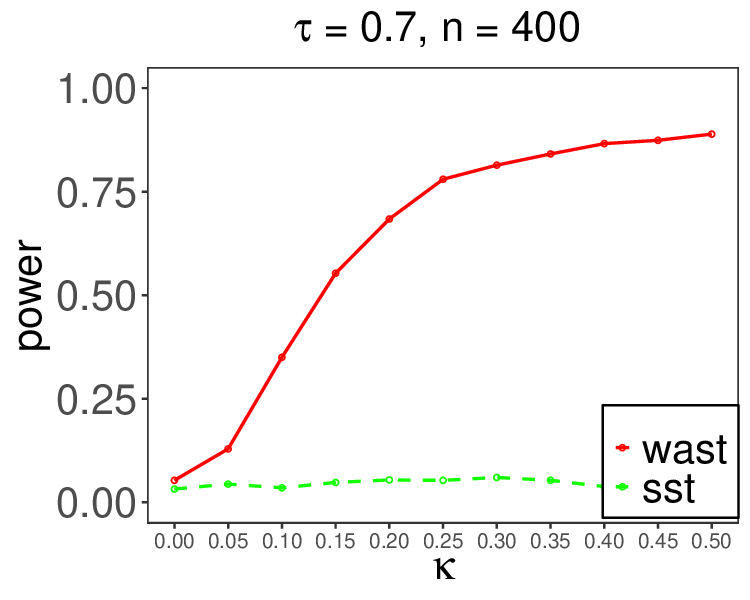}
		\includegraphics[scale=0.33]{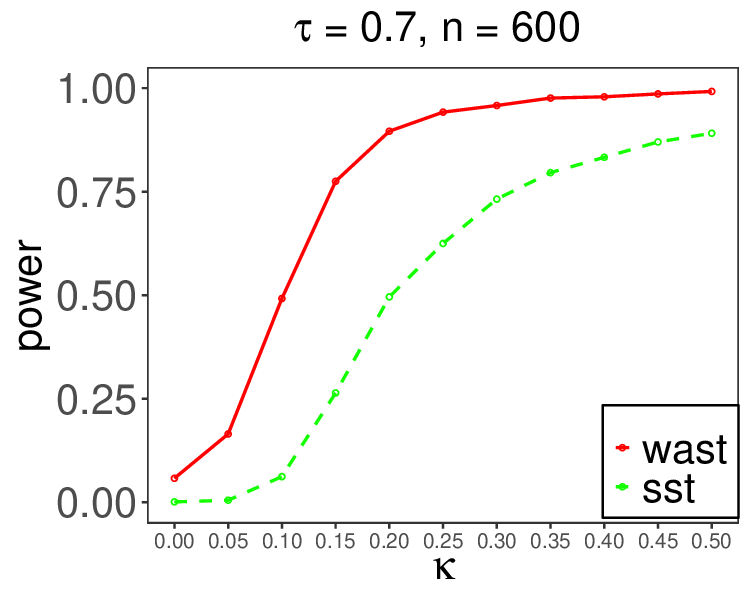}
		\caption{\it Powers of test statistic by the proposed WAST (red solid line) and SST (green dashed line) for $(p,q)=(50,20)$. From top to bottom, each row depicts the powers for probit model, semiparametric model, quantile regression with $\tau=0.2$, $\tau=0.5$ and $\tau=0.7$, respectively. Here the $\Gv$ $Z$ is generated from multivariate normal distribution with mean $\bzero$ and covariance $25I$.}
		\label{fig_qr5020_1}
	\end{center}
\end{figure}

Type \uppercase\expandafter{\romannumeral1} errors of the proposed test statistic are listed in Table \ref{table_size2}, and the power curves are depicted in Figure \ref{fig_qr13_1}-\ref{fig_qr5020_1}. From Table \ref{table_size2}, we can see that Type \uppercase\expandafter{\romannumeral1} errors of the proposed test statistic are all close to the nominal significant level 0.05, while that of the SST are much smaller the 0.05 except for the semiparametric model with B2 plus P2.
We can find from Figure \ref{fig_qr13_1}-\ref{fig_qr5020_1} that the power of the proposed test statistic increases much faster than that of SST.

\begin{table}[htp!]
	\def~{\hphantom{0}}
    \tiny
	\caption{Type \uppercase\expandafter{\romannumeral1} errors of the proposed WAST and SST based on resampling for probit regression model (ProbitRE), quantile regression (QuantRE) and semiparametric model (SPMoldel).
}
	\resizebox{\textwidth}{!}{
        \begin{threeparttable}
		\begin{tabular}{llcccccccc}
			\hline
			\multirow{2}{*}{Model}&\multirow{2}{*}{$(p,q)$}
			&\multicolumn{2}{c}{ $n=200$} && \multicolumn{2}{c}{ $n=400$} && \multicolumn{2}{c}{ $n=600$} \\
			\cline{3-4} \cline{6-7} \cline{9-10}
			&&   WAST & SST && WAST & SST && WAST & SST \\
			\cline{3-10}
			ProbitRE &$(1,3)$         & 0.062 & 0.003 && 0.050 & 0.012 && 0.047 & 0.010 \\
			&$(5,5)$                  & 0.044 & 0.000 && 0.030 & 0.000 && 0.041 & 0.001 \\
			& $(10,10)$               & 0.049 & 0.000 && 0.056 & 0.000 && 0.049 & 0.000 \\
			&$(50,20)$                & 0.054 & 0.006 && 0.058 & 0.037 && 0.043 & 0.000 \\
			[1 ex]
			SPModel &$(1,3)$          & 0.050 & 0.034 && 0.055 & 0.039 && 0.053 & 0.042 \\
			(B1+P1)&$(5,5)$           & 0.045 & 0.008 && 0.040 & 0.020 && 0.040 & 0.021 \\
			& $(10,10)$               & 0.046 & 0.001 && 0.041 & 0.007 && 0.041 & 0.018 \\
			&$(50,20)$                & 0.042 & 0.033 && 0.034 & 0.035 && 0.049 & 0.000 \\
			[1 ex]
			SPModel &$(1,3)$          & 0.057 & 0.044 && 0.049 & 0.048 && 0.047 & 0.053 \\
			(B2+P2)&$(5,5)$           & 0.056 & 0.007 && 0.062 & 0.013 && 0.034 & 0.013 \\
			& $(10,10)$               & 0.054 & 0.002 && 0.048 & 0.006 && 0.042 & 0.008 \\
			&$(50,20)$                & 0.069 & 0.031 && 0.054 & 0.032 && 0.040 & 0.000 \\
			[1 ex]
			QuantRE &$(1,3)$          & 0.038 & 0.010 && 0.054 & 0.036 && 0.043 & 0.034 \\
			($\tau=0.2$)&$(5,5)$      & 0.050 & 0.004 && 0.058 & 0.011 && 0.040 & 0.014 \\
			& $(10,10)$               & 0.050 & 0.000 && 0.053 & 0.002 && 0.053 & 0.002 \\
			&$(50,20)$                & 0.045 & 0.018 && 0.047 & 0.028 && 0.056 & 0.000 \\
			[1 ex]
			QuantRE &$(1,3)$          & 0.045 & 0.031 && 0.048 & 0.038 && 0.062 & 0.050 \\
			($\tau=0.5$)&$(5,5)$      & 0.055 & 0.013 && 0.055 & 0.028 && 0.044 & 0.037 \\
			& $(10,10)$               & 0.051 & 0.007 && 0.047 & 0.022 && 0.045 & 0.022 \\
			&$(50,20)$                & 0.053 & 0.040 && 0.064 & 0.033 && 0.040 & 0.002 \\
			[1 ex]
			QuantRE &$(1,3)$          & 0.035 & 0.026 && 0.041 & 0.032 && 0.066 & 0.040 \\
			($\tau=0.7$)&$(5,5)$      & 0.051 & 0.005 && 0.063 & 0.015 && 0.046 & 0.023 \\
			& $(10,10)$               & 0.064 & 0.000 && 0.056 & 0.008 && 0.037 & 0.009 \\
			&$(50,20)$                & 0.052 & 0.041 && 0.053 & 0.032 && 0.058 & 0.001 \\
			\hline
		\end{tabular}
\begin{tablenotes}
\item The nominal significant level is 0.05. Here the $\Gv$ $\bZ$ is generated from multivariate normal distribution with mean $\bzero$ and covariance $25\bI$.
\end{tablenotes}
\end{threeparttable}
	}
	\label{table_size2}
\end{table}

\subsection{Change plane analysis for quantile, probit and semiparametric models with \texorpdfstring{$Z$}{} from \texorpdfstring{$t_3$}{} distribution}\label{simulation_cpzt3}
We consider probit regression, quantile regression and semiparametric models with same settings as Section \ref{simulation_cpz} in current Supplementary Material but different $\Gv$ $\bZ$ and error's distribution of quantile regression model. We generate $Z_{ji}$ independently from $t_3$ distribution with 3 degrees of freedom, and set $\bZ_i=(1,Z_{1i},\cdots,Z_{(q-1)i})\trans$, where $j=1,\cdots,q-1$. For quantile regression, we generate error $\eps_i$ independently from $t_3$ distribution with degrees freedom 3.

\begin{figure}[!ht]
	\begin{center}
		\includegraphics[scale=0.33]{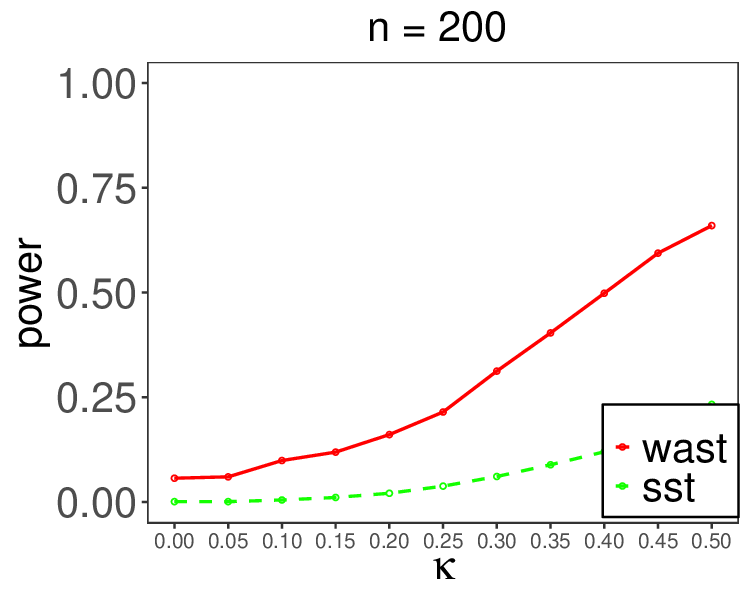}
		\includegraphics[scale=0.33]{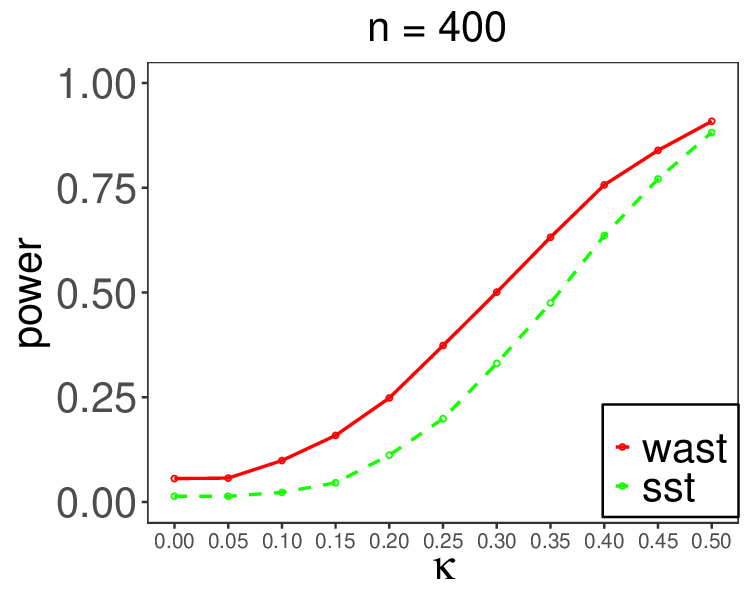}
		\includegraphics[scale=0.33]{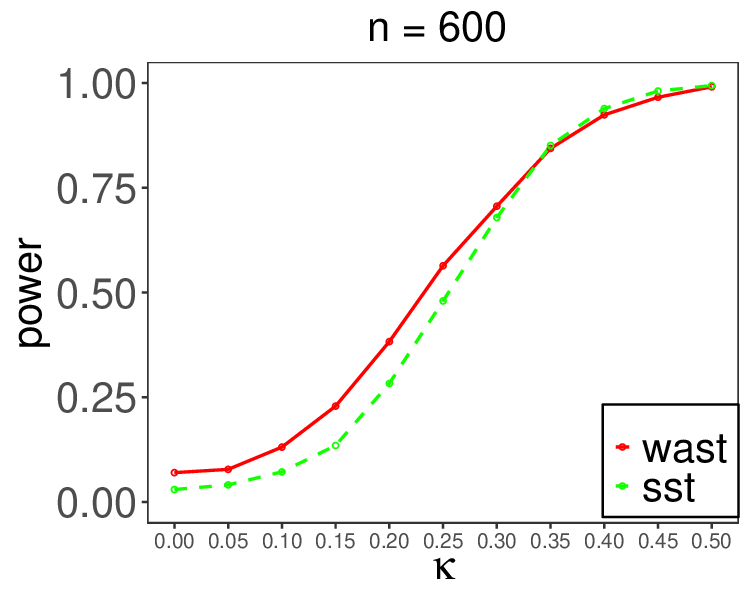}      \\
		\includegraphics[scale=0.33]{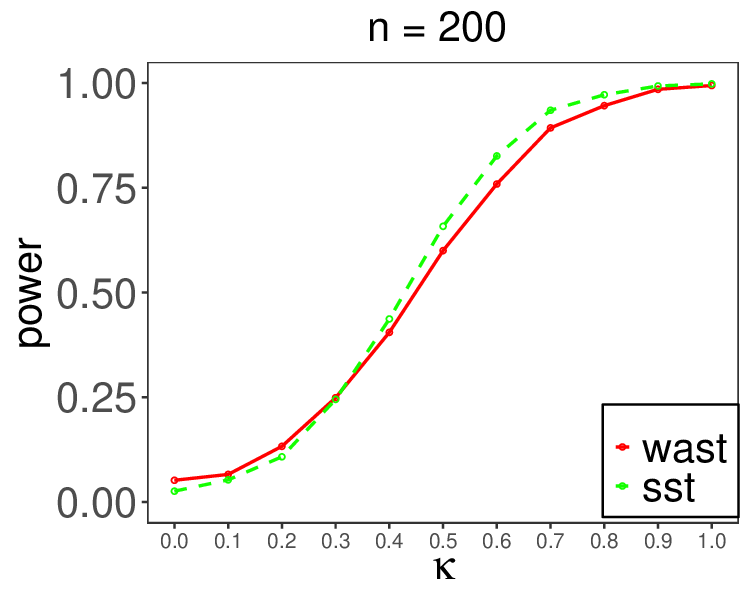}
		\includegraphics[scale=0.33]{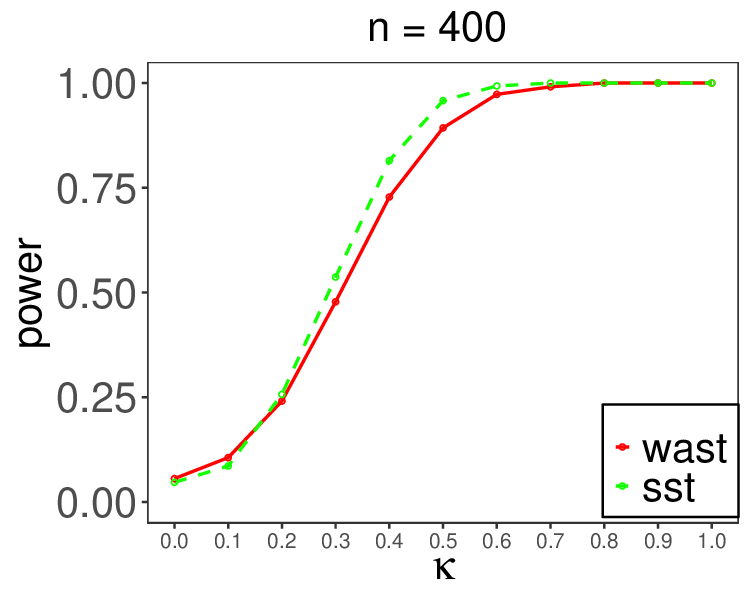}
		\includegraphics[scale=0.33]{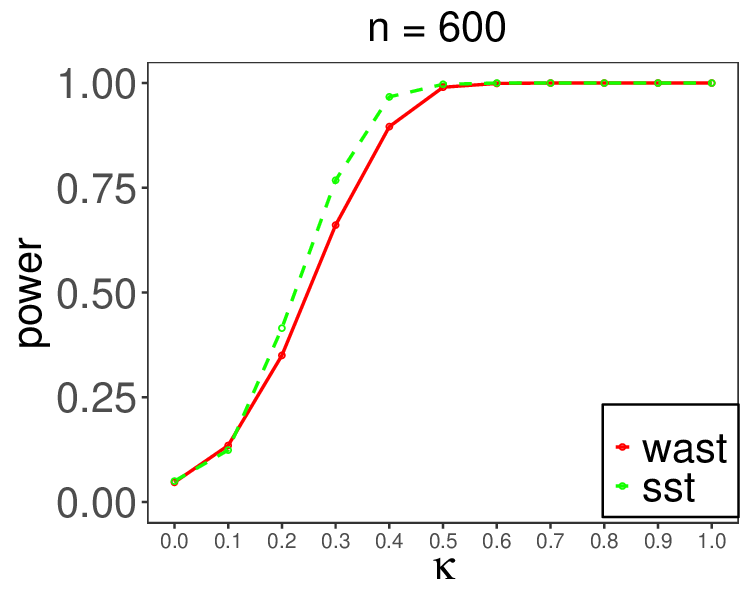}   \\
		\includegraphics[scale=0.33]{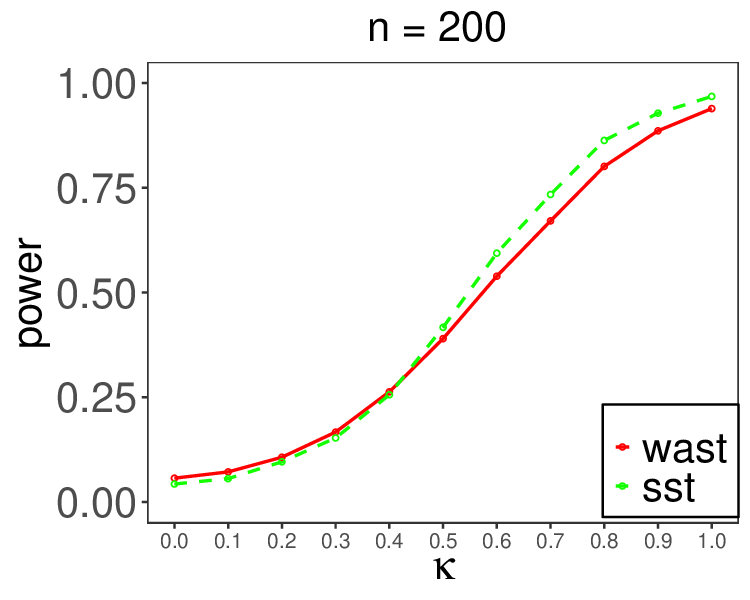}
		\includegraphics[scale=0.33]{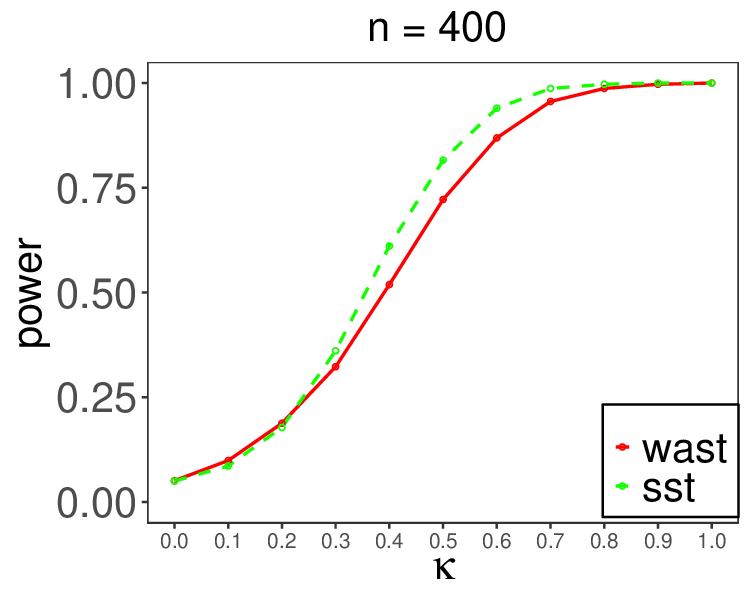}
		\includegraphics[scale=0.33]{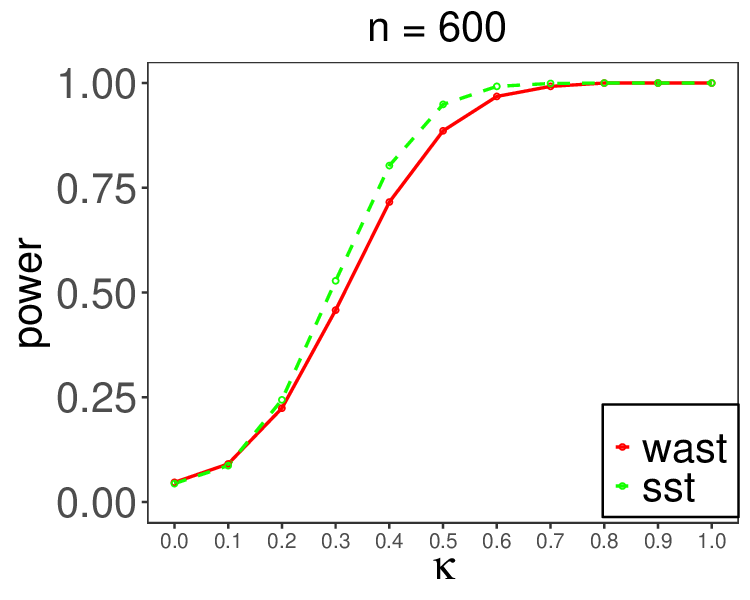}  \\
		\includegraphics[scale=0.33]{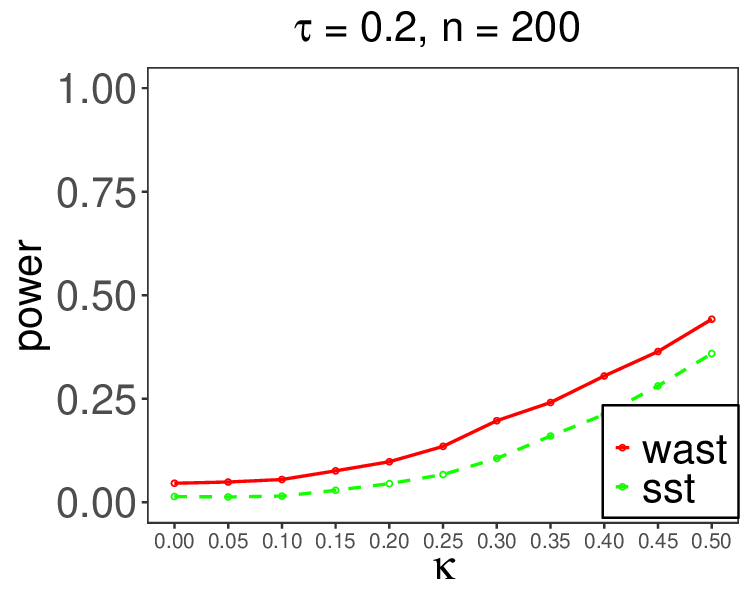}
		\includegraphics[scale=0.33]{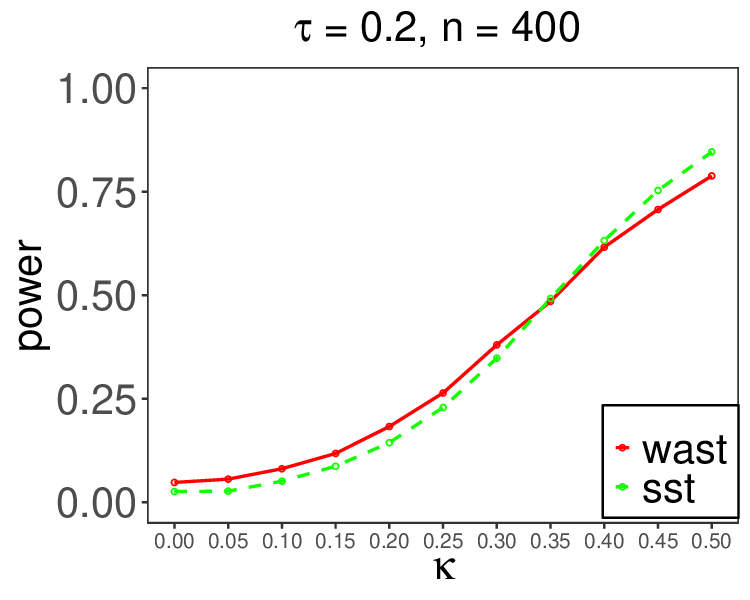}
		\includegraphics[scale=0.33]{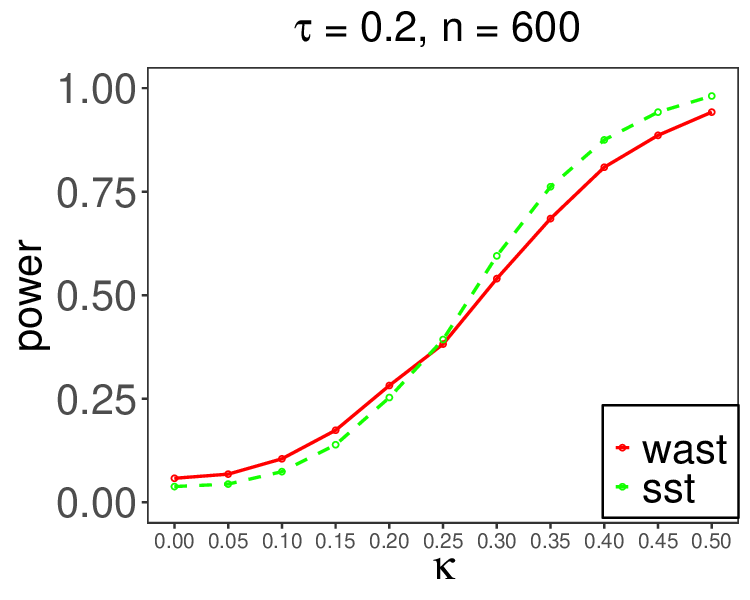}   \\
		\includegraphics[scale=0.33]{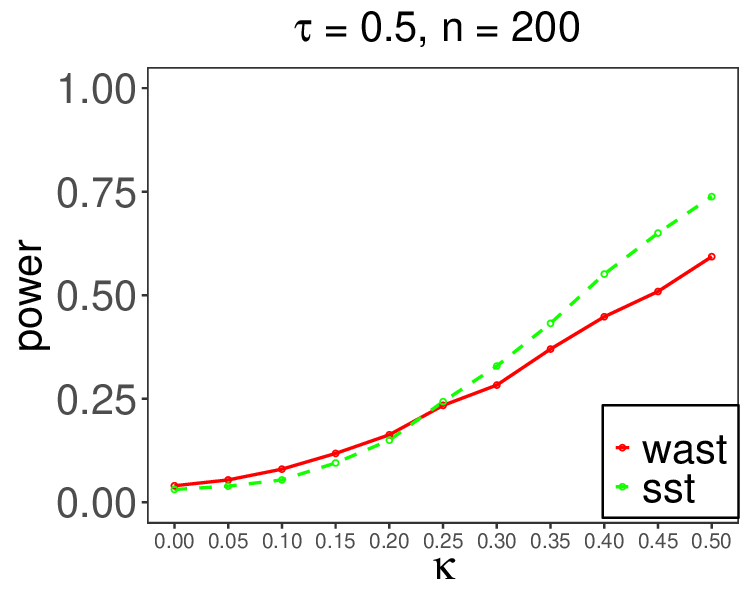}
		\includegraphics[scale=0.33]{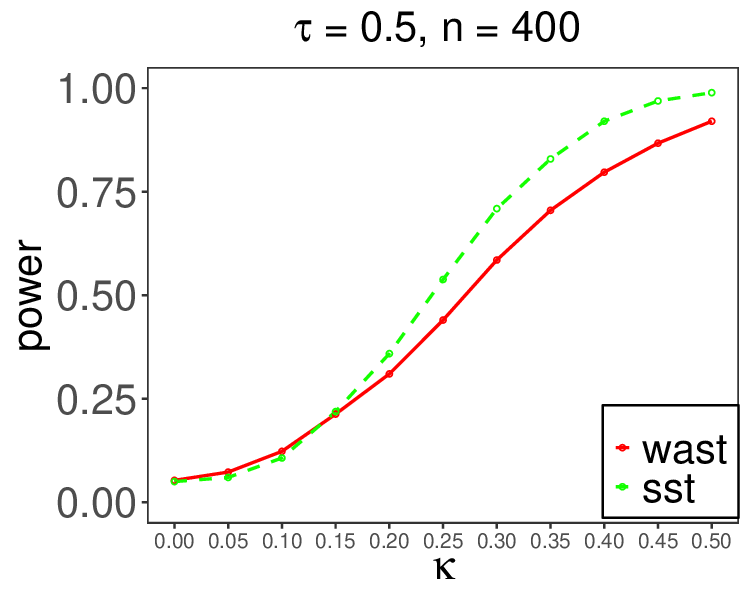}
		\includegraphics[scale=0.33]{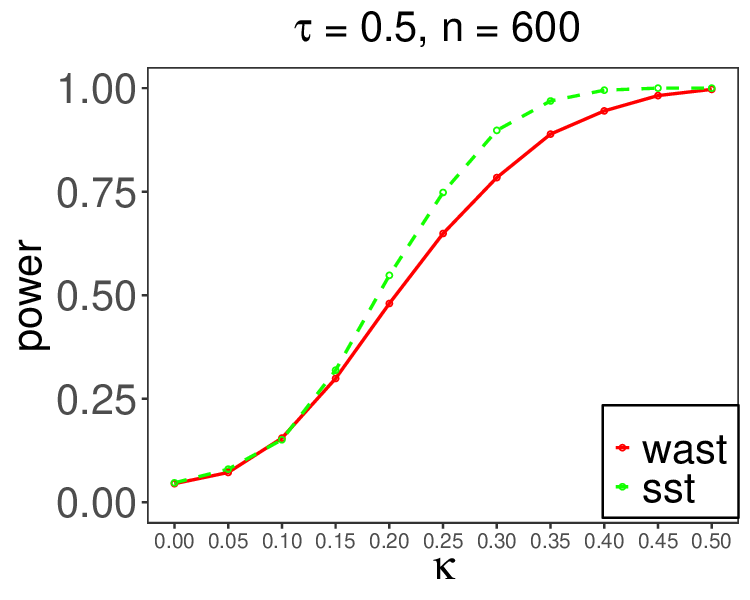}   \\
		\includegraphics[scale=0.33]{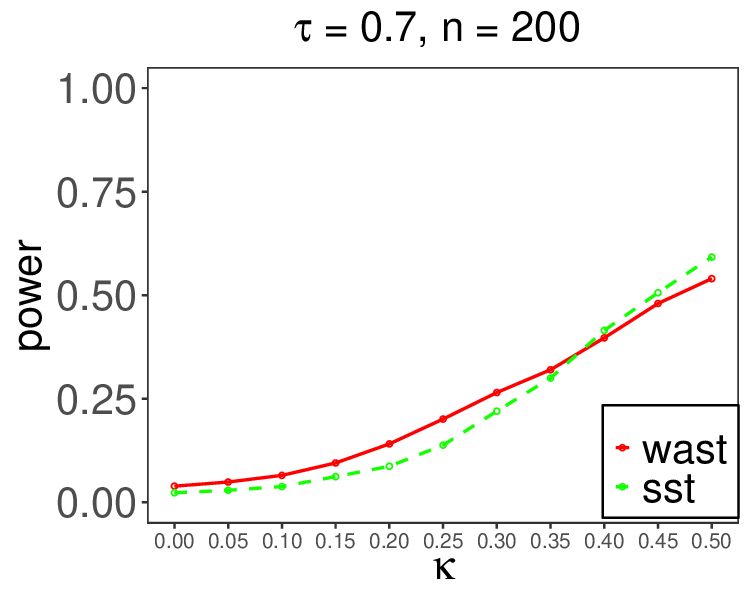}
		\includegraphics[scale=0.33]{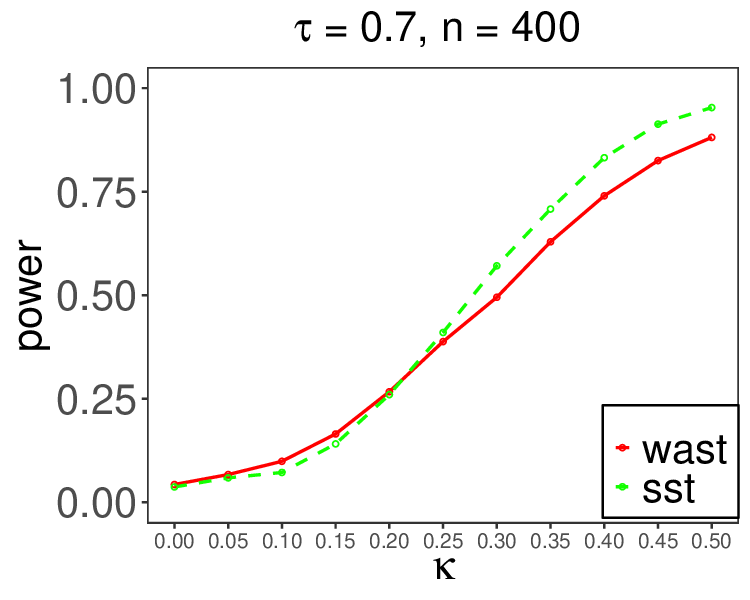}
		\includegraphics[scale=0.33]{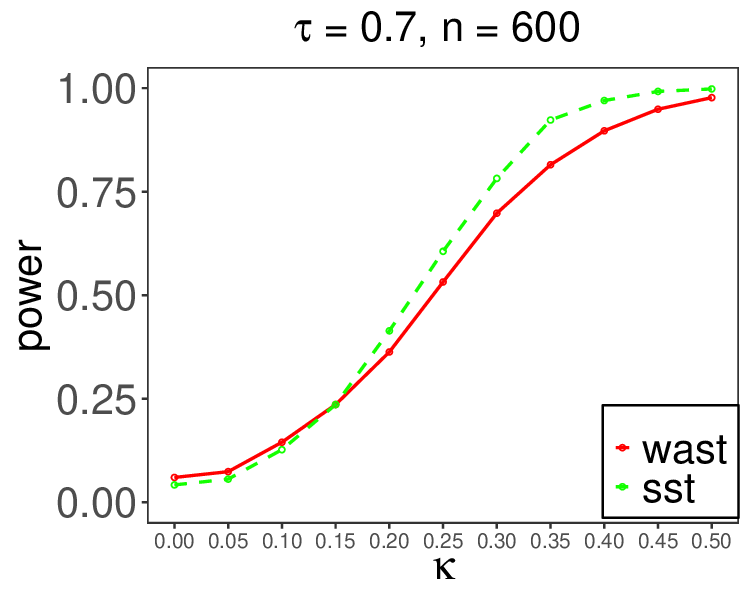}
		\caption{\it Powers of test statistic by the proposed WAST (red solid line) and SST (green dashed line) for $(p,q)=(1,3)$. From top to bottom, each row depicts the powers for probit model, semiparametric model, quantile regression with $\tau=0.2$, $\tau=0.5$ and $\tau=0.7$, respectively. Here the $\Gv$ $Z$ is generated from $t_3$ distribution.}
		\label{fig_qr13_zt3}
	\end{center}
\end{figure}

\begin{figure}[!ht]
	\begin{center}
		\includegraphics[scale=0.33]{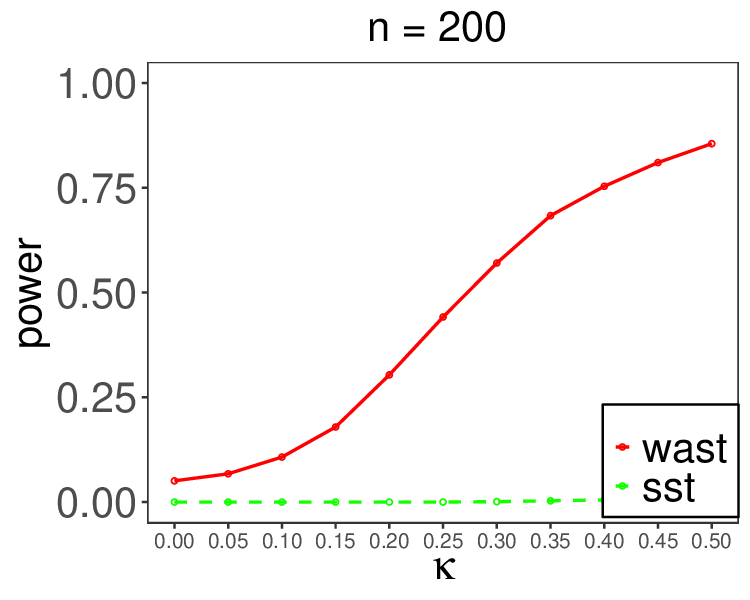}
		\includegraphics[scale=0.33]{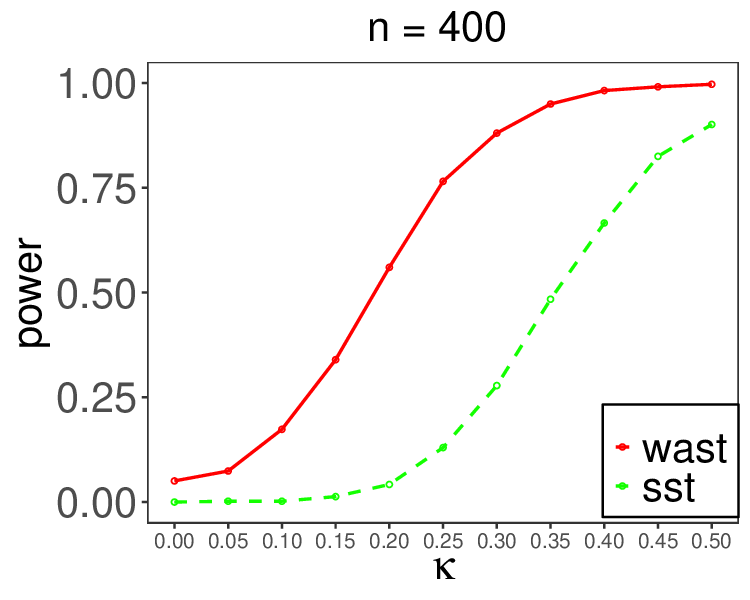}
		\includegraphics[scale=0.33]{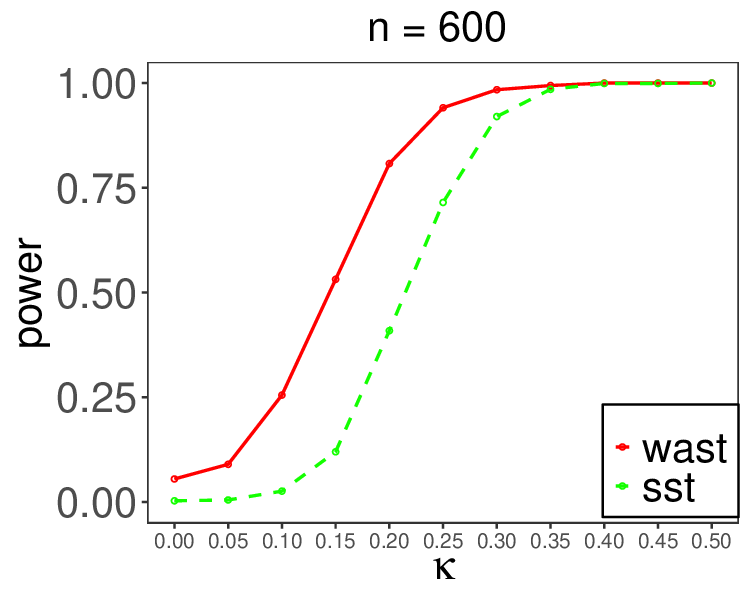}      \\
		\includegraphics[scale=0.33]{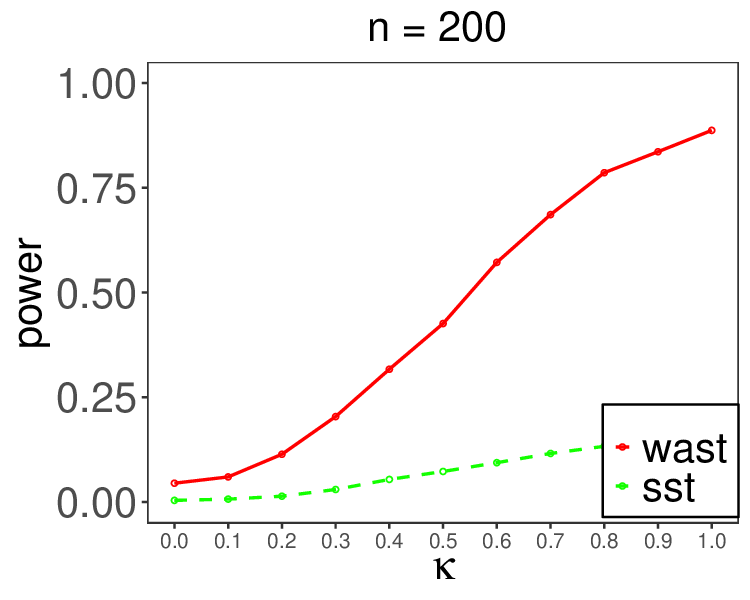}
		\includegraphics[scale=0.33]{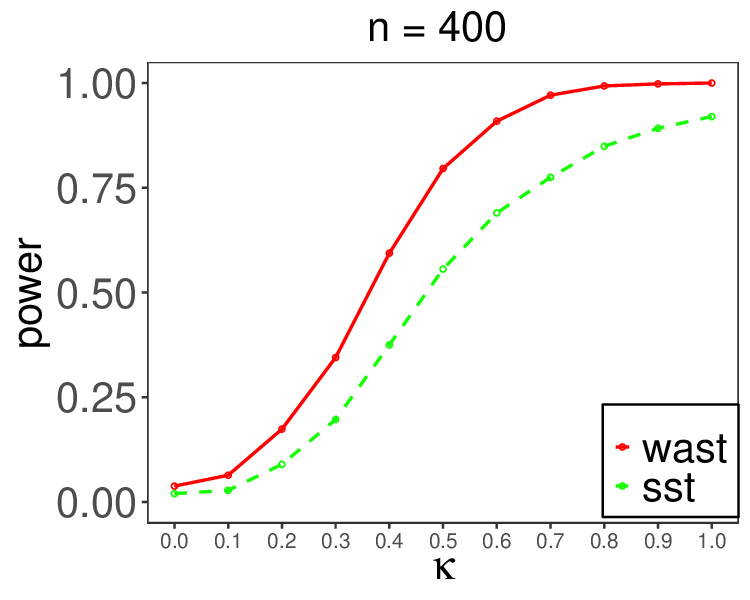}
		\includegraphics[scale=0.33]{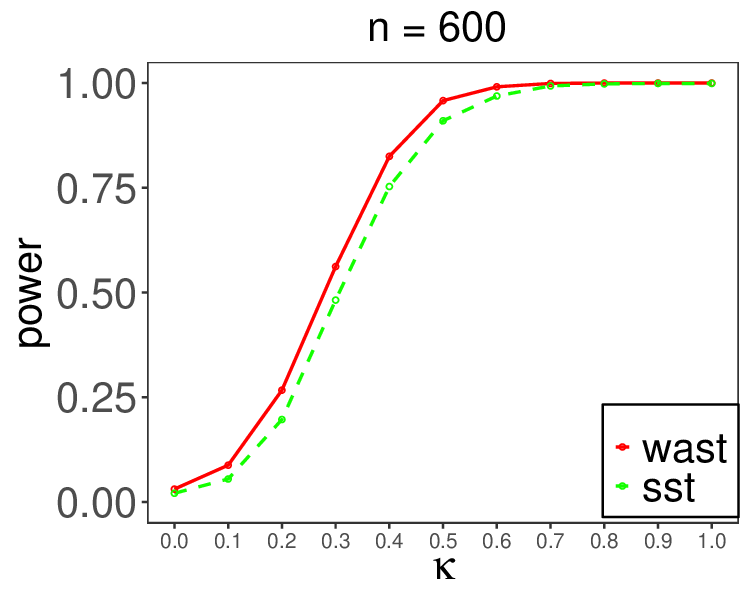}   \\
		\includegraphics[scale=0.33]{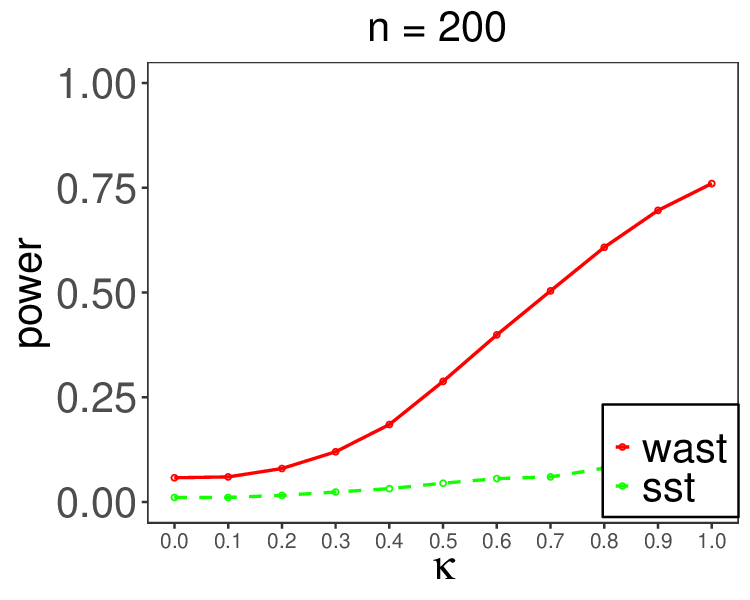}
		\includegraphics[scale=0.33]{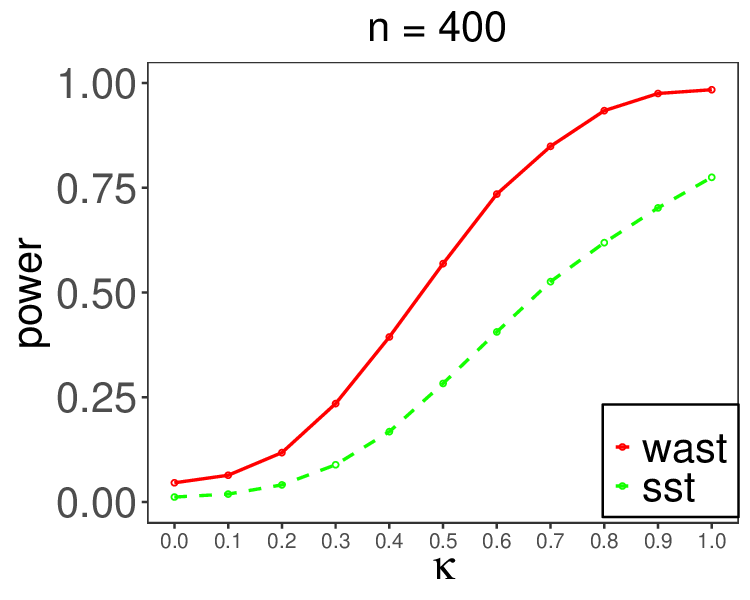}
		\includegraphics[scale=0.33]{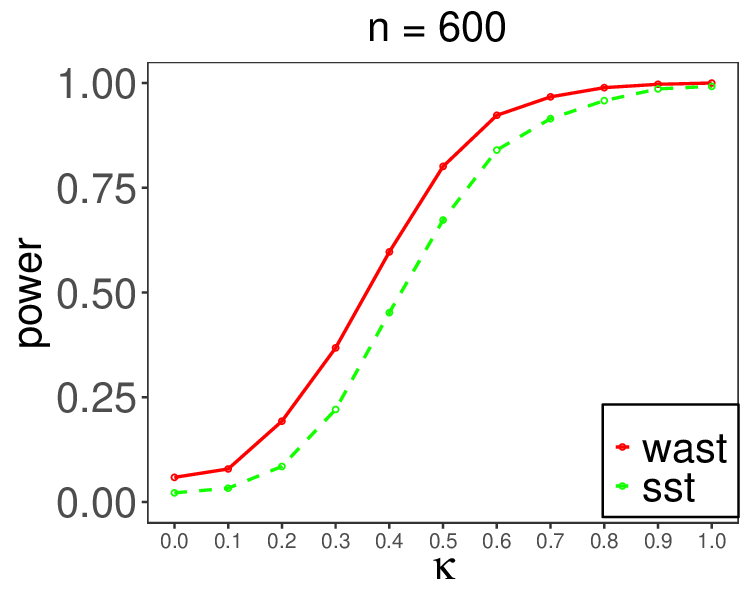}  \\
		\includegraphics[scale=0.33]{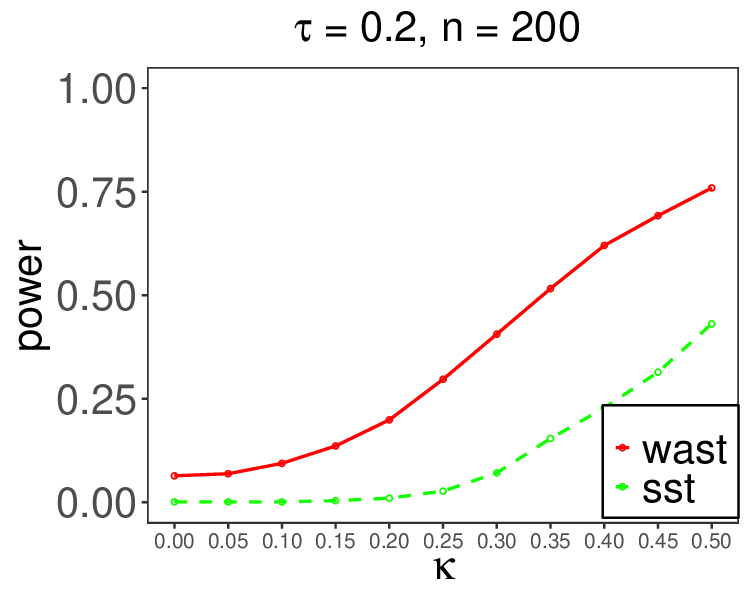}
		\includegraphics[scale=0.33]{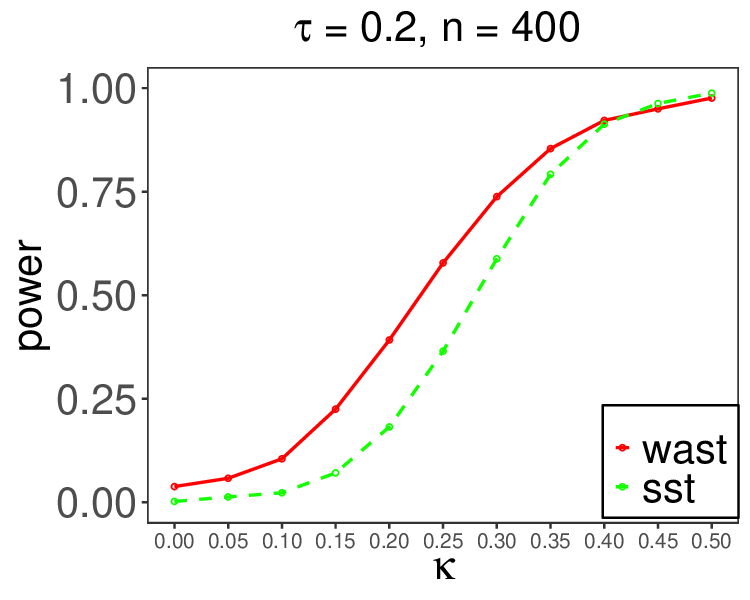}
		\includegraphics[scale=0.33]{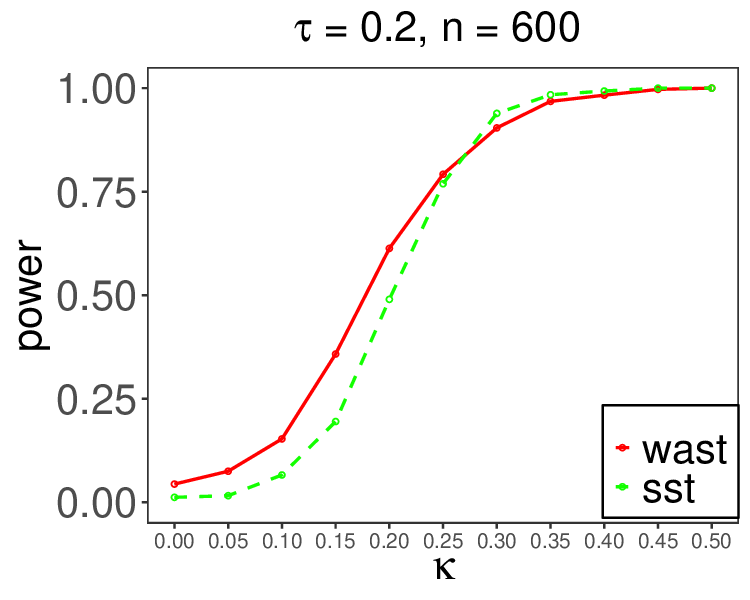}   \\
		\includegraphics[scale=0.33]{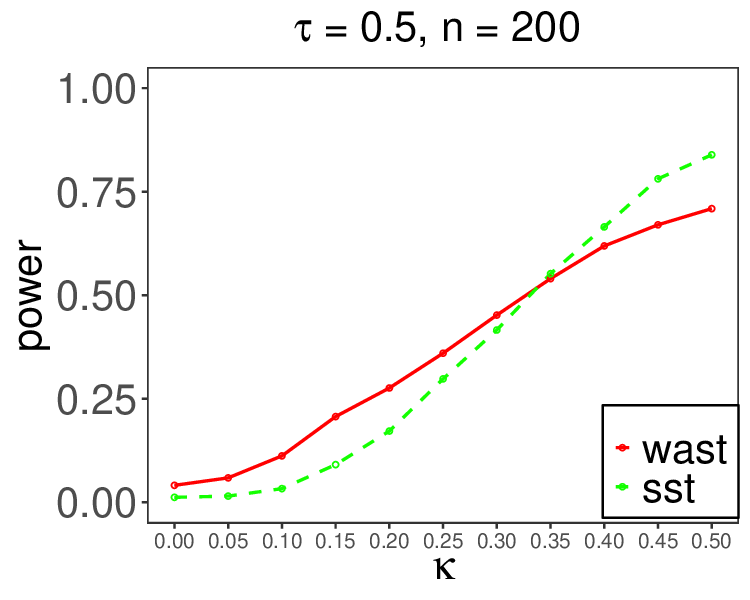}
		\includegraphics[scale=0.33]{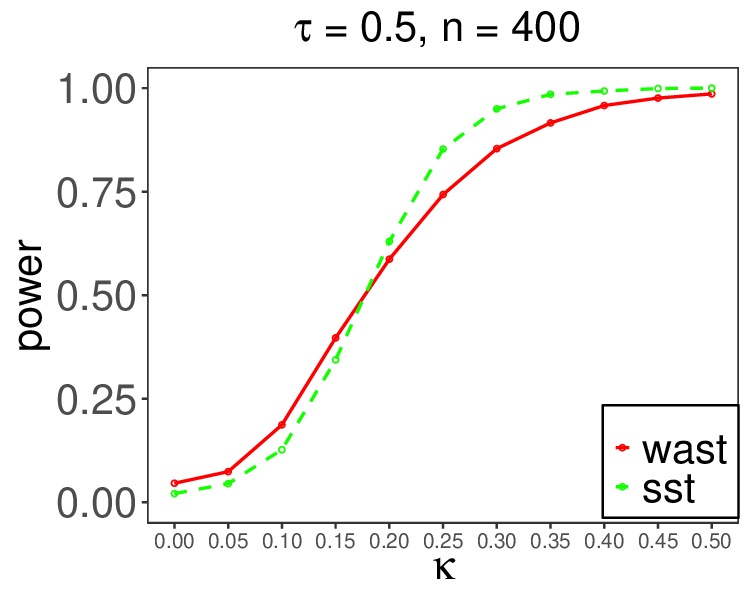}
		\includegraphics[scale=0.33]{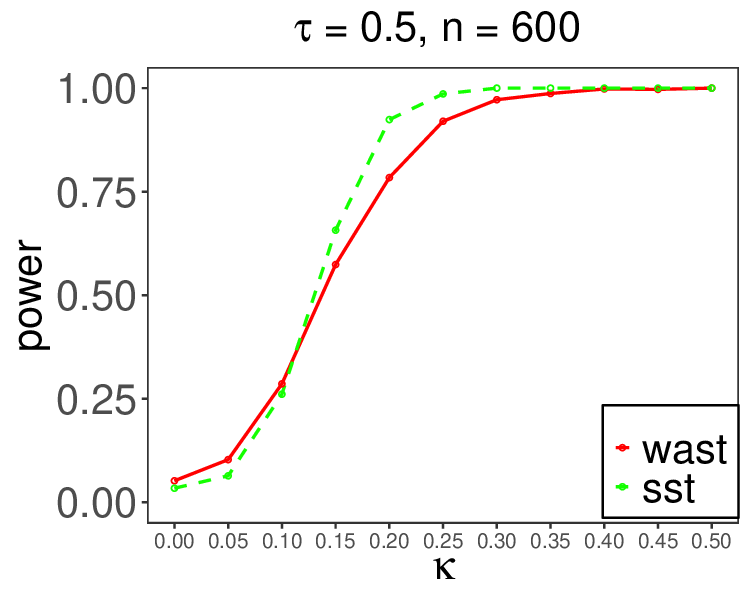}   \\
		\includegraphics[scale=0.33]{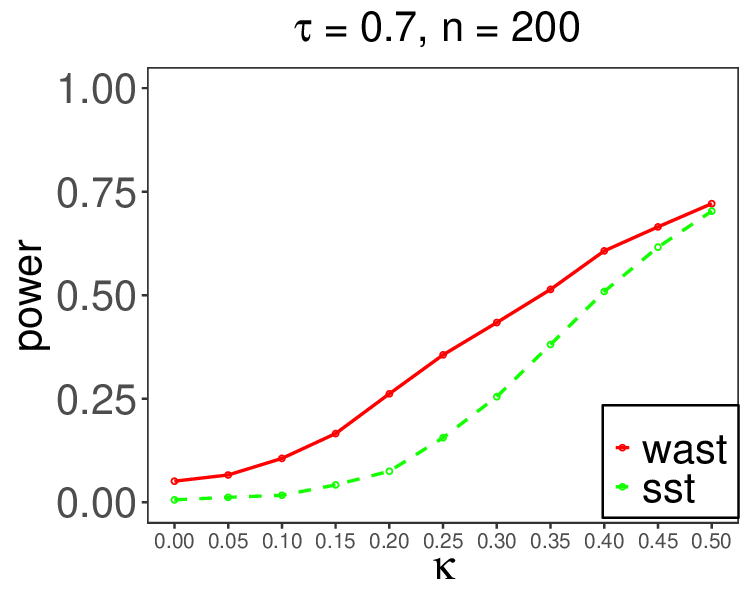}
		\includegraphics[scale=0.33]{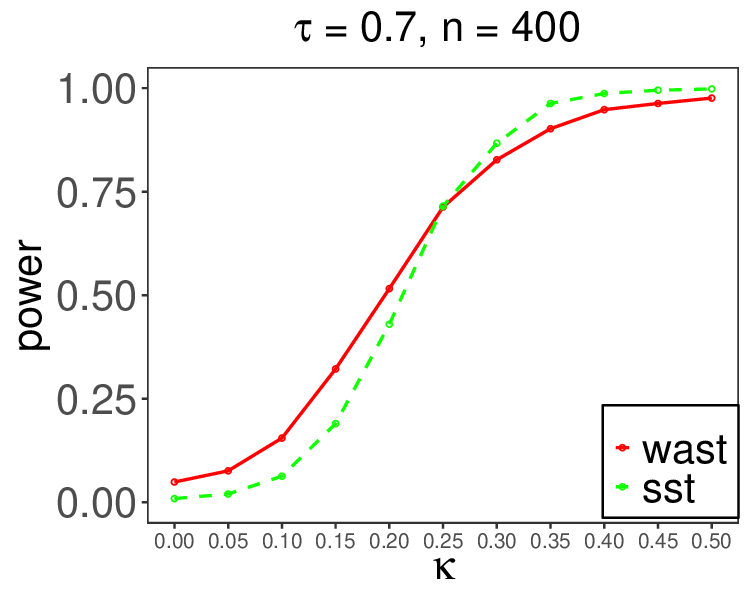}
		\includegraphics[scale=0.33]{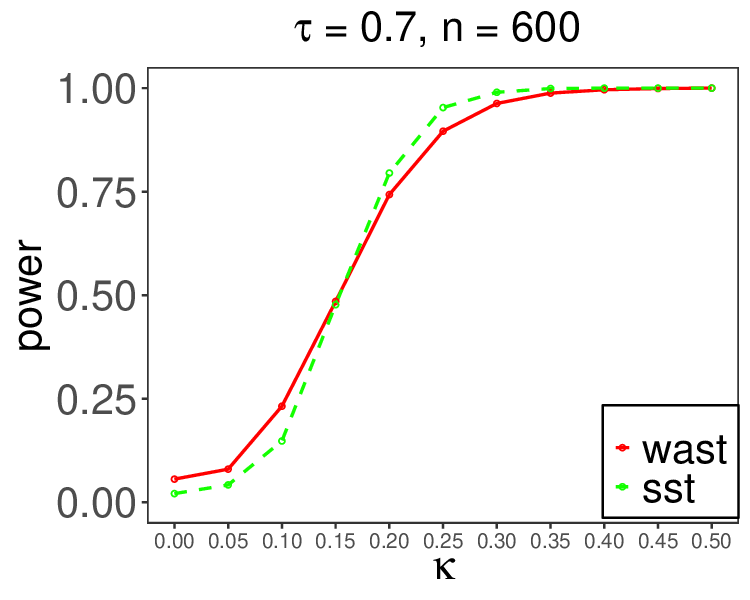}
		\caption{\it Powers of test statistic by the proposed WAST (red solid line) and SST (green dashed line) for $(p,q)=(5,5)$. From top to bottom, each row depicts the powers for probit model, semiparametric model, quantile regression with $\tau=0.2$, $\tau=0.5$ and $\tau=0.7$, respectively. Here the $\Gv$ $Z$ is generated from $t_3$ distribution.}
		\label{fig_qr55_zt3}
	\end{center}
\end{figure}

\begin{figure}[!ht]
	\begin{center}
		\includegraphics[scale=0.33]{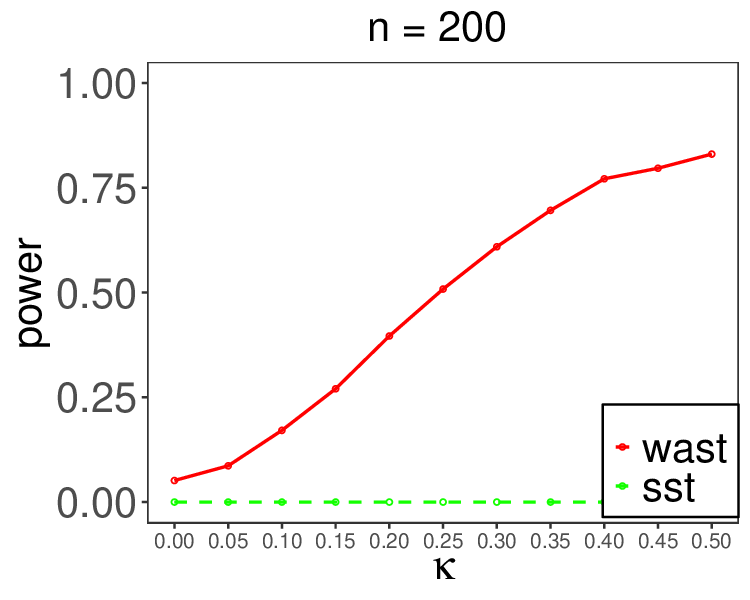}
		\includegraphics[scale=0.33]{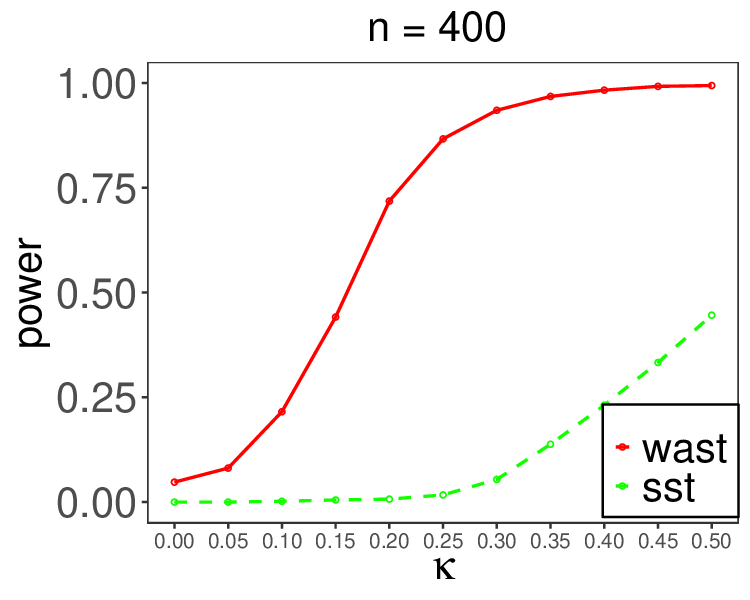}
		\includegraphics[scale=0.33]{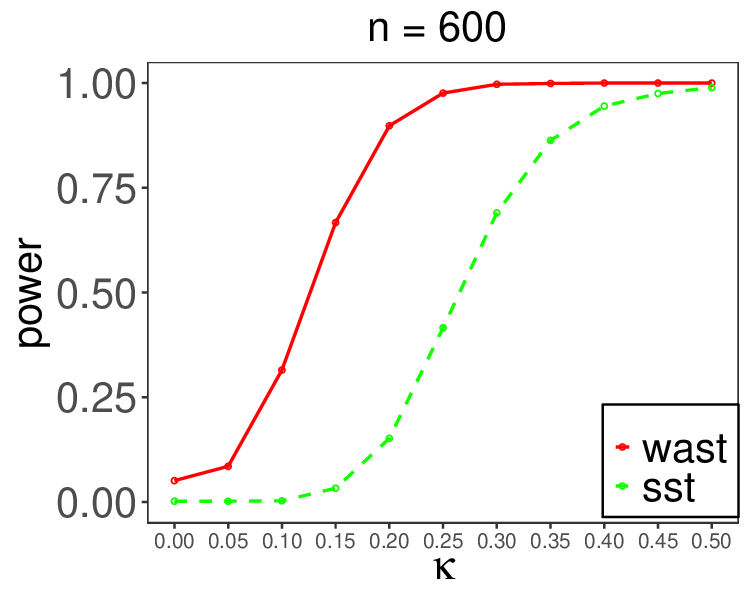}      \\
		\includegraphics[scale=0.33]{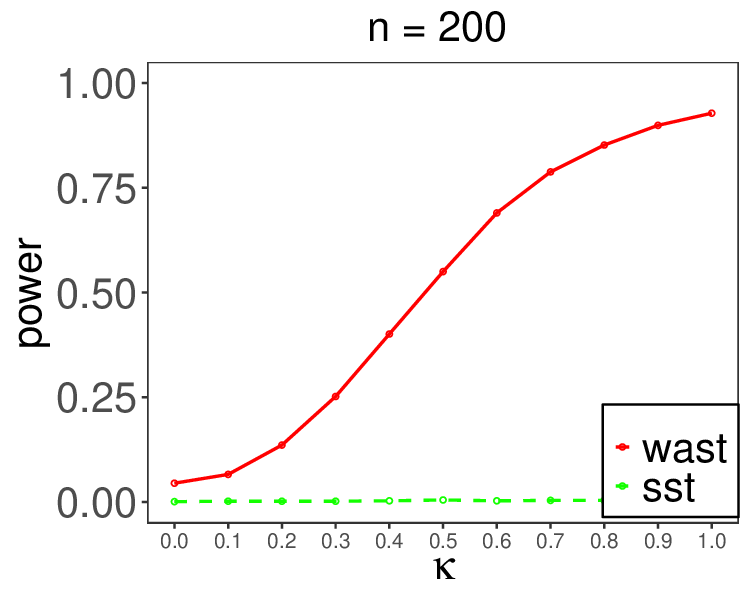}
		\includegraphics[scale=0.33]{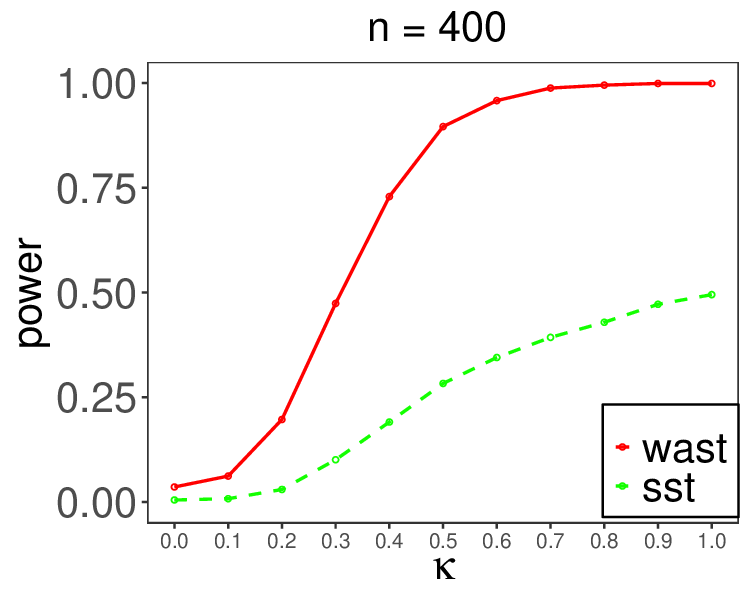}
		\includegraphics[scale=0.33]{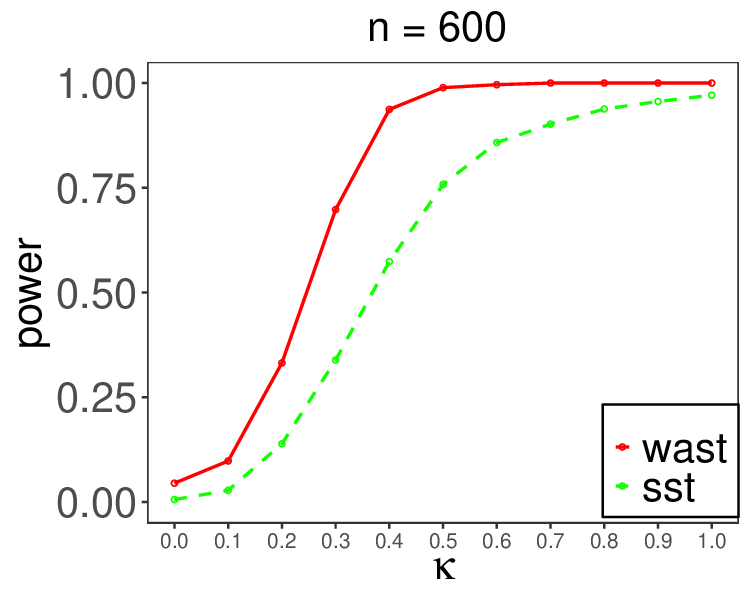}   \\
		\includegraphics[scale=0.33]{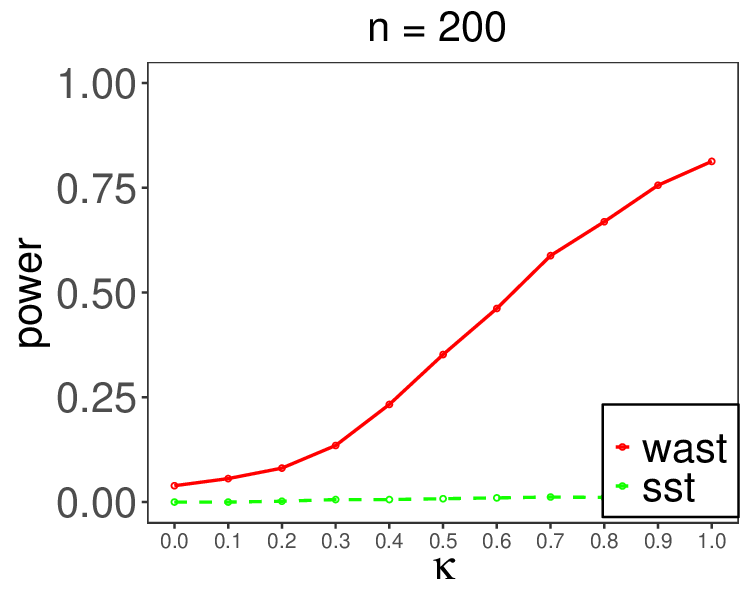}
		\includegraphics[scale=0.33]{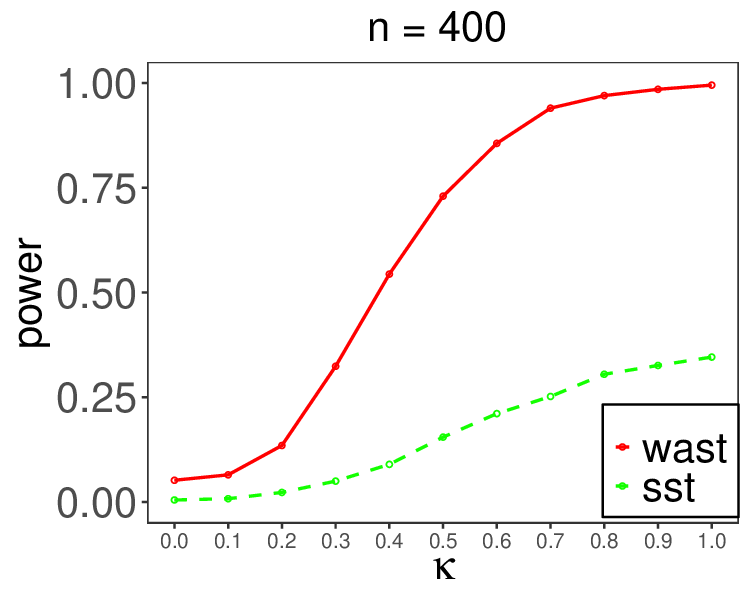}
		\includegraphics[scale=0.33]{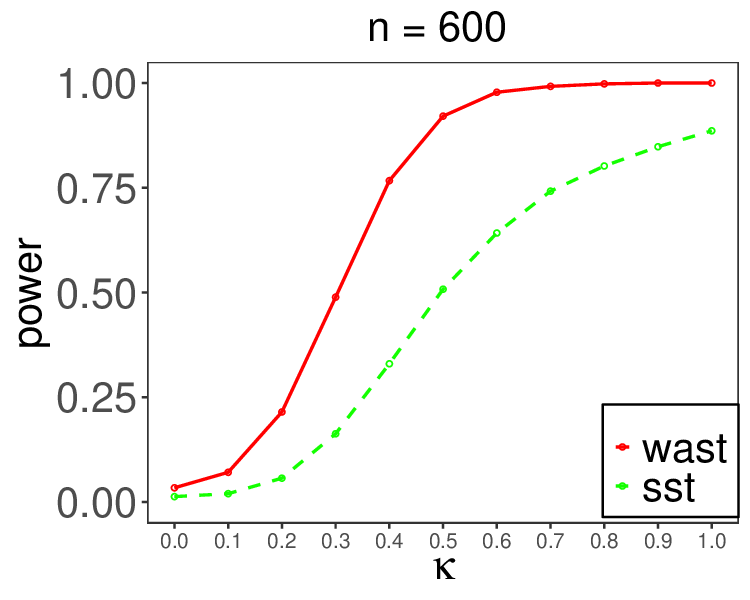}  \\
		\includegraphics[scale=0.33]{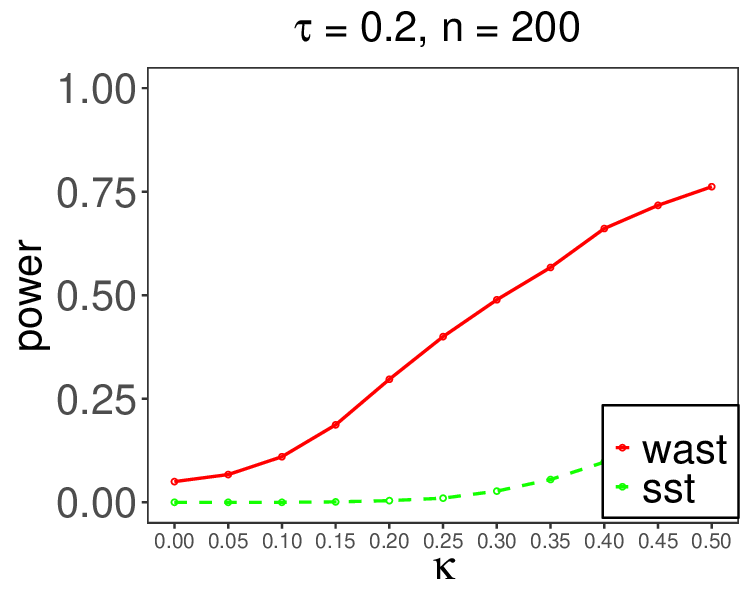}
		\includegraphics[scale=0.33]{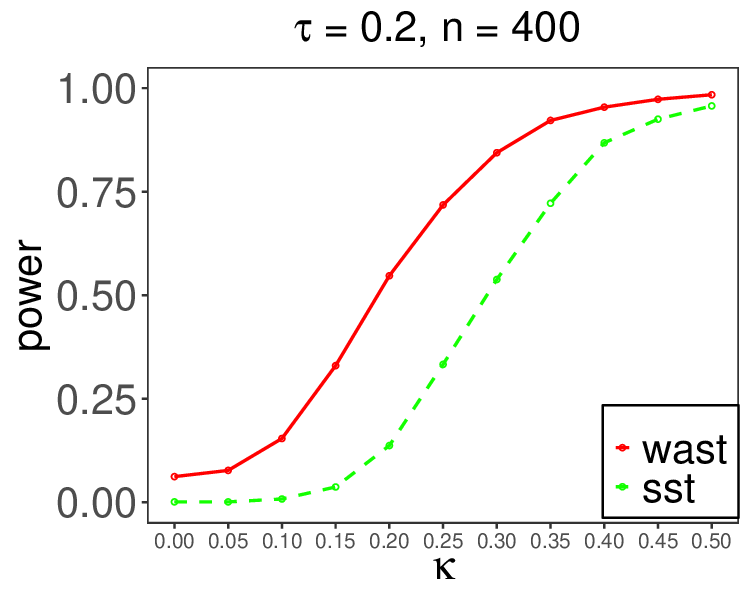}
		\includegraphics[scale=0.33]{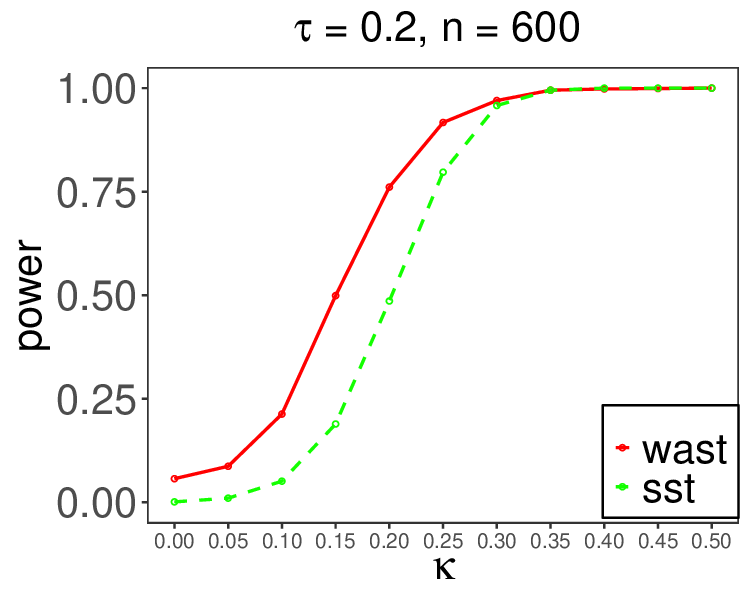}   \\
		\includegraphics[scale=0.33]{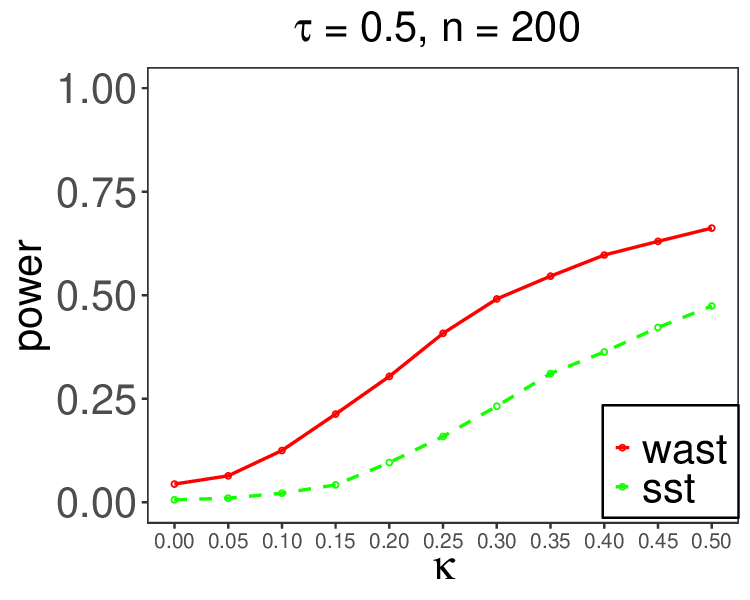}
		\includegraphics[scale=0.33]{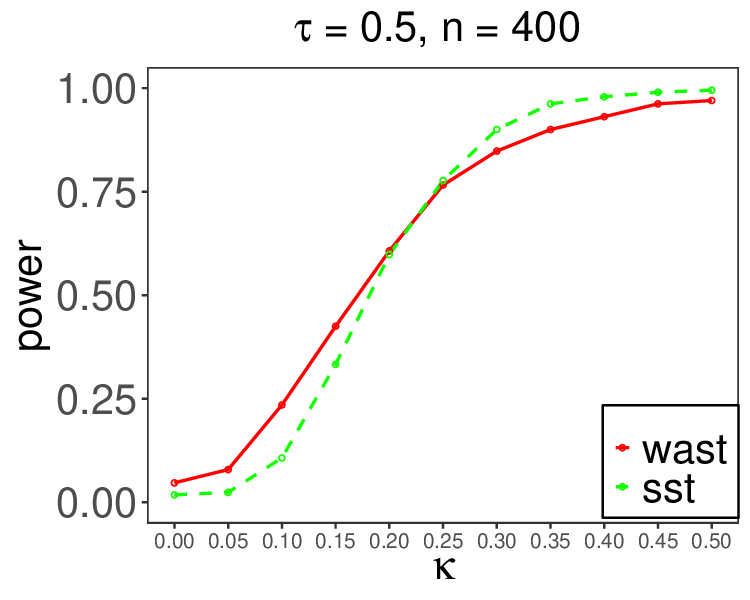}
		\includegraphics[scale=0.33]{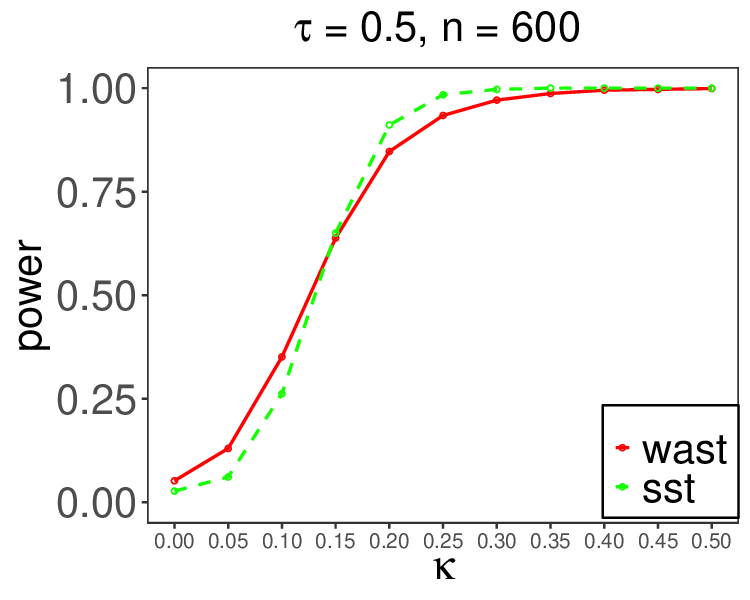}   \\
		\includegraphics[scale=0.33]{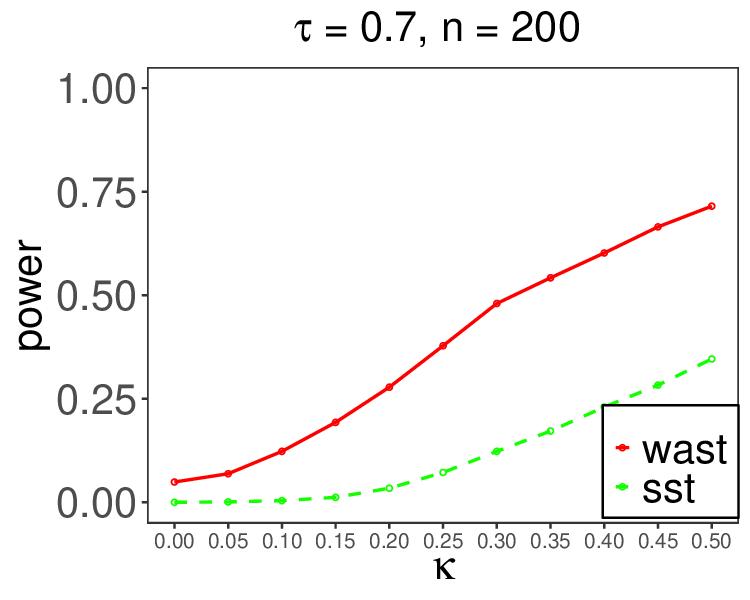}
		\includegraphics[scale=0.33]{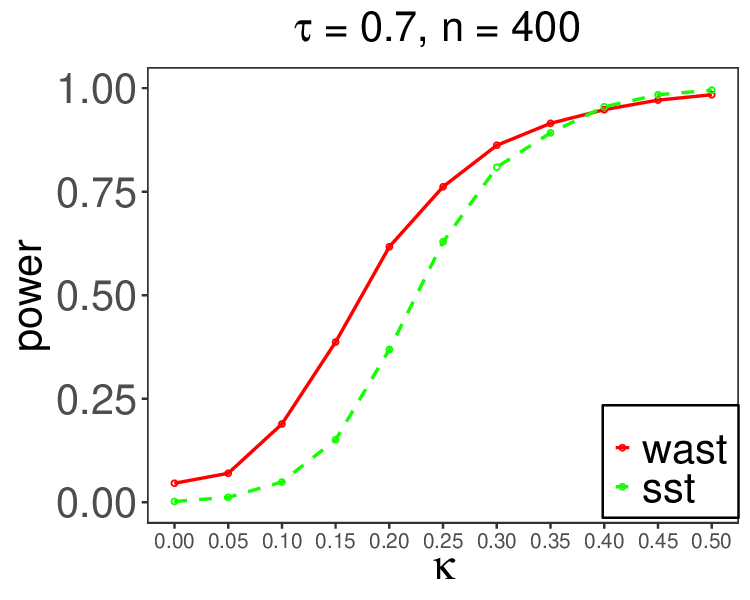}
		\includegraphics[scale=0.33]{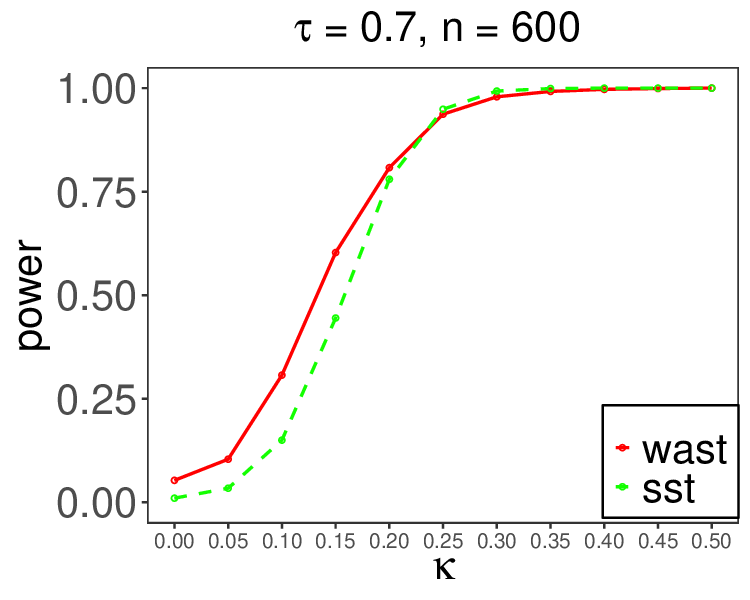}
		\caption{\it Powers of test statistic by the proposed WAST (red solid line) and SST (green dashed line) for $(p,q)=(10,10)$. From top to bottom, each row depicts the powers for probit model, semiparametric model, quantile regression with $\tau=0.2$, $\tau=0.5$ and $\tau=0.7$, respectively. Here the $\Gv$ $Z$ is generated from $t_3$ distribution.}
		\label{fig_qr1010_zt3}
	\end{center}
\end{figure}

\begin{figure}[!ht]
	\begin{center}
		\includegraphics[scale=0.33]{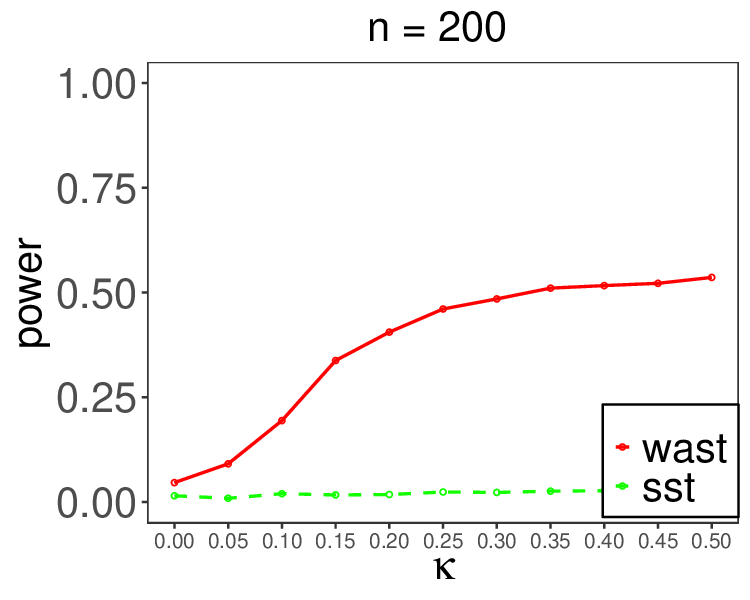}
		\includegraphics[scale=0.33]{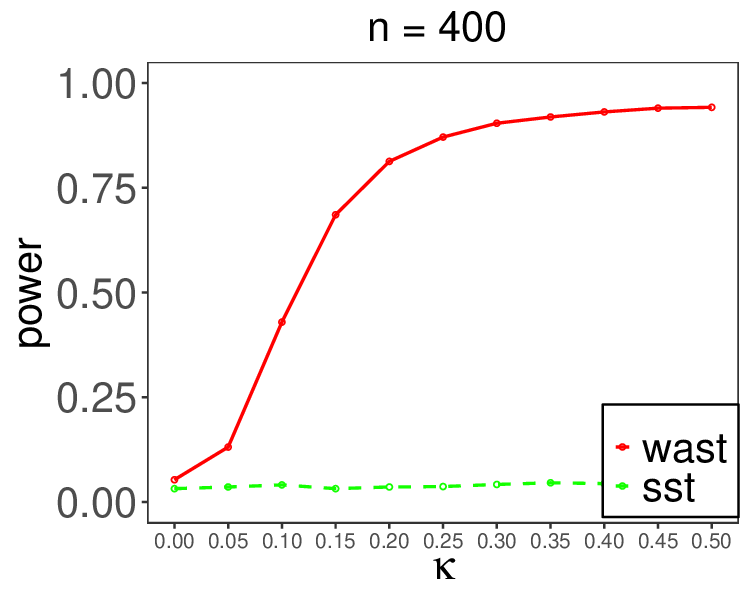}
		\includegraphics[scale=0.33]{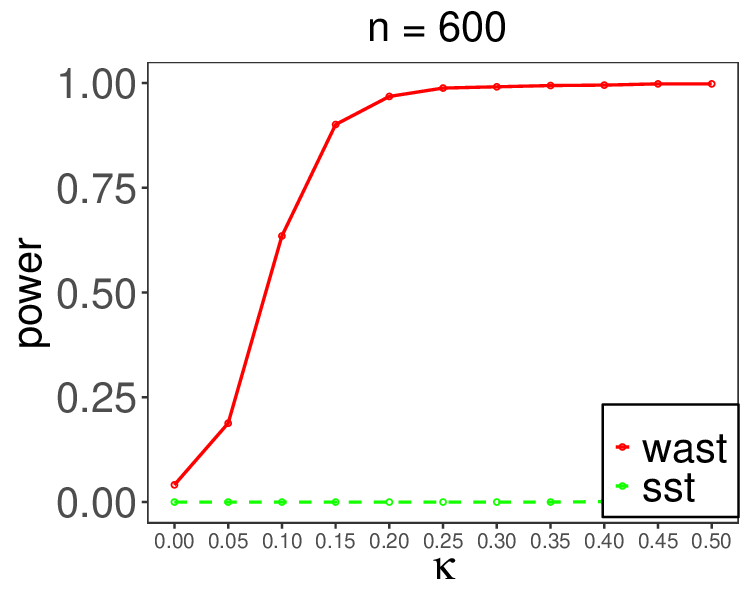}     \\
		\includegraphics[scale=0.33]{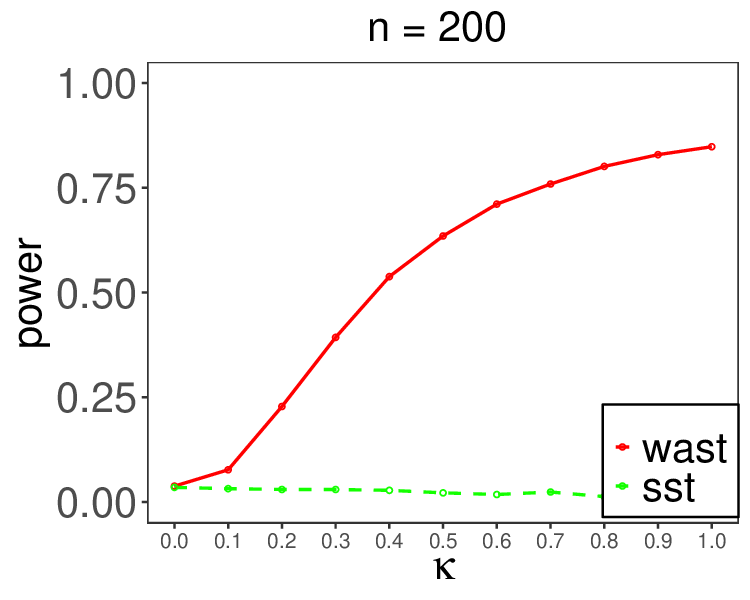}
		\includegraphics[scale=0.33]{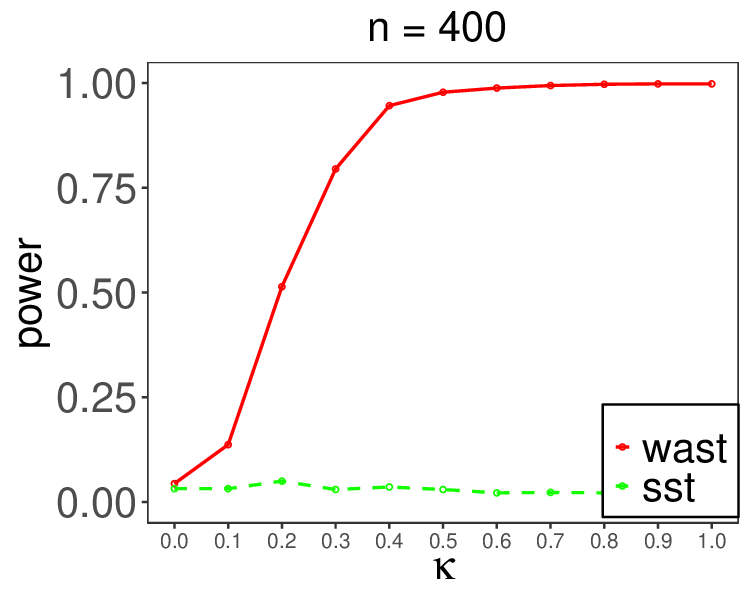}
		\includegraphics[scale=0.33]{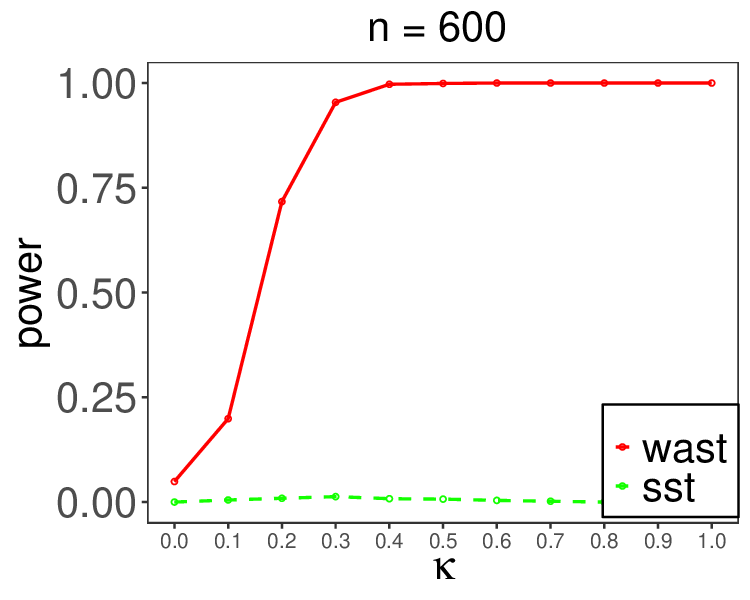}  \\
		\includegraphics[scale=0.33]{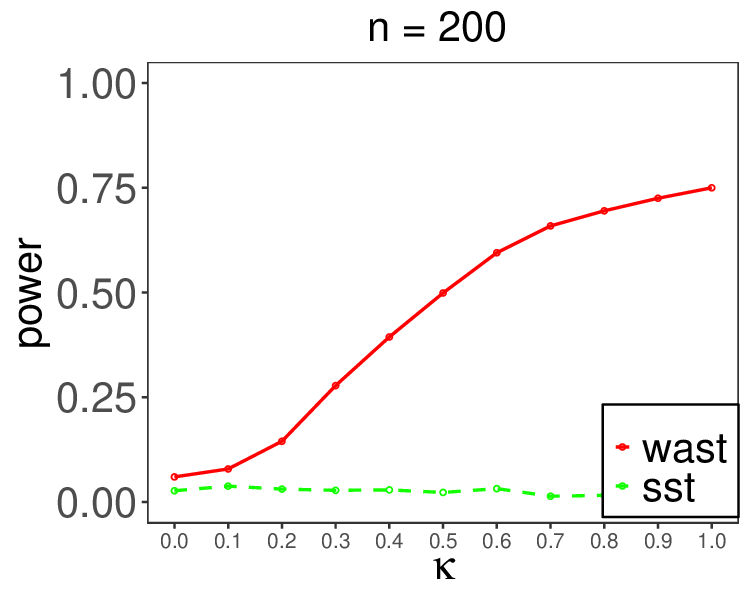}
		\includegraphics[scale=0.33]{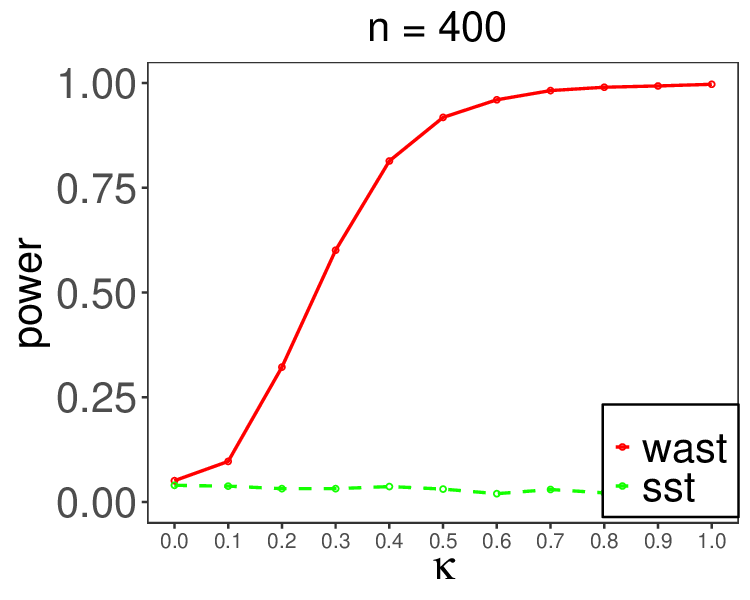}
		\includegraphics[scale=0.33]{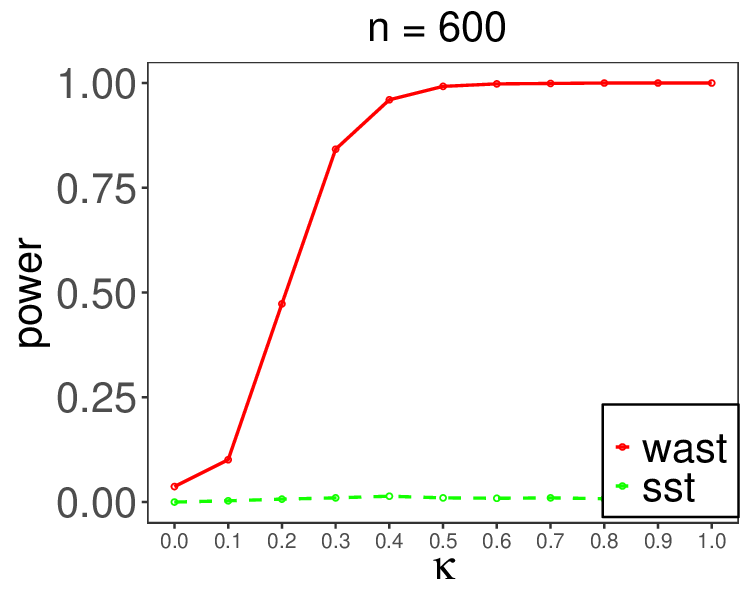} \\
		\includegraphics[scale=0.33]{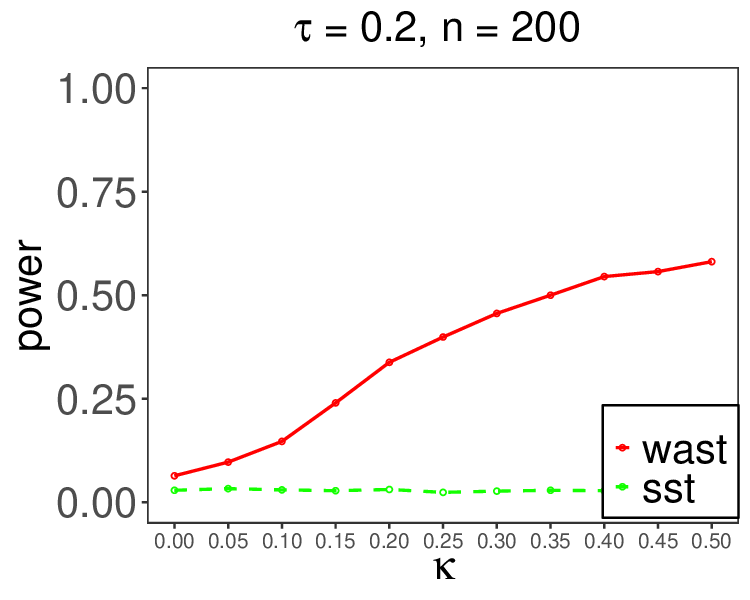}
		\includegraphics[scale=0.33]{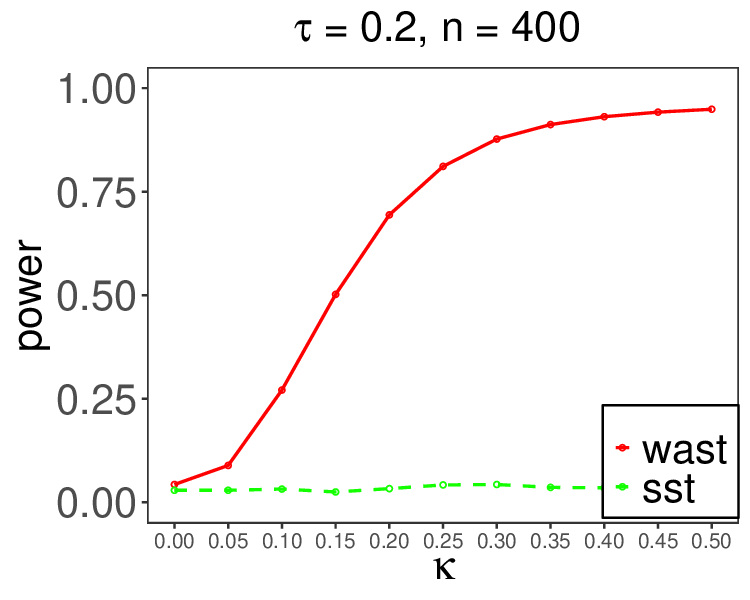}
		\includegraphics[scale=0.33]{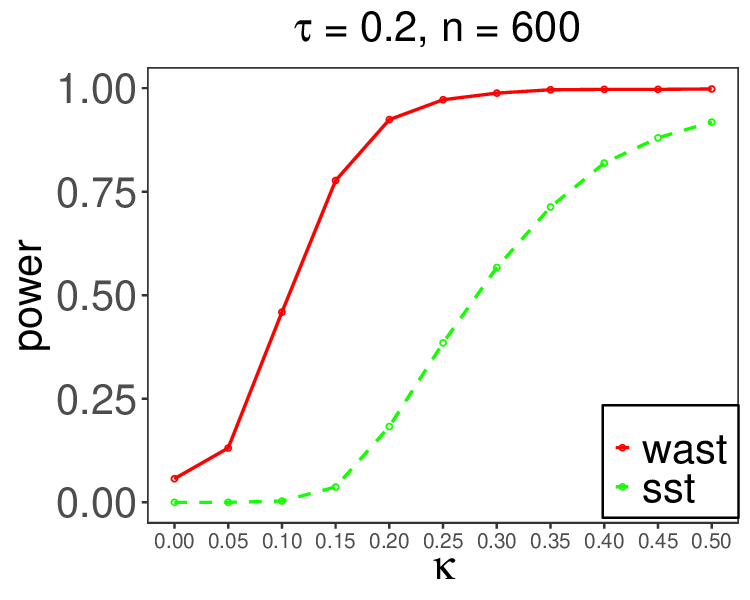}  \\
		\includegraphics[scale=0.33]{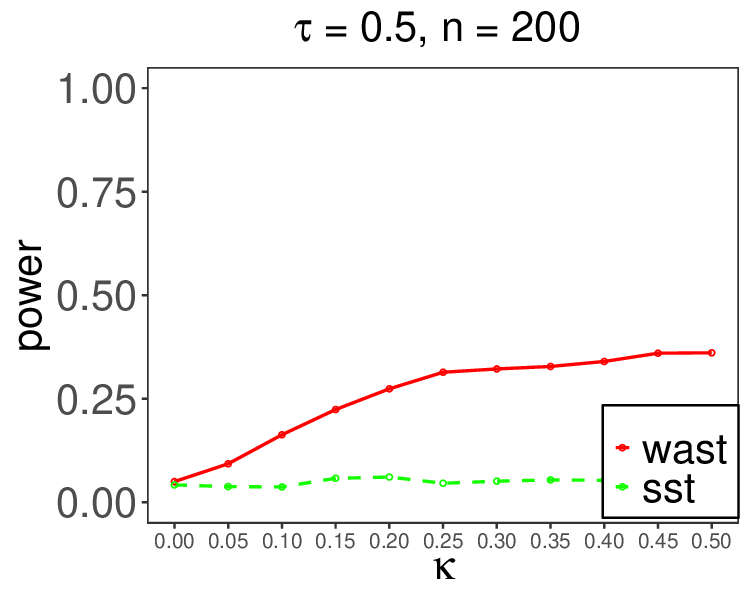}
		\includegraphics[scale=0.33]{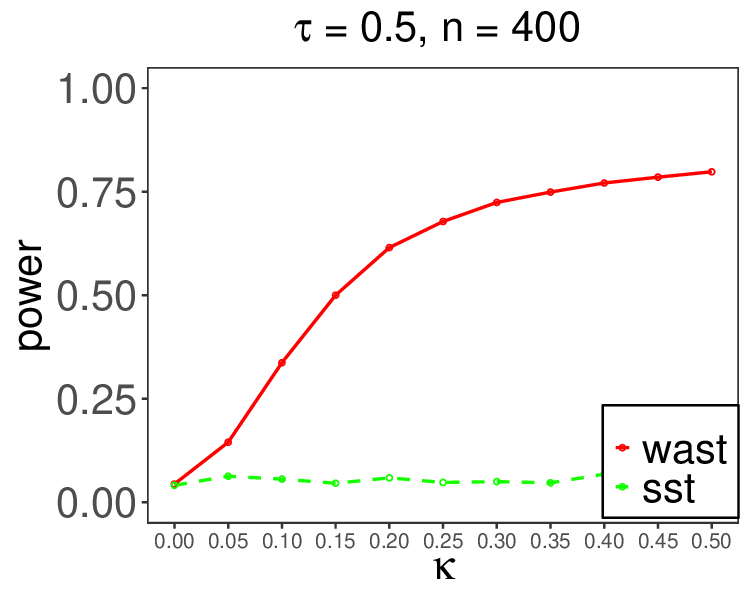}
		\includegraphics[scale=0.33]{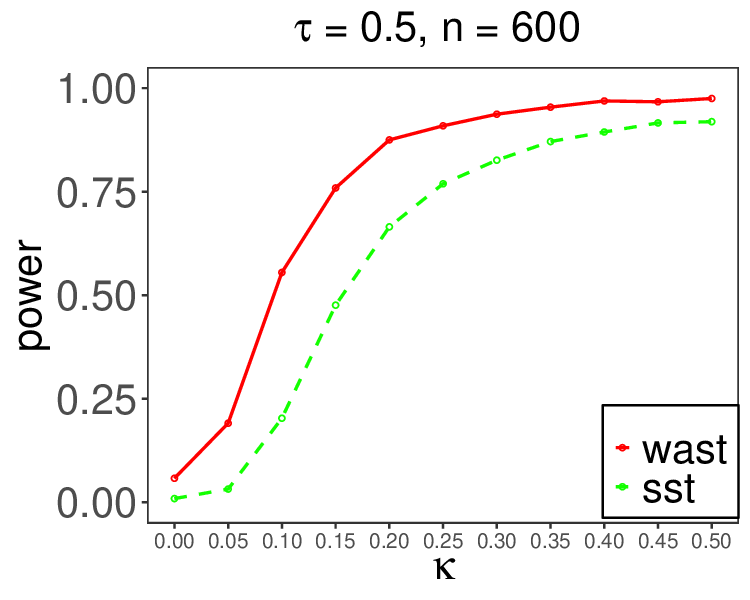}  \\
		\includegraphics[scale=0.33]{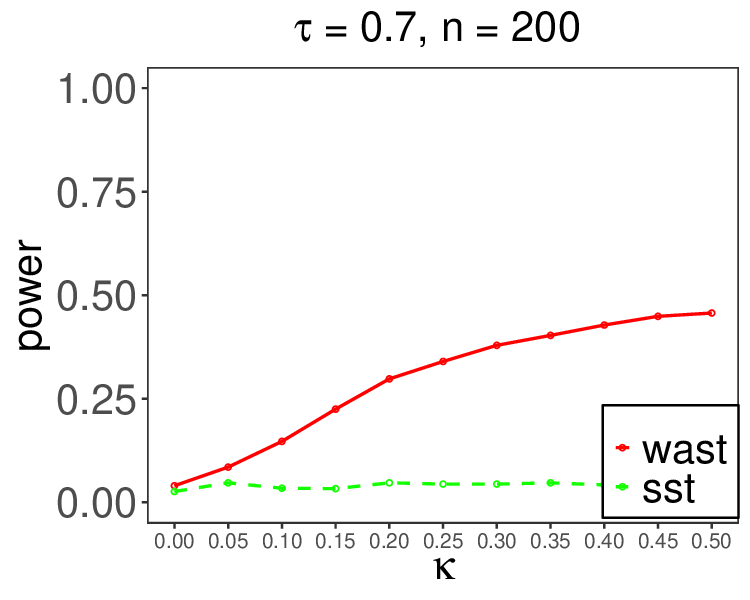}
		\includegraphics[scale=0.33]{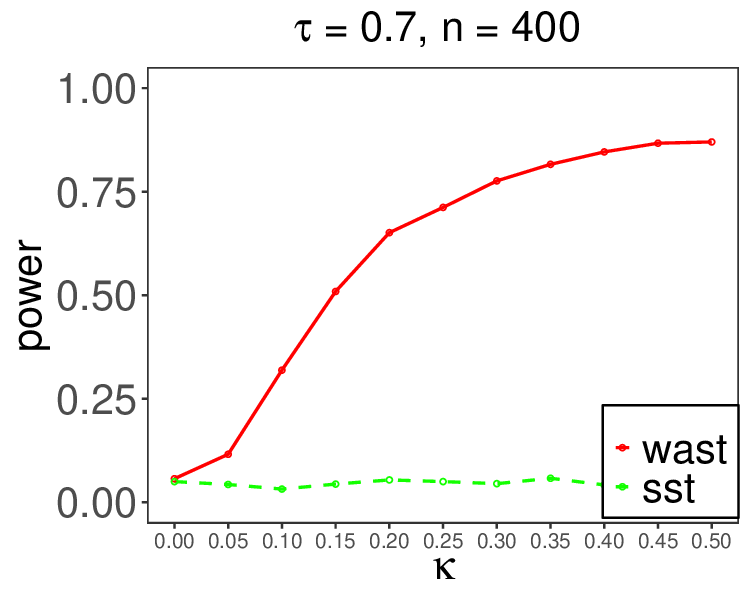}
		\includegraphics[scale=0.33]{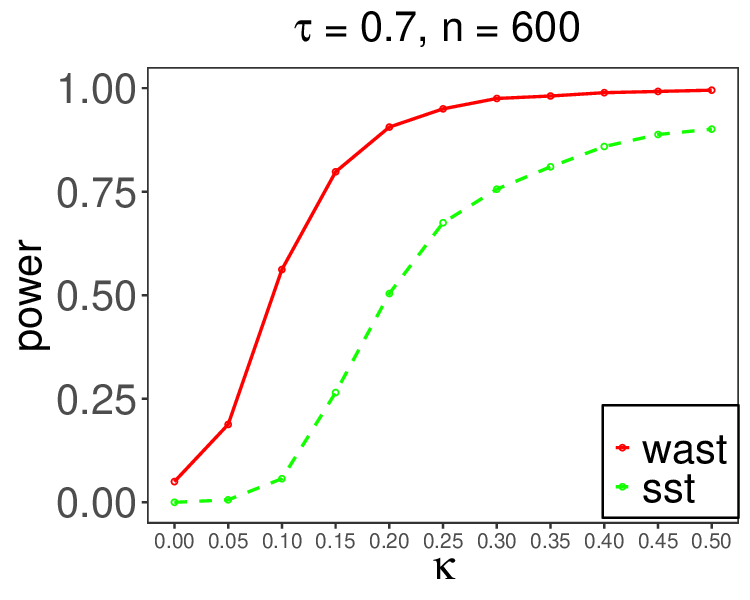}
		\caption{\it Powers of test statistic by the proposed WAST with (red solid line) and SST (green dashed line) for $(p,q)=(50,20)$. From top to bottom, each row depicts the powers for probit model, semiparametric model, quantile regression with $\tau=0.2$, $\tau=0.5$ and $\tau=0.7$, respectively. Here the $\Gv$ $Z$ is generated from $t_3$ distribution.}
		\label{fig_qr5020_zt3}
	\end{center}
\end{figure}

Type \uppercase\expandafter{\romannumeral1} errors of the proposed test statistic are listed in Table \ref{table_size_zt3}, and the power curves are depicted in Figure \ref{fig_qr13_zt3}-\ref{fig_qr5020_zt3}. We have the similar conclusion from Table \ref{table_size_zt3} and Figure \ref{fig_qr13_zt3}-\ref{fig_qr5020_zt3}.

\begin{table}[htp!]
	\def~{\hphantom{0}}
\tiny
	\caption{Type \uppercase\expandafter{\romannumeral1} errors of the proposed WAST and SST based on resampling for probit regression model (ProbitRE), quantile regression (QuantRE) and semiparametric model (SPMoldel).
}
	\resizebox{\textwidth}{!}{
        \begin{threeparttable}
		\begin{tabular}{llcccccccc}
			\hline
			\multirow{2}{*}{Model}&\multirow{2}{*}{$(p,q)$}
			&\multicolumn{2}{c}{ $n=200$} && \multicolumn{2}{c}{ $n=400$} && \multicolumn{2}{c}{ $n=600$} \\
			\cline{3-4} \cline{6-7} \cline{9-10}
			&&   WAST & SST && WAST & SST && WAST & SST \\
			\cline{3-10}
			ProbitRE &$(1,3)$         & 0.057  & 0.001 && 0.056  & 0.014 && 0.070  & 0.030 \\
			&$(5,5)$                  & 0.051  & 0.000 && 0.050  & 0.000 && 0.055  & 0.003 \\
			& $(10,10)$               & 0.051  & 0.000 && 0.048  & 0.000 && 0.051  & 0.002 \\
			&$(50,20)$                & 0.046  & 0.015 && 0.053  & 0.032 && 0.041  & 0.000 \\
			[1 ex]
			SPModel &$(1,3)$          & 0.052  & 0.026 && 0.056  & 0.047 && 0.047  & 0.050 \\
			(B1+P1)&$(5,5)$           & 0.045  & 0.004 && 0.038  & 0.020 && 0.031  & 0.021 \\
			& $(10,10)$               & 0.045  & 0.001 && 0.036  & 0.005 && 0.045  & 0.006 \\
			&$(50,20)$                & 0.038  & 0.035 && 0.044  & 0.032 && 0.049  & 0.000 \\
			[1 ex]
			SPModel &$(1,3)$          & 0.057  & 0.043 && 0.051  & 0.050 && 0.047  & 0.044 \\
			(B2+P2)&$(5,5)$           & 0.058  & 0.011 && 0.046  & 0.012 && 0.059  & 0.022 \\
			& $(10,10)$               & 0.039  & 0.000 && 0.052  & 0.005 && 0.034  & 0.013 \\
			&$(50,20)$                & 0.060  & 0.027 && 0.051  & 0.040 && 0.037  & 0.000 \\
			[1 ex]
			QuantRE &$(1,3)$          & 0.046  & 0.014 && 0.048  & 0.026 && 0.058  & 0.038 \\
			($\tau=0.2$)&$(5,5)$      & 0.064  & 0.001 && 0.038  & 0.002 && 0.044  & 0.012 \\
			& $(10,10)$               & 0.050  & 0.000 && 0.062  & 0.001 && 0.057  & 0.001 \\
			&$(50,20)$                & 0.064  & 0.029 && 0.043  & 0.029 && 0.057  & 0.000 \\
			[1 ex]
			QuantRE &$(1,3)$          & 0.040  & 0.031 && 0.053  & 0.050 && 0.045  & 0.047 \\
			($\tau=0.5$)&$(5,5)$      & 0.041  & 0.012 && 0.046  & 0.021 && 0.052  & 0.034 \\
			& $(10,10)$               & 0.044  & 0.006 && 0.047  & 0.018 && 0.052  & 0.027 \\
			&$(50,20)$                & 0.050  & 0.042 && 0.044  & 0.041 && 0.058  & 0.009 \\
			[1 ex]
			QuantRE &$(1,3)$          & 0.039  & 0.023 && 0.043  & 0.037 && 0.060  & 0.042 \\
			($\tau=0.7$)&$(5,5)$      & 0.051  & 0.006 && 0.049  & 0.009 && 0.056  & 0.021 \\
			& $(10,10)$               & 0.049  & 0.000 && 0.046  & 0.002 && 0.053  & 0.010 \\
			&$(50,20)$                & 0.040  & 0.026 && 0.057  & 0.050 && 0.050  & 0.000 \\
			\hline
		\end{tabular}
\begin{tablenotes}
\item The nominal significant level is 0.05. Here the $\Gv$ $\bZ$ is generated from $t_3$ distribution.
\end{tablenotes}
\end{threeparttable}
	}
	\label{table_size_zt3}
\end{table}

\subsection{Quantile regression with Cauchy distribution}\label{simulation_cpz_cauchy}
In this Section, we consider the quantile regression with error from standard Cauchy distribution. Other settings are same as Section \ref{simulation_cpzt3} in current Supplementary Material.

\begin{figure}[!ht]
	\begin{center}
		\includegraphics[scale=0.33]{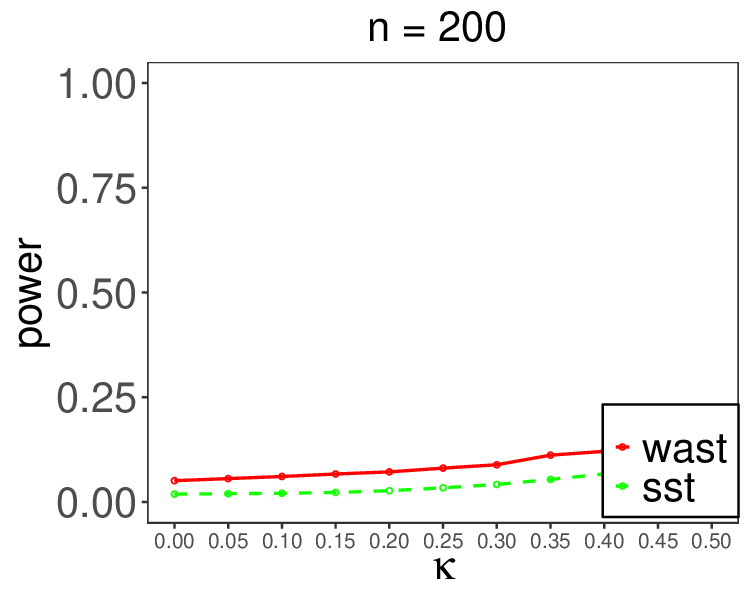}
		\includegraphics[scale=0.33]{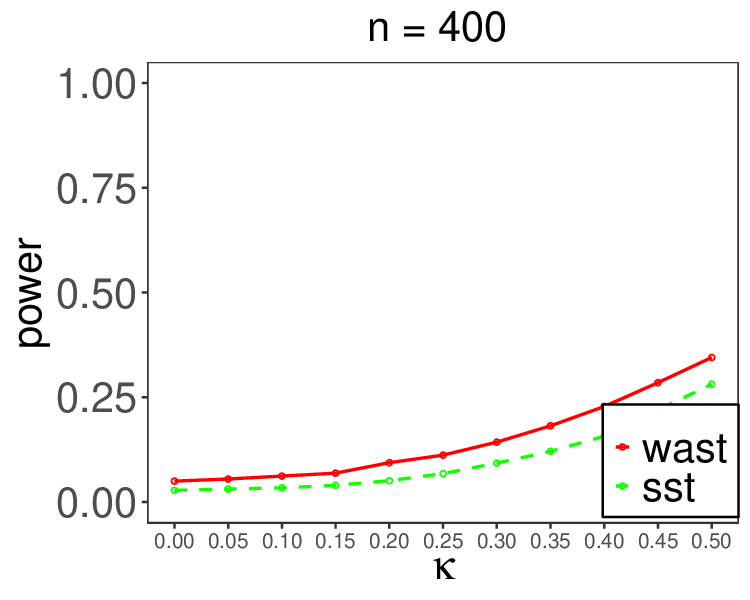}
		\includegraphics[scale=0.33]{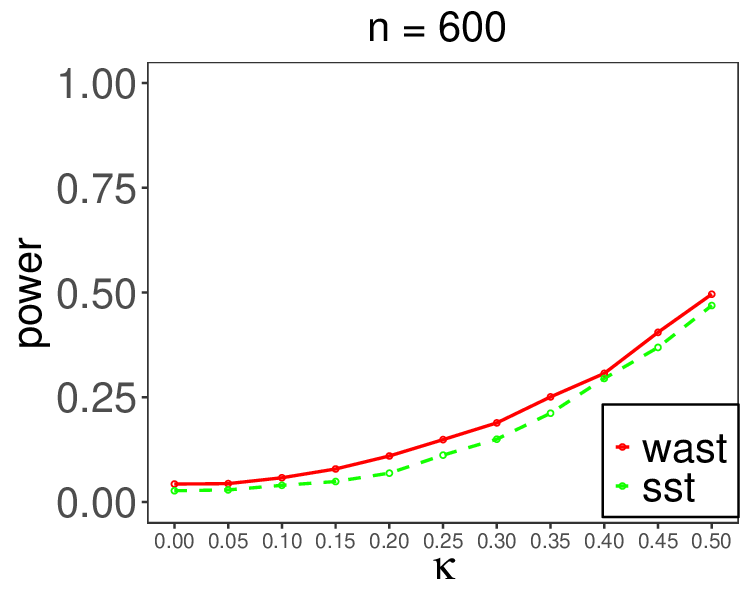}   \\
		\includegraphics[scale=0.33]{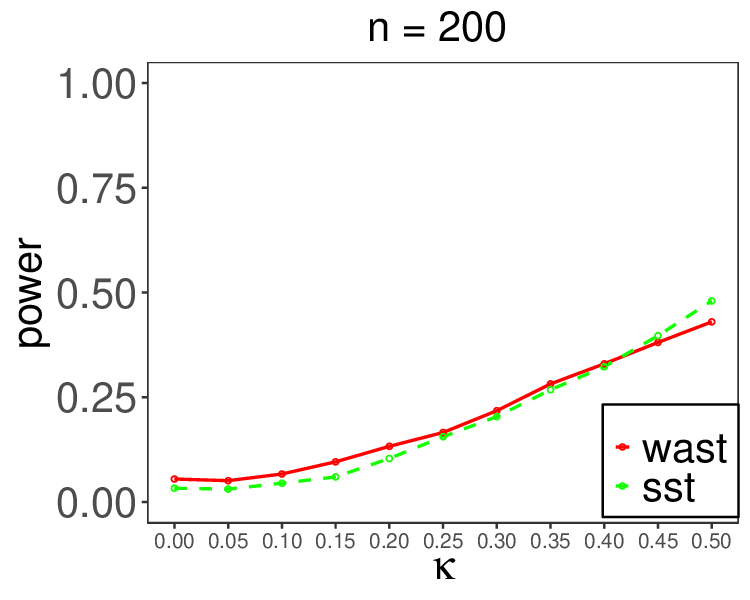}
		\includegraphics[scale=0.33]{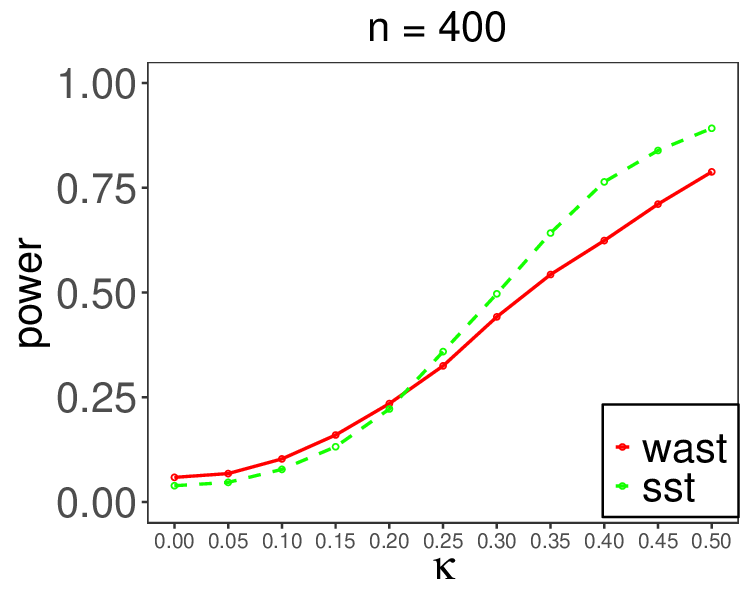}
		\includegraphics[scale=0.33]{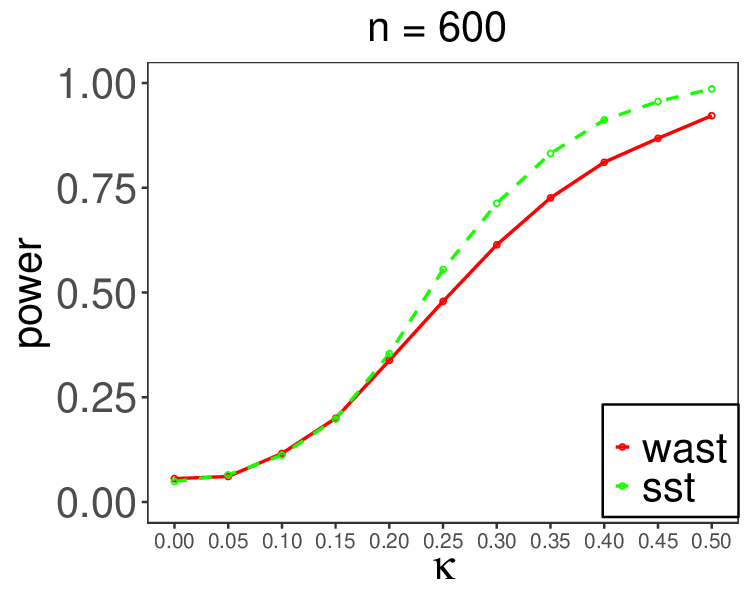}   \\
		\includegraphics[scale=0.33]{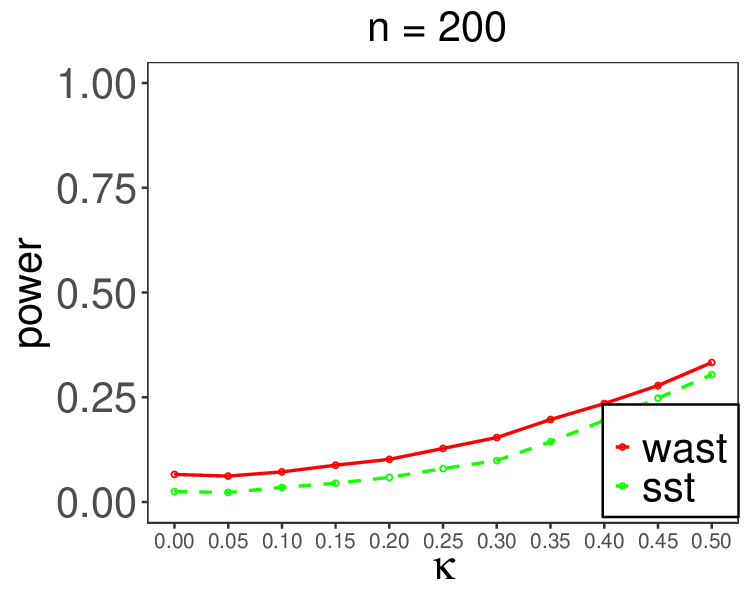}
		\includegraphics[scale=0.33]{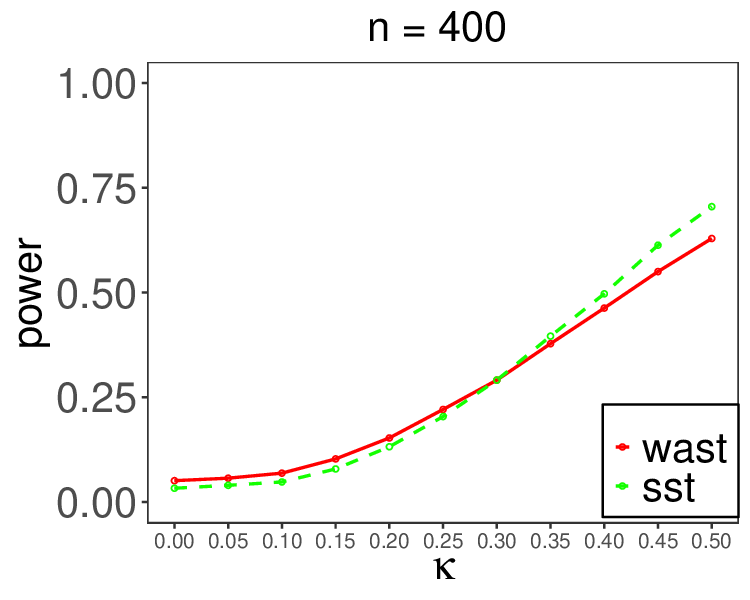}
		\includegraphics[scale=0.33]{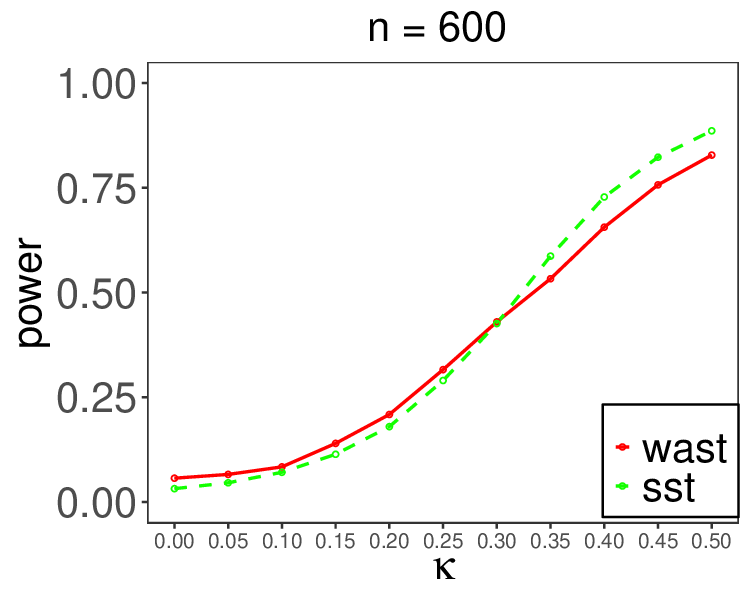}
		\caption{\it Powers of test statistic by the proposed WAST (red solid line) and SST (green dashed line) for $(p,q)=(1,3)$. From top to bottom, each row depicts the powers for quantile regression with $\tau=0.2$, $\tau=0.5$ and $\tau=0.7$, respectively. Here the $\Gv$ $Z$ is generated from $t_3$ distribution, and error from Cauchy distribution.}
		\label{fig_qr13_cauchy}
	\end{center}
\end{figure}

\begin{figure}[!ht]
	\begin{center}
		\includegraphics[scale=0.33]{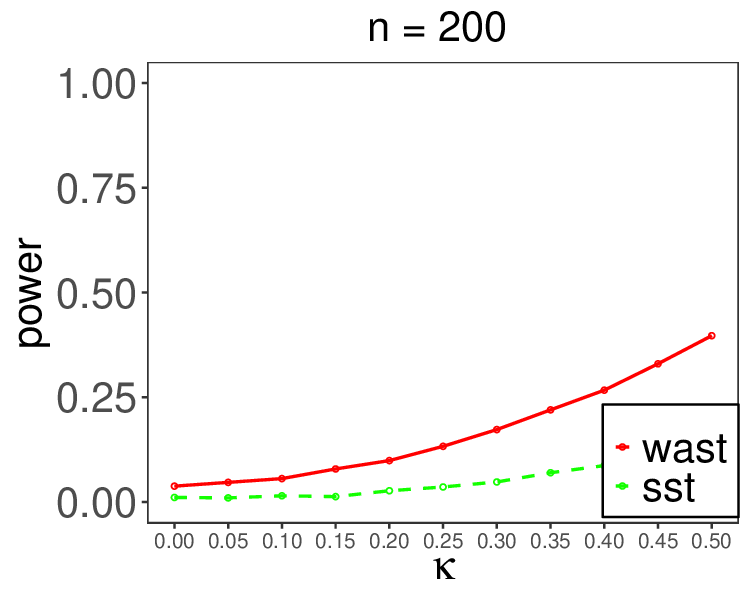}
		\includegraphics[scale=0.33]{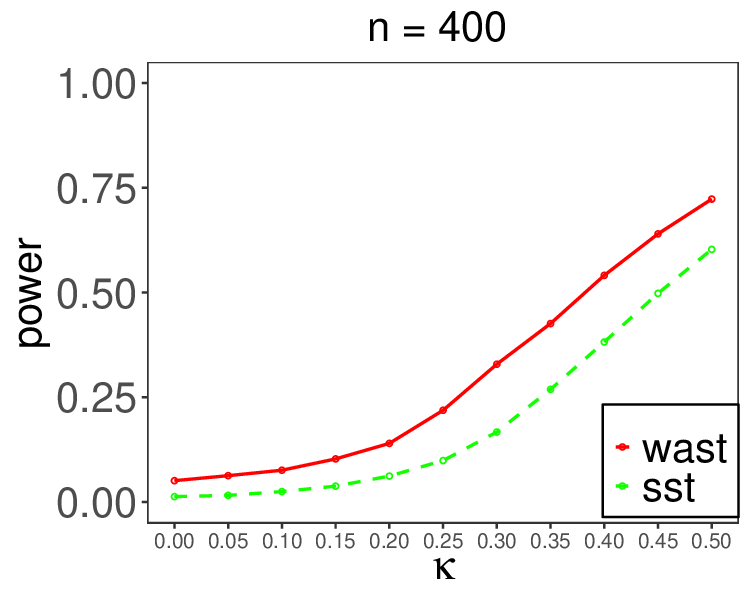}
		\includegraphics[scale=0.33]{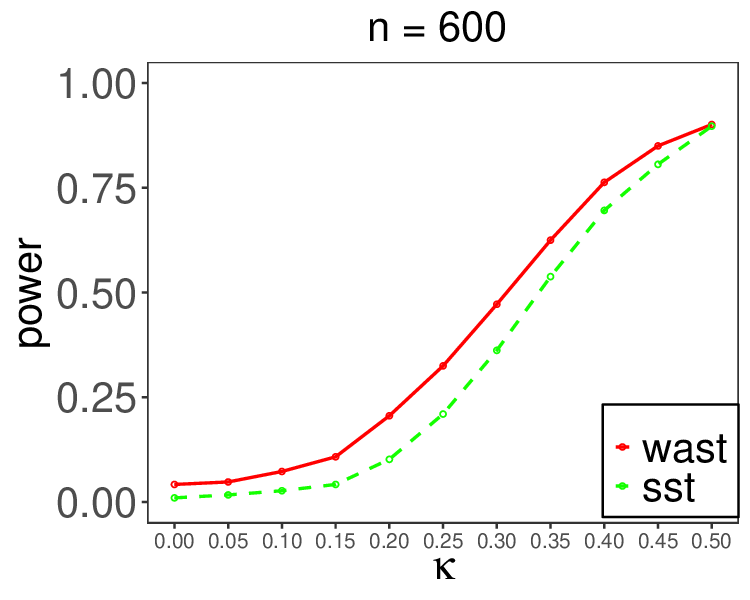}   \\
		\includegraphics[scale=0.33]{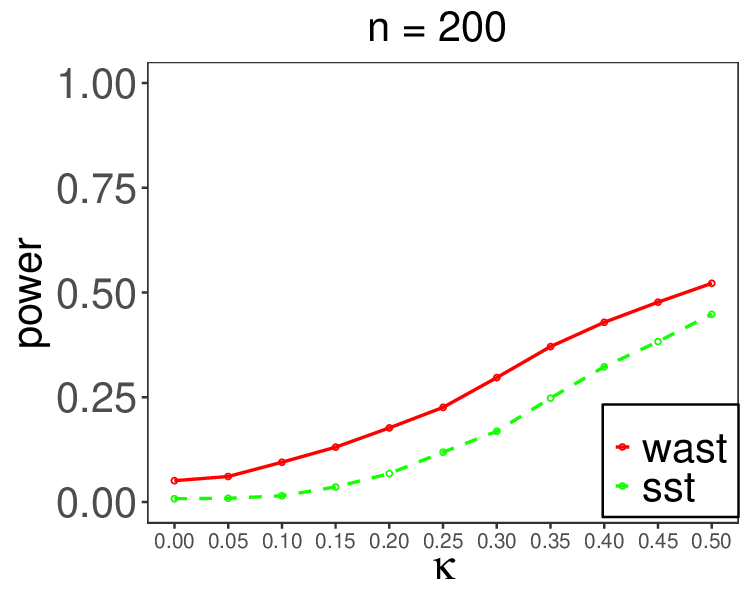}
		\includegraphics[scale=0.33]{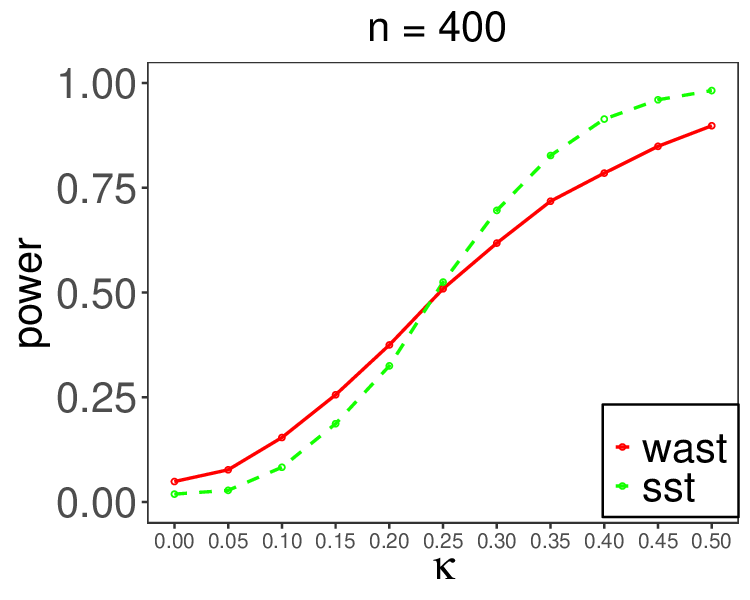}
		\includegraphics[scale=0.33]{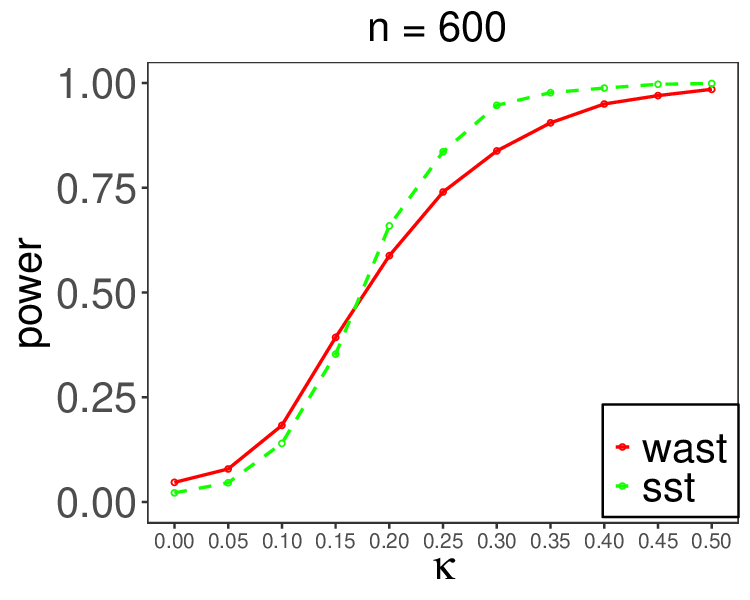}   \\
		\includegraphics[scale=0.33]{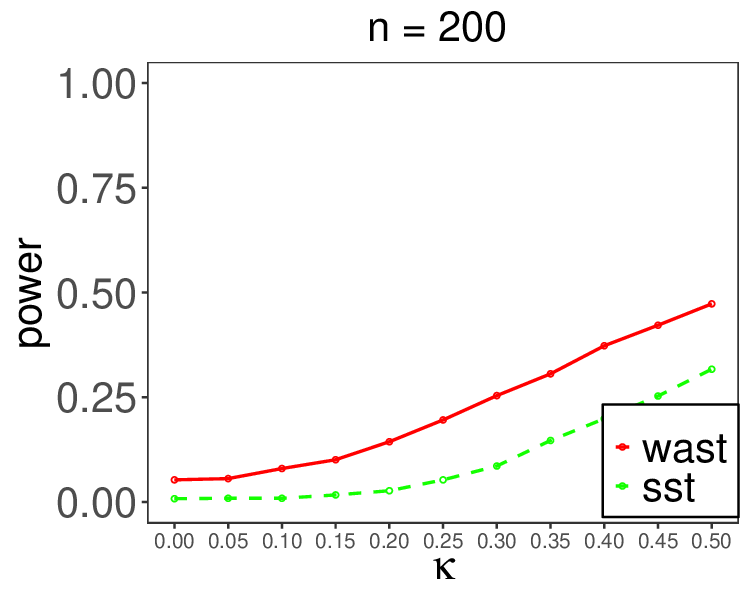}
		\includegraphics[scale=0.33]{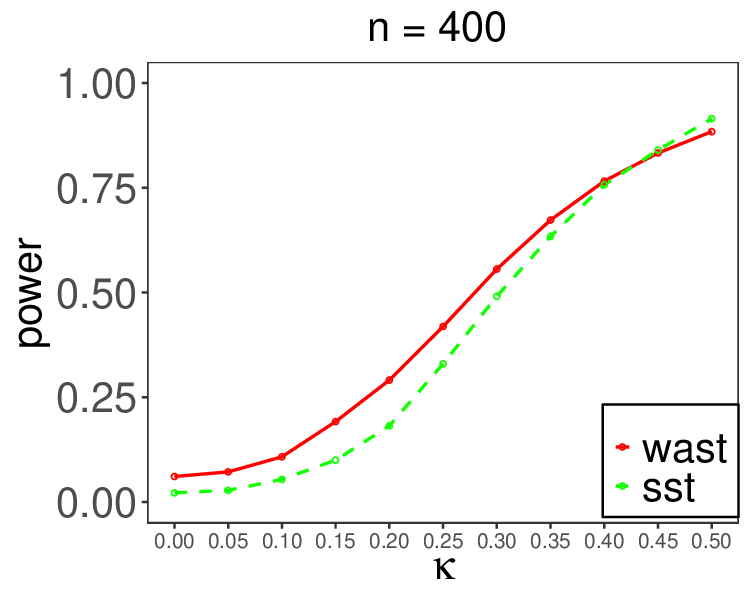}
		\includegraphics[scale=0.33]{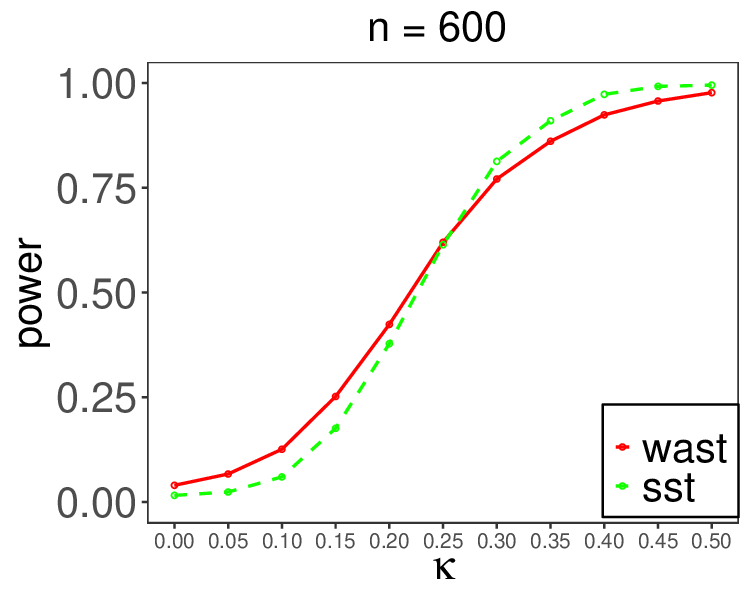}
		\caption{\it Powers of test statistic by the proposed WAST (red solid line) and SST (green dashed line) for $(p,q)=(5,5)$. From top to bottom, each row depicts the powers for quantile regression with $\tau=0.2$, $\tau=0.5$ and $\tau=0.7$, respectively. Here the $\Gv$ $Z$ is generated from $t_3$ distribution, and error from standard Cauchy distribution.}
		\label{fig_qr55_cauchy}
	\end{center}
\end{figure}

\begin{figure}[!ht]
	\begin{center}
		\includegraphics[scale=0.33]{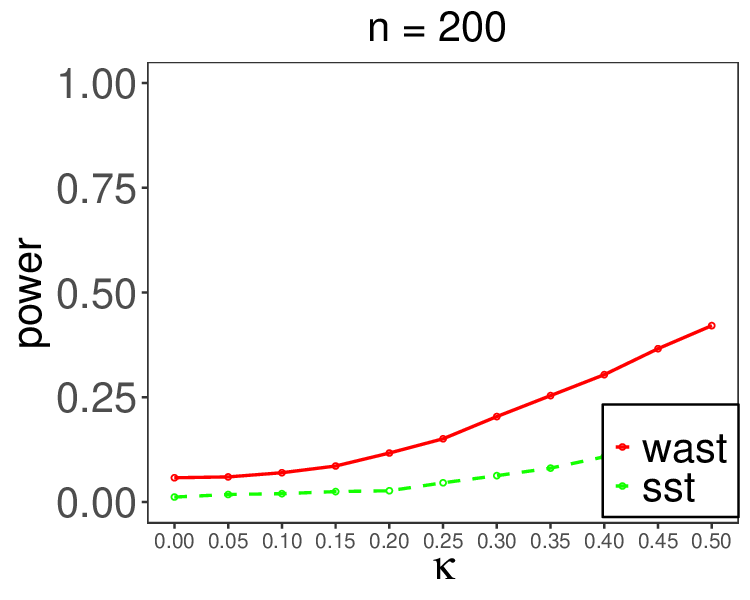}
		\includegraphics[scale=0.33]{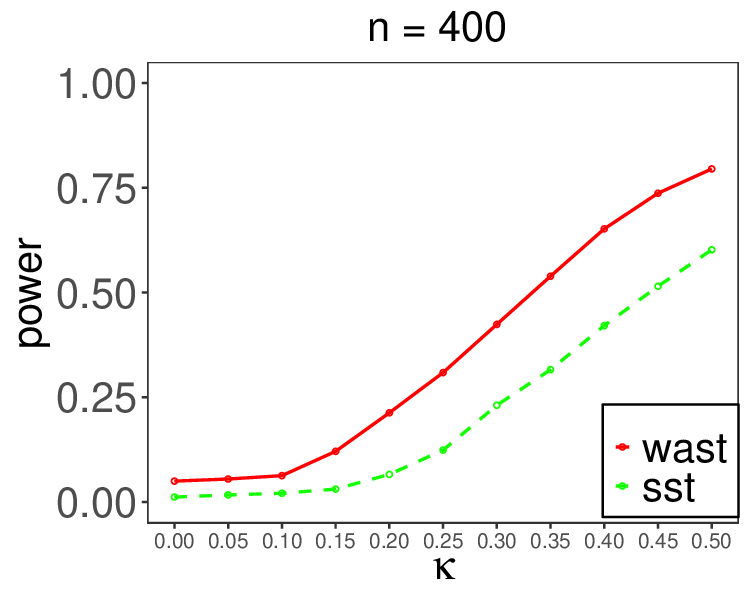}
		\includegraphics[scale=0.33]{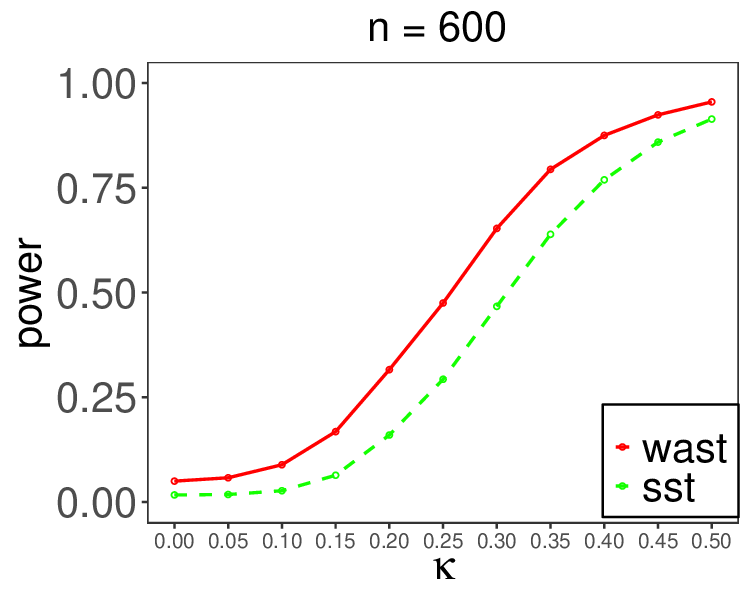}   \\
		\includegraphics[scale=0.33]{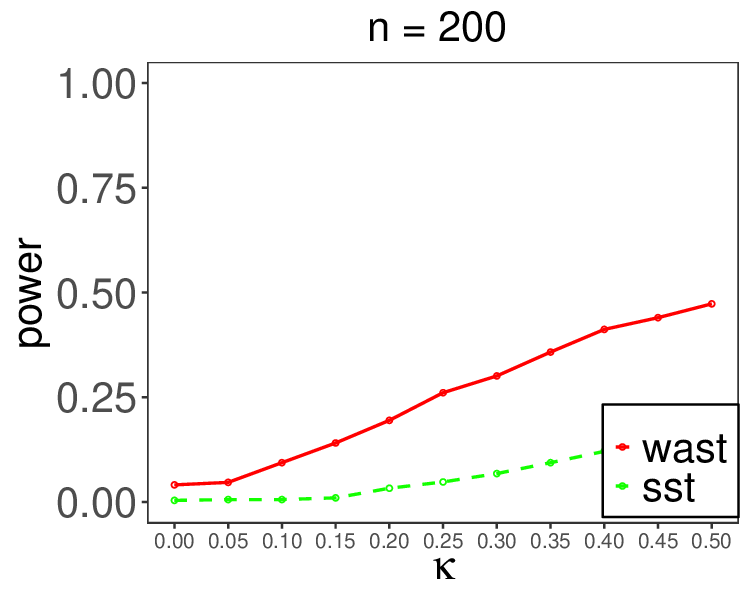}
		\includegraphics[scale=0.33]{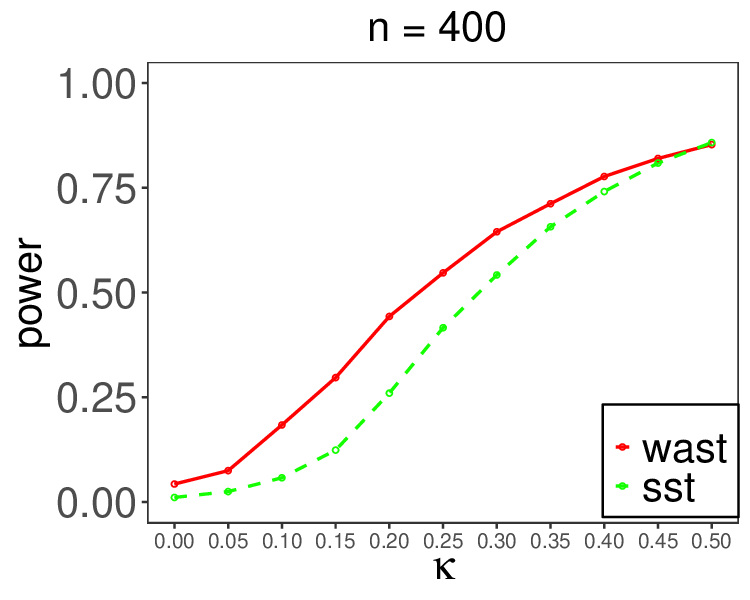}
		\includegraphics[scale=0.33]{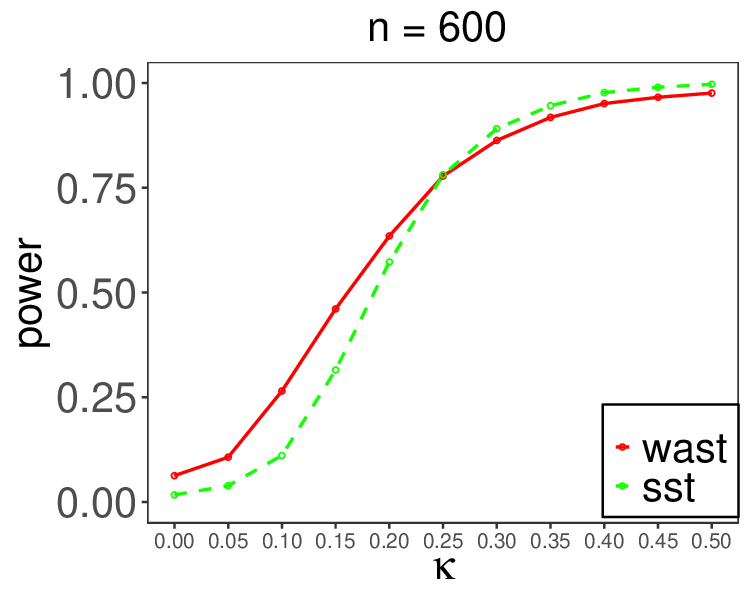}   \\
		\includegraphics[scale=0.33]{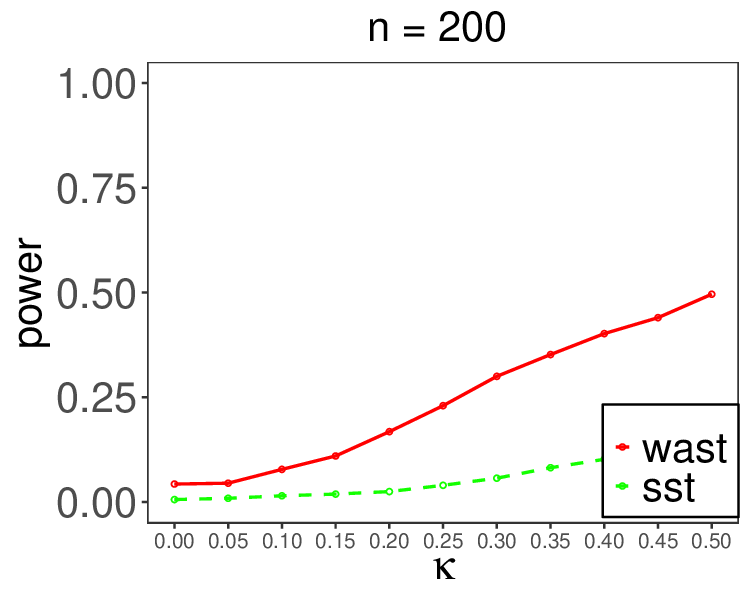}
		\includegraphics[scale=0.33]{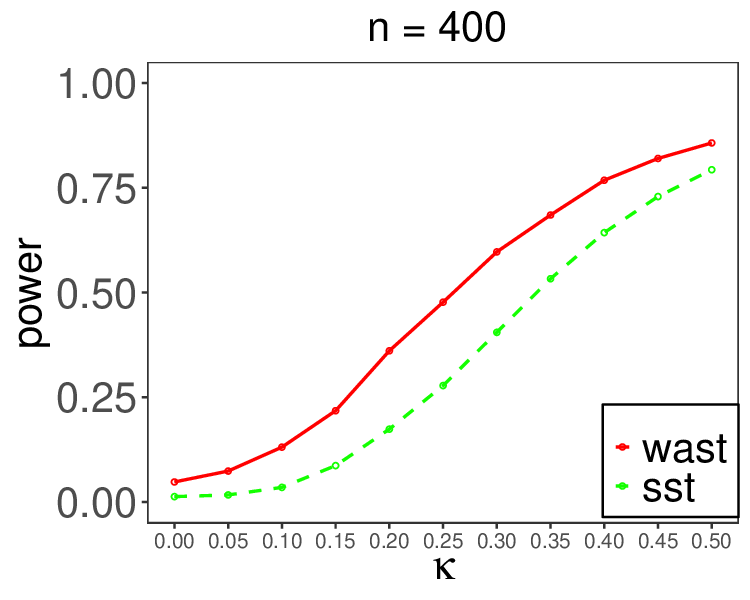}
		\includegraphics[scale=0.33]{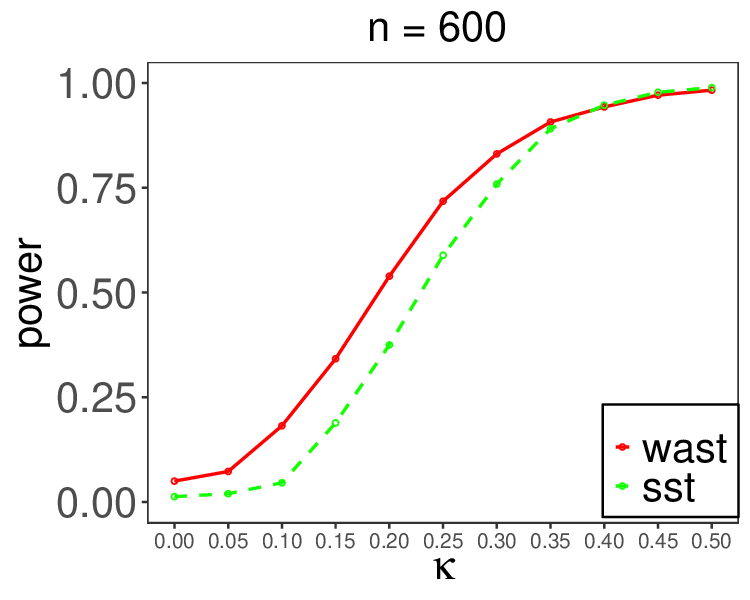}
		\caption{\it Powers of test statistic by the proposed WAST (red solid line) and SST (green dashed line) for $(p,q)=(10,10)$. From top to bottom, each row depicts the powers for quantile regression with $\tau=0.2$, $\tau=0.5$ and $\tau=0.7$, respectively. Here the $\Gv$ $Z$ is generated from $t_3$ distribution, and error from standard Cauchy distribution.}
		\label{fig_qr1010_cauchy}
	\end{center}
\end{figure}

\begin{figure}[!ht]
	\begin{center}
		\includegraphics[scale=0.33]{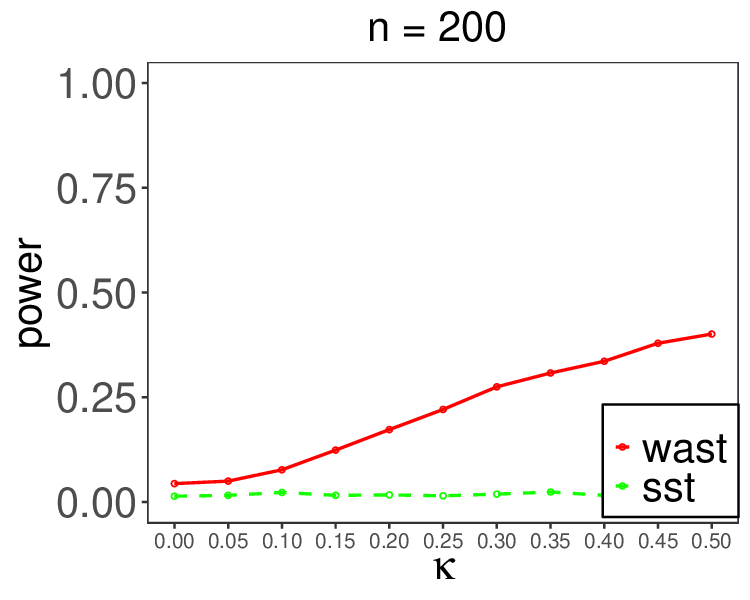}
		\includegraphics[scale=0.33]{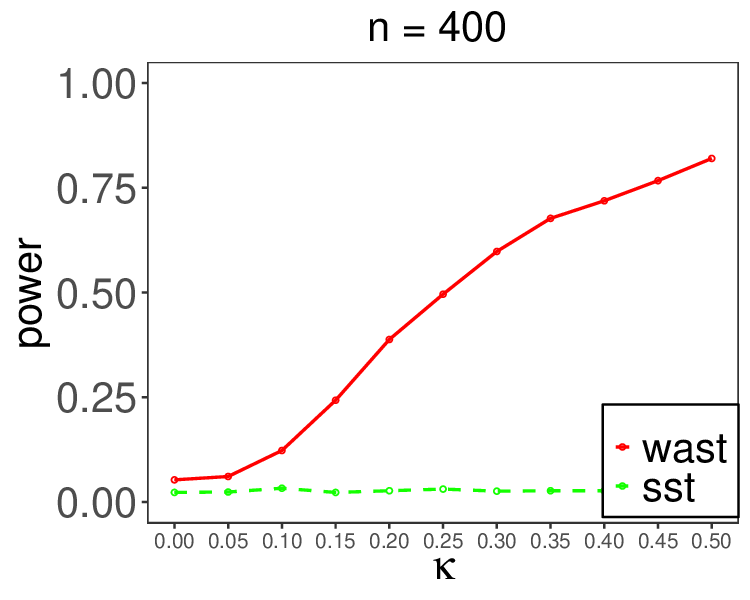}
		\includegraphics[scale=0.33]{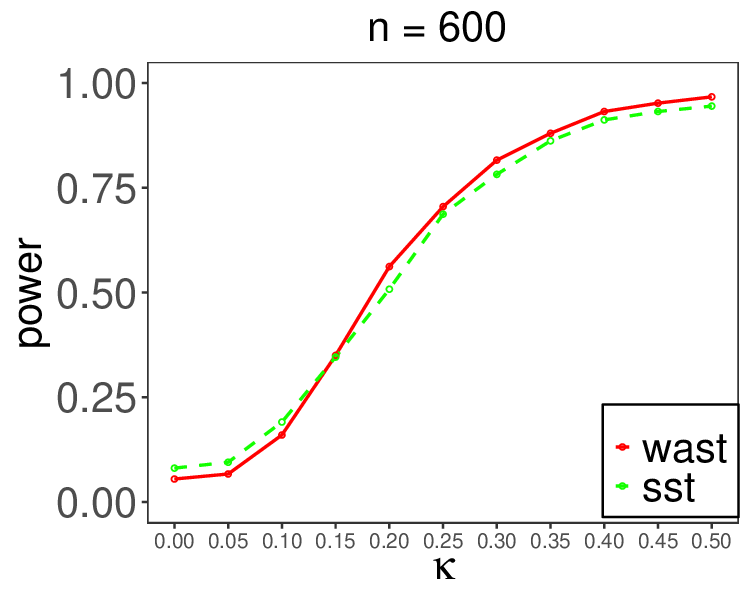}  \\
		\includegraphics[scale=0.33]{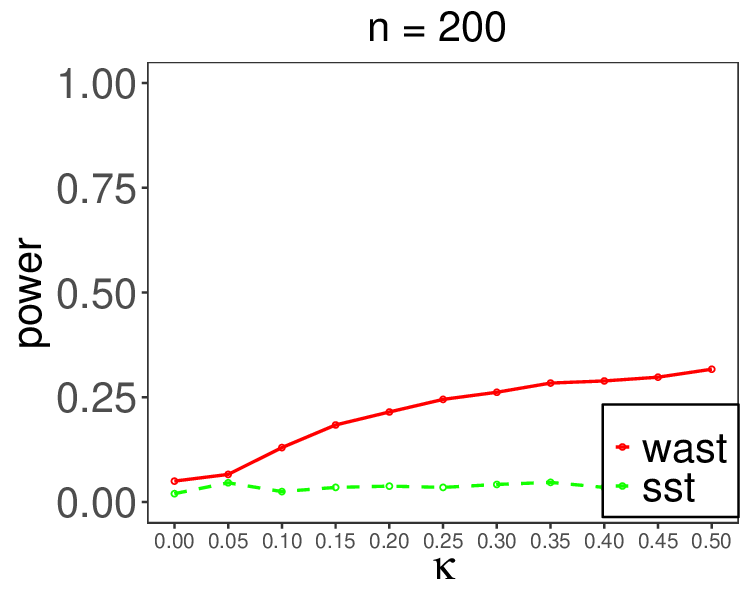}
		\includegraphics[scale=0.33]{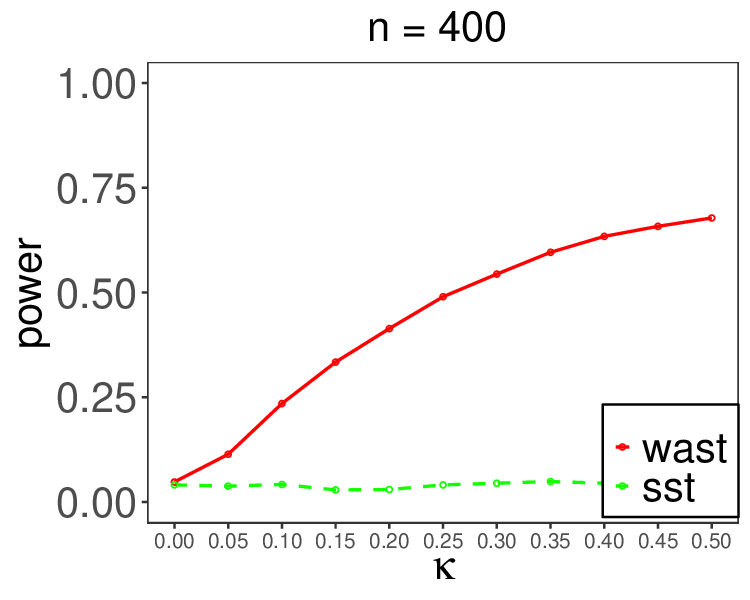}
		\includegraphics[scale=0.33]{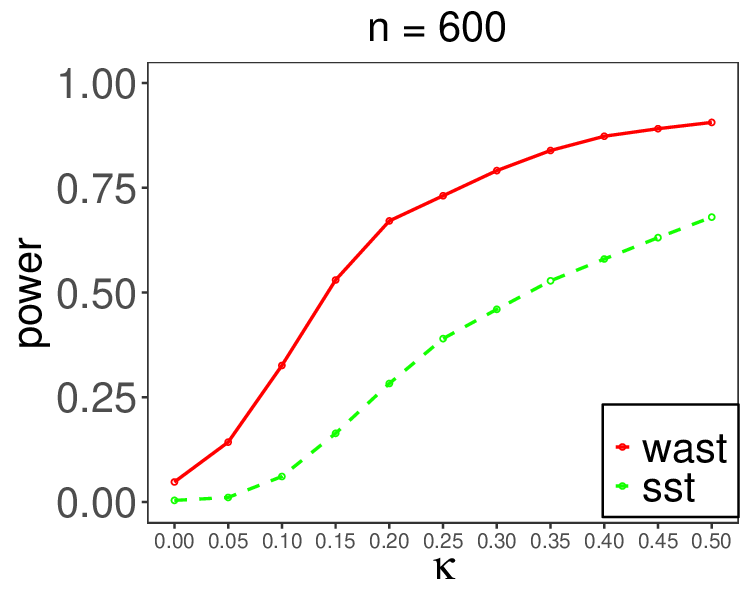}  \\
		\includegraphics[scale=0.33]{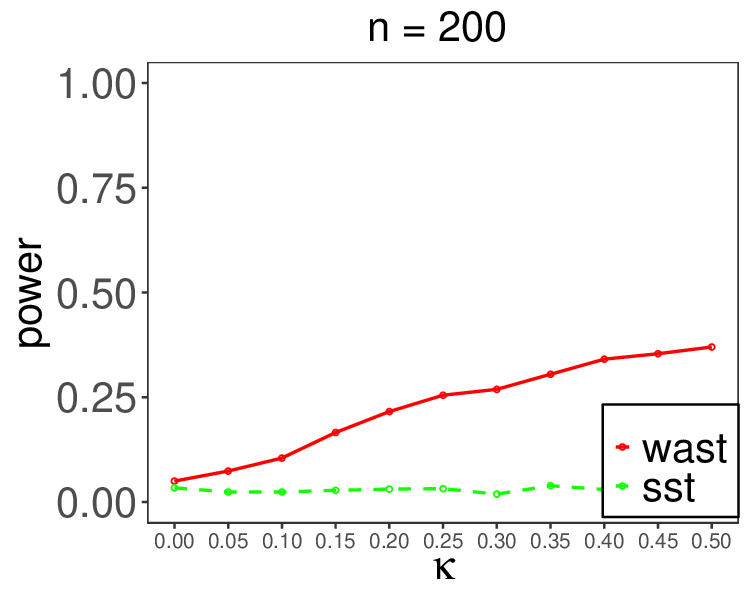}
		\includegraphics[scale=0.33]{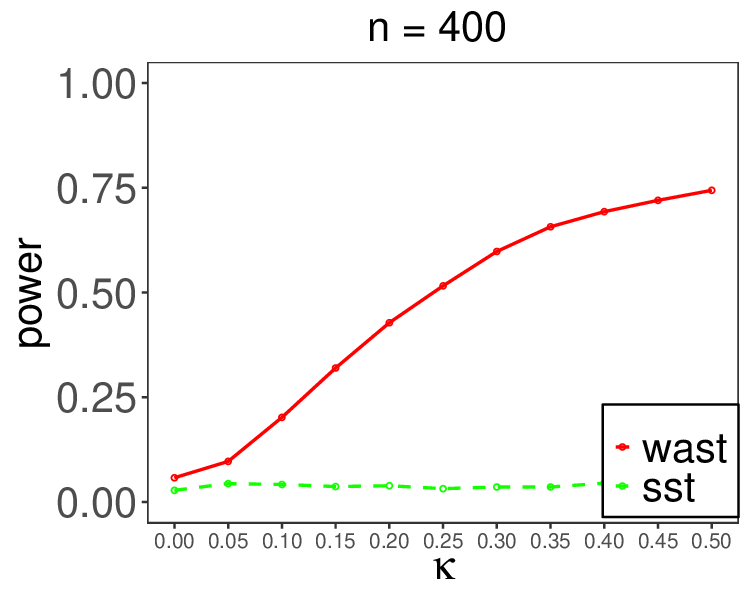}
		\includegraphics[scale=0.33]{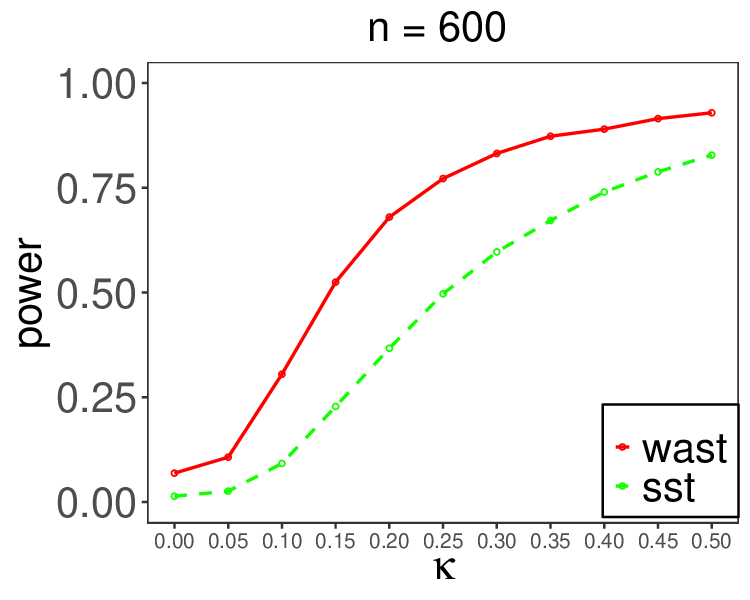}
		\caption{\it Powers of test statistic by the proposed WAST (red solid line) and SST (green dashed line) for $(p,q)=(50,20)$. From top to bottom, each row depicts the powers for quantile regression with $\tau=0.2$, $\tau=0.5$ and $\tau=0.7$, respectively. Here the $\Gv$ $Z$ is generated from $t_3$ distribution, and error from standard Cauchy distribution.}
		\label{fig_qr5020_cauchy}
	\end{center}
\end{figure}

Type \uppercase\expandafter{\romannumeral1} errors of the proposed test statistic are listed in Table \ref{table_size_cauchy}, and the power curves are depicted in Figure \ref{fig_qr13_cauchy}-\ref{fig_qr5020_cauchy}. We have the similar conclusion from Table \ref{table_size_cauchy} and Figure \ref{fig_qr13_cauchy}-\ref{fig_qr5020_cauchy}.

\begin{table}[htp!]
	\def~{\hphantom{0}}
    \tiny
	\caption{Type \uppercase\expandafter{\romannumeral1} errors of the proposed WAST 
    and SST based on resampling for quantile regression (QuantRE) with error from standard Cauchy distribution. 
    }
	\resizebox{\textwidth}{!}{
        \begin{threeparttable}
		\begin{tabular}{llcccccccc}
			\hline
			\multirow{2}{*}{Model}&\multirow{2}{*}{$(p,q)$}
			&\multicolumn{2}{c}{ $n=200$} && \multicolumn{2}{c}{ $n=400$} && \multicolumn{2}{c}{ $n=600$} \\
			\cline{3-4} \cline{6-7} \cline{9-10}
			&&   WAST & SST && WAST & SST && WAST & SST \\
			\cline{3-10}
			QuantRE &$(1,3)$          & 0.051 & 0.019 && 0.050 & 0.028 && 0.043 & 0.027 \\
			($\tau=0.2$)&$(5,5)$      & 0.038 & 0.011 && 0.051 & 0.013 && 0.042 & 0.010 \\
			& $(10,10)$               & 0.058 & 0.012 && 0.050 & 0.012 && 0.050 & 0.017 \\
			&$(50,20)$                & 0.044 & 0.014 && 0.053 & 0.023 && 0.055 & 0.081 \\
			[1 ex]
			QuantRE &$(1,3)$          & 0.055 & 0.033 && 0.059 & 0.039 && 0.056 & 0.049 \\
			($\tau=0.5$)&$(5,5)$      & 0.051 & 0.008 && 0.049 & 0.019 && 0.047 & 0.022 \\
			& $(10,10)$               & 0.041 & 0.004 && 0.043 & 0.011 && 0.063 & 0.017 \\
			&$(50,20)$                & 0.050 & 0.020 && 0.048 & 0.041 && 0.048 & 0.004 \\
			[1 ex]
			QuantRE &$(1,3)$          & 0.066 & 0.025 && 0.051 & 0.033 && 0.057 & 0.032 \\
			($\tau=0.7$)&$(5,5)$      & 0.053 & 0.008 && 0.061 & 0.022 && 0.040 & 0.016 \\
			& $(10,10)$               & 0.043 & 0.006 && 0.048 & 0.013 && 0.050 & 0.013 \\
			&$(50,20)$                & 0.050 & 0.034 && 0.058 & 0.028 && 0.069 & 0.014 \\
			\hline
		\end{tabular}
\begin{tablenotes}
\item The nominal significant level is 0.05. Here the $\Gv$ $\bZ$ is generated from $t_3$ distribution.
\end{tablenotes}
\end{threeparttable}
	}
	\label{table_size_cauchy}
\end{table}

\section{Case studies}\label{sec:case_study}
\setcounter{figure}{0}
\def\thefigure{D.\arabic{figure}}

\subsection{A histogram plot}
The histogram plot of the support incomes for the elderly in Shanxi Province of China is presented in this section. From this plot, we can observe that the support income has a heavy-tailed distribution.

\begin{figure}
	\centering
	\includegraphics[scale=1]{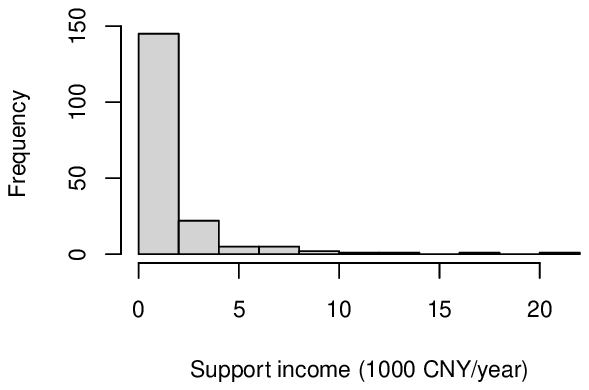}
	\caption{Histogram of support incomes in Shanxi Province of China.}
	\label{fig_hist}
\end{figure}

\subsection{Boston housing dataset}
We apply the proposed test to the Boston housing dataset from \citet{harrison1978hedonic}, which includes 14 variables that were measured across 506 census tracts in the Boston area. 
The response variable is the logarithm of the median value (MV) of the houses in those census tracts of the Boston Standard Metropolitan Statistical Area. 
The meaning of the other 13 social-economic, structural, and regional variables is as follows:
(1) Crime rate by town (CRIM);
(2) Proportion of residential land with a floor area of over 25,000 square feet (ZN);
(3) Proportion of nonretail business in each town (INDUS);
(4) Charles River dummy variable coding 1 if tract bounds the river and 0 otherwise (CHAS);
(5) Nitrogen oxide concentration in parts per 10 million (NOX);
(6) Average number of rooms per dwelling (RM);
(7) Proportion of self-occupied units built before 1940 (AGE);
(8) Weighted distances to five employment centers in Boston (DIS);
(9) Index of accessibility to radial highways (RAD);
(10) Full-value property-tax rate per \$10,000 (TAX);
(11) Pupil/teacher ratio by town (PTRATIO);
(12) $B = 1000( Bk - 0.63)^2$, where $Bk$ is the proportion of black people by town (B); and
(13) Percentage of people with lower status (LSTAT);

As explored in \citet{harrison1978hedonic}, the interest was to analyze the impact of air quality on the median housing costs. 
Thus, we take NOX as $\bX$ for our analysis. Furthermore, based on the regression tree depicted in Figure 1 of \citet{shannon2001tree},  LSTAT was identified as the root node with two descendent nodes CRIM and RM. Therefore, we take CRIM, RM and LSTAT as the grouping variables $\bZ$. 
The baseline covariate $\tbX$ includes seven variables: NOX, ZN, CHAS, DIS, RAD, TAX, and B. 
We consider the linear change-plane regression model, i.e.,
\begin{align*}
	Y_i = \tbX_i\trans\balpha + \bX_i\trans\bbeta\bone(\bZ_i\trans\btheta\geq0)+\eps_i, \quad i=1,\cdots,506,
\end{align*}
with $\tbX_i=(1,NOX_i,ZN_i,CHAS_i,DIS_i,RAD_i,TAX_i,B_i)\trans$, $\bZ_i=(1,CRIM_i, RM_i, LSTAT_i)\trans$, and $\bX_i=NOX_i$.

The p-value of the WAST is less than $0.0001$, which is calculated based on 5000 bootstrap samples. 
Therefore, this gives a strong evidence for rejecting the null hypothesis, that is, there exists a subgroup in which the air quality has different effects on the median value of the houses.

We, then, estimate parameters $(\balpha,\bbeta,\btheta)$ using a similar way as in the previous subsection. 
The estimators are $\hat{\balpha}=(3.463,-0.836,0.001,0.157,-0.025,0.004,-0.001,0.001)$, 
 $\hat{\bbeta}=1.002$, and $\hat{\btheta}=(-0.986,-0.009,0.166,-0.012).$
This gives one of subgroups, denoted by $G^+=\{i:\bZ^{T}_{i}\hat{\btheta}\geq 0\}$, $$G^+=\{i: -0.009*CRIM+0.166*RM-0.012*LSTAT\ge 0.986\}.$$ 
There are 132 census tracts in the subgroup $G^+$. From the estimation results, it is observed that the NOX variable has a positive influence on the MV of the houses in the subgroup $G^+$, while conversely, it exerts a negative impact on MV in another subgroup. 

\begin{figure}[!ht]
    \centering
    \subfloat{\includegraphics[scale=0.482]{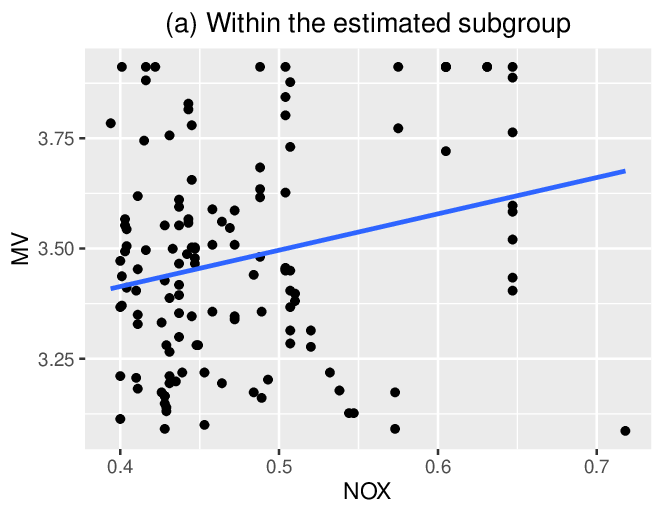}\label{fig:a}}
    \subfloat{\includegraphics[scale=0.482]{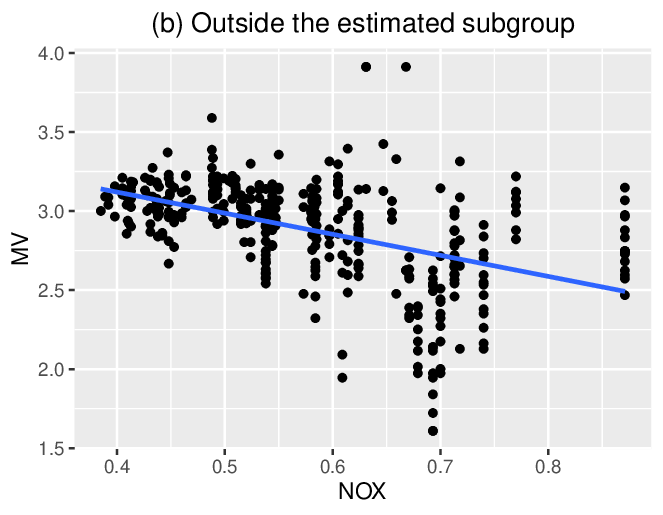}\label{fig:b}}
    \subfloat{\includegraphics[scale=0.482]{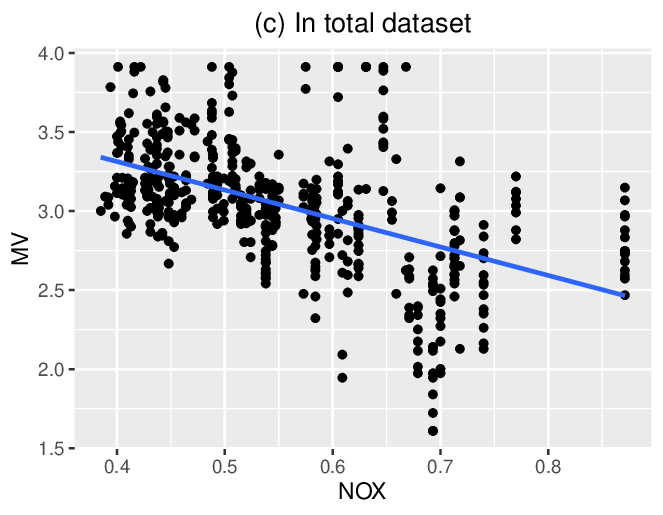}\label{fig:c}}
    \caption{Regression lines of the linear relationship of MV on NOX fitted on three different groups: (a) the subgroup $G^+$; (b) another subgroup $G^-=\{i:\bZ^{T}_{i}\hat{\btheta}< 0\}$; and (c) the total dataset.}\label{fig_cor}
\end{figure}

To demonstrate the impact of NOX on MV in different groups, we plot the regression lines in \autoref{fig_cor}. 
Based on \autoref{fig_cor} (a), it is evident that in the subgroup $G^+$, the median value of houses increases as the air quality decreases, while as shown in \autoref{fig_cor} (b), in another subgroup $G^-=\{i:\bZ^{T}_{i}\hat{\btheta}< 0\}$, the air quality has a positive impact on the median value of houses. This implies that in census tracts characterized by lower crime rates, higher room counts per dwelling, and smaller numbers of people with lower income, the median house value will still rise even in the presence of declining air quality. Thereby, combining \autoref{fig_cor} (a) and (b), the existence of significant  heterogeneity further supports our testing findings that there are subgroups.

\bibliographystyle{apalike}
\bibliography{papers}

\end{document}